\documentclass[11pt]{amsart}
\usepackage{amssymb}
\usepackage{mathabx}
\usepackage{xypic}
\usepackage{hhline}
\usepackage{dcpic}
\usepackage{hyperref}

\newtheorem{theorem}{Theorem}[section]
\newtheorem{lemma}[theorem]{Lemma}
\newtheorem{rem}[theorem]{Remark}
\newtheorem{hypothesis}[theorem]{Hypothesis}

\newtheorem{proposition}[theorem]{Proposition}
\newtheorem{corollary}[theorem]{Corollary}
\newtheorem{assumption}[theorem]{Assumption}

\newtheorem{definition}[theorem]{Definition}
\newtheorem{example}[theorem]{Example}

\theoremstyle{remark}

\numberwithin{theorem}{subsection}
\numberwithin{equation}{subsection}

\begin{document}


\makeatletter
\def\Ddots{\mathinner{\mkern1mu\raise\p@
\vbox{\kern7\p@\hbox{.}}\mkern2mu
\raise4\p@\hbox{.}\mkern2mu\raise7\p@\hbox{.}\mkern1mu}}
\makeatother

\newcommand{\OP}[1]{\operatorname{#1}}
\newcommand{\GO}{\OP{GO}}
\newcommand{\leftexp}[2]{{\vphantom{#2}}^{#1}{#2}}
\newcommand{\leftsub}[2]{{\vphantom{#2}}_{#1}{#2}}
\newcommand{\rightexp}[2]{{{#1}}^{#2}}
\newcommand{\rightsub}[2]{{{#1}}_{#2}}
\newcommand{\AI}{\OP{AI}}
\newcommand{\gen}{\OP{gen}}
\newcommand{\triv}{\OP{triv}}
\newcommand{\Asai}{\OP{Asai}}
\newcommand{\Image}{\OP{Im}}
\newcommand{\Res}{\OP{Res}}
\newcommand{\Ad}{\OP{Ad}}
\newcommand{\Out}{\OP{Out}}
\newcommand{\tr}{\OP{tr}}
\newcommand{\spec}{\OP{spec}}
\newcommand{\scopy}{\OP{end}}
\newcommand{\ord}{\OP{ord}}
\newcommand{\Cent}{\OP{Cent}}
\newcommand{\wellip}{\OP{w-ell}}
\newcommand{\Nrd}{\OP{Nrd}}
\newcommand{\Std}{\OP{Std}}
\newcommand{\JL}{\OP{JL}}
\newcommand{\temp}{\OP{temp}}
\newcommand{\BC}{\OP{BC}}
\newcommand{\sgn}{\OP{sgn}}
\newcommand{\SU}{\OP{SU}}
\newcommand{\Hom}{\OP{Hom}}
\newcommand{\Inter}{\OP{Int}}
\newcommand{\diag}{\OP{diag}}
\newcommand{\Sym}{\OP{Sym}}
\newcommand{\GSp}{\OP{GSp}}
\newcommand{\GL}{\OP{GL}}
\newcommand{\GSO}{\OP{GSO}}
\newcommand{\bdd}{\OP{bdd}}
\newcommand{\Int}{\OP{Int}}
\newcommand{\art}{\OP{art}}
\newcommand{\vol}{\OP{vol}}
\newcommand{\cusp}{\OP{cusp}}
\newcommand{\un}{\OP{un}}
\newcommand{\ellip}{\OP{ell}}
\newcommand{\sph}{\OP{sph}}
\newcommand{\gsimp}{\OP{sim-gen}}
\newcommand{\Aut}{\OP{Aut}}
\newcommand{\disc}{\OP{disc}}
\newcommand{\sdisc}{\OP{s-disc}}
\newcommand{\aut}{\OP{aut}}
\newcommand{\End}{\OP{End}}
\newcommand{\barQ}{\OP{\overline{\mathbf{Q}}}}
\newcommand{\barQp}{\OP{\overline{\mathbf{Q}}_{\it p}}}
\newcommand{\Gal}{\OP{Gal}}
\newcommand{\simp}{\OP{sim}}
\newcommand{\odd}{\OP{odd}}
\newcommand{\Normal}{\OP{Norm}}
\newcommand{\Ind}{\OP{Ind}}
\newcommand{\St}{\OP{St}}
\newcommand{\unit}{\OP{unit}}
\newcommand{\reg}{\OP{reg}}
\newcommand{\SL}{\OP{SL}}
\newcommand{\Frob}{\OP{Frob}}
\newcommand{\Id}{\OP{Id}}
\newcommand{\GSpin}{\OP{GSpin}}
\newcommand{\Norm}{\OP{Norm}}

\title[Endoscopic Classification ...]{Endoscopic Classification of representations of Quasi-Split Unitary Groups}

\author{Chung Pang Mok}

\address{Hamilton Hall, McMaster University, Hamilton, Ontario, L8S 4K1}

\email{cpmok@math.mcmaster.ca}

\begin{abstract}
In this paper we establish the endoscopic classification of tempered representations of quasi-split unitary groups over local fields, and the endoscopic classification of the discrete automorphic spectrum of quasi-split unitary groups over global number fields. The method is analogous to the work of Arthur \cite{A1} on orthogonal and symplectic groups, based on the theory of endoscopy and the comparison of trace formulas on unitary groups and general linear groups.
\end{abstract}

\maketitle

\tableofcontents

{\section{\textbf{Introduction}}

In the work \cite{A1} Arthur established, for symplectic and orthogonal groups, the endoscopic classification of tempered irreducible admissible representations over local fields, namely the construction of packets of representations, and the endoscopic classification of the discrete automorphic spectrum of these groups over global number fields. The method is based on the theory of endoscopy in both the standard and twisted case, together with the comparison of trace formula on the classical groups involved and the twisted trace formula for general linear groups (modulo the stabilization of the twisted trace formula for general linear groups, which is work in progress of Waldspurger and others, {\it c.f.} \cite{W7,W8,W9,MW3}). In this paper, we use Arthur's methods to establish analogous results for quasi-split unitary groups, generalizing the results of Rogawski \cite{R}. 

As in \cite{A1} the method of proof is a long induction argument that establish both the local and global theorems simultaneously; the induction argument for the proof of the local and global classification is completed at the end of the paper. We refer the reader to section two below for the description of the main classification results in the local and global setting. 

\bigskip

\noindent The content of the paper is as follows:

In section two we give the formalism of parameters and state the main classification theorems in both the local and global setting. In stating the global theorems we follow Arthur \cite{A1} and use a formal version of global parameter, defined in terms of conjugate self-dual cuspidal automorphic representations on general linear groups, to avoid the reference to the conjectural automorphic Langlands group. At the end of section two we also give a review of previous results obtained by earlier authors.  

Section three states the local character identities that characterizes the representations in a packet. We also give the statement of the local intertwining relation. The local intertwining relation is the crucial ingredient for the analysis of the spectral terms of trace formula, and in reducing the construction of local packets of representations to the case where the packets correspond to square integrable parameters. 

Section four gives some preliminary comparison of the trace formulas of unitary groups and general linear groups, which in particular gives the existence of ``weak base change" for discrete automorphic representations on unitary groups. The comparison of trace formulas in section four forms the background for the more sophisticated comparsion in section five and six.

In section five we state the stable multiplicity formula, one of the main global results, and which when combined with the global intertwining relation, allows a term by term comparison in the spectral and endoscopic expansion of the discrete part of the trace formula. The proof of the global intertwining relation (which follows from the corresponding local result) is completed in section eight, while the stable multiplicity formula is completed only in section nine of the paper. However, in section five we will be able to establish them for a large class of ``degenerate" parameters, based on the induction hypothesis. 

In section six we extend the analysis of section five to the parameters that are square-integrable, but satisfy rather stringent local constraints at certain archimedean places. The fact that we are able to establish the stable multiplicity formula for these parameters will form the input to the proof of the local classification theorems to be carried out in section seven.

In section seven we construct the packets associated to generic local parameters, and obtain the local Langlands classification for tempered representations of quasi-split unitary groups. Among the technical results to be established is the local intertwining relation, which reduces the constructions of packets of tempered representations to the case of discrete series representations. The method of proof is global, drawing on the results from trace formulas comparison in section 5 and 6. This is based on the standard technique of embedding a discrete series representation of a local group as a local component of an automorphic representation.

Section eight constructs the packets associated to general parameters. The arguments are similar to that of section seven, combined with results from section seven and the duality operator of Aubert-Schneider-Stuhler. The induction arguments concerning the local theorems are finished in section seven.    

In the final section nine we complete the induction arguments of all the global theorems, and obtains the endoscopic classification of automorphic representations for quasi-split unitary groups. The completion of the global induction argument follows that of \cite{A1}, by considering auxiliary parameters. As a corollary, when combined with results on global automorphic descent of Ginzburg-Rallis-Soudry, we obtain the local and global generic packet conjecture.

\bigskip

Finally we follow the convention of \cite{A1}: {\it Theorems} stated in the paper are understood to be proved by the long induction argument only to be completed at the end of the paper. On the other hand, results that we prove along the argument will either be stated as {\it Propositions} or {\it Lemmas}.

\section*{Acknowledgement}
The author would like to thank Professor Arthur, besides his encouragements, for showing how to use the trace formula. He would also like to thank the participants of the trace formula working seminar: Florian Herzig, Zhifeng Peng, Shuichiro Takeda, Kam Fai Tam, Fucheng Tan, Patrick Walls, Bin Xu. The discussions during the seminars have been very helpful throughout the work. He is very grateful to Professor Moeglin for comments and pointing out inaccuracies in the initial draft of the work, and also to the referee for numerous helpful suggestions and corrections. The author is partially supported by the NSERC Discovery Grant, and the Early Researcher Award from the Ontario Ministry of Research and Innovation.

\newpage

\section*{Notation}

Throughout the paper we denote by $F$ a local or global field, depending on the context (always of characteristic zero). And throughout we fix an algebraic closure $\overline{F}$ of $F$. When $F$ is a global field, we generally denote by $v$ a place of $F$, and denote by $F_v$ the completion of $F$ at $v$, and we fix embeddings $\overline{F} \hookrightarrow \overline{F}_v$. In both the local and global context we denote by $\Gamma_F$ the absolute Galois group of $F$ (or just $\Gamma$ when the context of $F$ is clear), and by $W_F$ the Weil group of $F$. 

Throughout the paper, when we refer to a quadratic field extension $E$ of $F$, we always intend $E$ as a subfield of $\overline{F}$, or equivalently a field extension $E$ of $F$ when a specified embedding into $\overline{F}$. We fix the identification:
\[
E \otimes_F \overline{F} = \overline{F} \times \overline{F}
\]  
such that the projection to the left $\overline{F}$ factor correspond to the specified embedding of $E$ into $\overline{F}$.

Similarly, when $F$ and hence $E$ are global, then for each place $v$ of $F$, the composite of the embeddings $E \hookrightarrow \overline{F} \hookrightarrow \overline{F}_v$ specifies the corresponding embedding of $E$ into $\overline{F}_v$, which in the case where $v$ splits in $E$ singles out the particular place $w$ of $E$ above $v$, with $\overline{w}$ being the other place of $E$ above $v$. In this case where $v$ splits in $F$ we fix the identification:
\[
E \otimes _F F_v = E_w \times E_{\overline{w}}=F_v \times F_v
\]
where the projection to the left $F_v$ factor corresponds to the place $w$ of $E$ above $v$.

For $E$ a quadratic extension of $F$ as above, and $N \geq 1$, we denote by $U_{E/F}(N)$ the quasi-split unitary group in $N$ variables over $F$, whose group of $F$-points is given by
\begin{eqnarray*}
U_{E/F}(N)(F)=\{ g \in \GL_N(E) \,\ | \,\  \leftexp{t}{c(g)} J g =J    \}
\end{eqnarray*}
here $J$ is the anti-diagonal matrix 
\begin{equation*}
J=\left( \begin{array}{rrr}
& & 1 \\
&  \Ddots & \\
1 & &
\end{array} \right)
\end{equation*}
and $z \mapsto c(z)$ is the Galois conjugation of $E/F$. When the context of $E/F$ is clear we will just write $U(N)$ for $U_{E/F}(N)$.

We identify the centre of $U_{E/F}(N)$ as $U_{E/F}(1)$ (consisting of scalar matrices).

When discussing unitary groups, it will be convenient to include the case where $E$ is a split quadratic separable extension of $F$ (as an $F$-algebra), i.e. $E = F \times F$, with the conjugation of $E$ over $F$ being given by the interchange of the two factors. We then have $U_{E/F}(N) \cong \GL_{N/F}$. More precisely, we denote by $\iota_1$ (resp. $\iota_2$) the isomorphism $\iota_i: U_{E/F}(N) \stackrel{\simeq}{\rightarrow} \GL_{N/F}$ induced by the projection of $E$ to the left $F$ factor (resp. the right $F$ factor). Then the map 
\[
\GL_{N/F} \stackrel{\iota_1^{-1}}{\rightarrow} U_{E/F}(N) \stackrel{\iota_2}{\rightarrow} \GL_{N/F}
\]
is given by $g \mapsto J \leftexp{t}{g}^{-1} J^{-1}$.

We denote by $G_{E/F}(N)$ the algebraic group over $F$ given by
\[
G_{E/F}(N) = \Res_{E/F} \GL_{N/E}
\] 
where $\Res_{E/F}$ is Weil restriction of scalars of $E/F$ (this includes both the case where $E$ is a field or a split $F$-algebra). Thus the group of $F$-points of $G_{E/F}(N)$ is given by $\GL_N(E)$. When the context of $E/F$ is clear the group $G_{E/F}(N)$ is abbreviated as $G(N)$. We denote by $\theta$ the $F$-automorphism of $G_{E/F}(N)$, whose action on $F$-points is given by 
\begin{eqnarray}
& & \\
& & \theta(g) = \Phi_N \leftexp{t}{c(g)}^{-1} \Phi_N^{-1} \mbox{ for }
g \in G_{E/F}(N)(F)=\GL_N(E). \nonumber
\end{eqnarray}
here $\Phi_N$ is the anti-diagonal matrix with alternating $\pm 1$ entries:
\begin{equation}
\Phi_N = \left( \begin{array}{rrrr}
& & & 1 \\
& & -1 & \\
&  \Ddots & & \\
(-1)^{N-1} & & &
\end{array} \right).
\end{equation}
Then $\theta$ preserves the standard $F$-splitting of $G_{E/F}(N)$.

\section{\textbf{Statement of the main theorems}}

\subsection{$L$-groups and $L$-embeddings}

We first recall the description of the dual groups and $L$-groups involved. As usual if $G$ is a connected reductive algebraic group over $F$, we denote by $\widehat{G}$ the dual group of $G$ (as a complex algebraic group). For the $L$-group we always use the Weil form:
\[
\leftexp{L}{G} = \widehat{G} \rtimes W_F
\]
with the action of $W_F$ on $\widehat{G}$ factors through $\Gal(F^{\prime}/F)$, where $F^{\prime}$ is any finite Galois extension of $F$ over which $G$ splits.

For $G=U_{E/F}(N)$ we have 
\[
\widehat{U}_{E/F}(N) = \GL_N(\mathbf{C}).
\]
We denote by $\alpha$ the automorphism of $\widehat{U}_{E/F}(N) = \GL_N(\mathbf{C})$ given by
\begin{eqnarray}
\alpha(g) = \Phi_N \leftexp{t}{g}^{-1} \Phi_N^{-1}, \,\ g \in \GL_N(\mathbf{C})
\end{eqnarray}
with $\Phi_N$ as in (1.0.2).

\noindent The matrix $\Phi_N$ satisfies
\[
\leftexp{t}{\Phi_N} = (-1)^{N-1}\Phi_N, \,\ \Phi_N^2 = (-1)^{N-1}.
\]
Thus $\alpha$ is of order two. As usual $\alpha$ is the unique automorphism in its inner class that fixes the standard splitting of $\GL_N(\mathbf{C})$.

We have
\[
\leftexp{L}{U_{E/F}(N)} = \GL_N(\mathbf{C}) \rtimes W_F
\]
with the action of $W_F$ on $\GL_N(\mathbf{C})$ factors through $\Gal(E/F)$, and if $w_c \in W_F \smallsetminus W_E$, then $w_c$ acts as the automorphism $\alpha$.

For our purpose it is important to include the case where $E=F \times F$ is the split quadratic separable extension of $F$ as an $F$-algebra. In this case we set 
\[
\widehat{U}_{E/F}(N) = \{(g,\leftexp{t}{g}^{-1}), \,\ g \in \GL_N(\mathbf{C})  \} \subset \GL_N(\mathbf{C}) \times \GL_N(\mathbf{C}).
\]
The projection of $E$ to the left $F$-factor (resp. to the right $F$-factor) corresponds to the isomorphism $\widehat{\iota}_1: \widehat{U}_{E/F}(N) \stackrel{\simeq}{\rightarrow} \GL_N(\mathbf{C})$ given by the projection to the left $\GL_N(\mathbf{C})$ factor (resp. the isomorphism $\widehat{\iota}_2: \widehat{U}_{E/F}(N) \stackrel{\simeq}{\rightarrow} \GL_N(\mathbf{C})$ given by the projection to the right $\GL_N(\mathbf{C})$ factor). We also put in this case:
\[
\leftexp{L}{U_{E/F}(N)} = \widehat{U}_{E/F}(N) \times W_F
\]
(in other words, the $L$-action is trivial). Then $\widehat{\iota}_1$ (resp. $\widehat{\iota}_2$) induce the isomorphism $\leftexp{L}{\iota}_1 : \leftexp{L}{U_{E/F}(N)} \stackrel{\simeq}{\rightarrow} \leftexp{L}{\GL_{N/F}}$ (resp. $\leftexp{L}{\iota}_2: \leftexp{L}{U_{E/F}(N)} \stackrel{\simeq}{\rightarrow} \leftexp{L}{\GL_{N/F}}$).

The main example is the case where $F$ is a global field, and $E$ is a quadratic field extension of $F$, and $v$ a prime of $F$. If $v$ does not split in $E$, then we have (still denoting the prime of $E$ above $v$ as $v$):
\[
\widehat{U}_{E_v/F_v}(N) = \widehat{U}_{E/F}(N) = \GL_N(\mathbf{C})
\]
and the fixed embedding $\overline{F} \hookrightarrow \overline{F}_v$ induces $W_{F_v} \hookrightarrow W_F$, and hence the embedding $\leftexp{L}U_{E_v/F_v}(N) \hookrightarrow \leftexp{L}{U_{E/F}(N)}$. On the other hand, suppose that $v$ splits in $E$ into two primes $w,\overline{w}$. Then we have $E_v = E_w \times E_{\overline{w}}$, and we may identify $E_w=E_{\overline{w}}=F_v$. In this case we rename the maps $\widehat{\iota}_1, \leftexp{L}{\iota}_1$ (resp. $\widehat{\iota}_2,\leftexp{L}{\iota}_2$) as $\widehat{\iota}_w, \leftexp{L}{\iota}_w$ (resp. $\widehat{\iota}_{\overline{w}}, \leftexp{L}{\iota}_{\overline{w}}$). The fixed embedding $\overline{F} \hookrightarrow \overline{F}_v$ (which also correspond to embedding $W_{F_v} \hookrightarrow W_F$. ) singles out the prime between $w,\overline{w}$, which we assume without loss of generality to be $w$. We then define the corresponding embedding of $L$-groups:
\begin{eqnarray*}
& & \leftexp{L}{U_{E_v/F_v}(N)} \rightarrow \leftexp{L}{U_{E/F}(N)}\\
& & h \times \sigma \mapsto    \widehat{\iota}_{w} (h)    \times   \sigma, \,\ h \in \widehat{U}_{E_v/F_v}(N)
\end{eqnarray*}
(note that the image of any element of $W_{F_v}$ in $\Gal(E/F)$ is trivial).

\bigskip

For the group $G_{E/F}(N) = \Res_{E/F} \GL_{N/E}$, we have 
\[
\widehat{G}_{E/F}(N) = \GL_N(\mathbf{C}) \times \GL_N(\mathbf{C}),
\]
and
\[
\leftexp{L}{G_{E/F}(N)} = (\GL_N(\mathbf{C}) \times \GL_N(\mathbf{C})) \rtimes W_F
\] 
with the action of $W_F$ on $\GL_N(\mathbf{C}) \times \GL_N(\mathbf{C})$ factors through $\Gal(E/F)$, and if $w_c \in W_F \smallsetminus W_E$ as before, then $w_c$ acts on $\GL_N(\mathbf{C}) \times \GL_N(\mathbf{C})$ as the automorphsim $\beta$ that interchanges the two factors, i.e.
\[
\beta((g,h)) =(h,g) \mbox{ for } (g,h) \in \GL_N(\mathbf{C}) \times \GL_N(\mathbf{C}).
\]

Finally we denote by $\widehat{\theta}$ the following automorphism of $\widehat{G}_{E/F}(N)$:
\begin{eqnarray}
\widehat{\theta}((g,h)) &=& (\alpha(h),\alpha(g)) \\
&=& (\Phi_N \leftexp{t}{h}^{-1} \Phi_N^{-1},\Phi_N\leftexp{t}{g}^{-1}\Phi_N^{-1}) \nonumber
\end{eqnarray}
for $(g,h) \in \widehat{G}_{E/F}(N) = \GL_N(\mathbf{C}) \times \GL_N(\mathbf{C})$. The inner class of the automorphism $\widehat{\theta}$ on $\widehat{G}_{E/F}(N)$ is dual to that of $\theta$ on $G_{E/F}(N)$, and is the unique automorphism in its inner class that preserves the standard $\Gamma_F$-splitting of $\widehat{G}_{E/F}(N)$.

\bigskip

We record some embeddings of $L$-groups which will be used throughout the paper. 

Thus let $E$ be a quadratic field extension of $F$ as before. We denote by $c: z \mapsto c(z)$ the Galois conjugation of $E$ over $F$. In the local (resp. global) setting, denote by $\omega_{E/F}$ the quadratic character of $F^{\times}$ (resp. $\mathbf{A}_F^{\times}/F^{\times}$) that corresponds to the quadratic extension $E/F$ under local (resp. global) class field theory. We put
\begin{eqnarray}
\mathcal{Z}_E = \{ \chi : E^{\times} \rightarrow \mathbf{C}^{\times} \mbox{ unitary}, \,\ \chi \circ c  = \chi^{-1}   \}
\end{eqnarray} 
in the local case, and similarly
\begin{eqnarray}
\mathcal{Z}_E = \{ \chi :  \mathbf{A}_E^{\times}/E^{\times} \rightarrow \mathbf{C}^{\times} \mbox{ unitary}, \,\ \chi \circ c  = \chi^{-1}   \}
\end{eqnarray} 
in the global case. In other words $\mathcal{Z}_E$ is the set of conjugate self-dual characters of $E^{\times}$ in the local setting, and of $ \mathbf{A}_E^{\times}/E^{\times}$ in the global setting. We then have a partition:
\begin{eqnarray*}
\mathcal{Z}_E = \mathcal{Z}_E^+ \coprod \mathcal{Z}_E^-
\end{eqnarray*} 
where
\begin{eqnarray}
& & \\
& & \mathcal{Z}_E^+ = \{\chi \in \mathcal{Z}_E, \,\ \chi|_{F^{\times}} =1 \}, \,\ \mathcal{Z}_E^- = \{\chi \in \mathcal{Z}_E, \,\ \chi|_{F^{\times}} =\omega_{E/F} \} \nonumber
\end{eqnarray}
in the local case, and similarly 
\begin{eqnarray}
& & \\
& & \mathcal{Z}_E^+ = \{\chi \in \mathcal{Z}_E, \,\ \chi|_{\mathbf{A}_F^{\times}} =1 \}, \,\ \mathcal{Z}_E^- = \{\chi \in \mathcal{Z}_E, \,\ \chi|_{\mathbf{A}_F^{\times}} =\omega_{E/F} \} \nonumber
\end{eqnarray}
in the global case (note that the condition $\chi|_{\mathbf{A}_F^{\times}} =1$ or $\omega_{E/F}$ implies the conjugate self-duality condition $\chi \circ c = \chi^{-1}$; similarly in the local case). For the definitions of the $L$-embeddings we need to choose characters $\chi_{\kappa} \in \mathcal{Z}_E^{\kappa}$ for $\kappa = \pm1$. For $\kappa=+1$ it is of course natural to just take $\chi_+ =1$, but to allow flexibility in induction arguments we will work in the more general setting.

We often identify such a character $\chi_{\kappa} \in \mathcal{Z}_E^{\kappa}$ for $\kappa = \pm 1$ as a character on $W_E$ under local or global class field theory. Then $\chi_{\kappa}$ satisfies the following: for $w_c \in W_F \smallsetminus W_E$, we have 
\begin{eqnarray}
\chi_{\kappa}(w_c \sigma w_c^{-1}) = \chi_{\kappa}(\sigma)^{-1} \mbox{ for } \sigma \in W_E 
\end{eqnarray}
\begin{eqnarray}
\chi_{\kappa}(w_c^2) = \kappa.
\end{eqnarray}

Given a sign $\kappa = \pm 1$, and $\chi_{\kappa} \in \mathcal{Z}_E^{\kappa}$, define the following embedding of $L$-groups:
\begin{eqnarray*}
\xi_{\chi_{\kappa}} : \leftexp{L}{U_{E/F}(N)} & \longrightarrow & \leftexp{L}{G_{E/F}(N)} 
\end{eqnarray*}
given by the rule: for a fixed choice of $w_c \in W_F \smallsetminus W_E$, 
\begin{eqnarray}
& & g \rtimes 1 \mapsto  (g, \leftexp{t}{g}^{-1}) \rtimes 1   \mbox{ for }  g \in \GL_N(\mathbf{C})\\ 
& & I_N \rtimes \sigma \mapsto (\chi_{\kappa}(\sigma) I_N, \chi^{-1}_{\kappa}(\sigma) I_N) \rtimes \sigma \mbox{ for } \sigma \in W_E \nonumber \\
& & I_N \rtimes w_c \mapsto (\kappa \Phi_N, \Phi_N^{-1}) \rtimes w_c \nonumber
\end{eqnarray}

\noindent The $\widehat{G}_{E/F}(N)$-conjugacy class of $\xi_{\chi_{\kappa}}$ is independent of the choice of $w_c \in W_F \smallsetminus W_E$. In general we only regard the $L$-embedding $\xi_{\chi_{\kappa}}$ only up to $\widehat{G}_{E/F}(N)$-conjugacy (we will be explicit when we do the otherwise). 

When $\kappa=+1$, and $\chi_+=1$, $L$-embedding $\xi_1$ is independent of the choice of $w_c$, and is usually referred to as the {\it standard base change $L$-embedding}. On the other hand when $\kappa =-1$, there is no canonical choice of $\chi_-$. In the literature the $L$-embedding $\xi_{\chi_{\kappa}}$ for $\kappa=-1$ is usually referred to as a {\it twisted base change $L$-embedding}. We will refrain from using such terminology whenever possible to avoid confusion. 

Slightly more generally, let $N_1, \cdots N_r$ be non-negative integers such that 
\[
N_1 + \cdots N_r =N.
\] 
Given $\kappa_i = \pm 1$ for each $i=1,\cdots,r$,  put $\underline{\kappa} = (\kappa_1,\cdots,\kappa_r)$. Given $\chi_i \in \mathcal{Z}_E^{\kappa_i}$ for $i=1,\cdots,r$, put $\underline{\chi}=(\chi_1,\cdots,\chi_r)$. We refer to $\underline{\kappa}$ as the {\it signature} of $\underline{\chi}$. Define the $L$-embedding 
\begin{eqnarray}
\xi_{\underline{\chi}} : \leftexp{L}{(U(N_1) \times \cdots U(N_r)   )} \longrightarrow \leftexp{L}{G(N)}
\end{eqnarray}
by composing the product of the $L$-embeddings:
\[
\prod_{i=1}^r \xi_{\chi_{\kappa_i}} : \leftexp{L}{(U(N_1) \times \cdots U(N_r)   )} \rightarrow \leftexp{L}{(  G(N_1) \times \cdots G(N_r)  )}
\]
with the obvious diagonal $L$-embedding
\begin{eqnarray}
 \leftexp{L}{(  G(N_1) \times \cdots G(N_r)  )} \rightarrow \leftexp{L}{G(N)}.
\end{eqnarray}

\bigskip

Next we consider $L$-embedding between the $L$-groups of unitary groups. As above given a partition $N_1 + \cdots + N_r = N$, put 
\begin{eqnarray}
& & \kappa_i =(-1)^{N-N_i} \\
& & \underline{\kappa} =(\kappa_1,\cdots,\kappa_r) \nonumber
\end{eqnarray}

Given $\underline{\chi}=(\chi_1,\cdots,\chi_r)$ with signature $\underline{\kappa}$, define the $L$-embedding (for a fixed choice of $w_c \in W_F \smallsetminus W_E$ as above):
\[
\zeta_{\underline{\chi}}: \leftexp{L}{(U(N_1) \times \cdots \times U(N_r)  )} \rightarrow \leftexp{L}{U(N)}
\]
given by the rule: 
\begin{eqnarray}
& & \\
& & (g_1,\cdots,g_r) \rtimes 1 \mapsto \diag(g_1,\cdots,g_r) \rtimes 1, \,\ g_i \in \GL_{N_i}(\mathbf{C}) \nonumber \\
 & & (I_{N_1}, \cdots,I_{N_r}) \rtimes \sigma \mapsto \diag(\chi_{\kappa_1}(\sigma) I_{N_1}, \cdots, \chi_{\kappa_r}(\sigma) I_{N_r} ) \rtimes \sigma \mbox{ for } \sigma \in W_E \nonumber \\
& & (I_{N_1}, \cdots I_{N_r}) \rtimes w_c \mapsto \diag(\kappa_1 \Phi_{N_1},\cdots, \kappa_r \Phi_{N_r}) \cdot \Phi_N^{-1}\rtimes w_c \nonumber
\end{eqnarray}
The $\widehat{U}(N)$-conjugacy class of the $L$-embedding $\leftexp{L}{(U(N_1) \times \cdots \times U(N_r)  )} \rightarrow \leftexp{L}{U(N)}$ is independent of the choice of $w_c$. 

It is immediate to check that we have the commutative diagram: for $\underline{\chi}$ as above and $\chi^{\prime} \in \mathcal{Z}_E^{\kappa^{\prime}}$, we have:
\begin{eqnarray}
\xymatrix{ \leftexp{L}{(U(N_1) \times \cdots \times U(N_r))} \ar[r]^-{\zeta_{\underline{\chi}}} \ar[rd]_{\zeta_{\widetilde{\underline{\chi}}}} & \leftexp{L}{U(N)} \ar[d]^{\xi_{\chi^{\prime}}} \\ & \leftexp{L}{G(N)} }
\end{eqnarray}
with $\widetilde{\underline{\chi} }=(\chi^{\prime} \chi_1,\cdots,\chi^{\prime}\chi_r)$. Thus the signature of $\widetilde{\underline{\chi}}$ is given by $\widetilde{\underline{\kappa}}=(\kappa^{\prime} \kappa_1,\cdots,\kappa^{\prime} \kappa_r  )$ (recall that $\kappa_i$ is defined as in (2.1.12)).

\subsection{Formalism of local parameters}
We first recall the formalism of local $L$-parameters. Thus $F$ is local. Denote by $L_F$ the local Langlands group of $F$. Thus
\begin{eqnarray}
 L_F = \left \{ \begin{array}{c}  W_F \mbox{ if } F \mbox{ is archimedean} \\   \,\  W_F \times \SU(2) \mbox{ if } F \mbox{ is non-archimedean}
\end{array} \right.
\end{eqnarray} 

In general, if $G$ is a connected reductive group over $F$, then an $L$-parameter for $G(F)$ is an admissible homomorphism:
\[
\phi:L_F \longrightarrow \leftexp{L}{G}.
\]
Two $L$-parameters are equivalent if they are conjugate by $\widehat{G}$. The set of equivalence classes of $L$-parameters of $G(F)$ is noted as $\Phi(G(F))$, or when the context of $F$ is clear, noted as $\Phi(G)$. 

We also define the subset $\Phi_{\bdd}(G) \subset \Phi(G)$ of bounded parameters, i.e. consisting of those parameters $\phi \in \Phi(G)$ whose image in $\widehat{G}$ is bounded.

For our purpose we will mostly be concerned with the case where $G=U_{E/F}(N)$, or when $G=G_{E/F}(N)$ for a quadratic extension $E/F$ (in fact we will need the twisted version of $G_{E/F}(N)$ as will be discussed in section 3). 

In the case where $G=G_{E/F}(N)=G(N)$, we recall that ({\it c.f.} \cite{R} p. 48) there is a natural bijection between the set $\Phi(G(N)(F))$ of equivalence classes of $L$-parameters of $G(N)(F)$, and the set $\Phi(\GL_N(E))$ of equivalence classes of $L$-parameters of $\GL_n(E)$. Here an $L$-parameter of $\Phi(\GL_N(E))$ is identified as an (equivalence class of) $N$-dimensional admissible representation of $L_E$.

To define the bijection fix $w_c \in W_F \smallsetminus W_E$ as before (the resulting bijection on equivalence classes is independent of the choice of $w_c$). Then for $\phi \in \Phi(\GL_n(E))$, the $L$-parameter $\phi^{\prime} \in \Phi(G(N)(F))$ corresponding to $\phi$ is defined as:

\begin{eqnarray}
& & \phi^{\prime} : L_F \longrightarrow \leftexp{L}{G(N)} \\
& & \phi^{\prime}(\sigma) = (\phi(\sigma),\phi^{c}(\sigma)) \rtimes \sigma \mbox{ for } \sigma \in L_E \nonumber \\
& & \phi^{\prime}(w_c) = (\phi(w_c^2), I_N) \rtimes w_c \nonumber
\end{eqnarray}
here 
\begin{eqnarray} \phi^c(\sigma) :=\phi(w_c^{-1} \sigma w_c).
\end{eqnarray}
We will henceforth identify the two sets $\Phi(\GL_N(E))$ and $\Phi(G(N)(F))$.

\bigskip

For $\kappa=\pm1$ and $\chi_{\kappa} \in \mathcal{Z}_E^{\kappa}$, we have maps of $L$-parameters 
\begin{eqnarray*}
\xi_{\chi_{\kappa},*}: \Phi(U(N)) &\rightarrow & \Phi(G(N)) \\
\phi &\mapsto & \xi_{\chi_{\kappa}} \circ \phi.
\end{eqnarray*}
If $\phi \in \Phi(U(N))$, then the $L$-parameter in $\Phi(\GL_N(E))$ that corresponds to $\xi_{\chi_{\kappa}} \circ \phi \in \Phi(G(N))$ is just $\phi|_{L_E} \otimes \chi_{\kappa}$. In particular if $\kappa=+1$ and if we choose $\chi_+=1$, then the corresponding $L$-parameter of $\Phi(\GL_N(E))$ is just given by $\phi|_{L_E}$ (this is usually known as the standard base change of $L$-parameters). 

The maps $\xi_{\chi_{\kappa}}$ give an injection of the set of (equivalence classes of) $L$-parameters of $U_{E/F}(N)(F)$ to $G_{E/F}(N)(F)$ (this can be proved as in \cite{GGP} theorem 8.1, part ii, or see lemma 2.2.1 below).

\bigskip

The image of $\Phi(U(N))$ in $\Phi(G(N)) \cong \Phi(\GL_N(E))$ under $\xi_{\chi_{\kappa},*}$ can be characterized as follows. In general, if 
\[
\rho:L_E \rightarrow \GL(V)
\]
is an admissible representation (with $V \cong \mathbf{C}^N$), then $\rho$ is called {\it conjugate self-dual} if $\rho^c \cong \rho^{\vee}$, where as in (2.2.3)
\[
\rho^c(\sigma) := \rho(w_c^{-1} \sigma w_c) \mbox{ for } \sigma \in L_E.
\]
and $\rho^{\vee}$ is the contragredient of $\rho$. In addition, $\rho$ is called {\it conjugate self-dual with parity $\eta$}, for $\eta = \pm1$, if there exists a non-degenerate bilinear form $B(\cdot,\cdot)$ on $V$, satisfying, for all $x,y \in V$:
\begin{eqnarray}
 B(\rho^c(\sigma)x, \rho(\sigma)y) =B(x,y) 
\end{eqnarray}
\begin{eqnarray}
 B(x,y) = \eta \cdot B(y, \rho(w_c^2)x)
\end{eqnarray}
This notion is independent of the choice of $w_c \in W_F \smallsetminus W_E$.

More concretely, if we represent $\rho$ as a homomorphism $\rho: L_E \rightarrow \GL_N(\mathbf{C})$, then the condition for $\rho$ being conjugate self-dual of parity $\eta$ is that there exists an $A \in \GL_N(\mathbf{C})$ such that
\begin{eqnarray}
& &  \leftexp{t}{\rho^c(\sigma)} A \rho(\sigma) =A \mbox{ for } \sigma \in L_E \\
& & \leftexp{t}{A} = \eta \cdot A \cdot  \rho(w_c^2) \nonumber
\end{eqnarray}

Following the terminology of section 3 of \cite{GGP}, a conjugate self-dual representation of parity +1 (resp. -1) will be called a conjugate orthogonal (resp. conjugate symplectic) representation.

\begin{lemma}
Let $\chi_{\kappa} \in \mathcal{Z}_E^{\kappa}$. The image of $\xi_{\chi_{\kappa},*} : \Phi(U_{E/F}(N)) \rightarrow \Phi(G_{E/F}(N)) \cong \Phi(\GL_N(E))$ is given by the set of parameters in $\Phi(\GL_N(E))$ that are conjugate self-dual with parity $\eta$, with $\eta = (-1)^{N-1} \cdot \kappa$.
\end{lemma}
\begin{proof}
This is given for instance in theorem 8.1 of \cite{GGP}. For the convenience of the reader we include some details here. Thus let
\[
\phi: L_F \longrightarrow \leftexp{L}{U(N)}
\]
be an element of $\Phi(U(N))$. As usual we identify $\phi|_{L_E}$ as a representation $\rho:L_E \rightarrow \GL_N(\mathbf{C})$. Write
\[
\phi(w_c) = C \rtimes w_c.
\]
Set $A:= \Phi_N C^{-1}$. Then by direct computation (using $\leftexp{t}{\Phi_N} = (-1)^{N-1} \Phi_N$) we have
\begin{eqnarray*}
\rho(w_c^2) \rtimes w_c^2=\phi(w_c^2) &=& \phi(w_c)^2\\ & =&  A^{-1} \cdot (-1)^{N-1} \cdot \leftexp{t}{A} \rtimes w_c^2
\end{eqnarray*}
i.e.
\begin{eqnarray}
\leftexp{t}{A} = (-1)^{N-1} A \cdot \rho(w_c^2).
\end{eqnarray}
Similarly, using the identity
\[
\phi(w_c^{-1} \sigma w_c) = \phi(w_c)^{-1} \phi(\sigma) \phi(w_c)
\]
we obtain
\begin{eqnarray}
 \leftexp{t}{\rho^{c}(\sigma)} A \rho(\sigma)=A.
\end{eqnarray}
Now the element of $\Phi(\GL_N(E))$ that corresponds to $\xi_{\chi_{\kappa},*} \circ \phi$ is given by $\rho^{\prime} :=\rho \otimes \chi_{\kappa} = \phi|_{L_E} \otimes \chi_{\kappa}$. By (2.1.8) $\chi_{\kappa}(w_c^2) = \kappa$. So from (2.2.7) we have
\[
\leftexp{t}{A} = (-1)^{N-1} \kappa A \cdot \rho^{\prime}(w_c^2).
\]
On the other hand, from (2.1.7) and (2.2.8) we obtain
\begin{eqnarray*}
\leftexp{t}{(\rho^{\prime})^c(\sigma) } A \rho^{\prime}(\sigma) = A \mbox{ for } \sigma \in L_E.
\end{eqnarray*}
Hence the assertion that $\rho^{\prime}$ is conjugate self-dual of parity $(-1)^{N-1} \kappa$. 

The other direction of the assertion is similar. Thus let $\rho^{\prime} \in \Phi(\GL_N(E))$ be an admissible representation $\rho^{\prime}:L_E \rightarrow \GL_N(\mathbf{C})$ that is conjuagte self-dual of parity $\eta= (-1)^{N-1} \kappa$, as in the situation of (2.2.6). Let $\phi^{\prime} \in \Phi(G(N))$ be the $L$-parameter of $G(N)$ corresponding to $\rho^{\prime}$, as given by (2.2.2). Thus
\begin{eqnarray*}
& & \phi^{\prime} : L_F \longrightarrow \leftexp{L}{G(N)} \\
& & \phi^{\prime}(\sigma) =(\rho^{\prime}(\sigma),(\rho^{\prime})^{c}(\sigma)) \rtimes \sigma \mbox{ for } \sigma \in L_E \\
& & \phi^{\prime}(w_c) =(\rho^{\prime}(w_c^2),I_N) \rtimes w_c.
\end{eqnarray*}  
Put 
\[
\phi^{\prime \prime} :=(I_N,\leftexp{t}{A}) \cdot \phi^{\prime} \cdot  (I_N,\leftexp{t}{A}^{-1}).
\]
Then $\phi^{\prime \prime}$ is $\widehat{G}(N)$-conjugate to $\phi^{\prime}$, and by (2.2.6) (with $\rho$ being replaced by $\rho^{\prime}$), we have
\begin{eqnarray*}
& & \phi^{\prime \prime} (\sigma) =(\rho^{\prime}(\sigma),\leftexp{t}{\rho^{\prime}(\sigma)}^{-1}) \rtimes \sigma \mbox{ for } \sigma \in L_E \\
& & \phi^{\prime \prime}(w_c) = ( \rho^{\prime}(w_c^2)  \leftexp{t}{A}^{-1}, \leftexp{t}{A}  )\rtimes w_c = (\eta  A^{-1}, \leftexp{t}{A}) \rtimes w_c.
\end{eqnarray*}
Hence if we put $C:=(-1)^{N-1} A^{-1} \Phi_N^{-1}=A^{-1} \Phi_N$, then 
\[
(\eta A^{-1}, \leftexp{t}{A}) = (\kappa C \Phi_N, \leftexp{t}{C}^{-1} \Phi_N^{-1}) = (C,\leftexp{t}{C}^{-1}) \cdot (\kappa \Phi_N,\Phi_N^{-1}).
\]
Thus $\phi^{\prime \prime} = \xi_{\chi_{\kappa},*} \circ \phi$, where $\phi \in \Phi(U(N))$ is given by:
\begin{eqnarray*}
& & \phi(\sigma) = (\rho^{\prime} \otimes \chi_{\kappa}^{-1})(\sigma) \rtimes \sigma \mbox{ for } \sigma \in L_E\\
& & \phi(w_c) = C \rtimes w_c.
\end{eqnarray*}
\end{proof}

As in section 7 of \cite{GGP}, a more invariant way to formulate the notion of conjugate self-duality with parity is as follows. Let $V$ be a finite dimensional vector space over $\mathbf{C}$, and denote by $\Std$ the standard representation of $\GL(V)$ acting on $V$. Put:
\begin{eqnarray*}
& & H^0 := \GL(V) \times \GL(V) \times W_E \\
& & H := ( \GL(V) \times \GL(V) ) \rtimes W_F
\end{eqnarray*}
here $W_F$ acts on $ \GL(V) \times \GL(V)$ through $\Gal(E/F)$, with the non-trivial element of $\Gal(E/F)$ acts by permuting the two $\GL(V)$ factors. If we make the identification $V = \mathbf{C}^N$ then we have $H=\leftexp{L}{G(N)}$. Now as on p.28 of \cite{GGP} there is a decomposition
\begin{eqnarray*}
\Ind^H_{H^0} (\Std \boxtimes \Std \boxtimes \triv) = \Asai^+  \oplus \Asai^- 
\end{eqnarray*}
with $\Asai^+ $ and $\Asai^- $ being representations of $H$ of dimension equal to $(\dim V)^2$. They are distinguished by:
\[
\tr(w_c|\Asai^+ )= \dim V, \,\ \tr(w_c|\Asai^-   )= -\dim V.
\]
Given $\rho:L_E \rightarrow \GL(V)$ be an admissible representation, we view $\rho$ as an admissible homomorphism $\widetilde{\rho}: L_F \rightarrow H$. Then we put
\[
\Asai^+ \rho := \Asai^+ \circ \widetilde{\rho},\,\ \Asai^- \rho := \Asai^- \circ \widetilde{\rho}
\]
which are admissible representations of $L_F$ of dimension equal to $(\dim V)^2$. Then $\rho$ is conjugate self-dual if and only if $\rho^c \otimes \rho$ contains a non-gegenerate vector (in the sense that it corresponds to a non-degenerate bilinear form on the underlying space of $\rho$) that is fixed under the action of $L_E$. Similarly $\rho$ is conjugate orthogonal (resp. conjugate symplectic) if and only if $\Asai^+\rho$ (resp. $\Asai^-\rho$) contains a non-degenerate vector that is fixed under the action of $L_F$.

\bigskip

\begin{rem}
\end{rem}
Even if $\rho$ is a conjugate self-dual representation of $L_E$ with a parity, its parity need not be unique. If $\rho$ is irreducible then the parity is unique, as clear from the Schur's lemma and characterization given above in terms of the Asai representation (2.2.9). More generally if $\rho$ is of the form:
\begin{eqnarray*}
\rho = \rho_1 \oplus \cdots \oplus \rho_r
\end{eqnarray*}
with the $\rho_i$ being mutually non-isomorphic, irreducible conjugate self-dual representation with the same parity $\eta$, then $\rho$ is also conjugate self-dual with parity $\eta$, and the parity of $\rho$ is also unique.
\bigskip

Following Arthur, we introduce the set of local $A$-parameters, which plays the role of local components of the global classification. In general if $G$ is connected reductive group over $F$, then a local $A$-parameter is an admissible homomorphism:
\begin{eqnarray}
\psi: L_F \times \SU(2) \longrightarrow \leftexp{L}{G}
\end{eqnarray}
such that $\psi|_{L_F}$ has bounded image in $\widehat{G}$ (since $\SU(2)$ is compact this is equivalent to saying that $\psi$ has bounded image). As in the case of $L$-parameters, the equivalence condition on $A$-parameters is defined by $\widehat{G}$-conjugacy, and the set of equivalence classes of $A$-parameters is noted as $\Psi(G(F))$ or $\Psi(G)$. As in \cite{A1}, we need to introduce a larger set $\Psi^+(G)$, whose elements are (equivalence classes of) all the admissible homomorphsim as above, without the boundedness condition on the restriction to $L_F$. 

An $A$-parameter $\psi \in \Psi(G)$ (or $\Psi^+(G)$) is called {\it generic} if $\psi$ is trivial on the $\SU(2)$-factor of $L_F \times \SU(2)$, in which case it can be identified as an element of $\Phi(G)$. As in the case of $L$-parameters, we will mostly be concerned with the case where $G= U_{E/F}(N)$ or $G_{E/F}(N)$ (and also the twisted version of $G_{E/F}(N)$). All the discussions above for the $L$-parameters of these groups, in particular the statement of lemma 2.2.1, can be formulated verbatim for $A$-parameters.

\bigskip

Recall the following (\cite{A1} p. 24): given an $A$-parameter $\psi \in \Psi^+(G)$, we can extend $\psi$ by analytic continuation:
\[
\psi: L_F \times \SL_2(\mathbf{C}) \longrightarrow \leftexp{L}{G}.
\]
then we define the $L$-parameter $\phi_{\psi} \in \Phi(G)$ associated to $\psi$ by the formula:
\begin{eqnarray}
\phi_{\psi}(\sigma) = \psi\Big(\sigma,   \begin{pmatrix}   |\sigma|^{1/2}  &  0 \\    0  &     |\sigma|^{-1/2}     \\ \end{pmatrix}  \Big) \mbox{ for } \sigma \in L_{F}.
\end{eqnarray}

\bigskip

Finally as in \cite{A1},  we introduce the following groups, which play a key role in the local classification:
\begin{eqnarray*}
& & S_{\psi} = \Cent(\Image \psi, \widehat{G}) \\
& & \overline{S}_{\psi} = S_{\psi}/Z(\widehat{G})^{\Gamma_F}\\
& & \mathcal{S}_{\psi} = \pi_0(\overline{S}_{\psi}).
\end{eqnarray*}

We also define the element:
\begin{eqnarray}
s_{\psi} := \psi\Big(1,   \begin{pmatrix}   -1  &  0 \\    0  &    -1     \\ \end{pmatrix}  \Big)
\end{eqnarray}
then $s_{\psi}$ is a central semi-simple element of $S_{\psi}$. The significance of the element $s_{\psi}$ will be clear from the local character relations for the representations in an $A$-packet ({\it c.f.} theorem 3.2.1).

\subsection{Formal global parameters}

In this subsection we define the global parameters. Thus $F$ now denotes a global number field, and $E$ a quadratic extension of $F$. In the ideal situation the global parameters can be defined as in the local case, with $L_F$ being the conjectural automorphic Langlands group. 

In the absence of the conjectural automorphic Langlands group, we can still define, following Arthur \cite{A1}, the global parameters for $U_{E/F}(N)$ in a formal way, in terms of cuspidal automorphic representation of general linear groups. This will be done in the next subsection, and this subsection is a preliminary to this.

To begin with we define:
\[
\Psi_{\simp}(G_{E/F}(N)) = \Psi_{\simp}(N)
\]
to be the set of formal tensor products:
\[
\mu \boxtimes \nu
\]
where $\mu$ is a unitary cuspidal automorphic representation on $\GL_m(\mathbf{A}_E)$, and $\nu$ is an algebraic representation of $\SL_2(\mathbf{C})$ of dimension $n$, such that $N = m \cdot n$. By the theorem of Moeglin-Waldspurger \cite{MW}, such an element $\psi^N = \mu \boxtimes \nu \in \Psi_{\simp}(N)$ corresponds to an irreducible unitary representation $\pi_{\psi^N}$ of $\GL_N(\mathbf{A}_E)$ that belongs to the discrete automorphic spectrum, namely $\pi_{\psi^N}$ is given by the isobaric sum:
\begin{eqnarray}
& & \\
& &
(\mu \otimes | \det|_{\mathbf{A}_E}^{\frac{n-1}{2}}) \boxplus ( \mu \otimes | \det|_{\mathbf{A}_E}^{\frac{n-3}{2}} )  \boxplus \cdots \boxplus (\mu \otimes | \det|_{\mathbf{A}_E}^{\frac{-(n-1)}{2}}). \nonumber
\end{eqnarray}
We also identify in the evident way the automorphic representations of $\GL_N(\mathbf{A}_E)$ and $G_{E/F}(N)(\mathbf{A}_F)$.

Define the following operation on the set of (unitary) cuspidal automorphic representation on $\GL_m(\mathbf{A}_E)$: given $\mu$ as above, put:
\begin{eqnarray}
\mu^{*} := (\mu^c)^{\vee}
\end{eqnarray}
where $\mu^c:= \mu \circ c$, with $c$ being the Galois conjugation of $E$ over $F$, and $(\mu^c)^{\vee}$ is the contragredient of $\mu^c$. Note that the operation $\mu \mapsto \mu^*$ is induced by the outer automorphism $g \mapsto \leftexp{t}{c(g)}^{-1}$ of $G_{E/F}(N)$. 

We say that $\mu$ is {\it conjugate self-dual} if we have:
\[
\mu^{*} = \mu.
\]
Extending this notation slightly, if $\psi^N = \mu \boxtimes \nu \in \Psi_{\simp}(N)$, then we put
\[
(\psi^N)^* : = \mu^* \boxtimes \nu.
\]
(Note that $\nu$, being a finite dimensional algebraic representation of $\SL_2(\mathbf{C})$, is self-dual.)

We now define the subset 
\[
\widetilde{\Psi}_{\simp}(N) \subset \Psi_{\simp}(N)
\]
consisting of conjugate self-dual elements, i.e. elements $\psi^N \in \Psi_{\simp}(N)$ satisfying $(\psi^N)^*=\psi^N$.

More generally we define the set $\Psi(N)$ consisting of objects given by (unordered) formal direct sums   
\begin{eqnarray}
\psi^N = l_1 \psi_{1}^{N_1} \boxplus \cdots \boxplus l_r \psi_{r}^{N_r}
\end{eqnarray}
where $l_i \geq 1$ are integers, and $\psi^{N_i}_i \in \Psi_{\simp}(N_i)$, mutually distinct, and subject to the condition
\[
N= l_1 \cdot N_1 + \cdots + l_r \cdot N_r.
\]
Elements of $\Psi(N)$ will be referred to as formal global parameters for $G_{E/F}(N)$.

Given $\psi^N \in \Psi(N)$ as in (2.3.3), let $P$ be the standard parabolic subgroup of $G_{E/F}(N) = \Res_{E/F} \GL_{N/E}$ corresponding to the partition:
\begin{eqnarray*}
(\underbrace{N_1,\cdots,N_1}_{l_1} ,\cdots,\underbrace{N_r,\cdots,N_r}_{l_r})
\end{eqnarray*}
we take $\pi_{\psi^N}$ to be the representation of $G_{E/F}(N)(\mathbf{A}_F)=\GL_N(\mathbf{A}_E)$ given by (normalized) induction:
\[
\pi_{\psi^N} = \mathcal{I}_P(\underbrace{\pi_{\psi_1^{N_1}} \boxtimes \cdots \boxtimes  \pi_{\psi_1^{N_1}} }_{l_1} \boxtimes \cdots \boxtimes \underbrace{\pi_{\psi_r^{N_r}} \boxtimes \cdots \boxtimes \pi_{\psi_r^{N_r}}}_{l_r})
\]
which is irreducible, by a theorem of Bernstein \cite{Be} and the fact that the representations $\pi_{\psi_i^{N_i}}$ are irreducible and unitary. The association $\psi^N \leftrightarrow \pi_{\psi^N}$ is then a bijection between the full $L^2$ automorphic spectrum of $G_{E/F}(N)(\mathbf{A}_F)=\GL_N(\mathbf{A}_E)$ and $\Psi(N)$. This follows from the theorem of Moeglin-Waldspurger already quoted above, together with 
Langlands' theory of Eisenstein series \cite{L,MW2}.

With $\psi^N \in \Psi(N)$ as before, we say that $\psi^N$ is conjugate self-dual, if there is an involution $i \leftrightarrow i^*$ of the indexing set $\{1,\cdots,r\}$, such that
\begin{eqnarray*}
 (\psi_i^{N_i})^* = \psi_{i^*}^{N_{i^*}}  
\end{eqnarray*} 
(hence $N_i=N_{i^*}$), and that
\[
l_i = l_{i^*}.
\]
The subset of conjugate self-dual parameters of $\Psi(N)$ will be noted as $\widetilde{\Psi}(N)$. If $\psi^N \in \widetilde{\Psi}(N)$ satisfies in addition the condition that $(\psi_i^{N_i})^* = \psi_i^{N_i}$ (i.e. $i^*=i$), and $l_i=1$ for all $i$, then $\psi^N$ is called {\it elliptic}. The subset of elliptic parameters is noted as $\widetilde{\Psi}_{\ellip}(N)$. We thus have a chain of inclusion:
\[
\widetilde{\Psi}_{\simp}(N) \subset \widetilde{\Psi}_{\ellip}(N) \subset \widetilde{\Psi}(N).
\]
A parameter $\psi^N$ is called {\it generic}, if for all simple components $\psi_i^{N_i} = \mu_i \boxtimes \nu_i$ of $\psi^N$, the factor $\nu_i$ is the trivial representation of $\SL_2(\mathbf{C})$. We will generally denote such a generic parameter as $\phi^N$. Denote by $\Phi(N)$ the set of global generic parameters, and by $\widetilde{\Phi}(N)$ the subset of generic parameters that are conjugate self-dual. The set $\widetilde{\Phi}_{\simp}(N)$ of conjugate self-dual simple generic parameters, namely those given by a conjugate self-dual cuspidal automorphic representation on $\GL_N(\mathbf{A}_E)$, plays a fundamental role in the global classification of the automorphic spectrum of $U_{E/F}(N)$.

We also follow the notational convention of \cite{A1}: given $\psi^N \in \widetilde{\Psi}(N)$ as in (2.3.3), write $K_{\psi^N}$ for the indexing set $\{1,\cdots,r\}$. We have a decomposition:
\begin{eqnarray}
K_{\psi^N} = I_{\psi^N} \coprod J_{\psi^N} \coprod (J_{\psi^N})^*
\end{eqnarray}
where $I_{\psi^N}$ consists of the set of indices that are fixed under the involution $i \leftrightarrow i^*$, while $J_{\psi_N}$ consists of the set of indices whose orbit under the involution contains exactly two elements. We can then write:
\begin{eqnarray}
\psi^N= \Big(\bigboxplus_{i \in I_{\psi^N}}   l_i \psi_i^{N_i}  \Big)   \boxplus \Big(  \bigboxplus_{j \in J_{\psi^N}} l_j(\psi_j^{N_j} \boxplus \psi_{j^*}^{N_{j^*}})   \Big).
\end{eqnarray} 

Recall the result of Jacquet-Shalika \cite{JS}. In general suppose that $G$ is a connected reductive group over $F$. For $S$ any finite set of primes of $F$ outside of which $G$ is unramified, denote by $\mathcal{C}_{\mathbf{A}_F}^S(G)$ the set of adelic families of semi-simple $\widehat{G}$-conjugacy classes $c^S$ of the form:
\[
c^S = (c_v)_{v \notin S}
\]
with $c_v $ a $\widehat{G}$-conjugacy class in $\leftexp{L}{G_v} = \widehat{G} \rtimes W_{F_v}$ represented by an element of the form $c_v = t_v \rtimes \Frob_v$, with $t_v$ a semi-simple element of $\widehat{G}$. If $\pi = \otimes^{\prime}_v \pi_v$ is an irreducible admissible representation of $G(\mathbf{A}_F)$, then by taking $S$ to be a finite set of primes outside of which both $G$ and $\pi$ are unramified, we obtain the family $c^S(\pi) \in \mathcal{C}^S_{\mathbf{A}_F}(G)$ by the Satake isomorphism. 

Given $S,S^{\prime}$, and $c_1^S \in \mathcal{C}^S_{\mathbf{A}_F}(G)$, $c_2^{S^{\prime}} \in \mathcal{C}^S_{\mathbf{A}_F}(G)$, define the relation
\[
c_1^{S} \sim  c_2^{S^{\prime}}
\]
if the $\widehat{G}$-conjugacy classes $c_{1,v}$ and $c_{2,v}$ are the same for almost all primes $v$. Put
\[
\mathcal{C}_{\mathbf{A}_F}(G) = \varinjlim_{S} \mathcal{C}^{S}_{\mathbf{A}_F}(G).
\]

In the case where $G=G_{E/F}(N)$, put 
\[
\mathcal{C}_{\mathbf{A}_F}(N) = \mathcal{C}_{\mathbf{A}_F}(G_{E/F}(N)). 
\]
Then the theorem of Jacquet-Shalika \cite{JS} (coupled with the theorem of Moeglin-Waldspurger \cite{MW} mentioned above) gives an injection from $\mathcal{A}(N)$, the set of automorphic representations of $G_{E/F}(N)(\mathbf{A}_F)=\GL_N(\mathbf{A}_E)$ that belong to the $L^2$ spectrum, to $\mathcal{C}_{\mathbf{A}_F}(N)$:
\begin{eqnarray}
\mathcal{A}(N) &\hookrightarrow & \mathcal{C}_{\mathbf{A}_F}(N) \\
\pi &\mapsto & c(\pi). \nonumber
\end{eqnarray}

Denote by $\mathcal{C}(N)=\mathcal{C}_{\aut}(N)$ the subset of $\mathcal{C}_{\mathbf{A}_F}(N)$ given by the image of $\mathcal{A}(N)$. Then we can identify the following three sets:
\[
\mathcal{A}(N) = \Psi(N) = \mathcal{C}_{\aut}(N).
\]

We also denote by $\widetilde{\mathcal{C}}(N)=\widetilde{\mathcal{C}}_{\aut}(N)$ the subset of conjugate self-dual families of $\mathcal{C}_{\aut}(N)$. Then we can identify:
\[
\widetilde{\Psi}(N) = \widetilde{\mathcal{C}}_{\aut}(N).
\]

\bigskip

Finally we discuss the localization of global parameters. 

In general if $v$ is a prime of $F$, we denote by the subscript $v$ for the localization of various objects at the prime $v$. Thus let $v$ be a prime of $F$ that does not split in $E$. Then the localization $E_v$ is a quadratic extension of $F_v$, and $G_{E_v/F_v}(N)$ is the localization of $G_{E/F}(N)$ at $v$. Given a simple generic parameter $\phi^N \in \Phi_{\simp}(N)$, in other words a unitary cuspidal automorphic representation $\mu$ on $\GL_N(\mathbf{A}_E)$, the localization $\mu_v$ is then an irreducible admissible representation of $\GL_N(E_v)$. Hence by the local Langlands correspondence for general linear groups as established by Harris-Taylor \cite{HT} and Henniart \cite{H1}, $\mu_v$ corresponds to a local $L$-parameter in $\Phi_v(N) := \Phi(G_{E_v/F_v}(N))=\Phi(\GL_N(E_v))$, which we note as $\phi^N_{v}$. This gives the localization $\phi^N \rightarrow \phi^N_{v}$ from $\Phi_{\simp}(N)$ to $\Phi_v(N)$, which takes $\widetilde{\Phi}_{\simp}(N)$ to $\widetilde{\Phi}_v(N)$. Similarly it follows that we have a localization map $\psi^N \rightarrow \psi^N_{v}$ from $\Psi(N)$ to $\Psi_v^+(N)$, which takes $\widetilde{\Psi}(N)$ to $\widetilde{\Psi}^+_v(N)$. We emphasize that we cannot say {\it a priori} that the localization map takes $\Psi(N)$ to $\Psi_v(N)$, since the generalized Ramanujan conjecture for cuspidal automorphic representation of general linear groups is not yet known.

If $v$ is a prime of $F$ that splits in $E$, then $E_v = E_w \times E_{\overline{w}}$, where $w$ and $\overline{w}$ are the primes of $E$ above $v$. If $\phi^N \in \widetilde{\Phi}_{\simp}(N)$ is a simple generic parameter represented by a conjugate self-dual cuspidal automorphic representation $\mu$ of $\GL_N(\mathbf{A}_E)$, we denote by $\phi^N_{w}$ (resp. $\phi^N_{\overline{w}}$) the $L$-parameter of $\Phi(\GL_N(E_w))$ (resp. the $L$-parameter of $\Phi(\GL_N(E_{\overline{w}}))$) corresponding to the localization $\mu_w$ as an irreducible admissible representation of $\GL_N(E_w)$ (resp. $\mu_{\overline{w}}$ as irreducible admissible representation of $\GL_N(E_{\overline{w}})$). If we identify $E_w = E_{\overline{w}} = F_v$, then we have $\mu_w \cong \mu_{\overline{w}}^{\vee}$, and hence $\phi^N_{w}=(\phi^N_{\overline{w}})^{\vee}$. More generally, if $\psi^N \in \widetilde{\Psi}(N)$, then the localizations $\psi^N_{w},\psi^N_{\overline{w}} \in \widetilde{\Psi}^+(\GL_N(F_v))$ satisfies $\psi^N_{w} = (\psi^N_{\overline{w}})^{\vee}$.

\subsection{Endoscopic data and parameters}

We briefly recall the notion of endoscopic data (both standard and twisted case) from chapter two of Kottwitz-Shelstad \cite{KS}. Thus $G^0$ is a connected reductive group over $F$, and $\theta$ is a $F$-rational semi-simple automorphism of $G^0$. Denote by $G: = G^0 \rtimes \theta$ the resulting $G^0$-bitorsor (or a twisted group). The particular case where $\theta$ is the identity automorphism (and thus $G=G^0$) is referred to as the standard case. If $\theta$ is of finite order, which is the case considered in our situations, the twisted group $G= G^0 \rtimes \theta$ is equal to a component of the non-connected reductive group $ G^+ = G^0 \rtimes \langle \theta \rangle$ generated by $G$. 

An endoscopic datum for $G$ is a quadruple $(G^{\prime},\mathcal{G}^{\prime},s,\xi^{\prime})$, subject to the conditions listed in section 2.1 of \cite{KS}. Here $s$ is a sem-simple element of $\widehat{G} = \widehat{G}^0 \rtimes \widehat{\theta}$ of the form $s = s^0 \rtimes \widehat{\theta}$ for a semi-simple element $s^0 \in \widehat{G}^0$, and $\widehat{\theta}$ is the automorphism of $\widehat{G}^0$ that is dual to $\theta$ and preserves a fixed $\Gamma_F$-splitting of $\widehat{G}^0$, with $s$ considered up to translation by $Z(\widehat{G}^0)^{\Gamma_F}$. The most important condition for endoscopic data is that $G^{\prime}$ is a connected quasi-split group over $F$, called the endoscopic group, whose dual group $\widehat{G}^{\prime}$ is the connected centralizer of $s$ in $\widehat{G}^0$, equivalently $\widehat{G}^{\prime}$ is the connected $\widehat{\theta}$-twisted centralizer of $s^0$ in $\widehat{G}^0$:
\begin{eqnarray*}
\widehat{G}^{\prime} = \Cent(s,\widehat{G}^0)^0
\end{eqnarray*}
with
\begin{eqnarray*}
\Cent(s,\widehat{G}^0) & =&  \{ g \in \widehat{G}^0, \,\ gs =sg      \}\\
&= &  \{ g \in \widehat{G}^0, \,\ g s^0 = s^0 \widehat{\theta}(g)       \}.
\end{eqnarray*}
The datum $\mathcal{G}^{\prime}$ is a split extension of $W_F$ by $\widehat{G}^{\prime}$, and $\xi^{\prime}:\mathcal{G}^{\prime} \rightarrow \leftexp{L}{G}$ is an $L$-embedding. We refer to section 2.1 of \cite{KS} for the definition of equivalence of endoscopic data, and the definition of the outer automorphism group $\Out_G(G^{\prime})$ of the endoscopic data.

The set of equivalence classes of endoscopic data of $G$ will be noted as $\mathcal{E}(G)$. An endoscopic data is called {\it elliptic}, if we have
\[
(Z(\widehat{G}^{\prime})^{\Gamma_F})^0 \subset (Z(\widehat{
G})^{\Gamma_F})
\]
where $Z(\widehat{G}) := Z(\widehat{G}^0)^{\widehat{\theta}}$. The set of equivalence classes of elliptic endoscopic data of $G$ is noted as $\mathcal{E}_{\ellip}(G)$. Given an endoscopic datum $(G^{\prime}, \mathcal{G}^{\prime},s,\xi^{\prime})$, a lot of times when the context is clear, we will abuse notation and will write $G^{\prime} \in \mathcal{E}(G)$ to mean that we are taking the endoscopic datum $G^{\prime}=(G^{\prime},\mathcal{G}^{\prime},s,\xi^{\prime})$ to be representative of the equivalence class that it defines. 

In all the situations that we consider, the datum $\mathcal{G}^{\prime}$ can always be taken as the $L$-group of $G^{\prime}$, i.e. we can take $\mathcal{G}^{\prime} = \leftexp{L}{G^{\prime}}$, and thus the datum $\mathcal{G}^{\prime}$ will be omitted from the notation. The datum $\xi^{\prime}$ is thus an $L$-embedding of $\leftexp{L}{G^{\prime}}$ to $\leftexp{L}{G}$ extending the identity inclusion of $\widehat{G}^{\prime} $ to $\widehat{G}$. In this case, a general endoscopic datum is of the form $(M,s_M,\xi^{\prime}_M)$, with $M$ a Levi subgroup of $G^{\prime}$ for some $(G^{\prime},s,\xi^{\prime}) \in \mathcal{E}_{\ellip}(G)$, with $\xi^{\prime}_M$ being the composition of $\xi^{\prime}$ with the $L$-embedding $\leftexp{L}{M} \rightarrow \leftexp{L}{G^{\prime}}$ that is dual to the Levi embedding $M \hookrightarrow G^{\prime}$.  

We will usually denote an (equivalence class of) endoscopic datum as $G^{\prime}$, or $(G^{\prime},\xi^{\prime})$ if we want to emphasize the $L$-embedding $\xi^{\prime}$ of the datum. 

For our purpose, we need the to consider the case of standard endoscopy with $G=U_{E/F}(N)$, and the case of twisted endoscopy, with $G=\widetilde{G}_{E/F}(N) := G_{E/F}(N) \rtimes \theta$, with $\theta$ as in (1.0.1), and $\widehat{\theta}$ as in (2.1.2). We first consider the standard case. 

Thus $G=U_{E/F}(N)$. The set $\mathcal{E}_{\ellip}(G)$ of (equivalence classes of) elliptic endoscopic data of $G=U_{E/F}(N)$ are determined by Rogawski (\cite{R} section 4.6). They are given by
\[
(G^{\prime},\xi^{\prime}) = (U_{E/F}(N_1) \times U_{E/F}(N_2), \zeta_{\underline{\chi}})
\]
\begin{eqnarray*}
& & N_1,N_2 \geq 0, \,\ N= N_1 + N_2 , \\
& & \underline{\chi} =(\chi_1,\chi_2) 
\end{eqnarray*}
with $\underline{\chi}$ having signature:
\begin{eqnarray*}
  \underline{\kappa} =(\kappa_1,\kappa_2)=((-1)^{N-N_1},(-1)^{N-N_2}) \mbox{ as in } (2.1.12)
\end{eqnarray*}
and $\zeta_{\underline{\chi}}$ as in (2.1.13). The equivalence class of the endoscopic data is uniquely determined by $N_1,N_2$ and is independent of the choice of $\chi_1 \in \mathcal{Z}_E^{\kappa_1},\chi_2 \in \mathcal{Z}_E^{\kappa_2}$. We also have
\begin{eqnarray}
\Out_G(G^{\prime}) = \left \{ \begin{array}{c}   \mbox{ trivial if } N_1 \neq N_2 \\  \,\ \,\ \mathbf{Z}/2\mathbf{Z}    \mbox{ if } N_1 = N_2.
\end{array} \right.
\end{eqnarray}
In the case where $N_1=N_2$ the non-trivial element of $\Out_G(G^{\prime})$ is represented by the automorphism of $G^{\prime}$ that switches the two $U_{E/F}(N_1)$ factors. Following \cite{A1} we also put $\mathcal{E}_{\simp}(G)$ to be the (unique up to equivalence) endoscopic datum represented by $G=U_{E/F}(N)$ itself.

Next consider the twisted case. Thus we take $G= \widetilde{G}_{E/F}(N)= G_{E/F}(N) \rtimes \theta$. The set of elliptic twisted endoscopic data $\mathcal{E}_{\ellip}(G)$ of the twisted group $G = \widetilde{G}_{E/F}(N)$, whose set of equivalence classes we denote as $\widetilde{\mathcal{E}}_{\ellip}(N)$, are given by the following ({\it c.f.} \cite{R}, section 4.7):

\begin{eqnarray*}
(G^{\prime},\xi^{\prime}) =(U_{E/F}(N_1) \times U_{E/F}(N_2),\xi_{\underline{\chi}_{\underline{\kappa}}} )
\end{eqnarray*}
subject to the condition:
\begin{eqnarray*}
& & N_1,N_2 \geq 0, \,\  N=N_1+N_2 \\
& & \underline{\chi}_{\underline{\kappa}}=(\chi_1,\chi_2) \mbox{ with signature } \underline{\kappa}, \mbox{ and } \\
& & \underline{\kappa} = \left \{ \begin{array}{c}  (1,-1) \mbox{ or } (-1,1) \mbox{ if } N_1 \equiv N_2 \bmod{2} \\  \,\ \,\ (1,1) \mbox{ or } (-1,-1)  \mbox{ if }  \mbox{ if } N_2 \notequiv N_2 \bmod{2}
\end{array} \right.
\end{eqnarray*}
with $\xi_{\underline{\chi}_{\underline{\kappa}} }$ as in (2.1.10). The equivalence class of the endoscopic datum $(G^{\prime},\xi^{\prime})$ is uniquely determined by the pair $(N_1,N_2)$ and $\underline{\kappa}$, and is independent of the choice of $\underline{\chi}$ having signature $\underline{\kappa}$.
When there is no confusion we abbreviate occasionally the endoscopic datum $(G^{\prime},\xi_{\underline{\chi}_{\underline{\kappa}}})$ as $(G^{\prime},\underline{\kappa})$. These endoscopic data are mutually inequivalent, except for the case where $N_1=N_2$, in which case the data $(U(N_1) \times U(N_1), (1,-1))$ and $(U(N_1) \times U(N_1),  (-1,1) )  $ are equivalent. We also denote by $\widetilde{\mathcal{E}}(N)$ the set of (equivalence class of) all endoscopic data of $\widetilde{G}_{E/F}(N)$. 

For $G^{\prime} \in \widetilde{\mathcal{E}}_{\ellip}(N)$ as above one has $\widetilde{\Out}_N(G^{\prime}) :=\Out_G(G^{\prime})$ (with $G=\widetilde{G}_{E/F}(N)$) being trivial.

Following Arthur \cite{A1}, we denote by $\widetilde{\mathcal{E}}_{\simp}(N) = \mathcal{E}_{\simp}(G)$ the set of (equivalence class of) simple twisted endoscopic data, namely given by $(U_{E/F}(N), \xi_{\chi_{\kappa}})$ for $\kappa = \pm 1$ (thus there are exactly two of them up to equivalence). The fact that $\widetilde{\Out}_N(U(N))$ is trivial makes the argument in the proof of the classification theorems in this paper to be simpler as compared to the symplectic-orthogonal situation considered in \cite{A1}.

For future reference, given a simple twisted endoscopic datum $(U_{E/F}(N),\xi_{\chi_{\kappa}})$ of $\widetilde{G}_{E/F}(N)$, we refer to the sign $(-1)^{N-1} \cdot \kappa$ as the {\it parity} of the datum.

\bigskip

\begin{rem}
\end{rem}
We note that when $N$ is odd, the two simple twisted endoscopic data $(U_{E/F}(N),\xi_{\chi_+})$ and $(U_{E/F}(N),\xi_{\chi_{-}})$ of $\widetilde{\mathcal{E}}_{\simp}(N)$, they have no equivalent Levi sub-data $(M,\xi_M) \in \widetilde{\mathcal{E}}(N)$ in common, with $M$ a Levi subgroup of $U_{E/F}(N)$. On the other hand, when $N$ is even then this property is not true, for example for the Siegel Levi $M \cong G_{E/F}(N/2)$, the two Levi sub-data $(M,\xi_{\chi_+})$ and $(M,\xi_{\chi_-})$ are equivalent as twisted endoscopic datum of $\widetilde{G}_{E/F}(N)$ (here we are still denoting by $\xi_{\chi_+}$ the composition of the dual Levi embedding $\leftexp{L}{M} \hookrightarrow \leftexp{L}{U_{E/F}(N)}$ with $\xi_{\chi_+}$; similarly for $\xi_{\chi_-}$); in fact any Levi $M$ with such property is a subgroup of the Siegel Levi. Thus the twisted endoscopy theory for unitary group is similar to the symplectic-orthogonal situation studied in \cite{A1}. The author is grateful to Moeglin for pointing this out.

Even though in the above situation with $N$ being even, the two twisted endoscopic data $(M,\xi_{\chi_+})$ and $(M,\xi_{\chi_-})$ are equivalent, it is still important to distinguish them notationally; thus in the following sections we will always use different notations for these two data.

\bigskip

For the rest of section 2.4 we fix $\chi_+ \in \mathcal{Z}_E^+$ and $\chi_{-} \in \mathcal{Z}_E^-$ (we can just take $\chi_+=1$ but this is not necessary). We can now state the first main result of the paper, which we view as a ``seed theorem". It will be proved by induction and comparison of trace formulas, but the complete proof will only be achieved at the end of section nine (but special cases of the theorem have to be proved along the way).

\begin{theorem}
Suppose that $\phi^N \in \widetilde{\Phi}_{\simp}(N)$ is a conjugate self-dual simple generic global parameter. Then there exists a unique (up to equivalence) twisted endoscopic data $(G,\xi_{\chi_{\underline{\kappa}}}) \in \widetilde{\mathcal{E}}_{\ellip}(N)$, such that 
\begin{eqnarray*}
c(\phi^N) = \xi_{\chi_{\underline{\kappa}}}(c(\pi))
\end{eqnarray*}
for some representation $\pi$ in the discrete automorphic spectrum of $G(\mathbf{A}_F)$. Furthermore, $(G,\xi_{\chi_{\underline{\kappa}}})$ is in fact simple, thus of the form $(U_{E/F}(N),\xi_{\chi_{\kappa}})$ for a unique $\kappa \in \{\pm 1 \}$. 
\end{theorem}
As we will see in theorem 2.5.4 the next subsection, the sign $\kappa$ is uniquely determined by the order of pole at $s=1$ of a certain Asai $L$-function associated to $\phi^N$.

For the moment suppose that $N$ is fixed, and assume theorem 2.4.2 is valid for any integer $m \leq N$. We can then define the fundamental objects needed for the statement of the global classification. 

Thus let $\psi^N \in \widetilde{\Psi}(N)$ be written in the standard form:
\begin{eqnarray}
\psi^N= \Big(\bigboxplus_{i \in I_{\psi^N}}   l_i \psi_i^{N_i}  \Big)   \boxplus \Big(  \bigboxplus_{j \in J_{\psi^N}} l_j(\psi_j^{N_j} \boxplus \psi_{j^*}^{N_{j^*}})   \Big).
\end{eqnarray}

For each $i \in I_{\psi^N}$, we can apply theorem 2.4.2 to the simple generic factor $\mu_i \in \widetilde{\Phi}_{\simp}(m_i)$ of $\psi_i^{N_i} = \mu_i \boxtimes \nu_i$. This gives a simple endoscopic datum $(U_{E/F}(m_i),\xi_{\chi_{\delta_i}}) \in \widetilde{\mathcal{E}}_{\simp}(m_i)$ for a unique $\delta_i \in \{ \pm 1 \}$. Put $H_i :=U_{E/F}(m_i)$. On the other hand, for an index $j \in J_{\psi_N}$, we just put $H_j:=G_{E/F}(m_j)$. 

Denote by $\{K_{\psi^N}\}$ the set of orbits of $K_{\psi^N}$ under the involution $i \leftrightarrow i^*$, which for our purpose can be identified as $I_{\psi^N} \coprod J_{\psi^N}$. Thus for $k \in K_{\psi^N}$ we have a connected reductive group $H_k$ defined over $F$. As usual we let $\leftexp{L}{H_k}$ be the Weil form of the $L$-group of $H_k$. Form the fibre product:

\begin{eqnarray}
\mathcal{L}_{\psi^N} = \prod_{k \in \{K_{\psi^N}\}} (\leftexp{L}{H_k} \rightarrow W_F)
\end{eqnarray}
which will serve as a substitute for the global Langlands group. We also define $\mathcal{L}_{\psi^N/E}$ to be the inverse image of $W_E$ in $\mathcal{L}_{\psi^N}$ under the projection map $\mathcal{L}_{\psi^N} \rightarrow W_F$. 

For an index $k=i \in I_{\psi^N}$, we have the embedding:
\[
\widetilde{\mu}_i := \xi_{\chi_{\delta_i}}: \leftexp{L}{U_{E/F}(m_i)} \longrightarrow \leftexp{L}{G_{E/F}(m_i)}.
\]
On the other hand, if $j \in J_{\psi^N}$, then we define the embedding:
\begin{eqnarray}
& & \\
& & \widetilde{\mu}_j : \leftexp{L}{G_{E/F}(m_j)} \longrightarrow   \leftexp{L}{(G_{E/F}(m_j) \times G_{E/F}(m_j))} \subset \leftexp{L}{G_{E/F}(2m_j)} \nonumber
\end{eqnarray}
by setting
\begin{eqnarray}
\widetilde{\mu}_j(h_j \rtimes \sigma) = (h_j \oplus \widehat{\theta}(h_j)) \rtimes \sigma,
\end{eqnarray}
 for $h_j \in \widehat{H}_j =\widehat{G}_{E/F}(m_j) = \GL_{m_j}(\mathbf{C}) \times \GL_{m_j}(\mathbf{C})$, $\sigma \in W_F$.

\begin{definition}
To the parameter $\psi^N \in \widetilde{\Psi}(N)$ as above, we associate the $L$-homomorphism:
\[
\widetilde{\psi}^N : \mathcal{L}_{\psi^N} \times \SL_2(\mathbf{C}) \rightarrow \leftexp{L}{G_{E/F}(N)}
\]
as the direct sum
\begin{eqnarray}
\widetilde{\psi}^N = \Big( \bigoplus_{i \in I_{\psi^N}} l_i(\widetilde{\mu_i} \boxtimes \nu_i) \Big) \oplus \Big( \bigoplus_{j \in J_{\psi^N}} l_j(\widetilde{\mu}_j \boxtimes \nu_j). \Big) 
\end{eqnarray}
\end{definition}

\noindent Here we have identified an $n$-dimensional representation $\nu:\SL_2(\mathbf{C}) \rightarrow \GL_n(\mathbf{C})$ as a homomorphism: 
\begin{eqnarray*}
& & \widetilde{\nu}: \SL_2(\mathbf{C}) \rightarrow \widehat{G}_{E/F}(n)=\GL_n(\mathbf{C}) \times \GL_n(\mathbf{C}) \\
& & \widetilde{\nu}(g) = (\nu(g),\nu(g)), \,\ g \in \SL_2(\mathbf{C})
\end{eqnarray*}
(note that any finite dimensional representation of $\SL_2(\mathbf{C})$ is self-dual). Henceforth we will just write $\widetilde{\nu}$ as $\nu$.

\bigskip

\begin{rem}
\end{rem}
Note that the group $\mathcal{L}_{\psi^N}$ can be defined from $\psi^N$ without using theorem 2.4.2. On the other hand, the definition of the $L$-embedding $\widetilde{\psi}^N$ depends crucially on the use of theorem 2.4.2.

\bigskip

\begin{definition}
Suppose that $(U_{E/F}(N),\xi_{\chi_{\kappa}}) \in \widetilde{\mathcal{E}}_{\simp}(N)$. Define $\Psi(U_{E/F}(N),\xi_{\chi_{\kappa}})$ to be the set consisting of pairs $\psi=(\psi^N, \widetilde{\psi})$, where $\psi^N \in \widetilde{\Psi}(N)$, and 
\[
\widetilde{\psi}: \mathcal{L}_{\psi_N} \times \SL_2(\mathbf{C}) \rightarrow \leftexp{L}{U_{E/F}(N)}
\] 
is an $L$-homomorphism (considered up to $\widehat{U}_{E/F}(N)$-conjugacy), such that
\begin{eqnarray}
\widetilde{\psi}^N = \xi_{\chi_{\kappa}} \circ \widetilde{\psi}.
\end{eqnarray}
\end{definition}

Note that $\widetilde{\psi}$ is determined by $\widetilde{\psi}^N$ and $\xi_{\chi_{\kappa}}$, i.e. $\Psi(U_{E/F}(N),\xi_{\chi_{\kappa}})$ can be defined to be consisting of the set of $\psi^N \in \widetilde{\Psi}(N)$ such that $\widetilde{\psi}^N$ factors through the $L$-embedding $\xi_{\chi_{\kappa}}$.

For $\psi=(\psi^N,\widetilde{\psi}) \in \Psi(U_{E/F}(N),\xi_{\chi_{\kappa}})$, we denote $\mathcal{L}_{\psi^N}$ as $\mathcal{L}_{\psi}$. We also denote the indexing sets $K_{\psi^N},I_{\psi^N}, J_{\psi^N}$ in (2.4.3) as $K_{\psi},I_{\psi},J_{\psi}$. As usual we denote by $\Phi(U_{E/F}(N),\xi_{\chi_{\kappa}})$ the subset of generic parameters, i.e. those that are trivial on the $\SL_2(\mathbf{C})$ factor.

\bigskip

\begin{rem}
\end{rem}
In the above notation, for the component $\mu_i$ of $\psi_i^{N_i} = \mu_i \boxtimes \nu_i \in \widetilde{\Psi}_{\simp}(N_i)$, we have $\widetilde{\mu}_i:\mathcal{L}_{\psi} \rightarrow \leftexp{L}{G_{E/F}(m_i)}$ factors through $\xi_{\chi_{\delta_i}}$ tautologically by construction. The analogue of lemma 2.2.1 in the current setting (which can be proved by exactly the same argument) shows that $\widetilde{\mu}_i|_{\mathcal{L}_{\psi/E}}$, identified as an admissible representation $\widetilde{\mu}_i|_{\mathcal{L}_{\psi_{/E}}} :\mathcal{L}_{\psi/E} \rightarrow \GL_{m_i}(\mathbf{C})$, is conjugate self-dual with parity given by the sign $\delta_i (-1)^{m_i-1}$. On the other hand the $n_i$-dimensional representation $\nu_i$ of $\SL_2(\mathbf{C})$ is orthogonal if $n_i$ is odd, and symplectic if $n_i$ is even. Hence the representation 
\[
\widetilde{\mu}_i|_{\mathcal{L}_{\psi/E}} \boxtimes \nu_i : \mathcal{L}_{\psi/E } \times \SL_2(\mathbf{C}) \rightarrow \GL_{N_i}(\mathbf{C}) \mbox{ (with } N_i = m_i n_i) 
\]
is conjugate self-dual with parity $\delta_i (-1)^{m_i - 1 + n_i -1}=\delta_i(-1)^{m_i+n_i}$. Thus if we define:
\begin{eqnarray}
\kappa_i := \delta_i (-1)^{N_i -m_i -n_1+1}
\end{eqnarray}
Then we have $\kappa_i (-1)^{N_i - 1} = \delta_i (-1)^{m_i+n_i}$. By the analogue of lemma 2.2.1 again the $L$-homomorphism
 \[
\widetilde{\mu}_i \boxtimes \nu_i : \mathcal{L}_{\psi} \times \SL_2(\mathbf{C}) \rightarrow \leftexp{L}{G_{E/F}(N_i)}
\]
factors through $\xi_{\chi_{\kappa_i}}$. Thus the simple parameter $\psi_i^{N_i} = \mu_i \boxtimes \nu_i \in \widetilde{\Psi}_{\simp}(N_i)$ defines a parameter $\psi_i = (\psi_i^{N_i}, \widetilde{\psi}_{i}  )$ in $\Psi(U_{E/F}(N_i),\xi_{\chi_{\kappa_i}} )$, with $\kappa_i$ as in (2.4.9), and $\widetilde{\psi}_{i}:=\widetilde{\mu}_i \boxtimes \nu_i$

\bigskip

More generally, if $(G,\xi) \in \widetilde{\mathcal{E}}(N)$ represents a general twisted endoscopic datum of $\widetilde{G}_{E/F}(N)$, then $\Psi(G,\xi)$ can again be identified as the set of pairs $(\psi^N,\widetilde{\psi})$, where $\psi^N \in \widetilde{\Psi}(N)$, and 
\[
\widetilde{\psi} : \mathcal{L}_{\psi} \times \SL_2(\mathbf{C}) \rightarrow \leftexp{L}{G}
\]
is an $L$-embedding, such that
\[
\widetilde{\psi}^N = \xi \circ \widetilde{\psi}.
\] 

For example if $(G,\xi)=(G,\xi_{\underline{\chi}
_{\underline{\kappa}}}) \in \widetilde{\mathcal{E}}_{\ellip}(N)$, we can write:
\[
G = G_1 \times G_2
\]
and with $\underline{\kappa}=(\kappa_1,\kappa_2)$, such that $(G_{i},\xi_{\chi_{\kappa_i}})$ represents a simple twisted endoscopic datum of the form $(U_{E/F}(N_{i}),\xi_{\chi_{\kappa_i}}) \in \widetilde{\mathcal{E}}_{\simp}(N_{i})$. Then we have
\[
\Psi(G,\xi_{\underline{\chi}_{\underline{\kappa}}}) =\Psi(G_1,\xi_{\chi_{\kappa_1}}) \times \Psi(G_2,\xi_{\chi_{\kappa_2}}).
\]

\bigskip

We note that, in contrast to the case where $G \in \widetilde{\mathcal{E}}_{\simp}(N)$, the projection $(\psi^N,\widetilde{\psi}) \rightarrow \psi^N$ is in general no longer injective.

\bigskip

Of particular importance will be the following:

\begin{definition}
Let $(G,\xi) \in \widetilde{\mathcal{E}}_{\ellip}(N)$ be an elliptic twisted endoscopic datum. Define $\Psi_2(G,\xi)$ to be the subset of $\Psi(G,\xi)$ consisting of $\psi=(\psi^N,\widetilde{\psi})$ such that $\psi^N \in \widetilde{\Psi}_{\ellip}(N)$. 
\end{definition}
$\Psi_2(G,\xi)$ is known as the set of square-integrable parameters of $G$, with respect to $\xi$. Thus for example if $\psi = (\psi^N,\widetilde{\psi}) \in \Psi_2(U_{E/F}(N),\xi_{\chi_{\kappa}})$, then we have:
\[
\psi^N = \psi_1^{N_1} \boxplus \cdots \boxplus \psi_r^{N_r}
\]
for mutually distinct $\psi_i^{N_i} \in \widetilde{\Psi}_{\simp}(N_i)$, such that $\kappa_i (-1)^{N_i-1} = \kappa (-1)^{N-1}$ for all $i=1,\cdots,r$. Here $\kappa_i$ is defined as in (2.4.9). 

\bigskip

For general $(G,\xi) \in \widetilde{\mathcal{E}}_{\ellip}(N)$, we see that the projection $\psi=(\psi^N,\widetilde{\psi}) \rightarrow \psi^N$ on $\Psi_2(G,\xi)$ is again injective.

\bigskip

Given $(G,\xi) \in \widetilde{\mathcal{E}}(N)$, define (in a formal manner analogous to the local situation) the map:
\begin{eqnarray*}
& & \xi_* : \Psi(G,\xi) \rightarrow \widetilde{\Psi}(N) \\
& & \psi = (\psi^N, \widetilde{\psi}) \rightarrow \psi^N.
\end{eqnarray*}

\bigskip

Then from the definitions it follows easily that
\begin{eqnarray*}
& & \widetilde{\Psi}_{\ellip}(N) \\ 
&=& \coprod_{(G,\xi) \in \widetilde{\mathcal{E}}_{\ellip}(N)} \xi_*(\Psi(G,\xi)).
\end{eqnarray*}

Indeed suppose that $\psi^N = \bigboxplus_{i \in I} \psi_i^{N_i}  \in \widetilde{\Psi}_{\ellip}(N)$ for $\psi_i^{N_i} \in \widetilde{\Psi}_{\simp}(N_i)$, with $\widetilde{\psi}_i^{N_i}$ factors through $\xi_{\chi_{\kappa_i}}$ ($\kappa_i$ being determined as in (2.4.9)). Then $(G,\xi)=(G,\underline{\xi}_{\underline{\chi}_{\underline{\kappa}}})$ is determined as follows. Write $\psi^N$ in the form
\[
\psi^N = \psi^N_{O} \boxplus \psi^N_{S}
\]
with 
\[
\psi^N_{O} = \bigboxplus_{i \in I_O} \psi_i^{N_i}, \,\ \psi_{N,S} = \bigboxplus_{i \in I_S} \psi_i^{N_i}
\]
being the ``conjugate orthogonal" and ``conjugate symplectic" part of $\psi^N$ respectively; here the indexing set $I_O$ and $I_S$ are given by:
\[
I_O = \{i \in I, \,\  \kappa_i (-1)^{N_i-1} =1  \}, \,\  I_S = \{i \in I, \,\  \kappa_i (-1)^{N_i-1} =-1  \}.
\]
If we put
\[
N_O = \sum_{i \in I_O} N_i, \,\ N_S = \sum_{i \in I_S} N_i
\]
then
\[
G = U_{E/F}(N_O) \times U_{E/F}(N_S), \,\ \underline{\kappa} = ( (-1)^{N_O-1},(-1)^{N_S})
\]

\noindent Note that the simple twisted endoscopic data $(U_{E/F}(N_O), (-1)^{N_O-1})$ and $(U_{E/F}(N_S),(-1)^{N_S})$ (of $\widetilde{G}_{E/F}(N_O)$ and $\widetilde{G}_{E/F}(N_S)$ respectively) have {\it opposite parity}.

\bigskip

\begin{definition}
Let $(G,\xi)=(U_{E/F}(N),\xi_{\chi_{\kappa}})$, and $\psi=(\psi^N,\widetilde{\psi}) \in \Psi(G,\xi_{\chi_{\kappa}})$. Define
\begin{eqnarray}
S_{\psi}(G) := \Cent(\Image \widetilde{\psi},\widehat{G})
\end{eqnarray}
the centralizer in $\widehat{G}$ of the image of $\widetilde{\psi}$. We also put
\begin{eqnarray*}
& & \overline{S}_{\psi} (G) := S_{\psi} /Z(\widehat{G})^{\Gamma_F}\\
& & \mathcal{S}_{\psi}(G) := \pi_0(\overline{S}_{\psi}).
\end{eqnarray*}
\end{definition}

Similar to the local case (2.2.12), we define the canonical central element $s_{\psi} \in S_{\psi}(G)$:
\begin{eqnarray}
s_{\psi} = \widetilde{\psi} \Big(1, \begin{pmatrix}   -1 &  0 \\    0  &    -1    \\ \end{pmatrix}     \Big).
\end{eqnarray}

The computation of $S_{\psi}(G)$ can easily be done. Thus let $\psi =(\psi^N,\widetilde{\psi}) \in \Psi(U_{E/F}(N),\xi_{\chi_{\kappa}})$. We can write $\psi^N$ in the form
\begin{eqnarray}
& & \\
& & \psi^N= \Big(\bigboxplus_{i \in I_{\psi}^+}   l_i \psi_i^{N_i}  \Big)   \boxplus \Big(\bigboxplus_{i \in I_{\psi}^-}   l_i \psi_i^{N_i}  \Big)  \boxplus \Big(  \bigboxplus_{j \in J_{\psi}} l_j(\psi_j^{N_j} \boxplus \psi_{j^*}^{N_{j^*}})   \Big) \nonumber
\end{eqnarray}
here $I_{\psi} = I_{\psi}^+ \coprod I_{\psi}^-$ is the partition determined as follows. Recall the signs $\kappa_i$ defined as in (2.4.9). We declare, given $i \in I_{\psi}$:

\begin{eqnarray}
 i \in \left \{ \begin{array}{c}  I_{\psi}^+ \mbox{ if }  \kappa_i = \kappa (-1)^{N-N_i} \\   I_{\psi}^- \mbox{ if } \kappa_i =  \kappa (-1)^{N-N_i-1}.
\end{array} \right.
\end{eqnarray} 

\begin{rem}
\end{rem}
By the analogue of lemma 2.2.1 in the present situation, the $L$-homomorphism:
\[
\widetilde{\psi}_i^{N_i}  : \mathcal{L}_{\psi_i} \times \SL_2(\mathbf{C}) \rightarrow \leftexp{L}{G_{E/F}(N_i)}
\]
corresponds to an $N_i$-dimensional representation of $\mathcal{L}_{\psi_i/E} \times \SL_2(\mathbf{C})$ that is conjugate self-dual with parity equal to $\kappa_i (-1)^{N_i-1}$, while the $L$-homomorphism
\[
\widetilde{\psi}^N : \mathcal{L}_{\psi} \rightarrow \leftexp{L}{G_{E/F}(N)}
\]
corresponds to an $N$-dimensional representation of $\mathcal{L}_{\psi/E} \times \SL_2(\mathbf{C})$ that is conjugate self-dual with parity $\kappa (-1)^{N-1}$. Thus the set $ I_{\psi}^+$ corresponds to the set of simple sub-parameters of $\psi^N$ that are conjugate self-dual of the same parity as $\psi^N$, while $I_{\psi}^-$ consists of those that are of opposite parity. 

\bigskip

It is then easy to see that, in order for $\widetilde{\psi}^N$ to factor through $\xi_{\chi_{\kappa}}: \leftexp{L}{U_{E/F}(N)} \rightarrow \leftexp{L}{G_{E/F}(N)}$, the set $I_{\psi}^-$ must be of even cardinality, and then we have ({\it c.f.} \cite{GGP} section 4):

\begin{eqnarray}
& & \\
& & S_{\psi}(U_{E/F}(N)) = \prod_{i \in I_{\psi}^+} O(l_i,\mathbf{C}) \times \prod_{i \in I_{\psi}^-} Sp(l_i,\mathbf{C}) \times \prod_{j \in J_{\psi}} \GL(l_j,\mathbf{C}). \nonumber
\end{eqnarray}
We have $Z(\widehat{U}_{E/F}(N))^{\Gamma_F}=\{\pm I_N\}$, identified in the natural way as an element of the right hand side of (2.4.14). It follows that
\begin{eqnarray}
| \mathcal{S}_{\psi} | = \left \{ \begin{array}{c}   2^{|I_{\psi}^+ |} \mbox{ if } l_i \mbox{ is even for all } i \in I_{\psi}^+ \\  2^{|I_{\psi}^+ |-1} \mbox{ otherwise. } \end{array} \right.
\end{eqnarray}

\bigskip

Since the twisted trace formula for the twisted group $\widetilde{G}_{E/F}(N)$ plays a key role in the proofs of the main results, we need to extend these notions to the twisted group $\widetilde{G}_{E/F}(N)$. Thus for $\psi^N \in \widetilde{\Psi}(N)$, we define
\begin{eqnarray}
\widetilde{S}_{\psi^N}(N) = S_{\psi^N}(\widetilde{G}_{E/F}(N)) = \Cent(\Image \widetilde{\psi}^N,\widehat{\widetilde{G}}_{E/F}(N))
\end{eqnarray}
and
\begin{eqnarray}
& & S^*_{\psi^N}(N)=S_{\psi^N}(\widetilde{G}^0_{E/F}(N)) \\
&=& S_{\psi^N}(G_{E/F}(N)) = \Cent(\Image \widetilde{\psi}^N,\widehat{G }(N)). \nonumber
\end{eqnarray}

Both $\widetilde{S}_{\psi^N}(N)$ and $S^*_{\psi^N}(N)$ are connected, with $\widetilde{S}_{\psi^N}(N)$ being a bi-torsor under the connected group $S^*_{\psi^N}(N)$. Since these are connected, they are much simpler objects. In particular the set of connected components
\[
\widetilde{\mathcal{S}}_{\psi^N}=\pi_0(\overline{\widetilde{S}}_{\psi^N}) = \pi_0( \widetilde{S}_{\psi^N}/ Z(\widehat{G}_{E/F}(N))^{\Gamma_F} )
\]
and 
\[
\mathcal{S}^*_{\psi^N}  = \pi_0(\overline{S}^*_{\psi^N}) = \pi_0(S_{\psi^N}^*/Z(\widehat{G}_{E/F})^{\Gamma_F} )
\]
are just singleton. However, they still play an important part of the proofs. These objects can also be formulated in the case of local parameters in the evident manner.

\bigskip

We now discuss the localizations of these objects. Thus suppose $v$ is a prime of $F$ that does not split in $E$. The localization of $G_{E/F}(N)$ at $v$ is $G_{E_v/F_v}(N)$, and similarly the localization of $U_{E/F}(N)$ at $v$ is $U_{E_v/F_v}(N)$. Recall from the last subsection that we have defined the localization mapping $\psi^N \rightarrow \psi^N_{v}$ from $\widetilde{\Psi}(N)$ to $\widetilde{\Psi}^+_v(N)$. If $\psi=(\psi^N,\widetilde{\psi}) \in \Psi(U_{E/F}(N),\xi_{\chi_{\kappa}})$ we would like to know that the localization $\psi^N_{v}$ also factors through $\xi_{\chi_{\kappa,v}}$ (here we denote by $\chi_{\kappa,v} \in \mathcal{Z}_{E_v}$ the localization of $\chi_{\kappa}$ at $E_v$). This is the content of the second ``seed" theorem:

\begin{theorem}
Suppose that $\phi =(\phi^N,\widetilde{\phi}) \in \Phi_{\simp}(U_{E/F}(N),\xi_{\chi_{\kappa}})$ is a simple generic global parameter. Then for any place $v$ of $F$ that does not split in $E$, the localization $\phi^N_{v}:L_{F_v} \rightarrow \leftexp{L}{G_{E_v/F_v}(N)}$ factors through $\xi_{\chi_{\kappa}}$. 
\end{theorem}
In other words, given $\phi=(\phi^N,\widetilde{\phi}) \in \Phi_{\simp}(U_{E/F}(N),\xi_{\chi_{\kappa}})$, and $v$ does not split in $E$, we can define the localization $\phi_v$ as the ($\widehat{U}_{E_v/F_v}(N)$-conjugacy class of) $L$-homomorphism:
\begin{eqnarray}
\phi_v: L_{F_v}  \rightarrow \leftexp{L}{U_{E_v/F_v}(N)}
\end{eqnarray}
such that $\phi^N_{v} = \xi_{\chi_{\kappa,v}} \circ \phi_v$, where $\xi_{\chi_{\kappa,v}}$
\[
\xi_{\chi_{\kappa,v}}: \leftexp{L}{U_{E_v/F_v}(N)} \longrightarrow \leftexp{L}{G_{E_v/F_v}(N)}.
\]
is as in (2.1.9) (in the local situation).

The proof of theorem 2.4.10 is not elementary, and the complete proof is achieved only in section nine by an elaborate induction argument (but again, special cases have to be proved along the way). In any case, for fixed $N$, assuming the validity of theorem 2.4.10 for for integers $m \leq N$, we obtain the immediate:

\begin{corollary}
Assume the validity of theorem 2.4.10 for $m \leq N$. Suppose that $\psi=(\psi^N,\widetilde{\psi}) \in \Psi(U_{E/F}(N),\xi_{\chi_{\kappa}})$. Then for for each prime $v$ of $F$ that does not split in $E$, the localization $\psi^N_{v}:L_{F_v} \rightarrow \leftexp{L}{G_{E_v/F_v}(N)}$ factors through $\xi_{\chi_{\kappa,v}}$. 
\end{corollary}
We can thus define the localization $\psi_v$ as the $\widehat{U}_{E_v/F_v}(N)$-conjugacy class of) $L$-homomorphism
\[
\psi_v:L_{F_v} \times \SU(2) \rightarrow \leftexp{L}{U_{E_v/F_v}(N)}
\]
such that $\psi^N_{v} =  \xi_{\chi_{\kappa}}     \circ \psi_v$.

\bigskip

The case where $v$ splits into two primes $w,\overline{w}$ in $E$ is elementary, and follows from the discussion at the end of section 2.3. Namely that we have $E_v = E_{w} \times E_{\overline{w}}$. We denote $U_{E_v/F_v}(N)$ as $U(N)_v$. Recall that we have the isomorphism $\iota_w: U(N)_v \stackrel{\simeq}{\rightarrow} \GL_{N/E_w}$ (resp. $\iota_{\overline{w}}: U(N)_v \stackrel{\simeq}{\rightarrow} \GL_{N/E_{\overline{w}}}$) corresponding to the projection of $E_v$ to $E_w$ (resp. projection to $E_{\overline{w}}$). If we identify $E_w =E_{\overline{w}}=F_v$ then the map $\iota_{\overline{w}} \circ \iota_w^{-1}$ is given by $g \mapsto J \leftexp{t}{g}^{-1} J^{-1}$. From the discussion at th end of the last subsection, given $\psi=(\psi^N,\widetilde{\psi}) \in \Psi(U_{E/F}(N),\xi_{\chi_{\kappa}})$, we can consider the localizations $\psi^N_{w}$ and $\psi^N_{\overline{w}}$ as elements of $\Psi^+(\GL_N(E_w))$ (resp. $\Psi^+(\GL_N(E_{\overline{w}})))$. 

We have the representation $\pi_{\psi^N_{w}}$ of $\GL_N(E_w)$ associated to $\psi^N_w$ as follows. First consider the case that $\psi^N_w \in \Psi(\GL_N(E_w))$. Then $\pi_{\psi^N_{w}}$ is the irreducible admissible representation of $\GL_N(E_w)$ whose $L$-parameter is given by $\phi_{\psi^N_{w}}$ as in (2.2.11), under the local Langlands classification for general linear groups. In general, given the parameter $\psi^N_w \in \Psi^+(\GL_N(E_w))$, there is a partition
\[
N=N_1 + \cdots + N_r
\]
and $\lambda_1,\cdots,\lambda_r \in \mathbf{R}$ satisfying $\lambda_1 > \cdots > \lambda_r$, and 
\[
\psi^{N_i}_w \in \Psi(\GL_{N_i}(E_w))
\]
such that
\[
\psi^N_w = \psi^{N_1}_{w} \otimes |\cdot |^{\lambda_1} \oplus \cdots \oplus \psi^{N_r}_{w} \otimes |\cdot |^{\lambda_r}. 
\]
Let $\pi_{\psi^{N_i}_w}$ be the irreducible admissible representation of $\GL_{N_i}(E_w)$ associated to $\psi^{N_i}_w$ as above. The partition $N=N_1 + \cdots + N_r$ defined the standard parabolic subgroup $P$ of $\GL_{N}(E_w)$, with Levi component $\GL_{N_1}(E_w) \times \cdots \times \GL_{N_r}(E_w)$. We then define:
\[
\pi_{\psi^N_w} :=\mathcal{I}_P \big(    ( \pi_{\psi^{N_1}_w} \otimes |\det|^{\lambda_1} ) \boxtimes \cdots \boxtimes   (\pi_{\psi^{N_r}_w} \otimes |\det|^{\lambda_r} ) \big).
\]
(thus is not irreducible in general because we are taking the full parabolic induction). In a similar way we define the representation $\pi_{\psi^N_{\overline{w}}}$ associated to $\psi^N_{\overline{w}}$. Note that in the case of generic parameters, then this correspondence amounts to the local Langlands classification for standard representations of general linear groups, instead of irreducible representations. 

From the conjugate self-duality of $\psi^N$, we have $\psi^N_{w} = (\psi^N_{\overline{w}})^{\vee}$ (if we make the identification $E_w =E_{\overline{w}} =F_v$), and without loss of generality, we may assume that $\psi^N_w= \leftexp{t}{(\psi^N_{\overline{w}})}^{-1}$ as actual homomorphism. Hence we have $\pi_{\psi^N_w} = (\pi_{\psi^N_{\overline{w}}})^{\vee}$. Thus we see that the pull-back of the representation $\pi_{\psi^N_{w}}$ to $U(N)_v$ via $\iota_{w}$ is isomorphic to the pull-back of the representation $\pi_{\psi^N_{\overline{w}}}$ to $U(N)_v$ via $\iota_{\overline{w}}$. We denote this representation of $U(N)_v$ as $\pi_{\psi_v}$, and we define $\psi_v \in \Psi^+(U(N)_v)$ to be the parameter obtained by composing $\psi^N_{w}: L_{F_v} \times \SU(2) \rightarrow \leftexp{L}{\GL_{N/E_w}}$ with the isomorphism $\leftexp{L}{\iota_w}^{-1}:\leftexp{L}{\GL_{N/E_w}} \stackrel{\simeq}{\rightarrow} \leftexp{L}{U_v(N)}$ (which is the same as that obtained by composing $\psi^N_{\overline{w}}$ with $\leftexp{L}{\iota_{\overline{w}}}^{-1}$). More concretely:
\[
\psi_v: L_{F_v} \times \SU(2) \rightarrow \leftexp{L}{U_{E_v/F_v}(N)}
\] 
is given by:
\[
\psi_v(\sigma) = (\psi^N_{w}(\sigma), \psi^N_{\overline{w}}(\sigma))
\]
Then we declare that $\pi_{\psi_v}$ corresponds to the parameter $\psi_v$.

\bigskip

Given theorem 2.4.10 (and corollary 2.4.11), we can now define the localization maps $L_{F_v} \rightarrow \mathcal{L}_{\psi}$ and for the component groups $\mathcal{S}_{\psi} \rightarrow \mathcal{S}_{\psi_v}$, which play a critical role in the global classification. 

First assume that $v$ does not split in $E$. Thus we assume the validity of theorem 2.4.10, and hence corollary 2.4.11, for all $m \leq N$. Let $\psi=(\psi^N,\widetilde{\psi}) \in \Psi(U_{E/F}(N),\xi_{\chi_{\kappa}})$. With notation as in (2.4.3) for $\psi^N$, consider an index $k=i$ that belongs to $ I_{\psi}$, with $\mu_k$ being the generic component of the simple parameter $\psi_k^{N_k} = \mu_k \boxtimes \nu_k \in \widetilde{\Psi}_{\simp}(N_k)$. Theorem 2.4.2 applied to $\mu_k \in \widetilde{\Phi}_{\simp}(m_k)$ gives a simple twisted endoscopic datum $(H_k,\xi_{\chi_{\delta_k}})$ for a unique sign $\delta_k \in \{ \pm 1\}$, and $H_i=U_{E/F}(m_k)$. Recall that we put $\widetilde{\mu}_k: \leftexp{L}{U_{E/F}(m_k)} \rightarrow \leftexp{L}{G_{E/F}(m_k)}$ to be given by $\xi_{\chi_{\delta_k}}$. The pair $(\mu_k,\widetilde{\mu}_k)$ thus define an element in $\widetilde{\Phi}_{\simp}(m_k,\xi_{\chi_{\delta_k}})$. Theorem 2.4.10 thus says that the localization $(\mu_k)_v$, which we identify as the $L$-parameter $(\mu_k)_v: L_{F_v} \rightarrow \leftexp{L}{G_{E_v/F_v}(m_k)}$, factors through $\xi_{\chi_{\delta_k,v}}$: 
\[
(\mu_k)_v: L_{F_v} \rightarrow \leftexp{L}{U_{E_v/F_v}(m_k)} \stackrel{\xi_{\chi_{\delta_k,v}}}{\longrightarrow} \leftexp{L}{G_{E_v/F_v}(m_k)}.
\]

By composing the map $L_{F_v} \rightarrow \leftexp{L}{U_{E_v/F_v}(m_k)}$ with the $L$-map:
\[
\leftexp{L}{U_{E_v/F_v}(m_k)} \rightarrow \leftexp{L}{U_{E/F}(m_k)} = \leftexp{L}{H_k}
\]
we obtain the commutative diagram of $L$-homomorphisms:
\begin{eqnarray}
\xymatrix{L_{F_v} \ar[r] \ar[d] & W_{F_v} \ar[d] \\  \leftexp{L}{H_k} \ar[r] & W_F }
\end{eqnarray}

If $k \in J_{\psi}$, then we have a similar diagram, which is elementary in this case and does not depend on theorem 2.4.10, namely that it is just given by the local Langlands classification for general linear groups. Taking fibre product over $\{K_{\psi}\}$ we obtain a commutative diagram

\begin{eqnarray}
\xymatrix{L_{F_v} \ar[r] \ar[d] & W_{F_v} \ar[d] \\  \mathcal{L}_{\psi} \ar[r] & W_F }
\end{eqnarray}
which fits into the larger commutative diagram:
\begin{eqnarray}
\xymatrix{L_{F_v} \times \SU_2 \ar[r]^{\psi_v } \ar[d] & \leftexp{L}{U_{E_v/F_v}(N)} \ar[r] \ar[d] & W_{F_v} \ar[d] \\  \mathcal{L}_{\psi} \times \SL_2(\mathbf{C}) \ar[r]^{\widetilde{\psi}} & \leftexp{L}{U_{E/F}(N)} \ar[r] & W_F }
\end{eqnarray}

We note that the same discussion can be carried out if $v$ splits into two primes $w,\overline{w}$ in $E$, and is elementary in the sense that it does not depend on using theorem 2.4.10. We just repeat the gist of the above discussion. From the discussion at the end of section 2.3, we have the localizations (with $\mu_k$ being the generic component of the simple parameter $\psi_k^{N_k}=\mu_k \boxtimes \nu_k$ as before):
\begin{eqnarray*}
& & (\mu_k)_w : L_{E_w} \rightarrow \leftexp{L}{\GL_{m_k/E_w}} = \GL_{m_k}(\mathbf{C}) \times W_{E_w} \\
& & (\mu_k)_{\overline{w}} : L_{E_{\overline{w}}} \rightarrow \leftexp{L}{\GL_{m_k/E_{\overline{w}}}} = \GL_{m_k}(\mathbf{C}) \times W_{E_{\overline{w}}}
\end{eqnarray*}
Making the identification $E_w=E_{\overline{w}}=F_v$, the two $L$-parameters $(\mu_k)_w$ and $(\mu_k)_{\overline{w}}^{\vee}$ are conjugate under $\GL_N(\mathbf{C})$. Without loss of generality assume that we have $(\mu_k)_{\overline{w}} = \leftexp{t}{(\mu_k)_w}^{-1}$ as actual homomorphism. We then define:
\[
 (\mu_k)_v : L_{F_v} \rightarrow \leftexp{L}{U_{E_v/F_v}(m_k)} 
\]
by the rule
\begin{eqnarray}
 (\mu_k)_v (\sigma) = (  (\mu_k)_w(\sigma), (\mu_k)_{\overline{w}}(\sigma) ). 
\end{eqnarray}

Now we have fixed the embedding $\overline{F} \hookrightarrow \overline{F}_v$ (which also correspond to the embedding $W_{F_v} \hookrightarrow W_F$), thus fixes the choice of a prime above $v$ for each finite extension of $F$. In particular without loss of generality suppose that $w$ is the prime of $E$ above $v$ singled out by this embedding.   Corresponding to this we have the $L$-embedding
\begin{eqnarray}
\leftexp{L}{U_{E_v/F_v}(m_k)} \rightarrow \leftexp{L}{U_{E/F}(m_k)}
\end{eqnarray}
induced by $\widehat{\iota}_{w^{\prime}}$ and the embedding $W_{F_v} \hookrightarrow W_F$ ({\it c.f.} the discussion in section 2.1). 

By composing (2.4.22) with (2.4.23) we obtain the diagram (2.4.19) in the case where $v$ splits in $E$. Similarly we obtain (2.4.20) and (2.4.21) in the case where $v$ splits in $E$.

\bigskip

Thus in all cases the commutative diagram (2.4.21) allows us to define localization maps: 
\begin{eqnarray}
& & S_{\psi} \rightarrow S_{\psi_v} \\
& & \overline{S}_{\psi} \rightarrow \overline{S}_{\psi_v} \nonumber \\
& & \mathcal{S}_{\psi} \rightarrow \mathcal{S}_{\psi_v} \nonumber
\end{eqnarray}
with the first two localization maps being injective.

\subsection{Statement of main results}
With the preparation of the previous subsections we can now state the main results to be proved in this paper. The proof will be a long induction argument. Additional theorems that are needed to be proved along the way will be stated in the following sections.

We begin with the local classification theorem. Thus $F$ is now a local field. In general if $G$ is a connected reductive group over $F$, denote by $\Pi(G) = \Pi(G(F))$ the set of irreducible admissible representations of $G(F)$, and by $\Pi_{\temp}(G)$ the subset of irreducible tempered representations of $G(F)$.

\begin{theorem}
Suppose $F$ is local, and $E$ a quadratic extension of $F$. 
\bigskip

\noindent (a) For any local parameter $\psi \in \Psi(U_{E/F}(N))$, there is a finite multi-set $\Pi_{\psi}$, equipped with a canonical mapping
\[
\pi \longrightarrow \langle \cdot, \pi \rangle, \,\ \pi \in \Pi_{\psi}
\]
from $\Pi_{\psi}$ to the character group $\widehat{\mathcal{S}}_{\psi}$ of $\mathcal{S}_{\psi}$. If both $U_{E/F}(N)$ and $\pi$ are unramified then $\langle \cdot, \pi \rangle =1$. All the representations in the packet $\Pi_{\psi}$ are irreducible unitary representations.

\bigskip

\noindent (b) If $\psi=\phi \in \Phi_{\bdd}(U_{E/F}(N))$ (hence a generic parameter of $\Psi(U_{E/F}(N))$), then $\Pi_{\phi}$ is multiplicity free, and all the representations in $\Pi_{\phi}$ are tempered representations. The mapping from $\Pi_{\phi}$ to $\widehat{\mathcal{S}}_{\phi}$ is injective, and $\Pi(U_{E/F}(N))$ is the disjoint union of the packets $\Pi_{\phi}$ for all $\phi \in \Phi(U_{E/F}(N))$. If $F$ is non-archimedean, then the map from $\Pi_{\phi}$ to $\widehat{\mathcal{S}}_{\phi}$ is bijective.

\bigskip

\noindent Furthermore $\Pi_{\temp}(U_{E/F}(N))$ is the disjoint union of the packets $\Pi_{\phi}$ for all $\phi \in \Phi_{\bdd}(U_{E/F}(N))$.

\bigskip

\noindent (c) Let $\xi_{\chi} : \leftexp{L}{U_{E/F}(N)} \hookrightarrow \leftexp{L}{G_{E/F}(N)}$ be an $L$-embedding defined by a character $\chi \in \mathcal{Z}_E$. For any local parameter $\psi \in \Psi(U_{E/F}(N))$, put $\psi^N = \xi_{\chi,*} \psi$, regarded as a representation $\psi^N: L_E \times \SU_2 \rightarrow \GL_N(\mathbf{C})$ (as in section 2.2). Consider $\det(\psi_N)$ as a one-dimensional character $\det(\psi^N) : W_E \rightarrow \mathbf{C}^{\times}$. Then for any representations in the packet $\Pi_{\psi}$, their central character has parameter given by $\det(\psi^N)$, corresponding to the $L$-embedding $\leftexp{L}{U_{E/F}(1)} \hookrightarrow \leftexp{L}{G_{E/F}(1)}$ defined by the character $\chi^{\otimes N}$. 

\end{theorem}

In particular, part (b) of theorem 2.5.1 gives the local Langlands classification of representations for the group $U_{E/F}(N)$.

\bigskip

As in \cite{A1}, before we turn to the global classification, we also need the local packets $\Pi_{\psi}$ associated to parameters $\psi \in \Psi^+(U_{E/F}(N))$, due to the potential failure of the generalized Ramanujan conjecture for unitary cuspidal automorphic representations on general linear groups. 

Given $\psi \in \Psi^+(U_{E/F}(N))$, we have a standard parabolic subgroup $P=MN_P$ of $U_{E/F}(N)$, a parameter $\psi_M \in \Psi(M)$, and a point $\lambda$ in the open chamber of $P$ in the real vector space
\[
\mathfrak{a}_M^* =X(M)_F \otimes \mathbf{R} 
\] 
such that $\psi$ is the composition of $\psi_{M,\lambda}$ under the $L$-embedding $\leftexp{L}{M} \hookrightarrow \leftexp{L}{U_{E/F}(N)}$; here $\psi_{M,\lambda} \in \Psi^+(M)$ is the twist of $\psi_M$ by $\lambda$, i.e. $\psi_{M,\lambda}$ is the product of $\psi_M$ with the central Langlands parameter of $M$ that is dual to the unramified quasi-character 
\[
\chi_{\lambda} :m \mapsto e^{\lambda H_M(m)}, \,\ m \in M(F).
\]
Where as usual $H_M:M(F) \rightarrow \mathfrak{a}_M=\Hom(X(M)_F,\mathbf{R})$ is the homomorphism defined by the condition that for $\chi \in X(M)_F$,
\[
e^{ \langle H_M(m), \chi \rangle } = |\chi(m)|, \,\ m \in M(F).
\]
 The Levi subgroup $M \subset U_{E/F}(N)$ is a product of several general linear factors $G_{E/F}(N_i)$, with a group $G_- = U_{E/F}(N_-) \in \widetilde{\mathcal{E}}_{\simp}(N_-)$ for $N_- \leq N$. Since the local Langlands classification for the general linear factors is already known, it thus follows from theorem 2.5.1 that we can construct the packet $\Pi_{\psi_M}$ and the pairing $\langle \cdot, \pi_M \rangle$ for $\pi_M \in \Pi_{\psi_M}$. Letting $\pi_{M,\lambda}$ be the twist of $\pi_M$ by the character $\chi_{\lambda}$, we then define the packet $\Pi_{\psi}$ associated to $\psi \in \Psi^+(U_{E/F}(N))$ by:
\begin{eqnarray*}
\Pi_{\psi} = \{ \pi = \mathcal{I}_P(\pi_{M,\lambda}): \,\ \pi_M \in \Pi_{\psi_M}   \}.
\end{eqnarray*} 

\noindent Note that since the representations $\pi$ in the packet $\Pi_{\psi}$ are defined as full parabolic induction $\pi = \mathcal{I}_P(\pi_{M,\lambda})$, they are in general neither irreducible nor unitary. In any case, it also follows easily from the definitions that $S_{\psi_M}=S_{\psi}$ and $\mathcal{S}_{\psi_M}=\mathcal{S}_{\psi}$, and hence we can define the pairing $\langle \cdot, \pi \rangle $ for $\pi \in \Pi_{\psi}$ by the rule: for $s \in \mathcal{S}_{\psi}=\mathcal{S}_{\psi_M}$, we set
\begin{eqnarray*}
\langle s,\pi \rangle =\langle s,\pi_M \rangle, \,\ \pi=\mathcal{I}_P(\pi_{M,\lambda})
\end{eqnarray*}

\noindent In theorem 3.2.1 of section 3.2, we will state the endoscopic character identities that characterizes the representations in the packet $\Pi_{\psi}$ associated to parameter $\psi \in \Psi(U_{E/F}(N))$. One of the reasons that we defined the packets $\Pi_{\psi}$ associated to $\psi \in \Psi^+(U_{E/F}(N))$ above as full parabolic induction of the representations in the packet $\Pi_{\psi_M}$, is that the endoscopic character identities are still valid for the packets $\Pi_{\psi}$ for $\psi \in \Psi^+(U_{E/F}(N))$, by analytic continuation from that of $\Pi_{\psi_M}$.

\bigskip
\bigskip

We now turn to the global classification. So $F$ is now a global field, $E$ a quadratic extension of $F$, and as before for any prime $v$ of $F$ we use subscript $v$ to denote localization of various objects at $v$. For $\kappa=\pm 1$, we fix $\chi_{\kappa} \in \mathcal{Z}_E^{\kappa}$ (for $\kappa=+1$ we can of course simply take $\chi_+=1$, but this is not necessary). 

Thus let $\psi=(\psi^N,\widetilde{\psi})) \in \Psi(U_{E/F}(N),\xi_{\chi_{\kappa}})$ be a global parameter. If $v$ is a place of $F$ that does not split in $E$, then by corollary 2.4.11, the localization $\psi^N_{v}$ as a parameter $L_{F_v} \times \SU_2 \rightarrow \leftexp{L}{G_{E_v/F_v}(N)}$ factors through $\xi_{\chi_{\kappa},v}$, and we defined $\psi_v$ to be the parameter in $\Psi^+(U_{E_v/F_v}(N))$ such that $\psi^N_{v} = \xi_{\chi_{\kappa},v} \circ \psi_v$. By the local classification theorem 2.5.1 (extended to parameters in $\Psi^+(U_{E_v/F_v}(N))$ as above), we have a local packet $\Pi_{\psi_v}$ corresponding to $\psi_v$. 

Suppose now $v$ splits in $E$. Recall that in this case we defined the localization $\psi_v \in \Psi^+(U_{E_v/F_v}(N))$ and the representation $\pi_{\psi_v}$ of $U_{E_v/F_v}(N)$. We denote the singleton packet $\{ \pi_{\psi_v}\}$ as $\Pi_{\psi_v}$. In this case the group $S_{\psi_v}$ (which is isomorphic to $S_{\psi_w}$ or $S_{\psi_{\overline{w}}}$), is connected, and hence $\mathcal{S}_{\psi_v}$ is trivial, and we simply define the pairing $\langle \cdot , \pi_{\psi_v} \rangle$ to be trivial. 

Define the global packet $\Pi_{\psi}$ associated to $\psi$ as the restricted tensor product of the local packets $\Pi_{\psi_v}$:
\begin{eqnarray}
& & \Pi_{\psi} = \otimes_v^{\prime}  \Pi_{\psi_v}  \\
&=&  \{\pi = \otimes^{\prime}_v \pi_v, \,\  \pi_v \in \Pi_{\psi_v}, \,\ \langle \cdot , \pi_v  \rangle =1 \mbox{ for almost all } v \} \nonumber
\end{eqnarray}

Recall that as in (2.4.24) we have a localization map $\mathcal{S}_{\psi} \rightarrow \mathcal{S}_{\psi_v}$ (if $v$ splits in $E$ then this is just the trivial map). If $x \in \mathcal{S}_{\psi}$, we denote by $x_v \in \mathcal{S}_{\psi_v}$ its localization. We then have a map
\begin{eqnarray}
\Pi_{\psi} &\rightarrow& \widehat{\mathcal{S}}_{\psi} \\
 \pi &\mapsto & \langle \cdot , \pi \rangle \nonumber \\
\langle x, \pi \rangle &=& \prod_v \,\  \langle x_v, \pi_v \rangle \nonumber
\end{eqnarray}

We briefly recall the definition of the canonical sign character $\epsilon_{\psi} \in \widehat{\mathcal{S}}_{\psi}$ (section 1.5 of \cite{A1}).Thus write $G=U_{E/F}(N)$, and denote by $\widehat{\mathfrak{g}}$ the Lie algebra of $\widehat{G}=\GL_N(\mathbf{C})$. Define the representation:
\[
\tau_{\psi} : S_{\psi} \times \mathcal{L}_{\psi} \times \SL_2(\mathbf{C}) \rightarrow \GL(\widehat{\mathfrak{g}})
\]
by setting
\begin{eqnarray}
\tau_{\psi}(s,g,h) = \Ad(s \cdot \widetilde{\psi} (g \times h))
\end{eqnarray}
where $\Ad$ is the adjoint representation of $\leftexp{L}{G}$ on $\widehat{\mathfrak{g}}$. The representation $\tau_{\psi}$ is orthogonal, hence self-dual, since it preserves the Killing form on $\widehat{\mathfrak{g}}$. Decompose $\tau_{\psi}$ as:
\begin{eqnarray}
\tau_{\psi} = \oplus_{\alpha} \tau_{\alpha} = \oplus_{\alpha} (\lambda_{\alpha} \otimes \mu_{\alpha} \otimes \nu_{\alpha})
\end{eqnarray}
into irreducible representations $\lambda_{\alpha},\mu_{\alpha},\nu_{\alpha}$ of $S_{\psi},\mathcal{L}_{\psi}$ and $\SL_2(\mathbf{C})$ respectively. Then the character $\epsilon_{\psi}(x)$ is defined, for $x \in \mathcal{S}_{\psi}$, as:
\begin{eqnarray}
\epsilon_{\psi}(x) = \rightexp{\prod_{\alpha}}{\prime} \det(\lambda_{\alpha}(s)  )
\end{eqnarray}
here $s$ is any element in $S_{\psi}$ whose image in $\mathcal{S}_{\psi}$ is equal to $x$, and $\prod^{\prime}_{\alpha}$ denotes the product over the indices $\alpha$ such that $\mu_{\alpha}$ is symplectic and such that
\begin{eqnarray}
\epsilon(1/2,\mu_{\alpha})=-1
\end{eqnarray}
The epsilon factor $\epsilon(s,\mu_{\alpha})$ is defined as the product $\prod_v \epsilon(s, \mu_{\alpha,v} , \psi_{F_v})$, where $\mu_{\alpha,v}$ is the representation of $L_{F_v}$ defined by considering the pull-back of the representation $\mu_{\alpha}$ on $\mathcal{L}_{\psi}$ to a representation $\mu_{\alpha,v}$ on $L_{F_v}$ via the map $L_{F_v} \rightarrow \mathcal{L}_{\psi}$ of (2.4.20), and $\epsilon(s,\mu_{\alpha,v},\psi_{F_v})$ being the local epsilon factor as defined in \cite{T}, with respect to the localization $\psi_{F_v}$ of a non-trivial additive character $\psi_F: \mathbf{A}_F/F \rightarrow \mathbf{C}^{\times}$ (since $\mu_{\alpha,v}$ is symplectic the central value $\epsilon(1/2, \mu_{\alpha,v},\psi_{F_v})$ is independent of the choice of $\psi_{F_v}$). In this paper we will always use the Langlands normalization of local $L$ and $\epsilon$ factors (as discussed in \cite{T}).

Note that one has $\epsilon_{\psi} =1$ if $\psi=\phi$ is a generic parameter. 

We put
\begin{eqnarray}
\Pi_{\psi}(\epsilon_{\psi}) = \{  \pi \in \Pi_{\psi}, \,\ \langle \cdot , \pi  \rangle = \epsilon_{\psi}  \}.
\end{eqnarray}

Now we can state the global classification result. Denote by $\mathcal{H}(U_{E/F}(N))$ the global Hecke algebra of $U_{E/F}(N)(\mathbf{A}_F)$.

\begin{theorem}
Fix $\kappa \in \{\pm 1\}$ and $\chi_{\kappa} \in \mathcal{Z}_E^{\kappa}$. We have a $\mathcal{H}(U_{E/F}(N))$-module decomposition of the $L^2$-discrete automorphic spectrum of $U_{E/F}(N)(\mathbf{A}_F)$:
\begin{eqnarray}
& & L^2_{\disc}(U_{E/F}(N)(F)  \backslash U_{E/F}(N)(\mathbf{A}_F)) \\
&=& \bigoplus_{\psi \in \Psi_2(U_{E/F}(N),\xi_{\chi_{\kappa}})} \bigoplus_{\pi \in \Pi_{\psi}(\epsilon_{\psi})} \,\ \pi. \nonumber
\end{eqnarray}
\end{theorem}

One can of course take $\kappa=+1$, and $\chi_+=1$. But this is not necessary (and in fact in order to carry out the induction arguments, we need to include the case where $\chi_{\kappa}$ is non-trivial).

\begin{rem}
\end{rem}
By virtue of part (b) of theorem 2.5.1, we see that theorem 2.5.2 implies the multiplicity one result for representations $\pi$ that belong to a global packet corresponding to a generic parameter. 

\bigskip

The final theorem we state in this subsection gives crucial information about signs. 

In general let $\phi^N \in \Phi_{\simp}(N)$, i.e. a (unitary) cuspidal automorphic representation of $\GL_N(\mathbf{A}_E)$. We then have the Rankin-Selberg $L$-function as defined by Jacquet-Piatetskii-Shapiro-Shalika \cite{JPSS}:
\[
L(s,\phi^N \times (\phi^N)^c)
\]
which has a factorization:
\begin{eqnarray}
L(s,\phi^N \times (\phi^N)^c) = L(s,\phi^N,\Asai^+) \cdot L(s,\phi^N,\Asai^-)
\end{eqnarray}
here $L(s,\phi^N,\Asai^{\pm})$ are the two Asai $L$-functions of $\phi^N$. The local factors of these Asai $L$-functions are studied in Goldberg \cite{G}, which are special cases of the $L$-functions studied by Shahidi \cite{S}. Note that in \cite{G}, the function $L(s,\phi^N,\Asai)$ of {\it loc. cit.} is being denoted as the ``+" Asai $L$-function $L(s,\phi^N,\Asai^+)$ here, while the function $L(s, \phi^N \otimes \chi,\Asai)$ with $\chi \in \mathcal{Z}_E^-$ of {\it loc. cit.} (which is independent of the choice of $\chi \in \mathcal{Z}_E^-$) is being denoted as the ``--" Asai $L$-function $L(s,\phi^N,\Asai^-)$ here. In particular for $\chi_- \in \mathcal{Z}_E^-$, we have $L(s,\phi^N,\Asai^-) = L(s,\phi^N \otimes \chi_-,\Asai^+)$.

\bigskip

When $\phi^N \in \widetilde{\Phi}_{\simp}(N)$, namely $\phi^N$ being conjugate self-dual, we have:
\[
L(s,\phi^N \times (\phi^N)^c) = L(s,\phi^N \times (\phi^N)^{\vee})
\] 
hence has a simple pole at $s=1$ by [JPSS]. Furthermore by Shahidi's theorem \cite{S}, both functions $L(s,\phi^N,\Asai^{\pm})$ are non-zero at $s=1$. Hence we see that exactly one of the functions $L(s,\phi^N,\Asai^{\pm})$ has a pole at $s=1$, which is a simple pole.

\bigskip

Finally if $\phi_1^{N_1} \in \Phi_{\simp}(N_1),\phi_2^{N_2} \in \Phi_{\simp}(N_2)$ let $\epsilon(s,\phi_1^{N_1} \times \phi_2^{N_2})$ be the Rankin-Selberg $\epsilon$-factor associated to the pair $\phi_1^{N_1} , \phi_2^{N_2}$ \cite{JPSS}.

\begin{theorem} As before $F$ is global.

(a) Suppose $\phi^N \in \widetilde{\Phi}_{\simp}(N)$ is a conjugate self-dual simple generic parameter. Let $\kappa$ be the sign associated to $\phi^N$ as in theorem 2.4.2. Then the Asai $L$-function
\[
L(s, \phi^N,\Asai^{\eta})
\]
has a pole at $s=1$, where $\eta= (-1)^{N-1} \cdot \kappa$.

We say that $\phi^N$ is conjugate orthogonal (resp. conjugate symplectic) if $\eta=1$ (resp. $\eta=-1$).

(b) Suppose $\phi_1^{N_1} \in \widetilde{\Phi}_{\simp}(N_1)$, $\phi_2^{N_2} \in \widetilde{\Phi}_{\simp}(N_2)$. Assume that $\phi_1^{N_1}$ and $\phi_2^{N_2}$ are conjuagte self-dual of the same parity (i.e. both conjugate orthogonal or both conjugate symplectic). Then 
\begin{eqnarray}
\epsilon(1/2,\phi_1^{N_1} \times (\phi_2^{N_2})^c) =1.
\end{eqnarray}
\end{theorem}

\begin{rem}
\end{rem}
Part (a) of theorem 2.5.4 can thus be phrased in more common terms: given a unitary cuspidal automorphic representation of $\Pi$ of $\GL_N(\mathbf{A}_E)$ that is conjugate self-dual, we have $\Pi$ arises as ``standard base change" from $U_{E/F}(N)$ if and only if $L(s,\Pi,\Asai^{(-1)^{N-1}})$ has a pole at $s=1$ (with a similar result for the ``twisted base change"). Thus this is the global automorphic analogue of lemma 2.2.1. 

\bigskip

\begin{rem}
\end{rem}
Although theorem 2.5.2 can be stated without invoking theorem 2.5.4, nevertheless theorem 2.5.4 plays a critical role in the proof of theorem 2.5.2. Part (b) of theorem 2.5.4 is the global automorphic analogue of Proposition 5.2, part 2, of \cite{GGP}, where such a result is proved in {\it loc. cit.} for local root number associated to conjugate orthogonal representations of local Langlands group. 

\bigskip

\begin{rem}
\end{rem}
Although part (b) of theorem 2.5.4 concerns simple generic parameters {\it a priori}, the place where it enters in the induction argument actually concerns the case of {\it non-generic} parameters, namely in the comparison of trace formulas in section five and six. Accordingly we formulate part (b) of theorem 2.5.4, as a statement in terms of the integer $N$, as follows (the reason for this formulation will be made clear in section 5.8, also {\it c.f.} proposition 6.1.5):

\bigskip

{\it For $i=1,2$ let $\psi_i^{N_i} = \mu_i \boxtimes \nu_i \in \widetilde{\Psi}_{\simp}(N_i)$, with $N_1 + N_2 \leq N$. Assume that $\mu_1$ and $\mu_2$ are of the same parity, and that $\nu_1 \otimes \nu_2$ is a direct sum of an {\it odd} number of irreducible {\it even} dimensional representations of $\SL_2(\mathbf{C})$. Then we have:}
\begin{eqnarray}
\epsilon(1/2,\mu_1 \times \mu_2^c)=1.
\end{eqnarray}

\bigskip

\begin{example}
\end{example}
We give an important class of examples where the condition on order of pole of Asai $L$-function follows automatically from local condition at an archimedean prime. Suppose that 
$\psi^N =\phi^N \in \widetilde{\Phi}_{\simp}(N)$ is a simple generic parameter, i.e. $\psi^N$ is given by a (unitary) conjugate self-dual cuspidal automorphic representation $\phi^N$ of $\GL_N(\mathbf{A}_E)$. Suppose there exists an archimedean place $v$ of $F$, such that $F_v = \mathbf{R}$ and $E_v = \mathbf{C}$, and such that $\phi^N_{v}$ arise as the ``standard base change" of a parameter associated to a representation on the {\it compact} unitary group $U_{\mathbf{C}/\mathbf{R}}(N,0)$. Then we claim that $L(s,\phi^N,\Asai^{(-1)^{N-1}})$ has a pole at $s=1$. In other words, the sign $\kappa$ associated to $\phi^N$ is equal to $+1$ in this case, and the pair $\phi:=(\phi^N,\widetilde{\phi})$ defines a parameter in $\Phi_{\simp}(U_{E/F}(N),\xi_{\chi_{+}})$ (i.e. $\phi^N$ arises from $U_{E/F}(N)$ via standard base change).

Indeed suppose that we had $\kappa =-1$. Then by theorem 2.4,10, we have $\phi^N_{v}$ factors through $\xi_{\chi_{-}}$ for this particular archimedean place $v$. Thus $\phi^N_{v}$ is conjugate self-dual with parity $(-1)^{N-1} \cdot (-1) = (-1)^N$, by lemma 2.2.1. On the other hand, the given condition on $\phi^N_{v}$ implies that it is an elliptic parameter, i.e. $\phi^N_{v} \in \widetilde{\Phi}_{\ellip,v}(N)$. Thus the parity of $\phi^N_{v}$ is unique by remark 2.2.2. The assertion that $\phi^N_{v}$ is elliptic follows immediately from the fact that the $L$-parameter of such a $\phi^N_{v}$ is of the form (e.g. \cite{BC} prop. 4.3.2):
\begin{eqnarray}
& & \\ 
& & 
z \in L_{E_v} = W_{E_v}= \mathbf{C}^{\times} \mapsto \diag( (z/\overline{z})^{a_1},\cdots,(z/\overline{z})^{a_N} ) \nonumber
\end{eqnarray}
with $a_1 > \cdots > a_N$ and $a_i \in \frac{N+1}{2} + \mathbf{Z}$ for all $i$ (and the usual interpretation of $(z/\overline{z})^{1/2}$, i.e. $(z/\overline{z})^{1/2} := z/|z \overline{z}|^{1/2} $). 

Now take $w_c =j \in L_{F_v} \smallsetminus L_{E_v} = W_{F_v} \smallsetminus W_{E_v}$ such that $j^2 =-1$. Then from the condition that $a_i \in \frac{N+1}{2} + \mathbf{Z}$ we have, with $\phi^N_{v}$ as in the above form (2.5.12):
\[
\phi^N_{v}(j^2) = \phi^N_{v}(-1) = (-1)^{N-1} I_N
\]
From equation (2.2.6) we see that this implies $\phi^N_{v}$ is conjuagte self-dual of parity $(-1)^{N-1}$ (taking the matrix $A$ there to be the identity matrix), a contradiction. Thus we conclude that $\kappa =1$. We formulate this as:

\begin{corollary}
Suppose that $\Pi$ is a conjugate self-dual (unitary) cuspidal automorphic representation on $\GL_N(\mathbf{A}_E)$. Suppose there exists a place $v$ of $F$ such that $F_v=\mathbf{R}$, $E_v=\mathbf{C}$, and that the $L$-parameter corresponding to the irreducible admissible representation $\Pi_v$ of $\GL_N(E_v)$ is given by (2.5.12). Then 
\[
L(s,\Pi,\Asai^{(-1)^{N-1}})
\]
has a pole at $s=1$, and $\Pi$ descends to $U_{E/F}(N)$ via ``standard base change", and does not descend to $U_{E/F}(N)$ via ``twisted base change".
\end{corollary}

We note that special case of such descent result to unitary groups have of course been considered by many authors before in the case where $E/F$ is a CM extension of a totally real field $F$, and $\phi^N$ being conjugate self-dual cuspidal automorphic representation on $GL_N(\mathbf{A}_E)$, that is cohomological at all the archimedean places of $E$ ({\it c.f.} the discussion in the next section). 

\subsection{Review of earlier results}

The problem on classification of automorphic representations on unitary groups has of course been studied by many authors before. Besides the complete results by Rogawski \cite{R} in the case of unitary groups in three variables, there is the work of Harris-Labesse \cite{HL} on base change for unitary groups, under the condition that the cuspidal automorphic representation on the unitary group is supercuspidal at two finite places of $F$ that split in $E$. Results obtained using the converse theorem in the globally generic case were obtained by Kim-Krishnamurthy \cite{KK1,KK2} and Codgell-Piatetski-Shapiro-Shahidi \cite{CPSS}. There are also results based on the method of automorphic descent due to Ginzburg-Rallis-Soudry \cite{GRS}. We will in fact combine our results with that of Ginzburg-Rallis-Soudry to deduce the generic packet conjecture in section nine.  

Under the assumption that $E/F$ is a CM extension over a totally real field, special cases of endoscopic transfer for unitary groups has been studied by many authors: e.g. Clozel-Harris-Labesse \cite{CHL}, S.Morel \cite{Mo}, S.W.Shin \cite{Shi}. More complete results on endoscopic classification of discrete automorphic representations of unitary groups have been obtained by Labesse \cite{La}, under suitable local assumptions on the automorphic representations involved. Most of these local assumptions are of the form: the discrete automorphic representation on the unitary group belongs to discrete series at all archimedean places of $F$, and at finite places of $F$ that does not split in $E$ it is spherical. In particular we refer to the work \cite{STF} for discussion of these results.

We of course need to point out that as in \cite{A1}, the results obtained in this paper is still conditional on the stabilization of the twisted trace formula, whereas the results obtained by the earlier authors as discussed above are unconditional (under the suitable local assumptions). 

The construction of local $L$-packets for unitary groups associated to square-integrable parameters in the non-archimedean case is also obtained by Moeglin-Tadic \cite{MT} and Moeglin \cite{Moe}. We would like to mention that the work \cite{Moe} yields information about supercuspidal representations that is {\it a priori} not available in the general framework of endoscopic classification.

\section{\textbf{Local character identities and the intertwining relation}}

In this entire section $F$ is a local field, and $E$ a quadratic extension of $F$. We lay down some of the local preliminaries needed for the comparison of the trace formulas in later parts of the paper. In particular we state the local character identities characterizing the representations in a packet, and the local intertwining relation (whose full proof will be completed only in section 8), which plays a crucial role in both the comparison of trace formulas in section five and six, and in reducing the construction of packets for general parameters to the case of square-integrable parameters.

\subsection{Local endoscopic transfer of test functions}

The local character identities that characterize the representations in a packet is based on the endoscopic transfer of test functions, which are defined with respect to the transfer factors of Langlands-Shelstad-Kottwitz. We refer back to section 2.4 for discussion of endoscopic data, and in particular pertaining to the case of unitary groups.

We briefly recall some of the discussion as in section 2.1 of \cite{A1} (in a simplified form, which is sufficient for our purpose; we refer to {\it loc. cit.} for the general setting). Thus 
\[
G = (G^0,\theta)
\]
is the data of a twisted group defined over $F$, and for which we also denote by the same symbol $G$ for the twisted group $G = G^0 \rtimes \theta$. If $\theta$ is the identity then of course $G=G^0$ is just a group. Given an element $\gamma \in G(F)$, we denote by $G_{\gamma}$ the centralizer of $\gamma$ in $G^0$:
\[
G_{\gamma}=\{x \in G^0, \,\ x^{-1} \gamma x = \gamma \}.
\]
A semi-simple element $\gamma$ is called strongly regular, if $G_{\gamma}$ is an abelian group. Given a strongly regular element $\gamma$, we form the invariant orbital integral:
\[
f_{G}(\gamma) = |D(\gamma)|^{1/2} \int_{G_{\gamma}(F) \backslash G^0(F)} f(x^{-1} \gamma x) \,\ dx
\] 
with $D(\gamma)$ being the Weyl discriminant
\[
D(\gamma) = \det(1-\Ad(\gamma)|_{\mathfrak{g}^0/\mathfrak{g}_{\gamma}})
\] 
with $\mathfrak{g}^0,\mathfrak{g}_{\gamma}$ being the Lie algebras of $G^0, G_{\gamma}$ respectively. We denote 
\[
\mathcal{I}(G) = \{f_G, \,\ f \in \mathcal{H}(G)   \}
\]
the space of functions on the set of strongly regular $G^0(F)$-conjugacy classes in $G(F)$ spanned by the invariant orbital integrals. Here $\mathcal{H}(G)$ is the Hecke algebra (or module in the twisted case) of smooth, compactly supported functions on $G(F)$ defined with respect to a maximal compact subgroup of $G^0(F)$. For the case $G=\widetilde{G}_{E/F}(N) = G_{E/F}(N) \rtimes \theta$ (with $\theta$ defined as in (1.0.1)), which is the only twisted group that we will consider, we denote $\mathcal{H}(\widetilde{G}_{E/F}(N))$ as $\widetilde{\mathcal{H}}(N)$.

The spectral interpretation of the space of orbital integrals $\mathcal{I}(G)$ is given by characters. Thus let $(\pi^0,V)$ be an irreducible unitary representation of $G^0(F)$, with an unitary extension $\pi$ of $\pi^0$ to $G(F)$, i.e. $\pi$ is a function from $G(F)$ to the space of unitary operators on $V$ satisfying
\[
\pi(x_1 x x_2) = \pi^0(x_1) \pi(x) \pi^0(x_2), \,\ x_1,x_2 \in G^0(F), x \in G(F).
\] 
Then we have the character:
\[
f_G(\pi) = \tr \pi(f) = \tr \Big( \int_{G(F)} f(x) \pi(x)  dx \Big), \,\ f \in \mathcal{H}(G).
\]
The two functions $\{f_G(\gamma)\}$ and $\{f_G(\pi)\}$ determine each other.

Suppose now that $G=G^0$ is a (connected) quasi-split group. Given $f \in \mathcal{H}(G)$, we have the function $f^G$ on the set of strongly regular stable conjugacy classes of $G(F)$:
\[
f^G(\delta) :=\sum_{\gamma \rightarrow \delta} f_G(\gamma)
\] 
with the sum indicates the set of conjugacy classes of $G(F)$ that belongs to the stable conjugacy class $\delta$. The function $f^G$ is known as the stable orbital integral of $f$. Denote
\[
\mathcal{S}(G) = \{ f^G, \,\ f \in \mathcal{H}(G)      \}
\] 
the space spanned by the stable orbital integrals of elements in $\mathcal{H}(G)$. We thus have a map
\begin{eqnarray*}
\mathcal{H}(G) & \rightarrow & \mathcal{S}(G) \\
f &\rightarrow & f^G.
\end{eqnarray*}

Back to the general case where $G$ is allowed to be a twisted group. Let $G^{\prime} = (G^{\prime},s,\xi^{\prime})$ be an endoscopic data for $G$ as in section 2.4 (in particular $G^{\prime}$ is a quasi-split connected group). By Langlands-Shelstad \cite{LS} (untwisted case) and Kottwitz-Shelstad \cite{KS}, we have the notion of transfer factor $\Delta(\delta,\gamma)$ with respect to $G$ and and the endoscopic datum $G^{\prime}=(G^{\prime},s,\xi^{\prime})$; here $\delta$ is a strongly $G$-regular stable conjugacy class in $G^{\prime}(F)$, and $\gamma$ is a stongly regular $G^0(F)$-conjugacy class in $G(F)$. In general $\Delta$ is well-determined up to a complex multiplicative constant of absolute value one. We remark that as in Arthur's papers, the factor $\Delta_{IV}$ as in \cite{KS} is already absorbed into the definition of orbital integrals.

Given a definition of the transfer factor $\Delta$, we can then define the Langlands-Shelstad-Kottwitz transfer of $f$ to a function on the set of strongly $G$-regular stable conjugacy classes in $G^{\prime}(F)$:
\begin{eqnarray}
f^{G^{\prime}} (\delta) := \sum_{\gamma} \Delta(\delta,\gamma) f_G(\gamma), \,\ f \in \mathcal{H}(G)
\end{eqnarray}
(the sum runs over all the strongly regular $G^0(F)$-conjugacy classes in $G(F)$). We caution that in the notation for $f^{G^{\prime}}$ the symbol $G^{\prime}$ denotes an endoscopic datum, i.e. the transfer depends on $G^{\prime}$ as an endoscopic datum (not just as an endoscopic group).

The Langlands-Shelstad-Kottwitz transfer conjecture, which is now a theorem, asserts that for any $f \in \mathcal{H}(G)$, the transfer $f^{G^{\prime}}$ belongs to $\mathcal{S}(G^{\prime})$ for any endoscopic data $G^{\prime}$ of $G$. In other words we have:
\begin{eqnarray}
f^{G^{\prime}} = \rightexp{f^{\prime}}{G^{\prime}}
\end{eqnarray} 
for some $f^{\prime} \in \mathcal{H}(G^{\prime})$. Thus we have a map:
\begin{eqnarray}
\mathcal{H}(G) & \rightarrow & \mathcal{S}(G^{\prime}) \\
f &\rightarrow & f^{G^{\prime}}. \nonumber
\end{eqnarray}

We refer to section 2.1 of \cite{A1} for a discussion of the history of this conjecture. In the case when $F$ is archimedean, the transfer conjecture is established by Shelstad \cite{Sh1,Sh4}. In the non-archimedean case, the transfer conjecture is established as a consequence of the works of Waldspurger \cite{W1,W2,W3}, and the works of Ngo \cite{N} on the fundamental lemma (in the case of unitary groups, which is our main concern here, the fundamental lemma was established earlier in \cite{LN}). The fundamental lemma asserts that when both $G$ and $G^{\prime}$ are unramified, and if $f$ is taken to be the characteristic function of $K \rtimes \theta$, with $K$ being a $\theta$-stable hyperspecial maximal compact subgroup of $G^0(F)$, then $f^{\prime}$ in (3.1.2) can be taken to be the characteristic function of a hyperspecial maximal compact subgroup of $G^{\prime}(F)$. See also section 4.2 for related discussion in the global situation.

In this paper we will only be concerned with the case $G=\widetilde{G}_{E/F}(N)$, and the case where $G$ is a product of quasi-split unitary groups $U_{E/F}(N)$. Since these groups have standard $\Gamma_F$-splitting, with the splitting being $\theta$-stable in the twisted case $G=\widetilde{G}_{E/F}(N)$, we can use the Whittaker normalization of the transfer factor (section 5.3 of \cite{KS}), which will be the normalization of the local transfer factor used in this paper. 

We fix the standard $\theta$-stable $\Gamma_F$-splitting 
\[
S=(B,T,\{ x_{\alpha}\}_{\alpha \in \Delta^+})
\]
of the group $G^0$ (here $T$ is the standard diagonal maximal torus of $G^0$, $B$ is the standard upper triangular Borel subgroup of $G^0$, and $\Delta^+$ is the set of positive simple roots of $(B,T)$, which is not to be confused with the transfer factor). Then for an endoscopic datum $G^{\prime}$ of $G$, it is defined in \cite{KS} section 5.3, the transfer factor $\Delta_S$ normalized with respect to the splitting $S$ of $G^0$. For a given choice of non-trivial additive character $\psi_F: F \rightarrow \mathbf{C}^{\times}$, we then have the Whittaker normalized transfer factor $\Delta$, which is related to $\Delta_S$ by:
\begin{eqnarray*}
\Delta = \frac{\epsilon(1/2, \tau_G , \psi_F   )}{\epsilon(1/2, \tau_{G^{\prime}},\psi_F )}  \Delta_S
\end{eqnarray*}

\noindent here $\tau_G$ (resp. $\tau_{G^{\prime}}$) is the Artin representation of $\Gamma_F$ on the $\mathbf{C}$-vector space $X^*(T)^{\theta} \otimes_{\mathbf{Z}} \mathbf{C}$ (resp. on $X^*(T^{\prime}) \otimes_{\mathbf{Z}} \mathbf{C}$, with $T^{\prime}$ being the standard diagonal maximal trous of $G^{\prime}$). The $\epsilon$ factors are defined as in \cite{T} in the Langlands' normalization. The Whittaker normalized transfer factor depends only on the corresponding Whittaker data, i.e. on the data   
\[
(B,\lambda)
\]   
with $\lambda:N_B \rightarrow \mathbf{C}^{\times}$ ($N_B$ being the unipotent radical of $B$) being determined by $\psi_F$ and the splitting $S$ by:
\[
\lambda(n) =\psi_F ( \sum_{\alpha \in \Delta^+} \alpha(n) ).
\]

\bigskip

Recall as in section 2.4 that for $G=U_{E/F}(N)$, an endoscopic datum $(G^{\prime},s,\xi^{\prime})$ of $G$ is uniquely determined up to equivalence by the endoscopic group $G^{\prime}$. On the other hand, for the twisted group $G=\widetilde{G}_{E/F}(N)$, the endoscopic group alone does not determine (the equivalence class of) the endoscopic datum. Thus when we want to emphasize the $L$-embedding in the endoscopic datum we write $(G^{\prime},\xi^{\prime})$ for the endoscopic datum, instead of just $G^{\prime}$ alone. So for instance given $\widetilde{f} \in \widetilde{\mathcal{H}}(N) = \mathcal{H}(\widetilde{G}_{E/F}(N))$, and $(G^{\prime},\xi^{\prime}) \in \widetilde{\mathcal{E}}(N)$ an endoscopic datum of $\widetilde{G}_{E/F}(N)$, we denote $\widetilde{f}^{G^{\prime}}$ as $\widetilde{f}^{(G^{\prime},\xi^{\prime}) }$ to emphasize that it is the Kottwitz-Shelstad transfer of $\widetilde{f}$ to $G^{\prime}$ with respect to the datum $(G^{\prime},\xi^{\prime})$.

\bigskip

Finally we need a property of the Langlands-Kottwitz-Shelstad transfer which plays a crucial role in the analysis of the trace formulas. 

First briefly recall the notion of compatible family of functions (section 2.1 of \cite{A1}) for $\widetilde{\mathcal{E}}(N)$. Switching notation, instead of writing an element of $\widetilde{\mathcal{E}}(N)$ as $(G^{\prime},\xi^{\prime})$, we denote such an element as $(G,\xi)$. Now given $(G,\xi)$ and $(M,\xi_M) \in \widetilde{\mathcal{E}}(N)$, we say that $(M,\xi_M)$ is a {\it Levi subdata} of $(G,\xi)$, if $M$ is a Levi subgroup of $G$, and $\xi_M$ is up to $\widehat{G}_{E/F}(N)$-conjugacy equal to the composition of the $L$-embedding $\leftexp{L}{M} \rightarrow \leftexp{L}{G}$ (that is dual to the inclusion $M \hookrightarrow G$) with $\xi$. Then given $f \in \mathcal{H}(G)$, we have the descent map $f_G \rightarrow f_M$ on the invariant orbital integrals, characterized by the condition: for any admissible representation $\pi_M$ of $M(F)$, denote by $\mathcal{I}_P(\pi_M)$ the normalized parabolic induction of $\pi_M$ to $G(F)$ (here $P$ is a parabolic subgroup of $G$ with Levi component $M$), then we have
\[
f_M(\pi_M) = f_G(\mathcal{I}_P(\pi_M)).
\]
Similarly we have the descent map for the stable orbital integrals $f^{G} \rightarrow f^{M}$, with $f^{M} :=( f_M)^{M}$.

Now consider a family of function indexed by the twisted endoscopic data of $\widetilde{G}_{E/F}(N)$:
\begin{eqnarray}
\mathcal{F} = \{f \in \mathcal{H}(G): \,\  (G,\xi) \mbox{ a twisted endoscopic datum}   \}.
\end{eqnarray}
Again we emphasize that it is the endoscopic data $(G,\xi)$ that forms the indexing set, and not the endoscopic groups themselves. Then the family (3.1.4) is called a {\it compatible family}, if, firstly for any $f \in \mathcal{H}(G)$ associated to $(G,\xi) \in \widetilde{\mathcal{E}}(N)$, and $h \in \mathcal{H}(M)$ associated to $(M,\xi_M) \in \widetilde{\mathcal{E}}(N)$, with $(M,\xi_M)$ a Levi subdata of $(G,\xi)$, we have:
\begin{eqnarray*}
f^{M} = h^{M}
\end{eqnarray*}
and secondly, the family is to be compatible with equivalence of twisted endoscopic data, in the natural sense. We give the condition in the two most significant cases: the first case is as follows: let $G_1=(U_{E/F}(N),\xi_1)$ and $G_2=(U_{E/F}(N),\xi_2)$ with $\xi_1=\xi_{\chi_1}$ and $\xi_2=\xi_{\chi_2}$ such that $\chi_1,\chi_2 \in \mathcal{Z}_E^{\kappa}$. Thus $G_1,G_2$ are equivalent as twisted endoscopic data. The character $\widetilde{a}:=\chi_2/\chi_1$ restricts to the trivial character on $F^{\times}$, and hence descends to a character $a$ on $U_{E/F}(1)$, i.e. such that 
\[
\widetilde{a}(u)=a(u/c(u)).
\]
Then the compatibility condition is
\[
f_2^{G_2} =(a \circ \det)\cdot f_1^{G_1}
\]
here $\det : U_{E/F}(N) \mapsto U_{E/F}(1)$ is the determinant map, and $f_1$ (resp. $f_2$) is the function in the compatible family (3.1.4) associated to $G_1$ (resp. $G_2$). Here $f_i^{G_i}$ is of course just the stable orbital integral $f_i^{U_{E/F}(N)}$, but we have denoted it as $f_i^{G_i}$ for clarity.

The second case is as follows: let $N$ be even, and $G=(U_{E/F}(N),\xi)$ and $G^{\vee}=(U_{E/F}(N),\xi^{\vee})$ be representatives of the two equivalence classes of simple twisted endoscopic data of $\widetilde{G}_{E/F}(N)$, with $\xi=\xi_{\chi_+}, \xi^{\vee}=\xi_{\chi_-}$ such that $\chi_{\pm} \in \mathcal{Z}_E^{\pm}$. As in remark 2.4.1, let $L \cong G_{E/F}(N/2)$ be the Siegel Levi of $U_{E/F}(N)$. We equip $L$ with the $L$-embedding $\xi$ and regard $(L,\xi)$ as a (non-elliptic) twisted endoscopic datum of $\widetilde{G}_{E/F}(N)$, which we still denote by $L$ for the moment; similarly we equip the Levi subgroup $L$ with the $L$-embedding $\xi^{\vee}$ and denote $L^{\vee}=(L,\xi^{\vee})$ the resulting twisted endoscopic datum of $\widetilde{G}_{E/F}(N)$. The datum $L$ (resp. $L^{\vee}$ is a Levi sub-datum of $G$ (resp. $G^{\vee}$). 

The datum $L$ and $L^{\vee}$ are equivalent as elements of $\widetilde{\mathcal{E}}(N)$. For $\widetilde{f} \in \widetilde{\mathcal{H}}(N)$. Then the compatibility condition is 
\begin{eqnarray}
f^{\vee, L^{\vee}} =  ((\chi_-/\chi_+)^{N/2} \circ \det) \cdot f^L.
\end{eqnarray}
where $f$ (resp. $f^{\vee}$) is the function associated to the twisted endoscopic datum $G$ (resp. $G^{\vee}$) in the family (3.1.4). The formulation of the compatibility with equivalence of endoscopic data in general is similar. 

The condition (3.1.5) is motivated as follows: for $\widetilde{f} \in \widetilde{\mathcal{H}}(N)$, the definition of Kottwitz-Shelstad transfer factor gives the following equality:
\begin{eqnarray*}
\widetilde{f}^{L^{\vee}} =  ((\chi_-/\chi_+)^{N/2} \circ \det) \cdot \widetilde{f}^L
\end{eqnarray*}
{\it c.f.} (3.1.8) below.
\bigskip

We note that since any $(M,\xi_M)$ is a Levi subdata of a certain elliptic twisted endoscopic data $(G,\xi)$, it follows that for a compatible family of functions $\mathcal{F}$, the stable orbital integrals $h^{M}$ (for $h \in \mathcal{H}(M)$ are determined by the stable orbital integrals $f^{G}$ for a set of equivalence classes of elliptic twisted endoscopic data $(G,\xi)$. In the following we will take this convention; namely we understood $\widetilde{\mathcal{E}}(N)$ and $\widetilde{\mathcal{E}}_{\ellip}(N)$ as a set of representatives of the equivalence classes of twisted endoscopic data, and specifies only the functions in a compatible family associated to the set of representatives of elliptic twisted endoscopic data.

Recall the space of invariant orbital integrals on $\mathcal{I}(\widetilde{G}_{E/F}(N))$, which we denote as $\widetilde{\mathcal{I}}(N)$. We similarly define the endoscopic version of this space
\[
\widetilde{\mathcal{I}}^{\mathcal{E}}(N)  \subset \bigoplus_{(G,\xi)\in \widetilde{\mathcal{E}}(N)}  \mathcal{S}(G),
\]
\begin{eqnarray}
& & \\
& & \widetilde{\mathcal{I}}^{\mathcal{E}}(N) = \{\bigoplus_{(G,\xi) \in \widetilde{\mathcal{E}}(N)  } f^{G} \,\ \big| \,\ f \mbox{ belongs to a compatible family } \mathcal{F}  \}. \nonumber
\end{eqnarray}

It follows from the basic properties of the transfer factors, and the transfer conjecture of Langlands-Shelstad-Kottwitz, that the image of the collective transfer mapping:
\begin{eqnarray}
& & \widetilde{\mathcal{I}}(N) \rightarrow  \bigoplus_{(G,\xi)\in \widetilde{\mathcal{E}}(N)}  \mathcal{S}(G) \\
& & \widetilde{f} \rightarrow \bigoplus_{(G,\xi)\in \widetilde{\mathcal{E}}(N)} \widetilde{f}^{(G,\xi)} \nonumber
\end{eqnarray}
lies in $\widetilde{\mathcal{I}}^{\mathcal{E}}(N)$. 

\begin{proposition} Notations as above.

\bigskip

(a) The map $\widetilde{\mathcal{I}}(N) \rightarrow \widetilde{\mathcal{I}}^{\mathcal{E}}(N)$ defined by (3.1.7) is an isomorphism. In particular, a family of functions $\mathcal{F}$ in (3.1.4) forms a compatible family, if and only if there exists $\widetilde{f} \in \widetilde{\mathcal{H}}(N)$, such that 
\begin{eqnarray}
\widetilde{f}^{(G,\xi)} = f^G
\end{eqnarray}
for any $(G,\xi) \in \widetilde{\mathcal{E}}(N)$ and $f \in \mathcal{F}$ is the function associated to $(G,\xi)$.

\bigskip

(b) For $(G,\xi) \in \widetilde{\mathcal{E}}_{\simp}(N)$ (thus $G=U_{E/F}(N))$ the transfer mapping
\begin{eqnarray}
\widetilde{\mathcal{H}}(N) & \rightarrow & \mathcal{S}(U_{E/F}(N))  \\
\widetilde{f} & \rightarrow & \widetilde{f}^{(G,\xi)} \nonumber
\end{eqnarray}
is surjective.

\end{proposition}

\noindent Part (a) of Proposition 3.1.1 is exactly the content of Proposition 2.1.1 of \cite{A1} (which applies to the general setup of twisted endoscopy). Similarly, the argument in the proof of Corollary 2.1.2 of \cite{A1} shows that the transfer mapping 
\begin{eqnarray*}
\widetilde{\mathcal{H}}(N) & \rightarrow & \mathcal{S}(U_{E/F}(N))^{\widetilde{\Out}_N(U(N))}  \\
\widetilde{f} & \rightarrow & \widetilde{f}^{(G,\xi)} \nonumber
\end{eqnarray*}
is surjective. Hence we are done since $\widetilde{\Out}_N(U(N))$ is trivial.

\subsection{Characterization of the local classification}

In this subsection we give the precise statement of the local character relations for the representations in a packet, which in particular characterize the local classification in term of endoscopic transfer (both the standard and the twisted case). 

Recall that the classification is based on transfer of conjugate self-dual representations on $G_{E/F}(N)(F)=\GL_N(E)$. If $\pi$ is an irreducible conjugate self-dual representation of $G_{E/F}(N)(F)$, there are two different ways we can extend $\pi$ to the group $\widetilde{G}_{E/F}^+(N)(F) = G_{E/F}(N) \rtimes \langle \theta \rangle $. We need to choose the extension given by Whittaker normalization. The discussion follows section 2.2 of \cite{A1}.

Suppose first that $\pi$ is a tempered irreducible admissible representation of $G_{E/F}(N)(F)$. Then $\pi$ has a $(B(N),\lambda)$-Whittaker functional $\omega$, where $(B(N),\lambda)$ is the standard Whittaker data of $G_{E/F}(N)$, with $\lambda$ being the character on $N_{B(N)}$ as in section 3.1 defined with respect to a fixed non-trivial additive character $\psi_F$ of $F$. The Whittaker functional $\omega$ is unique up to scalar multiple, and we have
\[
\omega(\pi(n) v) = \lambda(n)  \omega(v), \,\ n \in N_{B(N)}(F), v \in V_{\infty}
\] 
where $V_{\infty}$ is the underlying vector space of smooth vectors of the representation $(\pi,V)$. Since we are assuming that $\pi$ is conjugate self-dual, the representation $\pi \circ \theta$ is equivalent to $\pi$, hence there exists a non-zero intertwining operator $I$, which is unique up to scalar, from $\pi$ to $\pi \circ \theta$.

Since the Whittaker data $(B(N),\lambda)$ is stable under $\theta$, the linear form $\omega \circ \theta$ is again a non-zero $(B(N),\lambda)$-Whittaker functional of $\pi$, hence we have
\[
\omega \circ I = c \omega
\]
for a non-zero constant $c$. We define:
\begin{eqnarray}
\widetilde{\pi}(\theta) := c^{-1} I
\end{eqnarray}
then $\widetilde{\pi}(\theta)$ is the unique intertwining operator from $\pi$ to $\pi \circ \theta$ such that
\begin{eqnarray}
\omega = \omega \circ \widetilde{\pi}(\theta).
\end{eqnarray}
The intertwining operator $\widetilde{\pi}(\theta)$ thus provides a unitary extension $\widetilde{\pi}$ of $\pi$ to $\widetilde{G}_{E/F}^+(N)$, in particular to the twisted group $\widetilde{G}_{E/F}(N)$, satisfying
\begin{eqnarray}
\widetilde{\pi}(x_1 x x_2 ) =  \pi(x_1) \widetilde{\pi}(x) \pi(x_2)
\end{eqnarray}
for $x_1,x_2 \in \widetilde{G}_{E/F}^0(N)(F) = G_{E/F}(N)(F) = \GL_N(E)$, and $x \in \widetilde{G}_{E/F}(N)(F)$.

Next, suppose $\pi$ is replaced by a conjugate self-dual standard representation $\rho$ of $G_{E/F}(N)(F)=\GL_N(E)$ (not necessarily irreducible). Then $\rho$ is of the form given as an induced representation:
\[
\rho = \mathcal{I}_P(\sigma), \,\ \sigma = \pi_{M,\lambda}.
\]
Here $M$ is identified as the standard diagonal Levi subgroup:
\[
G_{E/F}(m_1) \times \cdots \times G_{E/F}(m_k)
\]
of $G_{E/F}(N)$, $P$ being the standard parabolic subgroup of $G_{E/F}(N)$ containing $M$, and $\pi_M$ is an irreducible tempered representation of $M(F)$, given as 
\[
\pi_M = \pi_1 \boxtimes \cdots \boxtimes \pi_k
\]
for irreducible tempered representation $\pi_i$ of $G_{E/F}(m_i)(F)=GL_{m_i}(E)$. While $\lambda = (\lambda_1,\cdots,\lambda_k) \in \mathbf{R}^k$ with $\lambda_1 > \cdots > \lambda_k$, such that $\sigma=\pi_{M,\lambda}$ is given by the twist of $\pi_{M}$ by $\lambda$: 
\begin{eqnarray}
& & \sigma(x) = \pi_{M,\lambda}(x) \\
&=&  \pi_1(x_1) |\det x_1|^{\lambda_1} \boxtimes \cdots \boxtimes \pi_k(x_k) | \det x_k|^{\lambda_k}, \nonumber \\
& & x = (x_1,\cdots,x_k) \in M(F)     \nonumber
\end{eqnarray}

\noindent From the conjugate self-duality of $\rho$, it follows that:
\begin{eqnarray*}
\pi_i^* &\cong & \pi_{k+1-i},  \\
    - \lambda_i &=& \lambda_{k+1-i}                   \,\ 1 \leq i \leq k 
\end{eqnarray*}

\noindent In particular $M$ is stable under $\theta$, and that $\pi_M$ and $\sigma = \pi_{M,\lambda}$ are conjugate self-dual, i.e. $\pi_M$ is isomorphic to $\pi_M \circ \theta$, and similarly for $\sigma$. An argument similar to above applied to $\pi_M$ yields the Whittaker normalized  extension $\widetilde{\pi}_M$ of $\pi_M$ to the group:
\[
M^+(F) = M(F) \rtimes \langle  \theta \rangle
\] 
and hence also an extension $\widetilde{\sigma}$ of $\sigma$ to $M^+(F)$. Denote by $P^+$ the normalizer of $P$ in $\widetilde{G}_{E/F}^+(N) = G_{E/F}(N) \times \langle \theta \rangle$ (note that $P$ is stable under $\theta$). We then obtain the extension $\widetilde{\rho}$ of $\rho$ to $\widetilde{G}^+_{E/F}(N)(F)$ given by:
\[
\widetilde{\rho} = \mathcal{I}_{P^+} (\widetilde{\sigma})
\] 
In particular the extension $\widetilde{\rho}$ of $\rho$ to the twisted group $\widetilde{G}_{E/F}(N)(F)$.

Finally if $\pi$ is a general irreducible admissible conjugate self-dual representation of $G_{E/F}(N)(F)=\GL_N(E)$, then $\pi$ is the Langlands quotient of a unique standard representation $\rho$ of $\GL_N(E)$. Now we have a bijection between irreducible representations and standard representations of $\GL_N(E)$, and by the local Langlands classification for $\GL_N(E)$, these are in bijection with the parameters $\Phi(N)=\Phi(\GL_N(E))$. Thus we see that $\rho$ is also conjugate self-dual. Hence by the above discussion $\rho$ has a Whittaker normalized extension $\widetilde{\rho}$ to $\widetilde{G}_{E/F}^+(N)(F)$. We take $\widetilde{\pi}$ to be the Langlands quotient of $\widetilde{\rho}$, which we define to be the Whittaker normalized extension of $\pi$ to $\widetilde{G}^+_{E/F}(N)$ and hence to the twisted group $\widetilde{G}_{E/F}(N)$. Thus we obtain the extension of $\pi$ to $\widetilde{G}_{E/F}(N)(F)$ in all cases.  

\bigskip

Suppose that $\psi^N \in \widetilde{\Psi}(N)=\Psi(\widetilde{G}_{E/F}(N))$ is a conjugate self-dual parameter of $G_{E/F}(N)$. As in (2.2.11) we denote by $\phi_{\psi^N}$ the parameter in $\widetilde{\Phi}(N)$ associated to $\psi^N$:
\begin{eqnarray}
\phi_{\psi^N}(\sigma) = \psi^N\Big(\sigma,   \begin{pmatrix}   |\sigma|^{1/2}  &  0 \\    0  &     |\sigma|^{-1/2}     \\ \end{pmatrix}  \Big), \mbox{ for } \sigma \in L_{F}.
\end{eqnarray}
Let $\rho_{\psi^N}= \rho_{\phi_{\psi^N}}$ be the standard representation of $G_{E/F}(N)(F)=\GL_N(E)$ associated to the $L$-parameter $\phi_{\psi^N}$, with Langlands quotient $\pi_{\psi^N}=\pi_{\phi_{\psi^N}}$. Then $\pi_{\psi^N}$ is an irreducible conjugate self-dual representation of $G_{E/F}(N)(F)=\GL_N(E)$. Furthermore $\pi_{\psi^N}$ is unitary, by the classification of irreducible unitary representations of general linear groups (Tadi\'c \cite{Ta} in the non-archimedean case and Vogan \cite{V} in the archimedean case). We have the linear form on $\widetilde{\mathcal{H}}(N)$:
\begin{eqnarray}
& & \widetilde{f} \mapsto \widetilde{f}_N(\psi^N), \,\ \widetilde{f} \in \widetilde{\mathcal{H}}(N)  \\
& & \widetilde{f}_N(\psi^N) := \tr \widetilde{\pi}_{\psi^N}(\widetilde{f}),  \nonumber
\end{eqnarray}
with $\widetilde{\pi}_{\psi^N}$ the canonical extension of $\pi_{\psi^N}$ as above. 

\bigskip
Recall from section 1.4 of \cite{A1} that we have a bijection:
\begin{eqnarray*}
(G^{\prime},\psi^{\prime}) \longleftrightarrow (\psi,s)
\end{eqnarray*}
(here we have taken $G$ to be $U_{E/F}(N)$ for simplicity; the bijection has a obvious variant for the twisted group $G=\widetilde{G}_{E/F}(N)$, for example in the analysis in section five). In this bijection the left hand side consists of pairs $(G^{\prime},\psi^{\prime})$, where $G^{\prime}=(G^{\prime},s^{\prime},\xi^{\prime})$ is an endoscopic datum of $G$. Here two such data $(G_1^{\prime},s_1^{\prime},\xi_1^{\prime}),(G_2^{\prime},s_2^{\prime},\xi_2^{\prime})$ are identified if the images $\xi_i^{\prime}(\leftexp{L}{G_i^{\prime}})$ of $\leftexp{L}{G_i^{\prime}}$ in $\leftexp{L}{G}$ are the same (for $i=1,2$), and $s^{\prime}_1,s^{\prime}_2$ are equal up to translation by an element in $Z(\widehat{G})^{\Gamma_F}$ (the identification relation is thus stronger than the equivalence relation on endoscopic data), and $\psi^{\prime}$ is an {\it actual} $L$-homomorphism $\psi^{\prime} : L_F \times \SU_2 \rightarrow \leftexp{L}{G}$ with image contained in $\xi^{\prime}(\leftexp{L}{G^{\prime}})$; while the right hand side consists of pairs where $\psi$ is an {\it actual} homomorphism $\psi:L_F \times \SU_2 \rightarrow \leftexp{L}{G}$, and $s$ is a semi-simple element of $\overline{S}_{\psi}$. We briefly recall this correspondence. Thus given the pair $(G^{\prime},\psi^{\prime})$, we take $\psi$ to be the $L$-homomorphism obtained by composing $\psi^{\prime}$ with the $L$-embedding $\xi^{\prime}:\leftexp{L}{G^{\prime}} \hookrightarrow \leftexp{L}{G}$, and $s$ being the image of the element $s^{\prime}$ associated to the endoscopic datum $G^{\prime}$ in the quotient $\overline{S}_{\psi}=S_{\psi}/Z(\widehat{G})^{\Gamma_F}$ (the element $s^{\prime}$ lies in $S_{\psi}$ precisely because $\psi$ factors through $\xi^{\prime}: \leftexp{L}{G^{\prime}} \rightarrow \leftexp{L}{G}$). 

Conversely given the pair $(\psi,s)$, the pair $(G^{\prime},\psi^{\prime})$ is given as follows. Put
\[
\widehat{G}^{\prime} := \Cent(s,\widehat{G})^0
\]
and the $L$-action of $W_F$ on $\widehat{G}^{\prime}$ being determined by $\psi|_{W_F}$. This determines the $L$-group $\leftexp{L}{G^{\prime}}$, and the endoscopic datum $G^{\prime}=(G^{\prime},s^{\prime},\xi^{\prime})$. By construction the $L$-homomorphism $\psi$ factors through the image of $\xi^{\prime}(\leftexp{L}{G^{\prime}})$ in $\leftexp{L}{G}$, which we take as $\psi^{\prime}$.

It follows that given a pair $(\psi,s)$ with $\psi \in \Psi(G)$, and $s \in \overline{S}_{\psi}$, we can choose an endoscopic datum $G^{\prime} =(G^{\prime},s^{\prime},\xi^{\prime})$ (i.e. an actual datum and not just as an equivalence class), and $\psi^{\prime} \in \Psi(G^{\prime})$, such that $\psi=\xi \circ \psi^{\prime}$. A more precise study of this correspondence is carried out in section five, also {\it c.f.} section 4.8 of \cite{A1}).

We can now state the main result concerning the character relations in the local packets, which made precise the statement of theorem 2.5.1.

\begin{theorem}
(a) Suppose that $G=G_1 \times G_2 = U_{E/F}(N_1) \times U_{E/F}(N_2)$ with $N_1+N_2=N$, and that $\psi \in \Psi(G)$. Then there is a unique stable linear form on $\mathcal{H}(G)$:
\begin{eqnarray}
f \mapsto f^G(\psi), \,\ f \in \mathcal{H}(G)
\end{eqnarray}
with the following property: for any $L$-embedding $\xi:\leftexp{L}{G} \rightarrow \leftexp{L}{G_{E/F}(N)}$ defining an endoscopic datum $(G,\xi)$ whose equivalence class is in $\widetilde{\mathcal{E}}_{\ellip}(N)$, we have
\begin{eqnarray}
\widetilde{f}^{(G,\xi)}(\psi) = \widetilde{f}_N( \xi_* \psi), \,\ \widetilde{f} \in \mathcal{H}(N),
\end{eqnarray}
together with a secondary property:
\begin{eqnarray}
f^G(\psi) = f^{G_1}(\psi_{1}) f^{G_2}(\psi_{2})
\end{eqnarray}
in case
\begin{eqnarray*}
 \psi = \psi_1 \times \psi_{2}, \,\ \psi_{i} \in \Psi(G_i),
\end{eqnarray*}
and
\[
f^G = f^{G_1} \times f^{G_2}, \,\ f^{G_i} \in \mathcal{S}(G_i)
\]
are composite.

(b) For every $\psi \in \Psi(U_{E/F}(N))$, there exists a finite multi-set $\Pi_{\psi}$ (i.e. a set with multiplicities), whose elements are irreducible unitary representations, together with a mapping:
\begin{eqnarray}
\Pi_{\psi} & \rightarrow & \widehat{\mathcal{S}}_{\psi} \\
\pi & \mapsto & \langle \cdot, \pi \rangle \nonumber
\end{eqnarray}
with the following property. If $s$ is a semi-simple element of the centralizer $S_{\psi}=S_{\psi}(U_{E/F}(N))$, and $(G^{\prime},\psi^{\prime})$ is the pair corresponding to $(\psi,s)$ as described above (here $G^{\prime}$ is understood as an endoscopic datum of $G$), then we have:
\begin{eqnarray}
& & \\
& & f^{G^{\prime}}(\psi^{\prime}) = \sum_{\pi \in \Pi_{\psi}} \langle s_{\psi} x,\pi\rangle f_{U(N)}(\pi), \,\ f \in \mathcal{H}(U_{E/F}(N)) \nonumber
\end{eqnarray}
where $x$ is the image of $s$ in $\mathcal{S}_{\psi}$, and $s_{\psi}$ is the element defined as in (2.2.12).

\end{theorem}

\bigskip

\begin{rem}
\end{rem}
From proposition 3.1.1(b), we have that the Kottwitz-Shelstad transfer mapping:
\begin{eqnarray}
& & \widetilde{\mathcal{H}}(N) \rightarrow \mathcal{S}(U_{E/F}(N)) \\
& & \widetilde{f} \mapsto \widetilde{f}^{(U(N),\xi)} \nonumber
\end{eqnarray}
is surjective for any $\xi:\leftexp{L}{U_{E/F}(N)} \hookrightarrow \leftexp{L}{G_{E/F}(N)}$. Hence equation (3.2.8) uniquely characterizes the stable linear form $f \mapsto f^{G}(\psi)$ on $\mathcal{H}(G)$ when $G=U_{E/F}(N)$. But we note that the stable linear form $f \mapsto f^G(\psi)$ depends only on $G$ and $\psi$, and {\it not} on $G$ as a twisted endoscopic datum of $\widetilde{G}_{E/F}(N)$. The secondary property (3.2.9) then characterizes the linear form $f^G(\psi)$ in the composite case.

\bigskip

\begin{rem}
\end{rem}
Theorem 3.2.1 make precise the local classification theorem stated as theorem 2.5.1. We note that part (c) of theorem 2.5.1 concerning central characters follows from part (a) of theorem 3.2.1. Indeed it is clear that for the assertion on central character given by theorem 2.5.1(c), it suffices to treat the case where $\xi=\xi_{\triv}: \leftexp{L}{G} \hookrightarrow \leftexp{L}{G_{E/F}(N)}$ is the ``standard base change" $L$-embedding (i.e. $\chi \in \mathcal{Z}_E^+$ is the trivial character). Then the assertion on central character follows from (3.2.8), by considering test functions on $G(F)$ and $\widetilde{G}_{E/F}(N)$ equivariant with respect to centre of $G$, together with the basic property of Kottwitz-Shelstad transfer mapping (for example \cite{KS} p.53 and p.112).

As with the main theorems in section 2.5, theorem 3.2.1 is proved by a long induction argument to be completed in section 8. The endoscopic character relations stated in part (b) of the theorem uniquely characterize the packet $\Pi_{\psi}$ for $\psi \in \Psi(U_{E/F}(N))$.

\bigskip

In order to establish theorem 3.2.1, it is necessary to augment it with the local intertwining relation. We state this in section 3.4, after some preparations on normalization of intertwining operators in the next subsection.

\subsection{Normalization of local intertwining operators}

As in \cite{A1}, the intertwining operators of induced representations play an important role in the analysis of the spectral terms of the trace formula. Following Arthur there will be three steps in the normalization of the local intertwining operators. Steps one and two will be done in this section. The first step is to define normalization factors for the local intertwining integrals, in terms of local $L$ and $\epsilon$ factors.

Thus we denote $G=U_{E/F}(N)$. Again it is convenient to allow the case where $E$ is split, i.e. $E=F \times F$. For the most part the discussion in the split case is simpler, since $U_{E/F}(N) \cong \GL_{N/F}$ when $E = F \times F$, and the essential case is the case where $E$ is a field extension of $F$.

For a fixed proper Levi subgroup $M$ of $G$, we denote by $\mathcal{P}(M)$ the set of parabolic subgroups of $G$ defined over $F$ with Levi component $M$. If $P \in \mathcal{P}(M)$ we denote by $N_P$ the unipotent radical of $P$.

Now $M$ has a decomposition of the form:
\begin{eqnarray}
M \cong G_{E/F}(N_1^{\prime}) \times \cdots \times G_{E/F}(N^{\prime}_{r^{\prime}}) \times G_-
\end{eqnarray}
with $G_- \cong U_{E/F}(N_-)$ for $N_- < N$, such that
\begin{eqnarray}
2 N_1^{\prime} + \cdots + 2 N_{r^{\prime}}^{\prime} + N_- = N.
\end{eqnarray}

In this section we are concerned with representations induced from $M(F)$. Thus {\it for this section only}, the symbols $\phi,\psi,\pi$ will be reserved with respect to $M$ instead of $G$. 

First suppose that $\phi \in \Phi(M)$ is a generic parameter. In general both the unnormalized local intertwining operators and the normalizing factors can have poles at particular points. To circumvent this we consider twists of $\phi$. More precisely, let $\lambda$ be a point in general position in the vector space:
\begin{eqnarray*}
 \mathfrak{a}^*_{M,\mathbf{C}} = \mathfrak{a}^*_{M} \otimes_{\mathbf{R}} \mathbf{C} 
\end{eqnarray*}
where
\[
\mathfrak{a}^*_M = X^*(M)_F \otimes_{\mathbf{Z}} \mathbf{R}
\]
we consider the twist 
\[
\phi_{\lambda}(w) = \phi(w) |w|^{\lambda}, \,\ w \in L_F
\]
here $|w|^{\lambda}$ is the point in
\[
A_{\widehat{M}} = (Z(\widehat{M})^{\Gamma_F})^0 = X^*(M)_F \otimes_{\mathbf{Z}} \mathbf{C}^{\times}
\]
associated to $|w|$ and $\lambda \in \mathfrak{a}^*_{M,\mathbf{C}} = X*(M)_F \otimes_{\mathbf{Z}} \mathbf{C}$. For instance if $\lambda \in X^*_F(M)$, then
\[
|w|^{\lambda} = \lambda^{\vee}(|w|)
\]
where $\lambda^{\vee} \in X_*(Z(\widehat{M}))$ is the element that correspond to $\lambda$ under the correspondence $X^*(M) \cong X_*(Z(\widehat{M}))$.

\noindent The map $w \mapsto |w|^{\lambda}$ is a central Langlands parameter for $M(F)$ that is the Langlands dual to the unramified character of $M(F)$:
\[
m \mapsto \exp(\lambda H_M(m)), \,\ m \in M(F).
\]
with $H_M: M(F) \rightarrow \mathfrak{a}_M = \Hom_{\mathbf{Z}}(X^*_F(M),\mathbf{R})$ being defined in the usual manner:
\[
\langle  H_M(m), \chi \rangle = \log |\chi(m)|, \,\ m \in M(F), \chi \in X^*_F(M).
\]

\bigskip

Given $P^{\prime},P \in \mathcal{P}(M)$, write $\rho_{P^{\prime}|P}$ for the adjoint representation of $\leftexp{L}{M}$ on the quotient:
\[
\widehat{\mathfrak{n}}_{P^{\prime}} / \widehat{\mathfrak{n}}_{P^{\prime}} \cap \widehat{\mathfrak{n}}_P
\]
with $\widehat{\mathfrak{n}}_P$ the Lie algebra of the unipotent radical of $\widehat{P}$. Denote by $\rho^{\vee}_{P^{\prime}|P}$ the contragredient of $\rho_{P^{\prime}|P}$. Then $\rho^{\vee}_{P^{\prime}|P} \circ \phi_{\lambda}$ is a finite dimensional representation of $L_F$. We then have the Artin local $L$ and $\epsilon$ factors defined as in \cite{T} (always with respect to the Langlands normalization):
\begin{eqnarray*}
& & L(s,\rho^{\vee}_{P^{\prime}|P} \circ \phi_{\lambda}) \\
& & \epsilon(s, \rho^{\vee}_{P^{\prime}|P} \circ \phi_{\lambda} ,\psi_F).
\end{eqnarray*}

\noindent We also define the $\delta$-factor

\begin{eqnarray*}
\delta(\rho^{\vee}_{P^{\prime}|P} \circ \phi_{\lambda},\psi_F) := \epsilon(1/2, \rho^{\vee}_{P^{\prime}|P} \circ \phi_{\lambda} ,\psi_F)^{-1} \epsilon(0, \rho^{\vee}_{P^{\prime}|P} \circ \phi_{\lambda} ,\psi_F).
\end{eqnarray*}

\noindent We then define the local normalizing factor 
\[
r_{P^{\prime}|P}(\phi_{\lambda}) = r_{P^{\prime}|P}(\phi_{\lambda},\psi_F)
\]
given by:
\begin{eqnarray}
& & \\
& & r_{P^{\prime}|P}(\phi_{\lambda},\psi_F)=
 \delta( \rho^{\vee}_{P^{\prime}|P} \circ \phi_{\lambda},\psi_F)^{-1} L(0, \rho^{\vee}_{P^{\prime}|P} \circ \phi_{\lambda})  L(1, \rho^{\vee}_{P^{\prime}|P} \circ \phi_{\lambda})^{-1}. \nonumber
\end{eqnarray}
Then $r_{P^{\prime}|P}(\phi_{\lambda})$ is a meromorphic function of $\lambda \in \mathfrak{a}^*_{M,\mathbf{C}}$.

\bigskip

\noindent Next suppose that $\psi \in \Psi(M)$. As in (2.2.11) we have the $L$-parameter $\phi_{\psi} \in \Phi(M)$ atatched to $\psi$. For $P^{\prime},P \in \mathcal{P}(M)$, put
\begin{eqnarray}
r_{P^{\prime}|P}(\psi_{\lambda})=r_{P^{\prime}|P}(\psi_{\lambda},\psi_F) := r_{P^{\prime}|P}(\phi_{\psi,\lambda},\psi_F).
\end{eqnarray}

\bigskip

Next recall the unnormalized local intertwining operators ({\it c.f.} \cite{A7} section 1). Thus let $\pi=(\pi,V)$ be an irreducible unitary representation of $M(F)$. For $\lambda \in \mathfrak{a}^*_{M,\mathbf{C}}$ we have the twist:
\[
\pi_{\lambda}(m) = \pi(m) \exp(\lambda H_M(m)).
\]
 
\noindent For $P \in \mathcal{P}(M)$ we have the parabolically induced representation $\mathcal{I}_P(\pi_{\lambda})$. We identify the underlying space of $\mathcal{I}_P(\pi_{\lambda})$ as the underlying space $\mathcal{H}_P(\pi)$ of $\mathcal{I}_P(\pi)$: namely fix a maximal compact subgroup $K$ of $G(F)$ in good position relative to $M$. Then we take $\mathcal{H}_P(\pi)$ to be the Hilbert space of functions $\phi:K \rightarrow V$, satisfying
\[
\phi(nmk) = \pi(m) \phi(k)
\]
for $k \in K$, $n \in K \cap N_P(F)$, $m \in K \cap M(F)$ (in particular the space $\mathcal{H}(\pi)$ is unchanged when we twist $\pi$ by $\lambda \in \mathfrak{a}^*_{M,\mathbf{C}}$). Then with repsect to $\mathcal{I}_P(\pi_{\lambda})$, such a function $\phi$ extends to a function on $G(F)$ by the rule: 
\[
\phi(x) = \pi(M_P(x)) \phi(K_P(x)) \cdot e^{ (\lambda + \rho_P) H_P(x) }
\]
where for $x \in G(F)$ we have written $x=N_P(x) M_P(x) K_P(x)$ as in {\it loc. cit.} As usual $\rho_P$ is half the sum of roots of the pair $(P,A_M)$ (with multiplicities), and $H_P:G(F) \rightarrow \mathfrak{a}_M$ is defined by:
\[
H_P(nmk) = H_M(m)
\]
for $n \in N_P(F), m \in M(F), k \in K$.

\bigskip

Then for $P^{\prime},P \in \mathcal{P}(M)$ we have the unnormalized local intertwining operator:
\[
J_{P^{\prime}|P}(\pi_{\lambda}) : \mathcal{H}_{P}(\pi) \longrightarrow \mathcal{H}_{P^{\prime}}(\pi)
\]
that intertwines the induced representations $\mathcal{I}_P(\pi_{\lambda})$ and $\mathcal{I}_{P^{\prime}}(\pi_{\lambda})$, which is defined for $Re(\lambda)$ in a certain sector of $\mathfrak{a}_M^*$ by an absolutely convergent integral:
\begin{eqnarray}
& & (J_{P^{\prime}|P}(\pi_{\lambda}))(h)(k) \\
&=& \int_{N_{P^{\prime}}(F) \cap N_P(F) \backslash N_{P^{\prime}}(F)} \pi(M_P(y) ) h( K_P(y) k ) e^{ (\lambda +\rho_P)(H_P(y) ) } \,\ dy \nonumber
\end{eqnarray}
for $h \in \mathcal{H}_P(\pi)$, $k \in K$, and we have written $y = N_P(y) M_P(y) K_P(y)$. The measure $dy$ is defined with respect to the Haar measure on $F$ that is self-dual with respect to the additive character $\psi_F$ ({\it c.f.} the discussion after equation (2.3.7) of \cite{A1}). 

For general $\lambda \in \mathfrak{a}_{M,\mathbf{C}}^*$ the operator $J_{P^{\prime}|P}(\pi_{\lambda})$ is obtained as meromorphic continuation.

\bigskip

The proof of the main theorems in this paper will be by induction on the integer $N$. The induction hypothesis will be made formally in section 5.3. In our present situation, since $M$ is a proper Levi subgroup of $G$, we can assume that the main theorems, both local and global, hold for $G_-=U_{E/F}(N_-)$. Since the local and global theorems are already known for the factors $G_{E/F}(N^{\prime}_i)(F) = \GL_{N^{\prime}_i}(E)$, we can thus assume that the main results hold for $M$. Thus given $\psi \in \Psi(M)$, we have the local packet $\Pi_{\psi}$, which in particular consists of irreducible unitary representations of $M(F)$. Define, for $P^{\prime},P \in \mathcal{P}(M)$ and $\lambda \in \mathfrak{a}^*_{M,\mathbf{C}}$, and $\pi \in \Pi_{\psi}$:
\begin{eqnarray}
R_{P^{\prime}|P}(\pi_{\lambda},\psi_{\lambda}) = r_{P^{\prime}|P}(\psi_{\lambda})^{-1} J_{P^{\prime}|P}(\pi_{\lambda}).
\end{eqnarray}
It is independent of $\psi_F$, as the dependence of $r_{P^{\prime}|P}(\psi_{\lambda})$ on $\psi_F$ is the same as that of $J_{P^{\prime}|P}(\pi_{\lambda})$.

\bigskip

In this subsection we fix an $L$-embedding $\xi =\xi_{\chi} : \leftexp{L}{G} = \leftexp{L}{U_{E/F}(N)} \hookrightarrow \leftexp{L}{G_{E/F}(N)}$, with $\chi \in \mathcal{Z}_E$, and we consider $G = (G,\xi)$ as an endoscopic datum in $\widetilde{\mathcal{E}}_{\simp}(N)$ (in the split case where $E=F \times F$ the choice of $\chi$ is of course unnecessary). We continue to denote by $\xi$ the $L$-embedding $\leftexp{L}{G_-} \hookrightarrow \leftexp{L}{G_{E/F}(N_-)}$ given by the restriction of $\xi$ to $\leftexp{L}{G_-}$, and similarly regard $G_- = (G_-,\xi)$ as an element of $\widetilde{\mathcal{E}}_{\simp}(N_-)$. 

\bigskip

\noindent With respect to the decomposition (3.3.1) put
\begin{eqnarray*}
& & \widetilde{M}^0 =  \\
& & G_{E/F}(N_1^{\prime}) \times \cdots \times G_{E/F}(N_r^{\prime}) \times G_{E/F}(N_-) \times G_{E/F}(N_r^{\prime}) \times \cdots \times G_{E/F}(N_1^{\prime}).
\end{eqnarray*}
Then $\widetilde{M}^0$ is a Levi subgroup of $G_{E/F}(N)$, and $\xi : \leftexp{L}{G} \hookrightarrow \leftexp{L}{G_{E/F}(N)}$ restricts to $\xi : \leftexp{L}{M} \hookrightarrow \leftexp{L}{\widetilde{M}^0}$ (with $\leftexp{L}{\widetilde{M}^0} \hookrightarrow \leftexp{L}{G_{E/F}(N)}$ being the Levi embedding dual to $\widetilde{M}^0 \hookrightarrow G_{E/F}(N)$). Then we again regard $M=(M,\xi)$ as an endoscopic datum (whose equivalence class defines an element in $\widetilde{\mathcal{E}}(N)$).  

\bigskip

\noindent We then note the following: the representation $\rho^{\vee}_{P^{\prime}|P}$ of $\leftexp{L}{M}$ factors as 
\[
\rho^{\vee}_{P^{\prime}|P} = \widetilde{\rho}^{\vee}_{P^{\prime} |P} \circ \xi
\]
here $\widetilde{\rho}^{\vee}_{P^{\prime} |P}$ is a representation of $\leftexp{L}{\widetilde{M}^0}$ whose simple constituents are given by Rankin-Selberg and $\Asai^{\pm}$ representations of the general linear factors of $\leftexp{L}{\widetilde{M}^0}$.

Now suppose that $ \phi \in \Phi(M)$ is a generic parameter. Consider the parameter $\phi^{\widetilde{M}^0} := \xi_* \phi \in \widetilde{\Phi}(\widetilde{M}^0)$. Then the representation $\rho^{\vee}_{P^{\prime}|P} \circ \phi_{\lambda}$ can be identified as $\widetilde{\rho}^{\vee}_{P^{\prime}|P} \circ \phi^{\widetilde{M}^0}_{\lambda}$. 

Thus we have 
\begin{eqnarray*}
& & L(s, \rho^{\vee}_{P^{\prime}|P} \circ \phi_{\lambda}) = L(s, \widetilde{\rho}^{\vee}_{P^{\prime}|P}  \circ \phi^{\widetilde{M}^0}_{\lambda} ) \\
& &  \epsilon(s, \rho^{\vee}_{P^{\prime}|P} \circ \phi_{\lambda},\psi_F) = \epsilon(s, \widetilde{\rho}^{\vee}_{P^{\prime}|P}  \circ \phi^{\widetilde{M}^0}_{\lambda} ,\psi_F). 
\end{eqnarray*}
On the other hand, if $\pi_{\phi^{\widetilde{M}^0}}$ is the irreducible admissible representation of $G_{E/F}(N)(F)=\GL_N(E)$ that corresponds to $\phi^{\widetilde{M}^0}$ under the local Langlands classification, then we also have the representation theoretic $L$ and $\epsilon$-factors:
\begin{eqnarray*}
& &  L(s,\pi_{\phi^{\widetilde{M}^0}_{\lambda}} ,\widetilde{\rho}^{\vee}_{P^{\prime}|P}  ) \\
& &   \epsilon(s, \pi_{\phi^{\widetilde{M}^0}_{\lambda}},\widetilde{\rho}^{\vee}_{P^{\prime}|P}  ,\psi_F) 
\end{eqnarray*}
which are given by a product of Rankin-Selberg factors on \linebreak $G_{E/F}(N^{\prime})(F) \times G_{E/F}(N^{\prime \prime})(F)$ (with $N^{\prime},N^{\prime \prime} \leq N$) and $\Asai^{\pm}$ factors on $G_{E/F}(N^{\prime})(F)$, as studied in Goldberg \cite{G}, and are special case of the $L$ and $\epsilon$ factors defined by Shahidi \cite{S}. For the Rankin-Selberg constituent of $\widetilde{\rho}^{\vee}_{P^{\prime}|P}$ the Artin and representation theoretic $L$ and $\epsilon$ factors are equal, by the local Langlands classification of \cite{HT,H1}. On the other hand, for the $\Asai^{\pm}$ constituent of $\widetilde{\rho}^{\vee}_{P^{\prime}|P}$, the result of Henniart \cite{H2} gives the equality of the Artin and representation theoretic $\Asai^{\pm}$ $L$-factors, while the Artin and representation theoretic $\Asai^{\pm}$ $\epsilon$-factors are equal up to a constant that is a root of unity. Thus we have:
\begin{eqnarray*}
& &  L(s, \rho^{\vee}_{P^{\prime}|P} \circ \phi_{\lambda}) = L(s, \widetilde{\rho}^{\vee}_{P^{\prime}|P}  \circ \phi^{\widetilde{M}^0}_{\lambda} )=L(s,\pi_{\phi^{\widetilde{M}^0}_{\lambda}} ,\widetilde{\rho}^{\vee}_{P^{\prime}|P}  )  \\
& & \epsilon(s, \rho^{\vee}_{P^{\prime}|P} \circ \phi_{\lambda},\psi_F) = \epsilon(s, \widetilde{\rho}^{\vee}_{P^{\prime}|P}  \circ \phi^{\widetilde{M}^0}_{\lambda} ,\psi_F) = \zeta  \cdot \epsilon(s, \pi_{\phi^{\widetilde{M}^0}_{\lambda}},\widetilde{\rho}^{\vee}_{P^{\prime}|P}  ,\psi_F) 
\end{eqnarray*}
for $\zeta$ a constant that is a root of unity. In particular the Artin $\delta$-factor 
\begin{eqnarray*}
& & \delta(\rho^{\vee}_{P^{\prime}|P} \circ \phi_{\lambda} ,\psi_F) =\delta( \widetilde{\rho}^{\vee}_{P^{\prime}|P} \circ \phi^{\widetilde{M}^0}_{\lambda},\psi_F) \\
&=& \epsilon(0, \widetilde{\rho}^{\vee}_{P^{\prime}|P} \circ \phi^{\widetilde{M}^0}_{\lambda},\psi_F) \cdot \epsilon(1/2, \widetilde{\rho}^{\vee}_{P^{\prime}|P} \circ \phi^{\widetilde{M}^0}_{\lambda},\psi_F)^{-1}  
\end{eqnarray*}
and the representation theoretic $\delta$-factor 
\begin{eqnarray*}
& & \delta(\pi_{\phi^{\widetilde{M}^0}_{\lambda}}, \widetilde{\rho}^{\vee}_{P^{\prime}|P} ,\psi_F) \\ &:=& \epsilon(0, \pi_{\phi^{\widetilde{M}^0}_{\lambda}}, \widetilde{\rho}^{\vee}_{P^{\prime}|P},\psi_F) \cdot \epsilon(1/2, \pi_{\phi^{\widetilde{M}^0}_{\lambda}}, \widetilde{\rho}^{\vee}_{P^{\prime}|P},\psi_F)^{-1}
\end{eqnarray*}
are equal.

\bigskip

For the case of a non-generic parameter $\psi \in \Psi(M)$ we similarly consider $\psi^{\widetilde{M}^0} :=\xi_* \psi \in \widetilde{\Psi}(\widetilde{M}^0)$, and the above discussion concerning the $L$, $\epsilon$ and $\delta$-factors applies to the associated generic parameter $\phi_{\psi}^{\widetilde{M}^0} := \phi_{\psi^{\widetilde{M}^0}} \in \widetilde{\Phi}(\widetilde{M}^0)$. In particular, if we put
\[
\pi_{\psi^{\widetilde{M}^0}_{\lambda}} := \pi_{\phi_{\psi,\lambda}^{\widetilde{M}^0}} 
\]
then we have
\begin{eqnarray*}
& & r_{P^{\prime}|P}(\psi_{\lambda})= r_{P^{\prime}|P}(\psi_{\lambda},\psi_F) \\
&=& \delta (\pi_{\psi_{\lambda}^{\widetilde{M}^0}},\widetilde{\rho}^{\vee}_{P^{\prime}|P},\psi_F)^{-1} L(0, \pi_{\psi_{\lambda}^{\widetilde{M}^0}},\widetilde{\rho}^{\vee}_{P^{\prime}|P}) L(1, \pi_{\psi_{\lambda}^{\widetilde{M}^0}},\widetilde{\rho}^{\vee}_{P^{\prime}|P})^{-1}.
\end{eqnarray*}

\bigskip

\begin{proposition}
Assume that the local and global theorems of section 2.5 are valid if $N$ is replaced by any integer $N_- < N$. 

\bigskip
\noindent (a) The operators (3.3.6) satisfy the following multiplicative property: for $P,P^{\prime},P^{\prime \prime} \in \mathcal{P}(M)$, we have
\begin{eqnarray}
R_{P^{\prime \prime}|P }(\pi_{\lambda},\psi_{\lambda}) = R_{P^{\prime \prime}|P^{\prime} }(\pi_{\lambda},\psi_{\lambda}) R_{P^{\prime}|P }(\pi_{\lambda},\psi_{\lambda}).
\end{eqnarray}

\noindent (b) We have the adjoint relation:
\begin{eqnarray}
R_{P^{\prime}|P}(\pi_{\lambda},\psi_{\lambda})^* = R_{P|P^{\prime}}(\pi_{-\overline{\lambda}},\psi_{-\overline{\lambda}})=R_{P^{\prime}|P}(\pi_{-\overline{\lambda}},\psi_{-\overline{\lambda}})^{-1}
\end{eqnarray}
In particular the operator $R_{P^{\prime}|P}(\pi_{\lambda},\psi_{\lambda})$ is uniatry and hence analytic if $\lambda$ is purely imaginary. Hence we can define the operator
\[
R_{P^{\prime}|P}(\pi,\psi) := R_{P^{\prime}|P}(\pi_0,\psi_0).
\]
\end{proposition}

\begin{proof}
Since the argument is basically the same as the proof of proposition 2.3.1 of \cite{A1}, we will only sketch the argument (the proof uses local-global arguments; that is why we need the validity of the global theorems also for $N_- < N$, even though the statement of the proposition by itself is a local statement).  

We first consider the case where $\psi=\phi \in \Phi_{\bdd}(M)$ is a generic parameter. Then $\Pi_{\psi}=\Pi_{\phi}$ is an $L$-packet consisting of irreducible tempered representations of $M(F)$ (and consequently $\phi$ is determined by $\pi$ by the local classification theorem for $M$). 

When $F$ is archimedean, then the assertions of the proposition is known from previous results of Arthur \cite{A7}, thus we can ssume that $F$ is non-archimedean. Furthermore the reduction procedure of \cite{A7} reduces the proof of the assertion to the case where $\pi \in \Pi_2(M)$, the set of representations of $M(F)$ that is square integrable modulo centre. Thus with respect to the decomposition:
\[
M= G_{E/F}(N_1^{\prime}) \times \cdots \times G_{E/F}(N_r^{\prime}) \times G_-
\] 
we have 
\[
\pi = \pi_1^{\prime} \boxtimes \cdots \boxtimes \pi_r^{\prime} \boxtimes \pi_-
\] 
with $\pi_- \in \Pi_2(G_-(F))$, and $\pi_i^{\prime} \in \Pi_2(G_{E/F}(N_i^{\prime})(F))=\Pi_2(\GL_{N_i^{\prime}}(E))$. Further reduction ({\it c.f.} proof of proposition 2.3.1 of \cite{A1}) allows the reduction to the case that the representations $\pi_i^{\prime}$ are supercuspidal representations of $G_{E/F}(N_i^{\prime})(F)=\GL_{N_i^{\prime}}(E)$.

If $E=F \times F$ is split, then $G_- = U_{E/F}(N_-) \cong \GL_{N_-/F}$, and the reduction procedure as above reduces to the case that $\pi$ is supercuspidal, and hence has a Whittaker model, and the results follow from Shahidi \cite{S2}. Thus we may assume that $E$ is a field extension of $F$. 

With these reductions, we now apply global methods. As in {\it loc. cit.}, we embed the representations $\pi_-$ and the $\pi_i^{\prime}$ 's as local components of discrete automorphic representations. In general, for any connected reductive group $\dot{G}$ over a number field $\dot{F}$, we denote by $\Pi_2(\dot{G})$ the set of irreducible admissible representations of $\dot{G}(\mathbf{A}_{\dot{F}})$ that belong to the discrete automorphic spectrum of $\dot{G}(\mathbf{A}_{\dot{F}})$. We first need

\begin{lemma}
Given local objects
\[
F,E, \pi_- \in \Pi_2(G_-) 
\]
there exists global objects 
\[
\dot{F}, \dot{E}, \dot{\pi}_- \in \Pi_2(\dot{G}_-), \mbox{ with } \dot{G}_- = U_{\dot{E}/\dot{F}}(N_-)
\]
and $L$-embedding
\[
\dot{\xi} = \dot{\xi}_{\dot{\chi}}: \leftexp{L}{\dot{G}_-} = \leftexp{L}{U_{\dot{E}/\dot{F}}(N_-)} \hookrightarrow \leftexp{L}{G_{\dot{E}/\dot{F}}(N_-)}
\]
together with a prime $u$ of $\dot{F}$, such that
\begin{enumerate}

\bigskip

\item $(F,E,\xi_{\chi},\pi_-)=(\dot{F}_u,\dot{E}_u,\dot{\xi}_{\dot{\chi},u} ,\dot{\pi}_{-,u})$

\bigskip

\item $\dot{\pi}_-$ is classified by a global generic parameter $\dot{\phi}_- \in \Phi(U_{\dot{E}/\dot{F}}(N_-),\dot{\xi}_{\dot{\chi}})$.

\bigskip

\item For any non-archimedean prime $v$ of $\dot{F}$ other than $u$, the representation $\dot{\pi_-}_v$ of $\dot{G}_-(\dot{F}_v)=U_{\dot{E}_v/\dot{F}_v}(N_-)(\dot{F}_v)$ has a non-zero vector fixed by a special maximal compact subgroup of $\dot{G}_-(\dot{F}_v)$.

\end{enumerate}
\end{lemma}

\bigskip

Lemma 3.3.2 is proved by an argument using the simple version of the invariant trace formula. We will establish a stronger form of the assertion in section 7.2, {\it c.f.} remark 7.2.5 (the argument there is of course independent of the propositions of this subsection), and for the moment take it for granted. 

\bigskip

For the general linear factors we similarly have:

\begin{lemma}
In the situation of lemma 3.3.2, we can in addition choose the number fields $\dot{E}/\dot{F}$ and the prime $u$ of $\dot{F}$ such that for $i=1,\cdots,r$, we have a cuspidal automorphic representation $\dot{\pi}_i^{\prime}$ of $G_{\dot{E}/\dot{F}}(N_i^{\prime})(\mathbf{A}_{\dot{F}})=\GL_{N_i^{\prime}}(\mathbf{A}_{\dot{E}})$, such that
\begin{enumerate}

\bigskip

\item $\dot{\pi}^{\prime}_{i,u} = \pi^{\prime}_i$

\bigskip

\item For each non-archimedean prime $v$ of $\dot{F}$ other than $u$, the representation $\dot{\pi}^{\prime}_i$ of $G_{\dot{E}_v/\dot{F}_v}(N_i^{\prime})(\dot{F}_v)=\GL_{N_i^{\prime}}(\dot{E}_v)$ has a non-zero vector fixed by a special maximal compact subgroup.
\end{enumerate} 
\end{lemma}

\noindent Lemma 3.3.3 follows from the results of Shahidi (\cite{S2}, proposition 5.1), which applies in the present situation, due to the fact that the representation $\pi^{\prime}_i$ of $G_{E/F}(N_i^{\prime})(F)=\GL_{N_i^{\prime}}(E)$ is supercuspidal, hence has Whittaker model and thus the result of \cite{S2} applies.

\bigskip

Back to the proof of proposition 3.3.1, we put
\begin{eqnarray*}
 \dot{M}= G_{\dot{E}/\dot{F}}(N_1^{\prime}) \times \cdots \times G_{\dot{E}/\dot{F}}(N_r^{\prime}) \times \dot{G}_- 
\end{eqnarray*}
and 
\begin{eqnarray*}
\dot{\pi}= \dot{\pi}^{\prime}_1 \boxtimes \cdots \boxtimes \dot{\pi}^{\prime}_r \boxtimes \dot{\pi}_-
\end{eqnarray*}
then $\dot{\pi}$ is an irreducible representation of $\dot{M}(\mathbf{A}_{\dot{F}})$ that lies in the discrete automorphic spectrum of $\dot{M}(\mathbf{A}_{\dot{F}})$. It is classified by a global generic parameter $\dot{\phi} =(\dot{\phi}^{\dot{\widetilde{M}}^0 },  \dot{\widetilde{\phi}}  )  \in \Phi(\dot{M},\dot{\xi})$ (here $\dot{\phi}^{\dot{\widetilde{M}}^0 } \in \Phi( \dot{\widetilde{M}}^0 )$), and satisifes the following: $\dot{\pi}_u = \pi$, $\dot{\phi}_u = \phi$, and for any non-archimedean prime $v$ of $\dot{F}$ other than $u$, the representation $\dot{\pi}_v$ has a non-zero vector fixed by a special maximal compact subgroup of $\dot{M}(\dot{F}_v)$. 

We can identify $\dot{M}$ as a Levi subgroup of $\dot{G} := U_{\dot{E}/\dot{F}}(N)$. Note that $\dot{G}_u = U_{E/F}(N)=G$. Similarly put
\begin{eqnarray*}
& & \dot{\widetilde{M}}^0 = \\
& & G_{\dot{E}/\dot{F}}(N_1^{\prime}) \times \cdots \times G_{\dot{E}/\dot{F}}(N_r^{\prime}) \times G_{\dot{E}/\dot{F}}(N_-) \times G_{\dot{E}/\dot{F}}(N_r^{\prime}) \times \cdots \times G_{\dot{E}/\dot{F}}(N_1^{\prime})
\end{eqnarray*}
then $\dot{\widetilde{M}}^0$ is a Levi subgroup of $G_{\dot{E}/\dot{F}}(N)$. Now the same character $\dot{\chi} \in \mathcal{Z}_{\dot{E}}$ defining $\dot{\xi} = \dot{\xi}_{\dot{\chi}} : \leftexp{L}{\dot{G}_-} \hookrightarrow \leftexp{L}{G_{\dot{E}/\dot{F}}(N_-)}$ defines the $L$-embedding $\leftexp{L}{\dot{G}} \hookrightarrow \leftexp{L}{G_{\dot{E}/\dot{F}}(N)}$ which we still denote as $\dot{\xi}$, and restricts to an $L$-embedding $\dot{\xi}: \leftexp{L}{\dot{M}} \hookrightarrow \leftexp{L}{\dot{\widetilde{M}}^0}$ (with the $L$-embedding $\leftexp{L}{\dot{\widetilde{M}}^0} \hookrightarrow \leftexp{L}{G_{\dot{E}/\dot{F}}(N)}$ being dual to the embedding $\dot{\widetilde{M}}^0 \hookrightarrow G_{\dot{E}/\dot{F}}(N)$).

For any $\dot{P} \in \mathcal{P}(\dot{M})$, we have the global induced representation $\mathcal{I}_{\dot{P}}(\dot{\pi})$, and for $\dot{P}^{\prime},\dot{P} \in \mathcal{P}(\dot{M})$ the global intertwining operator $M_{\dot{P}^{\prime}|\dot{P}}(\dot{\pi}_{\lambda}) : \mathcal{I}_{\dot{P}} (\dot{\pi}_{\lambda}) \rightarrow \mathcal{I}_{\dot{P}^{\prime}}(\dot{\pi}_{\lambda})$ (for $\lambda \in \mathfrak{a}^*_{\dot{M},\mathbf{C}}=\mathfrak{a}_{M,\mathbf{C}}^*$), defined by meromorphic continuation of the product:
\begin{eqnarray}
M_{\dot{P}^{\prime}|\dot{P}}(\dot{\pi}_{\lambda}) = \bigotimes_v J_{\dot{P}^{\prime}_v|\dot{P}_v} (\dot{\pi}_{v,\lambda}).
\end{eqnarray}
We then have Langlands' functional equation (section 1 of \cite{A8}): for $\dot{P},\dot{P}^{\prime},\dot{P}^{\prime \prime} \in \mathcal{P}(\dot{M})$, we have:
\begin{eqnarray}
M_{\dot{P}^{\prime \prime}|P}(\dot{\pi}_{\lambda}) = M_{\dot{P}^{\prime \prime}|\dot{P}^{\prime}}(\dot{\pi}_{\lambda}) M_{\dot{P}^{\prime}|P}(\dot{\pi}_{\lambda})
\end{eqnarray}

\noindent On the other hand, fix a non-trivial additive character $\psi_{\dot{F}} : \mathbf{A}_{\dot{F}}/\dot{F} \rightarrow \mathbf{C}^{\times}$ such that $\psi_{\dot{F},u}=\psi_F$. If we define the global normalizing factor:
\begin{eqnarray}
r_{\dot{P}^{\prime}|\dot{P}}(\dot{\phi}_{\lambda}) = \prod_v r_{\dot{P}^{\prime}_v|\dot{P}_v}(\dot{\phi}_{v,\lambda}) =  \prod_v r_{\dot{P}^{\prime}_v|\dot{P}_v}(\dot{\phi}_{v,\lambda}, \psi_{\dot{F},v}) 
\end{eqnarray}
as an Euler product, then we claim that we have the following:
\begin{enumerate}
\bigskip

\item $r_{\dot{P}^{\prime}|\dot{P}}(\dot{\phi}_{\lambda})$ has meromorphic continuation for all $\lambda \in \mathfrak{a}^*_{\dot{M},\mathbf{C}}$.

\bigskip

\item $r_{\dot{P}^{\prime \prime}|P}(\dot{\phi}_{\lambda}) = r_{\dot{P}^{\prime \prime}|\dot{P}^{\prime}}(\dot{\phi}_{\lambda}) r_{\dot{P}^{\prime}|P}(\dot{\phi}_{\lambda})$.
\end{enumerate}

\noindent Granting these two properties of the global normalizing factor $r_{\dot{P}^{\prime}|\dot{P}}(\dot{\phi}_{\lambda})$ for a moment, we then see that if we define
\begin{eqnarray*}
R_{\dot{P}^{\prime}|\dot{P}}(\dot{\pi}_{\lambda},\dot{\phi}_{\lambda}) = \bigotimes_v R_{\dot{P}^{\prime}_v|\dot{P}_v}(\dot{\pi}_{v,\lambda},\dot{\phi}_{\lambda})
\end{eqnarray*}
as a product, then $R_{\dot{P}^{\prime}|\dot{P}}(\dot{\pi}_{\lambda},\dot{\phi}_{\lambda})$ has an meromorphic continuation to all $\lambda \in \mathfrak{a}^*_{\dot{M},\mathbf{C}}$, and we have, by virtue of (3.3.9), (3.3.11) and the meromorphic continuation property, the following:
\begin{eqnarray}
M_{\dot{P}^{\prime}|\dot{P}}(\dot{\pi}_{\lambda}) = r_{\dot{P}^{\prime}|\dot{P}}(\dot{\phi}_{\lambda}) \cdot R_{\dot{P}^{\prime}|\dot{P}}(\dot{\pi}_{\lambda}).
\end{eqnarray}

\noindent Hence by (3.3.10) and item (2) above for $r_{\dot{P}^{\prime}|\dot{P}}(\dot{\phi}_{\lambda})$, we have
\begin{eqnarray}
R_{\dot{P}^{\prime \prime}|\dot{P}}(\dot{\pi}_{\lambda},\dot{\phi}_{\lambda}) =R_{\dot{P}^{\prime \prime}|\dot{P}^{\prime}}(\dot{\pi}_{\lambda},\dot{\phi}_{\lambda}) R_{\dot{P}^{\prime}|\dot{P}}(\dot{\pi}_{\lambda},\dot{\phi}_{\lambda}).
\end{eqnarray}

Now let $v$ be a prime of $\dot{F}$ other than $u$. If $v$ is an archimedean prime, then as we have quoted in the beginning of the proof, the identity
\begin{eqnarray}
& & \\
& &
R_{\dot{P}_v^{\prime \prime}|\dot{P}_v}(\dot{\pi}_{v,\lambda},\dot{\phi}_{v,\lambda}) =R_{\dot{P}_v^{\prime \prime}|\dot{P}_v^{\prime}}(\dot{\pi}_{v,\lambda},\dot{\phi}_{v,\lambda}) R_{\dot{P}_v^{\prime}|\dot{P}_v}(\dot{\pi}_{v,\lambda},\dot{\phi}_{v,\lambda}) \nonumber
\end{eqnarray}
follows from \cite{A7}. On the other hand, if $v$ is an non-archimedean prime othen than $u$, then from the construction of $\dot{\pi}$, the representation $\dot{\pi}_v$ has a non-zero vector fixed by a special maximal compact subgroup, hence by \cite{C,CS}, the representation $\dot{\pi}_v$ has a Whittaker model. Thus Shahidi's result \cite{S2} applies and hence the identity (3.3.14) also holds for non-archimedean $v$ other than $u$. On combining with (3.3.13) we then deduce:
\begin{eqnarray}
& &\\
& & R_{\dot{P}_u^{\prime \prime}|\dot{P}_u}(\dot{\pi}_{u,\lambda},\dot{\phi}_{u,\lambda}) =R_{\dot{P}_v^{\prime \prime}|\dot{P}_u^{\prime}}(\dot{\pi}_{u,\lambda},\dot{\phi}_{u,\lambda}) R_{\dot{P}_u^{\prime}|\dot{P}_u}(\dot{\pi}_{u,\lambda},\dot{\phi}_{u,\lambda}). \nonumber
\end{eqnarray}

\noindent Since $\dot{\pi}_u=\pi$ and $\dot{\phi}_u=\phi$, this proves (3.3.7) for the parameter $\phi$ and $\pi \in \Pi_{\phi}$, at least for $\pi$ satisfying the local conditions as in the beginning of the proof. But then by the reduction procedure of \cite{A7} as quoted in the beginning of the proof, this implies the validity of (3.3.7) for all $\phi \in \Phi_{\bdd}(M)$ and $\pi \in \Pi_{\phi}$.   

It remains to justify property (1) and (2) above for the global normalizing factor $r_{\dot{P}^{ \prime}|\dot{P}}(\dot{\phi}_{\lambda})$. Recall that we have $\dot{\phi} \in \Phi(\dot{M},\dot{\xi})$. By the discussion before the proof of the proposition, we have for each prime $v$ of $\dot{F}$:
\begin{eqnarray*}
& & L(s, \rho^{\vee}_{\dot{P}^{\prime}_v|\dot{P}_v} \circ \dot{\phi}_{v,\lambda} ) = L(s,  \pi_{\dot{\phi}^{\dot{\widetilde{M}}^0_v}_{v,\lambda}} , \widetilde{\rho}^{\vee}_{\dot{P}^{\prime}_v|\dot{P}_v} ) \\
& & \epsilon(s, \rho^{\vee}_{\dot{P}^{\prime}_v|\dot{P}_v} \circ \dot{\phi}_{v,\lambda} ,\psi_{\dot{F}_v}) = \zeta_v \cdot  \epsilon(s,  \pi_{\dot{\phi}^{\dot{\widetilde{M}}^0_v}_{v,\lambda}}, \widetilde{\rho}^{\vee}_{\dot{P}_v^{\prime}|\dot{P}_v} ,\psi_{\dot{F}_v})
\end{eqnarray*}  
where the $L$ and $\epsilon$-factors on the right hand side are given by a product of representation theoretic Rankin-Selberg factors on $G_{\dot{E}_v/\dot{F}_v}(N^{\prime})(\dot{F}_v) \times G_{\dot{E}_v/\dot{F}_v}(N^{\prime \prime})(\dot{F}_v)$ (with $N^{\prime},N^{\prime \prime} \leq N$), and $\Asai^{\pm}$ factors on $G_{\dot{E}_v/\dot{F}_v}(N^{\prime \prime \prime})(\dot{F}_v)$ (with $N^{\prime \prime \prime} \leq N$), and $\zeta_v$ is a constant equal to a root of unity (that is equal to one for almost all $v$). Note that by construction we have
\[
\dot{\phi}^{\dot{\widetilde{M}}^0} = \rightexp{\bigotimes_v}{\prime} \pi_{\dot{\phi}_v^{\dot{\widetilde{M}}_v^0 } }.
\]
Hence the result of \cite{JPSS} and \cite{S} gives the continuation of the global $L$ and $\epsilon$-functions, originally defined by Euler products:
\begin{eqnarray*}
& & L(s, \dot{\phi}^{\dot{\widetilde{M}}^0}_{\lambda}  , \widetilde{\rho}^{\vee}_{\dot{P}^{\prime}|\dot{P}} ) = \prod_v L(s,  \pi_{\dot{\phi}_{v,\lambda}^{\dot{\widetilde{M}}_v^0 } } , \widetilde{\rho}^{\vee}_{\dot{P}^{\prime}_v|\dot{P}_v} ) \\
& & \epsilon(s,  \dot{\phi}^{\dot{\widetilde{M}}^0}_{\lambda} , \widetilde{\rho}^{\vee}_{\dot{P}^{\prime}|\dot{P}} ,\psi_{\dot{F}_v}) = \prod_v \epsilon(s,  \pi_{\dot{\phi}_{v,\lambda}^{\dot{\widetilde{M}}_v^0 } }, \widetilde{\rho}^{\vee}_{\dot{P}_v^{\prime}|\dot{P}_v} ,\psi_{\dot{F}_v})
\end{eqnarray*}
as meromorphic functions for all $s$ and $\lambda$. In particular this gives the meromorphic continuation of the global normalizing factor:
\begin{eqnarray}
r_{\dot{P}^{\prime}|\dot{P}}(\dot{\phi}_{\lambda}) &=&\epsilon(1/2, \dot{\phi}^{\dot{\widetilde{M}}^0}_{\lambda}, \widetilde{\rho}^{\vee}_{\dot{P}^{\prime}|\dot{P}} ) \epsilon(0,\dot{\phi}^{\dot{\widetilde{M}}^0}_{\lambda} , \widetilde{\rho}^{\vee}_{\dot{P}^{\prime}|\dot{P}} )^{-1} \\
& & \times   L(0,\dot{\phi}^{\dot{\widetilde{M}}^0}_{\lambda} , \widetilde{\rho}^{\vee}_{\dot{P}^{\prime}|\dot{P}}  ) L(1,\dot{\phi}^{\dot{\widetilde{M}}^0}_{\lambda} , \widetilde{\rho}^{\vee}_{\dot{P}^{\prime}|\dot{P}}  )^{-1}.  \nonumber
\end{eqnarray}

\noindent Thus we have property (1). To establish property (2) for the global normalizing factor, it suffices, again by the reduction procedure in \cite{A7}, to treat the case where $\dot{M}$ is a maximal Levi subgroup of $\dot{G}$, in which case there are only two parabolic subgroup of $\dot{G}$ with $\dot{M}$ as Levi component, which we name as $\dot{P}$ and $\dot{\overline{P}}$ (thus $\dot{\overline{P}}$ is the opposite parabolic of $\dot{P}$), and it suffices to take $\dot{P}^{\prime}=\dot{\overline{P}}, \dot{P}^{\prime \prime} =\dot{P}$. Then the identity (2) that we want to establish boils down to:
\begin{eqnarray}
1=    r_{\dot{P}|\dot{\overline{P}} }(\dot{\phi}_{\lambda}) \cdot r_{\dot{\overline{P}}|\dot{P} }(\dot{\phi}_{\lambda}).
\end{eqnarray}

\noindent To establish (3.3.17) we use the global function equation (which follows from \cite{JPSS} and \cite{S}):
\begin{eqnarray}
& & \\
& & L(s, \dot{\phi}^{\dot{\widetilde{M}}^0}_{\lambda} , \widetilde{\rho}^{\vee}_{\dot{\overline{P}}  |\dot{P}}  ) = \epsilon(s, \dot{\phi}^{\dot{\widetilde{M}}^0}_{\lambda} , \widetilde{\rho}^{\vee}_{\dot{\overline{P}}  |\dot{P}}  ) L(1-s, \dot{\phi}^{\dot{\widetilde{M}}^0}_{\lambda} , \widetilde{\rho}_{\dot{\overline{P}}  |\dot{P}}  ). \nonumber
\end{eqnarray} 
Then (3.3.16) becomes
\begin{eqnarray}
& & \\
& & r_{\dot{P}^{\prime}|\dot{P}}(\dot{\phi}_{\lambda}) = \epsilon(1/2, \dot{\phi}^{\dot{\widetilde{M}}^0}_{\lambda} , \widetilde{\rho}^{\vee}_{\dot{P}^{\prime}|\dot{P}} ) L(1, \dot{\phi}^{\dot{\widetilde{M}}^0}_{\lambda} , \widetilde{\rho}_{\dot{\overline{P}}  |\dot{P}}  )  L(1, \dot{\phi}^{\dot{\widetilde{M}}^0}_{\lambda} , \widetilde{\rho}^{\vee}_{\dot{\overline{P}}  |\dot{P}}  )^{-1}. \nonumber
\end{eqnarray}
Since $\widetilde{\rho}^{\vee}_{\dot{\overline{P}}  |\dot{P}} = \widetilde{\rho}_{\dot{P}|\dot{\overline{P}}}$, it follows from (3.3.19) that
\begin{eqnarray}
& & \\
& &  r_{\dot{P}| \dot{\overline{P}}}(\dot{\phi}_{\lambda})  r_{\dot{\overline{P}}|\dot{P}}(\dot{\phi}_{\lambda}) =  \epsilon(1/2, \dot{\phi}^{\dot{\widetilde{M}}^0}_{\lambda} , \widetilde{\rho}^{\vee}_{\dot{P}|\dot{\overline{P}}} )  \epsilon(1/2, \dot{\phi}^{\dot{\widetilde{M}}^0}_{\lambda} , \widetilde{\rho}^{\vee}_{\dot{\overline{P}}|\dot{P}} ) \nonumber
\end{eqnarray}
but the functional equation (3.3.18) (applied two times) implies that the right hand side of (3.3.20) is equal to one, thus we obtain (3.3.17).

\bigskip

Finally, the deduction of (3.3.7) for the case of general parameters $\psi \in \Psi(M)$ and the adjoint relation (3.3.8), from the results already established above for the parameters $\psi = \phi \in \Phi_{\bdd}(M)$, can be done exactly as in the proof of proposition 2.3.1 of \cite{A1}.

\end{proof}

\bigskip

This finishes step one in the construction of normalized local intertwining operator. The second step is to convert the operators $R_{P^{\prime}|P}(\pi,\psi)$ into certain  intertwining operators of induced representations attached to $w \in W(M)$. 

Thus let $P \in \mathcal{P}(M)$, and $w \in W(M)$. As before the quasi-split group $G=U_{E/F}(N)$ has the standard $F$-splitting
\[
S=(T,B,\{x_{\alpha} \})
\]
we assume that $M$ and $P$ are standard in the sense that they contain $T$ and $B$ respectively. Then as in \cite{LS} we can choose a lifting $\widetilde{w} \in N(M)(F)$ of $w$ with respect to the splitting $S$ as follows. Let $w_T$ be the representative of $w$ in the rational Weyl group $W_F(G,T)$ that stabilizes the simple roots of $(B \cap M,T)$. Write $w_T=w_{\alpha_1} \cdots w_{\alpha_r}$ in its reduced decomposition relative to the roots $(B,T)$. Then
\[
\widetilde{w} = \widetilde{w}_{\alpha_1} \cdots \widetilde{w}_{\alpha_r}
\]
where 
\[
\widetilde{w}_{\alpha} = \exp(x_{\alpha}) \exp(-x_{-\alpha}) \exp(x_{\alpha}).
\]

\bigskip

Now as before let $\psi \in \Psi(M)$, and $\pi \in \Pi_{\psi}$; in particular $\pi$ is an irreducible unitary representation of $M(F)$. For $w \in W(M)$, define another representation:
\[
(w \pi)(m) = \pi(\widetilde{w}^{-1} m \widetilde{w}), \,\ m \in M(F)
\]  
of $M(F)$ on the underlying space $V_{\pi}$ of $\pi$. Put
\[
P^{\prime} = w^{-1} P := \widetilde{w}^{-1} P \widetilde{w}
\]
then we have the intertwining isomorphism
\begin{eqnarray}
l(\widetilde{w},\pi): \mathcal{H}_{P^{\prime}}(\pi) \longrightarrow \mathcal{H}_P(w \pi)
\end{eqnarray}
from $\mathcal{I}_{P^{\prime}}(\pi)$ to $\mathcal{I}_{P}(w \pi)$, by:
\begin{eqnarray}
(l(\widetilde{w},\pi) \phi^{\prime})(x) = \phi^{\prime}(\widetilde{w}^{-1} x), \,\ \phi^{\prime} \in \mathcal{H}_{P^{\prime}}(\pi), x \in G(F).
\end{eqnarray}

Since the mapping $w \rightarrow \widetilde{w}$ is not multiplicative in $w$, the map $w \mapsto l(\widetilde{w},\pi)$ does not satisfy the cocycle property. As in \cite{A1}, we remedy this by multiplying $l(\widetilde{\pi},\pi)$ with its own normalizing factors.

First with $\psi$ as above, we have the irreducible unitary representation $\pi_{\psi^{\widetilde{M}^0}}$ that corresponds to the parameter $\psi^{\widetilde{M}^0} :=\xi_{*} \psi \in \widetilde{\Psi}(\widetilde{M}^0(F))$. Define the $\epsilon$ factor associated to $w$:
\begin{eqnarray}
& & \\
& & \epsilon_P(w,\psi)=\epsilon(w,\psi,\psi_F)=\epsilon(1/2,\pi_{\psi^{\widetilde{M}^0}},\widetilde{\rho}^{\vee}_{w^{-1}P|P},\psi_F), \,\ w \in W(M). \nonumber
\end{eqnarray}
where as before $\epsilon(1/2,\pi_{\psi^{\widetilde{M}^0}},\widetilde{\rho}^{\vee}_{w^{-1}P|P},\psi_F)$ is the representation theoretic $\epsilon$ factor. 

\bigskip

Secondly also define the $\lambda$-factor:
\begin{eqnarray}
\lambda(w) = \lambda(w,\psi_F)  = \prod_{a} \lambda(F_a/F,\psi_F)
\end{eqnarray}
where the product ranges over the reduced roots of $(B,A_B)$ such that $w_T a < 0$ (here $A_B$ is the split component of $B$ or $T$). The extension $F_a$ is either the quadratic extension $E$ of $F$ or equal to $F$ itself, and is determined by the condition that $G_{a,sc} = \Res_{F_a/F} \SL_2$, with $G_a$ being the Levi subgroup of $G$ of semisimple rank one determined by the root $a$; equivalently each $a$ determines an orbit $\{\alpha\}$ of roots of $(B,T)$ under $\Gamma_F$ and $F_a$ is the extension that corresponds to the stabilizer of $\alpha$. Finally $\lambda(F_a/F,\psi_F)$ is the Langlands constant associated with this extension, given by
\begin{eqnarray*}
\lambda(F_a/F,\psi_F) = \frac{\epsilon(s, \Ind^F_{F_a} 1_{F_a},\psi_F)}{\epsilon(s,1_F,\psi_F)}
\end{eqnarray*}  
(which is a constant independent of $s$, and the $\epsilon$-factors being the Artin ones in the Langlands normalization as in \cite{T}). We obviously have $\lambda(F/F,\psi_F)=1$, and by equation (2.7) of \cite{KeS}, we have
\begin{eqnarray}
\lambda(E/F,\psi_F)^2 = \omega_{E/F}(-1)
\end{eqnarray}
(as before $\omega_{E/F}$ being the quadratic character of $F^{\times}$ that corresponds to the extension $E/F$ under local class field theory). To simplify the notation we omit the notational dependence on the additive character $\psi_F$ for the factors $\epsilon(w,\psi,\psi_F)$ and $\lambda(w,\psi_F)$ below.

\bigskip

\begin{lemma}
The product
\begin{eqnarray}
l(w,\pi,\psi) := \lambda(w)^{-1} \epsilon_P(w,\psi) l(\widetilde{w},\pi), \,\ w \in W(M)
\end{eqnarray}
satisfies the cocycle relation:
\begin{eqnarray}
& & l(w^{\prime} w,\pi,\psi) \\
& =&  l(w^{\prime},w \pi,w \psi) l(w,\pi,\psi), \,\ w^{\prime}, w \in W(M). \nonumber
\end{eqnarray}
\end{lemma}
\begin{proof}

Again since the proof is basically the same as lemma 2.3.4 of \cite{A1} we will only sketch the argument. Thus let $w^{\prime},w \in W(M)$. From the fact that both $\widetilde{w}^{\prime} \widetilde{w}$ and $\widetilde{w^{\prime} w}$ are representatives of $w^{\prime} w$ in $G(F)$ that preserve the same splitting of $M$, we have:
\[
\widetilde{w}^{\prime} \widetilde{w} = z(w^{\prime},w) \widetilde{w^{\prime} w} 
\]  
for $z(w^{\prime},w)$ in the centre of $M(F)$. It follows that
\[
l(\widetilde{w^{\prime} w},\pi) = \eta_{\pi}(z(w^{\prime},w)) l( \widetilde{w}^{\prime} ,w \pi) l(\widetilde{w},\pi)
\]
where $\eta_{\pi}$ is the central character of $\pi$. Moreover, from lemma 2.1.A of \cite{LS}, we have the following expression for $z(w^{\prime},w)$: 
\begin{eqnarray*}
z(w^{\prime},w) = (-1)^{\lambda(w^{\prime},w)}
\end{eqnarray*}
where $\lambda(w^{\prime},w) \in X_*(T)$ is the cocharacter of $T$ determined by: 
\[
\lambda(w^{\prime},w) = \sum_{\alpha \in \Sigma_B(w^{\prime}_T,w_T)} \alpha^{\vee}.
\]
In the summation $\Sigma_B(w^{\prime}_T,w_T)$ is the subset of roots $\Sigma_B$ of $(B,T)$ such that $w_T \alpha \notin \Sigma_B$ but $w^{\prime}_T w_T \alpha \in \Sigma_B$, and $\alpha^{\vee}$ being the coroot of $\alpha$. The argument in lemma 2.3.4 of \cite{A1} shows that the rational Weyl group $W_0^M$ of $(M,T)$ permutes the set $\Sigma_B(w^{\prime}_T,w_T)$, and hence $z(w^{\prime},w)$ is fixed by $W^M_0$. Since it is also fixed by $\Gamma_F$ (because $w^{\prime}_T,w_T$ are in the rational Weyl group $W_F(G,T)$), it follows that $\lambda(w^{\prime},w) \in X_*(A_M)$. In particular the element $z(w^{\prime},w) \in A_M(F)$. 

\bigskip

In general for any $\lambda \in X_*(A_M)$ denote by $\lambda^{\vee}$ the $\Gamma_F$-invariant character of $\widehat{M}$ that corresponds to $\lambda$, under the isomorphism:
\begin{eqnarray}
X_*(A_M) \cong X^*(\widehat{M})^{\Gamma_F}.
\end{eqnarray}
We also extend such an element $\lambda^{\vee}$ to a character on $\leftexp{L}{M}$ that is trivial on $W_F$. Then for any $u \in W_F$ with $x \in F^{\times}$ being the image of $u$ in $W_F^{ab} \cong F^{\times}$, then as a consequence of part (c) of theorem 2.5.1 applied to $G_-$ (the part on central character), we have
\begin{eqnarray}
\eta_{\pi}(x^{\lambda}) = \lambda^{\vee} \circ \phi_{\psi}(u).
\end{eqnarray}

\bigskip

\noindent Hence we have:
\begin{eqnarray}
\eta_{\pi}(z(w^{\prime},w) ) = \lambda^{\vee}(\phi_{\psi}(u))
\end{eqnarray}
for any $u \in W_F$ whose image in $W_F^{ab} \cong F^{\times}$ is equal to $-1$, and $\lambda^{\vee} = \lambda^{\vee}(w^{\prime},w):= (\lambda(w^{\prime},w))^{\vee}$.

\bigskip

Denote by $\widehat{\Sigma}^r_P$ be the set of reduced roots of $(\widehat{P},A_{\widehat{M}})$. We then have a decomposition:
\begin{eqnarray}
\lambda^{\vee}(w^{\prime},w) = \sum_{\beta \in \widehat{\Sigma}^r_P(w^{\prime},w)} \lambda^{\vee}_{\beta}(w^{\prime},w)
\end{eqnarray} 

\noindent with $\widehat{\Sigma}^r_P(w^{\prime},w)$ being the set
\[
\widehat{\Sigma}^r_P(w^{\prime},w) = \{\beta \in \widehat{\Sigma}_P^r, \,\ w \beta <0, w^{\prime} w \beta >0   \}
\]
and 
\[
\lambda^{\vee}_{\beta}=\lambda^{\vee}_{\beta}(w^{\prime},w) = \sum_{\alpha \in \Sigma_{\beta}} \alpha^{\vee}
\]
with $\Sigma_{\beta}$ the set of roots $\alpha$ of $(B,T)$ whose coroot $\alpha^{\vee}$ restricts to a positive multiple of $\beta$ on $A_{\widehat{M}}$. The same argument applied above to $\lambda^{\vee}(w^{\prime},w)$ can also be applied to $\lambda^{\vee}_{\beta}(w^{\prime},w)$ and shows that $\lambda_{\beta}^{\vee} (w^{\prime},w) \in X^*(\widehat{M})^{\Gamma_F}$. Hence we have by (3.3.30) and (3.3.31):
\begin{eqnarray}
\eta_{\pi}(z(w^{\prime},w))=\prod_{\beta \in \widehat{\Sigma}^r_P(w^{\prime},w)} (\lambda_{\beta}^{\vee} \circ \phi_{\psi})(u).
\end{eqnarray}

For $\beta \in \widehat{\Sigma}^r_P$, denote by $\rho_{\beta}$ the adjoint representation of $\leftexp{L}{M}$ on
\[
\widehat{\mathfrak{n}}_{\beta} := \bigoplus_{\alpha \in \Sigma_{\beta}} \widehat{\mathfrak{n}}_{\alpha} \subset \widehat{\mathfrak{n}}_P.
\]
We have $\rho_{\beta}^{\vee} = \rho_{-\beta}$, and similar to before, we have a factorization $\rho_{\beta} = \widetilde{\rho}_{\beta} \circ \xi$, where $\widetilde{\rho}_{\beta}$ is a representation of $\leftexp{L}{\widetilde{M}^0}$ whose simple constituents are Rankin-Selberg or $\Asai^{\pm}$ representations of the general linear factors of $\leftexp{L}{\widetilde{M}^0}$. We then have $\rho_{\beta} \circ \phi_{\psi} \cong \widetilde{\rho}_{\beta} \circ \phi_{\psi^{\widetilde{M}^0}}$. Then from the definition (3.3.23) we have
\begin{eqnarray}
\epsilon_P(w,\psi) = \prod_{\beta \in \widehat{\Sigma}^r_P, w \beta <0} \epsilon(1/2, \pi_{\psi^{\widetilde{M}^0}}, \widetilde{\rho}_{\beta}  ,\psi_F). 
\end{eqnarray}

Put 
\[
e(w^{\prime},w) = \epsilon_P(w^{\prime},w\psi) \epsilon_P(w,\psi) \epsilon_P(w^{\prime} w,\psi)^{-1}
\]
then by (3.3.33), we have the expression
\begin{eqnarray}
& & \\
 e(w^{\prime},w) &=& \prod_{\beta \in \widehat{\Sigma}^r_P(w^{\prime},w)} \epsilon(1/2,\pi_{\psi^{\widetilde{M}^0}}, \widetilde{\rho}_{\beta},\psi_F)  \epsilon(1/2,\pi_{\psi^{\widetilde{M}^0}}, \widetilde{\rho}_{-\beta},\psi_F) \nonumber \\
 &=& \prod_{\beta \in \widehat{\Sigma}^r_P(w^{\prime},w)} \epsilon(1/2,\pi_{\psi^{\widetilde{M}^0}}, \widetilde{\rho}_{\beta},\psi_F)  \epsilon(1/2,\pi_{\psi^{\widetilde{M}^0}}, \widetilde{\rho}^{\vee}_{\beta},\psi_F). \nonumber
\end{eqnarray}

For each $\beta$, the product of the two representation theoretic $\epsilon$-factors is equal to the product of the corresponding Artin $\epsilon$-factors ({\it c.f. \cite{H2}}):
\begin{eqnarray*}
& &  \epsilon(1/2,\pi_{\psi^{\widetilde{M}^0}}, \widetilde{\rho}_{\beta},\psi_F)  \epsilon(1/2,\pi_{\psi^{\widetilde{M}^0}}, \widetilde{\rho}^{\vee}_{\beta},\psi_F) \\ 
&=& \epsilon(1/2,   \widetilde{\rho}_{\beta} \circ \phi_{\psi^{\widetilde{M}^0}} ,\psi_F) \epsilon(1/2,\widetilde{\rho}^{\vee}_{\beta} \circ \phi_{\psi^{\widetilde{M}^0}},\psi_F).
\end{eqnarray*}

By equation (3.6.8) of \cite{T}, we have for $u \in W_F$ that maps to $-1$ under $W_F^{ab} \cong F^{\times}$ the following:
\[
\epsilon(1/2,   \widetilde{\rho}_{\beta} \circ \phi_{\psi^{\widetilde{M}^0}} ,\psi_F) \epsilon(1/2,\widetilde{\rho}^{\vee}_{\beta} \circ \phi_{\psi^{\widetilde{M}^0}},\psi_F) = \det (\widetilde{\rho}_{\beta} \circ \phi_{\psi^{\widetilde{M}^0}})(u).
\]

\noindent Hence we obtain:
\begin{eqnarray}
& & \\
& & e(w^{\prime},w) = \prod_{\beta \in \widehat{\Sigma}^r_P(w^{\prime},w)} \det(\widetilde{\rho}_{\beta} \circ \phi_{\psi^{\widetilde{M}^0}})(u) = \prod_{\beta \in \widehat{\Sigma}^r_P(w^{\prime},w)} \det(\rho_{\beta} \circ \phi_{\psi})(u). \nonumber
\end{eqnarray}

\noindent Similarly for the $\lambda$-factor, for $w^{\prime},w \in W(M)$ we have, from the definition (3.3.24) and equation (3.3.25):
\begin{eqnarray}
\lambda(w^{\prime} w) \lambda(w)^{-1} \lambda(w^{\prime})^{-1} = \prod_{\beta \in \widehat{\Sigma}^r_P(w^{\prime},w)} \omega_{E/F}(-1)^{n_{\beta}}
\end{eqnarray}
where $n_{\beta}$ is the number of $\Gamma_F$-orbits of roots $\{ \alpha \}$ associated to reduced roots $a$ of $(B,A_B)$ with $F_a=E$ and whose $\Gamma_F$-orbit $\{\alpha\}$ is in $\Sigma_{\beta}$. 

\noindent Hence by (3.3.32), (3.3.35) and (3.3.36), we see that to prove (3.3.27) we have to show (for $\beta \in \widehat{\Sigma}^r_P(w^{\prime},w)$):
\begin{eqnarray}
 \det \rho_{\beta}(\phi_{\psi}(u)) = \lambda^{\vee}_{\beta}(\phi_{\psi}(u)) \omega_{E/F}(-1)^{n_{\beta}}. 
\end{eqnarray}

\noindent Write 
\[
\phi_{\psi}(u) = m_{\psi}(u) \rtimes u, \,\ m_{\psi}(u) \in \widehat{M}
\]
then from the definition of $\rho_{\beta}$ we have $\det \rho_{\beta} = \lambda^{\vee}_{\beta}$ on $\widehat{M}$. Thus in particular
\[
\det \rho_{\beta}(m_{\psi}(u)) = \lambda_{\beta}^{\vee}(m_{\psi}(u)).
\]

\noindent Thus
\begin{eqnarray*}
& & \det \rho_{\beta}(\phi_{\psi}(u)) \\
& =& \det \rho_{\beta}(m_{\psi}(u)) \det \rho_{\beta}(1 \rtimes u) \\
&=& \lambda_{\beta}^{\vee}(m_{\psi}(u)) \det \rho_{\beta}(1 \rtimes u) \\
&=& \lambda_{\beta}^{\vee}(\phi_{\psi}(u)) \det \rho_{\beta}(1 \rtimes u).
\end{eqnarray*}

\noindent Hence to prove (3.3.37) it remains to establish
\begin{eqnarray}
\det \rho_{\beta}(1 \rtimes u) = \omega_{E/F}(-1)^{n_{\beta}}.
\end{eqnarray}
Since $u$ maps to $-1$ under $W_F^{ab} \cong F^{\times}$, we have $\omega_{E/F}(-1) =-1$ if and only if the image of $u$ in $\Gal(E/F)$ is nontrivial. Hence it suffices to show
\begin{eqnarray}
 \det \rho_{\beta}(1 \rtimes \sigma)   =(-1)^{n_{\beta}}
\end{eqnarray} 
for $\sigma \in W_F$ that maps to the nontrivial element of $\Gal(E/F)$. This can be verified by a direct computation.

As before, we have the standard decomposition
\[
M = G_{E/F}(N_1^{\prime}) \times \cdots \times G_{E/F}(N_r^{\prime}) \times U_{E/F}(N_-)
\]
and with $M$ being embedded in $G=U_{E/F}(N)$ in the standard way with respect to the partition:
\[
N = N_1^{\prime} + \cdots + N_r^{\prime} + N_- + N_r^{\prime} + \cdots + N_1^{\prime}.
\]

\noindent Then with respect to the standard splitting, the set of reduced roots of $(\widehat{P},A_{\widehat{M}})$ is given as follows. First under the decomposition of $M$ as above, we have
\[
A_{\widehat{M}} = \underbrace{\mathbf{C}^{\times} \times \cdots \times \mathbf{C}^{\times}}_{r}
\]
and where for $1 \leq i \leq r$, the factor $\mathbf{C}^{\times}$ of $A_{\widehat{M}}$ at the $i$-th spot corresponds to the factor $G_{E/F}(N_i^{\prime})$ of $M$. Denote by $s_i$, for $1 \leq i \leq r$ the character of $A_{\widehat{M}}$ given by projection to the $i$-th coordinate. Then we have:

\begin{eqnarray}
& & \\
& &  \widehat{\Sigma}^r_P =  \widehat{\Sigma}^r_P(1) \cup \widehat{\Sigma}^r_P(2) \cup \widehat{\Sigma}^r_P(3)  \nonumber \\
& &  \widehat{\Sigma}^r_P (1) =  \{\beta= \beta^{(1)}_{ij} = s_i s_j^{-1}, \,\  1 \leq i < j \leq r  \}, \nonumber \\
& &  \widehat{\Sigma}^r_P (2) =  \{\beta= \beta^{(2)}_{ij} = s_i s_j, \,\  1 \leq i < j \leq r  \}, \nonumber \\
& & \widehat{\Sigma}^r_P(3)= \{ \beta=\beta^{(3)}_i=s_i, \,\ 1 \leq i \leq r \}.  \nonumber
\end{eqnarray}

Given $\beta \in \widehat{\Sigma}^r_P$ as in (3.3.40), the set of roots $\alpha \in \Sigma_B$ that belongs to $\Sigma_{\beta}$ is given as follows. First denote 
\[
T_N:= \underbrace{G_{E/F}(1) \times \cdots \times G_{E/F}(1)}_{N}
\] 
the standard diagonal maximal torus of $G_{E/F}(N)$. For $1\leq u \leq N$, denote by $t_u$ the character of $T_{N/\overline{F}}$:
\[
t_u : T_{N/\overline{F}} \rightarrow \mathbf{G}_{m/\overline{F}}
\]
given by the composition of the projection to the $u$-th factor, followed by the projection
\[
(E \otimes_F \overline{F})^{\times} = \overline{F}^{\times} \times \overline{F}^{\times}
\]
to the left $\overline{F}^{\times}$-factor corresponding to the specified embedding of $E$ into $\overline{F}$ ({\it c.f.} the section on Notation). Then under the standard embedding of $T$ to $T_N$, we have:

\begin{eqnarray*}
& & \\
 & & \mbox{ For } \beta = \beta^{(1)}_{ij} \in \widehat{\Sigma}^r_P(1), \,\   \Sigma_{\beta} \nonumber \\
 & =& \Big\{\alpha^{(ij)}_{uv}  =  t_u t_v^{-1}, \,\     (\sum_{k=1}^{i-1} N_k^{\prime})+1  \leq  u \leq \sum_{k=1}^i N_k^{\prime}, \,\    (\sum_{k=1}^{j-1} N_k^{\prime})+1  \leq  v \leq \sum_{k=1}^j N_k^{\prime} \Big\}    \\
 & & \bigcup \Big\{\alpha^{(ij)}_{uv} = t_u t_v^{-1},  (\sum_{k=1}^{r} N_k^{\prime}) + N_- + (\sum_{k=j+1}^{r} N_k^{\prime} ) + 1  \leq  u \leq  (\sum_{k=1}^{r} N_k^{\prime}) + N_- + (\sum_{k=j}^{r} N_k^{\prime} ), \nonumber \\  & & \,\ \,\ \,\ \,\   (\sum_{k=1}^{r} N_k^{\prime}) + N_- + (\sum_{k=i+1}^{r} N_k^{\prime} ) + 1  \leq  v \leq  (\sum_{k=1}^{r} N_k^{\prime}) + N_- + (\sum_{k=i}^{r} N_k^{\prime} )     \Big\}. \nonumber
\end{eqnarray*}

\begin{eqnarray*}
& & \\
 & & \mbox{ For } \beta = \beta^{(2)}_{ij} \in \widehat{\Sigma}^r_P(2),  \nonumber \\
  & & \Sigma_{\beta}  = \Big\{\alpha^{(ij)}_{uv}  =  t_u t_v^{-1}, \,\     (\sum_{k=1}^{i-1} N_k^{\prime})+1  \leq  u \leq \sum_{k=1}^i N_k^{\prime}, \\  
& & \,\ \,\ \,\  (\sum_{k=1}^{r} N_k^{\prime})+N_- + (\sum_{k=j+1}^{r} N_k^{\prime}) +1  \leq  v \leq (\sum_{k=1}^r N_k^{\prime}) + N_- + (\sum_{k=j}^{r} N_k^{\prime}) \Big\}    \nonumber \\
 & & \bigcup \Big\{\alpha^{(ij)}_{uv} = t_u t_v^{-1}, \,\   (\sum_{k=1}^{j-1} N_k^{\prime})+1  \leq  u \leq \sum_{k=1}^j N_k^{\prime}, \\
& & \,\ \,\  (\sum_{k=1}^{r} N_k^{\prime})+N_- + (\sum_{k=i+1}^{r} N_k^{\prime}) +1  \leq  v \leq (\sum_{k=1}^r N_k^{\prime}) + N_- + (\sum_{k=i}^{r} N_k^{\prime})    \nonumber    \Big\}. \nonumber
\end{eqnarray*}

\begin{eqnarray*}
& & \\
 & & \mbox{ For } \beta = \beta^{(3)}_{i} \in \widehat{\Sigma}^r_P(3),   \,\ \Sigma_{\beta} \\
&=&  \Big\{ \alpha^{(i)}_{uv} = t_u t_v^{-1}, \,\   (\sum_{k=1}^{i-1} N_k^{\prime})+1   \leq u \leq \sum_{k=1}^i N_k^{\prime}, \,\ (\sum_{k=1}^r N_k^{\prime})+1  \leq v \leq (\sum_{k=1}^{r} N_k^{\prime}) + N_- \Big\} \\
& & \bigcup \Big\{ \alpha^{(i)}_{uv} = t_u t_v^{-1}, \,\  (\sum_{k=1}^r N_k^{\prime})+1  \leq u \leq (\sum_{k=1}^{r} N_k^{\prime}) + N_-  , \,\    (\sum_{k=1}^{i-1} N_k^{\prime})+1   \leq v \leq \sum_{k=1}^i N_k^{\prime}   \Big\} \\
& & \bigcup \Big\{  \alpha^{(i)}_{uv} = t_u t_v^{-1}, \,\    (\sum_{k=1}^{i-1} N_k^{\prime})+1   \leq u \leq    \sum_{k=1}^{i} N_k^{\prime}, \\                                   
& & \,\ \,\ \,\  (\sum_{k=1}^r N_k^{\prime}) +N_- + (\sum^r_{k=i+1} N_k^{\prime} )+1 \leq   v \leq (\sum_{k=1}^r N_k^{\prime}) + N_- + \sum^r_{k=i} N_k^{\prime}  \Big\}.
\end{eqnarray*}

Put 
\[
Y = \big\{ t_u t^{-1}_{N-u+1}, \,\    1 \leq u \leq \sum_{k=1}^r N_k^{\prime}       \big\}.
\]

Then we see that for $\beta \in \widehat{\Sigma}^r_P(1)$ or $\widehat{\Sigma}^r_P(2)$, and $\alpha \in \Sigma_{\beta}$, the $\Gamma_F$-orbit $\{\alpha\}$ has exactly two elements (of the form $\{t_u t_v^{-1}, t_{N-v+1} t_{N-u+1}^{-1} \}$) and is contained in $\Sigma_{\beta}$. It follows that $n_{\beta} = |\Sigma_{\beta}|/2$
\[
\det \rho_{\beta} (1 \rtimes u) = (-1)^{|\Sigma_{\beta}|/2} = (-1)^{n_{\beta}}.
\]  

Finally for $\beta \in \widehat{\Sigma}^r_P(3)$, and $\alpha \in \Sigma_{\beta}$, we similarly see that the $\Gamma_F$-orbit $\{\alpha\}$ has exactly two elements if and only if $\alpha \in \Sigma_{\beta} \backslash ( \Sigma_{\beta} \cap Y)$, in which case $\{\alpha\}$ is contained in $\Sigma_{\beta}$ (and the orbit is again of the form $\{t_u t_{v}^{-1}, t_{N-v+1} t_{N-u+1}^{-1} \}$). It follows that $n_{\beta} = |\Sigma_{\beta} \backslash ( \Sigma_{\beta} \cap Y)|/2$, and 
\[
\det \rho_{\beta}(1 \rtimes u) =(-1)^{|\Sigma_{\beta} \backslash (\Sigma_{\beta}  \cap Y)|/2 } =(-1)^{n_{\beta}}
\]
as required.

\end{proof}

\bigskip

With proposition 3.3.1 and lemma 3.3.4, we now define:
\begin{eqnarray}
R_P(w,\pi,\psi): \mathcal{H}_P(\pi) \longrightarrow \mathcal{H}_{w^{-1}P}(\pi) \longrightarrow \mathcal{H}_P(w \pi)
\end{eqnarray}
as the composition
\begin{eqnarray}
R_P(w,\pi,\psi) =  l(w,\pi,\psi)  \circ  R_{w^{-1}P|P}(\pi,\psi).
\end{eqnarray}

Then the operator $R_P(w,\pi,\psi)$ has a decomposition as follows: 
\begin{eqnarray}
R_P(w,\pi_{\lambda},\psi_{\lambda}) = r_P(w,\psi_{\lambda})^{-1} J_P(\widetilde{w},\pi_{\lambda})
\end{eqnarray} 
where 
\[
J_P(\widetilde{w},\pi_{\lambda}) = l(\widetilde{w},\pi_{\lambda}) J_{w^{-1}P|P}(\pi_{\lambda})
\]
is the unnormalized intertwining operator from $\mathcal{H}_P(\pi)$ to $\mathcal{H}_P(w \pi)$, and $r_P(w,\psi_{\lambda})$ is the normalizing factor:
\[
r_P(w,\psi_{\lambda}) = \lambda(w) \epsilon_P(w,\psi_{\lambda})^{-1} r_{w^{-1}P|P}(\psi_{\lambda}).
\]
We have
\begin{eqnarray}
\end{eqnarray}
\begin{eqnarray}
 r_P(w,\psi_{\lambda}) & = & \lambda(w)  \epsilon(0, \pi_{\psi^{\widetilde{M}^0},\lambda}, \widetilde{\rho}^{\vee}_{w^{-1}P|P},\psi_F)^{-1} \nonumber \\ &  \times &  L(0, \pi_{\psi^{\widetilde{M}^0},\lambda}, \widetilde{\rho}^{\vee}_{w^{-1}P|P})   L(1, \pi_{\psi^{\widetilde{M}^0},\lambda}, \widetilde{\rho}^{\vee}_{w^{-1}P|P})^{-1}. \nonumber 
\end{eqnarray}

Finally by proposition 3.3.1 and lemma 3.3.4 we have the cocycle property, which we state as:
\begin{proposition}
\begin{eqnarray}
R_P(w^{\prime} w ,\pi,\psi) = R_P(w^{\prime}, w \pi,w \psi) \circ R_P(w,\pi,\psi)
\end{eqnarray}
for $w^{\prime},w \in w(M)$.
\end{proposition}
\subsection{The local intertwining relation, part I}

In this subsection we give the statement of the local intertwining relation, following section 2.4 of \cite{A1}, and is a supplement of theorem 3.2.1. Among other things, the local intertwining relation reduces the construction of the packets to the case of square integrable parameters.

We revert back to the notation before section 3.3. Thus $\pi, \phi,\psi$ refer to objects for $G=U_{E/F}(N)$, while we denote by $\pi_M,\phi_M,\psi_M$ the objects for a Levi subgroup $M$ of $G$. As in section 3.3, we fix $P \in \mathcal{P}(M)$. The choice of $P$ allows us to identify
\begin{eqnarray}
W(M) \cong W(\widehat{M}).
\end{eqnarray}

As in section 3.3 we assume that $M$ is a proper Levi subgroup of $G$. Given $\psi_M \in \Psi(M)$, we can assume by induction as in section 3.3 that the packet $\Pi_{\psi_M}$ is already defined. We denote by $\psi \in \Psi(G)$ the parameter of $G$ obtained by composing $\psi_M$ with the $L$-embedding $\leftexp{L}{M} \rightarrow \leftexp{L}{G}$ that is dual to $M \hookrightarrow G$.

Recall the groups
\begin{eqnarray*}
& & S_{\psi}(G)=S_{\psi} \\
& & \overline{S}_{\psi}(G)=\overline{S}_{\psi} \\
& & \mathcal{S}_{\psi}(G) = \mathcal{S}_{\psi}
\end{eqnarray*}
associated to $\psi \in \Psi(G)$ at the end of section 2.2. Similarly we have
\begin{eqnarray*}
& & S_{\psi_M}(M) = S_{\psi_M} \\
& & \overline{S}_{\psi_M}(M) = \overline{S}_{\psi_M} \\
& & \mathcal{S}_{\psi_M}(M)=\mathcal{S}_{\psi_M}.
\end{eqnarray*}

Put 
\[
A_{\widehat{M}} = (Z(\widehat{M})^{\Gamma_F})^0.
\]
Then we have
\[
S_{\psi_M}(M) = \Cent(A_{\widehat{M}}, S_{\psi}(G) )
\]
and $\mathcal{S}_{\psi_M}(M)$ is the image of $S_{\psi_M}(M)$ in $\mathcal{S}_{\psi}(G)$. 

Recall some other finite groups as in section 2.4 of \cite{A1}. Put
\begin{eqnarray*}
& & N_{\psi}(G,M) = \Norm(A_{\widehat{M}},S_{\psi}) \\
& & \overline{N}_{\psi}(G,M) =N_{\psi}(G,M)/Z(\widehat{G})^{\Gamma_F} \\
& & \mathfrak{N}_{\psi}(G,M) = \pi_0(\overline{N}_{\psi}(G,M)).
\end{eqnarray*}
Then $S_{\psi_M}(M) = \Cent(A_{\widehat{M}},S_{\psi}(G))$ is a normal subgroup of $N_{\psi}(G,M)$, and similarly $\mathcal{S}_{\psi_M}(M)$ is a normal subgroup of $\mathfrak{N}_{\psi}(G,M)$. Put
\[
W_{\psi}(G,M) = W(S_{\psi},A_{\widehat{M}})
\]
the Weyl group of the pair $(S_{\psi},A_{\widehat{M}})$, i.e. the group of automorphisms of $A_{\widehat{M}}$ induced by conjugation by $S_{\psi}$. Then we have
\begin{eqnarray*}
W_{\psi}(G,M) &=&  \Norm(A_{\widehat{M}},S_{\psi})/\Cent(A_{\widehat{M}},S_{\psi}) \\
& =& \pi_0 \big( \Norm(A_{\widehat{M}},S_{\psi}) / \Cent(A_{\widehat{M}},S_{\psi}) \big) \\
&=& \mathfrak{N}_{\psi}(G,M)/\mathcal{S}_{\psi_M}.
\end{eqnarray*}
Note that we have $w \psi_M = \psi_M$ for $w \in W_{\psi}(G,M)$. Also from the definition the elements of $W_{\psi}(G,M)$ stabilize $A_{\widehat{M}}$, hence normalize $\widehat{M}$. Thus $W_{\psi}(G,M)$ can be regarded naturally as a subgroup of $W(\widehat{M})$, hence a subgroup of $W(M)$ by the identification (3.4.1) given by $P$. 

\noindent We also put 
\[
W^0_{\psi}(G,M) = W(S_{\psi}^0,A_{\widehat{M}}) 
\]
the Weyl group of the pair $(S_{\psi}^0,A_{\widehat{M}} )$. Then $W_{\psi}^0(G,M)$ is naturally a normal subgroup of $\mathfrak{N}_{\psi}(G,M)$, and also of $W_{\psi}(G,M)$ (i.e. the map $W_{\psi}^0(G,M) \rightarrow \mathfrak{N}_{\psi}(G,M) \rightarrow W_{\psi}(G,M)$ is injective). Put
\begin{eqnarray*}
& & R_{\psi}(G,M) = W_{\psi}(G,M)/W_{\psi}^0(G,M) \\
& & \mathcal{S}_{\psi}(G,M) = \mathfrak{N}_{\psi}(G,M)/W_{\psi}^0(G,M).
\end{eqnarray*}
Then $\mathcal{S}_{\psi}(G,M)$ is a subgroup of $\mathcal{S}_{\psi}$, and $\mathcal{S}_{\psi_M}(M)$ is a normal subgroup of $\mathcal{S}_{\psi}(G,M)$, whose quotient is given by $R_{\psi}(G,M)$. The relations between these groups can be summarized by the commutative diagram of short exact sequences of finite groups, which plays a crucial role in the analysis of trace formulas in section five and six:

\begin{eqnarray}
\xymatrix{& & 1 \ar[d] & 1 \ar[d] \\
 & & W_{\psi}^0(G,M) \ar@{=}[r] \ar[d]  & W_{\psi}^0(G,M) \ar[d] \\
1  \ar[r] & \mathcal{S}_{\psi_M}(M) \ar[r] \ar@{=}[d] & \mathfrak{N}_{\psi}(G,M)  \ar[r] \ar[d]  \ar@<.8 ex>@{<--}[d] & W_{\psi}(G,M) \ar[r] \ar[d] \ar@<.8 ex>@{<--}[d] & 1 \\
1 \ar[r] &  \mathcal{S}_{\psi_M}(M) \ar[r] & \mathcal{S}_{\psi}(G,M) \ar[r] \ar[d] & R_{\psi}(G,M) \ar[r] \ar[d] & 1 \\ 
& & 1 & 1 
  }
\end{eqnarray}

The two vertical exact sequence splits, and the dotted arrows stand for splittings determined by the chamber of $\overline{P}_{\psi}$ in the Lie algebra of $A_{\widehat{M}}$, where $\overline{P}_{\psi}$ is the parabolic subgroup of $\overline{S}^0_{\psi}=S^0_{\psi}/Z(\widehat{G})^{\Gamma_F}$ given by
\[
\overline{P}_{\psi} =( \widehat{P} \cap S^0_{\psi})/Z(\widehat{G})^{\Gamma_F}.
\]
In other words, the image of the splitting to $\mathfrak{N}_{\psi}(G,M) \rightarrow \mathcal{S}_{\psi}(G,M)$ determined by $\overline{P}_{\psi}$ is characterized as the elements of $\mathfrak{N}_{\psi}(G,M)$ whose conjugation action on $\overline{S}_{\psi}$ preserves $\overline{P}_{\psi}$.

Given $u \in \mathfrak{N}_{\psi}(G,M)$, write $w_u$ and $x_u$ for the image of $u$ in $W_{\psi}(G,M)$ and $\mathcal{S}_{\psi}(G,M)$ under the (3.4.2) respectively.  

We have the following action of $\mathfrak{N}_{\psi}(G,M)$ on $M$: an element $u \in \mathfrak{N}_{\psi}(G,M)$ acts on $M$ via
\begin{eqnarray}
m \rightarrow \Inter(\widetilde{w}_u)(m)= \widetilde{w}_u m \widetilde{w}^{-1}_u
\end{eqnarray}
where as in section 3.3 $\widetilde{w}_u$ is the lift of the element $w_u$ (regarded as an element in $W(M)$) to $G(F)$ that preserves the splitting of $M$. The action of $\mathfrak{N}_{\psi}(G,M)$ on $M$ thus by definition factors through $W_{\psi}(G,M)$. In particular for $u \in \mathfrak{N}_{\psi}(G,M)$, we can form the twisted group:
\begin{eqnarray}
\widetilde{M}_u = M \rtimes \widetilde{w}_u.
\end{eqnarray} 

Define the twisted centralizer
\begin{eqnarray}
\widetilde{S}_{\psi_M,u} &=& S_{\psi_M}(\widetilde{M}_u) :=\Cent(\Image \psi_M, \widehat{\widetilde{M}}_u   )
\end{eqnarray}
with $\widehat{\widetilde{M}}_u : = \widehat{M} \rtimes \widetilde{w}_u$. Since $w_u \psi_M = \psi_M$, the twisted centralizer $\widetilde{S}_{\psi_M,u}$ is non-empty and hence a $S_{\psi_M}$ bi-torsor. Put
\[
\widetilde{\mathcal{S}}_{\psi_M,u} = \mathcal{S}_{\psi_M}(\widetilde{M}_u) :=\widetilde{S}_{\psi_M,u} / S^0_{\psi_M} Z(\widehat{M})^{\Gamma_F}.
\]
We can then naturally identify:
\[
\mathcal{S}_{\psi_M,u} \cong \mathfrak{N}_{\psi}(G,M)(w_u)
\]
where we have denoted by $\mathfrak{N}_{\psi}(G,M)(w_u)$ the fibre of $\mathfrak{N}_{\psi}(G,M)$ above $w_u$ under the middle horizontal short exact sequence of (3.4.2).
We denote by $\widetilde{u}$ the element of $\widetilde{\mathcal{S}}_{\psi_M,u}$ corresponding to $u \in \mathfrak{N}_{\psi}(G,M)(w_u)$ under this identification (i.e. $\widetilde{u}$ is the element $1 \rtimes \widetilde{w}_u$ of $\widetilde{\mathcal{S}}_{\psi_M,u}$).

Suppose that $\pi_M \in \Pi_{\psi_M}$ is stable under $\widetilde{w}_u$, i.e. $\pi_M$ extends to a representation $\widetilde{\pi}_M$ on the twisted group $\widetilde{M}_u$. Recall since we are assuming the existence of the packet $\Pi_{\psi_M}$ (as a consequence of the induction hypothesis) we have the pairing:
\begin{eqnarray}
\langle \cdot , \cdot \rangle:   \mathcal{S}_{\psi_M} \times \Pi_{\psi_M} \rightarrow \{\pm 1\}. 
\end{eqnarray}
We claim that we have a canonical extension of this pairing to 
\begin{eqnarray}
\langle \cdot , \cdot \rangle : \widetilde{\mathcal{S}}_{\psi_M,u} \times  \widetilde{\Pi}_{\psi_M,u} \rightarrow \{ \pm 1 \}     
\end{eqnarray}
with $\widetilde{\Pi}_{\psi_M,u} \subset \Pi_{\psi_M}$ being the subset consisting of those $\pi_M \in \Pi_{\psi_M}$ that are stable under $\widetilde{w}_u$.

\bigskip

\noindent We first note that any $\pi_M \in \widetilde{\Pi}_{\psi_M,u}$ has a {\it canonical} extension $\widetilde{\pi}_M$ to $\widetilde{M}_u$. Indeed, we may assume that $M$ has the form:
\begin{eqnarray*}
& & M= G_{E/F}(N_1^{\prime}) \times \cdots \times G_{E/F}(N_r^{\prime}) \times G_- \\
& & G_- = U_{E/F}(N_-), \,\ N_- < N
\end{eqnarray*}
with respect to which $\pi_M$ has a decomposition:
\begin{eqnarray*}
& & \pi_M = \pi_1^{\prime} \boxtimes \cdots \pi_r^{\prime} \boxtimes \pi_- \\
& & \psi_M = \psi_1^{\prime} \times \cdots \times \psi_r^{\prime} \times \psi_-.
\end{eqnarray*}

\bigskip

\noindent For any $u \in \mathfrak{N}_{\psi}(G,M)$, we see that the action of $\widetilde{w}_u$ on $M$ has to preserve the $U_{E/F}(N_-)$ factor of $M$. Furthermore, from the fact that
\begin{eqnarray*}
 \widetilde{\Out}_{N_-}(U_{E/F}(N_-))=1
\end{eqnarray*}
we see that the action of $\widetilde{w}_u$ on the $G_-$ factor is given by inner automorphism of $G_-$. In particular the the action of the lift $\widetilde{w}_u$ commutes with the factor $G_-$. Hence we only need to extend the representation $\pi_1^{\prime} \boxtimes \cdots \boxtimes \pi_r^{\prime}$ to the twisted group $\widetilde{G}_{E/F}(N_1^{\prime}) \times \cdots \times \widetilde{G}_{E/F}(N_r^{\prime})$. This can be done canonically using the Whittaker normalization as in section 3.2 by using the standard Whittaker data on the group $G_{E/F}(N_1^{\prime}) \times \cdots \times G_{E/F}(N_r^{\prime})$ (note that by construction $\widetilde{w}_u$ preserves the Whittaker data of $M$ and hence of $G_{E/F}(N_1^{\prime}) \times \cdots \times G_{E/F}(N_r^{\prime})$). In particular for $u \in \mathfrak{N}_{\psi}(G,M)$ we have an intertwining operator $\widetilde{\pi}_M(w_u): w_u \pi_M \rightarrow \pi_M$, thus giving the extension $\widetilde{\pi}_M$ of $\pi_M$ to $\widetilde{M}_u$. 

Thus we can identify $\widetilde{\Pi}_{\psi_M,u}$ as the packet of representations on the twisted group $\widetilde{M}_u$. We also note that the Whittaker normalization ensures that for $u_1,u_2 \in \mathfrak{N}_{\psi}(G,M)$ satisfying $w_{u_i} \pi_M \cong \pi_M$ ($i=1,2$), we have $\widetilde{\pi}_M(w_{u_1 u_2})=\widetilde{\pi}_M(w_{u_1}) \widetilde{\pi}_{M}(w_{u_2})$. 

\bigskip

\noindent It remains to extend the pairing (3.4.6) to (3.4.7). Again from the decomposition $M=G_{E/F}(N_1^{\prime}) \times \cdots \times G_{E/F}(N_r^{\prime}) \times U_{E/F}(N_-)$, and the fact that the twisted centralizer of $\psi_i^{\prime}$ for the general linear factors $G_{E/F}(N_i^{\prime})$ are connected, we can make the $\mathcal{S}_{\psi_M}=\mathcal{S}_{\psi_-}(U_{E/F}(N_-))$-equivariant identification $\widetilde{S}_{\psi_M,u}  \cong \mathfrak{N}_{\psi}(G,M)(w_u) \cong \mathcal{S}_{\psi_M}=\mathcal{S}_{\psi_-}(U_{E/F}(N_-))$ (again we are using the fact that the action of $\widetilde{w}_u$ on $G_-=U_{E/F}(N_-)$ is given by inner automorphism). This gives the pairing (3.4.7), namely if $\widetilde{u}^{\prime} \in \mathcal{S}_{\psi_M,u}$ corresponds to an element $ x_{u^{\prime}} \in \mathcal{S}_{\psi_-}(U_{E/F}(N_-))$ under this identification, then we have:
\[
\langle \widetilde{u}^{\prime},\widetilde{\pi}_M \rangle = \langle x_{u^{\prime}}  , \pi_- \rangle
\]
where the pairing on the right hand side is given by the packet $\Pi_{\psi_-}$ associated to $\psi_- \in \Psi(U_{E/F}(N_-))$.

\begin{rem}
\end{rem}
In \cite{A1} Arthur needs to consider the case where $G=SO(2n)$. For the case $G=SO(2n)$ then in general one has no canonical extension of $\pi_M$ to the twisted group $\widetilde{M}_u$, and thus in \cite{A1} one needs to consider the case of twisted endoscopy for the group $SO(2N_-)$ (here $N_- \leq N$), with the outer automorphism being given by conjugation by $O(2N_-)$. This technical problem does not arise for the groups $SO(2N+1),Sp(2N)$ or $U_{E/F}(N)$.

\bigskip

Thus given $u \in \mathfrak{N}_{\psi}(G,M)$, and $\pi_M \in \widetilde{\Pi}_{\psi_M,u}$, we have the intertwining operator $\widetilde{\pi}_M(w_u) :  w_u \pi_M \rightarrow \pi_M$, and we continue to denote by $\widetilde{\pi}_M(w_u)$ the intertwining operator $\mathcal{I}_P( w_u \pi_M) \rightarrow \mathcal{I}_P( \pi_M)$ induced by $\widetilde{\pi}_M$. Put:
\begin{eqnarray}
R_P(w_u,\widetilde{\pi}_M,\psi_M) :=  \widetilde{\pi}_M(w_u)  \circ R_P(w_u,\pi_M,\psi_M).
\end{eqnarray}

If we put 
\[
\mathfrak{N}_{\psi,\pi_M}(G,M) :=\{ u \in \mathfrak{N}_{\psi}(G,M), \,\ w_u \pi_M \cong \pi_M\}
\]
then by combining the discussion above with proposition 3.3.5, we see that the map:
\[
u \mapsto R_P(w_u,\widetilde{\pi}_M,\psi_M)
\]
is a group homomorphism from $\mathfrak{N}_{\psi,\pi_M}(G,M)$ to the group of unitary intertwining operators of $\mathcal{I}_P(\pi_M)$.

\bigskip

\begin{definition}
For $u \in \mathfrak{N}_{\psi}(G,M)$, we define the linear form:
\begin{eqnarray*}
 f \rightarrow f_G(\psi,u), \,\ f \in \mathcal{H}(G) 
\end{eqnarray*}
\begin{eqnarray}
& &  f_G(\psi,u) \\
& :=& \sum_{\pi_M \in \widetilde{\Pi}_{\psi_M,u}} \langle \widetilde{u} ,  \widetilde{\pi}_M \rangle \tr(R_P(w_u,\widetilde{\pi}_M,\psi_M) \mathcal{I}_P(\pi_M,f)    ). \nonumber
\end{eqnarray}
\end{definition}

Following the convention of \cite{A1}, if we interpret the term $\langle \widetilde{u},\widetilde{\pi}_M   \rangle$ as being equal to zero if $\pi_M$ is not stable under $w_u$, then we can just  write (3.4.9) as 
\begin{eqnarray}
& & \\
& & f_G(\psi,u) := \sum_{\pi_M \in \Pi_{\psi_M}} \langle \widetilde{u} ,  \widetilde{\pi}_M \rangle \tr(R_P(w_u,\widetilde{\pi}_M,\psi_M) \mathcal{I}_P(\pi_M,f)    ). \nonumber
\end{eqnarray}

With the previous notations, for $s \in S_{\psi}$ a semi-simple element, we let $(G^{\prime},\psi^{\prime})$ the pair corresponding to $(\psi,s)$ under the correspondence as described in section 3.2. Assume that part (a) of theorem 3.2.1 is valid for the pair $(G^{\prime},\psi^{\prime})$, i.e. assuming the existence of the stable linear form
\[
f^{\prime} \rightarrow \rightexp{f^{\prime}}{G^{\prime}}(\psi^{\prime}), \,\ f^{\prime} \in \mathcal{H}(G^{\prime})
\]
satisfying (3.2.8). Define the linear form
\begin{eqnarray}
& & f \rightarrow f_G^{\prime}(\psi,s), \,\ f \in \mathcal{H}(G) \\
& & f^{\prime}_G(\psi,s) =f^{\prime} (\psi^{\prime}) := f^{G^{\prime}} (\psi^{\prime}). \nonumber
\end{eqnarray}

We can now state:
\begin{theorem} (The local intertwining relation, {\it c.f.} theorem 2.4.1 of \cite{A1}) 

\bigskip

\noindent Given any $u \in \mathfrak{N}_{\psi}(G,M)$ and $s \in S_{\psi}$, such that the image of $s$ in $\mathcal{S}_{\psi}$ is equal to the image $x_u$ of $u$ in $\mathcal{S}_{\psi}(G,M)$, we have the identity
\begin{eqnarray}
f_G(\psi,u) = f^{\prime}_G(\psi,s_{\psi} s), \,\ f \in \mathcal{H}(G).
\end{eqnarray}
\end{theorem}
In particular the local intertwining relation implies that the linear form $f_G(\psi,u)$ depends only on the image $x_u$ of $u$ in $\mathcal{S}_{\psi}(G,M)$. As with theorem 3.2.1, the local intertwining relation (3.4.12) is to be proved by a long induction argument. The induction argument of the proof of the local intertwining relation will be completed in section eight. However, special cases of the local intertwining relation have to be proved along the way in section five, which in particular forms an important ingredient in the comparison of trace formulas in section five and six. 
\bigskip

The most important case is when $\psi_M \in \Psi_2(M)$, the subset of square integrable parameters in $\Psi(M)$, i.e. the set of $\psi_M$ such that $\overline{S}_{\psi_M}$ is finite. Then $M$ is determined by $\psi$ up to conjugation by $G$ ({\it c.f.} Proposition 8.6 of \cite{B}). In this case $T_{\psi}:=A_{\widehat{M}}$ is then a maximal torus in $S_{\psi}^0$, with $\widehat{M}=\Cent(T_{\psi},\widehat{G})$, and $\mathcal{S}_{\psi}(G,M)$ is equal the full group $\mathcal{S}_{\psi}$. 

We also make the following abbreviation (here we put $\overline{T}_{\psi} = T_{\psi}/Z(\widehat{G})^{\Gamma_F}$):
\begin{eqnarray*}
& & \mathfrak{N}_{\psi} = \mathfrak{N}_{\psi}(G,M) =  \Norm(\overline{T}_{\psi}, \overline{S}_{\psi})/\overline{T}_{\psi} \\
& & W_{\psi} = W_{\psi}(G,M) = \Norm(\overline{T}_{\psi},\overline{S}_{\psi})/\Cent(\overline{T}_{\psi},\overline{S}_{\psi})  \\
& & W_{\psi}^0 = W_{\psi}^0(G,M) = \Norm( \overline{T}_{\psi},\overline{S}_{\psi}^0)/\overline{T}_{\psi}\\
& & R_{\psi} = R_{\psi}(G,M) \\
& &  S_{\psi}^1 = S_{\psi_M} = \Cent(T_{\psi},S_{\psi}) \\
& & \mathcal{S}_{\psi}^1 = \mathcal{S}_{\psi_M}.
\end{eqnarray*}  

In this basic case where $\psi_M \in \Psi_2(M)$, we have that any $\pi_M \in \Pi_{\psi_M}$ is stable under the group $W_{\psi}$ and hence $\mathfrak{N}_{\psi}$. This follows immediately from the discussion after equation (3.4.7), due to the fact that the components of the representation $\pi_M$ associated to the general linear factors are determined by their parameters. Hence each $\pi_M$ has a canonical extension $\widetilde{\pi}_M$ to the group $M(F) \rtimes \mathfrak{N}_{\psi}$. Thus in this basic case we see that the map:
\[
u \mapsto R_P(w_u,\widetilde{\pi}_M,\psi_M)
\]
is a group homomorphism from $\mathfrak{N}_{\psi}$ to the group of unitary intertwining operators of $\mathcal{I}_P(\pi_M)$.

Furthermore the bottom exact sequence of (3.4.2) of abelian groups:
\begin{eqnarray}
1 \rightarrow \mathcal{S}^1_{\psi} \rightarrow \mathcal{S}_{\psi} \rightarrow R_{\psi} \rightarrow 1
\end{eqnarray}
splits canonically in this case. Indeed these groups can be explicated as follows. Following the local form of the notation (2.4.12), given a parameter $\psi \in \Psi(G)=\Psi(U_{E/F}(N))$, we denote by $\psi^N \in \widetilde{\Psi}(N) = \Psi(\widetilde{G}_{E/F}(N))$ the parameter $\psi^N = \xi \circ \psi$, where we recall as in the last subsection we chose $\xi = \xi_{\chi_{\kappa}} : \leftexp{L}{U_{E/F}(N)} \rightarrow \leftexp{L}{G_{E/F}(N)}$ so that $(U_{E/F}(N),\xi_{\chi_{\kappa}})$ is a simple twisted endoscopic datum of $\widetilde{G}_{E/F}(N)$. By lemma 2.2.1 $\psi^N$ is conjugate self dual with parity $ (-1)^{N-1} \cdot \kappa$. We can then write
\begin{eqnarray*}
\psi^N = \bigoplus_{i \in I_{\psi}^+} l_i \psi_i^{N_i} \,\ \oplus \,\ \bigoplus_{i \in I_{\psi}^-}  l_i \psi_i^{N_i} \,\ \oplus \,\ \bigoplus_{j \in J_{\psi}} l_j (\psi_j^{N_j} \oplus \psi_{j^*}^{N_{j^*}} )
\end{eqnarray*}
with $\psi_i^{N_i} \in \widetilde{\Psi}_{\simp}(N_i),\psi^{N_j}_j \in \widetilde{\Psi}_{\simp}(N_j)$ being simple parameters (i.e. irreducible as representations of $L_E \times \SU(2)$), and where $I_{\psi}^+$ is the set of indices $i$ such that $\psi_i^{N_i}$ is conjugate self-dual of the same parity as $\psi^N$ (i.e. with parity $ (-1)^{N-1} \cdot \kappa$), $I_{\psi}^-$ is the set of indices $i$ such that $\psi_i^{N_i}$ is conjugate self-dual with parity different from $\psi^N$ (i.e. has parity $ (-1)^N \cdot \kappa$), and $J_{\psi}$ is the set of indices $j$ such that $\psi_j^{N_j}$ is not conjugate self-dual. 

As in (2.4.14) we have
\begin{eqnarray}
& & \\
& & S_{\psi} = \Big( \prod_{i \in I_{\psi}^+} O(l_i,\mathbf{C}) \Big) \times \Big( \prod_{i \in I^-_{\psi}} Sp(l_i,\mathbf{C}) \Big) \times  \Big(  \prod_{j \in J_{\psi}}  \GL(l_j,\mathbf{C})  \Big). \nonumber
\end{eqnarray}
From this it follows that we can identify $\pi_0(S_{\psi})$ as the group $\Sigma$ consisting of functions
\[
\sigma:I_{\psi}^+ \rightarrow  \mathbf{Z}/2 \mathbf{Z}.  
\]
Since $Z(\widehat{G})^{\Gamma_F} = \{ \pm I_N\}$, with $I_N$ being identified in the obvious way as an element on the right hand side of (3.4.14), we have
\begin{eqnarray}
\mathcal{S}_{\psi} &=& \pi_0(\overline{S}_{\psi}) \\
&=& \overline{\Sigma} := \Sigma  /   \langle \overline{\sigma} \rangle \nonumber
\end{eqnarray}
where $\overline{\sigma} \in \Sigma$ is the function such that $\overline{\sigma}(i) =-1$ for $i \in I_{\psi,o}^+$, and $\overline{\sigma}(i)=1$ for $i \in I_{\psi,e}^+$. Here 
\begin{eqnarray}
& & I_{\psi,o}^+ = \{i \in I_{\psi}^+, \,\ l_i \mbox{ odd} \} \\
& & l_{\psi,e}^+ = \{ i \in I_{\psi}^+,\,\ l_i \mbox{ even} \}. \nonumber
\end{eqnarray}  

Recall that in the basic case $\psi_M \in \Psi_2(M)$ the group $A_{\widehat{M}}$ is a maximal torus in $S_{\psi}$. Hence we can take $T_{\psi} = A_{\widehat{M}}$ to be a maximal torus of the product of groups on the right hand side of (3.4.14). Since $S_{\psi_M} = \Cent(A_{\widehat{M}},S_{\psi})$, we see that $\pi_0(S_{\psi_M})$ corresponds to the subgroup $\Sigma^{1} \subset \Sigma$ consisting of functions $\sigma \in \Sigma$ that are supported on $I_{\psi,o}^+$. Hence if we denote by $\overline{\Sigma}^1$ the image of $\Sigma^1$ in $\overline{\Sigma}$, then we have $\mathcal{S}^1_{\psi} = \mathcal{S}_{\psi_M}\cong \overline{\Sigma}^1$. The quotient of $\overline{\Sigma}$ by $\overline{\Sigma}^1$ is canonically isomorphic to the group $R$ consisting of functions
\[
\rho: I_{\psi,e}^+ \rightarrow \mathbf{Z}/2\mathbf{Z}.
\]
Hence $R_{\psi} \cong R$, and the exact sequence (3.4.13) is just given by:
\begin{eqnarray}
0 \rightarrow \overline{\Sigma}^1 \rightarrow \overline{\Sigma} \rightarrow R \rightarrow 0.
\end{eqnarray}
It is clear that the exact sequence (3.4.17) splits canonically, hence the exact sequence (3.4.13) splits canonically in the basic case that $\psi_M \in \Psi_2(M)$. 

\bigskip

\begin{proposition}
Let $\psi \in \Psi(G)$ be a parameter in the complement of $\Psi_2(G)$, and assume that the linear form $f^{\prime}(\psi,s)$ exists for any $s \in S_{\psi}$. In addition assume that for any proper Levi subgroup $M$ of $G$, theorem 3.4.3 holds for parameters in $\Psi_2(M)$. Then the packet $\Pi_{\psi}$ and the pairing $\langle x ,\pi \rangle$ of part (b) of theorem 3.2.1 exist, and satisfies (3.2.11).
\end{proposition}
\begin{proof}
Since the proof is the same as the proof of proposition 2.4.3 of \cite{A1} we only give a sketch. By assumption $\psi \notin \Psi_2(G)$, so we can choose a proper Levi subgroup $M$ of $G$, and parameter $\psi_M \in \Psi_2(M)$ such that $\psi_M$ maps to $\psi$. Since $M$ is proper we can assume by induction hypothesis as above that the packet $\Pi_{\psi_M}$ is already constructed. For each $\pi_M \in \Pi_M$, we denote by $\xi_M$ the character of $\mathcal{S}^1_{\psi}$ associated to $\pi_M$, namely $\xi_M(x_M) = \langle x_M,\pi_M\rangle$. Recall that in this basic case any $\pi_M \in \Pi_{\psi_M}$ is invariant under $\mathfrak{N}_{\psi}$.

Abusing notation slightly, we still denote by $\Pi_{\psi_M}$ the following (reducible) admissible unitary representation of $\mathcal{S}^1_{\psi} \times M(F)$:
\begin{eqnarray}
\Pi_{\psi_M} := \bigoplus_{(\xi_M,\pi_M)} (\xi_M \otimes \pi_M)
\end{eqnarray}
(here $\pi_M$ can occur with multiplicities; as we only know {\it a priori} that $\Pi_{\psi_M}$ is a multi-set).

We claim that the representation $\Pi_{\psi_M}$ of $\mathcal{S}^1_{\psi} \times M(F)$ can be extended to a representation of $M(F) \rtimes \mathfrak{N}_{\psi}$. It suffices to work with each component $\xi_M \otimes \pi_M$ of $\Pi_{\psi_M}$. Indeed by (3.4.7) we have an extension $\widetilde{\xi}_M$ of $\xi_M$ to $\mathfrak{N}_{\psi}$ given by $\widetilde{\xi}_{M}(u)=\langle  \widetilde{u}, \widetilde{\pi}_M \rangle$. Hence it suffices to specify the action of $M(F) \rtimes \mathfrak{N}_{\psi}$ on $\pi_M \in \Pi_{\psi_M}$. This action then comes from the extension $\widetilde{\pi}_M$ of $\pi_M$ to $M(F) \rtimes \mathfrak{N}_{\psi}$ in the discussion after (3.4.7). 

Consider the parabolic induction of $\Pi_{\psi_M}$ from $\mathcal{S}^1_{\psi} \times M(F)$ to $\mathcal{S}^1_{\psi} \times G(F)$:
\begin{eqnarray}
& & \,\ \,\ \Pi_{\psi} := \mathcal{I}_P(\Pi_{\psi_M})\\
&=& \bigoplus_{(\xi_M,\pi_M)} (\xi_M \otimes \mathcal{I}_P(\pi_M)) \nonumber
\end{eqnarray}
which is a unitary representation of $\mathcal{S}^1 \times G(F)$. In particular we have the unitary action of $G(F)$ on $\Pi_{\psi}$ (with the action of $G(F)$ on the $\xi_M$'s being trivial). On the other hand, for each component $\pi_M$ of $\Pi_{\psi_M}$, we have the map
\[
u \mapsto R_P(w_u,\widetilde{\pi}_M,\psi_M) = \widetilde{\pi}(w_u) \circ R_P(w_u,\pi_M,\psi_M)
\]
which is a group homomorphism from $\mathfrak{N}_{\psi}$ (that factors through $W_{\psi}$) to the group of unitary intertwining operators on $\mathcal{I}_P(\pi_M)$. Hence we have an action of $\mathfrak{N}_{\psi}$ on $\Pi_{\psi}$: given $u \in \mathfrak{N}_{\psi}$ its action on the representation $\Pi_{\psi} = \bigoplus_{(\xi_M,\pi_M)} (\xi_M \otimes \mathcal{I}_P(\pi_M))$ is given by the operator:
\begin{eqnarray}
R_P(u,\widetilde{\Pi}_{\psi_M},\psi_M) = \bigoplus_{(\xi_M,\pi_M)} \widetilde{\xi}_M(u) \otimes R_P(w_u,\widetilde{\pi}_M,\psi_M).
\end{eqnarray}

By definition of intertwining operator this action of $\mathfrak{N}_{\psi}$ on $\Pi_{\psi}$ commutes with that of $G(F)$, hence we have an unitary action of $\mathfrak{N}_{\psi} \times G(F)$ on $\Pi_{\psi}$:
\begin{eqnarray}
& & \\
& & \Pi_{\psi}(u,g) =  R_P(u, \widetilde{\Pi}_{\psi_M},\psi_M) \circ \Pi_{\psi}(g)  , \,\ u \in \mathfrak{N}_{\psi}, g \in G(F). \nonumber
\end{eqnarray}

Now consider the restriction of this representation to the subgroup $\mathcal{S}_{\psi} \times G(F)$ of $\mathfrak{N}_{\psi} \times G(F)$ (under the splitting of the middle vertical exact sequence of (3.4.2) as determined by $P$). Then we have a decomposition:
\begin{eqnarray}
\Pi_{\psi} = \bigoplus_{(\xi,\pi)} (\xi \otimes \pi)
\end{eqnarray}
where $\xi$ is a character of $\mathcal{S}_{\psi}$, and $\pi$ is an irreducible unitary representation of $G(F)$ (which might occur with multiplicities). We then declare that the packet $\Pi_{\psi}$ associated to $\psi \in \Psi(G)$ is give by the $\pi$'s that occur in the decomposition (3.4.22) (with possible multiplicities), and the pairing
\[
\langle \cdot , \cdot \rangle : \mathcal{S}_{\psi} \times \Pi_{\psi} \rightarrow \{ \pm1 \}
\] 
being given by $\langle x,\pi \rangle = \xi(x)$, with $\xi$ being the character associated to $\pi$ in the decomposition $(3.2.23)$. It remains to verify the character relation (3.4.22). 

For this we finally use the local intertwining relation. First from the definition (3.4.20) and (3.4.21), we have for $f \in \mathcal{H}(G)$:
\begin{eqnarray}
& & \tr \Pi(u,f) \\
&=&   \sum_{\pi_M \in \Pi_{\psi_M}} \widetilde{\xi}_M(u) \tr( R_P(w_u,\widetilde{\pi}_M,\psi_M)   \mathcal{I}_P(\pi_M,f)    )    \nonumber \\
&=& \sum_{\pi_M \in \Pi_{\psi_M}} \langle \widetilde{u},\widetilde{\pi}_M  \rangle \tr( R_P(w_u,\widetilde{\pi}_M,\psi_M)   \mathcal{I}_P(\pi_M,f)    )    \nonumber \\
&=&  f_G(\psi,u). \nonumber
\end{eqnarray}
On the other hand, by the decomposition (3.4.22), we have
\begin{eqnarray}
\tr \Pi_{\psi}(x_u,f) = \sum_{\pi \in \Pi_{\psi}} \langle x_u, \pi \rangle f_G(\pi)
\end{eqnarray}

\noindent with $x_u$ being the image of $u$ in $\mathcal{S}_{\psi}$. The local intertwining relation implies that (3.4.23) depends only on the image $x_u$ of $u$ in $\mathcal{S}_{\psi}$. Hence for any semi-simple element $s \in S_{\psi}$ with image $x$ in $\mathcal{S}_{\psi}$, we have by the local intertwining relation (taking any element $u \in \mathfrak{N}_{\psi}$ mapping to $x$ in $\mathcal{S}_{\psi}$):
\begin{eqnarray}
 \sum_{\pi \in \Pi_{\psi}} \langle x,\pi \rangle f_G(\pi) =f_G(\psi,x) = f_G(\psi,u) = f^{\prime}_G(\psi,s_{\psi} s). 
\end{eqnarray} 
Making the substitution $s \rightarrow s_{\psi}^{-1} s = s_{\psi} s$, we have
\begin{eqnarray}
 \sum_{\pi \in \Pi_{\psi}} \langle s_{\psi} x,\pi \rangle f_G(\pi) =  f^{\prime}_G(\psi, s). 
\end{eqnarray} 
Since $f^{\prime}(\psi,s) = f^{G^{\prime} }(\psi^{\prime}) $ by definition we obtain the character relation (3.2.11) as required. 
\end{proof}

\bigskip

Finally we add the following lemma from \cite{A1}, whose proof applies verbatim:
\begin{lemma} (lemma 2.4.2 of \cite{A1})
Suppose that for any proper Levi subgroup $M$ of $G$, the local intertwining relation (3.4.12) is valid whenever $\psi \in \Psi_2(M)$. Then it holds for any $\psi \in \Psi(M)$.
\end{lemma}

\bigskip

\subsection{The local intertwining relation, part II}

In the previous section we have formulated the local intertwining relation for the case where $G=U_{E/F}(N)$ (more precisely, as in section 3.3 and 3.4, we regard $G=(G,\xi)$ as a simple twisted endoscopic datum of $\widetilde{G}_{E/F}(N)$, for a choice of $L$-embedding $\xi: \leftexp{L}{G} \hookrightarrow \leftexp{L}{G_{E/F}(N)}$). For the analysis of the spectral terms of the twisted trace formula for $G_{E/F}(N)$, it is necessary to formulate the corresponding intertwining relation for the twisted group $\widetilde{G}_{E/F}(N)$, which we do in this subsection. At the same time, we use the theory of Whittaker models for the induced representations to obtain information about the normalized intertwining operators. The reader can find a discussion of the notion of Whittaker model and Whittaker functionals in the beginning of section 2.5 of \cite{A1}, and so will not be repeated here. We also recall the fact that any irreducible tempered representation of $G_{E/F}(N)(F)=\GL_N(E)$, or more generally a product of $G_{E/F}(N_i)(F)$'s, is generic, i.e. possess a non-zero Whittaker functional.

Thus suppose that $M$ is a standard Levi subgroup of $G_{E/F}(N)$, (which is a product of $G_{E/F}(N_i)$'s), and that
\[
\psi^M : L_F \times \SU(2) \rightarrow \leftexp{L}{M}
\]
is a parameter in $\Psi(M)$. Then as in section 3.2 the parameter $\psi^M$ corresponds to an irreducible unitary representation $\pi_{\psi^M}$ of $M(F)$, which is the Langlands quotient of a standard representation $\rho_{\psi^M}$ (associated to the parameter $\phi_{\psi^M} \in \Phi(M)$). Similar to the previous subsection, we denote by $\psi^N \in \Psi(G_{E/F}(N))$ the parameter of $G_{E/F}(N)$ obtained by composing $\psi^M$ with the $L$-embedding $\leftexp{L}{M} \hookrightarrow \leftexp{L}{G_{E/F}(N)}$. We assume as usual that $\psi^N \in \widetilde{\Psi}(N)$. 

One has a natural formulation of the commutative diagram (3.4.2) in the present case, for both $G=G_{E/F}(N)$ and the twisted group $G=\widetilde{G}_{E/F}(N)$. For the case of $G=G_{E/F}(N)$ the discussion is the same as before, while for the the twisted case $G=\widetilde{G}_{E/F}(N)$, the four terms in the lower right corner have to be formulated in terms of $\widetilde{G}_{E/F}(N)$, by considering the twisted centralizer $\overline{\widetilde{S}}_{\psi^N}$, while the rest of the terms in (3.4.2) are formulated in terms of the usual centralizer $\overline{S}_{\psi^N}$ (recall that $\overline{\widetilde{S}}_{\psi^N}$ is a bi-torsor under $\overline{S}_{\psi^N}$). 

In fact since the centralizers $\overline{S}_{\psi^N}$ and $\overline{\widetilde{S}}_{\psi^N}$ are connected, it follows that all the terms in the last row of (3.4.2) are singleton, for both the case $G=G_{E/F}(N)$ or $G=\widetilde{G}_{E/F}(N)$.  In particular, for $G=\widetilde{G}_{E/F}(N)$, we have $\mathfrak{N}_{\psi^N}(G,M) = W_{\psi^N}(G,M)$, and both are torsors under the group $W^0_{\psi^N}(G,M) = W_{\psi^N}(G_{E/F}(N),M)$.

More generally denote by $\widetilde{W}(M)$ the Weyl set of outer automorphisms of $M$ induced from the component $\widetilde{G}_{E/F}(N)$, and $\widetilde{W}_{\psi^N}(M)$ for the stabilizer of $\psi^M$ in $\widetilde{W}(M)$. Then we have $\widetilde{W}_{\psi^N}(M) = W_{\psi^N}(\widetilde{G}_{E/F}(N),M)$, and any element of $\widetilde{W}_{\psi^N}(M)$ stabilizes $\pi_{\psi^M}$ or $\rho_{\psi^M}$. 

Given $w \in \widetilde{W}_{\psi^N}(M)$, we can write
\[
w = \theta(N) \circ w^0
\]
where $\theta(N) = \theta$ as in (1.0.1), and $w^0 \in W(M,M^{\prime})$. Here $W(M,M^{\prime})$ is the Weyl set of elements that conjugate $M$ to the other standard Levi subgroup $M^{\prime}$ of $G_{E/F}(N)$, and in the present context $M^{\prime}$ is the Levi subgroup that is paired with $M$ under the involution $\theta$ of $G_{E/F}(N)$. The representative 
\[
\widetilde{w} = \theta(N) \circ \widetilde{w}^0
\]
then preserves the standard Whittaker datum $(B_M,\chi_M)$ for $M$. 

\bigskip

Consider first the untwisted case $W_{\psi^N}(M) = W^0_{\psi^N}(G_{E/F}(N),M)$. Then for $w \in W_{\psi^N}(M)$, we can define the normalized intertwining operator
\[
R_P(w, \pi_{\psi^M}) = R_P(w, \pi_{\psi^M},\psi^N): \mathcal{I}_P(\pi_{\psi^M}) \rightarrow \mathcal{I}_P(w \pi_{\psi^M})      
\]
as before (here since $\psi^M$ and hence $\psi^N$ is determined by $\pi_{\psi^M}$ we suppress $\psi^N$ from the notation of the intertwining operator). We note that since $\pi_{\psi^M}$ is irreducible and unitary, the induced representation $\mathcal{I}_P(\pi_{\psi^M})$ of $G_{E/F}(N)(F)=\GL_N(E)$ is irreducible by the theorem of Bernstein \cite{Be}. In fact we see that $\mathcal{I}_P(\pi_{\psi^M})$ is the Langlands quotient corresponding to the parameter $\psi^N \in \widetilde{\Psi}(N)$. 

Going back to the twisted case, suppose that $w \in \widetilde{W}_{\psi^N}(M)$. We would like to define, for $P \in \mathcal{P}(M)$, the normalized twisted intertwining operator:
\begin{eqnarray}
\widetilde{R}_P(w,\pi_{\psi^M}) : \mathcal{H}_P(\pi_{\psi^M}) \rightarrow \mathcal{H}_P(w \pi_{\psi^M})
\end{eqnarray}
by a variant of (3.3.41) as follows. Namely we define:
\begin{eqnarray}
\widetilde{l}(\widetilde{w},\pi_{\psi^M}) : \mathcal{H}_{w^{-1}P}(\pi_{\psi^M}) \rightarrow \mathcal{H}_P(w \pi_{\psi^M})
\end{eqnarray}
by the rule (which is a modification of (3.3.22)):
\begin{eqnarray}
(\widetilde{l}(\widetilde{w},\pi_{\psi^M}) \phi)(x) = \phi(\widetilde{w}^{-1}  x \theta(N) )
\end{eqnarray}
for $\phi \in \mathcal{H}_{w^{-1} P}(\pi_{\psi^M})$ and $x \in G_{E/F}(N)(F)$. Then we put, similar to (3.3.26):
\begin{eqnarray}
& & \widetilde{l}(w,\pi_{\psi^M}) : \mathcal{H}_{w^{-1}P}(\pi_{\psi^M}) \rightarrow \mathcal{H}_P(w \pi_{\psi^M}) \\
& &  \widetilde{l}(w,\pi_{\psi^M}) =  \lambda(w)^{-1}  \epsilon_P(w,\psi^N) \widetilde{l}(\widetilde{w},\pi_{\psi^M}).    \nonumber
\end{eqnarray}
here the factors $\epsilon_P(w,\psi^N)$ and $\lambda(w)$ are defined by exactly the same formula as in (3.3.23) and (3.3.24).

\noindent We then define the twisted intertwining operator (3.5.1) by the composition:
\begin{eqnarray}
\widetilde{R}_P(w,\pi_{\psi^M}) :=  \widetilde{l}(w,\pi_{\psi^M}) \circ R_{w^{-1}P|P}(\pi_{\psi^M}).
\end{eqnarray}

\noindent Next, similar to the discussion in section 3.2 by consideration of Whittaker functional, the standard representation $\rho_{\psi^M}$ has a canonical intertwining operator
\[
\widetilde{\rho}_{\psi^M} : w \rho_{\psi^M} \rightarrow \rho_{\psi^M}.
\]

\noindent The intertwining operator $\widetilde{\rho}_{\psi^M}$ then in turn defines the intertwining operator $\widetilde{\pi}_{\psi^M}(w): w \pi_{\psi^M} \rightarrow \pi_{\psi^M}$ on the Langlands quotient. We put
\[
\widetilde{R}_P(w, \widetilde{\pi}_{\psi^M}) := \widetilde{\pi}_{\psi^M} \circ \widetilde{R}_P(w, \pi_{\psi^M})
\]
which is then a self-intertwining operator on $\mathcal{H}_P(\pi_{\psi^M})$. By construction we have:
\begin{eqnarray}
\widetilde{R}_P(w, \widetilde{\pi}_{\psi^M}) : \mathcal{I}_P(\pi_{\psi^M}) \rightarrow \mathcal{I}_P(\pi_{\psi^M}) \circ \theta(N)
\end{eqnarray}
i.e. the operator $\widetilde{R}_P(w, \widetilde{\pi}_{\psi^M})$ intertwines the two representations $ \mathcal{I}_P(\pi_{\psi^M})$ and $\mathcal{I}_P(\pi_{\psi^M}) \circ \theta(N)$. In the untwisted case the situation is of course similar (but without the occurence of $\theta(N)$ in (3.5.6)).

On the other hand, since $\mathcal{I}_P(\pi_{\psi^M})$ is the Langlands quotient corresponding to $\psi^N \in \widetilde{\Psi}(N)$, we can attach as in section 3.2 the intertwining operator on the representation $\mathcal{I}_P(\pi_{\psi^M})$ (with respect to Whittaker normalization):

\begin{eqnarray}
& & \\
& & \widetilde{\mathcal{I}}_P(\pi_{\psi^M},N) =\widetilde{\mathcal{I}}_P(\pi_{\psi^M},\theta(N))   : \mathcal{I}_P(\pi_{\psi^M}) \rightarrow \mathcal{I}_P(\pi_{\psi^M}) \circ \theta(N). \nonumber
\end{eqnarray}

We then have the following:

\begin{proposition} \cite{A12}
\bigskip

\noindent (a) For $w^0 \in W_{\psi^N}(M)$ we have
\begin{eqnarray}
R_P(w^0,\widetilde{\pi}_{\psi^M}) \equiv 1.
\end{eqnarray}

\bigskip

\noindent (b) For $w \in \widetilde{W}_{\psi^N}(M)$ we have
\begin{eqnarray}
\widetilde{R}_P(w, \widetilde{\pi}_{\psi^M}) = \widetilde{\mathcal{I}}_P(\pi_{\psi^M},N).
\end{eqnarray}
\end{proposition}

We can now formulate the local intertwining relation for the twisted group $G=\widetilde{G}_{E/F}(N)$. For this we inflate the induced representation $\mathcal{I}_P(\pi_{\psi^M})$ to $\widetilde{G}^+_{E/F}(N)$. We can then regard $\mathcal{I}_P(\pi_{\psi^M})$ as acting on the Hilbert space
\[
\widetilde{H}^+_P(\pi_{\psi^M}) = \mathcal{H}_P(\pi_{\psi^M}) \oplus \widetilde{\mathcal{H}}_P(\pi_{\psi^M})
\]
where $\widetilde{\mathcal{H}}_P(\pi_{\psi^M})$ is the space of functions supported on the component $\widetilde{G}_{E/F}(N)$. For any $w \in \widetilde{W}_{\psi^N}(M)$ we then obtain the linear transformation
\begin{eqnarray}
R_P(w,\widetilde{\pi}_{\psi^M}) : \mathcal{H}_P(\pi_{\psi^M}) \rightarrow \widetilde{\mathcal{H}}_P(\pi_{\psi^M})
\end{eqnarray}
by setting
\[
(R_P(w,\widetilde{\pi}_{\psi^M}) \phi)(x) = (\widetilde{R}_P(w,\widetilde{\pi}_{\psi^M}) \phi) (x \theta(N)^{-1})
\]
for $\phi \in \mathcal{H}_P(\pi_{\psi^M})$ and $x \in \widetilde{G}_{E/F}(N)(F)$. Note that from the definition we have:
\begin{eqnarray*}
R_P(w,\widetilde{\pi}_{\psi^M}) =  \widetilde{\pi}_{\psi^M} \circ R_P(w,\pi_{\psi^M}).
\end{eqnarray*}
Here
\begin{eqnarray}
& & \\
& & R_P(w,\pi_{\psi^M}) = l(\widetilde{w},\pi_{\psi^M}) \circ ( \lambda(w)^{-1} \epsilon_P(w,\psi^N ) R_{w^{-1}P | P}(\pi_{\psi^M})) \nonumber \\
&=& l(\widetilde{w},\pi_{\psi^M}) \circ ( r_P( w ,\psi^N)^{-1}  J_{w^{-1}P | P}(\pi_{\psi^M})) \nonumber
\end{eqnarray}
with $ l(\widetilde{w},\pi_{\psi^M}) : \mathcal{H}_{w^{-1}P}(\pi_{\psi^M}) \rightarrow \widetilde{\mathcal{H}}_P(\pi_{w \psi^M})$ given by
\begin{eqnarray}
( l(\widetilde{w},\pi_{\psi^M}) \phi)(x) = \phi(\widetilde{w}^{-1} x)
\end{eqnarray}
for $\phi \in \mathcal{H}_{w^{-1}P}(\pi_{\psi^M}) $ and $x \in \widetilde{G}_{E/F}(N)(F)$, and
\begin{eqnarray*}
r_P( w ,\psi^N) = \lambda(w) \epsilon_P(w,\psi^N)^{-1} r_{w^{-1}P|P}(\psi^M ).
\end{eqnarray*}
In other words, the definition of $R_P(w,\pi_{\psi^M})$ is given by exactly the same formula as in (3.3.43).

\bigskip

Given $\widetilde{f} \in \widetilde{\mathcal{H}}(N)$, the integration of $\widetilde{f}$ against the representation $\mathcal{I}_P(\pi_{\psi^M})$ gives the linear transformation:
\begin{eqnarray*}
\mathcal{I}_P(\pi_{\psi^M}, \widetilde{f}) : \widetilde{\mathcal{H}}_P(\pi_{\psi^M}) \rightarrow \mathcal{H}_P(\pi_{\psi^M}).
\end{eqnarray*}

For $u=w$ in the set $\mathfrak{N}_{\psi^N}(G,M) = W_{\psi^N}(G,M)$ (here $G=\widetilde{G}_{E/F}(N)$), define:
\begin{eqnarray}
\widetilde{f}_N(\psi^N,u) := \tr( R_P(w, \widetilde{\pi}_{\psi^M}) \circ \mathcal{I}_P(\pi_{\psi^M}, \widetilde{f}) ).
\end{eqnarray}

\bigskip

\noindent On the other hand, the induced representation $\mathcal{I}_P(\pi_{\psi^M})$ of $G_{E/F}(N)(F)= \widetilde{G}^0_{E/F}(N)(F)$ also extends to the bi-torsor, provided by the operator
\[
\widetilde{\mathcal{I}}_P(\pi_{\psi^M},N) = \widetilde{\mathcal{I}}_P(\pi_{\psi^M},\theta(N)).
\]
Thus as in section 3.1 we have the linear form given by the twisted character
\begin{eqnarray}
\widetilde{f}_N(\psi^N) = \tr \widetilde{\mathcal{I}}_P(\pi_{\psi^M}, \widetilde{f}), \,\ \widetilde{f} \in \widetilde{\mathcal{H}}(N).
\end{eqnarray}

If we define the linear transformation:
\begin{eqnarray}
& & \mathcal{I}_P(\pi_{\psi^M},N) : \mathcal{H}_P(\pi_{\psi^M}) \rightarrow \widetilde{\mathcal{H}}_P(\pi_{\psi^M}) \\
& & (\mathcal{I}_P(\pi_{\psi^M},N) \phi)(x) := (\widetilde{\mathcal{I}}_P(\pi_{\psi^M},N) \phi)(x \theta(N)^{-1})       \nonumber
\end{eqnarray}
then we have
\begin{eqnarray}
\widetilde{f}_N(\psi^N) &=& \tr \widetilde{\mathcal{I}}_P( \pi_{\psi^M} , \widetilde{f} ) \\
&=& \tr(\mathcal{I}_P(\pi_{\psi^M},N) \circ \mathcal{I}_P(\pi_{\psi^M}, \widetilde{f})). \nonumber
\end{eqnarray}

\noindent But by (3.5.9) we have
\begin{eqnarray}
R_P(w, \widetilde{\pi}_{\psi^M}) = \mathcal{I}_P(\pi_{\psi^M} ,N )
\end{eqnarray}
hence we obtain:
\begin{eqnarray}
\widetilde{f}_N(\psi^N,u) = \widetilde{f}_N(\psi^N), \,\ \widetilde{f} \in \widetilde{\mathcal{H}}(N).
\end{eqnarray}

\noindent In particular the linear form $\widetilde{f}_N(\psi^N,u)$ is independent of $u$, which is to be expected since the twisted centralizer $\overline{\widetilde{S}}_{\psi^N}$ is connected.

The endoscopic counterpart of the distribution $\widetilde{f}_N(\psi^N,u)$ can be defined as in the previous subsection. Thus let $s \in \overline{\widetilde{S}}_{\psi^N}$, and let $(G^{\prime},\psi^{\prime})$ be the pair corresponding to $(\psi^N,s)$ under (twisted version of) the correspondence described in section 3.2. Here $G^{\prime} \in \widetilde{\mathcal{E}}(N)$ and $\psi^{\prime} \in \Psi(G^{\prime})$. Assuming the validity of part (a) of theorem 3.2.1 for the pair $(G^{\prime},\psi^{\prime})$ (if $G^{\prime} \notin \widetilde{\mathcal{E}}_{\ellip}(N)$ then the assertion can easily be reduced to the case of a Levi subgroup of $G^{\prime}$), in particular the existence of the stable linear form $f^{G^{\prime}}(\psi^{\prime})$ for $f \in \mathcal{H}(G^{\prime})$. We can define the linear form:
\begin{eqnarray}
\widetilde{f}^{\prime}_N(\psi^N,s)=\widetilde{f}^{G^{\prime}}_N(\psi^N,s) := \widetilde{f}^{G^{\prime}}(\psi^{\prime}).
\end{eqnarray}

\noindent On the other hand, part (a) of theorem 3.2.1 asserts that 
\[
\widetilde{f}^{G^{\prime}}(\psi^{\prime}) = \widetilde{f}_N(\psi^N)
\]
hence we have:
\begin{eqnarray}
\widetilde{f}^{\prime}_N(\psi^N,s) = \widetilde{f}^{G^{\prime}}(\psi^{\prime}) = \widetilde{f}_N(\psi^N).
\end{eqnarray}

\noindent In particular the linear form $\widetilde{f}^{\prime}_N(\psi^N,s)$ is independent of $s$. Furthermore, combining (3.5.18) and (3.5.20), we have:
\begin{eqnarray}
\widetilde{f}^{\prime}_N(\psi^N,  s) = \widetilde{f}_N(\psi^N,u).
\end{eqnarray}

\noindent Now if we replace $s$ by $s_{\psi^N} s$, then (3.5.21) also holds with $s$ replaced by $s_{\psi^N} s$. We record this as:

\begin{corollary} (of proposition 3.5.1)
Assume that part (a) of theorem 3.2.1 holds for any pair
\[
(G,\psi), \,\ G \in \widetilde{\mathcal{E}}(N), \psi \in \Psi(G).
\]
Then for any $u \in \mathfrak{N}_{\psi^N}(\widetilde{G}_{E/F}(N),M)$ and $s \in \overline{\widetilde{S}}_{\psi^N}$ as above, we have
\begin{eqnarray}
\widetilde{f}_N(\psi^N,u) = \widetilde{f}^{\prime}_N(\psi^N, s_{\psi^N} s) , \,\ \widetilde{f} \in \widetilde{\mathcal{H}}(N).
\end{eqnarray}
\end{corollary}

\bigskip

Corollary 3.5.2 is thus the local intertwining relation for the twisted group $\widetilde{G}_{E/F}(N)$. Together with theorem 3.4.3 (i.e. the local intertwining relation for $G \in \widetilde{\mathcal{E}}_{\simp}(N)$) they form an important component in the comparison of trace formulas in section 5 and 6.

\bigskip

To conclude this subsection, we record some results on intertwining operators that can be deduced from the theory of Whittaker models; more precisely, from the works of Shahidi \cite{S,S1,S2}. Again we refer to section 2.5 of \cite{A1} for more detailed discussions.

For this discussion $G$ is a general connected quasi-split group over $F$. For our purpose we only need to apply the results in the case where $G=U_{E/F}(N)$ or $G_{E/F}(N)$. We let $(B,\chi)$ be a Whittaker datum for $G$. For $M$ a standard Levi subgroup of $G$, we denote by $(B_M,\chi_M)$ the corresponding Whittaker datum for $M$. 

We consider the particular case where $\pi_M$ is an irreducible tempered representation of $M(F)$ that is {\it generic}, i.e. possess a non-zero Whittaker functional. Let $P$ be the standrad parabolic subgroup of $G$ with Levi component $M$ (as usual standard means $P \supset B$), we again consider the induced representation $\mathcal{I}_P(\pi_M)$ of $G(F)$. As before denote by $W(\pi_M)$ the subgroup of $W(M)$ consisting of those Weyl elements $w$ such that $w \pi_M \cong \pi_M$. In the case of generic representations, the construction of the normalized intertwining operator for $w \in W(M)$:
\begin{eqnarray}
R_P(w,\pi_M) : \mathcal{I}_P(\pi_M) \rightarrow \mathcal{I}_P( w \pi_M)
\end{eqnarray}
is already given in \cite{S,S1,S2}. We have already used Shahidi's results of {\it loc. cit.} in the proof of proposition 3.3.1. Again, in the context of the local classification, since $\pi_M$ is tempered, the $L$-parameter of $M$ classifying of $\pi_M$ is determined by $\pi_M$, and so can be omitted from the notation. 

From the results of {\it loc. cit.} the normalized intertwining operators $R_P(w,\pi_M)$ are unitary, and satisfy the cocycle relation
\begin{eqnarray}
& & \\
& & R_P(w^{\prime} w,\pi_M) = R_P(w^{\prime},w \pi_M) \circ R_P(w,\pi_M), \,\ w^{\prime},w \in W(M). \nonumber
\end{eqnarray}

Thus in our context (i.e. for $G=U_{E/F}(N)$ or $G_{E/F}(N)$) if we are in the case where $\pi_M$ is a generic representation, then proposition 3.3.1 and 3.3.5 already follows from \cite{S,S1,S2}. 

In addition, when $w \in W(\pi_M)$, there is a {\it canonical} choice of intertwining operator $\widetilde{\pi}_M(w) : w \pi_M \rightarrow \pi_M$. More precisely, let $\omega$ be a Whittaker functional for $\pi_M$ with respect to the Whittaker datum $(B_M,\chi_M)$. Then $\widetilde{\pi}_M$ is uniquely determined by the condition
\[
\omega \circ \widetilde{\pi}_M = \omega.
\] 
We can then define the canonical self-intertwining operator
\begin{eqnarray}
& & R_P(w,\widetilde{\pi}_M) : \mathcal{I}_P(\pi_M) \rightarrow \mathcal{I}_P(\pi_M) \\
& &  R_P(w,\widetilde{\pi}_M) = \widetilde{\pi}_M(w) \circ R_P(w,\pi_M). \nonumber
\end{eqnarray}
Note that this is consistent with the construction of the self-intertwining operator (3.4.8) in section 3.4 (in the context where $G=U_{E/F}(N)$).

Furthermore, the induced representation $\mathcal{I}_P(\pi_M)$ inherits the corresponding $(B_M,\chi_M)$-Whittaker functional $\omega$ of $\pi_M$, as follows. Denote by $N=N_P$ the unipotent radical of $P$. The on the space $\mathcal{H}_{P,\infty}(\pi_M)$ of smooth functions in $\mathcal{H}_P(\pi_M)$, we have the Whittaker integral, for $h \in \mathcal{H}_{P}(\pi_M)$ and $\lambda \in \mathfrak{a}^*_{M,\mathbf{C}}$:
\begin{eqnarray}
& & \\
& & W_{\chi,\omega}(x,h,\pi_{M,\lambda}) = \int_{N_*(F)} \omega(h_{\pi_M,\lambda}(w_*^{-1} n_*  x)) \chi(n_*)^{-1} \,\ dn_*, \,\ x \in G(F). \nonumber
\end{eqnarray}   

\noindent Here $h_{\pi_M,\lambda}$ is the function on $G(F)$ with respect to the induced representation $\mathcal{I}_P(\pi_{M,\lambda})$, as in the discussion of section 3.3, i.e.
\[
h(x) = \pi_M(M_P(x)) h(K_P(x)) e^{(\lambda +\rho_P)(H_P(x))}, \,\ x \in G(F)
\] 
and $N_* = N_{P_*}$ is the unipotent radical of the standard parabolic subgroup $P_*=M_* N_*$ that is adjoint to $P=M N$, i.e. we have the condition
\[
M_* = w_* M w_*^{-1}
\] 
where $w_*=w_l w_l^M$, with $w_l$ and $w_l^M$ being the longest elements in the restricted Weyl groups of $G$ and $M$ respectively. 

\noindent The Whittaker integral (3.5.26) converges absolutely for $Re(\lambda)$ in a certain affine chamber, and extends to an entire function of $\lambda$. The linear functional $\Omega_{\chi,\omega}(\pi_M)$ on $\mathcal{I}_P(\pi_M)$ defined by:
\begin{eqnarray}
\Omega_{\chi,\omega}(\pi_M)(h) := W_{\chi,\omega}(1,h,\pi_M), \,\ h \in \mathcal{H}_{P,\infty}(\pi_M)
\end{eqnarray}
is then a (non-zero) $(B,\chi)$-Whittaker functional for $\mathcal{I}_P(\pi_M)$.

We then have the following result from the works of Shahidi \cite{S,S1,S2}:

\begin{proposition} ({\it c.f.} statement of Theorem 2.5.1 and Corollary 2.5.2 of \cite{A1})

\noindent (a) For $w \in W(\pi_M)$ we have
\begin{eqnarray}
\Omega_{\chi,\omega}(\pi_M) \circ R_P(w,\widetilde{\pi}_M) = \Omega_{\chi,\omega}(\pi_M).
\end{eqnarray}

\bigskip

\noindent (b) Let $(\Pi,\mathcal{V})$ be the unique irreducible $(B,\chi)$-generic subrepresentation of $\mathcal{I}_P(\pi_M)$. Then for $w \in W(\pi_M)$ we have:
\begin{eqnarray}
R_P(w, \widetilde{\pi}_M) \phi = \phi, \,\ \phi \in \mathcal{V}_{\infty}.
\end{eqnarray}
\end{proposition}

\noindent Note that part (b) is an immediate consequence of part (a). Indeed the restriction of $R_P(w,\widetilde{\pi}_M)$ to the irreducible subspace $\mathcal{V}$ is a non-zero scalar. Since the Whittaker functional $\Omega_{\chi,\omega}(\pi_M)$ on $\mathcal{H}_{P,\infty}(\pi_M)$ is supported on $\mathcal{V}_{\infty}$, it follows from (3.5.28) that this scalar must be equal to one.

\bigskip

\noindent Proposition 3.5.3 applies to those $\pi_M$ that are irreducible tempered generic representations of $M(F)$. We will apply proposition 3.5.3 to the following cases, where the condition that $\pi_M$ being irreducible tempered representation of $M(F)$ already implies genericity:

\begin{enumerate}
\item $M$ is a product of $G_{E/F}(N_i)$'s.

\item $M=T$ is minimal.

\item $F=\mathbf{C}$ and $G,M$ arbitrary.
\end{enumerate}

\bigskip

This concludes section 2. Starting from the next section,we will be concerned with the global situation. We will come back to the local study in section 7, based on the global results from section 4 to 6.

\section{\textbf{Trace formulas and their stabilization}}

In this section we return to the global setting. Thus $F$ is a global field. In this section we carry out some preliminary comparison of trace formulas. The two trace formulas that we need are the trace formula for the unitary group $U_{E/F}(N)$, and the twisted trace formula for the twisted group $\widetilde{G}_{E/F}(N)$. For the comparison we also need their stabilization.

Section 4 depends only on the discussion of section 3.1, and is independent of the rest of section 3.

\subsection{Discrete part of trace formula}

We begin with the discrete part of the trace formula for unitary group (for the discussion in the more general context see section 3.1 of \cite{A1}). Thus we put $G=U_{E/F}(N)$. Write $\mathcal{H}(G)$ for the global adelic Hecke algebra of $U_{E/F}(N)(\mathbf{A}_F)$, consisting of smooth, compactly supported, complex-valued functions of $G(\mathbf{A}_F)$ that are $K$-finite with respect to left and right action of a maximal compact subgroup $K$ of $G(\mathbf{A}_F)$. 

Fix the minimal Levi subgroup $M_0$ of $G$ to be the standard diagonal one (which is also the maximal torus $T$ of $G$). Denote by $\mathcal{L}$ the set of Levi subgroups of $G$ containing $M_0$, which we call the standard Levi subgroups of $G$ (in the literature the elements of $\mathcal{L}$ are also called semi-standard). 

As in \cite{A1} we let $t \geq 0$ be the parameter that controls the norm of the imaginary part of the archimedean infinitesimal characters of representations. The discrete part of the trace formula for $G$ is a linear form $I_{\disc,t}^G$ on $\mathcal{H}(G)$ given by:
\begin{eqnarray}
& &   \\ 
& & I_{\disc,t}^G(f) = \sum_{\{M\}} |W(M)|^{-1} \sum_{w \in W(M)_{\reg}} |\det(w-1)_{\mathfrak{a}^G_M}|^{-1} \tr ( M_{P,t}(w) \mathcal{I}_{P,t}(f)), \nonumber \\
& & \,\ \,\ \,\ \,\ \,\ \,\ f \in \mathcal{H}(G). \nonumber
\end{eqnarray}

The explanation of these terms is as follows:

\begin{itemize}

\item $\{M\}$ is the set of conjugacy classes of standard Levi subgroups $M \in \mathcal{L}$ under the action of the Weyl group $W^G_0 = W^G(M_0)$ of $G$ with respect to $M_0$. Equivalently the sum can be taken over the set of $G$-conjugacy classes of all Levi subgroup of $G$.

\item $W(M) = \Normal(A_M,G)/M$ is the relative Weyl group of $G$ with respect to $M$. Here $A_M$ is the maximal $F$-split component of the centre of $M$. 

\item $\mathfrak{a}_M$ (resp. $\mathfrak{a}_G$) is the real vector space $\Hom_{\mathbf{Z}}(X^*_F(M),\mathbf{R})$ (resp. $\Hom_{\mathbf{Z}}(X^*_F(G),\mathbf{R})$). The canonical complement of $\mathfrak{a}_G$ in $\mathfrak{a}_M$ is noted as $\mathfrak{a}_M^G$. In our case where $G=U_{E/F}(N)$, we have $\mathfrak{a}_G=0$ and hence $\mathfrak{a}_M^G=\mathfrak{a}_M$, but we will still use the same notation as in the general case.

\item $W(M)_{\reg}$ is the set of regular elements of $W(M)$, i.e. the set of elements $w \in W(M)$ such that $\det(w-1)|_{\mathfrak{a}_M^G} \neq 0$.

\item $P \in \mathcal{P}(M)$ is a parabolic subgroup of $G$ with $M$ as Levi component.

\item $\mathcal{I}_P$ is the representation of $G(\mathbf{A}_F)$ parabolically induced from the relative discrete spectrum 
\[
L^2_{\disc}(M(F) A_{M,\infty}^+ \backslash M(\mathbf{A}_F))
\] 
of $M(\mathbf{A}_F)$. Here we denote by $A_{M,\infty}^+$ the identity component of \linebreak $A_{M_{\mathbf{Q}}}(\mathbf{R})$ (with $M_{\mathbf{Q}} = \Res_{F/\mathbf{Q}}M$ and $A_{M_{\mathbf{Q}}}$ is the maximal $\mathbf{Q}$-split component of its centre). We write $\mathcal{H}_P$ for the underlying space of this representation, realized as the space of left $N_P(\mathbf{A}_F)$-invariant functions $\phi$ on $G(\mathbf{A}_F)$ such that the function $\phi(mk)$ on $M(\mathbf{A}_F) \times K$ belongs to the space
\[
L^2_{\disc}(M(F) A_{M,\infty}^+ \backslash M(\mathbf{A}_F)) \otimes L^2(K).
\]

\bigskip

\noindent As usual $\mathcal{I}_P(f)$ is the corresponding operator associated to $f \in \mathcal{H}(G)$ acting on $\mathcal{H}_P$ given by integration:
\[
\mathcal{I}_P(f)\phi= \int_{G(\mathbf{A}_F)} f(g) \mathcal{I}_P(g)\phi \,\ dg, \,\ \phi \in \mathcal{H}_P
\]
(see section 3.1 of \cite{A1} for a breif discussion on the choice of measure).

\item We have a decomposition
\[
\mathcal{I}_P=\bigoplus_{t \geq 0} \mathcal{I}_{P,t}
\]
where $\mathcal{I}_{P,t}$ is the subrepresentation whose irreducible constituents have the property that the imaginary part of their archimedean infinitesimal character has norm $t$. Denote by $\mathcal{H}_{P,t}$ the underlying space of $\mathcal{I}_{P,t}$.

\item $M_P(w) = l(w) \circ M_{P^{\prime}|P}$. Here $P^{\prime} = w^{-1} P$, and $M_{P^{\prime}|P}: \mathcal{H}_{P,t} \rightarrow \mathcal{H}_{P^{\prime},t}$ is the global intertwining operator, defined to be the analytic continuation at $\lambda =0$ of the following intertwining integral: in general for $P,P^{\prime} \in \mathcal{P}(M)$, and $\lambda \in (\mathfrak{a}_M^G)^*_{\mathbf{C}}$ with $Re(\lambda)$ in a certain affine chamber, it is given by:
\begin{eqnarray*}
& & (M_{P^{\prime}|P} \phi)(x) \\
&:=& \int_{ N_P(\mathbf{A}_F) \cap N_{P^{\prime}}(\mathbf{A}_F)   \backslash N_{P^{\prime}}(\mathbf{A}_F)  } \phi(n x ) e^{(\lambda + \rho_P)(H_P(nx))} \,\ dn \cdot e^{-(\lambda+ \rho_{P^{\prime}})H_{P^{\prime}}(x) },
\end{eqnarray*} 
\bigskip
\noindent and $l(w): \mathcal{H}_{P^{\prime},t} \rightarrow \mathcal{H}_{P,t}$ given by left translation by $\widetilde{w}^{-1}$:
\[
(l(w)\phi)(y) = \phi(\widetilde{w}^{-1} y) \mbox{ for } y \in G(\mathbf{A}_F)
\] 
defined by any representative $\widetilde{w}$ of $w$ in $G(F)$. The subrepresentations $\mathcal{I}_{P,t}$ are stable under $M_P(w)$, and we denote by $M_{P,t}(w): \mathcal{H}_{P,t} \rightarrow \mathcal{H}_{P,t}$ the restriction of $M_P(w)$ to $\mathcal{I}_{P,t}$.

\end{itemize}

We will denote by $R_{\disc}^G$ for the representation of $G(\mathbf{A}_F)$ on the discrete spectrum $L^2_{\disc}(G(F)  \backslash G(\mathbf{A}_F))$, and similarly $R^G_{\disc,t}$ for the subrepresentation of $R^G_{\disc,t}$ whose irreducible constituents have the property that the imaginary part of their archimedean infinitesimal characters has norm $t$ .

\bigskip

Next we turn to the twisted trace formula. Thus with notations as in section 2.4, let $G = \widetilde{G}_{E/F}(N) = G_{E/F}(N) \rtimes \theta$ be our twisted group. We put $G^0 = \widetilde{G}^0_{E/F}(N) = G_{E/F}(N)$. We denote by $\mathcal{H}(G) = \widetilde{\mathcal{H}}(N)$ the global Hecke module consisting of smooth, compactly supported functions on $G(\mathbf{A}_F)$, which are $K$-finite with respect to the left and right action of a maximal compact subgroup of $G^0(\mathbf{A}_F)$. 

Given $y_1 = x_1 \rtimes \theta, y_2=x_2 \rtimes \theta \in G = G^0 \rtimes \theta$, we have the usual convention for the interpretation of the algebraic operation:
\begin{eqnarray}
& & y_1^{-1} \cdot y_2 = \theta^{-1}(x_1^{-1} \cdot x_2) \in G^0 \\
& & y_1 \cdot y_2^{-1} = x_1 \cdot x_2^{-1} \in G^0. \nonumber
\end{eqnarray}

Let $M_0$ be the diagonal minimal Levi subgroup of $G^0$ (which is the maxiaml diagonal torus of $G^0 = G_{E/F}(N)$). Again denote by $\mathcal{L}$ the set of Levi subgroups of $G^0$ containing $M_0$. The discrete part of the twisted trace formula for $G$ is a linear form on $\mathcal{H}(G)$ which bears a formal resemblance to the untwisted case:
\begin{eqnarray}
& &   \\ 
& & I_{\disc,t}^G(f) = \sum_{\{M\}} |W(M)|^{-1} \sum_{w \in W(M)_{\reg}} |\det(w-1)_{\mathfrak{a}^G_M}|^{-1} \tr ( M_{P,t}(w) \mathcal{I}_{P,t}(f)), \nonumber \\
& & \,\ \,\ \,\ \,\ \,\ \,\ f \in \mathcal{H}(G). \nonumber
\end{eqnarray}

The explanation of these terms is as follows:
\begin{itemize}

\item $\{M\}$ is the set of conjugacy classes of standard Levi subgroups $M \in \mathcal{L}$ under the action of the Weyl group $W_0^{G^0} = W^{G^0}(M_0)$ of $G^0$ with respect to $M_0$. Equivalently the sum can be taken over the set of $G^0$-conjugacy classes of all Levi subgroup of $G^0$.

\item $W(M) = \Normal(A_M,G)/M$ is the relative Weyl set of $G$ with respect to $M$.

\item $\mathfrak{a}_M$ is as before, while $\mathfrak{a}_G :=\mathfrak{a}_{G^0}^{\theta}$ is the subspace of $\theta$-invariants of $\mathfrak{a}_{G^0}$, which also occur as the space of $\theta$-coinvariants of $\mathfrak{a}_{G^0}$. Then $\mathfrak{a}_M^G$ is the kernel of the map:
\[
\mathfrak{a}_M \rightarrow \mathfrak{a}_{G^0} \rightarrow \mathfrak{a}_G
\]
In the case $G=\widetilde{G}_{E/F}(N)$ we have $\dim \mathfrak{a}_{G^0}=1$ and $\mathfrak{a}_G=0$, so $\mathfrak{a}_M^G=\mathfrak{a}_M$. But we will still use the same notation in the general case.

\item $W(M)_{\reg}$ is the set of regular elements, i.e. consisting of elements $w \in W(M)$ such that $\det(w-1)|_{\mathfrak{a}_M^G} \neq 0$.

\item $\mathcal{I}_{P,t} : \mathcal{H}_{P,t} \rightarrow \mathcal{H}_{P,t}^0$ is the map defined as follows. First $\mathcal{H}_{P,t}^0$ is the underlying space of the induced representation defined for $G^0$ as above, and $\mathcal{H}_{P,t}$ is the space of functions $\phi$ on $G(\mathbf{A}_F)$ such that, for every $y \in G(\mathbf{A}_F)$, the function
\[
x \mapsto \phi(xy) \mbox{ is in } \mathcal{H}^0_{P,t}. 
\]
Then for every $y \in G(\mathbf{A}_F)$, the linear map $\mathcal{I}_{P,t}(y): \mathcal{H}_{P,t} \rightarrow \mathcal{H}_{P,t}^0$ is defined just by right translation by $y$. The operator $\mathcal{I}_{P,t}(f):\mathcal{H}_{P,t} \rightarrow \mathcal{H}_{P,t}^0$ for $f \in \mathcal{H}(G)$ is then defined by integration as in the untwisted case.

\item $M_{P,t}(w) : \mathcal{H}_{P,t}^0 \rightarrow \mathcal{H}_{P,t}$ is defined by $M_{P,t} = l(w) \circ M_{P^{\prime}|P}$, with $P^{\prime} = w^{-1}P$, and $M_{P^{\prime}|P}: \mathcal{H}_{P,t}^0 \rightarrow \mathcal{H}^0_{P^{\prime},t}$ is {\it the same} global intertwining operator defined as above. While $l(w): \mathcal{H}^0_{P^{\prime},t} \rightarrow \mathcal{H}_{P,t}$ is defined analogous to the untwisted case by:
\begin{eqnarray}
(l(w) \phi)(y) = \phi(\widetilde{w}^{-1} y), \,\ y \in G(\mathbf{A}_F), \phi \in \mathcal{H}^0_{P^{\prime},t}
\end{eqnarray}
for any representative $\widetilde{w} \in G(F)$ of $w$, with the operation $\widetilde{w}^{-1} y$ understood as in (4.1.2).
\end{itemize}

\bigskip

For our purpose we will only need the twisted trace formula for the group $G=\widetilde{G}_{E/F}(N)$, in which case we note $I_{\disc,t}^G$ as $\widetilde{I}^N_{\disc,t}$.

The invariant trace formula, in both its standard and twisted version, is established in \cite{A2,A3,LW}. The proof of the main theorems of this paper, in particular the theorems stated in section 2.5, is based on the comparison of the two trace formulas via their stabilized version.

\subsection{Stabilization of trace formula}

The stabilization of trace formulas is crucial for their comparison. In the case of connected groups, the stabilization is established by Arthur \cite{A4,A5,A6}, combined with the work of Waldspruger \cite{W1,W2,W5}, Ngo \cite{N} , and Chaudouard-Laumon \cite{CL1,CL2} on the fundammental lemma (both standard and weighted version) and the Langlands-Shelstad transfer conjecture (in the case of unitary groups, which is our situations here, the standard fundamental lemma was established earlier in \cite{LN}). We also need the stable trace formula for twisted groups, which is currently being carried out by Waldspurger (the necessary fundamental lemmas and transfer conjecture in the twisted case were known from the works of Waldspurger \cite{W3,W4} and Ngo \cite{N}). We refer to section 3.2 of \cite{A1} for a more detailed discussion. Here we just recall the formalism.

Thus we denote by $G$ either the unitary group $U_{E/F}(N)$, or the twisted group $\widetilde{G}_{E/F}(N)$. In the previous subsection we defined the linear form $I_{\disc,t}^G$ on $\mathcal{H}(G)$ (for $t \geq 0$). The stabilization of $I_{\disc,t}^G$ refers to a decomposition:

\begin{eqnarray}
I^G_{\disc,t} (f) = \sum_{G^{\prime} \in \mathcal{E}_{\ellip}(G)} \iota(G,G^{\prime}) \cdot  \widehat{S}_{\disc,t}^{G^{\prime}}(f^{G^{\prime}}).
\end{eqnarray} 

Here in the summand we write $G^{\prime} \in \mathcal{E}_{\ellip}(G)$, but we caution that in the summand $G^{\prime}$ always denote an (equivalence class of) endoscopic data, not just the endoscopic group itself. As in the local case in section 3, we denote by $f^{G^{\prime}}$ the Langlands-Kottwitz-Shelstad transfer of $f$ to $G^{\prime}$ with respect to the $L$-embedding $\leftexp{L}{G^{\prime}} \rightarrow \leftexp{L}{G}$ that is part of the endoscopic data, with the global transfer factor $\Delta(\delta,\gamma)$ being normalized by the fixed standard splitting of $G$. Thus with the global analogue of the notations in section 3:

\begin{eqnarray}
f^{G^{\prime}}(\delta) = \sum_{\gamma} \Delta(\delta,\gamma) f_G(\gamma)
\end{eqnarray}
here $\delta$ is a strongly $G$-regular stable conjugacy class in $G^{\prime}(F)$, and $\gamma$ runs over the strongly regular $G^0(F)$-conjugacy classes in $G(F)$. The global transfer factor is the one defined by the fixed global standard splitting $S$ of $G$ (in the twisted case, the splitting has to be $\theta$-stable, which is the case for the standard splitting): $\Delta=\Delta_S$. From section 7.3 of \cite{KS}, the product formula is valid:
\[
\Delta_S(\delta,\gamma) = \prod_v \Delta_{S,v}(\delta_v,\gamma_v)
\] 
where $\Delta_{S,v}$ is the local transfer factor defined with repsect to the localization of the global splitting $S$. On the other hand, the local transfer factor $\Delta_v$ of section 3 is defined with respect to choice of a Whittaker normalization. The two are related by:
\begin{eqnarray*}
\Delta_v &=& \frac{\epsilon(1/2, \tau_{G,v} ,\psi_{F_v} )}{\epsilon(1/2, \tau_{G^{\prime},v}   , \psi_{F_v})} \Delta_{S,v} \\
&=&  \epsilon(1/2,  r_v , \psi_{F_v}) \Delta_{S,v}.
\end{eqnarray*}
Here $\tau_G$ (resp. $\tau_{G^{\prime}}$) denotes the Artin representation of $\Gamma_F$ on the $\mathbf{C}$-vector space $X^*(T_G)^{\theta} \otimes_{\mathbf{Z}} \mathbf{C}$ (resp. $X^*(T_{G^{\prime}}) \otimes_{\mathbf{Z}} \mathbf{C}$), where $T_G$, $T_{G^{\prime}}$ are maximal torus of $G$, $G^{\prime}$, and $\tau_{G,v}$, $\tau_{G^{\prime},v}$ are the restrictions of the Artin representation $\tau_G,\tau_{G^{\prime}}$ to the decomposition group $\Gamma_{F_v}$. We put $r : = \tau_G - \tau_{G^{\prime}}$ the corresponding virtual representation of $\Gamma_F$, with $r_v$ its restriction to $\Gamma_{F_v}$. By section 5.3 of \cite{KS}, the virtual representation $r$ is orthogonal, hence by \cite{FQ} the global epsilon factor 
\[
\epsilon(1/2,r)= \prod_{v} \epsilon(1/2,r_v,\psi_{F_v})  
\] 
is equal to one. Thus the product formula:
\[
\Delta(\delta,\gamma) = \Delta_S(\delta,\gamma) = \prod_v \Delta_v(\delta_v,\gamma_v)
\] 
is valid. 

It follows that if $f = \otimes_v^{\prime} f_v \in \mathcal{H}(G) =\otimes_v^{\prime} \mathcal{H}(G_v)$, then  
\[
f^{G^{\prime}} = \otimes_v^{\prime} f^{G^{\prime}}_v \in \mathcal{S}(G^{\prime}) = \otimes_v^{\prime} \mathcal{S}(G^{\prime}_v).
\]

The term $\widehat{S}_{\disc,t}^{G^{\prime}}$ is a linear form on $\mathcal{S}(G^{\prime})$, which depends only on $G^{\prime}$ as an endoscopic group, and not on the endoscopic datum associated to $G^{\prime}$. On the other hand, the transfer map $f \rightarrow f^{G^{\prime}}$ certainly depends on $G^{\prime}$ as an endoscopic datum of $G$. Thus the linear form $f \rightarrow \widehat{S}_{\disc,t}^{G^{\prime}}(f^{G^{\prime}})$ depends on $G^{\prime}$ as an endoscopic datum.

Analogous to the local situation of section 3.1, we have the map:
\begin{eqnarray}
& & \mathcal{H}(G^{\prime}) \rightarrow \mathcal{S}(G^{\prime}) \\
& & f \rightarrow f^{G^{\prime}}, \,\ f \in \mathcal{H}(G^{\prime}) \nonumber
\end{eqnarray}  
given by stable orbital integrals (notation as in the local situation of section 3.1):
\[
f^{G^{\prime}}(\delta) = \sum_{\gamma \rightarrow \delta} f_{G^{\prime}}(\gamma) 
\]
and we denote by $S_{\disc,t}^{G^{\prime}}$ the pull-back of the linear form $\widehat{S}_{\disc,t}^{G^{\prime}}$ to $\mathcal{H}(G^{\prime})$ via (4.2.3).

Finally the coefficients $\iota(G,G^{\prime})$ are given by, referring to the notations of  section 3.2 of \cite{A1}:
\begin{eqnarray}
& & \\
& & \iota(G,G^{\prime}) = | \pi_0(\kappa_G)|^{-1} k(G,G^{\prime}) | \overline{Z}( \widehat{G}^{\prime} )^{\Gamma_F}    |^{-1} | \Out_G(G^{\prime}) |^{-1}. \nonumber
\end{eqnarray}
Here 
\begin{eqnarray}
\overline{Z}(\widehat{G}^{\prime})^{\Gamma_F} = Z(\widehat{G}^{\prime} )^{\Gamma_F} / (Z(\widehat{G}^{\prime})^{\Gamma_F} \cap Z(\widehat{G})^{\Gamma_F}).
\end{eqnarray}
Recall that $Z(\widehat{G}) := Z(\widehat{G}^0)^{\theta}$. While $\pi_0(\kappa_G)$ is the group of connected components of the group
\[
\kappa_G := Z(\widehat{G})^{\Gamma_F} \cap ( Z(\widehat{G}^0)^{\Gamma_F})^0.
\]
In our case the group $\kappa_G$ is already finite, with
\[
|\kappa_{U_{E/F}(N)}|=1, \,\ |\kappa_{\widetilde{G}_{E/F}(N)}|=2.
\]

The coefficient $k(G,G^{\prime})$ is given by:
\[
k(G,G^{\prime}) = |\ker^1(F,Z(\widehat{G}^0))|^{-1} |\ker^1(F,Z(\widehat{G}^{\prime})) |.
\]

For both $G= U_{E/F}(N)$ or $\widetilde{G}_{E/F}(N)$, the coefficient $k(G,G^{\prime})$ is equal to one (lemma 3.5.1 of \cite{R}). 

Hence for $G=U_{E/F}(N)$, we have:
\begin{eqnarray}
\iota(G,G^{\prime}) = | \overline{Z}(\widehat{G}^{\prime})^{\Gamma_F} |^{-1} |\Out_G(G^{\prime})|^{-1}.
\end{eqnarray}

For $G=\widetilde{G}_{E/F}(N)$, we have (now replacing $G$ by $\widetilde{G}_{E/F}(N)$ and $G^{\prime}$ by $G$):
\begin{eqnarray}
\widetilde{\iota}(N,G) := \iota(\widetilde{G}_{E/F}(N),G)= \frac{1}{2}  |\overline{Z}(\widehat{G})^{\Gamma_F}|^{-1} |\widetilde{\Out}_N(G)|^{-1}.
\end{eqnarray}

For instance one has:
\begin{eqnarray*}
 \iota(U(N), U(N_1) \times U(N_2))  = \left \{ \begin{array}{c}  1 \mbox{ if } N_1 =0 \mbox{ or } N_2 =0 \\   1/2 \mbox{ if } N_1 ,N_2 \neq 0, N_1 \neq N_2 \\ 1/4 \mbox{ if } N_1=N_2 \neq 0.
\end{array} \right. 
\end{eqnarray*}

Similarly, consider $G \in \widetilde{\mathcal{E}}_{\ellip}(N)$. Then we have:
\begin{eqnarray*}
 \widetilde{\iota}(N, G)  = \left \{ \begin{array}{c}  1/2 \mbox{ if } G \in \widetilde{\mathcal{E}}_{\simp}(N)  \\   1/4 \mbox{ otherwise. } 
\end{array} \right. 
\end{eqnarray*}

\bigskip

As we have recalled in the beginning of this subsection, the stabilization of the discrete part of the trace formula $I_{\disc}^G$ is now known when $G$ is a connected group. The case where $G$ is a twisted group is the work in progress of Waldspurger, {\it c.f.} \cite{W7,W8}. For our purpose we only need the case of the group $\widetilde{G}_{E/F}(N)$, which we formulate as hypothesis, and which we assume for the rest of the paper:

\begin{hypothesis}
\end{hypothesis}
The stabilization (4.2.1) of the discrete part of the trace formula $I_{\disc}^G$ is valid for the twisted group $\widetilde{G}_{E/F}(N)$.

\subsection{Preliminary comparison}

In this subsection we carry out a preliminary comparison of the trace formulas for $U_{E/F}(N)$ and $\widetilde{G}_{E/F}(N)$, and in particular derive the existence of ``weak base change". This forms the background of the more elaborate comparison in the next two sections.

Thus as in previous subsections, we denote by $G$ for either the group $U_{E/F}(N)$, or the twisted group $\widetilde{G}_{E/F}(N)$. We will denote by $S$ any finite set of primes containing all the archimedean primes, and such that $G$ (as a group over $F$) is unramified outside $S$. Recall that the global Hecke algebra (or module in the twisted case) $\mathcal{H}(G)$ is defined with respect to a maximal compact subgroup $K=\prod_v K_v$ of $G^0(\mathbf{A}_F)$. Put $K^S = \prod_{v \notin S} K_v$. Then $K_v$ is a hyperspecial maximal compact subgroup of $G(F_v)$ for $v \notin S$. Define $\mathcal{H}(G,K^S) \subset \mathcal{H}(G)$ to be the subalgebra (or submodule in the twisted case) consisting of functions that are left and right bi-invariant with respect to the action of $K^S \subset G^0(\mathbf{A}_F^S)$ on $G(\mathbf{A}_F^S)$. Also define the unramified Hecke algebra at primes outside $S$:
\[
\mathcal{H}^S_{\un} = \mathcal{H}^S_{\un}(G^0) = C^{\infty}_c( K^S\backslash  G^0(\mathbf{A}_F^S)  /K^S ).
\]

For $h \in \mathcal{H}^S_{\un}$, there is an action on $\mathcal{H}(G,K^S)$:
\[
f \rightarrow f_h, \,\ f \in \mathcal{H}(G,K^S)
\]
defined in terms of spectral multipliers. Thus if $\pi^0$ is an irreducible admissible unitary representation of $G^0(\mathbf{A}_F)$ that is unramified outside $S$, and has an extension $\pi$ to $G(\mathbf{A}_F)$, then 
\[
\tr \pi(f_h) = \widehat{h}(c^S(\pi)) \tr \pi(f)
\]
here $\widehat{h}$ is the Satake transform of $h$. We have $c^S(\pi) \in \mathcal{C}_{\mathbf{A}_F}^S(G)$, where in the case $G = G^0 \rtimes \theta$ is a twisted group, $\mathcal{C}_{\mathbf{A}_F}^S(G)$ is defined to be the set of families $c^S \in \mathcal{C}^S(G^0)$ such that
\[
\widehat{\theta}(c_v) = c_v \mbox{ for } v \notin S.
\]
In our case where $G=\widetilde{G}_{E/F}(N)$, this is the same as the definition of $\widetilde{\mathcal{C}}_{\mathbf{A}_F}(N)$ in section 2.3.

For given $f \in \mathcal{H}(G)$, let $S$ be such that $f \in \mathcal{H}(G,K^S)$. We then have a decomposition:
\begin{eqnarray}
I_{\disc,t}^G(f) = \sum_{c^S} I^G_{\disc,c^S,t}(f)
\end{eqnarray}
here $c^S$ runs over the subset $\mathcal{C}^S(G) \subset \mathcal{C}^S_{\aut}(G)$ consisting of those families $c^S(\pi)$ for $\pi$ an automorphic representation of $G(\mathbf{A}_F)$, and $I^G_{\disc,c^S,t}$ is the $c^S$ eigen-component of $I^G_{\disc,t}$ with respect to the action $f \rightarrow f_h$, i.e. we have
\begin{eqnarray}
I^G_{\disc,c^S,t}(f_h) =\widehat{h}(c^S) I^G_{\disc,c^S,t}(f), \,\ h \in \mathcal{H}^S_{\un}.
\end{eqnarray}

In the sum (4.3.1) the summand $I^G_{\disc,c^S,t}$ vanishes for all $c^S$ outside a finite subet of $\mathcal{C}_{\mathbf{A}_F}^S(G)$ that depends on $f$ only through a choice of its Hecke type ({\it c.f.} section 3.1 of \cite{A1} for a discussion of Hecke type).

We put
\begin{eqnarray*}
& & \mathcal{C}_{\mathbf{A}_F}(G) = \varinjlim_S \mathcal{C}_{\mathbf{A}_F}^S(G) \\
& & \mathcal{C}(G) = \varinjlim_S \mathcal{C}^S(G).
\end{eqnarray*}
For $c \in \mathcal{C}(G)$, we then define:
\begin{eqnarray}
I_{\disc,c,t}^G(f) = \sum_{c^S \rightarrow c} I^G_{\disc,c^S,t}(f).
\end{eqnarray}
If $c \in \mathcal{C}_{\mathbf{A}_F}(G) \smallsetminus \mathcal{C}(G)$, we simply put $I^G_{\disc,c,t}(f) =0$; this definition is consistent with (4.3.3). 

Suppose first that $G$ is a connected group. In the definition (4.1.1) of the spectral terms of $I_{\disc}^G(f)$ involving the Levi subgroup $M$, write, for $c \in \mathcal{C}_{\mathbf{A}_F}(G)$: 
\[
\mathcal{I}_{P,c,t}(f) = \bigoplus_{c(\pi)=c } \mathcal{I}_{P,t,\pi}(f) 
\]
here $\pi$ runs over representations of $G(\mathbf{A}_F)$ of the form $\pi = \mathcal{I}_P(\pi_M)$, with irreducible $\pi_M \subset L^2_{\disc}(M(F) A^+_{M,\infty}\backslash M(\mathbf{A}_F))$, such that the norm of the imaginary part of the archimedean infinitesimal characters of the irreducible constituents of $\pi$ is equal to $t$, and $\mathcal{I}_{P,t,\pi}(f)$ is the restriction of $\mathcal{I}_{P,t}(f)$ to $\pi$. We denote by $\mathcal{H}_{P,c,t}$ the underlying space of $\mathcal{I}_{P,c,t}$, and $M_{P,c,t}(w)$ is the restriction of $M_{P,t}(w)$ to the subspace $\mathcal{H}_{P,c,t}$. We then have for $c \in \mathcal{C}_{\mathbf{A}_F}(G)$:
\begin{eqnarray}
& & \\
& & I^G_{\disc,c,t}(f) = \sum_{\{M\}} |W(M)|^{-1} \sum_{w \in W(M)_{\reg}} |\det(w-1)_{\mathfrak{a}_M^G}|^{-1} \tr(M_{P,c,t}(w) \mathcal{I}_{P,c,t}(f) )\nonumber 
\end{eqnarray} 

Furthermore the expansion (4.3.4) is also valid in the case where $G$ is a twisted group by making the appropriate interpretation of the terms $\mathcal{I}_{P,c,t}(f), M_{P,c,t}(w)$. We need the twisted version only for the group $G=\widetilde{G}_{E/F}(N)$. In this case we write the distribution $I_{\disc,c,t}^G$ as $\widetilde{I}^N_{\disc,c,t}$ for $c \in \widetilde{\mathcal{C}}_{\mathbf{A}_F}(N)$. 

\bigskip

\begin{example}
\end{example}

Recall that for $G= \widetilde{G}_{E/F}(N)$, the theorem of Moeglin-Waldspurger \cite{MW} and Jacquet-Shalika \cite{JS}  implies that we have a bijection:
\begin{eqnarray}
\widetilde{\Psi}(N)  & \stackrel{\simeq}{\rightarrow} & \widetilde{\mathcal{C}}(N) \\
\psi^N & \mapsto & c(\psi^N). \nonumber
\end{eqnarray}
Thus if $c = c(\psi^N)$ for $\psi^N \in \widetilde{\Psi}(N)$, with $t=t(\psi^N)$ the norm of the imaginary part of its archimedean infinitesimal character, we write $\widetilde{I}^N_{\disc,c(\psi^N),t(\psi^N)}$ as $\widetilde{I}^N_{\disc,\psi^N}$. Then we have the decomposition:
\begin{eqnarray}
\widetilde{I}^N_{\disc,t} (\widetilde{f}) = \sum_{\substack{ \psi^N \in \widetilde{\Psi}(N), \\ t(\psi^N)=t }} \widetilde{I}^N_{\disc,\psi^N}(\widetilde{f})     , \,\ \widetilde{f} \in \widetilde{\mathcal{H}}(N).
\end{eqnarray}

\bigskip

Next we need a similar decomposition in the endoscopic expansion (4.2.1).

\begin{lemma} (Lemma 3.3.1 of \cite{A1})

(a) Suppose that $G$ equals $G^0$ and is quasi-split. Then there is a decomposition
\begin{eqnarray}
S_{\disc,t}^G(f) = \sum_{c \in \mathcal{C}_{\mathbf{A}_F }(G)} S_{\disc,c,t}^G(f), \,\ f \in \mathcal{H}(G)
\end{eqnarray}
for stable linear forms $S^G_{\disc,c,t}$ that satisfies the analogue of (4.3.2) and (4.3.3). Furthermore the summand $S_{\disc,c,t}^G(f)$ vanishes for all $c$ outside a finite subset of $\mathcal{C}_{\mathbf{A}_F}(G)$ that depends on $f$ only through a choice of Hecke type.

(b) Suppose that $G$ is either connected or is the twisted group $\widetilde{G}_{E/F}(N)$. Then for any $c \in \mathcal{C}_{\mathbf{A}_F}(G)$, there is a decomposition
\begin{eqnarray}
& & \\
& & I_{\disc,c,t}(f) = \sum_{G^{\prime} \in \mathcal{E}_{\ellip}(G)} \iota(G,G^{\prime}) \sum_{c^{\prime} \rightarrow c} \widehat{S}^{G^{\prime}}_{\disc,c^{\prime},t} (f^{G^{\prime}}), \,\ f  \in \mathcal{H}(G). \nonumber
\end{eqnarray}
Here the sum $\sum_{c^{\prime} \rightarrow c}$ runs over the elements $c^{\prime} \in \mathcal{C}_{\mathbf{A}_F}(G^{\prime})$ such that $c^{\prime}$ maps to $c$ under the $L$-embedding $\leftexp{L}{G^{\prime}} \rightarrow \leftexp{L}{G}$ that is part of the endoscopic datum $G^{\prime} \in \mathcal{E}_{\ellip}(G)$.
\end{lemma}

\bigskip

\begin{rem}
\end{rem}
In the case where $G$ is a connected quasi-split group the assertion is proved in \cite{A1}, lemma 3.3.1. For the twisted group $\widetilde{G}_{E/F}(N)$, the proof in {\it loc. cit.} carry through, under our hypothesis 4.2.1 that the stabilization of $I_{\disc,t}^G = \widetilde{I}_{\disc,t}^N$ is valid. 

\bigskip

In accordance with \cite{A1}, in the context of Lemma 4.3.2, for $c \in \mathcal{C}_{\mathbf{A}_F}(G)$, we write: 
\begin{eqnarray}
S^{G^{\prime}}_{\disc,c,t} (f^{\prime}) := \sum_{c^{\prime} \rightarrow c} S_{\disc,c^{\prime},t}^{G^{\prime}}(f^{\prime}), \,\ f^{\prime} \in \mathcal{H}(G^{\prime})
\end{eqnarray}
with the sum $\sum_{c^{\prime} \rightarrow c}$ is understood as above. Then we can write (4.3.8) as:
\begin{eqnarray}
I_{\disc,c,t}^G(f) = \sum_{G^{\prime} \in \mathcal{E}_{\ellip}(G)} \iota(G,G^{\prime})  \widehat{S}^{G^{\prime}}_{\disc,c,t} (f^{G^{\prime}}).
\end{eqnarray}

\bigskip

Again we caution that in the sum $\sum_{G^{\prime} \in \mathcal{E}_{\ellip}(G)}$, $G^{\prime}$ really denote an endoscopic datum of $G$, instead of just an endoscopic group of $G$. Indeed for $c \in \mathcal{C}_{\mathbf{A}_F}(G)$, the stable linear form $S_{\disc,c,t}^{G^{\prime}}$ depends on $G^{\prime}$ {\it as an endoscopic datum of $G$}, even thoug the stable linear form $S_{\disc,t}^{G^{\prime}}$ depends only on $G^{\prime}$ as an endoscopic group.

In the case where $G=U_{E/F}(N)$, this distinction is essentially harmless, because as we have seen as in section 2.4 the equivalence class of an endoscopic data in $\mathcal{E}_{\ellip}(U_{E/F}(N))$ is uniquely determined by the endoscopic group. On the other hand, for $G= \widetilde{G}_{E/F}(N)$, this distinction is important, as for instance we have two non-equivalent simple twisted endoscopic data $(U_{E/F}(N),\xi_{\chi_+})$ and $(U_{E/F}(N),\xi_{\chi_-})$ of $\widetilde{\mathcal{E}}_{\simp}(N)$ with the same twisted endoscopic group $U_{E/F}(N)$ (here as in section 2.4 he characters $\chi_+ \in \mathcal{Z}_E^+$ and $\chi_- \in \mathcal{Z}_E^-$ are the ones used to define the $L$-embeddings).

Changing the context of the notations, with $G$ replaced by $\widetilde{G}_{E/F}(N)$ and $G^{\prime}$ replaced by $G$. If $(G, \xi)$ is an elliptic twisted endoscopic datum of $\widetilde{G}_{E/F}(N)$, we define, parallel to (4.3.9), for $c^N \in \widetilde{\mathcal{C}}_{\mathbf{A}_F}(N)$ and $t \geq 0$:
\begin{eqnarray}
S_{\disc,c^N,t,\xi}^{G}(f) := \sum_{\xi (c) =c^N  } S_{\disc, c,t}^G(f), \,\ f \in \mathcal{H}(G)
\end{eqnarray}
with the sum runs over $c \in \mathcal{C}_{\mathbf{A}_F}(G)$ such that $\xi(c) = c^N$. In the case where $c^N =c(\psi^N)$ for $\psi^N \in \widetilde{\Psi}(N)$ and $t=t(\psi^N)$ we also write $S_{\disc,c^N, t,\xi}^{G}$ as $S_{\disc,\psi^N,\xi}^{G}$. Then parallel to (4.3.10) we have for $c^N \in \widetilde{\mathcal{C}}_{\mathbf{A}_F}(N)$ and $t \geq 0$:
\begin{eqnarray}
& & \\
& & \widetilde{I}^N_{\disc,c^N,t}(\widetilde{f}) =  \sum_{(G,\xi) \in \widetilde{\mathcal{E}}_{\ellip} (N)} \widetilde{\iota}(N,G) \widehat{S}_{\disc,c^N,t,\xi}^{G}(\widetilde{f}^{(G,\xi)}), \,\ \widetilde{f} \in \widetilde{\mathcal{H}}(N) \nonumber
\end{eqnarray}
(as in section 3.1 we have written the Kottwitz-Shelstad transfer $\widetilde{f}^{(G,\xi)}$ instead of $\widetilde{f}^G$ to emphasize the dependence on the endoscopic datum), and similarly, if $c^N=c(\psi^N)$ for $\psi^N \in \widetilde{\Psi}(N)$, we have
 \begin{eqnarray}
& & \\
& & \widetilde{I}^N_{\disc,\psi^N}(\widetilde{f}) =  \sum_{(G,\xi) \in \widetilde{\mathcal{E}}_{\ellip} (N)} \widetilde{\iota}(N,G) S_{\disc,\psi^N,\xi}^{G}(\widetilde{f}^{(G,\xi)} ), \,\ \widetilde{f} \in \widetilde{\mathcal{H}}(N). \nonumber
\end{eqnarray}

We similarly define for $(G,\xi) \in \widetilde{\mathcal{E}}_{\ellip}(N)$, and $\psi^N \in \widetilde{\Psi}(N)$ the following:

\begin{eqnarray}
I_{\disc,\psi^N,\xi}^{G}(f) =  \sum_{\xi(c) = c(\psi^N)  } I^G_{\disc,c,t(\psi^N)}   (f)                  ,  \,\  f \in \mathcal{H}(G).
\end{eqnarray}

\begin{eqnarray}
R_{\disc,\psi^N,\xi}^G = \bigoplus_{\xi(c) = c(\psi^N)} R_{\disc,c,t(\psi^N)}^G.
\end{eqnarray}

Then for $(G,\xi) \in \widetilde{\mathcal{E}}_{\simp}(N)$ (i.e. $G=U_{E/F}(N)$), and $c^N \in \widetilde{\mathcal{C}}_{\mathbf{A}_F}(N)$, we have the following decomposition parallel to (4.3.10), which follows immediately from the previous discussions:
\begin{eqnarray}
& & \\
& & I^{G}_{\disc,c^N,t,\xi}(f) = \sum_{(G^{\prime},\zeta^{\prime}) \in \mathcal{E}_{\ellip}(G)} \iota(G,G^{\prime}) \widehat{S}^{G^{\prime}}_{\disc,c^N,t,\xi \circ \zeta^{\prime}} (f^{G^{\prime}}), \,\ f \in \mathcal{H}(G). \nonumber
\end{eqnarray}
Similarly if $c^N=c(\psi^N)$ for $\psi^N \in \widetilde{\Psi}(N)$ we have:

\begin{eqnarray}
& & \\
& & I^{G}_{\disc,\psi^N,\xi}(f) = \sum_{(G^{\prime} ,\zeta^{\prime}) \in \mathcal{E}_{\ellip}(G)} \iota(G,G^{\prime}) \widehat{S}^{G^{\prime}}_{\disc,\psi^N,\xi \circ \zeta^{\prime}} (f^{G^{\prime}}), \,\ f \in \mathcal{H}(G). \nonumber
\end{eqnarray}

\bigskip

Before we prove the next proposition, we first formulate the global version of the notion of compatible family (3.1.4) from section 3.1. We formulate it as a family of functions:
\[
\mathcal{F}=\{f \in \mathcal{H}(G) \big| \,\ (G,\xi) \in \widetilde{\mathcal{E}}(N) \}
\] 
parametrized by equivalence class of global endoscopic data of $\widetilde{G}_{E/F}(N)$ (as noted in section 2.4 there are only finitely many equivalence classes of endoscopic data for $\widetilde{G}_{E/F}(N)$). We say that $\mathcal{F}$ is a decomposable compatible family, if for each $v$ there is a local compatible family of functions
\[
\mathcal{F}_v=\{ f_v \in \mathcal{H}(G_v)| \,\ (G_v,\xi_v) \in \widetilde{\mathcal{E}}_v(N) \}
\]
such that for $(G,\xi) \in \widetilde{\mathcal{E}}(N)$ the corresponding functions satisfy
\[
f = \prod_v f_v
\]
We then define more generally that the family $\mathcal{F}$ is called a compatible family, if there is a finite set of decomposable compatible families
\[
\mathcal{F}_i=\{f_i \in \mathcal{H}(G)| \,\ (G,\xi) \in \widetilde{\mathcal{E}}(N)  \}
\]
such that the the corresponding functions satisfiy
\[
f = \sum_i f_i.
\]
If $\mathcal{F}$ is a compatible family of functions as above, then the stable orbital integrals $\{f^G\}$ are determined by the functions $f \in \mathcal{H}(G)$ associated to $(G,\xi) \in \widetilde{\mathcal{E}}_{\ellip}(N)$. Hence we will often specify functions of a compatible family for the elliptic datum $\widetilde{\mathcal{E}}_{\ellip}(N)$.

It then follows as an immediate corollary of part (a) of proposition 3.1.1 that a family of functions $\mathcal{F}$ is a compatible family, if and only if there exists a function $\widetilde{f} \in \widetilde{\mathcal{H}}(N)$, such that for each datum $(G,\xi) \in \widetilde{\mathcal{E}}_{\ellip}(N)$ and the corresponding function $f \in \mathcal{H}(G)$ attached to $(G,\xi)$, we have:
\[
f^G = \widetilde{f}^{(G,\xi)}.
\]

\bigskip

We can now prove the following proposition, which is based on a preliminary comparison of the trace formulas. A more elaborate comparison is the subject of the next two sections.

\begin{proposition}
For $G = U_{E/F}(N)$, and $\xi: \leftexp{L}{U_{E/F}(N)} \hookrightarrow \leftexp{L}{G_{E/F}(N)}$ (thus $(G,\xi) \in \widetilde{\mathcal{E}}_{\simp}(N)$), we have 
\begin{eqnarray}
I^{G}_{\disc,c^N,t,\xi}(f)=0 = S^{G}_{\disc,c^N,t,\xi}(f), \,\ f \in \mathcal{H}(G)
\end{eqnarray}
unless $(c^N,t)=(c(\psi^N),t(\psi^N))$ for $\psi^N \in \widetilde{\Psi}(N)$.
\end{proposition}
\begin{proof}
We prove this by induction, assuming the proposition holds for $N_- < N$. Thus we assume that $(c^N,t)$ is not of the form $(c(\psi^N),t(\psi^N))$ for any $\psi^N \in \widetilde{\Psi}(N)$. We first analyze $I_{\disc,c^N,\xi}^G$ from the spectral expansion (4.3.4). If $M$ is a standard proper Levi subgroup of $G$, then we  have:
\[
M = G_{E/F}(N_1) \times \cdots G_{E/F}(N_r) \times U_{E/F}(N_-)
\]
with $2 N_1 \cdots+2 N_r +N_- =N$, and $N_- < N$. By induction the proposition holds for $U_{E/F}(N_-)$. Since the corresponding assertions holds for the factors $G_{E/F}(N_i)$ (as follows from the results of Moeglin-Waldspurger and Jacquet-Shalika), we deduce that the assertion holds for $M$ itself. Therefore it follows that $\mathcal{I}_{P,c^N,t,\xi} : = \oplus_{\xi(c) = c^N} \mathcal{I}_{P,c,t}$ is identifically zero. Hence on the right hand side of (4.3.4) we are reduced to the term $M=G$, i.e. we have:
\begin{eqnarray}
I^G_{\disc,c^N,t,\xi} = \tr R^G_{\disc,c^N,t,\xi}, \,\ f \in \mathcal{H}(G).
\end{eqnarray}

\noindent Next we look at the endoscopic expansion (4.3.16). If $(G^{\prime},\zeta^{\prime}) \in \mathcal{E}_{\ellip}(G) \smallsetminus \mathcal{E}_{\simp}(G)$, then $G^{\prime}$ is a proper product:
\[
G^{\prime} =G_1^{\prime} \times G_2^{\prime} =U_{E/F}(N_1^{\prime}) \times U_{E/F}(N_2^{\prime}), \,\ N_1^{\prime},N_2^{\prime} < N
\]
and we can write $\xi \circ \zeta^{\prime} =\xi^{\prime}_1 \times \xi^{\prime}_2$ as in (2.1.14), with $\xi^{\prime}_1:\leftexp{L}{U_{E/F}(N^{\prime}_1)} \hookrightarrow \leftexp{L}{G_{E/F}(N^{\prime}_1)}$ and $\xi^{\prime}_2:\leftexp{L}{U_{E/F}(N^{\prime}_2)} \hookrightarrow \leftexp{L}{G_{E/F}(N^{\prime}_2)}$. Consider any pair $c_1^{N_1} \in \widetilde{\mathcal{C}}_{\mathbf{A}_F}(N_1)$, $c_2^{N_2} \in \widetilde{\mathcal{C}}_{\mathbf{A}_F}(N_2)$ such that 
\[
c^N = c_1^{N_1} \times c_2^{N_2}.
\]
By assumption $(c^N,t)$ is not of the form $(c(\psi^N),t(\psi^N))$ for any $\psi^N \in \widetilde{\Psi}(N)$. It follows that for at least one $i \in \{1,2\}$, we have $(c_i^{N_i},t)$ is not of the form $(c(\psi_i^{N_i}),t(\psi_i^{N_i}))$ for any $\psi_i^{N_i} \in \widetilde{\Psi}(N_i)$.

Thus by induction one of the linear forms $S^{G_1^{\prime}}_{\disc,c_1^{N_1},t,\xi^{\prime}_1}$ or $S^{G_2^{\prime}}_{\disc,c_2^{N_2},t,\xi^{\prime}_2}$ is zero. It follows that the linear form $S_{\disc,c^N,t, \xi \circ \zeta^{\prime}}^{G^{\prime}}$ is also zero (indeed for instance if $f^{\prime} = f^{\prime}_1 \times f^{\prime}_2$ is a decomposable function of $G^{\prime}$ then we have
\[
S^{G^{\prime}}_{\disc,c^N,t,\xi \circ \zeta^{\prime}}(f^{\prime})=  \sum_{c^N = c_1^{N_1} \times c_2^{N_2}} S^{G_1^{\prime}}_{\disc,c_1^{N_1},t,\xi_1^{\prime}}(f_1^{\prime})  \cdot  S^{G_2^{\prime}}_{\disc,c_2^{N_2},t,\xi_2^{\prime}}(f_2^{\prime})  =0 ).
\]

\noindent Thus by (4.3.16) we have:
\begin{eqnarray}
I^G_{\disc,c^N,t,\xi}(f) = S^{G}_{\disc,c^N,t,\xi}(f), \,\ f \in \mathcal{H}(G).
\end{eqnarray}
And (4.3.19) and (4.3.20) gives
\begin{eqnarray}
& & \\
& & I^G_{\disc,c^N,t,\xi}(f) = S^{G}_{\disc,c^N,t,\xi}(f) = \tr R^G_{\disc,c^N,t,\xi}(f), \,\ f \in \mathcal{H}(G). \nonumber
\end{eqnarray}

Now we use (4.3.12). By the same reasoning, the only terms that survive on the right hand side of (4.3.12) are the terms coming from simple twisted endoscopic data. Hence we have:
\begin{eqnarray}
& & \\
& & \widetilde{I}^N_{\disc,c^N,t}(\widetilde{f}) = \sum_{(G,\xi) \in \widetilde{\mathcal{E}}_{\simp}(N) } \widetilde{\iota}(N,G) \widehat{S}^G_{\disc,c^N,t,\xi}(\widetilde{f}^{(G,\xi)}) \nonumber
\end{eqnarray}

On the other hand, since $(c^N,t)$ is not equal to $(c(\psi^N),t(\psi^N))$ for any $\psi_N \in \widetilde{\Psi}(N)$, the left hand side of (4.3.22) is zero by the theorem of Jacquet-Shalika \cite{JS} and Moeglin-Waldspurger \cite{MW}. Thus we have: 
\begin{eqnarray}
 \sum_{(G,\xi) \in \widetilde{\mathcal{E}}_{\simp}(N) } \widetilde{\iota}(N,G) \widehat{S}^G_{\disc,c^N,t,\xi}(\widetilde{f}^{(G,\xi)}) =0
\end{eqnarray}

We can now replace the summand $\widehat{S}^G_{\disc,c^N,t,\xi}(\widetilde{f}^{(G,\xi)})$ in (4.3.23) by $S^G_{\disc,c^N,t,\xi}(f)$, for any compatible family of functions:
\[
\mathcal{F}=\{ f \in \mathcal{H}(G): \,\ (G,\xi) \in  \widetilde{\mathcal{E}}_{\ellip}(N) \}.
\]
Thus we have
\begin{eqnarray}
 \sum_{(G,\xi) \in \widetilde{\mathcal{E}}_{\simp}(N) } \widetilde{\iota}(N,G) S^G_{\disc,c^N,t,\xi}(f) =0.
\end{eqnarray}
Hence combining with (4.3.21) we obtain:
\begin{eqnarray}
 \sum_{(G,\xi) \in \widetilde{\mathcal{E}}_{\simp}(N) } \widetilde{\iota}(N,G)  \tr R^G_{\disc,c^N,t,\xi}(f)    =0.
\end{eqnarray}
for any compatible family of functions $\mathcal{F}$ as above.

Now the coefficients $\widetilde{\iota}(N,G)$ are positive, and the term $R^G_{\disc,c^N,t,\xi}$ is a linear combination with non-negative integer coefficients of irreducible admissible representations on $G(\mathbf{A})$. Hence we may apply the result on vanishing of coefficients (see lemma 4.3.6 below) to (4.3.25). We then conclude the vanishing of $R^G_{\disc,c^N,t,\xi}$ for any $(G,\xi) \in \widetilde{\mathcal{E}}_{\simp}(N)$. Hence
\[
\tr R_{\disc,c^N,t,\xi}^G(f) = S_{\disc,c^N,t,\xi}^G(f) =0 , \,\ f \in \mathcal{H}(G).
\]
By (4.3.21) again, we have $I_{\disc,c^N,t,\xi}^G(f) =0$. Thus conclude the proof.
\end{proof}

\bigskip

\begin{rem}
\end{rem}
The result on vanishing of coefficient, which is stated as lemma 4.3.6 below, is proved as proposition 3.5.1 of \cite{A1}. The proof in {\it loc. cit.} applies to a general $G$ (twisted or not), even though for our purpose we only need this for the twisted group $G=\widetilde{G}_{E/F}(N)$. As in \cite{A1}, this result on vanishing of coefficients forms the basis of the proofs of the main theorems established in this paper.

\begin{lemma} (proposition 3.5.1 of \cite{A1})
Suppose we are given an identity:
\begin{eqnarray}
\sum_{G^{\prime} \in \mathcal{E}_{\ellip}(G)} \sum_{\pi \in \Pi(G^{\prime}(\mathbf{A}_F))} c_{G^{\prime}}(\pi) f_{G^{\prime}}(\pi) \equiv 0
\end{eqnarray}
for any compatible family of functions $\mathcal{F}=\{f \in \mathcal{H}(G^{\prime})| \,\ G^{\prime} \in \mathcal{E}(G)$, and such that the coefficients $c_{G^{\prime}}(\pi)$ are non-negative real numbers, and that the coefficients $c_{G^{\prime}}(\pi)$ vanish outside a finite set of $\pi$ that depends only on the choice of a Hecke type for $\mathcal{F}$. Then all the coefficients $c_{G^{\prime}}(\pi)$ vanish.
\end{lemma}

\bigskip

Finally we state two corollaries of proposition 4.3.4:

\begin{corollary}
For $G=U_{E/F}(N)$, and $\xi:\leftexp{L}{U_{E/F}(N)} \hookrightarrow \leftexp{L}{G_{E/F}(N)}$, we have
\begin{eqnarray}
I_{\disc,t}^G(f) = \sum_{\substack{ \psi^N \in \widetilde{\Psi}(N), \\ t(\psi^N)=t}} I_{\disc,\psi^N ,\xi}^G(f)
\end{eqnarray}
and
\begin{eqnarray}
S_{\disc,t}^G(f) = \sum_{\substack{ \psi^N \in \widetilde{\Psi}(N), \\ t(\psi^N)=t }   } S_{\disc,\psi^N,\xi}^G(f).
\end{eqnarray}
\end{corollary}

\begin{corollary}
As above $G=U_{E/F}(N)$, and $\xi: \leftexp{L}{U_{E/F}(N) \hookrightarrow \leftexp{L}{G_{E/F}(N)}}$, and $c^N \in \widetilde{\mathcal{C}}_{\mathbf{A}_F}(N)$, we have
\[
L^2_{\disc,c^N,t,\xi}(G(F) \backslash G(\mathbf{A})) = 0
\]
unless $(c^N,t)=(c(\psi^N),t(\psi^N))$ for some $\psi^N \in \widetilde{\Psi}(N)$. Thus we have a decomposition:
\begin{eqnarray}
& & \\
& & L^2_{\disc}(G(F) \backslash G(\mathbf{A})) = \bigoplus_{\psi^N \in \widetilde{\Psi}(N)} L^2_{\disc,\psi^N,\xi} (G(F)  \backslash G(\mathbf{A})). \nonumber
\end{eqnarray}
\end{corollary}

Corollary 4.3.7 follows immediately from proposition 4.3.4, while the proof of corollary 4.3.8 from proposition 4.3.4 is the same as the proof of corollary 3.4.3 of \cite{A1}.

\bigskip

We note that corollary 4.3.8 gives in particular the existence of ``weak base change" associated to discrete automorphic representations on $U_{E/F}(N)$ (with respect to the $L$-embedding $\xi$). We would also want to cut down the set of parameters $\widetilde{\Psi}(N)$ to $\Psi_2(U_{E/F}(N),\xi)$ (the latter has yet to be defined). This requires more elaborate comparison of the trace formulas, which is the subject of section five and six. 

\bigskip

To simplify the notations in the next two sections, we will often abbreviate a twisted endoscopic datum $(G,\xi)$ of $\widetilde{G}_{E/F}(N)$ just as $G$, when the context is clear. Similarly we will abbreviate the terms
\[
I^G_{\disc,\psi^N,\xi}, \,\  R^G_{\disc,\psi^N,\xi}, \,\ S^G_{\disc,\psi^N,\xi}
\]
etc, just as 
\[
I^G_{\disc,\psi^N}, \,\  R^G_{\disc,\psi^N}, \,\ S^G_{\disc,\psi^N}
\]
and always keep in mind that these are defined with respect to the $L$-embedding $\xi: \leftexp{L}{G} \hookrightarrow \leftexp{L}{G_{E/F}(N)}$ that is part of the datum for $G=(G,\xi)$.

\bigskip

In section 6 we will need a stronger version of lemma 4.3.6. To state this we need some more notion. Thus in general $G$ is a connected reductive group over a field $F$ that we take temporarily to be local. Denote by $T(G)$ the set of $G(F)$-orbits of triples:
\begin{eqnarray}
\tau = (M,\pi,r)
\end{eqnarray}
where $M \subset G$ is a Levi subgroup, $\sigma \in \Pi_2(M)$, an irreducible representation of $M(F)$ that is square-integrable modulo centre, and $r \in R(\sigma)$ is an element in the representation theoretic $R$-group associated to $\sigma$ (see section 3.5 of \cite{A1} for a detailed discussion). We refer to equation (3.5.3) of \cite{A1} for the definition of the virtual character $f_G(\tau)$ for $\tau \in T(G)$ (we remark that in the general case as treated in {\it loc. cit.} one has to consider an extension of the representation theoretic $R$-group in order to split the associated 2-cocycles; however in the present case where $G$ is a unitary group their cohomology classes are trivial). 

Going back to the situation where $F$ is a global field. We state:
\begin{proposition}
Suppose there exists a place $v$ of $F$, and a $G_1 \in \widetilde{\mathcal{E}}_{\simp}(N)$, such that the following holds:
\begin{eqnarray}
& & \sum_{G^{\prime} \in \mathcal{E}_{\ellip}(G)} \sum_{\pi \in \Pi(G^{\prime}(\mathbf{A}_F))} c_{G^{\prime}}(\pi) f_{G^{\prime}}(\pi) \\
&=&  \sum_{\tau_v \in T(G_{1,v})} d_1(\tau_v,f_1^{v} ) f_{1,v}(\tau_v) \nonumber
\end{eqnarray}
for every compatible family of functions $\mathcal{F}$, such that the function $f_1$ associated to $G_1$ is a product
\[
f_1 = f_{1,v} f_1^v, \,\ f_{1,v} \in \mathcal{H}(G_{1,v}), f_1^{v} \in \mathcal{H}(G_1(\mathbf{A}_F^v)).
\] 

\noindent We suppose that the coefficients $c_{G^{\prime}}(\pi)$ are as in proposition 4.3.9, in particular are non-negative real numbers. Also, we suppose that the coefficients $d_1(\tau_v,f_1^{v} )$, as a function of $\tau_v$, is supported on a finite set that depends only on the choice of a Hecke type for $f_v$, and vanishes for any $\tau_v$ of the form $(M_v,\sigma_v,1)$. Then all the coefficients $c_{G^{\prime}}(\pi)$ and $d_1(\tau_v,f_1^{v} )$ vanish. 
\end{proposition}

\noindent Proposition 4.3.9 corresponds to corollary 3.5.3 of \cite{A1}, where it is deduced from proposition 3.5.1 of {\it loc. cit.}, by using in addition corollary 2.1.2 of {\it loc. cit.} (which is on surjectivity of the Kottwistz-Shelstad transfer). In our present context, we can deduce proposition 4.3.9 by the same argument as in \cite{A1}, with corollary 2.1.2 of {\it loc. cit.} being replaced by part (b) of proposition 3.1.1.

\section{\textbf{The Standard model}}

In this section we begin to study the term by term comparison in the spectral and endoscopic expansions of the discrete part of the trace formula, a process that is coined by Arthur the ``standard model". Two key inputs in this comparison are the stable multiplicity formula, and the global intertwining relation. 

We will also begin in section 5.3 to carry out the formal part of the induction argument for the proof of the main theorems.

\subsection{Stable multiplicity formula}

In this subsection, we state the stable multiplicity formula. Besides the global theorems already stated in section 2.5, the stable multiplicity formula is the main global result to be proved in this paper, which at the same time is the main driving force in the proof of other global theorems.

\bigskip

We briefly recall some quantities defined as in \cite{A9} which are necessary for the statement of the stable multiplicity formula. Thus in general we let $S$ be a connected component of a reductive group over $\mathbf{C}$. Denote by $S^+$ the reductive group generated by the component $S$, and by $S^0$ the identity component of $S^+$. Then $S$ is a bi-torsor under $S^0$. For $T$ a maximal torus of $S^0$, denote by 
\[
W(S) = \Norm(T,S)/T
\]
the Weyl set of $S$ with respect to $T$. In particular 
\[
W^0 := W(S^0)
\]
is the Weyl group of $S^0$ with respect to $T$, and $W(S)$ is again a bitorsor under $W^0$ (hence $|W(S)|=|W^0|$). We define the sign function
\[
\sgn^0: W(S) \rightarrow \{ \pm 1 \}
\] 
to be $(-1)$ raised to the number of positive roots of $(S^0,T)$ that are mapped under $w$ to negative roots.

Define the set of regular elements:
\begin{eqnarray}
W_{\reg}(S) =\{w \in W(S),\,\  \det(w-1)_{\mathfrak{a}_T} \neq 0         \}
\end{eqnarray}
where $\mathfrak{a}_T$ is the $\mathbf{R}$-vector space:
\[
\mathfrak{a}_T:= \Hom(X^*(T),\mathbf{R})
\]
equivalently, $W_{\reg}(S)$ is the set of elements $w \in W(S)$ whose action on $T$ has only finite number of fixed points. Define the number 
\begin{eqnarray}
i(S) := \frac{1}{|W(S)|} \sum_{w \in W_{\reg}(S)} \frac{\sgn^0(w)}{|\det(w-1)_{\mathfrak{a}_T}|}.
\end{eqnarray}

We also follow \cite{A9} in the definition of some centralizers. Thus we define
\[
Z(S) = \Cent(S,S^0)
\]
for the centralizer of $S$ in $S^0$, and similarly for $s \in S$:
\begin{eqnarray}
& & S_s := \Cent(s,S^0) \\
& & S_s^0 := (S_s)^0 = \Cent(s,S^0)^0. \nonumber 
\end{eqnarray}

Denote by $S_{ss}$ the set of semi-simple elements of $S$, and in general if $\Sigma \subset S$ is invariant under conjugation by $S^0$, then we define by $\mathcal{E}(\Sigma)$ the set of equivalence classes of elements in 
\[
\Sigma_{ss} = \Sigma \cap S_{ss}
\] 
with the equivalence relation being defined as follows: if $s,s^{\prime} \in \Sigma_{ss}$, then $s^{\prime} \sim s$ if 
\begin{eqnarray}
s^{\prime} =  z^0 s^0 s (s^0)^{-1}
\end{eqnarray}
with $s^0 \in S^0$ and $z \in Z(S^0_s)^0$. For our purpose we only need to consider the set of elliptic elements:
\begin{eqnarray}
S_{\ellip}=\{s \in S_{ss}, \,\ |Z(S_s^0)| < \infty  \}
\end{eqnarray}
then the equivalence relation on $S_{\ellip}$ reduces to $S^0$-conjugacy, and we put
\begin{eqnarray}
\mathcal{E}_{\ellip}(S) := \mathcal{E}(S_{\ellip}).
\end{eqnarray}
The set $\mathcal{E}_{\ellip}(S)$ is finite. Then we have the following:

\begin{proposition} (theorem 8.1 of \cite{A9})
There exists unique constants $\sigma(S_1)$ assigned to connected complex reductive group $S_1$, such that the following property holds: 

\bigskip

\noindent (a) For any connected component $S$ of a complex reductive group as above, if we put
\begin{eqnarray}
e(S) := \sum_{s \in \mathcal{E}_{\ellip}(S)} \frac{1}{|\pi_0(S_s)|} \sigma(S_s^0)
\end{eqnarray}
then we have the equality
\begin{eqnarray}
i(S)=e(S).
\end{eqnarray} 

\bigskip

\noindent (b) For any central subgroup $Z_1$ of $S_1$, we have
\begin{eqnarray}
\sigma(S_1) = \sigma(S_1/Z_1)|Z_1|^{-1}
\end{eqnarray}
(thus in particular $\sigma(S_1)=0$ if $S_1$ has an infinite centre). 
\end{proposition}

\bigskip

We now specialize these constructions. Thus let $G=(G,\xi) \in \widetilde{\mathcal{E}}_{\simp}(N)$ (hence the underlying simple twisted endoscopic group is $U_{E/F}(N)$). Given $\psi^N \in \widetilde{\Psi}(N)$, we assume the validity of theorem 2.4.2 (the first ``seed" theorem) for all the simple generic constituents of $\psi^N$. We will begin the formal proof of the theorems in section 5.3 as an induction on $N$. Thus from the induction hypothesis one can assume the validity of theorem 2.4.2 for parameters in $\widetilde{\Phi}_{\simp}(N_-)$ for $N_- <N$, and the only situation not already covered by the induction hypothesis is the case where $\psi^N$ itself is a simple generic parameter in $\widetilde{\Phi}_{\simp}(N)$. In any case with this premise we can define the $L$-homomorphism $\widetilde{\psi}^N: \mathcal{L}_{\psi^N} \rightarrow \leftexp{L}{G_{E/F}(N)}$. Then $\widetilde{\psi}^N$ factors through $\xi$ if and only if $\psi^N$ defines a parameter $\psi = (\psi^N,\widetilde{\psi})$ in $\Psi(G,\xi)$, in which case we have $\psi^N \in \xi_* \Psi(G,\xi)$ (thus we have $\widetilde{\psi}^N$ factors through $\xi$ if and only if $\psi^N \in \xi_* \Psi(G,\xi)$). Define
\begin{eqnarray}
m_{\psi^N}^G= \left \{ \begin{array}{c}  1 \mbox{ if } \psi^N \in \xi_* \Psi(G,\xi)   \\    0 \mbox{ if } \psi^N \notin \xi_* \Psi(G,\xi).
\end{array} \right.
\end{eqnarray}

If $\psi \in \Psi(G,\xi)$, we have defined in section 2 the groups $S_{\psi},\overline{S}_{\psi},\mathcal{S}_{\psi}$. We have also defined the sign character
\[
\epsilon_{\psi}=\epsilon_{\psi}^G: \mathcal{S}_{\psi} \rightarrow \{ \pm 1 \}
\]
associated to $\psi$ in terms of global symplectic root numbers. 

The statement of the stable multiplicity formula is predicated on the validity of theorem 2.4.10 for the parameter $\psi$, and also part (a) of theorem 3.2.1 of applied to the localization of the parameter $\psi$ at each place $v$ of $F$ (in any case, the proof of the stable multiplicity formula, together with the proof of the other local and global theorems, are established simultaneously at the end of the induction argument in section 9. Thus assume the validity of these two theorems for $\psi=(\psi^N,\widetilde{\psi}) \in \Psi(G,\xi)$ as above. Then as in the discussion in section 2.4, for each prime $v$ of $F$, the localization $\psi^N_{v}$ (as a parameter of $L_{F_v} \times \SU(2)$) factors through $\xi_v$ (here $\xi_v$ being the localization of $\xi$ at $v$), i.e. 
\begin{eqnarray}
\psi^N_{v} = \xi_v \circ \psi_v
\end{eqnarray}
and $\psi_v \in \Psi_v^+(G)$ (note that $\psi_v$ is then uniquely determined up to $\widehat{G}$-conjugacy by $\psi^N_{v}$, since $\widetilde{\Out}_N(G)$ is trivial). In addition, part (a) of theorem 3.2.1 gives the stable linear form:
\begin{eqnarray*}
f_v \mapsto f_v^{G_v}(\psi_v), \,\ f_v \in \mathcal{H}(G_v).  
\end{eqnarray*} 
Note that if $v$ splits in $E$ then these assertions are elementary, namely that the corresponding assertions of theorem 2.4.10 and theorem 3.2.1 in this case is already known.

\noindent We can then define the global stable linear form 
\begin{eqnarray}
f \mapsto f^G(\psi):=\prod_v f_v^{G_v}(\psi_v), \,\ f = \prod_v f_v \in \mathcal{H}(G)
\end{eqnarray}
with almost all terms in the product being equal to one (and as usual the linear form being extended to non-decomposable functions by linearity).

\bigskip

We can now state:
\begin{theorem} (the stable multiplicity formula) For $\psi^N \in \widetilde{\Psi}(N)$, we have
\begin{eqnarray}
S_{\disc,\psi^N,\xi}^G(f) = \frac{1}{|\mathcal{S}_{\psi}|} \epsilon^G_{\psi}(s_{\psi}) \sigma(\overline{S}^0_{\psi}) f^G(\psi), \,\ f \in \mathcal{H}(G)
\end{eqnarray}
if $\psi^N \in \xi_* \Psi(G,\xi)$, i.e. that $\psi^N$ defines a parameter $\psi \in \Psi(G,\xi)$, and
\[
S_{\disc,\psi^N,\xi}^G(f) = 0, \,\ f \in \mathcal{H}(G)
\] 
if $\psi^N \notin \xi_* \Psi(G,\xi)$.
\end{theorem}

\bigskip

\begin{rem}
\end{rem}
We can state the stable multiplicity formula in both cases in the form:
\begin{eqnarray}
& & \\
& & S_{\disc,\psi^N,\xi}^G(f) = m_{\psi^N}^G|\mathcal{S}_{\psi}|^{-1} \epsilon^G_{\psi}(s_{\psi}) \sigma(\overline{S}^0_{\psi}) f^G(\psi), \,\ f \in \mathcal{H}(G). \nonumber
\end{eqnarray}
Namely that if $\psi^N \notin \xi_* \Psi(G,\xi)$, then $m_{\psi^N}^G=0$, and the right hand side of (5.1.14) is just interpreted as zero (even though the parameter $\psi$ in $\Psi(G,\xi)$ is not defined in this case). 

To simplify the notation we will also write $\epsilon^G_{\psi}(\psi)$ for the sign $\epsilon^G_{\psi}(s_{\psi})$.

\bigskip

\begin{rem}
\end{rem}

\noindent  An important special case of the stable multiplicity formula is the case where $\psi \in \Psi_2(G,\xi)$ is a square integrable parameter, in which case $\overline{S}_{\psi}$ is finite, i.e. $\overline{S}^0_{\psi}$ is trivial and hence $\sigma(\overline{S}_{\psi}^0)=1$.

\bigskip

With $G=(G,\xi) \in \widetilde{\mathcal{E}}_{\simp}(N)$ as before, we define the following chains of subsets of $\Psi(G)$, in terms of the centralizer group $S_{\psi}$ associated to a parameter $\psi \in \Psi(G)$ (here we allow ourselves to omit the explicit reference to the $L$-embedding in the notation for endoscopic datum and the set of global parameters):
\begin{eqnarray*}
& & \Psi_{\simp}(G) \subset \Psi_2(G) \subset \Psi_{\ellip}(G) \subset \Psi_{\disc}(G) \\
& & \Psi_{\simp}(G) \subset \Psi_2(G) \subset \Psi_{\sdisc}(G) \subset \Psi_{\disc}(G) 
\end{eqnarray*}
defined by the conditions:
\begin{eqnarray}
& & \Psi_{\simp}(G)=\{\psi \in \Psi(G), \,\ |\overline{S}_{\psi}|=1  \} \\
& & \Psi_2(G) = \{ \psi \in \Psi(G), \,\ |\overline{S}_{\psi}| < \infty   \} \nonumber \\
& & \Psi_{\ellip}(G) = \{\psi \in \Psi(G), \,\  |\overline{S}_{\psi,s}| < \infty \mbox{ for some } s \in \overline{S}_{\psi,ss}   \} \nonumber \\
& & \Psi_{\sdisc}(G) =\{\psi \in \Psi(G), \,\  |Z(\overline{S}^0_{\psi})|   < \infty       \} \nonumber \\
& & \Psi_{\disc}(G)=\{ \psi \in \Psi(G), \,\  |Z(\overline{S}_{\psi})| < \infty          \}. \nonumber
\end{eqnarray}

\bigskip

Similar considerations apply to the twisted group $G=\widetilde{G}_{E/F}(N)$. Thus given $\psi^N \in \Psi(\widetilde{G}_{E/F}(N))=\widetilde{\Psi}(N)$. As in section 2.4 we have the centralizer
\begin{eqnarray}
& & \\
& & S_{\psi^N}^* = \Cent(\Image \widetilde{\psi}^N , \widehat{\widetilde{G}}^0_{E/F}(N)) =\Cent( \Image \widetilde{\psi}^N, \widehat{G}_{E/F}(N)) \nonumber
\end{eqnarray}
as before, and the twisted centralizer
\begin{eqnarray}
& & \\
& & \widetilde{S}_{\psi^N} = \Cent(\Image \widetilde{\psi}^N, \widehat{\widetilde{G}}_{E/F}(N)) = \Cent(\Image \widetilde{\psi}^N,  \widehat{G}_{E/F}(N) \rtimes \widehat{\theta} ).   \nonumber
\end{eqnarray}
Then $\widetilde{S}_{\psi^N}$ is a bi-torsor under $S_{\psi^N}^*$, hence $\widetilde{S}_{\psi^N}$ is connected. Similarly we have
\begin{eqnarray*}
& & \overline{\widetilde{S}}_{\psi^N} = \widetilde{S}_{\psi^N}/Z(\widehat{G}_{E/F}(N))^{\Gamma_F}, \,\ \widetilde{\mathcal{S}}_{\psi^N} = \pi_0(\overline{\widetilde{S}}_{\psi^N}) \\
& & \overline{S}^*_{\psi^N} = S^*_{\psi^N}/Z(\widehat{G}_{E/F}(N))^{\Gamma_F}, \,\ \mathcal{S}^*_{\psi^N} = \pi_0(\overline{S}^*_{\psi^N})
\end{eqnarray*}
with $\overline{\widetilde{S}}_{\psi^N}$ and $\widetilde{\mathcal{S}}_{\psi^N}$ being bi-torsor under $\overline{S}^*_{\psi^N}$ and $\mathcal{S}^*_{\psi^N}$ respectively. In particular $\widetilde{\mathcal{S}}_{\psi^N}$ is a singleton (however $\widetilde{\mathcal{S}}_{\psi^N}$ still plays an important role in the  twisted trace formula for $\widetilde{G}_{E/F}(N)$).

\bigskip

We can define the chain of subsets of parameters
\begin{eqnarray}
& & \widetilde{\Psi}_{\simp}(N) \subset \widetilde{\Psi}_2(N) \subset \widetilde{\Psi}_{\ellip}  \subset \widetilde{\Psi}_{\disc}(N) \\
& & \widetilde{\Psi}_{\simp}(N) \subset \widetilde{\Psi}_2(N) \subset \widetilde{\Psi}_{\sdisc}(N) \subset \widetilde{\Psi}_{\disc}(N) \nonumber
\end{eqnarray}
by the same conditions as (5.1.15) above, using the twisted centralizer $\overline{\widetilde{S}}_{\psi_N}$. In fact we have $\widetilde{\Psi}_{\simp}(N) = \widetilde{\Psi}_2(N)$ and $\widetilde{\Psi}_{\sdisc}(N)=\widetilde{\Psi}_{\disc}(N)$, since $\overline{\widetilde{S}}_{\psi_N}$ is connected. Thus the chain in (5.1.18) reduces to $\widetilde{\Psi}_{\simp}(N) \subset \widetilde{\Psi}_{\ellip}(N) \subset \widetilde{\Psi}_{\disc}(N)$. Also the definition of $\widetilde{\Psi}_{\simp}(N)$ and $\widetilde{\Psi}_{\ellip}(N)$ coincides with that defined as in section 2.3.

\subsection{The global intertwining relation, part I}

In \cite{A1}, Arthur termed the standard model the process of term by term comparison of the spectral and endoscopic expansion of the discrete part of the trace formula. In this term by term comparison, an important role is played by the global intertwining relation, which is a corollary of the local intertwining relation.

For our purpose we need to consider both the standard case where $G=(G,\xi) \in \widetilde{\mathcal{E}}_{\simp}(N)$ is a simple twisted endoscopic datum (thus either $(G,\xi)=(U_{E/F}(N),\xi)$), or $G$ is the twisted group $\widetilde{G}_{E/F}(N)$.

We first consider the case $G=(G,\xi) \in \widetilde{\mathcal{E}}_{\simp}(N)$. Given $\psi=(\psi^N,\widetilde{\psi}) \in \Psi(G,\xi)$ (a condition which in particular entails the validity of the ``seed" theorem 2.4.2 for the simple generic constituents of $\psi^N$), then as in the local situation of section 3.4, we can choose $M$ a Levi-subgroup of $G$, which is uniquely determined by $\psi$ up to conjugation by $G$, and $\psi_M \in \Psi_2(M,\xi)$, such that $\widetilde{\psi}$ is the composition of $\widetilde{\psi}_M$ with the $L$-embedding $\leftexp{L}{M} \hookrightarrow \leftexp{L}{G}$. More precisely, in this case 
\begin{eqnarray*}
A_{\widehat{M}} := (Z(\widehat{M})^{\Gamma_F})^0
\end{eqnarray*}
is a maximal torus of $S_{\psi}^0$, which we denote as $T_{\psi}$, in which case $\widehat{M}= \Cent(T_{\psi},\widehat{G})$. And
\begin{eqnarray*}
\overline{T}_{\psi} := A_{\widehat{M}}/ A_{\widehat{M}} \cap Z(\widehat{G})^{\Gamma_F}
\end{eqnarray*}
is a maximal torus in $\overline{S}_{\psi}$. In particular we can identify $\mathfrak{a}^*_{\overline{T}_{\psi}}$, the linear dual of $\mathfrak{a}_{\overline{T}_{\psi}}$, as:
\begin{eqnarray}
\mathfrak{a}_{\overline{T}_{\psi}}^* \cong \mathfrak{a}^G_M.
\end{eqnarray}
As in the local situation, we put $\mathcal{S}^1_{\psi}:= \mathcal{S}_{\psi_M}$, and define the groups $\mathfrak{N}_{\psi}, W_{\psi}^0,W_{\psi},R_{\psi}$ as in the local situation.

We then have the global version of (3.4.2) of the commutative diagram of short exact sequences, which plays a crucial role in the comparison of the spectral and endoscopic expansions of trace formula:

\begin{eqnarray}
\xymatrix{& & 1 \ar[d] & 1 \ar[d] \\
 & & W_{\psi}^0 \ar@{=}[r] \ar[d]  & W_{\psi}^0 \ar[d] \\
1  \ar[r] & \mathcal{S}^1_{\psi}\ar[r] \ar@{=}[d] & \mathfrak{N}_{\psi}  \ar[r] \ar[d]  \ar@<.8 ex>@{<--}[d] & W_{\psi} \ar[r] \ar[d] \ar@<.8 ex>@{<--}[d] & 1 \\
1 \ar[r] &  \mathcal{S}^1_{\psi} \ar[r] & \mathcal{S}_{\psi} \ar[r] \ar[d] & R_{\psi}\ar[r] \ar[d] & 1 \\ 
& & 1 & 1 
  }
\end{eqnarray}

As in the local case the splittings of the vertical short exact sequence are determined by the choice of a parabolic subgroup $P \in \mathcal{P}(M)$ of $G$. The choice of $P$ allows the identification $W(\widehat{M}) \cong W(M)$ with respect to which the identification (5.2.1) is equivariant with respect to the action of $W(\widehat{M})$ and $W(M)$. Given $u \in \mathfrak{N}_{\psi}$, we denote by $w_u$ and $x_u$ the image of $u$ in $W_{\psi}$ and $\mathcal{S}_{\psi}$ respectively. We form the twisted group
\[
\widetilde{M}_u = M \rtimes \widetilde{w}_u
\]
and we can identify the twisted centralizer $\widetilde{\mathcal{S}}_{\psi_M,u}=\mathcal{S}_{\psi_M}(\widetilde{M}_u)$ as the coset $\mathcal{S}^1_{\psi} u$ 
\begin{eqnarray*}
\widetilde{\mathcal{S}}_{\psi_M,u} \cong \mathcal{S}^1_{\psi} u = \mathfrak{N}_{\psi}(w_u)
\end{eqnarray*}
with $\mathfrak{N}_{\psi}(w_u)$ stands for the fibre of $\mathfrak{N}_{\psi}$ over $w_u$ under the second horizontal short exact sequence of (5.2.2). We also denote by $\widetilde{u}$ the element $u$ regarded as an element of $\widetilde{\mathcal{S}}_{\psi_M,u}$. 

\bigskip

We now define the global linear form
\[
f \mapsto f_G(\psi,u), \,\ f \in \mathcal{H}(G)
\]
for $\psi \in \Psi(G,\xi)$ and $u \in \mathfrak{N}_{\psi}$. To make the inductive assumptions clear, we first consider the case that $M \neq G$, i.e. $\psi \notin \Psi_2(G,\xi)$. Then we can assume as part of the induction hypothesis that all the local and global theorems of section 2.4 and 2.5 are valid for $M$. In particular, corollary 2.4.11 applied to $M$ allows us to define, for each prime $v$ of $F$, the localization $\psi_{M,v} \in \Psi^+_{v}(M)$ of $\psi_M$ at $v$, and that the local packet $\Pi_{\psi_{M,v}}$ is defined (and hence we can form the global packet $\Pi_{\psi_M}$). Denote by $G_v,M_v,P_v$ the localization of $G,M,P$ at $v$. Then the localization $\psi_v \in \Psi^+_v(G)$ is defined and is given by the composition of $\psi_{M,v}$ with the $L$-embedding $\leftexp{L}{M_v} \rightarrow \leftexp{L}{G_v}$. It is then immediate that we have a morphism of the global diagram (5.2.2) to the local diagram (3.4.2) for each prime $v$ of $F$, provided we note the following: recall that in the local discussion of section 3.3 and 3.4 we also need to fix a $L$-embedding $\leftexp{L}{U_{E_v/F_v}(N)} \hookrightarrow \leftexp{L}{G_{E_v/F_v}(N)}$; in the global to local context, we always choose the one given by the localization $\xi_v$ of $\xi$ at $v$.

We then have the local linear form:
\[
f_v \mapsto f_{v,G_v}(\psi_v,u_v), \,\ f_v \in \mathcal{H}(G_v)
\] 
as defined in (3.4.9); here $u_v$ is the imgae of $u \in \mathfrak{N}_{\psi}$ in $\mathfrak{N}_{\psi_v}(G_v,M_v)$. More precisely, this linear form is the one obtained by analytic continuation from the one defined as in (3.4.9) (for parameters in $\Psi_v(G)$) to parameters in $\Psi_v^+(G)$, {\it c.f.} the discussion on p.30. Hence we can define the global linear form
\begin{eqnarray}
f_G(\psi,u) = \prod_v f_{v,G_v}(\psi_v,u_v), \,\ f=\prod_v f_v \in \mathcal{H}(G).
\end{eqnarray}
Note that almost all factors in (5.2.3) are equal to one, by part (b) of proposition 3.5.3. In explicit form (5.2.3) is given by
\begin{eqnarray}
& & \\
& & f_G(\psi,u) =  \sum_{\pi_M \in \Pi_{\psi_M}} \langle \widetilde{u},\widetilde{\pi}_M \rangle \tr( R_P(w_u,\widetilde{\pi}_M,\psi_M) \mathcal{I}_P(\pi_M,f)) \nonumber
\end{eqnarray}
with
\begin{eqnarray}
 R_P(w_u,\widetilde{\pi}_M,\psi_M) = \bigotimes_v  R_{P_v}(w_{u_v},\widetilde{\pi}_{M,v},\psi_{M,v}) 
\end{eqnarray}
is the global normalized intertwining operator, and
\begin{eqnarray}
  \langle \widetilde{u}, \widetilde{\pi}_M  \rangle = \prod_v \langle \widetilde{u}_v,\widetilde{\pi}_{M,v}   \rangle 
\end{eqnarray}
with $\langle \widetilde{u}_v,\widetilde{\pi}_{M,v}   \rangle$ being the extension of the local pairing $\langle \cdot , \cdot \rangle$ on $\mathcal{S}_{\psi_{M,v}} \times \Pi_{\psi_{M,v}}$ as in section 3.4 (before remark 3.4.1). 

\bigskip

In the case where $M=G$, i.e. when $\psi \in \Psi_2(G,\xi)$, then we have $W_{\psi}^0,W_{\psi},R_{\psi}$ being trivial, and $\mathcal{S}^1_{\psi}=\mathcal{S}_{\psi}=\mathfrak{N}_{\psi}$, and the normalized global intertwining operators are trivial. Thus the definition of the linear form $f_G(\psi,u)$ reduces to the existence of the global packet $\Pi_{\psi}$ and the corresponding pairing on $\mathcal{S}_{\psi} \times \Pi_{\psi}$ (again implicit is the validity of corollary 2.4.11 for $\psi$). Thus assuming its existence for $\psi$, we put
\begin{eqnarray}
f_G(\psi,u) = \sum_{\pi \in \Pi_{\psi}} \langle u  ,\pi \rangle f_G(\pi), \,\ f \in \mathcal{H}(G), u \in \mathfrak{N}_{\psi} = \mathcal{S}_{\psi}.
\end{eqnarray}

\bigskip

We now define the endoscopic counterpart of the linear form (5.2.3), again assuming the validity of corollary 2.4.11 for $\psi$. As in the local situation we have a correspondence
\begin{eqnarray}
(G^{\prime},\psi^{\prime}) \leftrightarrow (\psi,s)
\end{eqnarray}
with the interpretation that $G^{\prime}=(G^{\prime},\zeta^{\prime})$ is an endoscopic datum of $G$, and $\psi^{\prime} \in \Psi(G^{\prime}, \xi \circ \zeta^{\prime})$ 
(here for example if $G^{\prime} = U(N_1) \times U(N_2)$ and if $\xi \circ \zeta^{\prime} = \xi_1^{\prime} \times \xi_2^{\prime}$ in accordance with (2.1.14), with $\xi^{\prime}_i : \leftexp{L}{U(N_i)} \hookrightarrow \leftexp{L}{G(N_i)}$, then $\Psi(G^{\prime},\xi \circ \zeta^{\prime}) = \Psi(U(N_1),\xi_1^{\prime}) \times \Psi(U(N_2),\xi_2^{\prime})$). Then corollary 2.4.11 is also valid for $\psi^{\prime}$.

\noindent We make the assumption that for any $s \in \overline{S}_{\psi}$ the global stable linear form 
\[
f^{\prime} \mapsto (f^{\prime})^{G^{\prime}}(\psi^{\prime})
\]
is already defined on $\mathcal{H}(G^{\prime})$. We denote the linear form
\[
f \mapsto f^{G^{\prime}}(\psi^{\prime}), \,\ f \in \mathcal{H}(G)
\]
as $f^{\prime}_G(\psi,s)$.

\bigskip

We can now state the following:

\begin{theorem} With the above notations, we have:

\noindent (a) For any $w^0 \in W^0_{\psi}$ we have the triviality of the normalized global intertwining operator:
\begin{eqnarray}
R_P(w^0,\widetilde{\pi}_M,\psi_M)=1.
\end{eqnarray}

\bigskip

\noindent (b) (The global intertwining relation) For any $u \in \mathfrak{N}_{\psi}$ and $s \in \overline{S}_{\psi,ss}$ such that the image of $s$ in $\mathcal{S}_{\psi}$ is equal to $x_u$, we have
\begin{eqnarray}
f_G(\psi,u) = f^{\prime}_G(\psi,s_{\psi}s).
\end{eqnarray} 

\end{theorem}

\bigskip

\begin{rem}
\end{rem}
The global intertwining relation is of course a direct corollary of the local intertwining relation (theorem 3.4.3). As in the local case the global intertwining relation implies that the linear form $f_G(\psi,u)$ depends only on the image of $u$ in $\mathcal{S}_{\psi}$, and this is consistent with part (a) of the theorem 5.2.1.   

\bigskip

As in the local situation of section 3.5, we also need to formulate the global intertwining relation for the twisted group $\widetilde{G}_{E/F}(N)$. In the context of the global diagram (5.2.2) (and similarly for the local diagram (3.4.2)) $M$ is a Levi subgroup of $\widetilde{G}^0_{E/F}(N)=G_{E/F}(N)$ (which is just a finite product of $G_{E/F}(N_i)$'s), such that $\psi^N$ is the image of a parameter $\psi^M \in \Psi_2(M)$ under the $L$-embedding $\leftexp{L}{M} \hookrightarrow \leftexp{L}{G_{E/F}(N)}$. We denote by $\pi_{\psi^M}$ the discrete automorphic representation of $M(\mathbf{A}_F)$ corresponding to $\psi^M$. Here $W_{\psi^N}^0$ and $\mathcal{S}_{\psi^N}^1$ are defined with respect to the identity component $G_{E/F}(N)$, while $\mathfrak{N}_{\psi^N},W_{\psi^N}, \mathcal{S}_{\psi^N},R_{\psi^N}$ are defined with respect to $\widetilde{G}_{E/F}(N)$. In fact $\mathcal{S}_{\psi^N}^1$ is trivial and $\mathcal{S}_{\psi^N}$ and $R_{\psi^N}$ are all singleton, by virtue of the connectedness of the (twisted) centralizer $\overline{\widetilde{S}}_{\psi^N}$, and thus $\mathfrak{N}_{\psi^N}=W_{\psi^N}$ and both are bi-torsors under $W_{\psi^N}^0$.

\noindent The global diagram (5.2.2) can of course be formulated for the untwisted group $G_{E/F}(N)$. In this case we denote the correspondng objects as $\mathfrak{N}_{\psi^N}^*,W^*_{\psi^N}$ etc. to emphasize that they are defined with respect to $G_{E/F}(N)$. Note that in fact we have we have $W_{\psi^N}^0 = W^*_{\psi^N} = (W^*_{\psi^N})^0$.  

In both the twisted case and untwisted case $\overline{T}_{\psi^N}$ is the maximal torus of $\overline{\widetilde{S}}_{\psi^N}^0  = \overline{S}_{\psi^N}^*$, and we have the identification
\begin{eqnarray}
\mathfrak{a}^*_{\overline{T}_{\psi^N}} \cong \mathfrak{a}_M^{\widetilde{G}(N)^0} = \mathfrak{a}_M^{G(N)}.
\end{eqnarray}

Back to the twisted case, as in (5.2.3), we can formulate the global spectral distribution, by taking the product of the corresponding local distributions (3.5.13): for $u =w \in \mathfrak{N}_{\psi^N} = W_{\psi^N}$, put
\begin{eqnarray}
\widetilde{f}_N(\psi^N,u)  &:=& \prod_v \widetilde{f}_{v,N}(\psi^N,u_v) \\
&=&  \tr ( R_P( w, \widetilde{\pi}_{\psi^M} ,\psi_M) \mathcal{I}_P(\pi_{\psi^M}  ,\widetilde{f}  ) )        \nonumber
\end{eqnarray}
here
\[
R_P( w, \widetilde{\pi}_{\psi^M} ,\psi_M)= R_P( w, \widetilde{\pi}_{\psi^M} ) = \bigotimes_v R_P( w_v, \widetilde{\pi}_{\psi^M_v} )
\]
with the local normalized twisted intertwining operator $R_P( w_v, \widetilde{\pi}_{\psi^M_v} )$ as in (3.5.10). We define the ``untwisted" normalized intertwining operator $R_P(w^0,\widetilde{\pi}_M)$ for $w^0 \in W_{\psi^N}^0=W_{\psi^N}^*$ in a similar manner.

\bigskip

\noindent In fact, by virtue of (3.5.18), applied to each of the local linear forms $\widetilde{f}_{v,N}(\psi^N_v,u_v)$, we see that $\widetilde{f}_N(\psi^N,u)$ is independent of $u$; more precisely we have 
\begin{eqnarray}
\widetilde{f}_N(\psi^N,u) &=& \widetilde{f}_N(\psi^N) =\tr (\mathcal{I}_P(\pi_{\psi^M},N)   \mathcal{I}_P(\pi_{\psi^M},\widetilde{f}) )   \nonumber
\end{eqnarray}
here
\begin{eqnarray}
& &   \mathcal{I}_{P}(\pi_{\psi^M},N) : \mathcal{H}_P(\pi_{\psi^M}) \rightarrow \widetilde{\mathcal{H}}_P(\pi_{\psi^M})  \\
& & \mathcal{I}_P(\pi_{\psi^M},N) = \bigotimes_v \mathcal{I}_{P_v}(\pi_{\psi^M_v},N) \nonumber
\end{eqnarray}
with $ \mathcal{I}_{P_v}(\pi_{\psi^M_v},N)$ being the local intertwining operator as defined in (3.5.15).

\bigskip

In the same manner we can define the global endoscopic distribution: suppose $s \in \overline{\widetilde{S}}_{\psi^N}$, and $(G^{\prime},\psi^{\prime})$ is the pair that corresponds to $(\psi^N,s)$ (here $G^{\prime} \in \widetilde{\mathcal{E}}(N)$ and $\psi^{\prime} \in \Psi(G^{\prime})$). Assume that the global stable linear form $f^{G^{\prime}}(\psi^{\prime})=\prod_v f_v^{G^{\prime}_v}(\psi^{\prime}_v)$ on $\mathcal{H}(G^{\prime})$ is defined. Put:
\begin{eqnarray}
\widetilde{f}^{\prime}_N(\psi^N,s) = \widetilde{f}^{G^{\prime}}_N(\psi^N,s)  :&=& \prod_{v} \widetilde{f}^{G^{\prime}_v}_{v,N}(\psi^N_v,s_v) = \widetilde{f}^{G^{\prime}}(\psi^{\prime}). 
\end{eqnarray}

\noindent With this setup we then have the following proposition, which is a direct corollary of proposition 3.5.1 and corollary 3.5.2

\begin{proposition} 

\noindent (a) We have
\begin{eqnarray}
R_P(w^0,\widetilde{\pi}_M,\psi_M) \equiv 1 \mbox{ for } w^0 \in W_{\psi^N}^0.
\end{eqnarray}

\bigskip

\noindent (b) (global intertwining relation) Assume the global stable linear form is defined for any pair $(G^{\prime},\psi^{\prime})$, with $G^{\prime} \in \widetilde{\mathcal{E}}(N)$, and $\psi^{\prime} \in \Psi(G^{\prime})$. Then for $u \in \mathfrak{N}_{\psi^N}$ and semi-simple $s \in \overline{\widetilde{S}}_{\psi^N}$, we have
\begin{eqnarray}
\widetilde{f}_N(\psi^N,u) = \widetilde{f}^{\prime}_N(\psi^N,s_{\psi^N} s) , \,\ \widetilde{f}_N \in \widetilde{\mathcal{H}}(N).
\end{eqnarray}
\end{proposition}

Note that theorem 5.2.1 is a theorem that is still to be established in section 8, while proposition 5.2.3, which is the analogue of theorem 5.2.1 for $\widetilde{G}_{E/F}(N)$, is already known to hold. It is also more convenient to use a uniform notation; so in the future when we refer to proposition 5.2.3, we will use the notation as in the statement of theorem 5.2.1, wtih the understanding that $G=\widetilde{G}_{E/F}(N)$ in this case.

\subsection{The global intertwining relation, part II}

We begin in this subsection the formal induction proof of our main theorems for this paper. Recall that these are the ``seed" theorems 2.4.2 and 2.4.10 (and its corollary 2.4.11); the local classification theorem 2.5.1 and global classification theorem 2.5.2; theorem 2.5.4 about signs (with part (b) of theorem 2.5.4 being interpreted as statement on $N$ as in remark 2.5.7); theorem 3.2.1 and the local intertwining relation theorem 3.4.3; the stable multiplicity formula theorem 5.1.2 and theorem 5.2.1 (we have of course already used the induction hypothesis in the discussions in section 3.3 and 3.4). We will prove these theorems simultaneously by induction (which is to be completed in section 9). Each of these theorems is stated in terms of the integer $N$. And we assume as induction hypothesis that all these theorems are valid for any $N_- <N$. 

In this subsection we begin to analyze theorem 5.2.1, in particular the global intertwining relation. Recall that the global intertwining relation (5.2.10) implies in particular that the linear form $f_{G}(\psi,u)$ depends only on the image of $u \in \mathfrak{N}_{\psi}$ in $\mathcal{S}_{\psi}$, and that $f_G^{\prime}(\psi,s)$ depends only on the image of $s \in \overline{S}_{\psi,ss}$ in $\mathcal{S}_{\psi}$. We begin to analyze this (in)dependence for these two linear forms.

We will denote by $G$ either an element of $\widetilde{\mathcal{E}}_{\simp}(N)$ or $\widetilde{G}_{E/F}(N)$. Thus in the case of $\widetilde{\mathcal{E}}_{\simp}(N)$ we are suppressing the $L$-embedding in the endoscopic datum, in order to simplify the notation. Thus for instance if $G=(G,\xi) \in \widetilde{\mathcal{E}}_{\simp}(N)$, then we denote the set of parameters $\Psi(G,\xi)$ just as $\Psi(G)$, when there is no confusion. 

Thus let $\psi \in \Psi(G)$ (as in the last subsection, in the case where $G \in \widetilde{\mathcal{E}}_{\simp}(N)$ and $\psi=(\psi^N,\widetilde{\psi})$, then this is predicated on the validity of theorem 2.4.2 for the simple generic constituents of $\psi^N$; this follows from the induction hypothesis unless $\psi^N$ is a simple generic parameter in $\widetilde{\Phi}_{\simp}(N)$). We begin with treating $f_G^{\prime}(\psi,s)$. Given $s \in \overline{S}_{\psi}$ assume that the linear form 
\[
f^{\prime}_G(\psi,s), \,\ f \in \mathcal{H}(G)
\]
is defined. We begin by showing:

\begin{lemma}
The linear form 
\[
f^{\prime}_G(\psi,s), \,\ f \in \mathcal{H}(G)
\]
depends only on the image $x=x_s$ of $s$ in $\mathcal{S}_{\psi}$; in other words suppose that $s_1,s_2 \in \overline{S}_{\psi}$ are such that both $f^{\prime}_G(\psi,s_1)$ and $f^{\prime}_G(\psi,s_2)$ are defined. Then $f^{\prime}_G(\psi,s_1)=f^{\prime}_G(\psi,s_2)$ if $s_1$ and $s_2$ have the same image in $\mathcal{S}_{\psi}$.
\end{lemma}
\begin{proof}
This is proved by a standard descent argument as in section 4.5 of \cite{A1}.  Recall that the linear form $f^{\prime}_G(\psi,s)$ is defined in terms of the correspondnece:
\[
(G^{\prime},\psi^{\prime}) \leftrightarrow (\psi,s)
\]
as $f^{\prime}_G(\psi,s):=f^{G^{\prime}}(\psi^{\prime})$. 

In general if $s$ is replaced by an element $s_1$ that is conjugate to $s$ under $\overline{S}_{\psi}^0$, then under this correspondence pair $(G^{\prime},\psi^{\prime})$ is replaced by a pair $(G_1^{\prime},\psi_1^{\prime})$ under equivalence of endoscopic data. Then from the definition it follows that
\[
f^{G^{\prime} }(\psi^{\prime}) = f^{G_1^{\prime}}(\psi^{\prime}_1), \,\ f \in \mathcal{H}(G).
\] 

\noindent Let $(\overline{T}_{\psi},\overline{B}_{\psi})$ be a fixed Borel pair for $\overline{S}_{\psi}^0$. Since the semi-simple automorphism $\Int(s)$ of the complex connected reductive group $\overline{S}_{\psi}^0$ has to stabilize a Borel pair of $\overline{S}^0_{\psi}$, we see by replacing $s$ by a $\overline{S}_{\psi}^0$-conjugate $s_1 $ that we can assume $\Int(s_1)$ stabilize the original pair $(\overline{T}_{\psi},\overline{B}_{\psi})$. This last condition then determines $s_1 := s_x$ in terms of $x$ up to translation by $\overline{T}_{\psi}$. In particular 
\begin{eqnarray}
\overline{T}_{\psi,x} := \Cent(s_x,\overline{T}_{\psi})^0
\end{eqnarray}
does not depend on the choice of the representative $s_1=s_x$.

We thus only need to show that the linear form $f^{\prime}_G(\psi,s_x)$ is invariant when $s_x$ is replaced by an element in $ \overline{T}_{\psi} s_x $.  Now since $s_x \in \overline{S}_{\psi,ss}$ it is not hard to see that any element of $\overline{T}_{\psi}$ is of the form $ (t^{-1} \cdot s_x t s_x^{-1}) \cdot t_1$ for some $t \in \overline{T}_{\psi}$ and $t_1 \in \overline{T}_{\psi,x}$. Hence any element in $ \overline{T}_{\psi} s_x$ can be written as:
\[
t^{-1} s_x t t_1 = t^{-1} s_x t_1 t
\]
for some $t \in \overline{T}_{\psi}$ and $t_1 \in \overline{T}_{\psi,x}$. Hence it suffices to show that $f_G^{\prime}(\psi,s_x)$ is invariant under translating $s_x$ by an element in $\overline{T}_{\psi,x}$.

\noindent Put
\[
\widehat{M}^0_x := \Cent(\overline{T}_{\psi,x},\widehat{G}^0)
\]
then we can form $\leftexp{L}{M^0_x}$ (with the $L$-action on $\widehat{M}_x^0$ coming from the $L$-action on $\widehat{G}$). Then we have $\leftexp{L}{M} \subset \leftexp{L}{M^0_x}$ (recall that in the context of the previous subsection $\widehat{M}=\Cent(\overline{T}_{\psi},\widehat{G}^0)$), and thus $\leftexp{L}{M_x^0}$ is the $L$-group of a Levi subgroup $M_x^0$ of $G^0$ containing $M$. Furthermore since $s_x$ centralizes $T_{\psi,x}$ by definition, it stabilizes not only $\widehat{M}_x^0$ but also some $\Gamma_F$-invariant parabolic subgroup $\widehat{P}^0_x \in \mathcal{P}(M_x^0)$. Hence we can form the $\Gamma_F$-invariant {\it Levi subset}
\[
\widehat{M}_x := \widehat{M}_x^0 \rtimes \Int(s_x)
\]
that is dual to a Levi subset $M_x$ of $G$. Then the pair $(\psi,s_x)$ is the image of a pair 
\[
(\psi_{M_x},s_{M_x}), \,\ \psi_{M_x} \in \Psi(M_x), s_{M_x} \in \overline{S}_{\psi_{M_x}} 
\]
under the $L$-embeddings $\leftexp{L}{M_x^0} \subset \leftexp{L}{G^0}$ and $\leftexp{L}{M_x} \subset \leftexp{L}{G}$. This pair is in turn the image of a pair $(M_x^{\prime},\psi^{\prime}_{M_x})$ for an endoscopic datum $M_x^{\prime}$ of $M_x$, which can be treated as a Levi sub-datum of $G^{\prime}$. 

\noindent Put
\[
f^{\prime}_{M_x} := (f^{\prime})_{M^{\prime}_x}, \,\ f^{\prime} \in \mathcal{H}(G^{\prime})
\]
the descent of $f^{\prime}$ to $M^{\prime}_x$. By the descent property of the Langlands-Kottwitz-Shelstad transfer we have
\[
f^{\prime}(\psi^{\prime}) = f^{\prime}_{M_x}(\psi^{\prime}_{M_x})
\]
(here we are writing $f^{\prime}(\psi^{\prime})$ for $(f^{\prime})^{G^{\prime}}(\psi^{\prime})$, and similarly for $f^{\prime}_{M_x}(\psi^{\prime}_{M_x})$). Thus finally since $f^{\prime}_{M_x}(\psi^{\prime}_{M_x})$ does not change when $s_x$ is replace by a translate in $\overline{T}_{\psi,x}$ (indeed the image is the same in $\overline{S}_{\psi_{M_x}}$), we see that
\begin{eqnarray}
f_G^{\prime}(\psi,s_x) = f^{\prime}(\psi^{\prime}) =  f^{\prime}_{M_x}(\psi^{\prime}_{M_x})
\end{eqnarray}
is also invariant when $s_x$ is replaced by a translate in $\overline{T}_{\psi}$.  
\end{proof}

\begin{rem}
\end{rem}
Suppose that we have $\dim \overline{T}_{\psi,x} \geq 1$ for all $x \in \mathcal{S}_{\psi}$. Then the descent argument above, together with our induction hypothesis, implies that the linear form $f^{\prime}_G(\psi,s)$ is defined for any $s \in \overline{S}_{\psi}$ (and depends only on the image of $s$ in $\mathcal{S}_{\psi}$). 
\bigskip

We now turn to the spectral distribution $f_G(\psi,u)$ for $u \in \mathfrak{N}_{\psi}$. First suppose that $G$ is an element of $\widetilde{\mathcal{E}}_{\simp}(N)$. From the global diagram 5.2.2 again, we see that the fibre $\mathfrak{N}_{\psi}(x_u)$ of the projection of $\mathfrak{N}_{\psi}$ to $\mathcal{S}_{\psi}$ over $x_u$ is a bi-torsor under the subgroup $W_{\psi}^0$ of $\mathfrak{N}_{\psi}^*$. From the local discussion in section 3.4, we have the following: first the pairing
\begin{eqnarray*}
\langle \widetilde{u},\widetilde{\pi}_M \rangle
\end{eqnarray*}  
is unchanged when $u$ is replaced by a translate in $W_{\psi}^0$; second, for any $w_0 \in W_{\psi}^0$ and $u \in \mathfrak{N}_{\psi}$, we have:
\begin{eqnarray*}
R_p(w^0 w_u  ,\widetilde{\pi}_M,\psi_M) = R_P(w^0,\widetilde{\pi}_M,\psi_M) \cdot R_P(w_u,\widetilde{\pi}_M,\psi_M).
\end{eqnarray*} 
Hence from equation (5.2.4) we see that the assertion given by part (a) of theorem 5.2.1 (namely equation (5.2.9)): 
\[
R_P(w^0,\widetilde{\pi}_M,\psi_M)  \equiv 1
\]
implies that the linear form $f_G(\psi,u)$ depends only on the image $x_u$ of $u \in \mathfrak{N}_{\psi}$ in $\mathcal{S}_{\psi}$. 

In the case where $G=\widetilde{G}_{E/F}(N)$ we have seen in section 5.2 that the linear form $f_G(\psi,u)$ is independent of $u$ (this is by virtue of (3.5.18); see the discussion before equation (5.2.13)), and in fact, we have the validity of (5.2.15). The case where $G \in \widetilde{\mathcal{E}}_{\simp}(N)$ is necessarily more subtle. However it is not hard to establish this from the induction hypothesis for a class of ``degenerate" parameters:

\begin{lemma}
Suppose that $G \in \widetilde{\mathcal{E}}_{\simp}(N)$, and $\psi \in \Psi(G)$ such that
\begin{eqnarray}
\dim \overline{T}_{\psi} \geq 2.
\end{eqnarray}
Then we have
\[
R_P(w^0,\widetilde{\pi}_M,\psi_M) \equiv 1
\]
for any $w^0 \in W_{\psi}^0$. In other words part (a) of theorem 5.2.1 is valid, and hence the linear form $f_G(\psi,u)$ depends only on the image of $u$ in $\mathcal{S}_{\psi}$. 
\end{lemma}
\begin{proof}
Since $W_{\psi}^0$ is the Weyl group of the pair $(\overline{S}_{\psi}^0,\overline{T}_{\psi})$ it is generated by simple reflections $\{ w^0_{\alpha}\}$; hence it suffices to prove the assertion for $w^0_{\alpha}$. Since $\dim \overline{T}_{\psi} \geq 2$ by hypothesis, the element $w^0_{\alpha}$ centralizes a torus of positive dimension in $\widehat{G}$; the centralizer of this torus in $\widehat{G}$ is a proper $\Gamma_F$-invariant Levi subgroup $\widehat{M}_{\alpha}$ of $\widehat{G}$ containing $\widehat{M}$, and hence dual to a proper Levi subgroup $M_{\alpha}$ of $G$ containing $M$. Denote by $\psi_{\alpha} \in \Psi(M_{\alpha})$ the image of $\psi_M$ under the $L$-embedding $\leftexp{L}{M} \hookrightarrow \leftexp{L}{M_{\alpha}}$. Then $w^0_{\alpha}$ can be identified as a Weyl element of the pair $(\overline{S}^0_{\psi_{\alpha}},\overline{T}_{\psi_{\alpha}})$ (with $\overline{T}_{\psi_{\alpha}}=\overline{T}_{\psi}$), and we also have:
\[
R_P(w^0_{\alpha},\widetilde{\pi}_M,\psi_M) = R_{P \cap M_{\alpha}}(w^0_{\alpha},\widetilde{\pi}_M,\psi_M).
\]  

\noindent Since $M_{\alpha}$ is a proper Levi subgroup of $G$ it is of the form
\begin{eqnarray*}
& & M= G_{E/F}(N_1^{\prime}) \times \cdots \times G_{E/F}(N_r^{\prime}) \times G_- \\
& & G_- = U_{E/F}(N_-), \,\ N_- < N
\end{eqnarray*}
and hence by the induction hypothesis that equation (5.2.9) is valid for the factor $G_-$ of $M$ (we already know its validity for the part involving the general linear factors). Hence we have $R_{P \cap M_{\alpha}}(w^0_{\alpha},\widetilde{\pi}_M,\psi_M) \equiv 1$.
\end{proof}

\begin{proposition}
As before $G$ is either an element of $\widetilde{\mathcal{E}}_{\simp}(N)$ or $\widetilde{G}_{E/F}(N)$, and $\psi \in \Psi(G)$. Suppose that we have
\begin{eqnarray}
\dim \overline{T}_{\psi,x} \geq 1
\end{eqnarray}
for all $x \in \mathcal{S}_{\psi}$, and if $G \in \widetilde{\mathcal{E}}_{\simp}(N)$ assumes also that
\begin{eqnarray}
\dim \overline{T}_{\psi} \geq 2
\end{eqnarray}
Then theorem 5.2.1 holds for $\psi$. In particular the global intertwining relation is valid for $\psi$.
\end{proposition}
\begin{proof}
We have already seen (remark 5.3.2) that condition (5.3.4) implies that the linear form $f^{\prime}_G(\psi,s)$ is defined for any $s \in \overline{S}_{\psi}$ and depends only on its image in $\mathcal{S}_{\psi}$. In the case where $G \in \widetilde{\mathcal{E}}_{\simp}(N)$ condition (5.3.5) also implies (5.2.9) and that the linear form $f_G(\psi,u)$ depends only on the image of $u \in \mathfrak{N}_{\psi}$ in $\mathcal{S}_{\psi}$ (the case for $G=\widetilde{G}_{E/F}(N)$ holds without this condition). Hence for $u \in \mathfrak{N}_{\psi}$ and $s \in \overline{S}_{\psi}$ having common image $x \in \mathcal{S}_{\psi}$, we have (with the notation in the proof of lemma 5.3.1):
\begin{eqnarray}
& & f_G(\psi,u) = f_G(\psi,x) = f_{M_x}(\psi_{M_x},x) \\
& & f^{\prime}_G(\psi,s_{\psi }s)=f^{\prime}_G(\psi,s_{\psi} x) = f^{\prime}_{M_x}(\psi_{M_x},s_{\psi}x). \nonumber
\end{eqnarray}  
Since (5.3.4) implies that $M_x$ is proper in $G$, we can apply the induction hypothesis, together with proposition 5.2.3 to deal with the general linear factors in $M_x^0$, to conclude that the global intertwining relation $f_{M_x}(\psi_{M_x},x)=f^{\prime}_{M_x}(\psi_{M_x},s_{\psi}x)$ holds. We then conclude from (5.3.6).
\end{proof}

\begin{rem}
\end{rem}
Similarly suppose that $S_{\psi}(G)$ has a central torus of positive dimension. Then the same descent argument as in the above proof shows that theorem 5.2.1 holds for $\psi$.

\subsection{The spectral expansion, part I}

Following Arthur, the discrete part of the trace formula can be given two expansions: the spectral and endoscopic expansions. The global intertwining relation allows the term by term comparison of these two expansions. In the current and the next subsection we consider the spectral expansion.

As before $G$ will denote both an element in $\widetilde{\mathcal{E}}_{\simp}(N)$ and $\widetilde{G}_{E/F}(N)$. Given a parameter $\psi^N \in \widetilde{\Psi}(N)$, we denote
\begin{eqnarray}
r^G_{\disc,\psi^N} := \frac{1}{|\kappa_G|}\tr R^G_{\disc,\psi^N}(f)
\end{eqnarray}

\noindent ({\it c.f.} (4.2.5) for the definition of $\kappa_G$). Note that in the case where $G=(G,\xi) \in \widetilde{\mathcal{E}}_{\simp}(N)$, then we are abbreviating the notation $R^G_{\disc,\psi^N,\xi}$ of (4.3.15) as $R^G_{\disc,\psi_N}$ here. Similar other notational abbreviations will be made below. The global classifications theorem asserts in particular that $\psi^N$ does not contribute to $R^G_{\disc}$ unless $\psi^N \in \xi_* \Psi(G)$. Considerable work needs to be done before we can prove this, and so we have to work with a general parameter $\psi^N$ in the present setting. 

Till the end of section five, we assume either that $\psi^N$ is not a simple generic parameter, so that from the induction hypothesis the seed theorems 2.4.2 and 2.4.10 are valid for the simple generic constituents of $\psi^N$, or that $\psi^N$ is simple generic and the seed theorems 2.4.2 and 2.4.10 are established already for $\psi^N$.

\noindent We then have
\begin{eqnarray}
& & I^G_{\disc,\psi^N}(f) - r^G_{\disc,\psi^N}(f) \\
&=&  \sum_{\substack{ \{M\} \\ M \neq G^0}} \sum_{w \in W(M)_{\reg}} \frac{1}{|W(M)|} \frac{1}{|\det(w-1)_{\mathfrak{a}_M^G}|} \tr (M_{P,\psi^N}(w) \mathcal{I}_{P,\psi^N}(f) ).             \nonumber
\end{eqnarray}

\noindent We thus need an expansion for the term:
\[
\tr (M_{P,\psi^N}(w) \mathcal{I}_{P,\psi^N}(f) )            
\]
for $M \neq G^0$, $w \in W(M)_{\reg}$, and $P \in \mathcal{P}(M)$ fixed.

As before $M$ is a product of groups $G_{E/F}(N^{\prime})$ for $N^{\prime} < N$, together with $G_-=U_{E/F}(N_-)$ for $N_- < N$ in the case $G \in \widetilde{\mathcal{E}}_{\simp}(N)$. Thus from the induction hypothesis the local and global classification theorems hold for $M$. In the case $G=(G,\xi) \in \widetilde{\mathcal{E}}_{\simp}(N)$ we similarly also denote by $M$ as an endoscopic datum in $\widetilde{\mathcal{E}}(N)$, with the $L$-embedding given  by $\xi$ (or more precisely the composition of the $L$-embedding $\leftexp{L}{M} \rightarrow \leftexp{L}{G}$ with $\xi$). 

Denote by $\Psi_2(M,\psi^N)$ the subset of parameters $\psi_M \in \Psi_2(M) \,\  (=\Psi_2(M,\xi))$ such that $\xi_* \psi_M =\psi^N$. Then for $\psi_M \in \Psi_2(M,\psi^N)$, we have the subspace:
\begin{eqnarray*}
L^2_{\disc, \psi_M}(M(F) \backslash M(\mathbf{A}_F)) \subset L^2_{\disc}(M(F) A_{M,\infty}^+\backslash M(\mathbf{A}_F))
\end{eqnarray*}

By the global classification theorem for $M$ we have a decomposition:
\begin{eqnarray*}
L^2_{\disc, \psi_M}(M(F) \backslash M(\mathbf{A}_F)) = \bigoplus_{\pi_M \in \Pi_{\psi_M}} m(\pi_M) \pi_M 
\end{eqnarray*} 
and hence

\begin{eqnarray}
\mathcal{I}_{P,\psi^N} = \bigoplus_{  \psi_M \in \Psi_2(M,\psi^N) } \bigoplus_{\pi_M \in \Pi_{\psi_M}} m(\pi_M) \mathcal{I}_P(\pi_M).
\end{eqnarray}
In the case where $G=(G,\xi) \in \widetilde{\mathcal{E}}_{\simp}(N)$, note that in order for $\Psi_2(M,\psi^N)$ to be non-empty, we must have $\psi^N \in \xi_* \Psi(G)$, i.e. $\psi^N$ defines a parameter $\psi \in \Psi(G)$ such that $\psi^N = \xi_* \psi$ (if $G=\widetilde{G}_{E/F}(N)$ then we interpret $\psi$ just as another name for $\psi^N$). In particular, if $\psi^N \notin \xi_* \Psi(G)$, then we interpret all the terms below that involve $\psi$ as being empty.

\noindent Hence denoting by $M_P(w,\pi_M)$ the restriction of $M_{P,\psi^N}(w)$ to $\mathcal{I}_P(\pi_M)$, we can write
\begin{eqnarray}
& & \tr(M_{P,\psi^N}(w) \mathcal{I}_{P,\psi^N}(f)) \\
&=&  \sum_{ \psi_M \in \Psi_2(M,\psi^N) } \sum_{\pi_M \in \Pi_{\psi_M}}  m(\pi_M) \tr(M_P(w,\pi_M) \mathcal{I}_P(\pi_M,f)). \nonumber
\end{eqnarray}

\noindent In the summation we can limit to the set of parameters $\psi_M \in \Psi_2(M,\psi^N)$ that satisfies $w \psi_M = \psi_M$, i.e. the set of parameters of $\Psi_2(M,\psi^N)$ that extends to the twisted group $\widetilde{M}_w$ (defined as in (3.4.4)), and we denote this set as $\Psi_2(\widetilde{M}_w,\psi^N)$. In particular the sum is empty unless $w \in W_{\psi}$. 

\noindent We can apply the spectral multiplicity formula for $\pi_M$ (with respect to $M$):
\begin{eqnarray*}
m(\pi_M)= \frac{1}{|\mathcal{S}_{\psi_M}|} \sum_{x_M \in \mathcal{S}_{\psi_M}} \epsilon^M_{\psi_M}(x_M) \langle x_M  , \pi_M  \rangle
\end{eqnarray*}

\noindent Now for $u \in \mathfrak{N}_{\psi}(w)$ denote by $\epsilon^M_{\psi_M}(\widetilde{u})$ the canonical extension of $\epsilon^M_{\psi_M}$ to $\widetilde{\mathcal{S}}_{\psi_M,u}$ (evaluated at the element $\widetilde{u}$). The canonical extension $\epsilon^M_{\psi_M}(\widetilde{u})$ is defined by the same considerations as in the discussion before remark 3.4.1, using the fact that the general linear factors of $\widetilde{M}_w$ has no contribution to $\epsilon^M_{\psi_M}$. Then we can write (5.4.4) as

\begin{eqnarray}
\end{eqnarray}
\begin{eqnarray}
& & \tr(M_{P,\psi^N}(w)  \mathcal{I}_{P,\psi^N}(f)) \nonumber \\
&=& \sum_{\psi_M \in \Psi_2(M_w,\psi^N)} \sum_{\pi_M \in \Pi_{\psi_M}} \sum_{x_M \in \mathcal{S}_{\psi_M}} \frac{1}{|\mathcal{S}_{\psi_M}|} \epsilon^M_{\psi_M}(x_M) \langle x_M, \pi_M \rangle \tr(M_P(w,\pi_M) \mathcal{I}_{P}(\pi_M,f)   )  \nonumber \\
&=&  \sum_{\psi_M \in \Psi_2(M_w,\psi^N)} \sum_{\pi_M \in \Pi_{\psi_M}} \sum_{u \in \mathfrak{N}_{\psi}(w)} \frac{1}{|\mathcal{S}_{\psi_M}|} \epsilon^M_{\psi_M}(\widetilde{u}) \langle \widetilde{u},\widetilde{\pi}_M \rangle \tr( M_P(w,\pi_M)  \mathcal{I}_P(\pi_M,f) ).       \nonumber
\end{eqnarray}

\noindent For the last equality we have used:
\begin{eqnarray}
\sum_{x_M \in \mathcal{S}_{\psi_M}} \epsilon^M_{\psi_M}(x_M) \langle x_M,\pi_M \rangle = \sum_{u \in \mathfrak{N}_{\psi}(w)} \epsilon^M_{\psi_M}(\widetilde{u}) \langle \widetilde{u}, \widetilde{\pi}_M\rangle 
\end{eqnarray}
for any $w \in W_{\psi}$. That (5.4.6) is valid follows from the definition of the canonical extension. In the particular case $G=\widetilde{G}_{E/F}(N)$ then $\mathcal{S}_{\psi_M}$ is of course trivial, both $\epsilon^M_{\psi_M}(x_M)$ and $\epsilon^M_{\psi_M}(\widetilde{u})$ are equal to one, and that the pairing $\langle x_M,\pi_M \rangle$ and $\langle \widetilde{u},\widetilde{\pi}_M \rangle$ are trivial.

\bigskip

\noindent We also have the equality:
\begin{eqnarray}
M_P(w,\pi_M)= r_P(w,\psi_M) R_P(w,\widetilde{\pi}_M,\psi_M)
\end{eqnarray}
here $r_P(w,\psi_M)$ is the product over all $v$ of the local normalizing factors:
\begin{eqnarray}
r_P(w,\psi_M) = \prod_v r_{P_v} (w,\psi_{M,v}).
\end{eqnarray}
Indeed in the case where $G \in \widetilde{\mathcal{E}}_{\simp}(N)$ (5.4.7) follows from the way the normalized intertwining operator $R_P(w,\widetilde{\pi}_M,\psi_M)$ is defined. In the case where $G= \widetilde{G}_{E/F}(N)$, one must in addition also use the analogue of lemma 4.2.3 of \cite{A1}, in our context of the twisted group $\widetilde{G}_{E/F}(N)=(G_{E/F}(N),\theta)$, in order to see that the extension of the global representation $R^N_{\disc}$ to the twised group $\widetilde{G}_{E/F}(N)(\mathbf{A}_F)$, defined by the automorphism $\theta$, is compatible with the product over all the places of $F$ of the local extensions defined as in section 3.2 in terms of local Whittaker models.

We note for future reference the following formula for the global normalizing factor
\begin{eqnarray*}
& & r_P(w,\psi_M) = \prod_v r_{P_v}(w,\psi_{M,v}) \\
&=& \prod_v \lambda_v(w) \epsilon_v(0, \pi_{\xi_* \psi_M},  \widetilde{\rho}^{\vee}_{w^{-1}P|P})^{-1} L_v(0, \pi_{\xi_* \psi_M}, \widetilde{\rho}^{\vee}_{w^{-1}P|P})    L_v(1, \pi_{\xi_* \psi_M}, \widetilde{\rho}^{\vee}_{w^{-1}P|P})^{-1}   \\
&=& \epsilon(0, \pi_{\xi_* \psi_M},  \widetilde{\rho}^{\vee}_{w^{-1}P|P})^{-1} L(0, \pi_{\xi_* \psi_M}, \widetilde{\rho}^{\vee}_{w^{-1}P|P})    L(1, \pi_{\xi_* \psi_M}, \widetilde{\rho}^{\vee}_{w^{-1}P|P})^{-1}   
\end{eqnarray*}
with the notation being the global analogue of that of section 3.3 ({\it c.f.} (3.3.44)). Here we are using the fact that the global $\lambda(w)=\prod_v \lambda_v(w)$ factor is equal to one. 

To conclude the above analysis, we have the formula:
\begin{eqnarray}
& & \tr(M_{P,\psi^N}(w) \mathcal{I}_{P,\psi^N}(f))  \\
&=& \sum_{\psi_M \in \Psi_2(M_w,\psi^N)} \sum_{\pi_M \in \Pi_{\psi_M}} \sum_{u \in \mathfrak{N}_{\psi}(w)}  \nonumber \\
& & \,\ \,\ \,\ \,\ \,\  \frac{1}{|\mathcal{S}_{\psi_M}|} r_P(w,\psi_M) \epsilon^M_{\psi_M}(\widetilde{u}) \langle \widetilde{u},\widetilde{\pi}_M \rangle \tr(R_P(w,\widetilde{\pi}_M,\psi_M) \mathcal{I}_P(\pi_M,f) ). \nonumber
\end{eqnarray}

\subsection{The spectral expansion, part II}

We now analyze the terms on the right hand side of (5.4.2). By applying (5.4.9), we see that (5.4.2) is a five-fold sum:
\begin{eqnarray}
\sum_{\substack{\{M\}  \\ M \neq G^0}} \sum_{w \in W(M)_{\reg}} \sum_{\psi_M \in \Psi_2(\widetilde{M}_w,\psi^N)} \sum_{u \in \mathfrak{N}_{\psi}(w)} \sum_{\pi_M \in \Pi_{\psi_M}} 
\end{eqnarray}
of the summand
\begin{eqnarray}
& & \frac{1}{|W(M)| \cdot |\mathcal{S}_{\psi_M}| \cdot |\det(w-1)_{\mathfrak{a}_M^G}| } r_P(w,\psi_M) \epsilon^M_{\psi_M}(\widetilde{u})   \\
& & \,\ \,\ \,\ \,\ \,\ \times \langle \widetilde{u}, \widetilde{\pi}_M\rangle \tr(R_P(w,\widetilde{\pi}_M,\psi_M) \mathcal{I}_P(\pi_M,f)  ) \nonumber
\end{eqnarray}

\noindent which can be written as the fourfold sum
\begin{eqnarray}
\sum_{\substack{\{M\}  \\ M \neq G^0}} \sum_{w \in W(M)_{\reg}} \sum_{\psi_M \in \Psi_2(\widetilde{M}_w,\psi^N)} \sum_{u \in \mathfrak{N}_{\psi}(w)} 
\end{eqnarray}
of the summand
\begin{eqnarray}
& & \\
& &  \frac{1}{|W(M)| \cdot |\mathcal{S}_{\psi_M}| \cdot |\det(w-1)_{\mathfrak{a}_M^G}| } r_P(w,\psi_M) \epsilon^M_{\psi_M}(\widetilde{u})  f_G(\psi,u). \nonumber
\end{eqnarray}

\noindent The double sum in (5.5.3)
\begin{eqnarray}
\sum_{w \in W(M)_{\reg}} \sum_{\psi_M \in \Psi_2(\widetilde{M}_w,\psi^N)}
\end{eqnarray}
can be replaced by a simple sum over the set
\begin{eqnarray}
V_{\psi} =\{ (\psi_M, w) \in \Psi_2(M,\psi^N) \times W(M)_{\reg}| \,\  w \psi_M = \psi_M          \}.
\end{eqnarray}

\noindent Then the projection $(\psi_M,w) \rightarrow w$ gives a fibration:
\begin{eqnarray}
V_{\psi} \rightarrow \{W_{\psi,\reg} \}
\end{eqnarray}
where $\{W_{\psi,\reg}\}$ denotes the set of $W_{\psi}$-conjugacy classes in $W_{\psi,\reg}$. 

We have a natural conjugation action of $W(M)$ on the set $V_{\psi}$, and $W(M)$ acts transitively on the fibres of the map (5.5.7). The stabilizer in $W(M)$ of any parameter in $\Psi_2(M,\psi^N)$ is isomorphic to $W_{\psi}$. Hence we have
\begin{eqnarray}
|\Psi_2(M,\psi^N)| =m_{\psi^N}^G |W(M)| |W_{\psi} |^{-1}
\end{eqnarray}
where we recall that $m^G_{\psi^N}$ is defined as in (5.1.10) (if $G=(G,\xi) \in \widetilde{\mathcal{E}}_{\simp}(N)$ and $\psi^N \notin \xi_* \Psi(G)$ then $m_{\psi^N}^G=0$, and the right hand side of (5.5.8) is interpreted just as being equal to zero, which is consisitent with $\Psi_2(M,\psi^N)$ being empty). 

\bigskip

\noindent Furthermore the summand (5.5.4) does not depend on the choice of $\psi_M \in \Psi_2(M,\psi^N)$, hence the triple sum in (5.5.3) 
\begin{eqnarray}
\sum_{w \in W(M)_{\reg}} \sum_{\psi_M \in \Psi_2(\widetilde{M}_w,\psi^N)} \sum_{u \in \mathfrak{N}_{\psi}(w)}
\end{eqnarray}
reduces to a sum
\begin{eqnarray}
\sum_{c \in \{ W_{\psi,\reg} \}} \sum_{w \in c} \sum_{u \in \mathfrak{N}_{\psi}(w)} = \sum_{w \in W_{\psi,\reg} } \sum_{u \in \mathfrak{N}_{\psi}}
\end{eqnarray}
of the summand

\begin{eqnarray}
& &   \frac{|W(M)|}{|W_{\psi}|}  \frac{m_{\psi^N}^G}{|W(M)| \cdot |\mathcal{S}_{\psi_M}| \cdot |\det(w-1)_{\mathfrak{a}_M^G}| }  \\ 
& & \,\ \,\ \,\ \,\ \times r_P(w,\psi_M) \epsilon^M_{\psi_M}(\widetilde{u}) f_G(\psi,u) \nonumber \\
&= &   \frac{m_{\psi^N}^G}{|\mathcal{S}_{\psi}|}  \frac{1}{ \cdot |W^0_{\psi}| \cdot |\det(w-1)_{\mathfrak{a}_M^G}| }  \nonumber \\ 
& & \,\ \,\ \,\ \,\ \times r_P(w,\psi_M) \epsilon^M_{\psi_M}(\widetilde{u}) f_G(\psi,u) \nonumber
\end{eqnarray} 
for a fixed choice of $\psi_M \in \Psi_2(M,\psi^N)$. For the last equality we have used the fact that $|W_{\psi}| |\mathcal{S}_{\psi_M}|=|\mathfrak{N}_{\psi}|=|\mathcal{S}_{\psi}| |W_{\psi}^0|$ which follows from the commutative diagram (5.2.2).

Now the $G^0$-conjugacy class of $M$ is determined by the parameter $\psi_M$. We only need a marker to rule out the case $M=G^0$ to account for the summation 
\[
\sum_{\substack{\{M\} \\ M \neq G^0}}
\]
in (5.5.3). This can be done by introducing
\begin{eqnarray}
& & \\
W_{\psi,\reg}^{\prime} &=& W_{\psi,\reg} \smallsetminus \{1\} \nonumber \\
&=& \left \{ \begin{array}{c}  \emptyset \mbox{ if } \overline{S}_{\psi} \mbox{ is finite, i.e. } \psi \in \Psi_2(G) \\   \,\  W_{\psi,\reg} \mbox{ otherwise. } 
\end{array} \right. \nonumber
\end{eqnarray}

Denote by $\mathfrak{N}_{\psi,\reg}$ and $\mathfrak{N}_{\psi,\reg}^{\prime}$ the inverse image of $W_{\psi,\reg}$ and $W_{\psi,\reg}^{\prime}$ in $\mathfrak{N}_{\psi}$. At this point we need the crucial {\it spectral sign lemma}. Before stating this, we introduce one piece of terminology.

We called a parameter $\psi^N \in \widetilde{\Psi}(N)$ to be an $\epsilon$-{\it parameter}, if it is of the form
\begin{eqnarray*}
& & \psi^N=\psi_1^{N_1} \boxplus \psi_2^{N_2} \\
& & \psi_i^{N_i} \in \widetilde{\Psi}_{\simp}(N_i), \,\  N_i < N,\,\  N_1 + N_2 =N \\
& & \psi_i^{N_i} = \mu_i \boxtimes \nu_i
\end{eqnarray*}
such that $\mu_1$ and $\mu_2$ are of the {\it same parity}, i.e. $\mu_1$ and $\mu_2$ are either both conjugate orthogonal or both conjugate symplectic (in the sense of part (a) of theorem 2.5.4), and that $\nu_1 \otimes \nu_2$ is a direct sum of {\it odd} number of irreducible representations of $\SL_2(\mathbf{C})$ of {\it even} dimension (in particular $\nu_1$ and $\nu_2$ are of different parity). Thus in particular an $\epsilon$-parameter is an elliptic parameter in $\widetilde{\Psi}_{\ellip}(N)$, and $\psi^N \in \xi^{\prime}_* \Psi_2(G^{\prime})$ for $G^{\prime}=(G^{\prime},\xi^{\prime}) \in \widetilde{\mathcal{E}}_{\ellip}(N) \smallsetminus \widetilde{\mathcal{E}}_{\simp}(N)$.

\begin{lemma} (The spectral sign lemma)

\noindent Suppose that $G \in \widetilde{\mathcal{E}}_{\simp}(N)$. For any $u \in \mathfrak{N}_{\psi}$, we have the identity
\begin{eqnarray}
r_P(w_u,\psi_M) \epsilon^M_{\psi_M}(\widetilde{u}) = \sgn^0(w_u) \epsilon^G_{\psi}(x_u).
\end{eqnarray}

\noindent If $G=\widetilde{G}_{E/F}(N)$ and $\psi^N \in \widetilde{\Psi}(N)$ is not an $\epsilon$-parameter then (5.5.13) again holds; if $\psi^N$ is an $\epsilon$-parameter then (5.5.13) is equivalent to the assertion $\epsilon(1/2, \mu_1 \times \mu_2^c)=1$ (as asserted by part (b) of theorem 2.5.4). 
\end{lemma}
We will establish the spectral sign lemma in section 5.8, and we will take it for granted at the present moment. Here we just note that all the quantities on the left hand side of (5.5.13) depend only on $\psi$, and thus we denote 
\begin{eqnarray*}
& & r^G_{\psi}(w) := r_P(w,\psi_M), \,\ w \in W_{\psi} \\
& & \epsilon^1_{\psi}(u) := \epsilon^M_{\psi_M}(\widetilde{u}), \,\ u \in \mathfrak{N}_{\psi}.
\end{eqnarray*}

In the rest of this subsection, we assume, in the case $G=\widetilde{G}_{E/F}(N)$, that either $\psi^N$ is not an $\epsilon$-parameter, or that (5.5.13) is valid for $\psi^N$. Then the sum (5.5.3) becomes:
\begin{eqnarray}
& & \frac{m_{\psi^N}^G}{|\mathcal{S}_{\psi}|}\sum_{u \in \mathfrak{N}_{\psi,\reg}^{\prime}} \epsilon^G_{\psi}(x_u) \frac{1}{ |W_{\psi}^0 |} \frac{\sgn^0(w_u)}{|\det(w_u-1)_{\mathfrak{a}_M^G}|} f_G(\psi,u) \\
&=& |\kappa_G|^{-1} \frac{m_{\psi^N}^G}{|\mathcal{S}_{\psi}|}\sum_{u \in \mathfrak{N}_{\psi,\reg}^{\prime}} \epsilon^G_{\psi}(x_u) \frac{1}{ |W_{\psi}^0 |} \frac{\sgn^0(w_u)}{|\det(w_u-1)_{\mathfrak{a}_{\overline{T}_{\psi}}}|} f_G(\psi,u) \nonumber
\end{eqnarray}
here for the last equality we have used (for any $w \in W_{\psi}$):
\begin{eqnarray*}
& & |\det(w-1)_{\mathfrak{a}^G_M}|=|\det(w-1)_{\mathfrak{a}^G_{G^0}}| \cdot |\det(w-1)_{\mathfrak{a}_M^{G^0}}| \\
&=& |\det(\theta-1)_{\mathfrak{a}^G_{G^0}}| \cdot |\det(w-1)_{\mathfrak{a}_M^{G^0}}| \\
&=& |\kappa_G| \cdot |\det(w-1)_{\mathfrak{a}_M^{G^0}}|
\end{eqnarray*}
and the identification $\mathfrak{a}_M^{G^0} \cong \mathfrak{a}^*_{\overline{T}_{\psi}}$ (which is equivariant with respect to the identification $W(M) \cong W(\widehat{M})$).

\bigskip

We can then fibre the sum (5.5.14) over $\mathcal{S}_{\psi}$ (again with respect to the commutative diagram (5.2.2)): for $x \in \mathcal{S}_{\psi}$ define $\mathfrak{N}_{\psi}(x)$ to be the inverse image of $x$ in $\mathfrak{N}_{\psi}$, and put
\begin{eqnarray*}
& & \mathfrak{N}_{\psi,\reg}(x) = \mathfrak{N}_{\psi,\reg} \cap \mathfrak{N}_{\psi}(x) \\
& & \mathfrak{N}_{\psi,\reg}^{\prime}(x) = \mathfrak{N}^{\prime}_{\psi,\reg} \cap \mathfrak{N}_{\psi}(x)
\end{eqnarray*}
and let $W_{\psi,\reg}(x)$ and $W_{\psi,\reg}^{\prime}(x)$ be the bijective image of $\mathfrak{N}_{\psi,\reg}(x)$ and $\mathfrak{N}_{\psi,\reg}^{\prime}(x)$ under the map $\mathfrak{N}_{\psi} \rightarrow W_{\psi}$, then the sum (5.5.14), which is equal to (5.4.2), can be written as
\begin{eqnarray}
& &  I^G_{\disc, \psi^N}(f) - r_{\disc,\psi^N}^G(f) \\
&= & |\kappa_G|^{-1} \frac{m_{\psi^N}^G}{|\mathcal{S}_{\psi}|} \sum_{x \in \mathcal{S}_{\psi}}\sum_{u \in \mathfrak{N}_{\psi,\reg}^{\prime}(x)} \epsilon^G_{\psi}(x_u) \frac{1}{ |W_{\psi}^0 |} \frac{\sgn^0(w_u)}{|\det(w_u-1)_{\mathfrak{a}_{\overline{T}_{\psi}}}|} f_G(\psi,u). \nonumber
\end{eqnarray}

\bigskip

As noted above in (5.5.12) if $\psi \notin \Psi_2(G)$ then $\mathfrak{N}_{\reg}^{\prime}=\mathfrak{N}_{\reg}$. On the other hand if $\psi \in \Psi_2(G)$ then $\mathfrak{N}_{\reg}^{\prime} = \emptyset$, so that the sum on the right hand side of (5.5.15) is empty. In this case, we have $W_{\psi}^0$ being trivial and $\mathfrak{N}_{\psi}=\mathcal{S}_{\psi}$. The existence of the linear form $f_G(\psi,u)=f_G(\psi,x)$ ($x \in \mathcal{S}_{\psi}$) amounts to the existence of the packet $\Pi_{\psi}$ and the associated pairing $\langle \cdot, \cdot \rangle$ on $\mathcal{S}_{\psi} \times \Pi_{\psi}$, in which case the linear form $f_G(\psi,x)$ is given by:
\[
f_G(\psi,x) = \sum_{\pi \in \Pi_{\psi}} \langle x , \pi \rangle f_G(\pi).
\]

\bigskip

We can give a uniform treatment as follows: if $\psi \in \Psi_2(G)$ and if the linear form $f_G(\psi,x)$ is defined for $x \in \mathcal{S}_{\psi}$, put:
\begin{eqnarray}
& & \\
& & \leftexp{0}{r^G_{\disc,\psi^N}(f)} :=  r^G_{\disc,\psi^N}(f) -  |\kappa_G|^{-1} \frac{m_{\psi^N}^G}{|\mathcal{S}_{\psi}|} \sum_{x \in \mathcal{S}_{\psi}} \epsilon^G_{\psi}(x)  f_G(\psi,x) \nonumber \\
&=& r^G_{\disc,\psi^N}(f) -  |\kappa_G|^{-1} \frac{m_{\psi^N}^G}{|\mathcal{S}_{\psi}|} \sum_{x \in \mathcal{S}_{\psi}} \sum_{\pi \in \Pi_{\psi}} \epsilon^G_{\psi}(x) \langle x,\pi \rangle f_G(\pi) \nonumber \\
&=& r^G_{\disc,\psi^N}(f) - |\kappa_G|^{-1} \sum_{\pi \in \Pi_{\psi}} m(\pi) f_G(\pi) \nonumber
\end{eqnarray}
here
\[
m(\pi):=\frac{m_{\psi^N}^G}{|\mathcal{S}_{\psi}|} \sum_{x \in \mathcal{S}_{\psi}} \epsilon^G_{\psi}(x) \langle x,\pi \rangle, 
\]
and on the other hand, if $\psi \notin \Psi_2(G)$, then we just put
\begin{eqnarray}
\leftexp{0}{r^G_{\disc,\psi^N}(f)} := r^G_{\disc,\psi^N}(f)
\end{eqnarray}

\noindent then in all cases we have
\begin{eqnarray}
& &  I^G_{\disc, \psi^N}(f) - \leftexp{0}{r_{\disc,\psi^N}^G(f)} \\
&= &  |\kappa_G|^{-1} \frac{m_{\psi^N}^G}{|\mathcal{S}_{\psi}|} \sum_{x \in \mathcal{S}_{\psi}}\sum_{u \in \mathfrak{N}_{\psi,\reg}(x)} \epsilon^G_{\psi}(x_u) \frac{1}{ |W_{\psi}^0 |} \frac{\sgn^0(w_u)}{|\det(w_u-1)_{\mathfrak{a}_{\overline{T}_{\psi}}}|} f_G(\psi,u) \nonumber
\end{eqnarray}
(the symbol $\leftexp{0}{r^G_{\disc,\psi^N}}$ is of course for the reason that the global classification theorem asserts that it is identifically equal to zero).

If in addition we happen to know that the linear form $f_G(\psi,u)$ depends only on the image of $u \in \mathfrak{N}_{\psi}$ in $\mathcal{S}_{\psi}$, the the sum (5.5.18) further simplifies. We record this as
\begin{proposition}
Suppose that $f_G(\psi,u)$ is defined for all $u \in \mathfrak{N}_{\psi,\reg}$, then we have:
\begin{eqnarray}
& &  I^G_{\disc, \psi^N}(f) - \leftexp{0}{r_{\disc,\psi^N}^G(f)} \\
&= &  |\kappa_G|^{-1} \frac{m_{\psi^N}^G}{|\mathcal{S}_{\psi}|} \sum_{x \in \mathcal{S}_{\psi}}\sum_{u \in \mathfrak{N}_{\psi,\reg}(x)} \epsilon^G_{\psi}(x_u) \frac{1}{ |W_{\psi}^0 |} \frac{\sgn^0(w_u)}{|\det(w_u-1)_{\mathfrak{a}_{\overline{T}_{\psi}}}|} f_G(\psi,u). \nonumber
\end{eqnarray}
If in addition, the linear form $f_G(\psi,u)=f_G(\psi,x)$ depends only on the image $x$ of $u$ in $\mathcal{S}_{\psi}$, then we have
\begin{eqnarray}
& & \\
& &  I^G_{\disc, \psi^N}(f) - \leftexp{0}{r_{\disc,\psi^N}^G(f)} = |\kappa_G|^{-1} \frac{m_{\psi^N}^G}{|\mathcal{S}_{\psi}|} \sum_{x \in \mathcal{S}_{\psi}} i_{\psi}(x) \epsilon^G_{\psi}(x)  f_G(\psi,x) \nonumber
\end{eqnarray}
here
\begin{eqnarray}
i_{\psi}(x) := \sum_{w \in W_{\psi,\reg}(x)}  \frac{1}{ |W_{\psi}^0 |} \frac{\sgn^0(w)}{|\det(w-1)_{\mathfrak{a}_{\overline{T}_{\psi}}}|}
\end{eqnarray}
is the number as defined in (5.1.2) attached to the component of $\overline{S}_{\psi}$ indexed by $x$.
\end{proposition}

\subsection{The endoscopic expansion}

In this subsection we consider the endoscopic expansion, and is parallel to the two previous subsections. 

As in the last subsection we will denote by $G$ either an element of $\widetilde{\mathcal{E}}_{\simp}(N)$ or $\widetilde{G}_{E/F}(N)$. Again fix a parameter $\psi^N \in \widetilde{\Psi}(N)$. Parallel to the term $r^G_{\disc,\psi^N}$ in the spectral expansion, we form the following sum:

\begin{eqnarray}
& & \\
& & s^G_{\disc,\psi^N}(f) = \sum_{G^{\prime} \in \mathcal{E}_{\simp}(G) } \iota(G,G^{\prime}) \widehat{S}^{G^{\prime}}_{\disc,\psi^N}(f^{G^{\prime}}), \,\ f \in \mathcal{H}(G). \nonumber
\end{eqnarray}

\noindent In the case where $G \in \widetilde{\mathcal{E}}_{\simp}(N)$ we of course have $s^G_{\disc,\psi^N}=S^G_{\disc,\psi^N}$.

We would like to refine the endoscopic expansion
\begin{eqnarray}
& & I^{G}_{\disc,\psi^N}(f) - s^G_{\disc,\psi^N}(f) \\
&=&  \sum_{G^{\prime} \in \widetilde{\mathcal{E}}_{\ellip}(G) \smallsetminus \widetilde{\mathcal{E}}_{\simp}(G) } \iota(G,G^{\prime}) \widehat{S}^{G^{\prime}}_{\disc,\psi^N}(f^{G^{\prime}})  \nonumber
\end{eqnarray}
that can be compared term by term with the spectral expansion (5.5.19) or (5.5.20).

For each $G^{\prime} \in \mathcal{E}_{\ellip}(G) \smallsetminus \mathcal{E}_{\simp}(G)$, denote by $\Psi(G^{\prime},\psi^N)$ the set of parameters $\psi^{\prime} \in \Psi(G^{\prime})$ such that $\psi^{\prime}$ maps to $\psi^N$ under the embeddings $\leftexp{L}{G^{\prime}} \hookrightarrow \leftexp{L}{G^0} \hookrightarrow \leftexp{L}{G_{E/F}(N)}$ (with the $L$-embeddings as being part of the endoscopic data; if $G=\widetilde{G}_{E/F}(N)$ then $\leftexp{L}{G^0} \rightarrow \leftexp{L}{G_{E/F}(N)}$ is of course just the identity map). Again in the case $G=(G,\xi) \in \widetilde{\mathcal{E}}_{\simp}(N)$, $\Psi(G^{\prime},\psi^N)$ is empty unless $\psi^N \in \xi_* \Psi(G)$, and as in the last subsection we denote by $\psi \in \Psi(G)$ the parameter in $\Psi(G)$ defined by $\psi^N$.  

For $G^{\prime} = (G^{\prime},\zeta^{\prime}) \in \mathcal{E}_{\ellip}(G) \smallsetminus \mathcal{E}_{\simp}(G)$, we see that $G^{\prime}$ is a proper product:
\[
G^{\prime} = G^{\prime}_1 \times G^{\prime}_2, \,\ G^{\prime}_i =(G^{\prime}_i,\xi^{\prime}_i) \in \widetilde{\mathcal{E}}_{\simp}(N_i), \,\ N_i <N
\]
(in other words $\xi \circ \zeta^{\prime} = \xi_1^{\prime} \times \xi_2^{\prime}$, where $\xi: \leftexp{L}{G} \rightarrow \leftexp{L}{G_{E/F}(N)}$ is the $L$-embedding that is part of the endoscopic data of $G$ in case $G \in \widetilde{\mathcal{E}}_{\simp}(N)$; in case $G=\widetilde{G}_{E/F}(N)$ we simply interpret $\xi$ as the identity map $\leftexp{L}{\widetilde{G}_{E/F}(N)^0} = \leftexp{L}{G_{E/F}(N)} \rightarrow \leftexp{L}{G_{E/F}(N)}$). 

\bigskip

\noindent Hence from the induction hypothesis we have the validity of the stable multiplicity formula for each simple factor $G^{\prime}_i$. Hence we also have the stable multiplicity formula for $G^{\prime}$. More precisely, suppose that for $i=1,2$ we have parameters $\psi_i^{N_i} \in \widetilde{\Psi}(N_i)$. Denote by $\psi^{\prime}_i \in \Psi(G^{\prime}_i)$ the parameter of $G^{\prime}_i$ defined by $\psi_i^{N_i}$, if $\psi_i^{N_i} \in (\xi^{\prime}_i)_* \Psi(G^{\prime}_i)$ (i.e. $\psi_i^{N_i} = (\xi^{\prime}_i)_*\psi^{\prime}_i$). Then 
\begin{eqnarray*}
S^{G^{\prime}_i}_{\disc,\psi_i^{N_i}}(f^{\prime}_i) = \frac{m_{\psi_i^{N_i}}^{G^{\prime}_i}}{|\mathcal{S}_{\psi^{\prime}_i}|} \epsilon^{G^{\prime}_i}(\psi^{\prime}_i)  \sigma(\overline{S}^0_{\psi^{\prime}_i}) (f^{\prime}_i)^{G^{\prime}_i}(\psi^{\prime}_i), \,\ f^{\prime}_i \in \mathcal{H}(G^{\prime}_i).
\end{eqnarray*}

\noindent Define
\begin{eqnarray*}
& & \Psi(G^{\prime},\psi^N) \\
&=& \{\psi^{\prime} = \psi_1^{\prime} \times \psi_2^{\prime} \in \Psi(G^{\prime})| \,\  \psi_i^{\prime} \in \Psi(G_i^{\prime}) \text{ for } i=1,2, \,\  \psi^N =   \psi_1^{N_1} \boxplus \psi_2^{N_2}   \}.
\end{eqnarray*}

\noindent Then for $f^{\prime} = f^{\prime}_1 \times f^{\prime}_2 \in \mathcal{H}(G^{\prime})$, we have
\begin{eqnarray*}
& & S^{G^{\prime}}_{\disc,\psi^N}(f^{\prime}) =\sum_{\substack{ \psi^{\prime} = \psi^{\prime}_1 \times \psi^{\prime}_2 \\ \psi^{\prime} \in \Psi(G^{\prime},\psi^N) }} S^{G^{\prime}_1}_{\disc,\psi_1^{N_1}} (f^{\prime}_1) \cdot S^{G^{\prime}_2}_{\disc,\psi_2^{N_2}} (f^{\prime}_2) \\
&=& \sum_{\substack{ \psi^{\prime} = \psi^{\prime}_1 \times \psi^{\prime}_2 \\ \psi^{\prime} \in \Psi(G^{\prime},\psi^N) }} \prod_{i=1}^2 \frac{1}{|\mathcal{S}_{\psi^{\prime}_i}|} \epsilon^{G^{\prime}_i}(\psi^{\prime}_i)  \sigma(\overline{S}^0_{\psi^{\prime}_i}) (f^{\prime}_i)^{G^{\prime}_i}(\psi^{\prime}_i) \\
&=&  \sum_{\substack{ \psi^{\prime} = \psi^{\prime}_1 \times \psi^{\prime}_2 \\ \psi^{\prime} \in \Psi(G^{\prime},\psi^N) }} \Big( \prod_{i=1}^2 \frac{1}{|\mathcal{S}_{\psi^{\prime}_i}|} \epsilon^{G^{\prime}_i}(\psi^{\prime}_i)  \sigma(\overline{S}^0_{\psi^{\prime}_i})  \Big) (f^{\prime})^{G^{\prime}}(\psi^{\prime}).
\end{eqnarray*}

\noindent Thus the above formula holds for all $f^{\prime} \in \mathcal{H}(G)$ (i.e. not just those $f^{\prime}$ that are pure product). By the multiplicativity of various quantities we thus see that the stable multiplicity formula in the composite case takes the following form:
\begin{eqnarray}
& & \\
& & S^{G^{\prime}}_{\disc,\psi^N}(f^{\prime}) = \sum_{\psi^{\prime} \in \Psi(G^{\prime},\psi^N)} \frac{1}{|\mathcal{S}_{\psi^{\prime}}|} \epsilon^{\prime}(\psi^{\prime}) \sigma(\overline{S}^0_{\psi^{\prime}}) (f^{\prime})^{G^{\prime}}(\psi^{\prime}), \,\ f^{\prime} \in \mathcal{H}(G^{\prime}) \nonumber
\end{eqnarray}
here we have denoted $\epsilon^{G^{\prime}}(\psi^{\prime})$ as $\epsilon^{\prime}(\psi^{\prime})$. Note in particular that the stable linear form $(f^{\prime})^{G^{\prime}}(\psi^{\prime})$ is defined on $\mathcal{H}(G^{\prime})$.

Hence from (5.6.2) and (5.6.3) we see that the difference
\begin{eqnarray*}
 I^G_{\disc,\psi^N}(f) - s^G_{\disc,\psi^N}(f) 
\end{eqnarray*}
is given by a double sum
\begin{eqnarray}
\sum_{G^{\prime} \in \mathcal{E}_{\ellip}(G) \backslash \mathcal{E}_{\simp}(G)} \sum_{\psi^{\prime} \in \Psi(G^{\prime},\psi^N)} 
\end{eqnarray}
with the summand given by the product of
\begin{eqnarray}
\iota(G,G^{\prime}) = |\kappa_G|^{-1} |\overline{Z}(\widehat{G}^{\prime})^{\Gamma}|^{-1}  |\Out_G(G^{\prime})|^{-1}
\end{eqnarray}
and
\begin{eqnarray}
 \frac{1}{|\mathcal{S}_{\psi^{\prime}}|} \epsilon^{\prime}(\psi^{\prime}) \sigma(\overline{S}^0_{\psi^{\prime}}) f^{\prime}(\psi^{\prime})
\end{eqnarray}
(we have made the abbreviaation $(f^{\prime})^{G^{\prime}}(\psi^{\prime})=f^{\prime}(\psi^{\prime})$).

Note that by the property (5.1.9) we have $\sigma(\overline{S}^0_{\psi^{\prime}})=0$ if $\psi^{\prime} \notin \Psi_{\sdisc}(G^{\prime})$ (i.e. if $|Z(\overline{S}^0_{\psi^{\prime}})|=\infty$). Hence the sum (5.6.4) can be restricted to the double sum
\begin{eqnarray*}
\sum_{G^{\prime} \in \mathcal{E}_{\ellip}(G) \backslash \mathcal{E}_{\simp}(G)} \sum_{\psi^{\prime} \in  \Psi_{\sdisc}(G^{\prime}) \cap \Psi(G^{\prime},\psi^N)}. 
\end{eqnarray*}

The next step is to transform the double sum (5.6.4). The argument is as in section 4.4 of \cite{A1} so we will be brief.

As in the local setting described in section 3.2 we have a bijective correspondence:
\begin{eqnarray}
 (G^{\prime},\psi^{\prime}) \leftrightarrow  (\psi_G,s) 
\end{eqnarray}
where on the right hand side $\psi_G= \Psi(G)$ with $\widetilde{\psi}_G$ being an {\it actual} $L$-homomorphism $\widetilde{\psi}_G: \mathcal{L}_{\psi} \times \SL_2(\mathbf{C}) \rightarrow \leftexp{L}{G^0}$, and $s$ is an element of $\overline{S}_{\psi,ss}$; for the left hand side $G^{\prime}$ is an element in the set $E(G)$ consisting of endoscopic datum $G^{\prime}=(G^{\prime},\zeta^{\prime})$ for $G$, taken up to the image of $\zeta^{\prime}(\mathcal{G}^{\prime})$ in $\leftexp{L}{G^0}$, and up to translation of the associated semi-simple element $s \in \widehat{G}$ by $Z(G^0)^{\Gamma}$ (thus elements of  are {\it not} considered up to equivalence of endoscopic data), and $\widetilde{\psi}^{\prime}$ is similarly an actual $L$-homomorphism to $\leftexp{L}{G^{\prime}}$ such that $\widetilde{\psi}_G = \zeta^{\prime} \circ \widetilde{\psi}^{\prime}$. 

Note that $\widehat{G}^0$ acts on $E(G)$ by conjugation, and we have $\mathcal{E}(G) = \widehat{G}^0 \backslash \backslash E(G)$ (here and below we use ``$ \backslash \backslash$" to denote the set of orbits). We similarly define $E_{\ellip}(G),E_{\simp}(G)$, so that $\mathcal{E}_{\ellip}(G) = \widehat{G}^0 \backslash \backslash E_{\ellip}(G)$ (and similarly for $\mathcal{E}_{\simp}(G)$).  

Denote  
\begin{eqnarray}
Y_{\sdisc,\psi^N}(G) = \{y=(\psi_G,s_G)      \}
\end{eqnarray}
with $\psi_G$ as above in (5.6.7) whose $\widehat{G}^0$-conjugacy class belongs to $\Psi(G,\psi^N)$, and $s_G \in \overline{S}_{\psi_G,\ellip}$. We also define the subset
\[
Y^{\prime}_{\sdisc,\psi^N}(G) \subset Y_{\sdisc,\psi^N}(G)
\]
consisting of $y=(\psi_G,s_G)$ such that $s_G \in \overline{S}_{\psi_G,\ellip}^{\prime}$, where
\begin{eqnarray}
\overline{S}_{\psi_G,\ellip}^{\prime} := \{ s \in \overline{S}_{\psi_G,\ellip}| \,\ G^{\prime}_s \notin E_{\simp}(G)               \}
\end{eqnarray}
where $G^{\prime}_s$ is the endoscopic datum of $G$ corresponding to the pair $(\psi_G,s)$ (in the case where $G=G^0$ then $\overline{S}^{\prime}_{\psi_G,\ellip}= \overline{S}_{\psi_G,\ellip} \smallsetminus \{1\}$). 

\bigskip

\noindent The bijection (5.6.7) restricts to a bijection between the set $Y_{\sdisc,\psi_N}(G)$, and the set of $(G^{\prime},\psi^{\prime})$ such that $G^{\prime} \in E_{\ellip}(G)$, and that the $\widehat{G}^{\prime}$-conjugacy class of $\psi^{\prime}$ belongs to $\Psi_{\sdisc}(G^{\prime}) \cap \Psi(G,\psi^N)$; similarly it restricts to a bijection between the set $Y^{\prime}_{\sdisc,\psi^N}(G)$ and the set of $(G^{\prime},\psi^{\prime})$, such that in addition $G^{\prime} \in E_{\ellip}(G) \smallsetminus E_{\simp}(G)$ ({\it c.f.} p.45--46 of \cite{A9}; note that in {\it loc. cit.} the term ``weakly elliptic parameter" is used instead of the term ``stably discrete parameter" in \cite{A1} and here). 

The group $\widehat{G}^0$ acts on $Y^{\prime}_{\sdisc,\psi_N}(G)$ by conjugation on each component. The summand in (5.6.4) (namely the product of (5.6.5) and (5.6.6)) is constant on any $(G^{\prime},\psi^{\prime})$ that lies in the same $\widehat{G}^0$-orbit. However in the sum (5.6.4) we regard $G^{\prime}$ as an element of $\mathcal{E}_{\ellip}(G) = \widehat{G}^0 \backslash \backslash E_{\ellip}(G)$, and $\psi^{\prime}$ as an element of $\Psi(G^{\prime})$ (i.e. $\psi^{\prime}$ is considered up to $\widehat{G}^{\prime}$-conjugacy).

\bigskip

In terms of the bijection (5.6.7) we have the obvious fibration
\begin{eqnarray}
\widehat{G}^0 \backslash \backslash Y^{\prime}_{\sdisc,\psi^N} \longrightarrow \mathcal{E}_{\ellip}(G) \smallsetminus \mathcal{E}_{\simp}(G).
\end{eqnarray}

\bigskip

The group $\widehat{G}^0$ acts transitively on the fibres of (5.6.10). Its fibres can be described as follows. Let $y=(\psi_G,s_G) \in Y^{\prime}_{\disc,\psi^N}$ which corresponds to the pair $(G^{\prime},\psi^{\prime})$ under (5.6.7). First note that the stabilizer of $y$ in $\widehat{G}^0$ is given by the group
\begin{eqnarray}
S_y^+ = S^+_{\psi_G,s_G} := \Cent(s_G,S^*_{\psi_G}).
\end{eqnarray} 

\noindent On the other hand the stabilizer in $\widehat{G}^0$ of $G^{\prime}$ as an element of $\mathcal{E}_{\ellip}(G) \smallsetminus \mathcal{E}_{\simp}(G)$ is given by the automorphism group $\Aut_G(G^{\prime})$. Thus the fibre of (5.6.10) over $G^{\prime}$ (as an element of $\mathcal{E}_{\ellip}(G) \smallsetminus \mathcal{E}_{\simp}(G)$) is in bijection with the quotient
\[
\Aut_G(G^{\prime})/S_y^+.
\]
Since $\Psi(G^{\prime},\psi^N)$ is the set of parameters up to $\widehat{G}^{\prime}$-conjugacy, we see that the number of parameters in $\Psi(G^{\prime},\psi^N)$ that correspond to the fibre of (5.6.10) above $G^{\prime}$ (as an element of $\mathcal{E}_{\ellip}(G) \smallsetminus \mathcal{E}_{\simp}(G)$) is given by the cardinality of the quotient:
\begin{eqnarray}
\Aut_G(G^{\prime})/ \Int_G(G^{\prime}) S_y^+
\end{eqnarray} 
where $\Int_G(G^{\prime}):= \widehat{G}^{\prime} Z(\widehat{G}^0)^{\Gamma}$. Since $\Aut_G(G^{\prime})/\Int_G(G^{\prime}) =\Out_G(G^{\prime})$, we see that the cardinality of (5.6.12) is given by the cardinality of:
\begin{eqnarray}
\Out_G(G^{\prime}) / (\Int_G(G^{\prime}) S_y^+/\Int_G(G^{\prime})).
\end{eqnarray}

\noindent Since
\[
\Int_G(G^{\prime}) S_y^+/\Int_G(G^{\prime}) \cong S_y^+ / (S_y^+ \cap  \widehat{G}^{\prime} Z(\widehat{G}^0)^{\Gamma}    )
\]
the cardinality of (5.6.13) is given by
\begin{eqnarray}
|\Out_G(G^{\prime})|  |S_y^+/(S_y^+ \cap \widehat{G}^{\prime} Z(\widehat{G}^0)^{\Gamma} )  |^{-1}.
\end{eqnarray}

\bigskip

Thus the double sum (5.6.4) can be written as a simple sum over the set $\widehat{G}^0 \backslash \backslash Y^{\prime}_{\sdisc,\psi^N}(G)$, provided we multiply each summand by the number (5.6.14). Now for any $y=(\psi_G,s_G) \in Y^{\prime}_{\sdisc,\psi^N}(G)$, the $\psi_G$ are all conjugate under $\widehat{G}^0$ (because the $\widehat{G}^0$-conjugacy class of $\psi_G$ belongs to $\Psi(G,\psi^N)$, and the fact that $\widetilde{\Out}_N(G)$ is trivial). Since the stabilizer of any such $\psi_G$ in $\widehat{G}^0$ is given by $S^*_{\psi_G}$, we see that, provided $\Psi(G,\psi^N)$ is not empty (i.e. $m^G_{\psi^N}=1$), the set $\widehat{G}^0 \backslash \backslash Y^{\prime}_{\sdisc,\psi^N}(G)$ is in bijection with the quotient  
\begin{eqnarray}
 S^*_{\psi} \backslash \backslash \overline{S}_{\psi,\ellip}^{\prime} = \overline{S}^*_{\psi} \backslash \backslash \overline{S}_{\psi,\ellip}^{\prime}. 
\end{eqnarray}

\noindent And the sum over $\widehat{G}^0 \backslash \backslash Y^{\prime}_{\sdisc,\psi^N}$ can be replaced by the sum over (5.6.15) (with the sum over $\psi_G$ collapsed and each $\psi_G$ can be identified as $\psi$), provided we multiply the sum over (5.6.15) with the number
\begin{eqnarray}
m^G_{\psi^N} =|\Psi(G,\psi^N)|.
\end{eqnarray}

In order to utilize the quotient set $\mathcal{E}^{\prime}_{\psi,\ellip}  :=\mathcal{E}(\overline{S}^{\prime}_{\psi,\ellip})$ as defined in (5.1.6), we consider the quotient of $\overline{S}_{\psi,\ellip}^{\prime}$ by $\overline{S}_{\psi}^0 = (\overline{S}^*_{\psi})^0$ instead of by $\overline{S}^*_{\psi}$ in (5.6.15). Given $s \in \overline{S}^{\prime}_{\psi,\ellip}$, the centralizer of $s$ in $\overline{S}^*_{\psi}$ is given by 
\[
\overline{S}^+_{\psi,s} := \Cent(s,\overline{S}^*_{\psi})
\]
hence the orbit of $s$ under the conjugation action of $\overline{S}^*_{\psi}$ is in bijection with 
\[
\overline{S}^*_{\psi} / \overline{S}^+_{\psi,s}. 
\]
On the other hand the centralizer of $s$ in $\overline{S}^0_{\psi}$ is given by
\[
\overline{S}_{\psi,s} := (\overline{S}^0_{\psi})_s = \Cent(s,\overline{S}_{\psi}^0)
\]
hence the orbit of $s$ under the conjugation action of $\overline{S}_{\psi}^0$ is in bijection with
\[
\overline{S}^0_{\psi}  / \overline{S}_{\psi,s} .
\]
Hence the sum over (5.6.15) can be replaced by the sum over
\begin{eqnarray}
\overline{S}^0_{\psi} \backslash \backslash \overline{S}_{\psi,\ellip}^{\prime} = \mathcal{E}_{\psi,\ellip}^{\prime}
\end{eqnarray}
provided that we multiply the summand by the rescaling constant
\begin{eqnarray}
& & \Big|   ( \overline{S}^*_{\psi} / \overline{S}_{\psi,s}^+) \Big/ (\overline{S}^0_{\psi} / \overline{S}_{\psi,s})   \Big|^{-1} \\
&=&   | \overline{S}^+_{\psi,s} / \overline{S}_{\psi,s} |  | \overline{S}^*_{\psi} / \overline{S}^0_{\psi}   |^{-1}  .               \nonumber
\end{eqnarray}

\noindent To conclude the above analysis, we see that the double sum (5.6.4) can be replaced by the sum over $\mathcal{E}_{\psi,\ellip}^{\prime}$, with the summand being given by the product of (5.6.5), (5.6.6), (5.6.14), (5.6.16) and (5.6.18), i.e. is the product of the following two expressions:
\begin{eqnarray}
|S^+_{\psi,s}/(S^+_{\psi,s} \cap    \widehat{G}^{\prime}  Z(\widehat{G}^0)^{\Gamma}   ) |^{-1}  |\mathcal{S}_{\psi^{\prime}}|^{-1} |\overline{Z}(\widehat{G}^{\prime})^{\Gamma}|^{-1} |\overline{S}^+_{\psi,s}/\overline{S}_{\psi,s} |
\end{eqnarray}
\begin{eqnarray}
|\kappa_G|^{-1} m^G_{\psi^N} |\overline{S}^*_{\psi}/ \overline{S}^0_{\psi} |^{-1} \sigma(\overline{S}^0_{\psi^{\prime}}) \epsilon^{\prime}(\psi^{\prime}) f^{\prime}(\psi^{\prime}).
\end{eqnarray}

To compute this product, we first note that the term 
\[
|S_{\psi,s}^+ / (S_{\psi,s}^+ \cap \widehat{G}^{\prime}  Z(\widehat{G}^0)^{\Gamma})  |
\]
occuring in (5.6.19) can be written as
\begin{eqnarray}
|\overline{S}_{\psi,s}^+/ (\overline{S}^+_{\psi,s} \cap \overline{\widehat{G}^{\prime}} ) |
\end{eqnarray}
where 
\begin{eqnarray*}
& & \overline{\widehat{G}^{\prime}} := \widehat{G}^{\prime} Z(\widehat{G}^0)^{\Gamma} / Z(\widehat{G}^0)^{\Gamma} \cong \widehat{G}^{\prime}/(\widehat{G}^{\prime} \cap Z(\widehat{G}^0)^{\Gamma}).
\end{eqnarray*}

Now we have (using that $(\overline{S}_{\psi,s}^+)^0 =\Cent(s,\overline{S}_{\psi})^0 =\Cent(s,\overline{S}_{\psi}^0)^0 = \overline{S}^0_{\psi,s}$):
\begin{eqnarray*}
& & | \overline{S}^+_{\psi,s} / (\overline{S}^+_{\psi,s} \cap \overline{\widehat{G}^{\prime}} )|^{-1} |\mathcal{S}_{\psi^{\prime}}  |^{-1} \\
&=&  | \overline{S}^+_{\psi,s} / (\overline{S}^+_{\psi,s} \cap \overline{\widehat{G}^{\prime}} )|^{-1} | (\overline{S}^+_{\psi,s} \cap \overline{\widehat{G}^{\prime}}   ) / ( \overline{S}^0_{\psi,s} \overline{Z}(\widehat{G}^{\prime})^{\Gamma}      )   |^{-1} \\
&=& | \overline{S}^+_{\psi,s}     / \overline{S}^0_{\psi,s} \overline{Z}(\widehat{G}^{\prime})^{\Gamma}  |^{-1}.
\end{eqnarray*}
Hence (5.6.19) is equal to
\begin{eqnarray}
& & |\overline{S}_{\psi,s}/\overline{S}^0_{\psi,s}  \overline{Z}(\widehat{G}^{\prime})^{\Gamma}|^{-1} |\overline{Z}(\widehat{G}^{\prime})^{\Gamma} |^{-1}\\
&=&   |\overline{S}_{\psi,s}/\overline{S}^0_{\psi,s} |^{-1} |\overline{S}^0_{\psi,s} \cap \overline{Z}(\widehat{G}^{\prime})^{\Gamma} |^{-1}     \nonumber\\
&=& |\pi_0(\overline{S}_{\psi,s})|^{-1} |\overline{S}^0_{\psi,s} \cap \overline{Z}(\widehat{G}^{\prime})^{\Gamma} |^{-1}.     \nonumber
\end{eqnarray}

\noindent For the expression (5.6.20), we have $|\overline{S}^*_{\psi}/\overline{S}^0_{\psi} |= |\mathcal{S}_{\psi}|$. For the number $\sigma(\overline{S}^0_{\psi^{\prime}})$, we use:
\begin{eqnarray*}
\overline{S}^0_{\psi^{\prime}} \cong \overline{S}^0_{\psi,s}/(\overline{S}^0_{\psi,s} \cap  \overline{Z}(\widehat{G}^{\prime})^{\Gamma}).
\end{eqnarray*}
Hence from the property (5.1.9) we have
\[
\sigma(\overline{S}^0_{\psi^{\prime}}) = \sigma(\overline{S}^0_{\psi,s}) | \overline{S}^0_{\psi,s} \cap \overline{Z}(\widehat{G}^{\prime})^{\Gamma}  |.
\]
Thus the product of (5.6.19) and (5.6.20) is
\begin{eqnarray}
 & & |\kappa_G|^{-1} \frac{m^G_{\psi^N}}{ |\mathcal{S}_{\psi}|} |\pi_0(\overline{S}_{\psi,s})|^{-1} \epsilon^{\prime}(\psi^{\prime}) \sigma(\overline{S}^0_{\psi,s}) f^{\prime}(\psi^{\prime}) \\
&=&  |\kappa_G|^{-1} \frac{m^G_{\psi^N}}{ |\mathcal{S}_{\psi}|} |\pi_0(\overline{S}_{\psi,s})|^{-1} \epsilon^{\prime}(\psi^{\prime}) \sigma(\overline{S}^0_{\psi,s}) f^{\prime}_G(\psi,s). \nonumber
\end{eqnarray}

It remains to treat the factor $\epsilon^{\prime}(\psi^{\prime})$ in (5.6.23). For this we need the:

\begin{lemma} (The endoscopic sign lemma)
For any $s \in \overline{S}_{\psi}$, we have
\begin{eqnarray}
\epsilon^{\prime}(\psi^{\prime}) = \epsilon^G_{\psi}(s_{\psi} x_s).
\end{eqnarray}
\end{lemma}

\noindent The endoscopic sign lemma will be established in section 5.8 together with the spectral sign lemma. Take this for granted, we thus obtain the following expression for the difference $I^G_{\disc,\psi^N}(f) - s^G_{\disc,\psi^N}(f)$:

\begin{eqnarray}
& & I^G_{\disc,\psi^N}(f) - s^G_{\disc,\psi^N}(f) \\
&=& |\kappa_G|^{-1} \frac{m^G_{\psi^N}}{|\mathcal{S}_{\psi}|} \sum_{s \in \mathcal{E}^{\prime}_{\psi,\ellip}} \epsilon^G_{\psi}(s_{\psi} x_s)  |\pi_0(\overline{S}_{\psi,s})|^{-1}   \sigma(\overline{S}^0_{\psi,s})  f^{\prime}_G(\psi,s). \nonumber
\end{eqnarray}

\noindent We now make a change of variable $s \mapsto s_{\psi}^{-1} s = s_{\psi} s $. We have $\overline{S}_{\psi,s}= \overline{S}_{\psi,s_{\psi}s}$ since $s_{\psi}$ is a central element of $\overline{S}_{\psi}$. Then we have:
\begin{eqnarray}
& & I^G_{\disc,\psi^N}(f) - s^G_{\disc,\psi^N}(f) \\
&=& |\kappa_G|^{-1} \frac{m^G_{\psi^N}}{|\mathcal{S}_{\psi}|} \sum_{s \in \mathcal{E}^{\prime}_{\psi,\ellip}} \epsilon^G_{\psi}( x_s)  |\pi_0(\overline{S}_{\psi,s})|^{-1}   \sigma(\overline{S}^0_{\psi,s})  f^{\prime}_G(\psi,s_{\psi} s) \nonumber \\
&=&  |\kappa_G|^{-1} \frac{m^G_{\psi^N}}{|\mathcal{S}_{\psi}|} \sum_{x \in \mathcal{S}_{\psi}} \sum_{s \in \mathcal{E}^{\prime}_{\psi,\ellip}(x)} \epsilon^G_{\psi}( x)  |\pi_0(\overline{S}_{\psi,s})|^{-1}   \sigma(\overline{S}^0_{\psi,s})  f^{\prime}_G(\psi,s_{\psi} x) \nonumber
\end{eqnarray}
here we have denoted by
\[
\mathcal{E}^{\prime}_{\psi,\ellip}(x)
\]
the fibre of $\mathcal{E}^{\prime}_{\psi,\ellip}$ above $x$ under the natual projection $\mathcal{E}_{\psi,\ellip}^{\prime} \rightarrow \mathcal{S}_{\psi}$.

\bigskip

As noted above for $s \in \overline{S}_{\psi,\ellip} \smallsetminus \overline{S}^{\prime}_{\psi,\ellip}$ the pair $(G^{\prime},\psi^{\prime})$ corresponding to $(\psi,s)$ satisfies $G^{\prime} \in \mathcal{E}_{\simp}(G)$. Suppose we happen to know that the linear form $f^{\prime}_G(\psi,s) =f^{\prime}(\psi^{\prime})$ is defined for all $s \in \overline{S}_{\psi,\ellip}$.  Then parallel to the spectral side of (5.5.16) and (5.5.17) we define
\begin{eqnarray}
& & \\
& & \leftexp{0}{s^G_{\disc,\psi^N}(f)} := \sum_{G^{\prime} \in \mathcal{E}_{\simp}(G)} \iota(G,G^{\prime}) \leftexp{0}{\widehat{S}^{G^{\prime}}_{\disc,\psi^N} (f^{G^{\prime}} ) }, \,\ f \in \mathcal{H}(G) \nonumber
\end{eqnarray}
where
\begin{eqnarray}
\leftexp{0}{S^{G^{\prime}}_{\disc,\psi^N}(f^{\prime})} := S^{G^{\prime}}_{\disc,\psi^N}(f^{\prime}) - \frac{m^{G^{\prime}}_{\psi^N}}{|\mathcal{S}_{\psi^{\prime}}|} \epsilon^{\prime}(\psi^{\prime}) \sigma(\overline{S}^0_{\psi^{\prime}}) f^{\prime}(\psi^{\prime}).
\end{eqnarray}

\bigskip

\noindent In other words, we have $\leftexp{0}{S^{G^{\prime}}_{\disc,\psi^N}}$ vanishes if and only if the stable multiplicity formula holds for the distribution $S^{G^{\prime}}_{\disc,\psi^N}$. We can then carry out the same argument as above, but with the set $\overline{S}_{\psi,\ellip},Y_{\sdisc,\psi^N}(G), \mathcal{E}_{\psi,\ellip} := \mathcal{E}(\overline{S}_{\psi,\ellip})$ and $\mathcal{E}_{\psi,\ellip}(x)$, in place of $\overline{S}^{\prime}_{\psi,\ellip},Y^{\prime}_{\sdisc,\psi^N}(G) ,\mathcal{E}^{\prime}_{\psi,\ellip}, \mathcal{E}_{\psi,\ellip}(x)$, we obtain 

\begin{eqnarray}
& & I^G_{\disc,\psi^N}(f) - \leftexp{0}{s^G_{\disc,\psi^N}(f)}\\
& =& |\kappa_G|^{-1}\frac{m^G_{\psi^N}}{|\mathcal{S}_{\psi}|} \sum_{x \in \mathcal{S}_{\psi}} \sum_{s \in \mathcal{E}_{\psi,\ellip}(x)}  \epsilon^G_{\psi}(x) |\pi_0(\overline{S}_{\psi,s})|^{-1} \sigma(\overline{S}_{\psi,s}^0) f^{\prime}_G(\psi,s_{\psi}x) \nonumber \\
&=& |\kappa_G|^{-1}\frac{m^G_{\psi^N}}{|\mathcal{S}_{\psi}|} \sum_{x \in \mathcal{S}_{\psi}}  e_{\psi}(x) \epsilon^G_{\psi}(x)  f^{\prime}_G(\psi,s_{\psi}x) \nonumber
\end{eqnarray}
 here 
\begin{eqnarray}
e_{\psi}(x):= \sum_{s \in \mathcal{E}_{\psi,\ellip}(x)} |\pi_0(\overline{S}_{\psi,s})|^{-1} \sigma(\overline{S}^0_{\psi,s})
\end{eqnarray}
is the number defined in (5.1.7) for the component of $\overline{S}_{\psi}$ indexed by $x$. We record this as:

\begin{proposition}
We have:
\begin{eqnarray}
& & I^G_{\disc,\psi^N}(f) - s^G_{\disc,\psi^N}(f) \\
&=& |\kappa_G|^{-1} \frac{m^G_{\psi^N}}{|\mathcal{S}_{\psi}|} \sum_{s \in \mathcal{E}^{\prime}_{\psi,\ellip}} \epsilon^G_{\psi}(s_{\psi} x_s)  |\pi_0(\overline{S}_{\psi,s})|^{-1}   \sigma(\overline{S}^0_{\psi,s})  f^{\prime}_G(\psi,x). \nonumber
\end{eqnarray}

If in addition the linear form $f^{\prime}_G(\psi,s)$ is defined for all $s \in \overline{S}_{\psi,\ellip}$ then we have
\begin{eqnarray}
& & I^G_{\disc,\psi^N}(f) - \leftexp{0}{s^G_{\disc,\psi^N}(f)}\\
&=& |\kappa_G|^{-1}\frac{m^G_{\psi^N}}{|\mathcal{S}_{\psi}|} \sum_{x \in \mathcal{S}_{\psi}}  e_{\psi}(x) \epsilon^G_{\psi}(x)  f^{\prime}_G(\psi,s_{\psi}x). \nonumber
\end{eqnarray}
\end{proposition}

\subsection{The comparison}

We can now begin the first step in the term by term comparison of the spectral expansion (5.5.20) and the endoscopic expansion (5.6.32) (under the appropriate hypothesis for their validity). As in the previous subsections $G$ will denote either an element of $\widetilde{\mathcal{E}}_{\simp}(N)$ or the twisted group $\widetilde{G}_{E/F}(N)$.

We fix a parameter $\psi^N \in \widetilde{\Psi}(N)$. In the case where $G=(G,\xi) \in \widetilde{\mathcal{E}}_{\simp}(N)$ we will, to ease notation, abbreviate the condition $\psi^N \in \xi_* \psi(G)$ just as $\psi^N \in \Psi(G)$ (similar remarks apply to $\Psi_2(G),\Psi_{\ellip}(G)$, etc), in which case we will just denote by $\psi$ or $\psi_G$ the parameter in $\Psi(G)$ defined by $\psi^N$, i.e. such that $\psi=(\psi^N,\widetilde{\psi})$ and hence $\psi^N = \xi_* \psi$ (in the case where $G=\widetilde{G}_{E/F}(N)$ we interpret the condition $\psi^N \in \Psi(G)$ as automaticallly satisfied, in which case $\psi^N$ and $\psi$ refer to the same parameter).

A key ingredient in the comparison is the global intertwining relation for $\psi \in \Psi(G)$:
\begin{eqnarray}
f_G(\psi,u) =f^{\prime}_G(\psi,s_{\psi} s), \,\ u \in \mathfrak{N}_{\psi}, s \in \overline{S}_{\psi}
\end{eqnarray}
for $u$ and $s$ having the same image in $\mathcal{S}_{\psi}$. This of course still needs to be proved, along with the other main global theorems.

\bigskip

\begin{proposition}
With the above notations, suppose that either one of the following conditions hold:
\begin{enumerate}

\item $G \in \widetilde{\mathcal{E}}_{\simp}(N)$ and $\psi^N \notin \Psi(G)$.

\item $\psi^N \in \Psi(G)$, the linear forms $f_G(\psi,u)$ and $f_G^{\prime}(\psi,s)$ are defined for all $u \in \mathfrak{N}_{\psi}$ and $s \in \overline{S}_{\psi}$, and the global intertwining relation (5.7.1) is satisifed.

\end{enumerate}

\noindent Then we have 
\begin{eqnarray}
\leftexp{0}{r^G_{\disc,\psi^N}(f)} = \leftexp{0}{s^G_{\disc,\psi^N}(f)}, \,\ f \in \mathcal{H}(G).
\end{eqnarray}
\noindent In fact in case (1) we have
\begin{eqnarray}
\tr R_{\disc,\psi^N}^G(f) = \leftexp{0}{r^G_{\disc,\psi^N}(f)}= \leftexp{0}{S^G_{\disc,\psi^N}(f)} = S^G_{\disc,\psi^N}(f).
\end{eqnarray}
\end{proposition}
\begin{proof}
We first consider case (1). Then we have $m^G_{\psi^N}=0$, so the spectral expansion (5.5.19) and the endoscopic expansion (5.6.31) gives:
\begin{eqnarray}
& & \\
& & I^G_{\disc,\psi^N}(f) - \leftexp{0}{r^G_{\disc,\psi^N}(f)} = 0 = I^G_{\disc,\psi^N}(f) - s^G_{\disc,\psi^N}(f) \nonumber
\end{eqnarray}
and hence the assertions, upon noting that in the case where $\psi^N \notin \Psi(G)$, we have, from the definitions (5.5.17), (5.6.27) and (5.6.28):
\begin{eqnarray*}
& & \leftexp{0}{r^G_{\disc,\psi^N}(f)} = r^G_{\disc,\psi^N}(f) = \frac{1}{|\kappa_G|} \tr R^G_{\disc,\psi^N}(f) = \tr R^G_{\disc,\psi^N}(f) \\
& & \leftexp{0}{s^G_{\disc,\psi^N}(f)} =\leftexp{0}{S^G_{\disc,\psi^N}(f)} = S^G_{\disc,\psi^N}(f) 
\end{eqnarray*}
(since $|\kappa_G|=1$ for $G \in \widetilde{\mathcal{E}}_{\simp}(N)$). 

\bigskip

\noindent In case (2) we apply the same arguments (we have $m^G_{\psi^N}=1$ in this case). Namely that from the given assumption we have the validity of the expansions (5.5.20) and (5.6.32):
\begin{eqnarray*}
& & I^G_{\disc,\psi^N}(f) - \leftexp{0}{r^G_{\disc,\psi^N}(f)} = |\kappa_G|^{-1} \frac{1}{|\mathcal{S}_{\psi}|}\sum_{x \in \mathcal{S}_{\psi}} \epsilon_{\psi}(x) i_{\psi}(x) f_{G}(\psi,x) \\
& & I^G_{\disc,\psi^N}(f) - \leftexp{0}{s^G_{\disc,\psi^N}(f)} = |\kappa_G|^{-1} \frac{1}{|\mathcal{S}_{\psi}|}\sum_{x \in \mathcal{S}_{\psi}} \epsilon_{\psi}(x) e_{\psi}(x) f_{G}^{\prime}(\psi,s_{\psi} x).
\end{eqnarray*}
Hence from the equality $i_{\psi}(x) = e_{\psi}(x)$ for all $x \in \mathcal{S}_{\psi}$ of equation (5.1.8), and the global intertwining relation (5.7.1),  we obtain a term by term identification of the two expansions. Hence we obtain:
\[
 I^G_{\disc,\psi^N}(f) - \leftexp{0}{r^G_{\disc,\psi^N}(f)}  = I^G_{\disc,\psi^N}(f) - \leftexp{0}{s^G_{\disc,\psi^N}(f)}
\]
and the assertion follow. 

\end{proof}

We would in fact like to prove that the distributions in (5.7.2) actually vanish, i.e. both the spectral multiplicity and the stable multiplicity formulas are valid (in case (1) the spectral multiplicity formula is interpreted as the condition that $\psi^N$ does not contribute to the discrete spectrum of $G=(G,\xi)$ with respect to the $L$-embedding $\xi$). This proof of this will be complete only at the end of the induction argument in section 9. In this subsection however, we will be able to treat the case a class of ``degenerate" parameters $\psi^N$ from the induction hypothesis.

\begin{lemma}
Suppose that $\psi^N \notin \widetilde{\Psi}_{\ellip}(N)$. The for $G=\widetilde{G}_{E/F}(N)$ or any $G \in \widetilde{\mathcal{E}}_{\simp}(N)$ such that $\psi^N \in \Psi(G)$, the linear forms $f_G(\psi,u)$ and $f_G^{\prime}(\psi,s)$ are defined for all $u \in \mathfrak{N}_{\psi}$ and $s \in \overline{S}_{\psi}$. Furthermore, in the case where $G=\widetilde{G}_{E/F}(N)$ the global intertwining relation (5.7.1) is valid for $\psi^N$.
\end{lemma}
\begin{proof}
Since $\psi^N \notin \widetilde{\Psi}_{\ellip}(N)$, we see that $M$ is proper in $G^0$, and hence the spectral linear form $f_G(\psi,u)$ is defined. 

As for the endoscopic distribution, note that for any $s \in \overline{S}_{\psi}$, if $(G^{\prime},\psi^{\prime})$ is the pair corresponding to $(\psi,s)$ (thus here $G^{\prime} \in \mathcal{E}(G)$ and $\psi^{\prime} \in \Psi(G^{\prime})$), then the condition that $\psi^N \notin \widetilde{\Psi}_{\ellip}(N)$ implies that $\psi^{\prime} \notin \Psi_2(G^{\prime})$. Thus $\psi^{\prime}$ factors through the $L$-group of a proper Levi subgroup of $G^{\prime}$. It then follows from the induction hypothesis, together with descent argument as before, that the stable linear form $(f^{\prime})^{G^{\prime}}$ is defind for $f^{\prime} \in \mathcal{H}(G^{\prime})$. Hence $f^{\prime}(\psi,s) = f^{G^{\prime}}(\psi^{\prime})$ (here $f \in \mathcal{H}(G)$) is defined.

Finally we already know from proposition 5.2.3 the validity of the global intertwining relation for $\widetilde{G}_{E/F}(N)$ in this case, once the endoscopic distributions are defined. It is thus valid in this case.
\end{proof}

Thus we assume that $\psi_N \notin \widetilde{\Psi}_{\ellip}(N)$. We first consider the case where $G=\widetilde{G}_{E/F}(N)$, in which case as in section 4 we denote the distributions $r^G_{\disc,\psi^N}$ and $s^G_{\disc,\psi^N}$ as $\widetilde{r}^N_{\disc,\psi^N}$ and $\widetilde{s}^N_{\disc,\psi^N}$, etc. Since $\psi_N$ does not lie in $\widetilde{\Psi}_{\ellip}(N)$, hence not in $\widetilde{\Psi}_{\simp}(N)$, we have
\[
\leftexp{0}{\widetilde{r}^N_{\disc,\psi^N}(\widetilde{f})} = \widetilde{r}^N_{\disc,\psi^N}(\widetilde{f}) = \frac{1}{2}\tr \widetilde{R}^N_{\disc,\psi^N}(\widetilde{f}) =0, \,\ \widetilde{f} \in \widetilde{\mathcal{H}}(N)
\] 
by the theorem of Moeglin-Waldspurger \cite{MW} and Jacquet-Shalika \cite{JS} (recall also that $|\kappa_{\widetilde{G}(N)}|=2$). Thus by case (2) of proposition 5.7.1 and lemma 5.7.2 (applied to the case $G=\widetilde{G}_{E/F}(N)$), we have $\leftexp{0}{\widetilde{s}^N_{\disc,\psi^N}(\widetilde{f})}=0$. In other words
\begin{eqnarray}
\sum_{G \in \widetilde{\mathcal{E}}_{\simp}(N)} \widetilde{\iota}(N,G) \leftexp{0}{\widehat{S}^G_{\disc,\psi^N}(\widetilde{f}^G) }  =0.
\end{eqnarray} 

Now for any compatible family of function $\mathcal{F}=\{f \in \mathcal{H}(G)| \,\ G \in \widetilde{\mathcal{E}}_{\simp}(N) \}$ (i.e. for $G \in \widetilde{\mathcal{E}}_{\ellip}(N) \smallsetminus \widetilde{\mathcal{E}}_{\simp}(N)$ the function associated to $G$ is identifically zero), there exists, by consequence of proposition 3.3.1, a function $\widetilde{f} \in \widetilde{\mathcal{H}}(N)$ such that $\widetilde{f}^G=f^G$ for any $G \in \widetilde{\mathcal{E}}_{\simp}(N)$ (with $f \in \mathcal{F}$ being the function associated to $G \in \widetilde{\mathcal{E}}_{\simp}(N)$). We can then replace the term $\leftexp{0}{\widehat{S}^G_{\disc,\psi^N}(\widetilde{f}^G)}$ in (5.7.5) by $\leftexp{0}{S^G_{\disc,\psi^N}(f)}$ for each $G \in \widetilde{\mathcal{E}}_{\simp}(N)$. We thus obtain:
\begin{corollary}
Suppose that $\psi^N \notin \widetilde{\Psi}_{\ellip}(N)$. Then for any compatible family $\mathcal{F}$ as above we have
\begin{eqnarray*}
\sum_{G \in \widetilde{\mathcal{E}}_{\simp}(N)} \widetilde{\iota}(N,G) \leftexp{0}{S^G_{\disc,\psi^N}(f)}  =0.
\end{eqnarray*} 
\end{corollary}
(Of course there are only two terms in the above sum, corresponding to the two equivalence of datum $(U_{E/F}(N),\xi) \in \widetilde{\mathcal{E}}_{\simp}(N)$.)
\bigskip

In order to combine the information from corollary 5.7.3 and proposition 5.7.1 (applied to the case $G \in \widetilde{\mathcal{E}}_{\simp}(N)$) we must establish the global intertwining relation for $G \in \widetilde{\mathcal{E}}_{\simp}(N)$ and our fixed $\psi^N \notin \widetilde{\Psi}_{\ellip}(N)$. We will see that, with the exception of some cases, this can be proved from the induction hypothesis: 

\begin{proposition}
With $\psi^N \notin \widetilde{\Psi}_{\ellip}(N)$ as above, suppose in addition that $\psi^N \notin \Psi_{\ellip}(G)$ for any $G \in \widetilde{\mathcal{E}}_{\simp}(N)$. Then for any $G \in \widetilde{\mathcal{E}}_{\simp}(N)$ we have 
\begin{eqnarray}
 \tr R^G_{\disc,\psi^N}(f) =0= \leftexp{0}{S^G_{\disc,\psi^N}}(f) , \,\ f \in \mathcal{H}(G).
\end{eqnarray}
In particular $\psi^N$ does not contribute to the discrete spectrum of any $G=(G,\xi) \in \widetilde{\mathcal{E}}_{\simp}(N)$ (with respect to the $L$-embedding $\xi$), and the stable multiplicity formula is valid for $\psi^N$ with respect to any $G=(G,\xi) \in \widetilde{\mathcal{E}}_{\simp}(N)$.
\end{proposition}
\begin{proof}
Thus let $\psi^N \notin \widetilde{\Psi}_{\ellip}(N)$ be given. For any $G \in \widetilde{\mathcal{E}}_{\simp}(N)$, if $\psi_N \notin \Psi(G)$, then by case (1) of proposition 5.7.1 we have: 
\begin{eqnarray}
\tr R_{\disc,\psi^N}^G(f) = \leftexp{0}{S^G_{\disc,\psi^N}}(f), \,\ f \in \mathcal{H}(G).
\end{eqnarray}

Thus suppose that $\psi^N \in \Psi(G)$. Proposition 5.3.4 asserts that, if the following two conditions hold: 
\begin{eqnarray}
& & \dim \overline{T}_{\psi} \geq 2 \\
& & \dim \overline{T}_{\psi,x} \geq 1 \mbox{ for all } x \in \mathcal{S}_{\psi}
\end{eqnarray}
then theorem 5.2.1, and in particular the global intertwining relation (5.7.1), is valid for $\psi$ (with respect to $G$), in which case we again conclude from case (2) of proposition 5.7.1 that 
\begin{eqnarray}
\leftexp{0}{r^G_{\disc,\psi^N}(f)} = \leftexp{0}{S^G_{\disc,\psi^N}(f)}.
\end{eqnarray}
But (5.7.10) is exactly (5.7.7), namely that $\psi^N \notin \widetilde{\Psi}_{\ellip}(N)$ implies that $\psi \notin \Psi_2(G)$, and hence $\leftexp{0}{r^G_{\disc,\psi^N}(f)} = r^G_{\disc,\psi^N}(f) = \tr R^G_{\disc,\psi^N}(f)$.

Hence if (5.7.8) and (5.7.9) hold for any $G \in \widetilde{\mathcal{E}}_{\simp}(N)$ such that $\psi^N \in \Psi(G)$, then we conclude that (5.7.7) is valid for {\it any} $G \in \widetilde{\mathcal{E}}_{\simp}(N)$. Thus corollary (5.7.3) gives for any compatible family of functions $\mathcal{F}=\{f \in \mathcal{G}, \,\ G \in \widetilde{\mathcal{E}}_{\simp}(N) \}$:
\begin{eqnarray}
 \sum_{G \in \widetilde{\mathcal{E}}_{\simp}(N)} \widetilde{\iota}(N,G) \tr R_{\disc,\psi^N}(f) =0. 
\end{eqnarray}

\noindent Thus since all the coefficients $\widetilde{\iota}(N,G)$ are positive, and $\tr R^G_{\disc,\psi^N}$ can be written as a linear combination with positive coefficients of characters of irreducible representations on $G(\mathbf{A}_F)$, we see that the left hand side of (5.7.11) is of the form given as in the left hand side of (4.3.26). We can thus apply lemma 4.3.6 to conclude:
\[
\tr R^G_{\disc, \psi^N}(f) =0, \,\ f \in \mathcal{H}(G)
\] 
for all $G \in \widetilde{\mathcal{E}}_{\simp}(N)$, hence also the vanishing of $\leftexp{0}{S^G_{\disc,\psi^N}(f)}$ by (5.7.7) again.

\bigskip

Thus we must now deal with parameters $\psi^N$ for which (5.7.8) or (5.7.9) fail for some $G \in \widetilde{\mathcal{E}}_{\simp}(N)$ (such that $\psi^N \in \Psi(G)$). This requires a more detailed comparison than the argument above.

\bigskip

Thus let $G \in \widetilde{\mathcal{E}}_{\simp}(N)$ be fixed, such that $\psi^N \in \Psi(G)$, but (5.7.8) or (5.7.9) fail for $G$. Recall from (2.4.14) that we have a general description of the centralizer group $S_{\psi}$: if $\psi^N$ has a decomposition:
\[
\psi^N = l_1 \psi^{N_1}_{1} \boxplus \cdots \boxplus l_r \psi^{N_r}_{r}
\]
then we have
\begin{eqnarray*}
S_{\psi}= \prod_{i \in I_{\psi}^+(G)} O(l_i,\mathbf{C}) \times \prod_{i \in I^-_{\psi}(G)} Sp(l_i,\mathbf{C}) \times \prod_{j \in J_{\psi}(G)} \GL(l_j,\mathbf{C})
\end{eqnarray*}
with the meaning of $I_{\psi}^+(G),I_{\psi}^-(G)$ and $J_{\psi}(G)$ as in (2.4.12), (2.4.13). Now $\overline{T}_{\psi} = T_{\psi}/Z(\widehat{G})^{\Gamma}$, with $Z(\widehat{G})^{\Gamma}$ finite (of order two), and $T_{\psi}$ is a maximal torus of $S_{\psi}$. Thus we can easily enumerate all the possibility where (5.7.8) or (5.7.9) fails, as follows.

Suppose first that (5.7.8) fails. Then we see that the set $J_{\psi}$ is either empty or we just have one single factor $\GL(1,\mathbf{C}) = \mathbf{C}^{\times}$ for the contribution of $J_{\psi}$ to $S_{\psi}$. In the latter case $S_{\psi}$ would have a central torus of positive dimension, hence by remark 5.3.2, the global intertwining relation holds for $\psi$ with respect to $G$ by reduction to a proper Levi subgroup of $G$, and thus (5.7.7) holds. So we may assume that $J_{\psi}$ is empty. 

Next we consider the contribution of the set $I_{\psi}^-$ to $S_{\psi}$. Then we see that either $I_{\psi}^-$ is empty, or we have only one factor $Sp(2,\mathbf{C})$ in the contribution of $I_{\psi}^-$ to $S_{\psi}$.

Finally we consider the case where both $I_{\psi}^-$ and $J_{\psi}$ are empty, i.e. we only have the contribution of $O(l_i,\mathbf{C})$ for $i \in I^+_{\psi}$ to $S_{\psi}$. If $l_i \leq 2$ for all such $i$, then $\overline{S}_{\psi,\ellip}$ is non-empty, for instance by taking $s=(s_i)_i$ such that $\det s_i = -1$ for all $i \in I_{\psi}^+$ with $l_i=2$. Thus $\psi^N \in \Psi_{\ellip}(G)$ in this case, which is ruled out by our hypothesis. From this we see that the only remaining case is that $l_i=3$ for exactly one $i$, and the remaining $l_i$'s are equal to one.

By similar analysis, we see that if (5.7.8) holds but (5.7.9) fails, then the only possibility is that $I^-_{\psi}$ and $J_{\psi}$ are empty, and that $l_i \leq 2$ for all $i$, which is ruled out by the hypothesis that $\psi^N \notin \Psi_{\ellip}(G)$.

\bigskip

\noindent To summarize, we have isolated the following two cases:

\begin{eqnarray}
& & \psi^N= 2 \psi^{N_1}_1 \boxplus \psi^{N_2}_2 \boxplus \cdots \boxplus \psi^{N_r}_r \\
& & S_{\psi}(G) = Sp(2,\mathbf{C}) \times (\mathbf{Z}/2\mathbf{Z})^{r-1} \nonumber
\end{eqnarray}

\noindent and 

\begin{eqnarray}
& & \psi^N = 3\psi^{N_1}_1 \boxplus \psi^{N_2}_2 \boxplus \cdots \boxplus \psi^{N_r}_r \\
& & S_{\psi}(G) = O(3,\mathbf{C} ) \times (\mathbf{Z}/2\mathbf{Z})^{r-1}  \nonumber \\ & & \,\ \,\ \,\ \,\ \,\ \,\ = SO(3,\mathbf{C}) \times (\mathbf{Z}/2\mathbf{Z})^r. \nonumber
\end{eqnarray}

\noindent If $\psi^N$ does not belong to either of these two cases, then as we have seen the global intertwining relation holds for $\psi^N$ with respect to any $G \in \widetilde{\mathcal{E}}_{\simp}(N)$, and hence (5.7.7) is valid for all such $G$. Thus the argument in the first part of the proof can be applied to yield the assertions of the proposition in this case. 

So we must now treat the two cases (5.7.12) and (5.7.13). The treatment of these cases will be similar, so it suffices to treat the case of (5.7.12). Thus assume that $G \in \widetilde{\mathcal{E}}_{\simp}(N)$ such that $\psi^N \in \Psi(G)$ and $S_{\psi}(G)$ is an in (5.7.12). Denoting by $G^{\vee}$ the element of $\widetilde{\mathcal{E}}_{\simp}(N)$ other than $G \in \widetilde{\mathcal{E}}_{\simp}(N)$ (in other words $G$ and $G^{\vee}$ have the same underlying endoscopic group $U_{E/F}(N)$ but with different equivalence of endoscopic datum, namely the $L$-embedding to $\leftexp{L}{G_{E/F}(N)}$). We have $\psi^N \in \Psi(G^{\vee})$ if and only if $r=1$, in which case we have 
\[
S_{\psi_{G^{\vee}}}(G^{\vee}) = O(2,\mathbf{C}).
\]
But this would imply that $\psi^N \in \Psi(G^{\vee})_{\ellip}$, contradicting the hypothesis of the proposition. Thus we have $r >1$, i.e. that $G$ is the only element of $\widetilde{\mathcal{E}}_{\simp}(N)$ such that $\psi^N \in \Psi(G)$. Since $\psi^N \notin \Psi(G^{\vee})$, we know that (5.7.7) holds for $\psi^N$ with respect to $G^{\vee}$.

\bigskip

\noindent In any case, we know that the endoscopic expansion (5.6.32) is valid:
\[
I^G_{\disc,\psi^N}(f) - \leftexp{0}{S^G_{\disc,\psi^N}(f)} = \frac{1}{|\mathcal{S}_{\psi}|}\sum_{x \in \mathcal{S}_{\psi}} \epsilon^G_{\psi}(x) e_{\psi}(x) f^{\prime}_G(\psi,s_{\psi}x)  , \,\ f \in \mathcal{H}(G)
\]
(recall that by lemma 5.7.2 and our hypothesis that $\psi^N \notin \widetilde{\Psi}_{\ellip}(N)$, the distribution $f^{\prime}_G(\psi,s_{\psi} s)$ for semi-simple $s \in \overline{S}_{\psi}$ is defined, and depends only on the image $x$ of $s$ in $\mathcal{S}_{\psi}$, by lemma 5.3.1). On the spectral side, {\it a priori} we only have the expansion (5.5.19):
\begin{eqnarray*}
& & I^G_{\disc,\psi^N}(f) - \tr R^G_{\disc,\psi^N}(f)\\
 &=& \frac{1}{|\mathcal{S}_{\psi}|} \sum_{x \in \mathcal{S}_{\psi}} \epsilon^G_{\psi}(x) \Big( \frac{1}{|W_{\psi}^0|} \sum_{u \in \mathfrak{N}_{\psi,\reg}(x)} \frac{\sgn^0(w_u)}{|\det(w_u-1)|} \Big)       f_G(\psi,u).
\end{eqnarray*}

\noindent However, in the case of (5.7.12), the set $W_{\psi,\reg}=\{w\}$ is singleton, namely the unique non-trivial Weyl element of $Sp(2,\mathbf{C})$, hence $\mathfrak{N}_{\psi,\reg}(x)$ contains exactly one element for any $x \in \mathcal{S}_{\psi}$. We denote this unique element of $\mathfrak{N}_{\psi,\reg}(x)$ as $u_x$, and hence the spectral expansion just becomes:
\[
I^G_{\disc,\psi^N}(f) - \tr R^G_{\disc,\psi^N}(f) = \frac{1}{|\mathcal{S}_{\psi}|}\sum_{x \in \mathcal{S}_{\psi}} \epsilon^G_{\psi}(x) i_{\psi}(x) f_G(\psi,x)  , \,\ f \in \mathcal{H}(G)
\]
where we have taken the liberty of setting $f_G(\psi,x):=f_G(\psi,u_x)$ in this case. Combining the two expansions, and upon using the equality (5.1.8) $i_{\psi}(x) = e_{\psi}(x)$ again,  we obtain:
\begin{eqnarray}
& & \\
& & 
\leftexp{0}{S^G_{\disc,\psi^N}(f)} - \tr R^G_{\disc,\psi^N}(f) = \frac{1}{|\mathcal{S}_{\psi}|} \sum_{x \in \mathcal{S}_{\psi}} \epsilon^G_{\psi}(x) i_{\psi}(x) (f_G(\psi,x) - f^{\prime}_G(\psi,s_{\psi}x)). \nonumber
\end{eqnarray}

\noindent The number $i_{\psi}(x)$ is easily computed:
\begin{eqnarray}
& & \\
& & i_{\psi}(x) = \frac{1}{|W_{\psi}^0|} \sgn^0(w) |\det(w-1)|^{-1} = -1/4 \,\ \mbox{ for any } x \in \mathcal{S}_{\psi} \nonumber
\end{eqnarray}
since $|W_{\psi}^0|=2$, $|\det(w-1)|=2$, and $\sgn^0(w)=-1$, which is the crucial minus sign. On the other hand, we also note from the form of $S_{\psi}(G)$ in (5.7.12) that $W_{\psi}^0=W_{\psi}$. Hence the $R$-group $R_{\psi}$ is trivial, and thus we can identify $\mathcal{S}_{\psi_M} = \mathcal{S}_{\psi}$, and under this identification, if $x \in \mathcal{S}_{\psi}$ corresponds to $x_M \in \mathcal{S}_{\psi^M}$, then $\epsilon_{\psi_M}^M(x_M) = \epsilon^G_{\psi}(x)$. 

It follows, upon substituting the expression (5.2.4) of the linear form $f_G(\psi,x)=f_G(\psi,u_x)$ that
\begin{eqnarray}
& & \\
& & \frac{1}{|\mathcal{S}_{\psi}|} \sum_{x \in \mathcal{S}_{\psi}} \epsilon^G_{\psi}(x) i_{\psi}(x) f_G(\psi,x) \nonumber \\
&=&  -\frac{1}{4|\mathcal{S}_{\psi_M}|} \sum_{x_M \in \mathcal{S}_{\psi_M}} \epsilon^M_{\psi_M}(x_M) \sum_{\pi_M \in \Pi_{\psi_M}} \langle x_M ,\pi_M \rangle \tr( R_P(w,\widetilde{\pi}_M,\psi_M) \mathcal{I}_P(\pi_M,f) ) \nonumber \\
&=& -\frac{1}{4} \sum_{\pi_M \in \Pi_{\psi_M}} m(\pi_M) \tr( R_P(w,\widetilde{\pi}_M,\psi_M) \mathcal{I}_P(\pi_M,f) ) \nonumber
\end{eqnarray}
where 
\[
m(\pi_M) = \frac{1}{|\mathcal{S}_{\psi_M}|} \sum_{x_M \in \mathcal{S}_{\psi_M}} \epsilon^M_{\psi_M}(x_M) \langle x_M,\pi_M \rangle
\]
which by our induction hypothesis applied to $M \neq G$, is the multiplicity of $\pi_M$ appearing in the (relative) dsicrete spectrum of $M$.

For the endoscopic distribution $f^{\prime}_G(\psi,x)$, we have by the descent argument in the proof of lemma 5.3.1 that for any $x \in \mathcal{S}_{\psi}$:
\begin{eqnarray}
f^{\prime}_G(\psi,x)=f^{\prime}(\psi^{\prime})=f^{\prime}_M(\psi_M^{\prime})
\end{eqnarray} 
where the notations are as follows: if we denote by $x_M \in \mathcal{S}_{\psi_M}$ the element that correspond to $x \in \mathcal{S}_{\psi}$, then $(M^{\prime},\psi_M^{\prime})$ is the pair that correspond to $(\psi_M,x_M)$ (note that in the notations of the proof of lemma 5.3.1 we have $M_x=M$ in the present case). Our induction hypothesis that the local theorem 3.2.1 holds for $M$ and the parameter $\psi_M$ (more precisly the localizations of $\psi_M$ at all the places of $F$) yields the equality:
\begin{eqnarray}
f_M^{\prime}(\psi^{\prime}_M)= \sum_{\pi_M \in \Pi_{\psi_M}} \langle s_{\psi} x_M, \pi_M \rangle f_M(\pi_M)
\end{eqnarray}
(note that $s_{\psi}$ and $s_{\psi_M}$ gives the same element). Hence we have
\begin{eqnarray}
& & \\
& & \frac{1}{|\mathcal{S}_{\psi}|} \sum_{x \in \mathcal{S}_{\psi}} \epsilon^G_{\psi}(x) i_{\psi}(x) f_G^{\prime}(\psi,s_{\psi}x) \nonumber \\
&=& -\frac{1}{4 |\mathcal{S}_{\psi_M}|} \sum_{x_M \in \mathcal{S}_{\psi_M}} \epsilon^M_{\psi_M}(x_M) \sum_{\pi_M \in \Pi_{\psi_M}} \langle x_M,\pi_M \rangle f_M(\pi_M) \nonumber\\
&=& - \frac{1}{4} \sum_{\pi_M \in \Pi_{\psi_M}} m(\pi_M) f_M(\pi_M) \nonumber \\
&=& - \frac{1}{4} \sum_{\pi_M \in \Pi_{\psi_M}} m(\pi_M) \tr \mathcal{I}_P(\pi_M,f) \nonumber
\end{eqnarray}
with the last equality follows from the adjunction $f_M(\pi_M)=\tr\mathcal{I}_P(\pi_M,f)$.

Combining (5.7.14), (5.7.16) and (5.7.19), we obatin:
\begin{eqnarray}
& & \leftexp{0}{S^G_{\disc,\psi_N}(f)} - \tr R^G_{\disc,\psi_N}(f) \\
&=& \frac{1}{4} \sum_{\pi_M \in \Pi_{\psi_M}} m(\pi_M) \tr(1-R_P(w,\widetilde{\pi}_M,\psi_M) \mathcal{I}_P(\pi_M,f)), \,\ f \in \mathcal{H}(G). \nonumber
\end{eqnarray}

\bigskip

\noindent Now by corollary 5.7.3, for any compatible family of functions $\mathcal{F}=\{f_{G^*} \in \mathcal{H}(G^*), \,\ G^* \in \widetilde{\mathcal{E}}_{\simp}(N) \}$, we have
\[
\sum_{G^* \in \widetilde{\mathcal{E}}_{\simp}(N)} \widetilde{\iota}(N,G^*) \leftexp{0}{S^{G^*}_{\disc,\psi_N}(f_{G^*}) }=0.
\]
On the other hand, we already know that for $G^{\vee}$, equation (5.7.7) is valid for $G^{\vee}$. Hence subsituting (5.7.20) we obtain
\begin{eqnarray}
& & \sum_{G^* \in \widetilde{\mathcal{E}}_{\simp}(N)} \widetilde{\iota}(N,G^*) \tr R^{G^*}_{\disc,\psi^N}(f_{G^*}) \\
& & \,\ \,\ +  \frac{1}{4} \sum_{\pi_M \in \Pi_{\psi_M}} m(\pi_M) \tr(1-R_P(w,\widetilde{\pi}_M,\psi_M) \mathcal{I}_P(\pi_M,f_G)) =0. \nonumber
\end{eqnarray}

\noindent To conclude the proof we observe that since the element $w \in W_{\psi}^0 = W_{\psi}$ is of order two, the eigenvalues of the intertwining operator $R_P(w,\widetilde{\pi}_M,\psi_M)$ can only be $+1$ or $-1$. In particular it follows that we can write (5.7.21) as a linear combinations with non-negative coefficients of irreducible characters of representations on $G^* \in \widetilde{\mathcal{E}}_{\simp}(N)$, with repsect to the compatible family of functions $\mathcal{F}$. Thus (5.7.21) is of the form (4.3.26). So by lemma 4.3.6 all the coefficients must vanish. Thus in particular 
\begin{eqnarray}
& & \tr R_{\disc,\psi^N}^G(f) =0 , \,\ f \in \mathcal{H}(G) \\
& & R_P(w,\widetilde{\pi}_M,\psi_M) \equiv 1 \nonumber
\end{eqnarray}
and then going back to (5.7.20) we see that $\leftexp{0}{S^G_{\disc,\psi^N}(f)}=0$, as required.

\end{proof}

\bigskip

\begin{corollary} (of proof of proposition 5.7.4)
For $\psi^N$ as in proposition 5.7.4, and $G \in \widetilde{\mathcal{E}}_{\simp}(N)$ such that $\psi^N \in \Psi(G)$, we have the validity of part (a) of theorem 5.2.1. The global intertwining relation (namely part (b) of theorem 5.2.1) is valid for $\psi^N$ with respect to $G$, except for the cases (5.7.12) and (5.7.13), in which case we only have the weaker equality:
\begin{eqnarray}
\sum_{x \in \mathcal{S}_{\psi}} \epsilon_{\psi}^G(x)(f_G(\psi,x) - f^{\prime}_G(\psi,s_{\psi}x) )= 0.
\end{eqnarray}
\end{corollary}
\begin{proof}
We have already seen in the above proof that the hypothesis on $\psi^N$ allows us to establish both parts (a) and (b) of theorem 5.2.1 for $\psi^N$ with respect to any $G \in \widetilde{\mathcal{E}}_{\simp}(N)$ such that $\psi^N \in \Psi(G)$, except for the cases (5.7.12) and (5.7.13). However in these cases, we have established part (a) of theorem 5.2.1 in equation (5.7.22) (there is only one-nontrivial element in $W^0_{\psi}$ in these two cases). For (5.7.23) it follows from (5.7.14) together with the vanishing of the distributions $\tr R_{\disc,\psi_N}^G$ and $\leftexp{0}{S^G_{\disc,\psi_N}}$ just proved.
\end{proof}

\bigskip
\bigskip

In the previous analysis of this subsection we are mainly concerned with the case that $\psi^N \notin \widetilde{\Psi}_{\ellip}(N)$ is a non-elliptic parameter, and we saw by proposition 5.7.4 that we can obtain the stable multiplicity formula for $\psi^N$ with respect to $G \in \widetilde{\mathcal{E}}_{\simp}(N)$ in most cases (namely that if in addition $\psi^N \notin \Psi_{\ellip}(G)$ for all $G \in \widetilde{\mathcal{E}}_{\simp}(N)$). The treatment of the case where $\psi^N \in \widetilde{\Psi}_{\ellip}(N)$ is necessarily much more difficult. Here we only state a lemma for future use in the study of such parameters.

Thus let $G \in \widetilde{\mathcal{E}}_{\simp}(N)$ as before. Suppose that $\psi^N \in \widetilde{\Psi}_{\ellip}(N)$. We assume that the seed theorems 2.4.2 and 2.4.10 are valid for the simple generic constituents of $\psi^N$, and that the local classifications theorems are valid for the localizations of $\psi^N$ (in the case $\psi^N \in \Psi(G)$). Thus in particular, if $\psi^N \in \Psi(G)$, we can define the parameter $\psi \in \Psi_2(G)$ and stable linear form $f^G(\psi)$, and also the packet $\Pi_{\psi}$. 

\begin{lemma}
With hypothesis as above, suppose that the stable multiplicity formula is valid for $\psi^N$ with respect to $G \in \widetilde{\mathcal{E}}_{\simp}(N)$. Then the spectral multiplicity formula is also valid with respect to $G$.
\end{lemma}
\begin{proof}
The proof of the lemma is in fact implicit in the proof of proposition 5.7.4. However, if we are {\it given} the validity of the stable multiplicity formula, then the discussion is fairly simple.

For any proper Levi subgroup $M$ of $G$, the set $\Psi_2(M,\psi^N)$ is empty because $\psi^N \in \widetilde{\Psi}_{\ellip}(N)$. Hence by the discussion in section 5.4, we see that in the spectral expression for $I^G_{\disc,\psi^N}$, all the terms involving $M \neq G$ vanish. Thus 
\begin{eqnarray}
I^G_{\disc,\psi^N}(f) =  \tr R^G_{\disc,\psi^N}(f), \,\ f \in \mathcal{H}(G).
\end{eqnarray} 

On the other hand, since we are assuming the validity of the stable multiplicity formula for $\psi^N$ with respect to $G$, we have $\leftexp{0}{S^G_{\disc,\psi^N}(f)} =0$. Thus the endoscopic expansion (5.6.32) becomes:
\begin{eqnarray}
I^G_{\disc,\psi^N}(f) &=& \frac{m^G_{\psi^N}}{|\mathcal{S}_{\psi}|} \sum_{x \in \mathcal{S}_{\psi}} e_{\psi}(x) \epsilon^G_{\psi}(x) f^{\prime}_G(\psi,s_{\psi}x) \\
&=& \frac{m^G_{\psi^N}}{|\mathcal{S}_{\psi}|} \sum_{x \in \mathcal{S}_{\psi}}  \epsilon^G_{\psi}(x) f^{\prime}_G(\psi,s_{\psi} x) \nonumber
\end{eqnarray}
where the last equality follows from the fact, in case $m^G_{\psi^N} \neq 0$ which implies that $\psi^N \in \Psi_2(G)$, that the centralizer group $\overline{S}_{\psi}$ is finite, from which it follows that $e_{\psi}(x)=1$.  

From (5.7.25) we can already see that $I_{\disc,\psi^N}^G =0$ if $m^G_{\psi^N}=0$, i.e. when $\psi^N \notin \Psi(G)$, and so by (5.7.24) $\tr R^G_{\disc,\psi^N}(f)=0$. In other words $\psi^N$ does not contribute to the discrete spectrum of $G$. Thus we may now assume that $m^G_{\psi^N}=1$, i.e. $\psi^N \in \Psi_2(G)$.

For $x \in \mathcal{S}_{\psi} = \overline{S}_{\psi}$, let $(G^{\prime},\psi^{\prime})$ be the pair that corresponds to $(\psi,x)$. Then we have
\begin{eqnarray*}
f^{\prime}_G(\psi,x) = f^{\prime}(\psi^{\prime}) = \sum_{\pi \in \Pi_{\psi}} \langle s_{\psi} x,\pi \rangle f_G(\pi)
\end{eqnarray*}
and hence making the substitution $x \mapsto s_{\psi} x$ we obtain:
\begin{eqnarray}
f^{\prime}_G(\psi,s_{\psi} x) = \sum_{\pi \in \Pi_{\psi}} \langle x,\pi \rangle f_G(\pi)
\end{eqnarray}
for $x \in \mathcal{S}_{\psi}$.

\noindent Combining (5.7.24), (5.7.25) and (5.7.26) we obtain
\begin{eqnarray}
& & \tr R^G_{\disc,\psi^N}(f)= I^G_{\disc,\psi^N}(f) \\
&=& \frac{1}{|\mathcal{S}_{\psi}|} \sum_{x \in \mathcal{S}_{\psi}} \epsilon_{\psi}(x) \sum_{\pi \in \Pi_{\psi}} \langle x,\pi\rangle f_G(\pi) \nonumber \\
&=&  \sum_{\pi \in \Pi_{\psi}}  m(\pi) f_G(\pi)  \nonumber
\end{eqnarray}
where 
\[
m(\pi) = \frac{1}{|\mathcal{S}_{\psi}|} \sum_{\pi \in \Pi_{\psi}} \epsilon_{\psi}(x) \langle x,\pi \rangle 
\]
and this gives exactly the spectral multiplicity formula. 
\end{proof}

\subsection{The two sign lemmas}

Lemma 5.5.1 and 5.6.1 (the spectral and endoscopic sign lemmas respectively) play an important role in the term by term comparison of the spectral and endoscopic expansions of trace formulas, as seen in the previous subsections. The arguments given in section 4.6 of \cite{A1}, of the proofs of these two sign lemmas in the setting of symplectic and orthogonal groups, apply in the case of unitary groups also. However, it is necessary to first establish some preliminaries in our current setting of unitary groups, before we can refer to the proofs given in {\it loc. cit.}

As in the previous subsections, we denote by $G$ either an element of $\widetilde{\mathcal{E}}_{\simp}(N)$, or the twisted group $\widetilde{G}_{E/F}(N)$. Let $\psi \in \Psi(G)$. If $G=(G,\xi) \in \widetilde{\mathcal{E}}_{\simp}(N)$ denote $\psi^N = \xi_* \psi \in \widetilde{\Psi}(N)$; otherwise if $G=\widetilde{G}_{E/F}(N)$ then we let $\psi^N$ just to be another name for $\psi$. 

We refer to the notation of section 2.4. Thus we write
\begin{eqnarray*}
\psi^N &=& \bigboxplus_{k \in K_{\psi^N}} l_k \psi_k^{N_k} \\
&=& \bigboxplus_{i \in I_{\psi^N}} l_i \psi_i^{N_i} \boxplus \bigboxplus_{j \in J_{\psi^N}} l_j  (\psi_j^{N_j} \boxplus \psi_{j^*}^{N_{j^*}})
\end{eqnarray*}
here $k \leftrightarrow k^*$ is the involution on the set $K_{\psi^N}$ defined by the automorphism $\theta=\theta(N)$ of $G_{E/F}(N)$, with $I_{\psi^N}$ being the set of fixed points of $K_{\psi^N}$ under this involution, and $J_{\psi^N}$ a set of representatives of the order two orbits under this involution. We have
\[
K_{\psi^N} = I_{\psi^N} \coprod J_{\psi^N} \coprod (J_{\psi^N})^*
\]
and we identify $\{K_{\psi^N}\}$, the set of orbits of $K_{\psi^N}$ under the involution, as $I_{\psi^N} \coprod J_{\psi^N}$. In the context of section 5.5 and 5.6, we need to apply the spectral and endoscopic sign lemmas, only when $\psi^N \notin \widetilde{\Psi}_{\simp}(N)$. Thus we may assume that $\psi^N \notin \widetilde{\Psi}_{\simp}(N)$, and hence from the induction hypothesis, we can assume that theorem 2.4.2 and theorem 2.4.10 hold for the simple generic constituents of $\psi^N$. In particular the group $\mathcal{L}_{\psi}$ and the localization maps $L_{F_v} \rightarrow \mathcal{L}_{\psi}$ are defined for each place $v$ of $F$.

The first task is to establish a description of the sign character $\epsilon_{\psi}$ on $\mathcal{S}_{\psi} = \pi_0(\overline{S}_{\psi})$. As in {\it loc. cit.}, in the case of the twisted group $\widetilde{G}_{E/F}(N)$, the centralizer set $\overline{S}_{\psi}$ is a torsor under $\overline{S}^*_{\psi}$, with $\overline{S}^+_{\psi}$ being the (non-connected) group generated by $\overline{S}_{\psi}$; then $\epsilon_{\psi}$ is actually defined on the quotient $\mathcal{S}^+_{\psi} = \pi_0(\overline{S}^+_{\psi})$ of $\overline{S}_{\psi}^+$.

Recall that, as in section 2.5, we have the representation
\begin{eqnarray*}
& & \tau_{\psi}: \overline{S}_{\psi} \times \mathcal{L}_{\psi} \times \SL_2(\mathbf{C}) \rightarrow \GL(\widehat{\mathfrak{g}}) \\
& & \tau_{\psi}(s,g,h) = \Ad_G (s \cdot \widetilde{\psi}(g,h)  )
\end{eqnarray*}
and the sign character $\epsilon_{\psi}$ on $\mathcal{S}_{\psi}$ is defined in terms of the decomposition:
\begin{eqnarray}
\tau_{\psi} = \bigoplus_{\alpha} \tau_{\alpha} = \bigoplus_{\alpha} \lambda_{\alpha} \otimes \mu_{\alpha} \otimes \nu_{\alpha}.
\end{eqnarray} 

\noindent In order to give another description of the sign character $\epsilon_{\psi}$, we first give another decomposition of $\tau_{\psi}$. 

First recall that as in section 2.4, for each $k \in K_{\psi^N}$, we have the $L$-homomorphism associated to the simple generic constituent $\mu_k$ of $\psi_k^{N_k} = \mu_k \boxtimes \nu_k$ (here $\mu_k \in \widetilde{\Phi}_{\simp}(m_k)$ such that $N_k = m_k n_k$ with $n_k = \dim \nu_k$):
\begin{eqnarray}
\widetilde{\mu}_k : \mathcal{L}_{\psi} \longrightarrow \leftexp{L}{H_k} \hookrightarrow \leftexp{L}{G_{E/F}(m_k)}
\end{eqnarray}
(here $H_k \in \widetilde{\mathcal{E}}_{\simp}(m_k)$ or $\widetilde{G}^0_{E/F}(m_k)=G_{E/F}(m_k)$ itself). As in lemma 2.2.1, fix a choice of $w_c \in W_F \smallsetminus W_E$, then we can identify the $L$-homomorphism $\mathcal{L}_{\psi} \rightarrow \leftexp{L}{G_{E/F}(m_k)}$ of (5.8.2) as an irreducible representation
\begin{eqnarray}
\widetilde{\mu}_k\big|_{\mathcal{L}_{\psi/E}} : \mathcal{L}_{\psi/E} \longrightarrow \GL_{m_k}(\mathbf{C}).
\end{eqnarray}

\noindent For $g \in \mathcal{L}_{\psi/E}$ denote by $g^c \in \mathcal{L}_{\psi/E}$ the element $w_c g w_c^{-1}$; on the other hand, for $k \in K_{\psi^N}$, denote by $g_k$ the image of $g$ under (5.8.3). The for $k,k^{\prime} \in K_{\psi^N}$, denote by
\begin{eqnarray*}
R_{E,kk^{\prime}}(g) X = g_k \cdot X \cdot \leftexp{t}{g^c_{k^{\prime}}}
\end{eqnarray*}
the representation of $\mathcal{L}_{\psi/E}$ on the space of complex $m_k \times m_{k^{\prime}}$ matrices. Put
\begin{eqnarray}
R_{k k^{\prime}} = \Ind^{\mathcal{L}_{\psi}}_{\mathcal{L}_{\psi/E}} R_{E,k k ^{\prime}}
\end{eqnarray}
we have the relations:
\begin{eqnarray*}
& & R_{k k^{\prime}} \cong R_{k^{\prime} k} \\
& & R_{k k^{\prime}}^{\vee} \cong R_{k^* (k^{\prime})^{*}}
\end{eqnarray*}
and their associated $L$-function, defined with respect to the localization maps $L_{F_v} \rightarrow \mathcal{L}_{\psi}$, is the Rankin-Selberg $L$-function
\begin{eqnarray}
L(s,R_{kk^{\prime}}) = L(s,\mu_k \times \mu^c_{k^{\prime}}).
\end{eqnarray}

In the particular case where $k^{\prime}=k \in K_{\psi^N}$ we have the decomposition:
\begin{eqnarray}
R_{k k} = \Asai^+_k \oplus \Asai^-_k.
\end{eqnarray}
The representations $\Asai_k^{\pm}$ satisfy the relation
\begin{eqnarray*}
(\Asai_k^{\pm})^{\vee} \cong \Asai^{\pm}_{k^*}
\end{eqnarray*}
and whose associated $L$-functions are the Asai $L$-functions
\begin{eqnarray}
L(s,\Asai_k^{\pm}) = L(s,\mu_k,\Asai^{\pm}).
\end{eqnarray}

\noindent We refer to the family of representations $R_{k k^{\prime}}$, and $\Asai_k^{\pm}$ for $k,k^{\prime} \in K_{\psi^N}$ as {\it standard} representations of $\mathcal{L}_{\psi}$. For $\sigma$ a standard representation of $\mathcal{L}_{\psi}$ we have the global $L$-function $L(s,\sigma)$ with analytic continuation and functional equation:
\begin{eqnarray}
L(s,\sigma) = \epsilon(s,\sigma) L(1-s,\sigma^{\vee}).
\end{eqnarray}  
The epsilon factor $\epsilon(s,\sigma)$ is always being taken as automorphic $\epsilon$-factors. For $\sigma = R_{k k ^{\prime}}$ the automorphic $\epsilon$-factors are the same as the arithmetic $\epsilon$-factors, by the local Langlands classification for general linear groups \cite{HT,H1}.   

\noindent In addition, we refer to standard representations of the form:
\begin{eqnarray*}
& & \{ R_{kk}: \,\ k \in K_{\psi^N} \} \\ 
& & \{ R_{k k^*}: \,\ k \in K_{\psi^N}  \} \\
& & \{ \Asai_k^{\pm}: \,\ k \in K_{\psi^N}  \}
\end{eqnarray*}
as {\it diagonal} standard representations. We also note that the set of standard representations that are self-dual consists of the following:
\begin{eqnarray*}
& & \{R_{k k^{\prime}}: \,\  k,k^{\prime} \in I_{\psi^N} \} \\
& & \{ R_{k k^*}: \,\ k \in K_{\psi^N}     \} \\
& & \{ \Asai_k^{\pm}: \,\ k \in I_{\psi^N}       \}.
\end{eqnarray*}

\bigskip

\noindent For $\sigma$ a self-dual standard representation, it is either symplectic type or orthogonal type; more precisely, $\sigma$ is of symplectic type if $\sigma=R_{k k^{\prime}}$, where $k,k^{\prime} \in I_{\psi^N}$ that are conjugate self-dual of opposite parity (i.e. $R_{k k^{\prime}} = \Ind^{\mathcal{L}_{\psi}}_{\mathcal{L}_{\psi/E}} R_{E,k k^{\prime}}$ with $R_{E,k k^{\prime}}$ being conjugate symplectic), while the others are of orthogonal type. In particular we see that if $\sigma$ is a standard representation of symplectic type, then it must be irreducible of Rankin-Selberg type; hence in this case the arithmetic and automorphic $\epsilon$-factors defined by $\sigma$ coincides.

\begin{lemma}
There is a decomposition 
\begin{eqnarray}
\tau_{\psi} =\bigoplus_{\kappa} \tau_{\kappa} = \bigoplus_{\kappa} (\lambda_{\kappa} \otimes \sigma_{\kappa} \otimes \nu_{\kappa}), \,\ \kappa \in \mathcal{K}_{\psi^N}
\end{eqnarray}
for standard representations $\sigma_{\kappa}$ of $\mathcal{L}_{\psi}$, and irreducible representations $\lambda_{\kappa}$ and $\nu_{\kappa}$ of $\overline{S}_{\psi}$ and $\SL_2(\mathbf{C})$ respectively. The indexing set $\mathcal{K}_{\psi^N}$ has an involution
\[
\kappa \leftrightarrow \kappa^*
\] 
such that
\[
\tau_{\kappa}^{\vee} \cong \tau_{\kappa^*}.
\]
Furthermore, if a constituent $\sigma_{\kappa}$ is a diagonal standard representation, then the corresponding representation $\nu_{\kappa}$ is odd dimensional.
\end{lemma}
\begin{proof}
This can be proved as in lemma 4.6.1 of \cite{A1}. Namely that we first work with the product 
\[
\mathcal{A}_{\psi} = \mathcal{L}_{\psi} \times \SL_2(\mathbf{C}).
\]
One has the obvious variant of the definition of standard representation for $\mathcal{A}_{\psi}$. Then one can establish a decomposition
\[
\tau_{\psi} =\bigoplus_{\iota} \tau_{\iota} = \bigoplus_{\iota} \lambda_{\iota} \otimes \rho_{\iota}
\]
where $\lambda_{\iota}$ is irreducible representation of $\overline{S}_{\psi}$, and $\rho_{\iota}$ is a standard representation of $\mathcal{A}_{\psi}$, together with an involution $\iota \leftrightarrow \iota^*$ such that $\tau_{\iota}^{\vee} \cong \tau_{\iota^*}$. One then obtain the required decomposition (5.8.9), by observing that any standard representation $\rho_{\iota}$ of $\mathcal{A}_{\psi}$ has a decomposition of the form:
\[
\rho_{\iota} = \bigoplus_{\kappa } \sigma_{\kappa} \otimes \nu_{\kappa}
\]
for standard representations $\sigma_{\kappa}$ of $\mathcal{L}_{\psi}$ and irreducible representations $\nu_{\kappa}$ of $\SL_2(\mathbf{C})$. For instance if $\rho_{\iota}$ is of Rankin-Selberg type $R_{k k^{\prime}}$ (with respect to $\mathcal{A}_{\psi}$), then $\sigma_{\kappa}$ is the same Rankin-Selberg type representation for $\mathcal{L}_{\psi}$, while $\nu_{\kappa}$ ranges over the irreducible constituents of the tensor product:
\begin{eqnarray*}
\nu_{k} \otimes \nu_{k^{\prime}}.
\end{eqnarray*}
Similarly if $\rho_{\iota}$ is of the type $\Asai_k^{\eta}$ for $\eta = \pm 1$ (with respect to $\mathcal{A}_{\psi}$), then $\sigma_{\kappa}$ is again a standard representation $\Asai^{c \cdot \eta}_{k}$ for $\mathcal{L}_{\psi}$, with $c=\pm 1$ being determined by the condition that $c=+1$ if $\nu_k$ is odd dimensional (hence orthogonal), while $c=-1$ if $\nu_k$ is even dimensional (hence symplectic). Finally $\nu_{\kappa}$ ranges over appropriate irreducible consitutents of the tensor product:
\[
\nu_k \otimes \nu_k.
\]
Thus the decomposition (5.8.9) follows. For the last claim of the lemma, it suffices to note that in the case where a constituent $\sigma_{\kappa}$ is a diagonal standard representation, then the irreducible representation $\nu_{\kappa}$ of $\SL_2(\mathbf{C})$ must arise from the irreducible constituents of the tensor product $\nu_k \otimes \nu_k$ as above; such irreducible constituents must be odd dimensional.
\end{proof}

It is in the context of the decomposition (5.8.9), instead of the general decomposition (5.8.1), that we will establish the two sign lemmas. In particular, the same argument as in section 4.6 of \cite{A1} shows that we have the following expression for the sign character $\epsilon_{\psi}=\epsilon^G_{\psi}$:
\begin{eqnarray}
\epsilon_{\psi}(x)=\epsilon_{\psi}^G(x) = \prod_{\kappa \in \mathcal{K}_{\psi^N}^-} \det (\lambda_{\kappa}(s)), \,\ s \in \overline{S}_{\psi}
\end{eqnarray} 
where $x=x_s$ is the image of $s$ in $\mathcal{S}_{\psi}$, and the index set $\mathcal{K}_{\psi^N}^- \subset \mathcal{K}_{\psi^N}$ consists of those $\kappa \in \mathcal{K}_{\psi^N}$ such that
\begin{enumerate}

\item $\sigma_{\kappa}$ is symplectic ,

\item $\epsilon(1/2,\sigma_{\kappa})=-1$,

\item $\nu_{\kappa}(s_{\psi})=-1$; equivalently $\nu_{\kappa}$ is even dimentional, and hence symplectic.

\end{enumerate}

\bigskip

\noindent Following {\it loc. cit.} the next task is to relate the sign character $\epsilon_{\psi}^G$ with the corresponding sign character $\epsilon^M_{\psi_M}$ with respect to a Levi subgroup $M$ of $G^0$. Thus for a moment we consider the restriction
\[
\tau_{\psi,1} = \Ad_G \circ \widetilde{\psi}
\] 
of $\tau_{\psi}$ to the subgroup
\[
\mathcal{A}_{\psi} :=\mathcal{L}_{\psi} \times \SL_2(\mathbf{C}).
\]
As in the context of section 5.4 and 5.5, fix a Levi subgroup $M \subset G^0$, with dual Levi subgroup $\leftexp{L}{M} \subset \leftexp{L}{G^0}$, such that the set $\Psi_2(M,\psi^N)$ is non-empty, and let $\psi_M \in \Psi_2(M,\psi^N)$ that maps to $\psi$ under the $L$-embedding $\leftexp{L}{M} \subset \leftexp{L}{G^0}$. Denote by $\Ad_{G,M}$ the restriction of $\Ad_G$ to $\leftexp{L}{M}$. Thus
\[
\tau_{\psi,1} = \Ad_G \circ \widetilde{\psi} = \Ad_{G,M} \circ \widetilde{\psi}_M.
\]

\noindent We have the root space decomposition
\begin{eqnarray}
\widehat{\mathfrak{g}} = \bigoplus_{\xi} \widehat{\mathfrak{g}}_{\xi} = \widehat{\mathfrak{m}} \oplus \Big( \bigoplus_{\alpha \in \widehat{\Sigma}_M } \widehat{\mathfrak{g}}_{\alpha} \Big)
\end{eqnarray}
of $\widehat{\mathfrak{g}}$ with respect to the $\Gamma_F$-split component 
\[
A_{\widehat{M}} = (Z(\widehat{M})^{\Gamma_F})^0
\]
of the centre of $\widehat{M}$. Here $\widehat{\mathfrak{m}}$ is the Lie algebra of $\widehat{M}$, while $\xi$ ranges over the roots of $(\widehat{G}^0,A_{\widehat{M}})$ (including $0$), and $\widehat{\Sigma}_M$ denotes the set of non-zero roots. 

We then fix a decomposition
\[
\Ad_{G,M} \circ \widetilde{\psi}_M = \bigoplus_{h \in H_{\psi^N}} \sigma_h \otimes \nu_h
\]
that is compatible with the decompositions (5.8.9) and (5.8.11), in the following sense: the first condition means that we can write
\[
\sigma_h \otimes \nu_h = \sigma_{\kappa} \otimes \nu_{\kappa}, \,\ h \in H_{\psi^N}
\]
for a surjective mapping $h \rightarrow \kappa$ from the indexing set $H_{\psi^N}$ to $\mathcal{K}_{\psi^N}$. The second condition means that $H_{\psi^N}$ is a disjoint union of sets $H_{\psi^N,\xi}$, such that if $\rho_{\xi}$ is the restriction of $\Ad_{G,M}$ to the subspace $\widehat{\mathfrak{g}}_{\xi}$ of $\widehat{\mathfrak{g}}$, then
\begin{eqnarray}
\rho_{\xi} \circ \widetilde{\psi}_M = \bigoplus_{h \in H_{\psi^N,\xi}} \sigma_h \otimes \nu_h.
\end{eqnarray}
We denote by $\widehat{\mathfrak{g}}_h$ the subspace of $\widehat{\mathfrak{g}}_{\xi}$ on which the representation $\sigma_h \otimes \nu_h$ acts.

\bigskip

\noindent Let $H_{\psi^N,\xi}^-$ be the pre-image of $\mathcal{K}_{\psi^N}^-$ in $H_{\psi^N,\xi}$, and set
\[
\widehat{\mathfrak{g}}^-_{\psi^N,\xi} = \bigoplus_{h \in H^-_{\psi^N,\xi}} \widehat{\mathfrak{g}}_h.
\]
Similarly set $H_{\psi^N}^-$ to be the pre-image of $\mathcal{K}_{\psi^N}^-$ in $H_{\psi^N}$. Then we have
\[
H_{\psi^N}^- = \coprod_{\xi} H^-_{\psi^N,\xi}
\]
and we set
\begin{eqnarray}
\widehat{\mathfrak{g}}^-_{\psi^N} := \bigoplus_{\xi} \widehat{\mathfrak{g}}_{\psi^N,\xi}^- = \bigoplus_{h \in H^-_{\psi^N}} \widehat{\mathfrak{g}}_h
\end{eqnarray}
then the product $\overline{S}_{\psi} \times \mathcal{A}_{\psi} = \overline{S}_{\psi} \times \mathcal{L}_{\psi} \times \SL_2(\mathbf{C})$ acts on $\widehat{\mathfrak{g}}^-_{\psi^N}$, and there is a $\overline{S}_{\psi}^*$-equivariant morphism from $\overline{S}_{\psi}$ into the semisimple algebra $\End_{\mathcal{A}_{\psi}}(\widehat{\mathfrak{g}}^-_{\psi^N})$. Hence we can write (5.8.10) in the form:
\begin{eqnarray}
\epsilon^G_{\psi}(x_s) = \det (s, \End_{\mathcal{A}_{\psi}}(\widehat{\mathfrak{g}}^-_{\psi^N})), \,\ s \in \overline{S}_{\psi}. 
\end{eqnarray} 

\bigskip

To relate with the corresponding sign character on $M$, we set
\[
\widehat{\mathfrak{m}}^-_{\psi^N} := \widehat{\mathfrak{g}}^-_{\psi^N,0} \subset \widehat{\mathfrak{m}}.
\]

\noindent In the notation of section 5.1 let $\overline{N}_{\psi}$ be the normalizer of $\overline{T}_{\psi}$ in $\overline{S}_{\psi}$. Then for $n \in \overline{N}_{\psi}$ having image $u$ in $\mathfrak{N}_{\psi}$, write
\begin{eqnarray}
\epsilon^M_{\psi}(u) :=\det (n,\End_{\mathcal{A}_{\psi}}(\widehat{\mathfrak{m}}^-_{\psi^N})   )
\end{eqnarray} 
then the formula (5.8.15) matches the earlier definition of the canonical extension $\epsilon^M_{\psi_M}(\widetilde{u})$ of the character $\epsilon^M_{\psi_M}$ to $\widetilde{S}_{\psi_M,u}$. We can therefore write:
\[
\epsilon^M_{\psi}(u) = \epsilon^M_{\psi_M}(\widetilde{u}) = \epsilon^1_{\psi}(u).
\]

\noindent We also set, for $n \in \overline{N}_{\psi}$ having image $u \in \mathfrak{N}_{\psi}$ as above:
\begin{eqnarray}
\epsilon^{G/M}_{\psi}(u) := \det (n, \End_{\mathcal{A}_{\psi}}(\widehat{\mathfrak{g}}^-_{\psi^N}/ \widehat{\mathfrak{m}}^-_{\psi^N} ))
\end{eqnarray}
then $\epsilon^{G/M}_{\psi}(u)$ depends only on the image $w_u$ of $u$ in $W_{\psi}$ ({\it c.f.} the discussion after equation (4.6.6) of \cite{A1}). Thus we can write:
\begin{eqnarray}
\epsilon^G_{\psi}(x_u) = \epsilon^1_{\psi}(u) \epsilon^{G/M}_{\psi}(w_u), \,\ u \in \mathfrak{N}_{\psi}. 
\end{eqnarray}

\bigskip

\noindent The next step is to analyze the global normalizing factor 
\[
r_{\psi}(w) = r^G_{\psi}(w)=r_P(w,\psi_M)
\]
given by the value at $\lambda=0$ of the quotient:
\begin{eqnarray}
& & \\
& &
L(0,\rho_{P,w} \circ \widetilde{\psi}_{M,\lambda}) \epsilon(0,\rho_{P,w} \circ \widetilde{\psi}_{M,\lambda})^{-1} L(1, \rho_{P,w} \circ \widetilde{\psi}_{M,\lambda})^{-1}, \,\ \lambda \in \mathfrak{a}^*_{M,\mathbf{C}} \nonumber
\end{eqnarray}
here $\rho_{P,w} = \rho^{\vee}_{w^{-1}P|P}$, with $\rho_{w^{-1}P|P}$ being the adjoint representation of $\leftexp{L}{M}$ on
\[
w^{-1} \widehat{\mathfrak{n}}_P w/ w^{-1} \widehat{\mathfrak{n}}_P w \cap \widehat{\mathfrak{n}}_P.
\]
({\it c.f.} the paragraph before equation (5.4.9)). Here we are using the Artin notation for the $L$ and $\epsilon$-factors, but we remember that for the $\epsilon$-factor it is always given the automorphic definition; in fact as the analysis below shows (by virtue of lemma 5.8.1), we do not need to consider epsilon factors associated to Asai representations, and only need to consider Rankin-Selberg epsilon factors, for which both the Artin theoretic and the automorphic definitions coincide.

\bigskip
\noindent Denote by $\widehat{\Sigma}_P \subset \widehat{\Sigma}_M$ the set of roots of $(\widehat{P},A_{\widehat{M}})$. Put
\begin{eqnarray*}
\widehat{\Sigma}_{P,w} := \{\alpha \in \widehat{\Sigma}_P: \,\ w \alpha \notin \widehat{\Sigma}_P   \}.
\end{eqnarray*}  
Then we have
\[
\rho_{P,w} = \bigoplus_{\alpha \in \widehat{\Sigma}_{P,w}} \rho^{\vee}_{- \alpha} \cong \bigoplus_{\alpha \in \widehat{\Sigma}_{P,w}} \rho_{\alpha}
\]
since we have the isomorphism $\rho^{\vee}_{- \alpha} \cong \rho_{\alpha}$ determined by the Killing form on $\widehat{\mathfrak{g}}$. From the decomposition (5.8.12), we have
\begin{eqnarray}
& & L(s, \rho_{\alpha} \circ \widetilde{\psi}_{M,\lambda}) = L(s + \alpha(\lambda), \rho_{\alpha} \circ \widetilde{\psi}_M) \\
&=& \prod_{h \in H_{\psi^N,\alpha}} L(s+\alpha(\lambda), \sigma_h \otimes \nu_h   ).      \nonumber
\end{eqnarray}

\noindent In particular the $L$-functions $ L(s, \rho_{\alpha} \circ \widetilde{\psi}_{M,\lambda})$, and hence $L(s,\rho_{P,w} \circ \widetilde{\psi}_{M,\lambda})$, has analytic continuation and functional equation. As in {\it loc. cit.} we have the following expression for $r_{\psi}(w)$:
\begin{eqnarray*}
& & \\
 r_{\psi}(w) &=& \lim_{\lambda \rightarrow 0} L(1, \rho^{\vee}_{P,w} \circ \widetilde{\psi}_{M,\lambda}) L(1,\rho_{P,w} \circ \widetilde{\psi}_{M,\lambda} )^{-1} \nonumber\\
&=& \lim_{\lambda \rightarrow 0} \prod_{\alpha \in \widehat{\Sigma}_{P,w}} L(1, \rho_{\alpha}^{\vee} \circ \widetilde{\psi}_{M,\lambda}) L(1,\rho_{\alpha} \circ \widetilde{\psi}_{M,\lambda})^{-1} \nonumber \\
&=& \lim_{\lambda \rightarrow 0} \prod_{\alpha \in \widehat{\Sigma}_{P,w}} L(1- \alpha(\lambda), \rho_{\alpha}^{\vee} \circ \widetilde{\psi}_M   ) L(1+ \alpha(\lambda),\rho_{\alpha} \circ \widetilde{\psi}_M )^{-1}. \nonumber 
\end{eqnarray*}

\noindent Now the argument on p.38 of \cite{A9} shows that the two sets of of representations of $\mathcal{L}_{\psi}$ (up to isomorphism):
\begin{eqnarray*}
& & \{ \rho_{\alpha} \circ \widetilde{\psi}_M  \}_{\alpha \in \widehat{\Sigma}_{P,w}} \\
& & \{ \rho_{\alpha}^{\vee} \circ \widetilde{\psi}_M  \}_{\alpha \in \widehat{\Sigma}_{P,w}} 
\end{eqnarray*}
are in bijection. We hence obtain the following expression for $r_{\psi}(w)$:
\begin{eqnarray}
 r_{\psi}(w)= \lim_{\lambda \rightarrow 0}  \prod_{\alpha \in \widehat{\Sigma}_{P,w}} ( (- \alpha(\lambda))^{a_{-\alpha}} \alpha(\lambda)^{-a_{\alpha}} ) 
\end{eqnarray}
here we have denoted 
\[
a_{\alpha}  = \ord_{s=1} L(s,\rho_{\alpha},\widetilde{\psi}_M). 
\]
From (5.8.19), it follows that if we similarly denote 
\[
a_h= \ord_{s=1} L(s,\sigma_h \otimes \nu_h), \,\ h \in H_{\psi^N,\alpha}, 
\]
then we have
\begin{eqnarray*}
a_{\alpha} = \sum_{h \in H_{\psi^N,\alpha}} a_h.
\end{eqnarray*} 

\noindent We have the adjoint relation for the $L$-function of unitary representations (with respect to standard representations as above):
\begin{eqnarray}
L(s,(\sigma_h \otimes \nu_h)^{\vee}) =\overline{L(\overline{s},\sigma_h \otimes \nu_h)}
\end{eqnarray}
from which it follows that 
\[
a_h = \ord_{s=1} L(s, (\sigma_h \otimes \nu_h)^{\vee})
\]
and hence
\[
a_{-\alpha} = a_{\alpha}.
\]
Thus (5.8.20) becomes
\begin{eqnarray}
r_{\psi}(w)=\prod_{\alpha \in \widehat{\Sigma}_{P,w}} \prod_{h \in H_{\psi^N,\alpha}} (-1)^{a_h} 
\end{eqnarray}
In fact, it follows from the adjoint relation (5.8.21) that in the product of (5.8.22) it suffices to consider only those $h \in H_{\psi^N,\alpha}$ such that $\sigma_h$ is self-dual: $\sigma_h^{\vee} \cong \sigma_h$.

The analysis is now divided into the two cases depending on the parity of $n_h =\dim \nu_h$. We have:
\begin{eqnarray}
L(s,\sigma_h \otimes \nu_h) = \prod_{i=1}^{n_h}L(s+ \frac{1}{2}(n_h-2i+1)  ,\sigma_h ).
\end{eqnarray}  

\noindent Since the $L$-functions $L(s,\sigma_h)$ are associated to unitary cuspidal automorphic representations, with respect to standard representations, they are non-zero whenever $Re(s)\geq 1$ or $Re(s) \leq 0$; furthermore, the only possible poles are at $s=0$ and $s=1$. Hence on the right hand side of (5.8.23), the only terms that would contribute to $a_h=\ord_{s=1} L(s,\sigma_h \otimes \nu_h)$ are given by those $i$ such that
\[
1 + \frac{1}{2}(n_h -2i +1) = 0,\frac{1}{2},1.
\]

\noindent Suppose first that $n_h = \dim \nu_h$ is even, i.e. $\nu_h$ is symplectic. Then there is exactly one $i$ such that $1+ 1/2(n_h -2i+1) = 1/2$, and thus
\[
a_h=\ord_{s=1} L(s,\sigma_h \otimes \nu_h) = \ord_{s=1/2} L(s,\sigma_h). 
\]  
Since $\sigma_h$ is self-dual, it satisfies the functional equation
\begin{eqnarray}
L(s,\sigma_h) = \epsilon(s,\sigma_h) L(1-s,\sigma_h)
\end{eqnarray}
hence we have
\begin{eqnarray}
(-1)^{a_h} = \epsilon(1/2,\sigma_h).
\end{eqnarray}

Now since $\dim \nu_h$ is even, we have by lemma 5.8.1 that the standard representation $\sigma_h$ cannot be diagonal, thus $\sigma_h$ is self-dual of Rankin-Selberg type $R_{kk^{\prime}}$ with $k,k^{\prime} \in I_{\psi^N}$ such that $k^{\prime} \neq k, k^*$. 

It is at this point that we crucially apply the induction hypothesis for theorem 2.5.4(b). If $\sigma_h$ is orthogonal, then we can apply the induction hypothesis for theorem 2.5.4(b) to conclude that $\epsilon(1/2,\sigma_h)=1$, hence the only terms with $\dim \nu_h$ being even that contribute to (5.8.22) are given by those $\sigma_h$ that are symplectic, i.e. $h \in H_{\psi^N,\alpha}^-$. It equals the product
\begin{eqnarray}
r^-_{\psi}(w) = \prod_{\alpha \in \widehat{\Sigma}_{P,w}} (-1)^{|H^-_{\psi^N,\alpha}|}
\end{eqnarray}

\begin{rem}
\end{rem}
At this point it is necessary to remark on the use of induction hypothesis for theorem 2.5.4(b) in the above argument. Following remark 2.5.7, we see that, under the induction hypothesis, the following is valid with respect to our fixed integer $N$: suppose that for $i=1,2$ we have $k_i \in K_{\psi^N}$ corresponding to parameters $\psi_i^{N_i} \in \widetilde{\Psi}_{\simp}(N_i)$ with $\psi_i^{N_i} = \mu_{i} \boxtimes \nu_{i}$, such that $N_1 + N_2 <N$. Then if the pair $\mu_{1}$ and $\mu_{2}$ is conjugate self-dual of the same parity, i.e. $R_{k_1 k_2}$ is orthogonal, then
\[
\epsilon(1/2, \mu_{1} \times \mu_{2}^c) =1.
\] 

\bigskip

\noindent Recall that we allow $G$ to be an element in $\widetilde{\mathcal{E}}_{\simp}(N)$ or the twisted group $\widetilde{G}_{E/F}(N)$. First consider the case that $G$ is an element in $\widetilde{\mathcal{E}}_{\simp}(N)$. Then from the fact that $M \neq G^0$ and that $\Psi_2(M,\psi^N)$ is non-empty, we see that we must have $N_1 + N _2 <N$; hence the use of induction hypothesis for theorem 2.5.4(b) is valid. 

\noindent Next we consider the case that $G=\widetilde{G}_{E/F}(N)$. Again we have $M \neq G^0$, from which we have $N_1,N_2 < N$. However, in the case where $G=\widetilde{G}_{E/F}(N)$ we can certainly have $N_1 + N_2 = N$ (if $N_1+N_2 < N$ then we just use the induction hypothesis for theorem 2.5.4(b)). From the proof of lemma 5.8.1, we see that the set of irreducible representations $\nu_{\kappa}$ of $\SL_2(\mathbf{C})$ that is paired with the standard representation $R_{k_1 k_2}$ arises from the set of all irreducible constituents of the tensor product representation $\nu_1 \otimes \nu_2$. It follows that, in the case where $N_1+N_2=N$, we still obtain the same conclusion as before, {\it unless} $R_{k_1 k_2}$ is orthogonal, {\it and} that the tensor product $\nu_1 \otimes \nu_2$ decomposes into an {\it odd} number of irreducible representations of $\SL_2(\mathbf{C})$ of {\it even} dimension, i.e. unless $\psi^N=\psi_1^{N_1} \boxplus \psi_2^{N_2}$ is an $\epsilon$-parameter (as defined in the paragraph after equation (5.5.12)). It follows that in the case where $G=\widetilde{G}_{E/F}(N)$, then provided that $\psi^N$ is not an $\epsilon$-parameter, then the above reasoning leading to equation (5.8.26) is valid.

\noindent For future reference (more precisely for the proof of proposition 6.1.5), we also note that, in the case where $G=\widetilde{G}_{E/F}(N)$ and $\psi^N$ is an $\epsilon$-parameter as above, then the validity of the reasoning leading to (5.8.26) is equivalent to $\epsilon(1/2, \mu_1 \times \mu_2^c)=1$.

\bigskip

Now we consider the case that $n_h = \dim \nu_h$ is odd, i.e. $\nu_h$ is orthogonal. Again the functional equation (5.8.24) is satisfied since $\sigma_h$ is self-dual. Now if $n_h > 1$, then there is exactly one $i$ such that
\begin{eqnarray*}
& & 1+ \frac{1}{2}(n_h-2i+1) =1\\
& & 1+ \frac{1}{2}(n_h-2(i+1)+1) =0
\end{eqnarray*}
hence there are two terms on the right hand side of (5.8.23) that contribute to $a_h$. By the functional equation (5.8.24), we have
\[
\ord_{s=1} L(s,\sigma_h) = \ord_{s=0}L(s,\sigma_h)
\]
hence their contributions cancel in the product (5.8.22). Thus we are reduced to considering the case $n_h=1$, i.e. $\nu_h$ is the trivial one-dimensional representation of $\SL_2(\mathbf{C})$. Then $a_h=\ord_{s=1}L(s,\sigma_h)$. Again, we can apply the induction hypothesis for theorem 2.5.4(a) to obtain the value of $\ord_{s=1} L(s,\sigma_h)$: in the case where $\sigma_h$ is of the form $\Asai_k^{\pm}$ for $k \in I_{\psi^N}$, then theorem 2.5.4(a) asserts that $L(s,\sigma_h)$ has a pole at $s=1$, necessarily simple, if and only if $\sigma_h$ contains the trivial representation of $\mathcal{L}_{\psi}$; otherwise $L(s,\sigma_h)$ is analytic and non-zero at $s=1$; the same assertion is of course valid in the case where $\sigma_h$ is of Rankin-Selberg type, by the results of \cite{JPSS}. Thus the contribution of those terms on the right hand side of the product (5.8.22) with $\dim \nu_h$ odd are given by those $h$ such that $\sigma_h \otimes \nu_h$ containes the one-dimensional trivial representation of $\mathcal{A}_{\psi}=\mathcal{L}_{\psi} \times \SL_2(\mathbf{C})$. As in {\it loc. cit.} this total contribution (i.e. over $\alpha \in \widehat{\Sigma}_{P,w}$ and over $h \in H_{\psi^N,\alpha}$ such that $\sigma_h$ is self-dual and $\sigma_h \otimes \nu_h$ contains the one-dimensional trivial representation) is seen to be equal to the sign character $s^0_{\psi}(w)$ on $W_{\psi}$.

\noindent Thus to conclude, we have the factorization:
\begin{eqnarray}
r_{\psi}(w) = r_{\psi}^-(w) s_{\psi}^0(w).
\end{eqnarray}

\noindent From the factorization (5.8.17) and (5.8.27), it follows that the spectral sign lemma is equivalent to:
\begin{eqnarray}
r^-_{\psi}(w)  = \epsilon^{G/M}_{\psi}(w), \,\ w \in W_{\psi}.
\end{eqnarray}  
which by (5.8.26) is equivalent to proving the equality:
\begin{eqnarray}
\epsilon^{G/M}_{\psi}(w)   = \prod_{\alpha \in \widehat{\Sigma}_{P,w}} (-1)^{|H^-_{\psi^N,\alpha}|}
\end{eqnarray}

\noindent The argument for the proof of (5.8.29) is the same as that given in section 4.6 of \cite{A1}, from which we conclude the proof of the spectral sign lemma 5.5.1. Finally, the proof of the endoscopic sign lemma given in {\it loc. cit.} applies without change in our context, and so we also conclude the proof of the endoscopic sign lemma 5.6.1.

\section{\textbf{Study of Critical Cases}}

\subsection{The case of square-integrable parameters}

In this section we treat a class of parameters, which are critical for establishing the local theorems in section 7 and 8. These are not directly amenable to the induction arguments of the previous subsection, and we need the arguments as given in chapter 5 of \cite{A1}.

For this purpose we introduce a set of parameters:
\begin{eqnarray}
& & \widetilde{\mathcal{F}} = \coprod_N^{\infty} \widetilde{\mathcal{F}}(N) \\
& & \widetilde{\mathcal{F}}(N) \subset \widetilde{\Psi}(N). \nonumber
\end{eqnarray} 

We denote 
\begin{eqnarray*}
& & \widetilde{\mathcal{F}}_{\simp}(N) := \widetilde{\mathcal{F}} \cap \widetilde{\Psi}_{\simp}(N) \\
& & \widetilde{\mathcal{F}}_{\ellip}(N) := \widetilde{\mathcal{F}} \cap \widetilde{\Psi}_{\ellip}(N)
\end{eqnarray*}
etc, and we make the assumption that $\widetilde{\mathcal{F}}$ is the graded semi-group generated by $\widetilde{\mathcal{F}}_{\simp}(N)$ with respect to the operation $\boxplus$. We define $\deg \psi^N := N$ for $\psi^N \in \widetilde{\Psi}(N)$.

As in the previous section, we fix the integer $N$, and from the induction hypothesis that all the global theorems are valid for parameters in $\widetilde{\Psi}$ of degree less than $N$. We assume that the simple parameters of the family $\widetilde{\mathcal{F}}$ have degree less than or equal to $N$. For the arguments of section 6 we can work exclusively within the family $\widetilde{\mathcal{F}}$.

In order to carry out the arguments in this section we need to impose additional hypotheses on the family $\widetilde{\mathcal{F}}$. We first set up some notations:
\begin{eqnarray*}
& & \widetilde{\mathcal{F}}^{\prime}(N) := \widetilde{\mathcal{F}}(N) \smallsetminus \widetilde{\Phi}_{\simp}(N) \\
& & \widetilde{\mathcal{F}}_{\simp}^{\prime}(N) := \widetilde{\mathcal{F}}_{\simp}(N) \smallsetminus \widetilde{\Phi}_{\simp}(N) \\
& & \widetilde{\mathcal{F}}^{\prime}_{\ellip}(N) := \widetilde{\mathcal{F}}_{\ellip}(N) \smallsetminus \widetilde{\Phi}_{\simp}(N) 
\end{eqnarray*}
etc (recall that $\widetilde{\Phi}_{\simp}(N)$ is the subset of $\widetilde{\Psi}(N)$ consisting of simple generic parameters). 

For each $G = (G,\xi) \in \widetilde{\mathcal{E}}_{\simp}(N)$, put:
\begin{eqnarray}
\widetilde{\mathcal{F}}^{\prime}(G) :=\{\psi =(\psi^N,\widetilde{\psi}) \in \Psi(G)| \,\ \psi^N \in \widetilde{\mathcal{F}}^{\prime}(N)   \}
\end{eqnarray}
and similarly for the definition of $\widetilde{\mathcal{F}}^{\prime}_{\simp}(G), \widetilde{\mathcal{F}}^{\prime}_2(G), \widetilde{\mathcal{F}}^{\prime}_{\ellip}(G)$, etc. As before for $\psi^N \notin \widetilde{\Phi}_{\simp}(N)$, the induction hypothesis implies that the seed theorems 2.4.2 and 2.4.10 are valid for the simple generic constituents of $\psi^N$; thus $\widetilde{\psi}^N$ is defined and the condition $\psi =(\psi^N,\widetilde{\psi}) \in \Psi(G)$ is well-defined and is equivalent to $\widetilde{\psi}^N$ factoring through $\xi$ with $\widetilde{\psi}^N=\xi \circ \widetilde{\psi}$ ($\widetilde{\psi}$ is of course uniquely determined by $\widetilde{\psi}^N$).

We would like to define the set $\widetilde{\mathcal{F}}(G)$ by the same definition as in (6.1.2) except that we only require $\psi^N \in \widetilde{\mathcal{F}}(N)$. However, if $\psi^N \in \widetilde{\mathcal{F}}(N) \smallsetminus \widetilde{\mathcal{F}}^{\prime}(N)$, then $\psi^N$ is a simple generic parameter, and the seed theorems has yet to be shown to be valid for $\psi^N$. We will thus need to work with a provisional definition for the case of simple generic parameters.

On the other hand in the case where $G \in \widetilde{\mathcal{E}}_{\ellip}(N) \smallsetminus \widetilde{\mathcal{E}}_{\simp}(N)$ is composite, we can write
\[
G = G_1 \times G_2, \,\ G_i \in \widetilde{\mathcal{E}}_{\simp}(N_i), \,\ N_i < N.
\]
and we can define 
\[
\widetilde{\mathcal{F}}(G) := \widetilde{\mathcal{F}}(G_1) \times \widetilde{\mathcal{F}}(G_2) 
\]
with
\[
\widetilde{\mathcal{F}}(G_i) :=\{\psi =(\psi^{N_i},\widetilde{\psi}) \in \Psi(G_i)| \,\ \psi^{N_i} \in \widetilde{\mathcal{F}}(N_i)   \}
\]
being well-defined since the seed theorems are valid for simple generic parameters of degree less than $N$. Similarly for the definition of $\widetilde{\mathcal{F}}_2(G), \widetilde{\mathcal{F}}_{\ellip}(G)$ etc. in the case of composite $G$. 

Back to the case of simple $G \in \widetilde{\mathcal{E}}_{\simp}(N)$. We define the set
\[
\widetilde{\mathcal{F}}_{\gsimp}(G)
\]
to be consisting of pairs $\psi=(G,\psi^N)$ where $\psi^N \in \widetilde{\mathcal{F}}_{\gsimp}(N) := \widetilde{\mathcal{F}}(N) \cap \widetilde{\Phi}_{\simp}(N)$, satisfying the following two conditions:
\begin{enumerate}
\item The stable linear form $S^G_{\disc,\psi^N}$ is not zero.

\item There exists a stable linear form 
\[
f \mapsto f^G(\psi), \,\ f \in \mathcal{H}(G)
\]
\noindent such that 
\[
\widetilde{f}^G(\psi) = \widetilde{f}_N( \psi^N), \,\ \widetilde{f} \in \widetilde{\mathcal{H}}(N).
\]
\end{enumerate}

\bigskip

For $G=(G,\xi) \in \widetilde{\mathcal{E}}_{\simp}(N)$ we can formally define a map:
\begin{eqnarray}
& & \xi_* :  \widetilde{\mathcal{F}}_{\gsimp}(G) \rightarrow \widetilde{\mathcal{F}}_{\gsimp}(N) \\
& & \psi =(G,\psi^N) \mapsto \psi^N. \nonumber
\end{eqnarray}

\noindent Then we have 
\begin{eqnarray}
  \bigcup_{(G,\xi) \in \widetilde{\mathcal{E}}_{\simp}(N) } \xi_* \widetilde{\mathcal{F}}_{\gsimp}(G)  \subset \widetilde{\mathcal{F}}_{\gsimp}(N)
\end{eqnarray}
however at this point we cannot yet conclude that the inclusion (6.1.4) is an equality nor that the union in (6.1.4) is disjoint. 

We then simply define for $G \in \widetilde{\mathcal{E}}_{\simp}(N)$:
\begin{eqnarray}
& & \widetilde{\mathcal{F}}(G) := \widetilde{\mathcal{F}}^{\prime}(G) \bigsqcup \widetilde{\mathcal{F}}_{\gsimp}(G) \\
& & \widetilde{\mathcal{F}}_{\simp}(G) := \widetilde{\mathcal{F}}_{\simp}^{\prime}(G) \bigsqcup \widetilde{\mathcal{F}}_{\gsimp}(G) \nonumber \\
& & \widetilde{\mathcal{F}}_2(G) := \widetilde{\mathcal{F}}_2^{\prime}(G) \bigsqcup \widetilde{\mathcal{F}}_{\gsimp}(G) \nonumber \\
& & \widetilde{\mathcal{F}}_{\ellip}(G) := \widetilde{\mathcal{F}}^{\prime}_{\ellip}(G) \bigsqcup \widetilde{\mathcal{F}}_{\gsimp}(G) \nonumber
\end{eqnarray}
etc.

We can now make the following hypothesis on the set of parameters in $\widetilde{\mathcal{F}}$, which we impose in the rest of the section:
\begin{hypothesis} With notations as above:

\bigskip

\noindent (a) Suppose that $G=(G,\xi) \in \widetilde{\mathcal{E}}_{\ellip}(N)$, and that $\psi \in \widetilde{\mathcal{F}}_2(G)$. Then there exists a unique stable linear form
\[
f \mapsto f^G(\psi), \,\ f \in \mathcal{H}(G)
\]
satisfying:

\bigskip
\noindent (1) For $\widetilde{f} \in \widetilde{\mathcal{H}}(N)$ we have
\begin{eqnarray}
\widetilde{f}^G(\psi) = \widetilde{f}_N(\xi_* \psi).
\end{eqnarray}

\noindent (2) In case $G=G_1 \times G_2$ and $\psi = \psi_1 \times \psi_2$ are composite then we have
\begin{eqnarray}
f^G(\psi) = f_1^{G_1}(\psi_1) \times f_2^{G_2}(\psi_2)
\end{eqnarray}
if 
\[
f^G = f_1^{G_1} \times f_2^{G_2}.
\]

\bigskip

\noindent (b) The inclusion (6.1.4) is an equality (this part of course concerns only the simple generic parameters of degree $N$).
\end{hypothesis}

\bigskip

\begin{rem}
\end{rem}
\noindent (a) Part (a) of Hypothesis 6.1.1 is of course a global version of the statement of part (a) of theorem 3.2.1. The statement in hypothesis 6.1.1 holds also for the parameters $\psi \in \widetilde{\mathcal{F}}(G) \smallsetminus \widetilde{\mathcal{F}}_2(G)$ (in other words for $\psi^N \notin \widetilde{\mathcal{F}}_{\ellip}(N)$), by applying the induction hypothesis to a proper Levi subgroup of $G$.  

\bigskip

\noindent (b) In the case $G \in \widetilde{\mathcal{E}}_{\simp}(N)$, the definition of the set $\widetilde{\mathcal{F}}_{\gsimp}(G)$ of simple generic parameters of $G$ may look a bit unnatural. In fact, in the absence of theorem 2.4.2 for simple generic parameters of degree $N$, it is more reasonable to define the set
\begin{eqnarray*}
\widetilde{\mathcal{F}}^{\#}_{\gsimp}(G) 
\end{eqnarray*}
just to be consisting of pairs $\psi=(G,\psi^N)$ such that $\psi^N$ is simple generic and $S^G_{\disc,\psi^N}$ is not identifically zero. Then $\widetilde{\mathcal{F}}_{\gsimp}(G) \subset \widetilde{\mathcal{F}}_{\gsimp}^{\#}(G)$ but at this point we do not know {\it a priori} that this inclusion is an equality. On the other hand, we do have the equality
\begin{eqnarray*}
\bigcup_{(G,\xi) \in \widetilde{\mathcal{E}}_{\simp}(N)} \xi_* \widetilde{\mathcal{F}}^{\#}_{\gsimp}(G) = \widetilde{\mathcal{F}}_{\gsimp}(N)
\end{eqnarray*}
(but at this point we do not know {\it a priori} that the union is disjoint). This can be seen as follows: if $\psi^N$ is simple generic, (5.6.2) just becomes
\begin{eqnarray*}
\widetilde{I}^N_{\disc,\psi^N}(\widetilde{f}) = \sum_{G \in \widetilde{\mathcal{E}}_{\simp}(N)} \widetilde{\iota}(N,G) \cdot \widehat{S}^G_{\disc,\psi^N}(\widetilde{f}^G), \,\ \widetilde{f} \in \widetilde{\mathcal{H}}(N)  
\end{eqnarray*} 
indeed the terms $S^G_{\disc,\psi^N}$ for $G \in \widetilde{\mathcal{E}}_{\ellip}(N) \smallsetminus \widetilde{\mathcal{E}}_{\simp}(N)$ all vanish (as follows from proposition 4.3.4). Hence we conclude that the distribution $S^G_{\disc,\psi^N}$ is non-zero for some $G = (G,\xi)\in \widetilde{\mathcal{E}}_{\simp}(N)$, i.e. $\psi^N \in \xi_* \widetilde{\mathcal{F}}^{\#}_{\gsimp}(G)$. 

\noindent We thus similarly define
\begin{eqnarray*}
& & \widetilde{\mathcal{F}}^{\#}(G) := \widetilde{\mathcal{F}}^{\prime}(G) \bigsqcup \widetilde{\mathcal{F}}^{\#}_{\gsimp}(G) \\
& & \widetilde{\mathcal{F}}^{\#}_{\simp}(G) := \widetilde{\mathcal{F}}^{\prime}_{\simp}(G) \bigsqcup \widetilde{\mathcal{F}}^{\#}_{\gsimp}(G) \\
& & \widetilde{\mathcal{F}}^{\#}_2(G) := \widetilde{\mathcal{F}}^{\prime}_2(G) \bigsqcup \widetilde{\mathcal{F}}^{\#}_{\gsimp}(G) \\
\end{eqnarray*}
etc.

\noindent We also note that if we can establish part (a) of hypothesis 6.1.1, but for all $\psi \in \widetilde{\mathcal{F}}^{\#}_{2}(G)$ (the difference of course concerns only the simple generic parameters), then we would have $\widetilde{\mathcal{F}}_{\gsimp}(G) = \widetilde{\mathcal{F}}_{\gsimp}^{\#}(G)$. In particular (6.1.4) is then an equality.

\bigskip

As in the arguments in section 5, comparison with the twisted trace formula for $\widetilde{G}_{E/F}(N)$ plays a crucial role. Thus consider $\psi^N \in \widetilde{\mathcal{F}}_{\ellip}(N)$ an elliptic parameter. Recall the notion of an $\epsilon$-parameter (defined in the paragraph before the statement of the spectral sign lemma 5.5.1). The case of an $\epsilon$-parameter will be treated below, thus we first consider the case $\psi^N$ is not an $\epsilon$-parameter. We would like to use the spectral and endoscopic expansions of proposition 5.5.2 and 5.6.2 respectively. Strictly speaking these expansions are derived under the tacit assumption that $\psi^N$ is not a simple generic parameter (in order to apply the induction hypothesis to have the validity of the seed theorems for the simple generic constituents of $\psi^N$). However, with the presence of hypothesis 6.1.1 in force, the case of simple generic parameters can also be treated in the same (actually simpler) manner. 

Thus we first assume that $\psi^N \in \widetilde{\mathcal{F}}_{\ellip}(N) \smallsetminus \widetilde{\mathcal{F}}_{\gsimp}(N)$ is not a simple generic parameter. By hypothesis 6.1.1, we see that the hypothesis for the endoscopic expansion (5.6.32) of proposition 5.6.2 is satisfied (in fact for the endoscopic expansion we do not need to make the assumption that $\psi^N$ is not an $\epsilon$-parameter). In the case of $\widetilde{G}_{E/F}(N)$ and $\psi^N \in \widetilde{\Psi}_{\ellip}(N)$ the twisted centralizer $S_{\psi}$ is a connected abelian torsor; in particular $\mathcal{S}_{\psi^N}$ is reduced to a singleton $\{x_1\}$ and $\mathfrak{N}_{\psi^N} \cong W_{\psi^N} \cong \mathcal{S}_{\psi^N}=\{x_1\}$. Thus the expansion (5.6.32) just becomes:
\begin{eqnarray}
\widetilde{I}^N_{\disc,\psi^N} (\widetilde{f}) - \leftexp{0}{\widetilde{s}^N_{\disc,\psi^N}}(\widetilde{f}) = \frac{1}{2}  e_{\psi^N}(x_1) \widetilde{\epsilon}^N_{\psi^N}(x_1) \widetilde{f}^{\prime}_N(\psi, x_1)
\end{eqnarray}
(the element $s_{\psi^N}$ does not matter in this case since $\overline{S}_{\psi^N}$ is connected).

\noindent Similarly the spectral expansion (5.5.20) is just:
\begin{eqnarray}
\widetilde{I}^N_{\disc,\psi^N}(\widetilde{f}) - \leftexp{0}{\widetilde{r}^N_{\disc,\psi^N} (\widetilde{f})} = \frac{1}{2}  i_{\psi^N}(x_1) \widetilde{\epsilon}^N_{\psi^N}(x_1) \widetilde{f}_N(\psi, x_1).
\end{eqnarray}

We note the global intertwining relation for $\psi^N$, namely that in this case it reduces to (6.1.6) of hypothesis 6.1.1 (the situation is similar to the local case in section 3.5, {\it c.f.} (3.5.20)-(3.5.21)). Thus
\begin{eqnarray}
\widetilde{f}_N^{\prime}(\psi^N, x_1 ) = \widetilde{f}_N(\psi^N, x_1), \,\ \widetilde{f} \in \widetilde{\mathcal{H}}(N).
\end{eqnarray}

The expansions (6.1.8) and (6.1.9), together with (6.1.10) and the equality $i_{\psi^N}(x_1)=e_{\psi^N}(x_1)$ (equation (5.1.8)) gives:
\begin{eqnarray}
\leftexp{0}{\widetilde{r}^N_{\disc,\psi^N} (\widetilde{f})} = \leftexp{0}{\widetilde{s}^N_{\disc,\psi^N} (\widetilde{f})}.
\end{eqnarray}

For the twisted group $\widetilde{G}_{E/F}(N)$ the distribution $\leftexp{0}{\widetilde{r}^N_{\disc,\psi^N}(\widetilde{f})}$ of course vanishes, by the theorem of Moegin-Waldspurger \cite{MW} and Jacquet-Shalika \cite{JS}. Hence $\leftexp{0}{\widetilde{s}^N_{\disc,\psi^N}(\widetilde{f})}$ also vanishes. As in the the arguments in section 4 and 5, we can work interchangeably with $\widetilde{f} \in \widetilde{\mathcal{H}}(N)$ and a compatible family of functions $\mathcal{F} = \{f \in \mathcal{H}(G^*)| \,\ G^* \in \widetilde{\mathcal{E}}_{\ellip}(N)   \}$, and we have:

\begin{eqnarray}
\sum_{G^* \widetilde{\mathcal{E}}_{\simp}(N)}  \widetilde{\iota}(N,G^*) \leftexp{0}{ S^{G^*}_{\disc,\psi^N}(f)}  =0, \,\ f \in \mathcal{F}.
\end{eqnarray}

\bigskip

Next consider $G=(G,\xi) \in \widetilde{\mathcal{E}}_{\simp}(N)$, and $\psi^N$ is a simple generic parameter such that $\psi^N \in \xi_* \widetilde{\mathcal{F}}_{\gsimp}(G)$ (recall that at this point we can still not assert that the inclusion (6.1.4) is an equality). Let $\psi=(G,\psi^N) \in \widetilde{\mathcal{F}}_{\gsimp}(G)$ be the parameter defined by $\psi^N$. In this case we simply put 
\[
\mathcal{S}_{\psi}=\mathcal{S}_{\psi}(G) = \{1\}
\]
and by hypothesis 6.1.1, we have the existence of the stable linear form $f^G(\psi)$, which allows us to define $\leftexp{0}{S^G_{\disc,\psi^N}}$ as in (5.6.28), which simply becomes:
\begin{eqnarray}
\leftexp{0}{S^G_{\disc,\psi^N}(f)} = S^G_{\disc,\psi^N}(f) - f^G(\psi), \,\ f \in \mathcal{H}(G).
\end{eqnarray}

\noindent On the other hand we have another element $G^{\vee} \in \widetilde{\mathcal{E}}_{\simp}(N)$, and we simply put
\begin{eqnarray}
\leftexp{0}{S^{G^{\vee}}_{\disc,\psi^N}(f)} = S^{G^{\vee}}_{\disc,\psi^N}(f), \,\ f \in \mathcal{H}(G^{\vee}).
\end{eqnarray} 

\noindent This is to be expected since $G \in \widetilde{\mathcal{E}}_{\simp}(N)$ is the unique element through which $\psi^N$ factors, and thus $\psi^N$ should make no contribution to the distribution $S^{G^{\vee}}_{\disc,\psi^N}$. Then with this convention the identity (6.1.12) is still valid; its derivation is in fact much simpler as follows: from the fact that $\psi^N$ is a simple parameter, the terms on the right hand side of (5.4.2) (with $G=\widetilde{G}_{E/F}(N)$ in the notation there) have to vanish, again by the theorems of \cite{MW} and \cite{JS}. Thus 
\[
\widetilde{I}^N_{\disc,\psi^N}(\widetilde{f}) = \widetilde{r}^N_{\disc,\psi^N}(\widetilde{f}) = \frac{1}{2} \tr \widetilde{R}^N_{\disc,\psi^N}(\widetilde{f})=\frac{1}{2} \widetilde{f}_N(\psi^N).
\]

Similarly for the expansion (5.6.2) (again with $G=\widetilde{G}_{E/F}(N)$ in the notation there), all the terms associated to $G^* \in \widetilde{\mathcal{E}}_{\ellip}(N) \smallsetminus \widetilde{\mathcal{E}}_{\simp}(N)$ vanish, by virtue of the stable multiplicity formula applied to such $G^*$ (which follows from the induction hypothesis as such $G^*$ is composite). Thus 
\begin{eqnarray*}
& & \widetilde{I}^N_{\disc,\psi^N}(\widetilde{f}) = \widetilde{s}^N_{\disc,\psi^N}(\widetilde{f}) =\sum_{G^* \in \widetilde{\mathcal{E}}_{\simp}(N)} \widetilde{\iota}(N,G^*) \widehat{S}^{G^*}_{\disc,\psi^N}(\widetilde{f}^{G^*}) \\
&=& \widetilde{\iota}(N,G) \widetilde{f}^G(\psi) + \sum_{G^* \in \widetilde{\mathcal{E}}_{\simp}(N)} \widetilde{\iota}(N,G^*) \leftexp{0}{\widehat{S}^{G^*}_{\disc,\psi^N} (\widetilde{f}^{G^*})} \\
&=& \frac{1}{2}  \widetilde{f}^G(\psi) + \sum_{G^* \in \widetilde{\mathcal{E}}_{\simp}(N)} \widetilde{\iota}(N,G^*) \leftexp{0}{\widehat{S}^{G^*}_{\disc,\psi^N} (\widetilde{f}^{G^*})}
\end{eqnarray*}
(recall that here we are considering a fixed $G \in \widetilde{\mathcal{E}}_{\simp}(N)$ such that $\psi^N \in \xi_* \widetilde{\mathcal{F}}_{\gsimp}(G)$).

\noindent Combining the two expansions and again using the equality $\widetilde{f}^G(\psi) = \widetilde{f}_N(\psi^N)$ from (6.1.6), we again obtain (6.1.12), upon replacing $\widetilde{f} \in \widetilde{\mathcal{H}}(N)$ by a compatible family $\mathcal{F}$.

Thus to conclude the identity (6.1.12) is valid for $\psi^N \in \xi_* \widetilde{\mathcal{F}}_2(G)$ for any $G \in \widetilde{\mathcal{E}}_{\ellip}(N)$ with $\psi^N$ not an $\epsilon$-parameter.

\begin{proposition}
Suppose that $G=(G,\xi) \in \widetilde{\mathcal{E}}_{\ellip}(N)$, and $\psi^N \in \xi_* \widetilde{\mathcal{F}}_2(G)$ not an $\epsilon$-parameter, with $\psi^N = \xi_* \psi$ for $\psi \in \widetilde{\mathcal{F}}_2(G)$. Then we have: 

\noindent (1) If $G\in \widetilde{\mathcal{E}}_{\simp}(N)$, then the stable multiplicity formula is valid for the pair $(G,\psi)$ if and only if we have
\begin{eqnarray}
\leftexp{0}{S^{G^*}_{\disc,\psi^N}} = S^{G^*}_{\disc,\psi^N} =0
\end{eqnarray}
for $G^* \in \widetilde{\mathcal{E}}_{\simp}(N)$ with $G^* \neq G$. In the case $N$ is odd (6.1.15) holds, and hence the stable multiplicity formula is valid for $(G,\psi)$

\noindent (2) In the case where $G \notin \widetilde{\mathcal{E}}_{\simp}(N)$ we also have:
\[
\leftexp{0}{S_{\disc,\psi^N}^{G^*}} = S_{\disc,\psi^N}^{G^*} \equiv 0,\,\  R^{G^*}_{\disc,\psi^N} \equiv 0
\] 
for every $G^* \in \widetilde{\mathcal{E}}_{\simp}(N)$.
\end{proposition}
\begin{proof}
First we consider the case that $G \in \widetilde{\mathcal{E}}_{\simp}(N)$. Then by (6.1.12) we have
\begin{eqnarray}
\sum_{G^* \widetilde{\mathcal{E}}_{\simp}(N)}  \widetilde{\iota}(N,G^*) \leftexp{0}{ S^{G^*}_{\disc,\psi^N}(f_{G^*})}  =0
\end{eqnarray} 
for any compatible family of functions $\mathcal{F}=\{f_{G^*} \in \mathcal{H}(G^*)| \,\  G^* \in \widetilde{\mathcal{E}}_{\simp}(N) \}$. Recall that besides $G$, there is only one other element $G^{\vee} \in \widetilde{\mathcal{E}}_{\simp}(N)$. Thus if (6.1.15) holds, we have
\[
\leftexp{0}{S^G_{\disc,\psi^N}}(\widetilde{f}^G) =0
\]
for any $\widetilde{f} \in \widetilde{\mathcal{H}}(N)$. By the surjectivity of the Kottwitz-Shelstad transfer $\widetilde{f} \mapsto \widetilde{f}^G$ (Proposition 3.1.1(b)), it follows that $\leftexp{0}{S^G_{\disc,\psi^N}} \equiv 0$, i.e. the stable multiplicity formula is valid for the pair $(G,\psi)$. The converse is similar. 

In the case $N$ is odd, then since the two data $G$ and $G^{\vee}$ has no Levi sub-data in common ({\it c.f.} remark 2.4.1), we can choose the compatible family $\{f_{G^*}\}$ so that $f_{G^{\vee}} \equiv 0$, while $f_G$ being arbitrary. Hence from (6.1.16) we obtain the vanishing of $S^G_{\disc,\psi^N}$. Similarly for the vanishing of $S^{G^{\vee}}_{\disc,\psi^N}$.

\bigskip

Next we consider the case that $G=(G,\xi) \in \widetilde{\mathcal{E}}_{\ellip}(N) \smallsetminus \widetilde{\mathcal{E}}_{\simp}(N)$, and that $\psi^N \in \xi_* \widetilde{\mathcal{F}}_2(G)$ is not an $\epsilon$-parameter (note that in this case $\psi^N$ cannot be a simple parameter). The since $G$ is composite, the stable multiplicity formula is valid for the pair $(G,\psi)$ (which follows from our induction hypothesis). 

Since $\psi^N \notin  \widetilde{\Psi}(G^*)$ for $G^* =(G^*,\xi^*) \in \widetilde{\mathcal{E}}_{\simp}(N)$, we have by part (a) of proposition 5.7.1:
\begin{eqnarray}
& & \\
& & \tr R^{G^*}_{\disc,\psi^N}(f_{G^*}) = \leftexp{0}{S^{G^*}_{\disc,\psi^N}(f_{G^*})} = S^{G^*}_{\disc,\psi^N}(f_{G^*}), \,\ f_{G^*} \in \mathcal{H}(G^*). \nonumber
\end{eqnarray}

\noindent By the assumption that $\psi^N$ is not an $\epsilon$-parameter, the identity (6.1.16) is again valid, hence substituting (6.1.18) in (6.1.16) we have:
\begin{eqnarray}
\sum_{G^* \widetilde{\mathcal{E}}_{\simp}(N)}  \widetilde{\iota}(N,G^*)  \tr R^{G^*}_{\disc,\psi^N}(f_{G^*}) =0.
\end{eqnarray}
By the now familiar use of lemma 4.3.6, we deduce 
\[
\tr R^{G^*}_{\disc,\psi^N} \equiv 0
\]
for every $G^* \in \widetilde{\mathcal{E}}_{\simp}(N)$. Hence we deduce the rest of the proposition from (6.1.18).
\end{proof}

\bigskip

The case where $G \in \widetilde{\mathcal{E}}_{\simp}(N)$ with $N$ being even is much more subtle, and the full resolution is only completed in chapter 9. We first setup some preliminary considerations. Thus for $N$ even, let $L \cong G_{E/F}(N/2)$ be the Siegel Levi of $U_{E/F}(N)$, equipped as a Levi sub-datum of $G$. Similarly for for the other element $G^{\vee}=(G^{\vee},\xi^{\vee}) \in \widetilde{\mathcal{E}}_{\simp}(N)$ denote by $L^{\vee}$ the Levi sub-datum of $G^{\vee}$ with the same underlying group as $L$. The datum $L$ and $L^{\vee}$ are equivalent, {\it c.f.} remark 2.4.1. Define two transfer mappings:
\begin{eqnarray*}
 f &\mapsto& f^L=f_L, \,\ f \in \mathcal{H}(G) = \mathcal{H}(U_{E/F}(N))\\
 f^{\vee} &\mapsto& f^{\vee,L^{\vee}}=f^{\vee}_{L^{\vee}}, \,\ f^{\vee} \in \mathcal{H}(G^{\vee}) = \mathcal{H}(U_{E/F}(N))
\end{eqnarray*}  
to the space $\mathcal{S}(L) =\mathcal{I}(L)$ (resp. to the space $\mathcal{S}(L^{\vee})=\mathcal{I}(L^{\vee})$; of course $L$ and $L^{\vee}$ have the same underlying group and the two transfer mappings are the same, but it is still important to distinguish them notationally). Denote by $\widetilde{\mathcal{S}}^0(L)$ (resp. $\widetilde{\mathcal{S}}^0(L^{\vee})$) the image of the transfer mapping.

\begin{proposition}
Suppose that $G$ and $\psi$ are as in proposition 6.1.3, with $G \in \widetilde{\mathcal{E}}_{\simp}(N)$ and $N$ being even. Then there exist linear forms
\begin{eqnarray*}
h^L &\mapsto & h^L(\Lambda), \,\ h^L \in \widetilde{\mathcal{S}}^0(L) \\
h^{\vee,L^{\vee}} & \mapsto & h^{\vee,L^{\vee}}(\Lambda^{\vee}), \,\ h^{\vee,L^{\vee}} \in \widetilde{\mathcal{S}}^0(L)
\end{eqnarray*}
on $\widetilde{\mathcal{S}}^0(L)$ and $\widetilde{\mathcal{S}}^0(L^{\vee})$, such that
\[
S^G_{\disc,\psi^N}(f) = |\mathcal{S}_{\psi}|^{-1} \epsilon^G(\psi) f^G(\psi) - f^L(\Lambda), \,\ f \in \mathcal{H}(G)
\]
and
\[
S^{G^{\vee}}_{\disc,\psi^N}(f^{\vee}) = f^{\vee,L^{\vee}}(\Lambda^{\vee}), \,\ f^{\vee} \in \mathcal{H}(G^{\vee}). 
\]

\noindent Furthermore, we have 
\[
f^L(\Lambda)=f^{\vee,L^{\vee}}(\Lambda^{\vee})
\]
with $f$ (resp. $f^{\vee}$) the the function associated to $G$ (resp. $G^{\vee}$) in a compatible family, and the linear form
\[
f^{\vee} \mapsto f^{\vee,L^{\vee}}(\Lambda^{\vee}), \,\ f^{\vee} \in \mathcal{H}(G^{\vee})
\]
is a unitary character on $G^{\vee}(\mathbf{A}_F)=U_{E/F}(N)(\mathbf{A}_F)$.
\end{proposition}
\begin{proof}
We use (6.1.12), which simplifies to (using $\widetilde{\iota}(N,G)=\widetilde{\iota}(N,G^{\vee})=1/2$)
\[
\leftexp{0}{S^G_{\disc,\psi^N}(f)} + \leftexp{0}{S^{G^{\vee}}_{\disc,\psi^N}(f^{\vee})} =0
\]
where $f$ and $f^{\vee}$ are functions associated to $G$, resp. $G^{\vee}$, in the compatible family occuring in (6.1.12). The only condition between $f$ and $f^{\vee}$ imposed by the compatible family is given by (3.1.5). From this it follows that $\leftexp{0}{S^G_{\disc,\psi^N}(f)}$, resp. $\leftexp{0}{S^{G^{\vee}}_{\disc,\psi^N}(f)}$ factors through a linear form on $\widetilde{\mathcal{S}}^0(L)$, resp. $\widetilde{\mathcal{S}}^0(L^{\vee})$. More precisely we have linear forms $\Lambda$, resp. $\Lambda^{\vee}$ on $\widetilde{\mathcal{S}}^0(L)$, resp. $\widetilde{\mathcal{S}}^0(L^{\vee})$, such that: 
\begin{eqnarray*}
& & \leftexp{0}{S^{G}_{\disc,\psi^N}(f)} = - f^L(\Lambda), \,\ f \in \mathcal{H}(G) \\
& & \leftexp{0}{S^{G^{\vee}}_{\disc,\psi^N}(f^{\vee})} = f^{\vee,L^{\vee}}(\Lambda^{\vee}), \,\ f^{\vee} \in \mathcal{H}(G^{\vee})
\end{eqnarray*}
and satisfying
\[
f^L(\Lambda) = f^{\vee,L^{\vee}}(\Lambda^{\vee})
\]
when $f$ resp. $f^{\vee}$ are functions associated to $G$ resp. $G^{\vee}$ that is part of a compatible family. 

For the final assertion we use the equality from (5.7.3):
\[
\tr R^{G^{\vee}}_{\disc,\psi^N}(f^{\vee}) = \leftexp{0}{S^{G^{\vee}}_{\disc,\psi^N}}(f^{\vee}),\,\ f^{\vee} \in \mathcal{H}(G^{\vee})
\]
since $\leftexp{0}{S^{G^{\vee}}_{\disc,\psi^N}}(f^{\vee}) = f^{\vee,L^{\vee}}(\Lambda^{\vee})$, it follows that $\Lambda^{\vee}$ is a unitary character (possibly zero). 

\end{proof}

\bigskip
Thus the linear form $\Lambda$ (and $\Lambda^{\vee}$) is the obstruction to the validity of the stable multiplicity formula in the situation of proposition 6.1.4. The full resolution, namely the proof of the vanishing of $\Lambda$, is achieved in chapter 9. For generic parameters with certain local constraints at archimedean places, we will establish the vanishing of the linear form $\Lambda$ for these parameters in section 6.4.

Note that in the situation proposition 6.1.3, the stable multiplicity formula for such a parameter $\psi^N$ is already shown to be valid in the case when $N$ is odd; we simply interpret this as the vanishing of the linear form $\Lambda$.
\bigskip

We now treat the case where $\psi^N$ is an $\epsilon$-parameter.

\begin{proposition}
Suppose that $\psi^N =\psi_1^{N_1} \boxplus \psi_2^{N_2} \in \widetilde{\mathcal{F}}_{\ellip}(N)$ is an $\epsilon$-parameter. Then for any $G^* \in \widetilde{\mathcal{E}}_{\simp}(N)$ we have:
\begin{eqnarray}
& & \\
& & \tr R^{G^*}_{\disc,\psi^N}(f) =0= \leftexp{0}{S^{G^*}_{\disc,\psi^N}(f)} = S^{G^*}_{\disc,\psi^N}(f),\,\  f \in \mathcal{H}(G^*). \nonumber
\end{eqnarray}
Furthermore, for the simple generic constituents $\mu_1$ and $\mu_2$ of $\psi_1^{N_1}$ and $\psi_2^{N_2}$ respectively, we have
\begin{eqnarray}
\epsilon(1/2,\mu_1 \times \mu_2^c)=1
\end{eqnarray}
in accordance with part (b) of theorem 2.5.4.
\end{proposition}
\begin{proof}
The main issue is that, in the case where $\psi^N$ is an $\epsilon$-parameter, we do not have the identity (6.1.12) {\it a priori}; this is because the derivation of (6.1.12) relies on the spectral expansion (6.1.9) (namely the expansion (5.5.20) in the present setting). As seen in section 5.5, the derivation of the spectral expansion (6.1.9) depends on the spectral sign lemma, whose validity in the case of $\epsilon$-parameter is equivalent to (6.1.21), which we still have not yet proven at this point.

Recall that we have four quantities:
\begin{eqnarray*}
\widetilde{r}^N_{\psi^N}(w_u), \,\ \epsilon^1_{\psi^N}(u), \,\  \sgn^{N,0}(w_u), \,\ \epsilon^N_{\psi^N}(x_u)
\end{eqnarray*}
occuring in the statement of the spectral sign lemma for $\widetilde{G}_{E/F}(N)$ (note that $x_u \in \widetilde{\mathcal{S}}_{\psi^N}$ is the element that is denoted as $x_1$ in (6.1.8) and (6.1.9)). Since $\widetilde{S}_{\psi^N}$ is an abelian torsor, we have $\epsilon^N_{\psi^N}(x_u) =\epsilon^1_{\psi^N}(u) =\sgn^{N,0}(w_u)=1$, and it follows from the discussion of section 5.8 (in particular remark 5.8.2) that
\[
\widetilde{r}^N_{\psi^N}(w_u)=\widetilde{r}^N_{\psi^N}(w) = \epsilon(1/2,\mu_1 \times \mu_2^c)
\]
($w$ being the unqiue element of $W_{\psi^N}$). Then instead of (6.1.9), the spectral expansion is {\it a priori} of the form:
\begin{eqnarray}
& & \widetilde{I}^N_{\disc,\psi^N}(\widetilde{f}) = \widetilde{I}^N_{\disc,\psi^N}(\widetilde{f}) - \leftexp{0}{\widetilde{r}^N_{\disc,\psi^N} (\widetilde{f})} \\
&=& \frac{1}{2}  \epsilon(1/2,\mu_1 \times \mu_2^c) i_{\psi^N} (x_1)  \widetilde{f}_N(\psi^N,x_1). \nonumber
\end{eqnarray}

\noindent The endoscopic expansion (6.1.8) remains valid:
\begin{eqnarray}
\widetilde{I}^N_{\disc,\psi^N}(\widetilde{f}) - \leftexp{0}{\widetilde{s}^N_{\disc,\psi^N} (\widetilde{f})} = \frac{1}{2}   e_{\psi^N} (x_1)  \widetilde{f}^{\prime}_N(\psi^N,x_1). 
\end{eqnarray}

\noindent Hence using the identities $i_{\psi^N}(x_1)=e_{\psi^N}(x_1)$ and $\widetilde{f}^{\prime}_N(\psi^N,x_1)= \widetilde{f}_N(\psi^N,x_1)$, we obtain from (6.1.22) and (6.1.23):
\begin{eqnarray}
& & \leftexp{0}{\widetilde{s}^N_{\disc,\psi^N}(\widetilde{f})} \\ 
&=& \frac{1}{2} (\epsilon(1/2,\mu_1 \times \mu_2^c)-1) i_{\psi_N}(x_1) \widetilde{f}^{\prime}_N(\psi^N,x_1). \nonumber
\end{eqnarray} 

\noindent Now since $\psi^N$ is an $\epsilon$-parameter we have $\psi^N \in \xi_* \Psi_2(G)$, for a unique $G=(G,\xi) \in \widetilde{\mathcal{E}}_{\ellip}(N) \smallsetminus \widetilde{\mathcal{E}}_{\simp}(N)$ and $\psi \in \Psi_2(G)$ such that $\psi^N = \xi_* \psi$. Namely that we have
\begin{eqnarray*}
G = G_1 \times G_2, \,\ G_i =(G_i,\xi_i) \in \widetilde{\mathcal{E}}_{\simp}(N_i), \,\ N_i < N 
\end{eqnarray*}
\begin{eqnarray*}
& &  \psi^{N_i}_i \in (\xi_i)_* \widetilde{\mathcal{F}}_{\simp}(G_i) \\
& & \psi^{N_i}_i = (\xi_i)_* \psi_i, \,\ \psi_i \in \widetilde{\mathcal{F}}_{\simp}(G_i)
\end{eqnarray*}

\noindent We have $\widetilde{f}_N^{\prime}(\psi^N,x_1) = \widetilde{f}^G(\psi)$. Furthermore, since $\psi_i \in \widetilde{\mathcal{F}}_{\simp}(G_i)$, we have, by the usual argument of considering the spectral and endoscopic expansion with respect to $G_i$ and $\psi_i^{N_i}$, that
\begin{eqnarray}
\tr R^{G_i}_{\disc,\psi_i^{N_i}} (f_i)= S^{G_i}_{\disc,\psi_i^{N_i}}(f_i) = f_i^{G_i}(\psi_i), \,\ f_i \in \mathcal{H}(G_i).
\end{eqnarray} 
(note that $\mathcal{S}_{\psi_i}(G_i)$ is trivial since $\psi_i \in \widetilde{\mathcal{F}}_{\simp}(G_i)$). We claim that the stable linear form $f^G(\psi)$ on $G$ is a sum, with non-negative coefficients, of irreducible characters of representations $\pi = \pi_1 \boxtimes \pi_2$ on $G=G_1 \times G_2$. Indeed it suffices to verify this in the case when $f = f_1 \times f_2$, in which case we have $f^G = f_1^{G_1} \times f_2^{G_2}$, in which case we have by (6.1.7)
\[
f^G(\psi) = f_1^{G_1}(\psi_1) \times f_2^{G_2}(\psi_2)
\]
and we conclude from (6.1.25).

\noindent We also note that the number $i_{\psi^N}(x_1)$ is positive, which follows immediately from the definition and the fact that $\overline{S}_{\psi^N}$ is an abelian torsor, namely that in fact we have $i_{\psi^N}(x_1)=1/4$. Thus finally we can write (6.1.24) in the form:
\begin{eqnarray}
& & \\
& & \sum_{G^* \in \widetilde{\mathcal{E}}_{\simp}(N)} \widetilde{\iota}(N,G^*) \leftexp{0}{\widehat{S}^{G^*}_{\disc,\psi^N}(\widetilde{f}^{G^*}) } + \frac{1}{8}(1-\epsilon(\frac{1}{2}, \mu_1 \times \mu_2^c) )  \widetilde{f}^G(\psi) =0. \nonumber
\end{eqnarray}

\noindent Now in (6.1.26) we can work interchangeably between $\widetilde{f} \in \widetilde{\mathcal{H}}(N)$ and a compatible family of functions $\mathcal{F} =\{f_{G^*} \in \mathcal{H}(G^*)| \,\ G^* \in \widetilde{\mathcal{E}}_{\ellip}(N) \}$. Denote by $f$ the function in the family $\mathcal{F}$ that is associated to $G$. We thus have:
\begin{eqnarray}
& & \\
& & \sum_{G^* \in \widetilde{\mathcal{E}}_{\simp}(N)} \widetilde{\iota}(N,G^*) \leftexp{0}{S^{G^*}_{\disc,\psi^N}(f_{G^*}) } + \frac{1}{8}(1-\epsilon(\frac{1}{2}, \mu_1 \times \mu_2^c) )  f^G(\psi) =0. \nonumber
\end{eqnarray}

\noindent By case (1) of proposition 5.7.1 we have
\begin{eqnarray}
\tr R^{G^*}_{\disc,\psi^N}=S^{G^*}_{\disc,\psi^N}=\leftexp{0}{S^{G^*}_{\disc,\psi^N}}.
\end{eqnarray} 
Hence we have: 
\begin{eqnarray}
& & \\
& & \sum_{G^* \in \widetilde{\mathcal{E}}_{\simp}(N)} \widetilde{\iota}(N,G^*) \tr R^{G^*}_{\disc,\psi^N}(f_{G^*})  + \frac{1}{8}(1-\epsilon(\frac{1}{2}, \mu_1 \times \mu_2^c) )  f^G(\psi) =0. \nonumber
\end{eqnarray} 
Recall that we have establsihed that $f^G(\psi)$ is a sum with positive coefficients of characters of irreducible representations on $G$. Since the coefficient $1-\epsilon(1/2,\mu_1 \times \mu_2)$ is non-negative (namely either equal to zero or two), we see that we can apply lemma 4.3.6 on vanishing of coefficients to (6.1.29). From the vanishing of coefficients we conclude that
\begin{eqnarray}
& &  \tr R^{G^*}_{\disc,\psi^N} \equiv 0 \mbox{ for } G^* \in \widetilde{\mathcal{E}}_{\simp}(N)     \\
& & \epsilon(1/2,\mu_1 \times \mu_2^c)=1 \nonumber
\end{eqnarray}
and this concludes the proof by (6.1.28) again.
\end{proof}

\bigskip

\begin{rem}
\end{rem}
The part of proposition 6.1.5 concerning the $\epsilon$-factor $\epsilon(1/2,\mu_1 \times \mu_2^c)$ forms the main induction step of the argument for the proof of part (b) of theorem 2.5.4 (in the present setting at least for parameters in the family $\widetilde{\mathcal{F}}$; we will complete the argument in section 9. Proposition 6.1.3 and 6.1.5 also completes the proof that parameters $\psi^N \in \xi_* \widetilde{\mathcal{F}}_2(G)$ for $G=(G,\xi) \in \widetilde{\mathcal{E}}_{\ellip}(N) \smallsetminus \widetilde{\mathcal{E}}_{\simp}(N)$ do not contribute to the discrete spectrum of $G^*$ for simple $G^* \in \widetilde{\mathcal{E}}_{\simp}(N)$.

\subsection{The case of elliptic parameters}

Recall that in the last subsection we have treated the case of square-integrable parameters $\psi \in \widetilde{\mathcal{F}}_2(G)$ for $G \in \widetilde{\mathcal{E}}_{\ellip}(N)$. In this section we treat the case of parameters $\psi \in \widetilde{\mathcal{F}}_{\ellip}(G)$ for $G \in \widetilde{\mathcal{E}}_{\simp}(N)$. Note that these two classes concerns exactly those parameters that are not covered by proposition 5.7.4. 

Thus for the rest of the subsection we let $G =(G,\xi) \in \widetilde{\mathcal{E}}_{\simp}(N)$. We denote:
\begin{eqnarray*}
& & \widetilde{\Psi}_{\ellip}^2(G) := \widetilde{\Psi}_{\ellip}(G) \smallsetminus \widetilde{\Psi}_2(G) \\
& & \widetilde{\mathcal{F}}^2_{\ellip}(G) := \widetilde{\mathcal{F}}(G) \cap  \widetilde{\Psi}_{\ellip}^2(G).
\end{eqnarray*}

\noindent In this subsection we will be concerned with the set of parameters $\widetilde{\mathcal{F}}_{\ellip}^2(G)$. In fact for parameters $\psi \in \widetilde{\Psi}_{\ellip}^2(G)$, we have $\xi_* \psi \notin \widetilde{\Psi}_{\ellip}(N)$, hence by remark 6.1.2 the conditions of hypothesis 6.1.1 already holds for such $\psi$ from the induction hypothesis, so in fact we can just work with the set $\widetilde{\Psi}_{\ellip}^2(G)$. But in accordance with the convention to work within the family $\widetilde{\mathcal{F}}$, we will state the results for the parameters $\widetilde{\mathcal{F}}_{\ellip}^2(G)$.

Thus let $\psi \in \widetilde{\mathcal{F}}_{\ellip}^2(G)$ and $\psi^N := \xi_* \psi$. Then $\psi^N$ is characterized by the condition:
\begin{eqnarray}
& & \psi^N = 2 \psi_1^{N_1} \boxplus \cdots \boxplus 2 \psi_q^{N_q} \boxplus \psi_{q+1}^{N_{q+1}} \boxplus \cdots \boxplus \psi_r^{N_r}, \,\ q \geq 1 \\
& & S_{\psi}(G) = O(2,\mathbf{C})^q \times O(1,\mathbf{C})^{r-q} \cong O(2,\mathbf{C})^q \times (\mathbf{Z}/2\mathbf{Z})^{r-q} \nonumber \\
& & |\mathcal{S}_{\psi}(G)|=2^{r- t} \nonumber
\end{eqnarray}
where $t=0$ if $r=q$ and $t=1$ if $r > q$. 

\noindent We denote by $\mathcal{S}_{\psi,\ellip}$ the set of components of $\overline{S}_{\psi}$ indexed by $x \in \mathcal{S}_{\psi}$ such that $\mathcal{E}_{\psi,\ellip}(x)$ (as defined in section 5.6) is non-empty. Then for each $x \in \mathcal{S}_{\psi,\ellip}$, there is a unique element $w_x \in W_{\psi,\reg}(x)$ (with $W_{\psi,\reg}(x)$ being empty if $x \notin \mathcal{S}_{\psi,\ellip}$).

One of the main tasks is to establish that such parameters $\psi \in \widetilde{\mathcal{F}}_{\ellip}^2(G)$ do not contribute to the discrete spectrum of $G$, and to establish the stable multiplicity formula for $\psi$, which reduces to the assertion that $S^G_{\disc,\psi^N} =0$ (note that from (6.2.1) we have that $\overline{S}_{\psi}^0$ contains a non-trivial central torus, so the number $\sigma(\overline{S}_{\psi})=0$ by property (5.1.9)). To establish this in general would require the global intertwining relation, which we have not establish at this point. In this subsection, we collect some results which would be needed in later subsections. 

Recall that in the present case for $x \in \mathcal{S}_{\psi,\ellip}$, the sets $\mathfrak{N}_{\psi,\reg}(x) \stackrel{\sim}{\rightarrow} W_{\psi,\reg}(x)$ are singleton. Thus denoting the unique element of $\mathfrak{N}_{\psi,\reg}(x)$ as $u_x$ (thus $u_x$ maps to $w_x \in W_{\psi,reg}(x)$), we can take the liberty of setting $f_G(\psi,x):=f_G(\psi,u_x)$ ($f \in \mathcal{H}(G)$) for $x \in \mathcal{S}_{\psi,\ellip}$.    

\begin{proposition}
For $\psi \in \widetilde{\mathcal{F}}_{\ellip}^2(G)$ as in (6.2.1), suppose that the index $r$ satisfies $r > 1$. Then there is a positive constant $c$ such that, for any compatible family $\mathcal{F}=\{f^* \in \mathcal{H}(G^*)| \,\ G^* \in \widetilde{\mathcal{E}}_{\simp}(N)  \}$, we have
\begin{eqnarray}
& & \sum_{G^* \in \widetilde{\mathcal{E}}_{\simp}(N)} \widetilde{\iota}(N,G^*) \tr R^G_{\disc,\psi^N}(f^*) \\
&=& c \sum_{x \in \mathcal{S}_{\psi,\ellip}}  \epsilon_{\psi}^G(x) \big( f_G^{\prime}(\psi,s_{\psi}x) - f_G(\psi,x) \big) \nonumber
\end{eqnarray}
(here as before $f$ is the function in the compatible family that is associated to $G$).
\end{proposition}  
\begin{proof}
Since in the case where $\psi \in \widetilde{\mathcal{F}}_{\ellip}^2(G)$ we have $\psi^N \notin \widetilde{\Psi}_{\ellip}(N)$. Hence by corollary 5.7.3 we again have the validity of the following identity:
\begin{eqnarray}
\sum_{G^* \in \widetilde{\mathcal{E}}_{\simp}(N)} \widetilde{\iota}(N,G^*) \leftexp{0}{S^{G^*}_{\disc,\psi^N}(f^*)} =0.
\end{eqnarray}

As before denote by $G^{\vee} =(G^{\vee},\xi^{\vee})$ the other element of $\widetilde{\mathcal{E}}_{\simp}(N)$ (i.e. other than $G$). Then we claim that
\begin{eqnarray}
\tr R^{G^{\vee}}_{\disc,\psi^N}(f^{\vee}) = \leftexp{0}{S^{G^{\vee}}_{\disc,\psi^N}(f^{\vee})  }, \,\ f^{\vee} \in \mathcal{H}(G^{\vee}).
\end{eqnarray} 

Indeed if $\psi^N \notin (\xi^{\vee})_* \Psi(G^{\vee})$ then this follows from case (1) of proposition 5.7.1. On the other hand if $\psi^N \in (\xi^{\vee})_* \Psi(G^{\vee})$, then from the form of $\psi^N$ in (6.2.1) we must have $r=q \geq 2$, and we have (denoting by $\psi^{\vee}$ the parameter in $\Psi(G^{\vee})$ defined by $\psi^N$):
\begin{eqnarray}
& & S_{\psi^{\vee}}(G^{\vee}) = Sp(2,\mathbf{C})^q \\
& & \mathcal{S}_{\psi^{\vee}} = \{1\}. \nonumber
\end{eqnarray}

\noindent It follows, with the notations of proposition 5.3.4, that
\begin{eqnarray*}
\dim \overline{T}_{\psi^{\vee}} = \dim \overline{T}_{\psi^{\vee},x} = q \geq 2
\end{eqnarray*}
(here $x$ is of course the unique trivial element of $\mathcal{S}_{\psi^{\vee}}(G^{\vee})$). Thus by proposition 5.3.4 the global intertwining relation is valid for the pair $(G^{\vee},\psi^{\vee})$, and so by case (2) of proposition 5.7.1, we see that (6.2.4) is again valid.

Thus it remains to treat the term for $G$ itself. We use the spectral expansion (5.5.20) and the endoscopic expansion (5.6.32) for the pair $(G,\psi)$. In both expansion the summand is empty if $x \notin \mathcal{S}_{\psi,\ellip}$, so we can limit the sum in both expansions to $x \in \mathcal{S}_{\psi,\ellip}$. The number $i_{\psi}(x) = e_{\psi}(x)$ is easy to compute: for $x \in \mathcal{S}_{\psi,\ellip}$, we have
\begin{eqnarray*}
i_{\psi}(x) = |W_{\psi}^0|^{-1} |\det(w_x -1 )|^{-1} = 1 \cdot \frac{1}{2^q} = \frac{1}{2^q}.
\end{eqnarray*}

\noindent Hence we have the expansions
\begin{eqnarray}
I^G_{\disc,\psi^N}(f)  -  \tr R^G_{\disc, \psi^N}(f)  = \frac{1}{2^{q+r-t}} \sum_{x \in \mathcal{S}_{\psi,\ellip}} \epsilon_{\psi}^G(x) f_G(\psi,x)
\end{eqnarray}
and
\begin{eqnarray}
& & \\
& & I^G_{\disc,\psi^N}(f)  -  \leftexp{0}{S^G_{\disc, \psi^N}(f)}  = \frac{1}{2^{q+r-t}} \sum_{x \in \mathcal{S}_{\psi,\ellip}} \epsilon_{\psi}^G(x) f^{\prime}_G(\psi,s_{\psi}x). \nonumber
\end{eqnarray}

\noindent Hence we have
\begin{eqnarray}
& &  \leftexp{0}{S^G_{\disc, \psi^N}(f)} \\
&=& \tr R^G_{\disc, \psi^N}(f) + \frac{1}{2^{q+r-t}} \sum_{x \in \mathcal{S}_{\psi,\ellip}} \epsilon_{\psi}^G(x) \big( f_G(\psi,x)  -  f^{\prime}_G(\psi,s_{\psi}x) \big). \nonumber
\end{eqnarray}

\noindent Substituting (6.2.4) and (6.2.8) into (6.2.3), we obtain the assertion (6.2.2), with the constant $c= 1/2^{q+r -t+1}$ (noting that $\widetilde{\iota}(N,G)=1/2$).

\end{proof}

\bigskip

We next conisder the case where $r=1$ in (6.2.1), in other words
\begin{eqnarray}
& & \psi^N = 2 \psi_1^{N_1} \\
& & S_{\psi}(G) = O(2,\mathbf{C}) \nonumber \\
& & |\mathcal{S}_{\psi}(G)|=2 \nonumber
\end{eqnarray} 
In particular $N=2 N_1$ has to be even. As above denote by $G^{\vee} = (G^{\vee},\xi^{\vee}) \in \widetilde{\mathcal{E}}_{\simp}(N)$ the other element not equal to $G$. Then we have $\psi^N \in (\xi^{\vee})_* \widetilde{\mathcal{F}}(G^{\vee})$, with $\psi^{\vee} \in \widetilde{\mathcal{F}}(G^{\vee})$ the parameter of $G^{\vee}$ defined by $\psi^N$. We have
\begin{eqnarray*}
& & S_{\psi^{\vee}}(G^{\vee}) = Sp(2,\mathbf{C}) \\
& & \mathcal{S}_{\psi^{\vee}}(G^{\vee}) = \{1\}.
\end{eqnarray*}

In the next proposition we denote by $x_1$, respectively $x_1^{\vee}$ the unique element of $\mathcal{S}_{\psi,\ellip}(G)$, and respectively $\mathcal{S}_{\psi^{\vee},\ellip}(G^{\vee})$ (of course $x_1$ is simply the non-trivial element of $\mathcal{S}_{\psi}$ and $x_1^{\vee}$ the identity element of $\mathcal{S}_{\psi^{\vee}}$). We also denote by $M^{\vee}$ the Levi component of the standard Siegel parabolic subgroup of $G^{\vee}$. Then we can identify $M^{\vee} = G_{E/F}(N/2)$. We regard $M^{\vee}=(M^{\vee},\xi^{\vee})$ as an endoscopic datum $\widetilde{\mathcal{E}}(N)$, and in particular as a Levi sub-datum of $G^{\vee} \in \widetilde{\mathcal{E}}_{\simp}(N)$. The parameter $\psi^{N}=2\psi_1^{N_1}$ then defines a parameter $\psi_{1,M^{\vee}} \in \widetilde{\Psi}(M^{\vee})$ such that $\psi^{N} = (\xi^{\vee})_* \psi_{1,M^{\vee}}$.   

\begin{proposition}
Let $\psi \in \widetilde{\mathcal{F}}^2_{\ellip}(G)$ such that $\psi^N$ is as in (6.2.9). Then for any compatible family $\mathcal{F}$ as above we have
\begin{eqnarray}
& & \sum_{G^* \in \widetilde{\mathcal{E}}_{\simp}(N)} \tr R^G_{\disc,\psi^N}(f^*)  + \frac{1}{8} \big(   (f^{\vee})^{M^{\vee}}(\psi_{1,M^{\vee}}) - f^{\vee}_{G^{\vee}}(\psi^{\vee},x_1^{\vee})           \big) \\
& & = \frac{1}{8} \big(  f^{\prime}_G(\psi,s_{\psi}x_1) - f_G(\psi,x_1)   \big). \nonumber
\end{eqnarray} 
\end{proposition}
\begin{proof}
The proof is similar to that of proposition 6.2.1; the only complication is that this time we do not know the global intertwining relation for the pair $(G^{\vee},\psi^{\vee})$, which is the reason for the extra term on the left hand side of (6.2.10).

Identity (6.2.3) remains valid in this case, and so is (6.2.8); in fact in the present case (6.2.8) simplies to:
\begin{eqnarray}
& &\\
& & \leftexp{0}{S^G_{\disc, \psi^N}(f)} = \tr R^G_{\disc, \psi^N}(f) +  \frac{1}{4}  \big( f_G(\psi,x_1)  -  f^{\prime}_G(\psi,s_{\psi}x_1) \big) \nonumber
\end{eqnarray}
here we note that the sign $\epsilon^G_{\psi}(x_1)=1$, as follows from the form of $\psi^N$ in (6.2.9), together with lemma 5.8.1 and the discussion after equation (5.8.10).

\noindent By the same argument we also have the following identity for the pair:
\begin{eqnarray}
& & \leftexp{0}{S^{G^{\vee}}_{\disc, \psi^N}(f^{\vee})} \\
&=& \tr R^{G^{\vee}}_{\disc, \psi^N}(f^{\vee}) +   i_{\psi^{\vee}}(x_1^{\vee}) \big( f^{\vee}_{G^{\vee}}(\psi^{\vee},x^{\vee}_1)  -  (f^{\vee})^{\prime}_{G^{\vee}}(\psi^{\vee}, x^{\vee}_1) \big) \nonumber
\end{eqnarray}
(in this case the element $s_{\psi^{\vee}}$ does not matter because $S_{\psi^{\vee}}$ is connected and hence $\mathcal{S}_{\psi^{\vee}}$ is trivial). The constant $ i_{\psi^{\vee}}$ is given by (here $w$ is the unique Weyl element of $Sp(2,\mathbf{C})$):
\begin{eqnarray}
 i_{\psi^{\vee}} = |W_{\psi^{\vee}}^0|^{-1} \frac{\sgn^0(w)}{|\det(w-1)|} = -1/4
\end{eqnarray}

\noindent and by descent we have
\begin{eqnarray*}
(f^{\vee})^{\prime}_{G^{\vee}}(\psi,x_1^{\vee}) = (f^{\vee})^{G^{\vee}}(\psi) = (f^{\vee})^{M^{\vee}}(\psi_{1,M^{\vee}}). 
\end{eqnarray*}

\noindent Thus we have
\begin{eqnarray}
& & \leftexp{0}{S^{G^{\vee}}_{\disc, \psi^N}(f^{\vee})} - \tr R^{G^{\vee}}_{\disc, \psi^N}(f^{\vee}) \\
&=& + \frac{1}{4} \big(   (f^{\vee})^{M^{\vee}}(\psi_{1,M^{\vee}}) - f^{\vee}_{G^{\vee}}(\psi^{\vee},x_1^{\vee})           \big). \nonumber
\end{eqnarray}

\noindent Substituting (6.2.11)-(6.2.14) into (6.2.3) (and again using $\widetilde{\iota}(N,G) = \widetilde{\iota}(N,G^{\vee})=1/2$) we obtain (6.2.10).
\end{proof}

\bigskip
\bigskip

In section 6.3 we will use the technique of this subsection to study the square-integrable parameters of section 6.1. In order to carry out the argument it is necessary to work with parameters in $\widetilde{\mathcal{F}}$ of degree larger than $N$. 

Thus let $\psi^{N_+}_+ \in \widetilde{\mathcal{F}}(N_+)$ (in the application we have in mind $N_+$ is typically larger than $N$), written in the usual form:
\begin{eqnarray}
\psi^{N_+}_+ = l_1 \psi_1^{N_1} \boxplus \cdots \boxplus l_r \psi_r^{N_r} 
\end{eqnarray}
with $\psi_i^{N_i} \in \widetilde{\mathcal{F}}_{\simp}(N_i)$, mutually distinct.

\noindent Put
\begin{eqnarray}
& & N_{+,-} := \sum_{l_i \text{odd}} N_i \\
& &  \psi^{N_{+,-}}_{+,-} = \bigboxplus_{l_i \text{odd}} \psi_i^{N_i}.        \nonumber
\end{eqnarray}

Suppose that $\psi^{N_+}_+ \in (\xi_+)_* \widetilde{\mathcal{F}}(G_+)$ for some $G_+=(G_+,\xi_+) \in \widetilde{\mathcal{E}}_{\simp}(N_+)$, with $\psi_+^{N_+} = (\xi_+)_* \psi_+$. Put:
\begin{eqnarray}
M_+ := G_{E/F}(N_1)^{l_1^{\prime}} \times \cdots \times  G_{E/F}(N_r)^{l_r^{\prime}} \times G_{+,-}
\end{eqnarray}
here $l_i^{\prime}$ is the greatest integer less than or equal to $l_i/2$, and $G_{+,-}=U_{E/F}(N_{+,-})$. Then $M_+$ is a Levi subgroup of $G_+$, but as before we regard $M_+ =(M_+,\xi_+)$ as an endoscopic datum $\widetilde{\mathcal{E}}(N)$, which is a Levi sub-datum of $G_+=(G_+,\xi_+) \in \widetilde{\mathcal{E}}_{\simp}(N_+)$. We also note that the $L$-embedding $\xi_+$ restricts to an $L$-embedding $\leftexp{L}{G_{+,-}} \hookrightarrow \leftexp{L}{G_{E/F}(N_{+,-})}$ which we denote as $\xi_{+,-}$. Thus we have the datum $G_{+,-}=(G_{+,-},\xi_{+,-}) \in \widetilde{\mathcal{E}}_{\simp}(N_{+,-})$.

Then we have the parameter $\psi_{+,M_+} \in \widetilde{\Psi}_2(M_+)$ of the form:
\begin{eqnarray}
\psi_{+,M_+} := \psi_1^{l_1^{\prime}} \times \cdots \times \psi_r^{l_r^{\prime}} \times \psi_{+,-}
\end{eqnarray}
with $\psi_{+,-} \in \widetilde{\mathcal{F}}_2(G_{+,-})$ and $\psi_i \in \widetilde{\mathcal{F}}_{\simp}(N_i)$ such that $\psi_{+,-}^{N_{+,-}} = \xi_{+,-,*} (\psi_{+,-})$ and $\psi_+^{N_+}=\xi_{+,*} (\psi_{+,M_+})$.
Note that $\psi_{+,M_+}$ maps to $\psi_+$ under the dual Levi embedding $\leftexp{L}{M_+} \hookrightarrow \leftexp{L}{G_+}$. 

Suppose that the condition $N_{+,-} \leq N$ is satisfied. Then the linear form $f^{\prime}_{G_{+,-}}(\psi_{+,-}, x_{+,-})$ is defined for any $x_{+,-} \in \mathcal{S}_{\psi_{+,-}}$, by hypothesis 6.1.1. It then follows that the linear form $f^{\prime}_{M_+}(\psi_{+,M_+},x_{M_+})$ is defined on $\mathcal{H}(M_+)$ for any $x_{M_+} \in \mathcal{S}_{\psi_{+,M_+}}$. The same descent argument as in lemma 5.7.2 then shows that the linear form $f^{\prime}_{G^+}(\psi_+,x_+)$ is defined for any $x_+ \in \mathcal{S}_{\psi_+}$. In particular, the stable linear form $f^{G_+}(\psi_+)$ is defined and we have
\[
f^{G_+}(\psi_+) = f^{M_+}(\psi_{+,M_+}).
\] 

\noindent It follows that we can form the distribution 
\[
\leftexp{0}{S^{G^*}_{\disc,\psi_+^{N_+}}}
\] 
on $\mathcal{H}(G^*)$, as the difference between $S^{G^*}_{\disc,\psi_+^{N_+}}$ and its expected value, for any $G^* =(G^*,\xi^*) \in \widetilde{\mathcal{E}}_{\ellip}(N_+)$. Namely that
\begin{eqnarray*}
& & \leftexp{0}{S^{G^*}_{\disc,\psi_+^{N_+}}(f^*)}\\
&:=& S^{G^*}_{\disc,\psi_+^{N_+}}(f^*) -   \sum_{\psi^* \in \widetilde{\mathcal{F}}(G^*,\psi_+^{N_+}) }   \frac{1}{|\mathcal{S}_{\psi^*}|} \epsilon_{\psi^*}(\psi^*) \sigma(\overline{S}^0_{\psi^*}) f^*(\psi^*)            \nonumber
\end{eqnarray*}
where $\widetilde{\mathcal{F}}(G^*,\psi_+^{N_+})$ is defined as follows: if $G^* = G^1 \times G^2$ with $G^i =(G^i,\xi^i) \in \widetilde{\mathcal{E}}_{\simp}(N^i)$ for $i=1,2$ (so that $ \xi^* = \xi^1 \times \xi^2   $), then
\begin{eqnarray*}
& & \widetilde{\mathcal{F}}(G^*,\psi_+^{N_+}) \\
&=& \{\psi^{*} = \psi^1 \times \psi^2 \in \Psi(G^*)| \,\  \psi^i \in \widetilde{\mathcal{F}}(G^i) \text{ for } i=1,2, \,\  \psi_+^{N_+} =   \psi^{1,N^1} \boxplus \psi^{2,N^2}   \}.
\end{eqnarray*}
here above $\psi^{i,N^i}:=(\xi^i)_* \psi^i \in \widetilde{\mathcal{F}}(N^i)$.

\begin{proposition}
Suppose that $\psi_+^{N_+} \in \widetilde{\mathcal{F}}(N_+)$ as above, such that
\begin{eqnarray}
N_1 + \cdots + N_r \leq N
\end{eqnarray}
Then for any compatible family of functions $\mathcal{F}^+ =\{f \in \mathcal{H}(G^*)| \,\ G^* \in \widetilde{\mathcal{E}}_{\ellip}(N_+)    \}$, we have
\begin{eqnarray}
\sum_{G^* \in \widetilde{\mathcal{E}}_{\ellip}(N_+)} \widetilde{\iota}(N_+,G^*) \leftexp{0}{S^{G^*}_{\disc,\psi_+^{N_+}}(f^*)}=0.
\end{eqnarray}
\end{proposition}
\begin{proof}
We first note the difference between (6.2.20) and (6.2.3): in (6.2.20) we have to sum over all $G^* \in \widetilde{\mathcal{E}}_{\ellip}(N_+)$, not just the simple endoscopic data; this is because in general we have $N_+ > N$, and so we cannot deduce {\it a priori} from the induction hypothesis the validity of the stable multiplicity formula for the composite $G^*$. 

To prove (6.2.20), we again use the spectral expansion (5.5.20) and the endoscopic expansion (5.6.32), applied to the pair $(\widetilde{G}_{E/F}(N_+),\psi_+^{N_+})$. Here for the endoscopic expansion (5.6.32), the term $\leftexp{0}{\widetilde{s}^{N_+}_{\disc,\psi_+^{N_+}}}$ needs to be defined instead to be:
\begin{eqnarray}
\leftexp{0}{\widetilde{s}^{N_+}_{\disc,\psi_+^{N_+}}} (\widetilde{f})= \sum_{G^* \in \widetilde{\mathcal{E}}_{\ellip}(N_+)} \widetilde{\iota}(N_+,G^*) \leftexp{0}{S^{G^*}_{\disc,\psi_+^{N_+}}(\widetilde{f}^{G^*})}.
\end{eqnarray} 

\noindent Under the assumption that (6.2.19), we have in particular $N_{+,-} \leq N$, so the linear forms $\widetilde{f}_{N_+}(\psi_+^{N_+},x_{+})$ and $\widetilde{f}^{\prime}_{N_+}(\psi_+^{N_+},x_+)$ ($\widetilde{f} \in \widetilde{\mathcal{H}}(N_+)$) that occur in both expansions are well-defined. Also recall the fact that for the twisted group $\widetilde{G}_{E/F}(N_+)$, we already know that the linear form $\widetilde{f}_{N_+}(\psi_+^{N_+},x_+)$ depends only on $x_+ \in \mathcal{S}_{\psi_+}$, by the discussion in the paragraph before lemma 5.3.3.

With the modification of $\leftexp{0}{\widetilde{s}^{N_+}_{\disc,\psi_+^{N_+}}}$ as above, then the same argument as in section 5.6 shows that the endoscopic expansion (5.6.32) remains valid in this case. As for the justification of the spectral expansion (5.5.20), the argument in section 5.4 and 5.5 leading up to (5.5.20) shows that, the spectral expansion is valid, provided that the spectral sign lemma (lemma 5.5.1) is valid for $\psi_+^{N_+}$. From the discussions in section 5.8, in particular the discussion in remark 5.8.2, we see that this is valid provided that for all $\epsilon$-sub-parameters of $\psi_+^{N_+}$ the condition (6.1.21) on $\epsilon$-factor is satisfied. But under condition (6.2.19), any $\epsilon$-sub-parameter of $\psi_+^{N_+}$ have degree at most $N$, and so (6.1.21) is indeed satisfied by proposition 6.1.5 (in the case where the sub-parameter has degree less than $N$ then the condition on $\epsilon$-factor follows of course from the induction hypothesis). 

Thus both the spectral and endoscopic expansions hold. Furthermore, the global intertwining relation for $(\widetilde{G}_{E/F}(N_+),\psi_+^{N_+})$ is also valid, by the same descent argument in proposition 5.3.4 (here we are of course using the fact that the global intertwining relation is valid for the pair $(\widetilde{G}_{E/F}(N_{+,-}),\psi_{+,-}^{N_{+,-}})$, by hypothesis 6.1.1 and the condition (6.2.19)). We thus conclude that
\[
 \leftexp{0}{\widetilde{r}^{N_+}_{\disc,\psi_+^{N_+}}}(\widetilde{f})= \leftexp{0}{\widetilde{s}^{N_+}_{\disc,\psi_+^{N_+}}}(\widetilde{f}).
\]
Since we already know that the left hand side vanishes this gives
\[
\leftexp{0}{\widetilde{s}^{N_+}_{\disc,\psi_+^{N_+}}}(\widetilde{f})=0
\]
and the usual argument of replacing $\widetilde{f} \in \widetilde{\mathcal{H}}(N_+)$ by a compatible family $\mathcal{F}^+$ concludes the proof.
\end{proof}

\subsection{Supplementary parameter}

We return to the study of square-integrable parameters in section 6.1, using the technique of ``enlarging" the original parameter. 

Thus we assume $G =(G,\xi) \in \widetilde{\mathcal{E}}_{\simp}(N)$ and $\psi \in \widetilde{\mathcal{F}}_2(G)$. One of the main issue is to resolve the provisional definition of the set of simple generic parameters $\widetilde{\mathcal{F}}_{\gsimp}(G)$ given in section 6.1, which is based on the stable multiplicity formula, with the definition in terms of the original construction in section 2.4, based on the seed theorem 2.4.2 (which of course has yet to be established for simple generic parameters of degree $N$). Similarly, part (a) of theorem 2.5.4 asserts another characterization in terms of the pole at $s=1$ of Asai $L$-functions. We will eventually establish the equivalence of all these three characterizations in section 9. 

Even though the complete proof of these three characterizations will be achieved only in section 9, we can establish these assertions in the next subsection for a particular set of parameters $\widetilde{\mathcal{F}}$ with certain local constraints. The argument is based on the results of this subsection.

As before write $\psi^N = \xi_* \psi \in \widetilde{\mathcal{F}}_{\ellip}(N)$ in the standard form:
\begin{eqnarray}
\psi^N = \psi_1^{N_1} \boxplus \cdots \boxplus \psi_r^{N_r}, \,\ \psi_i^{N_i} \in \widetilde{\mathcal{F}}_{\simp}(N_i).
\end{eqnarray}

We assume the $N_i$'s are arranged so that $N_1 \leq N_i$ for all $i$. Put $N_+ :=N_1 + N$, and form the supplementary parameter
\begin{eqnarray}
\psi_+^{N_+} := \psi_1^{N_1} \boxplus \psi^N = 2 \psi_1^{N_1} \boxplus \psi_2^{N_2} \boxplus \cdots \boxplus \psi_r^{N_r}.
\end{eqnarray}

\noindent Then $\psi_+^{N_+} \in \widetilde{\mathcal{F}}(N_+)$. Furthermore in the notation of (6.2.16) we have $N_{+,-} = N_2 + \cdots + N_r < N$. We denote by $G_+ = (G_+,\xi_+) \in \widetilde{\mathcal{E}}_{\simp}(N_+)$ the unique element of $\widetilde{\mathcal{E}}_{\simp}(N_+)$ that has the same parity as $G=(G,\xi) \in \widetilde{\mathcal{E}}_{\simp}(N)$. In other words, in the notation of section 2.4, if $\xi = \xi_{\chi}$ and $\xi_+ = \xi_{+,\chi_+}$, then with $\kappa_{\chi}$ (resp. $\kappa_{\chi_+}$) the sign such that $\chi \in \mathcal{Z}_E^{\kappa_{\chi}}$ (resp. $\chi_+ \in \mathcal{Z}_E^{\kappa_{\chi_+}}$), we have:
\[
\kappa_{\chi} (-1)^{N-1} = \kappa_{\chi_+} (-1)^{N_+ -1}.
\] 

\noindent Then we have $\psi_+^{N_+} \in (\xi_+)_* \widetilde{\mathcal{F}}_{\ellip}(G_+)$, and we denote $\psi_+ \in \widetilde{\mathcal{F}}_{\ellip}(G_+)$ the parameter of $G_+ = (G_+,\xi_+)$ defined by $\psi_+^{N_+}$. We also denote by $G_+^{\vee}$ the element of $\widetilde{\mathcal{E}}_{\simp}(N_+)$ other than $G_+$.

In the case where $N$ is even we also introduce the Levi subgroup $L_+ \cong G_1 \times G_{E/F}(N/2)$ of $G_+$. Here $G_1 =(G_1,\xi_1) \in \widetilde{\mathcal{E}}_{\simp}(N_1)$ is the datum such that $\psi_1^{N_1} = \xi_{1,*} \psi_1$. It is equipped with the structure of Levi sub-datum of $G_+$ as usual. It is adapted to the decomposition
\[
\psi_+^{N_+} = \psi_1^{N_1} \boxplus \psi^N
\]
and comes with the linear form:
\[
f^{L_+} \mapsto f^{L_+}(\psi_1 \times \Lambda), \,\ f \in \mathcal{H}(G_+).
\] 

As in the previous subsection, we consider separately the cases where $r > 1$ and $r=1$. We first treat the case that $r>1$. However, since we are now working with parameters of degree greater than $N$, we need to make a stronger induction hypothesis, namely that we need to assume that the stable multiplicity formula is valid for parameters in $\widetilde{\mathcal{F}}(N) \smallsetminus \widetilde{\mathcal{F}}_{\ellip}(N)$.

\begin{proposition}
Suppose that the stable multiplicity formula is valid for all the parameters in $\widetilde{\mathcal{F}}(N) \smallsetminus \widetilde{\mathcal{F}}_{\ellip}(N)$. If $r > 1$ in (6.3.1), then there are positive constants $b_+$ and $c$, such that for any compatible family of functions $\mathcal{F}^+ = \{ f^* \in \mathcal{H}(G^*)| \,\ G^* \in \widetilde{\mathcal{E}}_{\simp}(N_+)   \}$, we have the identity:
\begin{eqnarray}
& & \sum_{G^* \in \widetilde{\mathcal{E}}_{\simp}(N_+)} \widetilde{\iota}(N_+,G^*)  \tr R^{G^*}_{\disc,\psi_+^{N_+}}(f^*)  + b_+ f^{L_+}(\psi_1 \times \Lambda) \\
&=&  c \sum_{x_+ \in \mathcal{S}_{\psi_+,\ellip}} \epsilon_{\psi_+}^{G_+}(x_+) \big( f^{\prime}_{G_+}(\psi_+,s_{\psi_+} x_+) - f_{G_+}(\psi_+,x_+)   \big). \nonumber
\end{eqnarray}
(Here $f$ is the function in the compatible family associated to $G_+$.)
\end{proposition}
\begin{proof}
The statement (6.3.3) is formally the same as (6.2.2), and the strategy of the proof is similar. However, since we are now working with parameters of degree larger than $N$, we need to be more careful with our induction hypothesis. 

Since $N_1 + \cdots + N_r =N \leq N$, the hypothesis for proposition 6.2.3 is satisfied, and hence (6.2.20) is valid:
\begin{eqnarray}
\sum_{G^* \in \widetilde{\mathcal{E}}_{\ellip}(N_+)} \widetilde{\iota}(N_+,G_+) \leftexp{0}{S^{G^*}_{\disc,\psi_+^{N_+}}(f^*)}=0.
\end{eqnarray} 

We claim that for $G^* \in \widetilde{\mathcal{E}}_{\ellip}(N_+) \smallsetminus \widetilde{\mathcal{E}}_{\simp}(N_+)$, the stable multiplicity formula is valid for $\psi_+^{N_+}$ with respect to $G^*$, in other words, we have:
\begin{eqnarray}
\leftexp{0}{S^{G^*}_{\disc,\psi_+^{N_+}}(f^*)}=0
\end{eqnarray}
except (possibly) for the datum $G^* = G_1 \times G^{\vee}$. Here we recall
\begin{eqnarray}
& & \leftexp{0}{S^{G^*}_{\disc,\psi_+^{N_+}}(f^*)}\\
&:=& S^{G^*}_{\disc,\psi_+^{N_+}}(f^*) -   \sum_{\psi^* \in \widetilde{\mathcal{F}}(G^*,\psi_+^{N_+}) }   \frac{1}{|\mathcal{S}_{\psi^*}|} \epsilon_{\psi^*}(\psi^*) \sigma(\overline{S}^0_{\psi^*}) f^*(\psi^*).            \nonumber
\end{eqnarray}

\noindent In other words, the terms that correspond to $G^* \in \widetilde{\mathcal{E}}_{\ellip}(N_+) \smallsetminus \widetilde{\mathcal{E}}_{\simp}(N_+)$, except (possibly) for the datum $G^*=G_1 \times G^{\vee}$, do not contribute to (6.3.4). To establish (6.3.5) for composite $G^* \neq G_1 \times G^{\vee}$, we note that any $G^*  \in \widetilde{\mathcal{E}}_{\ellip}(N_+) \smallsetminus \widetilde{\mathcal{E}}_{\simp}(N_+)$ has a decomposition:
\begin{eqnarray}
& & G^* = G^1 \times G^2, \,\ G^k =(G^k,\xi^k) \in \widetilde{\mathcal{E}}_{\simp}(N^k) \\
& & N_+ = N^1 + N^2, \,\   0 < N^1 \leq N^2 <N \nonumber
\end{eqnarray}
and for any $\psi^* \in \widetilde{\mathcal{F}}(G^*,\psi_+^{N_+})$ we have a corresponding decomposition
\begin{eqnarray}
\psi^* = \psi^1 \times \psi^2, \,\ \psi^k \in \widetilde{\mathcal{F}}_{\ellip}(G^k)
\end{eqnarray}
such that, if we denote by $\psi^{k,N^k}$ the parameter $(\xi^k)_* \psi^k \in \widetilde{\mathcal{F}}(N^k)$, then we have
\begin{eqnarray}
\psi_+^{N_+} = \psi^{1,N^1} \boxplus \psi^{2,N^2}
\end{eqnarray}
and we can write (6.3.6) as
\begin{eqnarray}
& & \\
& & \leftexp{0}{S^{G^*}_{\disc,\psi_+^{N_+}}(f^*)} \nonumber \\
&=& \sum_{ \substack{ \psi^{1,N^1} \times \psi^{2,N^2} \\  \psi_+^{N_+} = \psi^{1,N^1} \boxplus \psi^{2,N^2}  } }  \Big( S^{G^*}_{\disc, \psi^{1,N^1} \times \psi^{2,N^2} }(f^*)  - \frac{m_{\psi^*}}{|\mathcal{S}_{\psi^*}|} \epsilon_{\psi^*}(\psi^*) \sigma(\overline{S}^0_{\psi^*}) f^*(\psi^*)  \Big)   \nonumber \\
&=& \sum_{\substack{ \psi^{1,N^1} \times \psi^{2,N^2} \\ \psi_+^{N_+} = \psi^{1,N^1} \boxplus \psi^{2,N^2}    }}   \leftexp{0}{ S^{G^*}_{\disc, \psi^{1,N^1} \times \psi^{2,N^2} }(f^*) } \nonumber
\end{eqnarray}
here $m_{\psi^*}=1$ if $\psi^{k,N^k} \in (\xi^k)_* \widetilde{\mathcal{F}}(G^k)$ for both $k=1,2$, and is equal to $0$ otherwise.

We first claim that (6.3.5) is valid if the partition $(N^1,N^2)$ is not equal to $(N_1,N)$. Indeed, suppose that $N^1 < N_1$. Then since $N_1 \leq N_i$ for all $i$, we see that the index of summation in (6.3.10) is empty, and thus (6.3.5) is valid. On the other hand, if $N^1 > N_1$, then we have $0< N^1 \leq N^2 < N$, so from the induction hypothesis the stable multiplicity formula is valid for each $\psi^{k,N^k}$ with respect to $G^k$ for $k=1,2$. Thus the stable multiplicity formula is valid for the parameter $\psi^{1,N^1} \times \psi^{2,N^2}$ with respect to $G^*$. Hence each of the summand in (6.3.10) vanishes, and again we have (6.3.5).

Thus we must analyze the case where $(N^1,N^2)=(N_1,N)$. The most significant case is 
\begin{eqnarray*}
& & G^* = G_1 \times G^{\vee} \\
& & (\psi^{1,N^1},\psi^{2,N^2}) = (\psi_1^{N_1},\psi^N) \nonumber
\end{eqnarray*}
here $G_1 =(G_1,\xi_1) \in \widetilde{\mathcal{E}}_{\simp}(N_1)$ is the unique element such that $\psi_1^{N_1} \in (\xi_1)_* \Psi(G_1)$, and $G^{\vee}$ is the element of $\widetilde{\mathcal{E}}_{\simp}(N)$ other than $G$.

Indeed, conisder any $G^*=G^1 \times G^2$ with the associated partition $(N^1,N^2)=(N_1,N)$, and consider a summand in the sum (6.3.10). Suppose that $\psi^{1,N^1} \neq \psi_1^{N_1}$. Then $\psi^{2,N^2}$ has to contain $2 \psi_1^{N_1}$ has a sub-parameter. It follows that $\psi^{2,N^2} \notin \widetilde{\mathcal{F}}_{\ellip}(N)$ (remember that $N^2=N$). Then by the hypothesis we made about non-elliptic parameters in $\widetilde{\mathcal{F}}(N)$, the stable multiplicity formula holds for $\psi^{2,N^2}$ with respect to $G^2$. As for $\psi^{1,N^1}$, the stable multiplicity formula holds for $\psi^{1,N^1}$ with respect to $G^1$ by the induction hypothesis, since $N^1=N_1 < N$. Thus the corresponding summand in (6.3.10) vanishes. Thus we may assume $\psi^{1,N^1}= \psi_1^{N_1}$ and $\psi^{2,N^2} = \psi^N$ in the sum (6.3.10) (in particular at most only one summand in the sum (6.3.10)). 

If $G^1 \neq G_1$ (as an element of $\widetilde{\mathcal{E}}_{\simp}(N_1)$), then by the induction hypothesis the stable multiplicity formula is valid for $\psi_1^{N_1} \in \widetilde{\mathcal{F}}_{\simp}(N_1)$ with respect to $G^1$ (again since $N_1 < N$), which in this case is the assertion of the vanishing:
\begin{eqnarray*}
S^{G^1}_{\disc,\psi_1^{N_1}}  =0
\end{eqnarray*}
and since we also have $m_{\psi^*}=0$ in the present case then it follows that both terms (of the single remaining summand) in the second line of (6.3.10) vanishes, and so does (6.3.10) itself. Hence we may assume $G^1=G_1$. Since $G^* \in \widetilde{\mathcal{E}}_{\ellip}(N_+)$ is composite, its two factors have to have opposite parity. It follows that we must have $G^2=G^{\vee}$ if $G^1=G_1$. Thus (6.3.5) is valid for any $G^* \in \widetilde{\mathcal{E}}_{\ellip}(N_+) \smallsetminus \widetilde{\mathcal{E}}_{\simp}(N_+)$ with $G^* \neq G_1 \times G^{\vee}$.

Thus we are reduced finally to the case $G^* = G_1 \times G^{\vee}$. But then by proposition 6.1.4, together with the validity (by induction hypothesis) of the stable multiplicity for the pair $(G_1,\psi_1)$, we have for $f^* = f_1 \times f_2^{\vee} \in \mathcal{H}(G_1 \times G^{\vee})$
\begin{eqnarray*}
& & S^{G_1 \times G^{\vee}}_{\disc,\psi_1^{N_1} \times \psi^N }(f^*) \\
&=& S^{G_1}_{\disc,\psi_1^{N_1}}(f_1) \cdot  S^{G^{\vee}}_{\disc,\psi^{N}}(f_2^{\vee}) \\
&=& f_1^{G_1}(\psi_1) \cdot f_2^{\vee,L^{\vee}}(\Lambda^{\vee}) 
\end{eqnarray*}
(note that the sign $\epsilon^{G_1}(\psi_1)=+1$ since $\psi_1$ is a simple parameter). Again $m_{\psi^*}=0$ and so 
\begin{eqnarray*}
& & \leftexp{0}{S^{G_1 \times G^{\vee}}_{\disc,\psi_1^{N_1} \times \psi^N  }(f^*)}=S^{G_1 \times G^{\vee}}_{\disc,\psi_1^{N_1} \times \psi^N  }(f^*) \\
&=& f_1^{G_1}(\psi_1) \cdot f_2^{\vee,L^{\vee}}(\Lambda^{\vee} )=  f^{L_+}(\psi_1 \times \Lambda).
\end{eqnarray*}
For the last equality we are using the fact that when $f^{L_+}= g_1 \times g_2$ with $g_1 \in \mathcal{H}(G_1)$ and $g_2 \in \mathcal{H}(G)$, then we have 
\[
f_1^{G_1} = g_1^{G_1}, \,\ f_2^{\vee,L^{\vee}} = g_2^{L}.
\]

We thus conclude that the only contribution to (6.3.4) coming from $G^* \in \widetilde{\mathcal{E}}_{\ellip}(N_+) \smallsetminus \widetilde{\mathcal{E}}_{\simp}(N_+)$ is:
\begin{eqnarray}
& & \sum_{G^* \in \widetilde{\mathcal{E}}_{\ellip}(N_+) \smallsetminus \widetilde{\mathcal{E}}_{\simp}(N_+)   }   \widetilde{\iota}(N_+,G^*) \leftexp{0}{S^{G^*}_{\disc,\psi_+^{N_+}}(f^*)  } \\
& & = \widetilde{\iota}(N_+,G_1 \times G^{\vee}) f^{L_+}(\psi_1 \times \Lambda). \nonumber
\end{eqnarray}

Hence it remains to analyze the summand in (6.3.4) associated to the two elements of $\widetilde{\mathcal{E}}_{\simp}(N_+)$, namely $G_+$ and $G_+^{\vee}$. We first consider $G_+^{\vee}$. We claim that
\begin{eqnarray}
\leftexp{0}{S^{G^{\vee}}_{\disc,\psi_+^{N_+}}} = \tr R^{G^{\vee}}_{\disc, \psi_+^{N_+}}.
\end{eqnarray}

\noindent We have $\psi_+^{N_+} \notin (\xi^{\vee})_* \Psi(G^{\vee})$, so
\[
\leftexp{0}{S^{G^{\vee}}_{\disc,\psi_+^{N_+}}} =S^{G^{\vee}}_{\disc,\psi_+^{N_+}}.
\] 
However, we cannot apply part (1) of proposition 5.7.1 to conclude (6.3.12) directly, since $N_+ > N$. However, it can still be established in the present case. Indeed, consider the expansion (5.4.2), applied to the pair $(G_+^{\vee},\psi_+^{N_+})$. The right hand side of (5.4.2) are terms that correspond to property Levi subgroups $M^{\vee}_+$ of $G_+^{\vee}$. However, we see that there is no proper $M_+^{\vee}$ such that $\psi_+^{N_+}$ contribute to the discrete spectrum $R^{M^{\vee}_+}_{\disc}$ of $M_+^{\vee}$. Indeed, any such $M_+^{\vee}$ is either attached to a partition of $N_+$ that is not compatible with $\psi_+^{N_+}$, in which case we apply corollary 4.3.8, or $M_+^{\vee}$ is a product of groups for which we can combine our induction hypothesis and the fact that $\widetilde{\Psi}_2(M_+^{\vee},\psi_+^{N_+})$ is empty. Thus we have:
\begin{eqnarray}
I^{G^{\vee}_+}_{\disc,\psi_+^{N_+}}(f^{\vee}) - \tr R^{G^{\vee}_+}_{\disc,\psi_+^{N_+}}(f^{\vee})=0, \,\ f^{\vee} \in \mathcal{H}(G_+^{\vee}).
\end{eqnarray}

\noindent Similarly consider the expansion (5.6.2) applied to $(G_+^{\vee},\psi_+^{N_+})$. The right hand side of (5.6.2) are terms that correspond to proper elliptic endoscopic data $(G_+^{\vee})^{\prime}$ of $G_+^{\vee}$. We similarly see that there is no proper $(G_+^{\vee})^{\prime}$ such that $\psi_+^{N_+}$ contribute to the stable distribution $S^{(G_+^{\vee})^{\prime}}_{\disc}$ of $(G_+^{\vee})^{\prime}$. Indeed, either the two factors of $(G_+^{\vee})^{\prime}$ are attached to a partition of $N_+$ that is incompatible with $\psi_+^{N_+}$, in which case we apply proposition 4.3.4, or the two factors of $(G_+^{\vee})^{\prime}$ are data for which we can apply the induction hypothesis, together with the fact that $\Psi( (G_+^{\vee})^{\prime} ,\psi_+^{N_+})$ is empty. Thus we have
\begin{eqnarray}
I^{G^{\vee}_+}_{\disc,\psi_+^{N_+}}(f^{\vee}) - S^{G^{\vee}_+}_{\disc,\psi_+^{N_+}}(f^{\vee})=0, \,\ f^{\vee} \in \mathcal{H}(G_+^{\vee})
\end{eqnarray}
thus we obtain (6.3.12) from (6.3.13) and (6.3.14).

Thus it remains to analyze the summand in (6.3.4) for $G_+$. We claim that we have:
\begin{eqnarray}
& & \\
& & \leftexp{0}{S^{G^{\vee}}_{\disc,\psi_+^{N_+}}}(f) - \tr R^{G^{\vee}}_{\disc, \psi_+^{N_+}}(f)  \nonumber \\
&=&  \frac{1}{2} \sum_{x_+ \in \mathcal{S}_{\psi_+,\ellip}}  \epsilon_{\psi_+}^{G_+}(x_+) \big(  f_{G_+}(\psi_+,x_+) - f^{\prime}_{G_+}(\psi_+,s_{\psi_+} x_+)  \big) , \,\ f \in \mathcal{H}(G_+). \nonumber \\
& & \,\ \,\ \,\ + f^{L_+}(\psi_1 \times \Lambda). \nonumber
\end{eqnarray}

\noindent In order to establish (6.3.15) we again use the spectral expansion (5.5.20) and the endoscopic expansion (5.6.32), applied to the pair $(G_+,\psi_+^{N_+})$. Since $N_+ > N$, we have to justify the validity of these two expansions. As for the spectral expansion (5.5.20), first note that, exactly the same as in (6.2.1), for $x_+ \in \mathcal{S}_{\psi_+,\ellip}$, the set $\mathfrak{N}_{\psi_+,\reg}(x_+)$ and $W_{\psi_+,\reg}(x_+)$ are singleton. And from the shape of $\psi_+^{N_+}$ in (6.3.2), we see that the (unique up to $G_+$-conjugacy) Levi subgroup $M_+$ of $G_+$ such that $\Psi_2(M_+,\psi_+^{N_+})$ is non-empty is the one given in (6.2.17), which in the present setting is given by:
\begin{eqnarray*}
& & M_+ = G_{E/F}(N_1) \times G_{+,-} \\
& & G_{+,-} = U_{E/F}(N_{+,-}) \\
& & N_{+,-} = N_2 + \cdots + N_r <N
\end{eqnarray*}
and as before we regard $M_+=(M_+,\xi_+) \in \widetilde{\mathcal{E}}(N_+)$ as an endoscopic datum, which is a Levi sub-datum of $G_+ \in \widetilde{\mathcal{E}}_{\simp}(N_+)$. By considering the derivation that leads to (5.5.20) (together with the discussion on the proof of the spectral sign lemma in remark 5.8.2), we see that it requires only the condition $N_{+,-} < N$ for its derivation. Hence (5.5.20) is again valid in the present setting.  

As for the endoscopic expansion (5.6.32), we see from the derivation leading up to (5.6.32) that, in order to have its validity, we need to have the stable multiplicity formula for $\psi_+^{N_+}$ with respect to any proper elliptic endoscopic data $G_+^{\prime}  \in \mathcal{E}_{\ellip}(G_+)$, with $(G_+)^{\prime} = G_1^{\prime} \times G_2^{\prime} $ and $G_k^{\prime}  \in \widetilde{\mathcal{E}}_{\simp}(N_k^{\prime})$ ($k=1,2$). In other words we need the stable multiplicity formula for
\begin{eqnarray*}
S^{G^{\prime}_+}_{\disc, \psi^{1,N_1^{\prime}} \times \psi^{2,N_2^{\prime}}  }(f^{\prime}) =  S^{G_1^{\prime}}_{\disc, \psi^{1,N_1^{\prime}} }(f_1^{\prime}) \cdot S^{G_2^{\prime}}_{\disc, \psi^{2,N_2^{\prime}} }(f_2^{\prime})
\end{eqnarray*}
attached to parameters $\psi^{k,N_k^{\prime}} \in \widetilde{\mathcal{F}}(N_k^{\prime})$ such that 
\begin{eqnarray*}
\psi_+^{N_+} = \psi^{1,N_1^{\prime}} \boxplus \psi^{2,N_2^{\prime}}
\end{eqnarray*}
and functions $f^{\prime} = f_1^{\prime} \times f_2^{\prime} \in \mathcal{H}(G_+^{\prime})$. However, we see that the only case that cannot be treated by the induction hypothesis is the case where
\begin{eqnarray*}
& & (\psi^{1,N_1^{\prime}} ,\psi^{2,N_2^{\prime}}) = (\psi_1^{N_1},\psi^N) \\
& & (G_+)^{\prime} = G_1 \times G.
\end{eqnarray*} 
Now in this case, we have to apply proposition 6.1.4, and the stable multiplicity formula for the pair $(G_1,\psi_1)$ (being valid again by induction hypothesis) to obtain:
\begin{eqnarray}
& & \leftexp{0}{S^{G_1 \times G}_{\disc,\psi_1^{N_1} \times \psi^N}(f^{\prime})} \\
&=& S^{G_1}_{\disc,\psi_1^{N_1}}(f_1^{\prime}) \cdot \leftexp{0}{ S^G_{\disc,\psi^N}(f_2^{\prime}) } \nonumber\\
&=&  (f_1^{\prime})^{G_1}(\psi_1) \cdot (f_2^{\prime})^L(\Lambda) \nonumber\\
&=& f^{L_+}(\psi_1 \times \Lambda). \nonumber
\end{eqnarray}
Thus (6.3.16) gives the extra correction term that must be added to the left hand side of (5.6.32) to have the valid endoscopic expansion. 

As in the proof of proposition 6.2.1, we can limit the summation in both expansions over $x_+ \in \mathcal{S}_{\psi_+,\ellip}$, and the number $i_{\psi_+}(x_+)=e_{\psi_+}(x_+)$ is given by $1/2$. We obtain (6.3.15) by combining the spectral and the endoscopic expansions.

Thus finally substituting (6.3.5), (6.3.11), (6.3.12) and (6.3.15) to (6.3.4) we obtain (6.3.3), with $c= \frac{1}{2} \cdot \widetilde{\iota}(N_+,G_+)$, and $b_+ = \widetilde{\iota}(N,G_1 \times G^{\vee})+\widetilde{\iota}(N,G_+)$.
\end{proof}

\bigskip

\begin{rem}
\end{rem}
In the case where $N$ is odd the linear form $\Lambda$ vanishes ({\it c.f.} proposition 6.1.3), and thus the term $b_+ f^{L_+}(\psi_1 \times \Lambda)$ in (6.3.3) can be omitted when $N$ is odd.
\bigskip

We now consider the case when $r=1$ in (6.3.1) and (6.3.2). We then have:
\begin{eqnarray}
& & \psi_+^{N_+} = 2 \psi^N \\
& & N_+ = 2 N \nonumber \\
& & N_{+,-} =0. \nonumber
\end{eqnarray}

\noindent Recall the definition of $G_+=(G_+,\xi_+)$ and $G_+^{\vee}=(G_+^{\vee},\xi_+^{\vee})$ as elements of $\widetilde{\mathcal{E}}_{\simp}(N_+)$ in the paragraph after equation (6.3.2). Similar to the situation of proposition 6.2.2 we denote by $M_+^{\vee}$ the element of $\widetilde{\mathcal{E}}(N_+)$ that is a Levi sub-datum of $G_+^{\vee}$, and whose underlying endoscopic group is the Levi component of the Siegel parabolic subgroup of $G_+^{\vee}$ (thus $M_+^{\vee} \cong G_{E/F}(N)$).

We have $\psi_+^{N_+} \in (\xi_+)_* \widetilde{\mathcal{F}}_{\ellip}(G_+)$ and $\psi_+^{N_+} \in (\xi^{\vee}_+)_* \widetilde{\mathcal{F}}_{\ellip}(G_+^{\vee})$. Denote by $\psi_+ \in \widetilde{\mathcal{F}}_{\ellip}(G_+)$ and $\psi_+^{\vee} \in \widetilde{\mathcal{F}}_{\ellip}(G_+^{\vee})$ the corresponding parameters. We have
\begin{eqnarray}
& & S_{\psi_+}(G_+) = O(2,\mathbf{C}) \\ 
& & S_{\psi_+^{\vee}}(G_+^{\vee}) = Sp(2,\mathbf{C}). \nonumber
\end{eqnarray}
We also denote by $\psi_{M_+^{\vee}} \in \widetilde{\Psi}(M_+^{\vee})$ the parameter of $M_+^{\vee}=(M_+^{\vee},\xi_+^{\vee})$ defined by $\psi^N$ (i.e. we have $\psi^N=\xi_{+,*}^{\vee} \psi_{M_+^{\vee}}$).

\noindent We again denote by $x_{+,1}$ the unique non-trivial element of $\mathcal{S}_{\psi_+}$, and by $x_{+,1}^{\vee} \in \mathcal{S}_{\psi_+^{\vee}}$ the unque element of $\mathcal{S}_{\psi_+^{\vee}}$ (which is just the trivial element). Then in both cases $x_{+,1}$ (resp. $x_{+,1}^{\vee}$) is the unique element of $\mathcal{S}_{\psi_+,\ellip}$ (resp. $\mathcal{S}_{\psi_+^{\vee},\ellip }$). We denote by $w_+$ and $u_+$ the unique element of $W_{\psi_+,\reg}=W_{\psi_+,\reg}(x_+)$ and $\mathfrak{N}_{\psi_+,\reg} =\mathfrak{N}_{\psi_+,\reg}(x_+)$ respectively, and similarly by $w_+^{\vee}$ and $u_+^{\vee}$ the unique element of $W_{\psi^{\vee}_+,\reg}=W_{\psi^{\vee}_+,\reg}(x^{\vee}_+)$ and $\mathfrak{N}_{\psi^{\vee}_+,\reg} =\mathfrak{N}_{\psi^{\vee}_+,\reg}(x^{\vee}_+)$ respectively.

Compared to the proof of proposition 6.2.2, the proof of proposition 6.3.3 below is complicated by the fact that the spectral sign lemma for the parameter $\psi_+^{N_+} \in \widetilde{\mathcal{F}}(N_+)$ with respect to $G_+$ or $G_+^{\vee}$, is unknown at this point, in the case where $\psi^N \in \widetilde{\mathcal{F}}_{\simp}(N)$ is a simple generic parameter, which in turn is closely related to part (a) of theorem 2.5.4 for the simple generic $\psi^N$. More precisely, define the sign $\delta_{\psi} = \pm 1$ by the following rule. We put $\delta_{\psi}=1$ if $\psi^N$ is not simple generic; on the other hand, if $\psi^N$ is simple generic, define $\delta_{\psi}$ as follows: in the notation of section 2.4 suppose that $\xi=\xi_{\chi}$, and $\kappa=\kappa_{\chi}$. Then we put $\delta_{\psi}=1$ if the Asai $L$-function
\begin{eqnarray*}
L(s,\psi^N, \Asai^{(-1)^{N-1} \kappa})
\end{eqnarray*} 
has a pole at $s=1$; otherwise we put $\delta_{\psi} =-1$ if the Asai $L$-function
\begin{eqnarray*}
L(s,\psi^N, \Asai^{(-1)^{N} \kappa})
\end{eqnarray*} 
has a pole at $s=1$. In other words $\delta_{\psi}=1$ if and only if part (a) of theorem 2.5.4 is valid for $\psi^N$ (in the case where $\psi^N$ is not simple generic its simple generic component has degree less than $N$ and so the result follows from induction hypothesis). 

From the discussion in section 5.8, we see that for the pair $(G_+,\psi_+)$ we have the following identity:
\begin{eqnarray*}
r_{\psi_+}(w_+) \epsilon^1_{\psi_+}(u_+) = \delta_{\psi} \sgn^0(w_+) \epsilon_{\psi_+}(x_+)
\end{eqnarray*}
and similarly for the pair $(G_+^{\vee},\psi_+^{\vee})$ we have:
\begin{eqnarray*}
r_{\psi^{\vee}_+}(w^{\vee}_+) \epsilon^1_{\psi^{\vee}_+}(u^{\vee}_+) = \delta_{\psi} \sgn^0(w^{\vee}_+) \epsilon_{\psi^{\vee}_+}(x^{\vee}_+).
\end{eqnarray*}

\noindent Thus the sign $\delta_{\psi}$ accounts for the deficiency of not knowing the spectral sign lemma for the pair $(G_+,\psi_+)$ and $(G_+^{\vee},\psi_+^{\vee})$ at this point.

Finally we follow the notation of \cite{A1} ({\it c.f.} the discussion before Lemma 5.3.2 of {\it loc. cit.}) and write the stable linear form $S^G_{\disc,\psi^N}$ on $\mathcal{H}(G)$ as $f \mapsto f^G(\Gamma)$. Thus proposition 6.1.4 give the identity:
\[
f^{G}(\psi) = f^G(\Gamma) + f^G(\Lambda)
\]
(noting that both $|\mathcal{S}_{\psi}|$ and $\epsilon^G(\psi)$ are equal to $+1$ since $\psi$ is a simple parameter).

As in the previous situation where $r >1$, we introduce the Levi subgroup
\[
L_+ \cong G \times G_{E/F}(N)
\]
of $G_+$ and equip $L_+$ the structure of Levi sub-datum of $G_+$. We then have the linear form:
\[
f \mapsto f^{L_+}(\Gamma \times \Lambda)
\]
on $\mathcal{H}(G_+)$.
 
\begin{proposition}
Suppose that $r = 1$ in (6.3.1), i.e. $\psi^N \in \widetilde{\mathcal{F}}_{\simp}(N)$, then for any compatible family of functions $\mathcal{F}^+ = \{ f^* \in \mathcal{H}(G^*)| \,\ G^* \in \widetilde{\mathcal{E}}_{\simp}(N_+)   \}$, we have the identity:
\begin{eqnarray}
& & \sum_{G^* \in \widetilde{\mathcal{E}}_{\simp}(N_+)} \widetilde{\iota}(N_+,G^*)  \tr R^{G^*}_{\disc,\psi_+^{N_+}}(f^*) \\
& & +     \frac{1}{8} \big( (f^{\vee})^{M_+^{\vee}}(\psi^{\vee}_{M_+^{\vee}}) - \delta_{\psi} f^{\vee}_{G^{\vee}_+}(\psi^{\vee}_+,x^{\vee}_{+,1})      \big)   + \frac{1}{2} f^{L_+}(\Gamma \times \Lambda)              \nonumber \\
&=&  \frac{1}{8} \big( f^{\prime}_{G_+}(\psi_+,s_{\psi_+} x_{+,1}) - \delta_{\psi} f_{G_+}(\psi_+,x_{+,1})      \big). \nonumber
\end{eqnarray}
(As in proposition 6.3.1 here $f$ is the function in the compatible family associated to $G_+$.)
\end{proposition}
\begin{proof}
Again the strategy of the proof is similar to that of proposition 6.2.2, with additional justification as in the proof of proposition 6.3.1, due to the fact that we are working with parameters of degree $N_+ = 2N > N$. Since the justification is similar, we will only indicate the additional complication in the present situation.

We again start with identity (6.2.20):
\begin{eqnarray}
\sum_{G^* \in \widetilde{\mathcal{E}}_{\ellip}(N_+)} \widetilde{\iota}(N_+,G_+) \leftexp{0}{S^{G^*}_{\disc,\psi_+^{N_+}}(f^*)}=0.
\end{eqnarray} 

We again claim that
\begin{eqnarray}
\leftexp{0}{S^{G^*}_{\disc,\psi_+^{N_+}}(f^*)}=0
\end{eqnarray}
for 
\begin{eqnarray}
G^* \in \widetilde{\mathcal{E}}_{\ellip}(N_+) \smallsetminus \widetilde{\mathcal{E}}_{\simp}(N_+),\,\ G^* \neq G \times G^{\vee}
\end{eqnarray}
Indeed, the same argument as in the proof of proposition 6.3.1 applies without change to this case also. To analyze the case $G^*=G \times G^{\vee}$, we apply proposition 6.1.4. Hence for $f^* = f_1 \times f_2^{\vee} \in \mathcal{H}(G^*)$, we have (noting that $\widetilde{\mathcal{F}}(G \times G^{\vee},\psi_+^{N_+})$ is empty):
\begin{eqnarray*}
& & \leftexp{0}{S^{G \times G^{\vee}}_{\disc, \psi_+^{N_+} }(f^*)}  = S^{G \times G^{\vee}}_{\disc, \psi_+^{N_+} } (f^*)     \\
&=& S^{G \times G^{\vee}}_{\disc, \psi^N \times \psi^N}(f^*) \nonumber \\
&=& S^{G}_{\disc,\psi^N}(f_1) \cdot S^{G^{\vee}}_{\disc,\psi^N}(f_2^{\vee})  \nonumber \\
&=& (f_1)^G(\Gamma) \cdot f_2^{\vee,L^{\vee}}(\Lambda) \nonumber \\
&=& f^{L_+}(\Gamma \times \Lambda).\nonumber
\end{eqnarray*}
Thus the only contribution to (6.3.20) from $G^* \in \widetilde{\mathcal{E}}_{\ellip}(N_+) \smallsetminus \widetilde{\mathcal{E}}_{\simp}(N_+)$ is from $G^* = G \times G^{\vee}$, given by:
\begin{eqnarray}
& & \widetilde{\iota}(N_+,G \times G^{\vee}) \leftexp{0}{S^{G \times G^{\vee}}_{\disc, \psi_+^{N_+} }(f^*)}\\
&=& \frac{1}{4} f^{L_+}(\Gamma \times \Lambda). \nonumber
\end{eqnarray} 

To obtain (6.3.19) it suffices to obtain the following two expressions:
\begin{eqnarray}
& & \leftexp{0}{S^{G_+}_{\disc, \psi_+^{N_+}}(f)} - \tr R^{G_+}_{\disc,\psi_+^{N_+}}(f) \\
& & \,\ + \frac{1}{4}(f^{G \times G})(\Gamma \times \Gamma) - f^{G \times G}(\psi \times \psi)) \nonumber  \\
&=& -\frac{1}{4} \big( f^{\prime}_{G_+}(\psi_+,s_{\psi_+} x_{+,1}) - \delta_{\psi} f_{G_+}(\psi_+,x_{+,1})      \big) \nonumber
\end{eqnarray}
and
\begin{eqnarray}
& & \leftexp{0}{S^{G_+^{\vee}}_{\disc, \psi_+^{N_+}}(f^{\vee})} - \tr R^{G_+^{\vee}}_{\disc,\psi_+^{N_+}}(f^{\vee})  \\
& & \,\ + \frac{1}{4} f^{L \times L}(\Lambda \times \Lambda) \nonumber \\
&=& \nonumber  \frac{1}{4} \big( (f^{\vee})^{M_+^{\vee}}(\psi^{\vee}_{M_+^{\vee}}) - \delta_{\psi} f^{\vee}_{G^{\vee}_+}(\psi^{\vee}_+,x^{\vee}_{+,1})      \big).           \nonumber     
\end{eqnarray}
Indeed since $\widetilde{\iota}(N_+,G_+)=\widetilde{\iota}(N_+,G_+^{\vee})=1/2$, we see that (6.3.19) would follow from substituting (6.3.21)-(6.3.25) to (6.3.20), together with the computation:

\begin{eqnarray*}
& & \frac{1}{4} f^{L_+}(\Gamma \times \Lambda) - \frac{1}{8}(f^{G \times G}(\Gamma \times \Gamma) - f^{G \times G}(\psi \times \psi)) -\frac{1}{8} f^{L \times L}(\Lambda \times \Lambda)\\
&=& \frac{1}{4} f^{L_+}(\Gamma \times \Lambda) - \frac{1}{8}(- 2 f^{L_+}(\Gamma \times \Lambda) - f^{L \times L}(\Lambda \times \Lambda))  -\frac{1}{8} f^{L \times L}(\Lambda \times \Lambda)\\
&=& \frac{1}{2} f^{L_+}(\Gamma \times \Lambda).
\end{eqnarray*}

\noindent The derivation of (6.3.24) and (6.3.25) are parallel to the derivation of (6.2.11) and (6.2.14) respectively. The two new phenomenon are: firstly there is the sign $\delta_{\psi}$ occuring as coefficients of the spectral distributions $f_{G_+}(\psi_+,x_+)$ and $f^{\vee}_{G_+^{\vee}}(\psi_+^{\vee},x_+^{\vee})$, which account for not knowing the validity of the spectral sign lemma for the pairs $(G_+,\psi_+)$ and $(G_+^{\vee},\psi_+^{\vee})$ at this point; secondly there are extra terms occuring on the left hand side of (6.3.24) and (6.3.25), due to not knowing the stable multiplicity formula for the distributions $S^{G \times G}_{\disc,\psi_+^{N_+}}$ and $S^{G^{\vee} \times G^{\vee}}_{\disc,\psi_+^{N_+}}$ at this point, and hence the occurence of extra terms in the endoscopic expansions for the distributions $I^{G_+}_{\disc,\psi_+^{N_+}}$ and $I^{G_+^{\vee}}_{\disc,\psi_+^{N_+}}$, {\it c.f.} the proof of proposition 6.3.1. In other words, the extra terms comes from:
\begin{eqnarray*}
& & \iota(G_+,G \times G) \leftexp{0}{S^{G \times G}_{\disc,\psi_+^{N_+}}(f^{G \times G})} \\
&=& \frac{1}{4} (S^{G \times G}_{\disc,\psi^N \times \psi^N}(f^{G \times G}) - f^{G \times G}(\psi \times \psi) )\\
&=& \frac{1}{4} (f^{G \times G}(\Gamma \times \Gamma) - f^{G \times G}(\psi \times \psi) )
\end{eqnarray*}
and (here $f^{\vee}$ is the function in the compatible family occuring in (6.3.20) that is associated to $G_+^{\vee}$):
\begin{eqnarray*}
& & \iota(G_+^{\vee},G^{\vee} \times G^{\vee}) \leftexp{0}{S^{G^{\vee} \times G^{\vee}}_{\disc,\psi_+^{N_+}}(f^{\vee, G^{\vee} \times G^{\vee}})} \\
&=& \frac{1}{4} S^{G^{\vee} \times G^{\vee}}_{\disc,\psi_+^{N_+}}(f^{\vee, G^{\vee} \times G^{\vee}})\\
&=& \frac{1}{4} f^{\vee,L^{\vee} \times L^{\vee}}(\Lambda^{\vee} \times \Lambda^{\vee}) = \frac{1}{4} f^{L \times L}(\Lambda \times \Lambda).
\end{eqnarray*}
\end{proof}

For later purpose (to be used in the completion of global induction arguments in chapter 9) we also record the following:

\begin{corollary}
In the situation of proposition 6.3.3, we have:
\begin{eqnarray}
& & \\
& & S^{G_+}_{\disc,\psi_+^{N_+}}(f) = \tr R^{G_+}_{\disc,\psi_+^{N_+}}(f) + \frac{1}{4}(\delta_{\psi} f_{G_+}(\psi_+,x_{+,1})-f^{G \times G}(\Gamma \times \Gamma)) \nonumber
\end{eqnarray}
and
\begin{eqnarray}
& & \\
& & S^{G_+^{\vee}}_{\disc,\psi_+^{N_+}}(f^{\vee}) = \tr R^{G_+^{\vee}}_{\disc,\psi_+^{N_+}}(f^{\vee}) - \frac{1}{4}(\delta_{\psi} f^{\vee}_{G_+^{\vee}}(\psi_+^{\vee},x_{+,1}^{\vee}) + f^{L \times L}(\Lambda \times \Lambda)). \nonumber
\end{eqnarray}
\end{corollary}
\begin{proof}
Same as corollary 5.3.3 of \cite{A1}. For example (6.3.26) follows from the endoscopic expansion:
\begin{eqnarray*}
& & I^{G_+}_{\disc,\psi_+^{N_+}}(f) = S^{G_+}_{\disc,\psi_+^{N_+}}(f) + \iota(G_+,G \times G) S^{G \times G}_{\disc,\psi_+^{N_+}}(f^{G \times G})\\
&=&  S^{G_+}_{\disc,\psi_+^{N_+}}(f) + \frac{1}{4} f^{G \times G}(\Gamma \times \Gamma)
\end{eqnarray*}
and the corresponding spectral expansion, which reads as:
\[
 I^{G_+}_{\disc,\psi_+^{N_+}}(f) =  \tr R^{G_+}_{\disc,\psi_+^{N_+}}(f) + \frac{1}{4}\delta_{\psi} f_{G_+}(\psi_+,x_{+,1})
\]
Similarly for the proof of (6.3.27). 

\end{proof}

\subsection{Generic parameters with local constraints}

In this final subsection, we refine the propositions established in the previous subsections, in the special case of particular families $\widetilde{\mathcal{F}}$ consisting of generic parameters with serious local constraints at the archimedean places. In particular all the parameters considered in this section are generic. 

As we have noted before, the global intertwining relation is the obstruction to obtaining complete information about the spectral multiplicity and the stable multiplicity formula in section 6.2 and 6.3. For the special class of families to be considered in this subsection, we can obtain complete information about the parameters treated in section 6.2 and 6.3, without knowing the validity of the global intertwining relation {\it a priori}. The global information obtained in this subsection will in turn be used to establish the local results in section 7.

As in section 6.1, we have a family of parameters $\widetilde{\mathcal{F}}$ that is the graded semi-group generated by the simple parameters $\widetilde{\mathcal{F}}_{\simp}$ (with the simple parameters of degree less than or equal to $N$). Consider the following condition on $\widetilde{\mathcal{F}}$:

\begin{assumption} There is a non-empty set $V=V(\widetilde{\mathcal{F}})$ of archimedean valuations of $F$, that does not split in $E$, for which the following three conditions hold:

\bigskip

\noindent (6.4.1)(a) Suppose that $\psi^N \in \xi_* \widetilde{\mathcal{F}}^{\#}_{\simp}(G)$ for some $G =(G,\xi) \in \widetilde{\mathcal{E}}_{\simp}(N)$. Then for any $v \in V$ we have $\psi^N_v \in (\xi_v)_* \Psi_2(G_v)$.

\bigskip

\noindent (6.4.1)(b) Suppose that $\psi^N \in \xi_* \widetilde{\mathcal{F}}_2^{\simp}(G)$, for some $G=(G,\xi) \in \widetilde{\mathcal{E}}_{\ellip}(N)$. Then there is a valuation $v \in V$, such that $\psi^N$ does not lies in $(\xi_v^*)_*\Psi^+(G^*_v)$ for any $G^*_v  = (G^*_v,\xi_v^*) \in \widetilde{\mathcal{E}}_{v,\simp}(N)$ with $G^*_v \neq G_v$ (as elements of $\widetilde{\mathcal{E}}_{\simp,v}(N)$).

\bigskip

\noindent (6.4.1)(c) Suppose that $\psi^N \in \xi_* \widetilde{\mathcal{F}}^2_{\disc}(G)$ for some $G=(G,\xi) \in \widetilde{\mathcal{E}}_{\ellip}(N)$, with $\psi \in \widetilde{\mathcal{F}}^2_{\disc}(G)$ the parameter of $G$ defined by $\psi^N$. Then there is a valuation $v \in V$ such that the kernel of the composition of mappings
\begin{eqnarray*}
\mathcal{S}_{\psi}(G) \rightarrow \mathcal{S}_{\psi_v}(G_v) \rightarrow R_{\psi_v}(G_v)
\end{eqnarray*}
contains no element whose image in the global $R$-group $R_{\psi}=R_{\psi}(G)$ is regular.
\end{assumption}
\noindent (Here in the above $\widetilde{\mathcal{F}}_{2}^{\simp}(G)$ and $\widetilde{\mathcal{F}}_{\disc}^2(G)$ are defined similarly to $\widetilde{\mathcal{F}}_{\ellip}^2(G)$ in section 6.2.)
\bigskip

Thus the parameters in $\widetilde{\mathcal{F}}$ satisfy rather serious local constraints at the archimedean places. We will construct such parameters in section 7, as an application of the simple version of the invariant trace formula. In fact, the construction shows that we can even assume that all the parameters in $\widetilde{\mathcal{F}}$ are generic parameters. This is sufficients for the purpose of establishing the propositions in this subsection, which be applied to the proof of the local classification of tempered representations in section 7. However, in accordance with the notations of the previous subsections, we will still use the notations $\psi^N,\psi$ for such parameters instead of the notations $\phi^N,\phi$.

Henceforth we fix such a family of generic parameters $\widetilde{\mathcal{F}}$ satisfying assumption 6.4.1 in this subsection. In accordance with the previous induction arguments we assume that all the local and global theorems are valid for parameters of degree less than $N$. The first thing is to show the validity of hypothesis 6.1.1:

\begin{lemma}
Hypothesis 6.1.1 holds for the family of parameters $\widetilde{\mathcal{F}}$. 
\end{lemma} 
\begin{proof}
By part (b) of remark 6.1.2, it suffices to consider part (a) of hypothesis 6.1.1, in the following setting where $G=(G,\xi) \in \widetilde{\mathcal{E}}_{\ellip}(N)$, and that $\psi \in \widetilde{\mathcal{F}}^{\#}_2(G)$. Given the pair $(G,\psi)$, write as usual write $\psi^N = \xi_* \psi$. We decompose the set of valuations of $F$ as a disjoint union
\[
V \coprod U \coprod V_{\un}
\]
with $V$ being the set of archimedean valuations for $\widetilde{\mathcal{F}}$ as in Assumption 6.4.1, $U$ being a finite set, and where $(G,\psi)$ is unramified at every places in $V_{\un}$. Consider a function:
\[
\widetilde{f} = \widetilde{f}_V \cdot \widetilde{f}_U \cdot \widetilde{f}_{\un}  , \,\ \widetilde{f} \in \widetilde{\mathcal{H}}(N)
\]
adapted to this decomposition. Without loss of generality assume $\widetilde{f}_{\un}$ is decomposable: 
and
\[
\widetilde{f}_{\un} = \prod_{v \in V_{\un}} \widetilde{f}_{v}
\]
with $\widetilde{f}_v \in \widetilde{\mathcal{H}}_v(N,\widetilde{K}_v(N))$, the spherical Hecke module of $\widetilde{\mathcal{H}}_v(N)$. Here $\widetilde{K}_v(N): =K_v(N) \rtimes \theta$, with $K_v(N)$ being the standard maximal compact subgroup of $G_{E_v/F_v}(N)$, which in particular is $\theta$-invariant.

In the first step we allow $\widetilde{f}_U \in \widetilde{\mathcal{H}}_U(N) :=\bigotimes_{v \in U} \widetilde{\mathcal{H}}_v(N)$ to vary, and fix a chocie of functions $\widetilde{f}_V$ and $\widetilde{f}_{\un}$.

For $\widetilde{f}_{\un} = \prod_{v \in V_{\un}} \widetilde{f}_{v}$ we can simply take $\widetilde{f}_v$ to be the characteristic function of $\widetilde{K}_v(N)$. Then in particular we have
\[
\widetilde{f}_{\un,N}(\psi^N_{\un}) = \prod_{v \in V_{\un}} \widetilde{f}_{v,N}(\psi^N_v)   \neq 0.
\]

As for the choice of $\widetilde{f}_V$, we use condition (6.4.1)(a,b). Recall that we are assuming that $\psi \in \widetilde{\mathcal{F}}^{\#}_2(G)$. It then follows from (6.4.1)(a) in the case $\psi \in \widetilde{\mathcal{F}}^{\#}_{\simp}(G)$, and (6.4.1)(b) in the case $\psi \in \widetilde{\mathcal{F}}_2^{\simp}(G)$, that there exists a valuation $w \in V$, such that $\psi^N_w$ does not lie in $(\xi_w^*)_*\Phi(G_w^*)$ for any $G_w^* =(G_w^*,\xi_w^*) \in \widetilde{\mathcal{E}}_{w,\simp}(N)$ with $G_w^* \neq G_w$. 

We can choose a function $\widetilde{f}_w \in \widetilde{\mathcal{H}}_w(N)$, such that $\widetilde{f}^{G^*_w}_w$ is identifically zero for any such $G^*_w \neq G_w$, but such that the value $\widetilde{f}_{w,N}(\psi^N_w)$ is non-zero. This follows from part (a) of proposition 3.1.1, and the twisted spectral transfer result of Mezo \cite{Me} and Shelstad \cite{Sh3}. At the other places $v \in V \smallsetminus \{w\}$, we can just choose $\widetilde{f}_v$ such that $\widetilde{f}_{v,N}(\psi^N_v)$ is non-zero. The function $\widetilde{f}_V$ then has the property that it value
\[
\widetilde{f}_{V,N}(\psi^N_V) = \prod_{v \in V} \widetilde{f}_{v,N}(\psi^N_v)
\]
is non-zero, while the twisted transfer
\[
\widetilde{f}_V^{G^*} = \prod_{v \in V} \widetilde{f}_v^{G^*_v}
\]
is identifically zero for $G^* \in \widetilde{\mathcal{E}}_{\simp}(N)$ with $G^* \neq G$ (note that since $w$ does not split in $E$, we have $G^* \neq G$ implies $G^*_w \neq G_w$).

We then apply the identity (4.3.13), which we recall here as:
\begin{eqnarray}
\widetilde{I}^N_{\disc,\psi^N}(\widetilde{f}) = \sum_{G^* \in \widetilde{\mathcal{E}}_{\ellip}(N)} \widetilde{\iota}(N,G^*) \widehat{S}^{G^*}_{\disc,\psi^N}(\widetilde{f}^{G^*})
\end{eqnarray}
to the function
\begin{eqnarray}
\widetilde{f} = \widetilde{f}_V \cdot \widetilde{f}_U \cdot \widetilde{f}_{\un}
\end{eqnarray}
 with $\widetilde{f}_U \in \widetilde{\mathcal{H}}_U(N)$ being allowed to vary. 

We claim that for any $G^* \in \widetilde{\mathcal{E}}_{\ellip}(N)$ with $G^* \neq G$, the term $\widehat{S}_{\disc,\psi^N}^{G^*}(\widetilde{f}^{G^*})$ vanishes. Indeed, first suppose that $G^* = (G^*,\xi^*) \in \widetilde {\mathcal{E}}_{\ellip}(N) \smallsetminus \widetilde{\mathcal{E}}_{\simp}(N)$. Then from induction hypothesis the stable multiplicity formula is valid for $\psi^N$ with respect to $G^*$. But since $\psi^N \in \xi_* \widetilde{\mathcal{F}}^{\#}_2(G)$, we know that  $\psi^N \notin (\xi^*)_* \Psi(G^*)$, hence we have the vanishing of $S_{\disc,\psi^N}^{G^*}$. 

Next suppose that $G^* \in \widetilde{\mathcal{E}}_{\simp}(N)$. Then by the chocie of $\widetilde{f}_V$ above we have $\widetilde{f}_V^{G^*} \equiv 0$, hence
\begin{eqnarray*}
& & \widehat{S}^{G^*}_{\disc,\psi^N}(\widetilde{f}^{G^*})\\
&=& \widehat{S}^{G^*}_{\disc,\psi^N}(\widetilde{f}_V^{G^*} \cdot \widetilde{f}_U^{G^*} \cdot \widetilde{f}_{\un}^{G^*}) =0
\end{eqnarray*}
as required.

Thus from (6.4.1) we obtain
\begin{eqnarray}
\widetilde{I}^N_{\disc,\psi^N}(\widetilde{f}) =  \widetilde{\iota}(N,G) \widehat{S}^{G}_{\disc,\psi^N}(\widetilde{f}^{G}).
\end{eqnarray}

The spectral distribution $\widetilde{I}_{\disc,\psi^N}$ can be analyzed as follows. Since $\psi^N \in \widetilde{\Psi}_{\ellip}(N)$, we can write $\psi^N$ in the standard way:
\begin{eqnarray}
\psi^N = \psi_1^{N_1} \boxplus \cdots \boxplus \psi_r^{N_r}
\end{eqnarray}
with $\psi_i^{N_i} \in \widetilde{\Psi}_{\simp}(N_i)$ mutually distinct. Then $\psi^N \in \Psi_2(\widetilde{M}^0)$, where $\widetilde{M}^0$ is the standard Levi subgroup of $G_{E/F}(N)$:
\begin{eqnarray}
\widetilde{M}^0 = G_{E/F}(N_1) \times \cdots \times G_{E/F}(N_r)
\end{eqnarray}
and there is a unique element $w \in  W_{\psi^N,\reg}$ such that $w$ induces the outer automorphism $\theta(N_i)$ of each general linear factor of $\widetilde{M}^0$. We then have (as a special case of the analysis in section 5.4-5.5):
\begin{eqnarray}
\widetilde{I}^N_{\disc,\psi^N}(\widetilde{f})  &=& |\det(w-1)_{\mathfrak{a}_{\widetilde{M}^0}^{\widetilde{G}(N)}  } |^{-1} \widetilde{r}^N_{\psi^N}(w) \widetilde{f}_N(\psi^N) \\
&=& \frac{1}{2^r} \widetilde{f}_N(\psi^N).    \nonumber 
\end{eqnarray}
for any $\widetilde{f} \in \widetilde{\mathcal{H}}(N)$. Here we are using the fact that
\[
 |\det(w-1)_{\mathfrak{a}_{\widetilde{M}^0}^{\widetilde{G}(N)}  } | = 2^r
\]
which is immediate from the shape of $\widetilde{M}^0$ in (6.4.5); and that
\[
\widetilde{r}^N_{\psi^N}(w) =1
\]
which again follows from the discussions in section 5.8.

Combining (6.4.3) and (6.4.6) we have
\begin{eqnarray}
\widehat{S}^G_{\disc,\psi^N}(\widetilde{f}^G) = \frac{1}{\widetilde{\iota}(N,G) 2^r} \widetilde{f}_N(\psi^N), \,\ \widetilde{f} \in \widetilde{\mathcal{H}}(N).
\end{eqnarray}

\noindent In particular with our current choice of $\widetilde{f} = \widetilde{f}_V \cdot \widetilde{f}_U \cdot \widetilde{f}_{\un}$, we have
\begin{eqnarray*}
& & \widehat{S}^G_{\disc,\psi^N}(\widetilde{f}_V^G \cdot \widetilde{f}_U^G \cdot \widetilde{f}_{\un}^G) \\
&=& \frac{1}{\widetilde{\iota}(N,G) 2^r} \widetilde{f}_{V,N}(\psi^N) \cdot  \widetilde{f}_{U,N}(\psi^N) \cdot  \widetilde{f}_{\un,N}(\psi^N)
\end{eqnarray*} 

\noindent In particular since $ \widetilde{f}_{V,N}(\psi^N)$ and $ \widetilde{f}_{\un,N}(\psi^N)$ are non-zero by the choice of $\widetilde{f}_V$, we see that $  \widetilde{f}_{U,N}(\psi^N)$ is the pull-back of a stable linear form, defined on the image of the twisted transfer
\begin{eqnarray}
\widetilde{\mathcal{H}}_U(N) \rightarrow \mathcal{S}_U(G) = \bigotimes_{v \in U} \mathcal{S}_v(G_v)
\end{eqnarray}
and it follows that we can write
\begin{eqnarray}
\widetilde{f}_{U,N}(\psi^N) = \widetilde{f}_U^G(\psi), \,\ \widetilde{f}_U \in \widetilde{\mathcal{H}}_U(N)
\end{eqnarray}
with $f_U \mapsto f_U^G(\psi)$ being a stable linear form defined on the image of (6.4.8).

\noindent On the other hand, since $V$ consists of archimedean places (that does not split in $E$), the twisted transfer results of \cite{Me} and \cite{Sh3} shows that (6.4.8) also holds with $U$ replaced by $V$, in particular the existence of the stable linear form $f_V^G(\psi)$ on $\widetilde{\mathcal{H}}_V(N)$. 

Since the set $U$ can be taken to be arbitrarily large, we see that there is a stable linear form $f \mapsto f^G(\psi)$, defined on the image of the twisted transfer
\begin{eqnarray}
\widetilde{\mathcal{H}}(N) \rightarrow \mathcal{S}(G)
\end{eqnarray} 
such that for any $\widetilde{f} \in \widetilde{\mathcal{H}}(N)$ we have
\begin{eqnarray}
\widetilde{f}_N(\psi^N) = \widetilde{f}^G(\psi).
\end{eqnarray}

\bigskip

Now first suppose that $G \in \widetilde{\mathcal{E}}_{\simp}(N)$. Then by proposition 3.1.1 (b), the twisted transfer mapping (6.4.10) is surjective, and hence the stable linear form $f^G(\psi)$ is defined on the whole of $\mathcal{S}(G)$, and it follows that hypothesis (6.1.1) holds for the pair $(G,\psi)$ (which concerns only the case of equation (6.1.6)).  

It remains to consider the case where $G \in \widetilde{\mathcal{E}}_{\ellip}(N) \smallsetminus \widetilde{\mathcal{E}}_{\simp}(N)$ is composite, thus we have
\[
G = G^1 \times G^2, \,\ G^i \in \widetilde{\mathcal{E}}_{\simp}(N^i), \,\ N^i < N
\] 
with $G^1=(G^1,\xi^1)$ and $G^2=(G^2,\xi^2)$ being of opposite parity as endoscopic data. We have $\psi = \psi^1 \times \psi^2$, with $\psi^i \in \Psi_2(G^i)$. In this case we define the linear form by the condition (6.1.7). Thus we define $f^G(\psi) = f^G(\psi^1 \times \psi^2)$ for $f \in \mathcal{H}(G)$ by the the following condition: if 
\[
f^G = f_1^{G^1} \times f_2^{G^2}
\]
is decomposable then
\[
f^G(\psi) = f^G(\psi^1 \times \psi^2)  := f_1^{G^1}(\psi^1) \times f_2^{G^2}(\psi^2).
\]

\noindent The main point in this case then is to show the validity of (6.1.6).

\bigskip

We first note that since $G$ is composite, we have $\widetilde{\iota}(N,G)=1/4$. Hence
\begin{eqnarray}
& & \widetilde{\iota}(N,G) \cdot 2^r \\
&=& 2^{r-2} = |\mathcal{S}_{\psi}| = |\mathcal{S}_{\psi_1}|\cdot |\mathcal{S}_{\psi_2}| \nonumber
\end{eqnarray}
as is easily verified.

\noindent On the other hand, by the induction hypothesis, the stable multiplicity formula is valid for $(G^i,\psi^i)$. Hence we have (denoting $\psi^{i,N^i}:= \xi^i_* \psi^i$):
\begin{eqnarray}
S^{G^i}_{\disc,\psi^{i,N^i}}(f_i) = \frac{1}{|\mathcal{S}_{\psi^i}|} f_i^{G^i}(\psi^i), \,\ f_i \in \mathcal{H}(G^i)
\end{eqnarray}
(recall that $\psi$ and hence $\psi^i$ are generic parameters, hence there involves no $\epsilon$-factor in the stable multiplicity formula). 

\bigskip

\noindent  By (6.4.13) and (6.4.12), we have, for $f \in \mathcal{H}(G)$ such that $f^G = f_1^{G^1} \times f_2^{G^2}$ is decomposable:
\begin{eqnarray}
& & S^G_{\disc,\psi^N}(f) = \widehat{S}^G_{\disc,\psi^N}(f^G) \\
&=& \widehat{S}^{G^1}_{\disc,\psi^{1,N^1}} (f_1^{G^1}) \cdot  \widehat{S}^{G^2}_{\disc,\psi^{2,N^2}} (f_2^{G^2})   \nonumber \\
&=& \frac{1}{|\mathcal{S}_{\psi^1}|} f_1^{G^1}(\psi^1) \cdot  \frac{1}{|\mathcal{S}_{\psi^2}|} f_2^{G^2}(\psi^2) \nonumber \\
&=&  (\widetilde{\iota}(N,G) 2^r)^{-1} \cdot f_1^{G^1}(\psi^1) \cdot  f_2^{G^2}(\psi^2) \nonumber \\
&=&  (\widetilde{\iota}(N,G) 2^r)^{-1} f^G (\psi^1 \times \psi^2) = (\widetilde{\iota}(N,G) 2^r)^{-1} f^G(\psi) \nonumber
\end{eqnarray}
and hence (6.4.14) also holds for all $f \in \mathcal{H}(G)$. 

\noindent  We thus obtain, for $\widetilde{f} \in \widetilde{\mathcal{H}}(N)$, by combining (6.4.14) with (6.4.7) that (replacing $f^G$ by $\widetilde{f}^G$ in (6.4.14))
\begin{eqnarray}
\widetilde{f}^G(\psi) = \widetilde{f}_N(\psi^N)
\end{eqnarray}
as required.
\end{proof}

\bigskip

\begin{rem}
\end{rem}
As discussed in remark 6.1.2, it follows from lemma 6.4.2 that the statement of hypothesis 6.1.1 actually holds for parameters $\psi \in \widetilde{\mathcal{F}}(G) \smallsetminus \widetilde{\mathcal{F}}_2(G)$ (for which only part (a) of hypothesis 6.1.1 is relevant).

\bigskip

We can now come back to the parameters treated in section 6.2:

\begin{proposition}
Suppose that
\[
(G,\psi), \,\ G = (G,\xi) \in \widetilde{\mathcal{E}}_{\simp}(N), \,\ \psi \in \widetilde{\mathcal{F}}^2_{\ellip}(G),
\]
is as in proposition 6.2.1, with $\widetilde{\mathcal{F}}$ being our family of generic parameters that satisfy Assumption 6.4.1. The we have (as usual with $\psi^N = \xi_* \psi$):
\begin{eqnarray}
\tr R^{G^*}_{\disc,\psi^N}(f^*) =0 = \leftexp{0}{S^{G^*}_{\disc,\psi^N}(f^*)}, \,\ f^* \in \mathcal{H}(G^*)
\end{eqnarray}
for every $G^* \in \widetilde{\mathcal{E}}_{\simp}(N)$, while the right hand side of (6.2.2) vanishes (note that since we are dealing with generic parameters here, the $\epsilon$-factors on the right hand side of (6.2.2) are all equal to one).
\end{proposition}
\begin{proof}
The proof is the same as the proof of lemma 5.4.3 of \cite{A1}, so again we will be brief with the argument. Let $v \in V=V(\widetilde{\mathcal{F}})$ be as in condition (6.4.1)(c), with respect to the $G \in \widetilde{\mathcal{E}}_{\simp}(N)$ that we are currently considering. Consider a decomposable function
\[
f = f_v f^v \in \mathcal{H}(G), \,\ f_v \in \mathcal{H}(G_v), f^v \in \mathcal{H}(G(\mathbf{A}_F^v))
\]
and similarly write the (generic) paramater $\psi \in \Psi(G)$ symbolically as
\[
\psi = \psi_v \psi^v
\]
then the terms on the right hand side of (6.2.2) can be factorized as:
\begin{eqnarray*}
& & f_G(\psi,x) = f_{v,G}(\psi_v,x_v) f^v_G(\psi^v,x^v) \\
& &  f^{\prime}_G(\psi,x) = f^{\prime}_{v,G}(\psi_v,x_v) (f^v)^{\prime}_G(\psi^v,x^v) 
\end{eqnarray*}
and thus the right hand side of (6.2.2) can be written as:
\begin{eqnarray}
& & \\
& & c \sum_{x \in \mathcal{S}_{\psi,\ellip}} \big( (f^v)^{\prime}_G(\psi^v,x^v) f^{\prime}_{v,G}(\psi_v,x_v) - f^v_G(\psi^v,x^v) f_{v,G}(\psi_v,x_v) \big). \nonumber
\end{eqnarray}

\noindent Arguing in the proof of lemma 5.4.3 of \cite{A1}, we can write (6.4.17) in the form
\begin{eqnarray}
\sum_{\tau_v \in T(G_v)} d(\tau_v,f^v) f_{v,G}(\tau_v)
\end{eqnarray}
with $T(G_v)$ being the set introduced in proposition 4.3.9. In the present context the coefficients $d(\tau_v,f^v)$ can be non-zero only for elements $\tau_v \in T(G_v)$ that are $W_0^G$-orbits of triples of the form:
\[
(M_v,\pi_v,r_v(x_v))
\]
where $M_v$ is the Levi subgroup of $G_v$ such that $\Psi_2(M_v,\psi^N_v)$ is non-empty ($M_v$ is unique up to conjugation, and as usual regarded as an endoscopic datum in $\widetilde{\mathcal{E}}(N)$ over $F_v$)). Note that $M_v$ is proper in $G_v$ but can be smaller than the localization at $v$ of the global Levi subgroup $M$ of $G$ associated to $\psi$. The representations $\pi_v$ runs over the local $L$-packet $\Pi_{\psi_v}(M_v)$, and $r_v(x)$ stands for the image of $x_v$ in the local $R$-group $R_{\psi_v}$ (with $x \in \mathcal{S}_{\psi,\ellip}$).

\noindent Furthermore, since $x \in \mathcal{S}_{\psi,\ellip}$ in (6.4.17), the image of $x$ in the global $R$-group $R_{\psi}=R_{\psi}(G)$ is a regular element. Hence by condition (6.4.1)(c) the image of $x_v$ in the local $R$-group $R_{\psi_v}$ is non-trivial. Hence it follows as in {\it loc. cit.} that the coefficient $f(\tau_v,f^v)$ in (6.4.18) is zero, for $\tau_v \in T(G_v)$ any element of the form $(M_v,\pi_v,1)$.  

On the other hand, the left hand side of (6.2.2) is, as usual, a linear combination with positive coefficients if irreducible characters of representations on $G(\mathbf{A}_F)$ for $G \in \widetilde{\mathcal{E}}_{\simp}(N)$. Thus proposition 4.3.9 applies. We conclude firstly that all the coefficients $d(\tau_v,f^v)$ vanish, and hence so is (6.4.18) and this the right hand side of (6.2.2). Secondly we also conclude that all the coefficients of the irreducible characters on the left hand side of (6.2.2) vanish. In particular we have
\begin{eqnarray}
R^{G^*}_{\disc,\psi^N} \equiv 0, \,\ G^* \in \widetilde{\mathcal{E}}_{\simp}(N).
\end{eqnarray}

Finally as we observe in the proof of proposition 6.2.1 that the difference
\begin{eqnarray}
\tr R^{G^*}_{\disc,\psi^N}(f^*) - \leftexp{0}{S^{G^*}_{\disc,\psi^N}(f^*)}, \,\ f^* \in \mathcal{H}(G^*)
\end{eqnarray}
vanishes unless $G^*=G$, in which case it is equal to the right hand side of (6.2.2). We have just seen that this expression vanishes. Thus (6.4.20) vanishes for all $G^* \in \widetilde{\mathcal{E}}_{\simp}(N)$, and the conclusion follows from (6.4.19).
\end{proof}

\begin{proposition}
Then $(G,\psi)$ be as in proposition 6.2.2. Thus $G \in \widetilde{\mathcal{E}}_{\simp}(N)$ such that $\psi \in \widetilde{\mathcal{F}}_{\ellip}^2(G)$ with $\psi^N$ as in (6.2.9). Then
\[
\tr R^{G^*}_{\disc,\psi^N} (f^*) = 0 = \leftexp{0}{S^{G^*}_{\disc,\psi^N}(f^*) } , \,\ G^* \in \widetilde{\mathcal{E}}_{\simp}(N).
\]
Furthermore the right hand side of (6.2.10) vanishes, and so does the expression
\begin{eqnarray}
\frac{1}{8}\big( (f^{\vee})^{M^{\vee}}(\psi_{1,M^{\vee}}) - f^{\vee}_{G^{\vee}}(\psi^{\vee},x_1^{\vee}) \big)
\end{eqnarray}
on the left hand side of (6.2.10).
\end{proposition} 
\begin{proof}
The proof is similar to the proof of proposition 6.4.4. The only difference is the appearance of the extra term (6.4.21) occuring on the left hand side of (6.2.10). But by descent as in the proof of proposition 5.7.4, the term (6.4.21) can be written as (in the evident notation):
\begin{eqnarray}
\frac{1}{8} \tr \big(  (1 - R_{P^{\vee}}(w_1^{\vee}, \widetilde{\pi}_{\psi_{1,M^{\vee}}} ,\psi_{1,M^{\vee}}  )  ) \mathcal{I}_{P^{\vee}} (\pi_{\psi_{1,M^{\vee}}}  ,f^{\vee}) \big).
\end{eqnarray}
Here referring to the notation as in the situation of proposition 6.2.2 the element $w_1^{\vee}$ is the unique regular element of $W_{\psi^{\vee}} = W_{\psi^{\vee}}^0$, in particular is an element of order two. Thus the eigenvalues of the intertwining operator $R_{P^{\vee}}(w_1^{\vee}, \widetilde{\pi}_{\psi_{1,M^{\vee}}} ,\psi_{1,M^{\vee}}  )$ are $\pm 1$. And it follows that the expression (6.4.22) is also a linear combination with positive coefficients of irreducible characters of representations of $G^{\vee}$. We conclude by the same argument as in proposition 6.4.4. As a bonus, we also see that the intertwining operator $R_{P^{\vee}}(w_1^{\vee}, \widetilde{\pi}_{\psi_{1,M^{\vee}}} ,\psi_{1,M^{\vee}}  )$ is identifically equal to one.
\end{proof}

\bigskip

We now apply the same method to the parameters treated in section 6.3.

\bigskip

\begin{proposition}
In the situation of proposition 6.3.1, we have for $G^* \in \widetilde{\mathcal{E}}_{\simp}(N_+)$
\begin{eqnarray}
\tr R^{G^*}_{\disc,\psi_+^{N_+}}(f^*) =0 = \leftexp{0}{ S^{G^*}_{\disc,\psi_+^{N_+}}(f^*) }, \,\ f^* \in \mathcal{H}(G^*).
\end{eqnarray}
Furthermore the linear form $\Lambda$ and the right hand side of (6.3.3) vanishes. In particular the stable multiplicity formula is valid for $\psi^N$.
\end{proposition}
\begin{proof}
As in the statement of proposition 6.3.1, the premise for the validity of equation (6.3.3) is the validity of the stable multiplicity formula for all the parameters in $\widetilde{\mathcal{F}}(N) \smallsetminus \widetilde{\mathcal{F}}_{\ellip}(N)$. But this is exactly the content of proposition 6.4.4 and 6.4.5, coupled with proposition 5.7.4. Thus equation equation (6.3.3) is valid.

We can apply essentially the same argument as in the proof of proposition 6.4.4; the only difference being the term $b_+f^{L_+}(\psi_1 \times \Lambda)$ occuring on the left hand side of (6.3.3) that needs to be taken care of. As in the proof of Lemma 5.4.5 of \cite{A1}, put
\[
G_1^{\vee} = G_1 \times G^{\vee} \in \widetilde{\mathcal{E}}_{\ellip}(N_+).
\]
and $L_1^{\vee}$ the Levi sub-datum of $G_1^{\vee}$ with the same underlying group as $L_+$. Furthermore for $f_1$ the function associated to $G_1^{\vee}$ in the compatible family occuring in (6.3.3), we have the equality:
\[
f_1^{L_1^{\vee}}(\psi_1 \times \Lambda^{\vee}) = f^{L_+}(\psi_1 \times \Lambda).
\]
From proposition 6.1.4, the linear form $\Lambda^{\vee}$ is a unitary character on $G^{\vee}$. On the other hand, the stable linear form defined by the simple parameter $\psi_1$ is also a unitary character; this follows from the stable multiplicity formula for $(G_1,\psi_1)$, together with the equality 
\[
S^{G_1}_{\disc,\psi_1^{N_1}} \equiv \tr R^{G_1}_{\disc,\psi_1^{N_1}}
\] 
which follows since $\psi_1$ is a simple parameter. Thus the linear form $f_1 \mapsto f_1^{L_1^{\vee}}(\psi_1 \times \Lambda^{\vee})$ is a unitary character on $G_1^{\vee}$, hence is a linear combination with non-negative coefficients of irreducible admissible representations on $G_1^{\vee}(\mathbf{A}_F)$. The same is thus true for the linear form $b_+ f_1^{L_1^{\vee}}(\psi_1 \times \Lambda^{\vee})$.

From this the same argument in the proof of proposition 6.4.4 can be applied, which also gives the vanishing of the linear form $\psi_1 \times \Lambda^{\vee}$. Since the linear form defined by $\psi_1$ does not vanish, it follows that $\Lambda^{\vee}$ and hence $\Lambda$ vanishes.

\end{proof}

\begin{proposition}
In the situation of proposition 6.3.3, we have for $G^* \in \widetilde{\mathcal{E}}_{\simp}(N_+)$
\begin{eqnarray}
\tr R^{G^*}_{\disc,\psi_+^{N_+}}(f^*) =0 = \leftexp{0}{ S^{G^*}_{\disc,\psi_+^{N_+}}(f^*) }, \,\ f^* \in \mathcal{H}(G^*).
\end{eqnarray}
Furthermore $\Lambda$ vanishes, i.e. the stable multiplicity formula is valid for $\psi^N$, and the terms occuring in (6.3.19):
\begin{eqnarray}
& &  \frac{1}{8} \big( (f^{\vee})^{M_+^{\vee}}(\psi^{\vee}_{M_+^{\vee}}) - \delta_{\psi} f^{\vee}_{G^{\vee}_+}(\psi^{\vee}_+,x^{\vee}_{+,1})      \big)                  \\
& &  \frac{1}{8} \big( f^{\prime}_{G_+}(\psi_+,s_{\psi_+} x_{+,1}) - \delta_{\psi} f_{G_+}(\psi_+,x_{+,1})      \big)
\end{eqnarray}
vanishes identically, and we have
\begin{eqnarray}
\delta_{\psi}=1.
\end{eqnarray}
\end{proposition}
\begin{proof}
We use similar argument as in the proofs of propositions 6.4.4 - 6.4.6. By descent as in the proof of proposition 6.4.5, the term (6.4.25) can be written as:
\begin{eqnarray}
\frac{1}{8} \tr \big(  (1 - \delta_{\psi} R_{P^{\vee}_+}(w_+^{\vee}, \widetilde{\pi}_{\psi^{\vee}_{M_+^{\vee}} } ,\psi^{\vee}_{M_+^{\vee}} )  ) \mathcal{I}_{P_+^{\vee}} (\pi_{\psi^{\vee}_{M_+^{\vee}} }  ,f^{\vee}) \big).
\end{eqnarray}

\noindent Again $w_+^{\vee}$ is an element of order two in $W_{\psi_+^{\vee}} = W_{\psi_+^{\vee}}^0$, and hence the eigenvalues of the intertwining operator $R_{P^{\vee}_+}(w_+^{\vee}, \widetilde{\pi}_{\psi^{\vee}_{M_+^{\vee}} } ,\psi^{\vee}_{M_+^{\vee}} )$ are $\pm 1$. Since $\delta_{\psi}$ is also a sign, it again follows that (6.4.28) is a linear combination with positive coefficients of irreducible characters of representations of $G^{\vee}_+$. 

Similar to the proof of proposition 6.4.6 (and with similar notations), the term $f^{L_+}(\Gamma \times \Lambda)$ occuring on the left hand side of (6.3.19) is equal to $f_1^{L_1^{\vee}}(\Gamma \times \Lambda^{\vee})$. It is a unitary character on $G_1^{\vee} = G \times G^{\vee} \in \widetilde{\mathcal{E}}_{\ellip}(N_+)$. Indeed $\Lambda^{\vee}$ is a unitary character on $G^{\vee}$ by proposition 6.1.4, while for $\Gamma$ we have:
\begin{eqnarray*}
h^G(\Gamma)&=&S^G_{\disc,\psi^N}(h)\\
&=& R^G_{\disc,\psi^N}(h), \,\ h \in \mathcal{H}(G).
\end{eqnarray*}
The second equality holds since $\psi^N$ is a simple parameter. Thus $\Gamma$ is a unitary character on $G$ also. 

Thus the left hand side of (6.3.19) is a linear combination with non-negative coefficients of irreducible characters of representations of $G^*$, for $G^* \in \widetilde{\mathcal{E}}_{\simp}(N_+)$. 

\noindent The right hand side of (6.3.19), namely the term (6.4.26), can be treated in exactly the same way as in proposition 6.4.4. Thus the same reasoning used in the proof of propositions 6.4.4 - 6.4.6 gives (6.4.24), the vanishing of the linear form $\Gamma \times \Lambda^{\vee}$, and the vanishing of (6.4.25) and (6.4.26). Furthermore, we also obtain the following:
\begin{eqnarray}
 \delta_{\psi} R_{P^{\vee}_+}(w_+^{\vee}, \widetilde{\pi}_{\psi^{\vee}_{M_+^{\vee}} } ,\psi^{\vee}_{M_+^{\vee}} ) \equiv 1
\end{eqnarray}
In particular the intertwining operator $R_{P^{\vee}_+}(w_+^{\vee}, \widetilde{\pi}_{\psi^{\vee}_{M_+^{\vee}} } ,\psi^{\vee}_{M_+^{\vee}} )$ is a scalar that is equal to $\pm 1$.

\bigskip

\noindent Finally we show (6.4.27); equivalently by (6.4.29) we need to show that $R_{P^{\vee}_+}(w_+^{\vee}, \widetilde{\pi}_{\psi^{\vee}_{M_+^{\vee}} } ,\psi^{\vee}_{M_+^{\vee}} )  \equiv 1$. We use the result from Whittaker models in section 3.5. Indeed since $\psi^{\vee}_{M_+^{\vee}}$ is a generic parameter, and $M^{\vee}_+ \cong G_{E/F}(N)$, we can apply proposition 3.5.3, and deduce that for each place $v$ of $F$, the local intertwining operator $ R_{P^{\vee}_{+,v}}(w_{+,v}^{\vee}, \widetilde{\pi}_{\psi^{\vee}_{M_+^{\vee},v} } ,\psi^{\vee}_{M_+^{\vee},v} )$ acts trivially on the generic subspace of $\mathcal{I}_{P_{+,v}^{\vee}}(\pi_{\psi^{\vee}_{M_+^{\vee},v}})$. Hence we see that the global intertwining operator    
\[
R_{P^{\vee}_+}(w_+^{\vee}, \widetilde{\pi}_{\psi^{\vee}_{M_+^{\vee}} } ,\psi^{\vee}_{M_+^{\vee}} ) = \bigotimes_v R_{P^{\vee}_{+,v}}(w_{+,v}^{\vee}, \widetilde{\pi}_{\psi^{\vee}_{M_+^{\vee},v} } ,\psi^{\vee}_{M_+^{\vee},v } )
\]
also acts trivially on the generic subspace of $\mathcal{I}_{P_+^{\vee}}(\pi_{\psi^{\vee}_{M_+^{\vee}}})$. Since we already know that it is a scalar, it must be identifically equal to one, as required.

Finally it remains to see that $\Lambda$ vanishes. From the vanishing of $\Gamma \times \Lambda^{\vee}$ it follows that either $\Gamma$ or $\Lambda^{\vee}$ vanishes. If it is $\Gamma$ that vanishes, it would follow that
\[
f^G(\psi) = f^G(\Gamma) + f^G(\Lambda) =f^G(\Lambda) = f^L(\Lambda), \,\ f \in \mathcal{H}(G)
\]
thus the linear form $f^G(\psi)$ is induced from $L(\mathbf{A}_F)$. This contradicts in particular the local condition 6.4.1(a). Hence it is $\Lambda^{\vee}$ and hence $\Lambda$ that vanishes.
\end{proof}

\bigskip

To conclude this, we see that for the specific class of families of parameters $\widetilde{\mathcal{F}}$ studied in this subsection (i.e. satisfying assumption 6.4.1), we have the following:

\begin{proposition}
Let $\psi^N \in \widetilde{\mathcal{F}}_{\simp}(N)$ be a simple generic parameter, with $\widetilde{F}$ being a family satisfying Assumption 6.4.1. Let $G=(G,\xi) \in \widetilde{\mathcal{E}}_{\simp}(N)$, with $\xi =\xi_{\chi_{\kappa}}$ for $\chi_{\kappa} \in \mathcal{Z}_E^{\kappa}$. Then the following are equivalent:
\bigskip

\noindent (a) The distribution $S^{G}_{\disc,\psi^N}$ is not identifically zero. 

\bigskip

\noindent (b) The distribution $\tr R^G_{\disc,\psi^N}$ is not identifically zero.

\bigskip

\noindent (c) The Asai $L$-function $L(s,\psi^N,\Asai^{(-1)^{N-1} \kappa}  )$ has a pole at $s=1$.

\end{proposition}
\begin{proof}
Since $\psi^N \in \widetilde{\mathcal{F}}(N)$ is simple generic, we have the identity:
\begin{eqnarray}
\tr R_{\disc,\psi^N}^G(f) = I^G_{\disc,\psi^N}(f)  = S^G_{\disc,\psi^N}(f), \,\ f \in \mathcal{H}(G)
\end{eqnarray}
(here we have used the induction hypothesis to deduce the vanishing of the right hand sides of (5.4.2) and (5.6.2) respectively for the distribution $I^G_{\disc,\psi^N}$).

\noindent We then deduce the equivalence of (a) and (b) from (6.4.30). On the other hand, since we have proved that $\widetilde{\mathcal{F}}_{\gsimp}(G) = \widetilde{\mathcal{F}}_{\gsimp}^{\#}(G)$ in lemma 6.4.2  (see part (b) of remark 6.1.2), if follows that the condition that $S^G_{\disc,\psi^N}$ is not identifically zero is equivalent to the condition that $\psi^N \in \xi_* \widetilde{\mathcal{F}}_{\gsimp}(G)$. Then by proposition 6.4.7, the condition that $\psi^N \in \xi_* \widetilde{\mathcal{F}}_{\gsimp}(G)$ implies $\delta_{\psi}=1$, i.e. part (c) is satisfied. Finally, since we know that eactly one of the two Asai $L$-functions $L(s,\psi^N,\Asai^{\pm})$ has a pole at $s=1$, we see that (a) and (c) are equivalent.

\end{proof}

\begin{rem}
\end{rem}

In particular we see from the proof of Proposition 6.4.8 that the definition of the set $\widetilde{\mathcal{F}}(G)_{\gsimp}$ of simple generic parameters is equivalent to the original definition given in section 2.4, and that the seed theorem 2.4.2, and theorem 2.5.3(a) are satisfied for simple generic parameters in $\widetilde{\mathcal{F}}(N)$. We also note that as a consequence, equation (6.1.4) for the families $\widetilde{\mathcal{F}}$ considered in this subsection (which we now know to be equality) is indeed a disjoint union.

\bigskip

The fact that the equivalences listed in proposition 6.4.8 should be valid for all simple generic parameters $\psi^N \in \widetilde{\Psi}_{\simp}(N)$ is of course a main part of the classification theorem. Although proposition 6.4.8 only applies to the class of parameters that satisfy rather stringent local conditions at the archimedean places (namely assumption 6.4.1), the results will be used in section 7 to establish the local classification theorems, which in turn would be used in the induction argument to establish the global results in section 9.

\section{\bf Local Classification}

In section seven we construct the packets associated to generic local parameters, and obtain the local Langlands classification for tempered representations of quasi-split unitary groups. Among the technical results to be established is the local intertwining relation, which reduces the constructions of packets of tempered representations to the case of discrete series representations. The method of proof is global, drawing on the results from trace formulas comparison in section 5 and 6. This is based on the standard technique of embedding a discrete series representation of a local group as a local component of an automorphic representation.

We need to remark about the use of induction in section 7. In this section we are going to establish the local theorems for generic parameters for {\it all} degrees $N$, by induction. Thus throughout the induction arguments of section 7, we fix the integer $N$ and we assume that all the local theorems for generic parameters hold for parameters of degree smaller than $N$. Besides, we also need to use global inputs from section 6.4; more precisely we will work with families $\dot{\widetilde{\mathcal{F}}}$ of global generic parameters to be introduced in section 7.3 (these are in particular global generic parameters with local constraints at archimedean places as in section 6.4); all the global arguments are to be carried out within families $\dot{\widetilde{\mathcal{F}}}$. Hence for the purpose of carrying out the induction arguments we also need to establish the global theorems for the families of global generic parameters $\dot{\widetilde{\mathcal{F}}}$ for {\it all} degrees in this section. Thus in the induction arguments of this section, we assume that all the global theorems are valid for the families of global generic parameters $\dot{\widetilde{\mathcal{F}}}$ of degree smaller than $N$.

\subsection{Resum\'e on local parameters and local packets}

In the entire section seven, we always denote by $F$ a local field, and $E$ a quadratic extension of $F$. In accordance with the previous notations, we generally denote by $G=(G,\xi) \in \widetilde{\mathcal{E}}_{\ellip}(N)$ an elliptic twisted endoscopic datum of the twisted group $\widetilde{G}_{E/F}(N)$ over $F$. In this subsection it suffices to consider the case $G \in \widetilde{\mathcal{E}}_{\simp}(N)$, i.e. the underlying endoscopic group $G$ is $U_{E/F}(N)$ (thus implicitly $E$ is a quadratic field extension of $F$, for otherwise if $E$ is the split extension then $G=\GL_{N/F}$ and the local classification results is of course known).

We will only need to consider generic parameters in section 7; the case of non-generic parameters being the subject of section 8. As before $\Phi(G)$ is the set of generic parameters of $G(F)$. The $L$-embedding $\xi : \leftexp{L}{G} \rightarrow \leftexp{L}{G_{E/F}(N)}$, which is part of the endoscopic datum for $G$, allows us to identify a parameter $\phi \in \Phi(G)$ with the $N$-dimensional representation $\phi^N$ associated to $\xi_* \phi$:
\[
\phi^N: L_E \rightarrow \GL_N(\mathbf{C})
\]
in accordance with lemma 2.2.1. Recall that if $\xi =\xi_{\chi}$ for some conjugate self-dual character $\chi \in \mathcal{Z}_E^{\kappa}$, $\kappa = \pm 1$ (notation as in (2.1.3) of {\it loc. cit.}), then $\phi^N$ is conjugate self-dual of parity $(-1)^{N-1} \kappa$.

We first consider the case that $F$ is non-archimedean. With $K=K_F$ being the standard open maximal compact subgroup of $G(F)$. Then $K$ is a special maximal compact subgroup of $G(F)$, and hyperspecial when $E/F$ is unramified. An irreducible admissible representation $\pi$ of $G(F)$ is said to be $K$-spherical if its restriction to $K$ contains the trivial representation. In particular $\pi$ is unramified when $E/F$ is unramified. A packet $\Pi_{\phi}$ is spherical if it contains a (unique) spherical representation. In this case the parameter $\phi$ can be chosen to factor through $\leftexp{L}{M_0}$, where $M_0$ is the minimal Levi subgroup of $G=U_{E/F}(N)$ given by the standard diagonal maximal torus. The resulting parameter $\phi_{M_0}:L_F \rightarrow \leftexp{L}{M_0}$, which factors through $W_F$ and corresponds to a spherical character $\pi_{M_0}$ of $M_0(F)$, under the Langlands correspondence for the torus $M_0$. If $\phi \in \Phi_{\bdd}(G)$, then $\pi_{M_0}$ is unitary, and the packet $\Pi_{\phi}$ is just the irreducible constituents of $\mathcal{I}_{P_0}(\pi_{M_0})$ ($P_0$ being the standard Borel subgroup of $G$). 

In addition when $E/F$ is unramified, then $\phi$ is an unramified parameter, and the packet $\Pi_{\phi}$ contains a unique unramified representation. More concretely, with $\xi =\xi_{\chi}$ for a character $\chi \in Z_E^{\kappa}$ as above, the representation $\phi^N = \xi_{\chi,*} \phi$ of $L_E$, which factors through $W_E$, decomposes as a sum of $N$ characters of $W_E$:
\begin{eqnarray}
\phi^N = \eta_1 \oplus \cdots \oplus \eta_N
\end{eqnarray} 
with the characters $\chi^{-1} \eta_i$ being unramified. Note that when $E/F$ is unramified, then we can always choose $\chi \in \mathcal{Z}_E^{\kappa}$ to be unramified.

It is immediate to check that if $E/F$ is unramified, then an unramified character $\eta$ of $W_E$ is conjugate self-dual (with respect to the extension $E/F$) if and only if it is self-dual, i.e. $\eta$ is a quadratic unramified character, and that $\eta$ is conjugate orthogonal (resp. conjugate symplectic) if and only if $\eta(\Frob_E)=+1$, i.e. $\eta$ is the trivial character (resp. $\eta(\Frob_E)=-1$, i.e. $\eta$ is the unique non-trivial unramified quadratic character). Recall that the representation $\phi^N$ has to be conjugate self-dual of parity $(-1)^{N-1} \kappa$. Thus the characters $\eta_i$ in (2.1.1) has to be subjected to condition as in (2.4.12-13). From this description it is also immediate to check that in the unramified case $\mathcal{S}_{\phi}$ is at most of order two (it can be of order two only when $N$ is even).

\bigskip

We next consider the case where $F$ is archimedean, i.e. $F \cong \mathbf{R}$ and $E \cong \mathbf{C}$. The archimedean case plays an important role in the global-local methods of section seven. We have $L_{\mathbf{C}}=W_{\mathbf{C}}=\mathbf{C}^{\times}$. The most important case being $\phi \in \Phi_2(G)$. In this case $\phi^N = \xi_* \phi$ (with $\xi=\xi_{\chi}$ for $\chi \in \mathcal{Z}_E^{\kappa}$) takes the form similar to (2.1.1):
\begin{eqnarray}
\phi^N = \eta_1 \oplus \cdots \oplus \eta_N
\end{eqnarray}
where $\eta_i$ are {\it distinct} conjugate self-dual characters of $\mathbf{C}^{\times}$ (with respect to $\mathbf{C}/\mathbf{R}$) of parity equal to $(-1)^{N-1} \kappa$. In general, a conjugate self-dual character $\eta$ of $\mathbf{C}^{\times}$ is of the form:
\begin{eqnarray}
\eta: z \mapsto (z/\overline{z})^{a}
\end{eqnarray} 
with $a \in \frac{1}{2}\mathbf{Z}$. The character $\eta$ is conjugate orthogonal (resp. conjugate symplectic) if and only if $a \in \mathbf{Z}$ (resp. $a \in \frac{1}{2} \mathbf{Z} - \mathbf{Z}$).

In (2.1.2), if $\eta_i(z)=(z/\overline{z})^{a_i}$ for $a_i \in \frac{1}{2} \mathbf{Z}$, then we identify the $N$-tuple:
\[
\mu_{\phi^N}:=(a_1,\cdots ,a_N) 
\]  
as the infinitesimal character of $\phi^N$. The infinitesimal character of $\phi$ is given by the shift:
\[
\mu_{\phi}=(b_1,\cdots,b_N)
\]
where $a_i=b_i+c$, for $c \in \frac{1}{2}\mathbf{Z}$ is such that $\chi(z)=(z/\overline{z})^c$. In particular for $\phi \in \Phi_2(G)$ we have $\mu_{\phi} \in (\mathbf{Z} + \frac{1}{2}(N-1))^N$.

We also define:
\begin{eqnarray}
d(\mu_{\phi}) = \inf( \min_i(b_i), \min_{i \neq j} (|b_i-b_j|)).
\end{eqnarray}
and similarly for $d(\mu_{\phi^N})$.

The description of general parameters $\Phi(G)$ is as in (2.4.12-13).

The archimedean case of the local theorems, stated as theorem 2.5.1, and which is stated in the more precise form as theorem 3.2.1, are already known in the case of (bounded) generic parameters, thanks to the works on Mezo and Shelstad. More precisely, part (b) of theorem 3.2.1 follows from the works of Shelstad \cite{Sh1,Sh2} on spectral transfer in standard endoscopy for real groups, while part (b) follows from the works of Mezo \cite{Me} and Shelstad \cite{Sh3} on spectral transfer in twisted endoscopy for real groups. Thus we have:
\begin{theorem}
Theorem 3.2.1 (and hence theorem 2.5.1) holds for bounded generic parameters in the archimedean case.
\end{theorem} 

The remaining local result for bounded generic parameters in the archimedean case to be established is the local intertwining relation, stated as theorem 3.4.3. We will need to use the following weaker version of the local intertwining relation, which follows from Shelstad's results \cite{Sh1,Sh2} ({\it c.f.} the discussion in section 6.1 of \cite{A1} concerning the results \cite{Sh1,Sh2}). First for $\phi \in \Phi_{\bdd}(G)$ and  $x \in \mathcal{S}_{\phi}$, we have the character identity from \cite{Sh1,Sh2}:
\[
f^{\prime}(\phi,x) = f^{\prime}(\phi^{\prime}) = \sum_{\pi \in \Pi_{\phi}} \langle x,\pi \rangle f_G(\pi), \,\ f \in \mathcal{H}(G).
\] 
Let $M$ be a Levi subgroup of $G$ such that $\Phi_2(M,\phi)$ (the set of square-integrable parameters of $M$ that maps to $\phi$) is non-empty, say $\phi_M \in \Phi_2(M,\phi)$. We then have the packet $\Pi_{\phi_M}$, and the packet $\Pi_{\phi}$ are constructed as irreducible constituents of $\mathcal{I}_P(\pi_M)$, with $\pi_M$ ranges over elements of $\Pi_{\phi_M}$. There is then a well-defined projection map $\Pi_{\phi} \rightarrow \Pi_{\phi_M}$ sending $\pi$ to $\pi_M$. 

We also have the $R$-group $R_{\phi}$ which is a quotient of $\mathcal{S}_{\phi}$ whose kernel is given by $\mathcal{S}_{\phi_M}$ (see the discussion in section 3.4). The results of Shelstad \cite{Sh1,Sh2} already gives the isomorphism between $R_{\phi}$ and the representation theoretic $R$-group $R(\pi_M)$ associated to any $\pi_M \in \Pi_{\phi_M}$. Then the weaker version of the local intertwining relation obtained from Shelstad's results that we will use is as follows:
\begin{proposition}
For every $\pi_M \in \Pi_{\phi_M}$, there exists a character $\epsilon_{\pi_M}$ on $\mathcal{S}_{\phi}$, which is the pull-back of a character on $R_{\phi}$, such that for every $x \in \mathcal{S}_{\phi}$, we have:
\begin{eqnarray}
f_G(\phi,x) = \sum_{\pi \in \Pi_{\phi}} \epsilon_{\pi_M}(x) \langle x,\pi \rangle f_G(\pi), \,\ f \in \mathcal{H}(G),
\end{eqnarray}
where in the sum occuring in (2.1.5), we have denoted, for $\pi \in \Pi_{\phi}$, its projection to $\Pi_{\phi_M}$ as $\pi_M$. 
\end{proposition}

Thus to establish the local intertwining relation for the case of archimedean generic parameters, we have to show that the characters $\epsilon_{\pi_M}$ are all trivial. We need the following result that is a consequence of Shelstad's results:
\begin{proposition} (\cite{Sh2}, theorem 11.5)
Let $(B,\omega)$ be Whittaker data for $G$, and $(B_M,\omega_M)$ be the corresponding Whittaker data for $M$ (recall that we normalize transfer factors according to Whittaker data, as in section 3.2). If $\pi$ (resp. $\pi_M$) is the unique generic representation in $\Pi_{\phi}$ (resp. $\Pi_{\phi_M}$), then we have $\langle \cdot, \rangle =1$ (resp. $\langle \cdot, \pi_M \rangle) =1$. Consequently the character $\epsilon_{\pi_M}$ is trivial for the generic $\pi_M$.
\end{proposition}

Finally we need the following result concerning archimedean packet for generic parameters. It is proved as lemma 6.1.2 of \cite{A1}. The proof in {\it loc. cit.} applies to any connected reductive groups over $\mathbf{R}$. 
\begin{proposition} (Lemma 6.1.2 of \cite{A1})
For $\phi \in \Phi_{\bdd}(G)$, the image of $\Pi_{\phi}$ in $\widehat{\mathcal{S}}_{\phi}$ generates $\widehat{\mathcal{S}}_{\phi}$ as a group.
\end{proposition}

\begin{rem}
\end{rem}
\noindent In the case where $F$ is non-archimedean and $\phi$ is a square-integrable parameter, the packet $\Pi_{\phi}$ associated to $\phi$ was constructed by Moeglin \cite{Moe} by a different method. The construction of Moeglin gives information about supercuspidal representations which are not available from the general endoscopic classification {\it apriori}.

\subsection{Construction of global representation}

As in \cite{A1}, the method of proof of the local theorems relies on global techniques. In order to globalize the local data, we need a series of preparations as follows. We follow the convention of Arthur \cite{A1} and put a $``\cdot"$ over a local object to denote a choice of globalization of the local object. 

\begin{lemma}
Let $F$ be either real or a $p$-adic field. For any positive integer $r_0$, there exists a totally real field $\dot{F}$ having at least $r_0$ real places, and a place $u$ of $\dot{F}$ such that $\dot{F}_u=F$. Furthermore if $E/F$ is a quadratic extension, then we can choose $\dot{F}$ and a totally imaginary quadratic extension $\dot{E}$ of $\dot{F}$, such that $u$ does not split in $\dot{E}$, and such that $\dot{E}_u=E$, and such that $\dot{E}$ is unramified over $\dot{F}$ at all finite places of $\dot{F}$ outside $u$ (In particular if $F=\mathbf{R}$ then $\dot{E}$ can be chosen to be unramified over $\dot{F}$ at all finite places).
\end{lemma}
This is well-known; for instance the first assertion is proved in \cite{A1}, lemma 6.2.1. The second assertion concerning quadratic extension is standard.

\begin{lemma}
Let $E/F$ and $\dot{E}/\dot{F}$ be as in lemma 7.2.1. Then for any $\chi \in \mathcal{Z}_E^{\kappa}$, there exists $\dot{\chi} \in \mathcal{Z}_{\dot{E}}^{\kappa}$ such that $\dot{\chi}_u =\chi$. 
\end{lemma}
This is again well-known and is in any case a simple exercise in Gr\"o{\ss}encharakter.

From lemma 7.2.1 and 7.2.2 it follows that if $G \in \widetilde{\mathcal{E}}(N)$, then there exists $\dot{G} \in \dot{\widetilde{\mathcal{E}}}(N)$ such that $\dot{G}_u=G$. Similar results of course hold for $\widetilde{\mathcal{E}}_{\ellip}(N)$ and $\widetilde{\mathcal{E}}_{\simp}(N)$.

Thus let $\dot{E}/\dot{F}$ be as above associated to $E/F$. We denote by $S_{\infty}$ the set of archimedean places of $\dot{F}$, and put
\[
S_{\infty}(u) :=S_{\infty}(u), \,\ S_{\infty}^u :=S_{\infty} - \{u\}.
\]
We may assume that the number of real places of $\dot{F}$ is large (for example for lemma 7.2.3 below, it suffices to have $|S_{\infty}(u)|\geq 2$).

\begin{lemma}
Given $G=(G,\xi) \in \widetilde{\mathcal{E}}_{\simp}(N)$, and a square-integrable representation $\pi \in \Pi_2(G)$ of $G(F)$, and $t \in \mathbf{N}$, there exists $\dot{G} = (\dot{G},\dot{\xi}) \in \dot{\widetilde{\mathcal{E}}}_{\simp}(N)$, and an automorphic representation $\dot{\pi}$ occuring in $L^2_{\disc,G(F) \backslash G(\mathbf{A})}$, with the following properties:
\begin{enumerate}
\item $(\dot{G}_u,\dot{\xi}_u)=(G,\xi)$ and $\dot{\pi}_u=\pi$.

\item For every place $v \notin S_{\infty}(u)$, the representation $\dot{\pi}_v$ is spherical at $v$.

\item For any $v \in S_{\infty}^u$, $\dot{\pi}_v$ is a square-integrable representation of $\dot{G}(\dot{F}_v)$, whose Langlands parameter $\phi_v \in \Phi_2(\dot{G}_v)$ satisfies $d(\mu_{\phi_v})>t$ (i.e. the parameter $\phi_v$ can be chosen to be in general position).
\end{enumerate}
\end{lemma}
\begin{proof}
Given lemma 7.2.1 and 7.2.2, the proof of lemma 7.2.3 is then exactly the same as that of lemma 6.2.2 of \cite{A1}, which is based on the simple version of the invariant trace formula. So we just give a sketch. We point out that the proof depends only on the results of section 4, and is thus independent of any induction hypothesis. 

Let 
\[
\dot{K}^{\infty,u} = \prod_{v \notin S_{\infty}(u)} \dot{K}_v
\] 
be the standard maximal open compact subgroup of $\dot{G}(\dot{\mathbf{A}}^{\infty,u})$. Put $\dot{f}^{\infty,u}$ to be equal to the characteristic function of $\dot{K}^{\infty,u}$. At the place $u$, we take $\dot{f}_u \in \mathcal{H}(G)$ to be a pseudo-coefficient $f_{\pi}$ of the representation $\pi$ (existence of $f_{\pi}$ is given by \cite{CD} and \cite{BDK}). At the places $S_{\infty}^u$, we arbitrarily choose Langlands parameters $\phi_v \in \Phi_2(\dot{G}_v)$ that satisfies $d(\mu_{\phi_v}) > t$, and we take $\dot{f}_v$ to be a stable pseudo-coefficient associated to $\phi_v$, i.e.
\[
\dot{f}_v = \sum_{\pi_v \in \Pi_{\phi_v}} f_{\pi_v}.
\]
Then put $\dot{f}(\dot{x}) := \dot{f}_u(\dot{x}_u) \dot{f}^u_{\infty}(\dot{x}^u_{\infty}) \dot{f}^{\infty,u}(\dot{x}^{\infty,u})$. The invariant trace formula of Arthur \cite{A3}, when applied to $\dot{f}$, simplifies and gives the identity (the simplification follows from \cite{A3} theorem 7.1, \cite{A10} theorem 5.1, and p.268 of \cite{A10}; {\it c.f.} the discussion in the proof of lemma 6.2.2 of \cite{A1}):

\begin{eqnarray}
\sum_{\dot{\pi}} m_{\disc}(\dot{\pi}) \dot{f}_{\dot{G}}(\dot{\pi}) = \sum_{\gamma} \vol(\dot{G}_{\gamma}(\dot{F})\backslash \dot{G}_{\gamma}(\dot{\mathbf{A}})) \dot{f}_{\dot{G}}(\gamma)
\end{eqnarray}
here on the left hand side $\dot{\pi}$ runs over the irreducible unitary representation of $\dot{G}(\dot{\mathbf{A}})$, and $m_{\disc}(\dot{\pi})$ is its multiplicity in $L^2_{\disc}(\dot{G}(\dot{F}) \backslash \dot{G}(\dot{\mathbf{A}})$. On the right hand side is a finite sum of nonzero terms where $\gamma$ runs over the semi-simple conjugacy classes in $\dot{G}(\dot{F})$ that are $\mathbf{R}$-elliptic at each place of $S_{\infty}$. 

Put $\dot{Z}^{\infty,u}=Z_{\dot{G}}(\dot{F}) \cap \dot{K}^{\infty,u}$ (here $Z_{\dot{G}}$ is the center of $\dot{G}$, which is given by $U_{\dot{E}/\dot{F}}(1)$). Note that $Z_{\dot{G}}(\dot{F}_v)$ is compact for each $v \in S_{\infty}(u)$ (since $v$ does not split in $\dot{E}$ for $v \in S_{\infty}(u)$). Hence $\dot{Z}^{\infty,u}$ is a finite cyclic group. We can then choose $\phi_v$ for $v \in S_{\infty}^u$ such that the product $\dot{f}_u \dot{f}_{\infty}^u$ is invariant under translation by $\dot{Z}^{\infty,u}$. Then for $\gamma \in Z_{\dot{G}}(\dot{F})$, we have $\dot{f}_{\dot{G}}(\gamma) \neq 0$ only for $\gamma \in \dot{Z}^{\infty,u}$, in which case we have
\begin{eqnarray}
\dot{f}_{\dot{G}}(\gamma) = \dot{f}_{\dot{G}}(1).
\end{eqnarray}
From the results of Harish-Chandra (lemma 23 of \cite{HC1}, lemma 17.4 and 17.5 of \cite{HC2}), the terms on the right hand side of (7.2.1) is dominated by the terms with $\gamma \in Z_{\dot{G}}(\dot{F})$, when the infinitesimal characters $\mu_{\phi_v}$ for $v \in S_{\infty}^u$ are in sufficient general position. Hence from (7.2.2), we can choose the $\phi_v$ so that the right hand side of (7.2.1) is non-zero.

Thus from the non-vanishing the of the left hand side of (7.2.1), there exists a $\dot{\pi}$ such that
\begin{eqnarray}
m_{\disc}(\dot{\pi}) \dot{f}_{\dot{G}}(\dot{\pi}) = m_{\disc}(\dot{\pi}) f_{\pi}(\dot{\pi}_u) \dot{f}^u_{\infty}(\dot{\pi}^u_{\infty}) \dot{f}^{\infty,u}(\dot{\pi}^{\infty,u})   \neq 0
\end{eqnarray}

From (7.2.3) it follows that $\dot{\pi}_v$ is spherical for $v \notin S_{\infty}(u)$. For $v \in S_{\infty}^u$, we can deduce that $\dot{\pi}_v$ has to be tempered, from the fact that $\dot{\pi}_v$ is unitary and its infinitesimal character is in general position. Thus since $\dot{f}_v$ is a stable pseudo-coefficient of $\Pi_{\phi_v}$ it follows that $\dot{\pi}_v \in \Pi_{\phi_v}$. For the remaining place $u$, since $f_{\pi}$ is the pseudo-coefficient of $\pi$, it suffices to show that $\dot{\pi}_u$ is tempered.

First from corollary 4.3.8, there exists a unique $\dot{\psi}^N \in \widetilde{\Psi}_{\ellip}(N)$, such that $\dot{\pi}$ belongs to
\[
L^2_{\disc,\dot{\psi}^N}(\dot{G}(\dot{F})\backslash \dot{G}(\dot{\mathbf{A}}))
\]
and as in section 4 we have the stabilization of the twisted trace formula for $\dot{\widetilde{G}}(N)$:
\begin{eqnarray}
\widetilde{I}^N_{\disc,\dot{\psi}^N}(\dot{\widetilde{f}}) = \sum_{\dot{G}^* \in \dot{\widetilde{\mathcal{E}}}_{\ellip}(N)} \widetilde{\iota}(N,\dot{G}^*) \widehat{S}^{\dot{G}^*}_{\disc,\dot{\psi}^N}(\dot{\widetilde{f}}^{\dot{G}^*}), \,\ \dot{\widetilde{f}} \in \dot{\widetilde{\mathcal{H}}}(N).
\end{eqnarray}

We choose decomposable functions:
\[
\dot{\widetilde{f}} = \dot{\widetilde{f}}_u \cdot \dot{\widetilde{f}}^u_{\infty} \cdot \dot{\widetilde{f}}^{\infty,u} \in \dot{\widetilde{\mathcal{H}}}(N)
\]
\[
\dot{f}= \dot{f}_u \cdot \dot{f}^u_{\infty} \cdot \dot{f}^{\infty,u} \in \dot{\mathcal{H}}(\dot{G})
\]
such that for each $v$, we have $\dot{\widetilde{f}}^{\dot{G}_v}_v=\dot{f}_v^{\dot{G}_v}$. For the place at $u$ we allow $\dot{f}_u$ to be a variable component, while for $v \neq u$, we let $\dot{f}_v$ be chosen as before. In addition, for $v \in S_{\infty}^u$ we require that $\dot{\widetilde{f}}_v^{G^*_v} =0$ for any $\dot{G}^* \in \dot{\widetilde{\mathcal{E}}}_{\ellip}(N)$ other than $\dot{G}$ itself (these choices are possible by proposition 3.1.1). Then (7.2.4) simplifies to 
\begin{eqnarray}
\widetilde{I}^N_{\disc,\dot{\psi}^N}(\dot{\widetilde{f}}) = \widetilde{\iota}(N,\dot{G}) S^{\dot{G}}_{\disc,\dot{\psi}^N}(\dot{f}).
\end{eqnarray}
On the other hand, under our choice of $\dot{f}$, we have the equality:
\[
\tr R^{\dot{G}}_{\disc,\dot{\psi}^N}(\dot{f})=I^{\dot{G}}_{\disc,\dot{\psi}^N}(\dot{f}) = S^{\dot{G}}_{\disc,\dot{\psi}^N}(\dot{f}).
\]
Thus we have
\begin{eqnarray}
\widetilde{I}^N_{\disc,\dot{\psi}^N}(\dot{\widetilde{f}}) = \iota(N,\dot{G}) \tr R^{\dot{G}}_{\disc,\dot{\psi}^N}(\dot{f}).
\end{eqnarray}
As a linear form of $\dot{f}_u$, the right hand side of (7.2.6) is non-zero (for example when evaluated at $\dot{f}_u = f_{\pi}$, by (7.2.3)). Hence the left hand side of (7.2.6) is non-zero as a linear form of $\dot{\widetilde{f}}_u$. Since the left hand side of (7.2.6) is a non-zero multiple of the linear form $\dot{\widetilde{f}}_N(\dot{\psi}^N)$, we see in particular that $\dot{\psi}^N_v$ is spherical for $v \notin S_{\infty}(u)$, and for $v \in S_{\infty}^u$, we have $\dot{\widetilde{f}}_v(\dot{\psi}^N_v) \neq 0$. 

Recall that for $v \in S^u_{\infty}$ we have $\dot{\widetilde{f}}_v^{\dot{G}_v}=\dot{f}_v^{G_v}$, and $\dot{\widetilde{f}}_v^{\dot{G}^*_v}=0$ for $\dot{G}^* \neq \dot{G}$. Hence given the non-vanishing of $\dot{\widetilde{f}}_v(\dot{\psi}^N_v) $, it follows from the spectral transfer results of \cite{Me} and \cite{Sh3} that $\phi_{\dot{\psi}_v} \in \dot{\xi}_{v,*} \Phi(\dot{G}_v)$, and the infinitesimal characters of $\phi_v$ and that of $\phi_{\dot{\psi}^N_v}$ corresponds (under the $L$-embedding $\dot{\xi}_v$). The infinitesimal character of $\phi_{\dot{\psi}^N_v}$, cannot be in general position if $\dot{\psi}^N_v$ is a non-generic parameter. It thus follows that $\dot{\psi}^N$ itself has to be a generic parameter $\dot{\phi}^N$. 

One first show that $\dot{\phi}^N_u \in \widetilde{\Phi}_{\bdd}(N)$. We have seen above that
\begin{eqnarray}
c^u \dot{\widetilde{f}}_{u,N}(\dot{\phi}^N_u)  = \tr R^{\dot{G}}_{\disc,\dot{\phi}^N}(\dot{f})
\end{eqnarray}
for a non-zero scalar $c^u$, and that (7.2.7) is non-zero when $\dot{f}_u=f_{\pi}$. Since $\dot{f}_{u,\dot{G}_u}$ is cuspidal, $\dot{\widetilde{f}}_u$ can be chosen so that $\dot{\widetilde{f}}_{u,\dot{G}_u}$ is cuspidal. Then from the non-vanishing of $\dot{\widetilde{f}}_u(\dot{\phi}_u)$ it follows that $\dot{\phi}_u \in \widetilde{\Phi}_{\ellip}(N)$ (this follows for example by considering descent to a proper Levi subset of $\widetilde{G}(N)$ if we were to have $\dot{\phi}_u \notin \widetilde{\Phi}_{\ellip}(N)$). In particular $\dot{\phi}_u \in \widetilde{\Phi}_{\bdd}(N)$, and hence $\dot{\phi}_u$ corresponds to a tempered representation of $\GL_N(\dot{E}_u)=\GL_N(E)$, under the local Langlands classification for general linear groups, and hence the linear form $\dot{\widetilde{f}}_{u,N}(\dot{\phi}_u)$ is tempered. As in proof of lemma 6.2.2 of \cite{A1}, one then deduce from this that the linear form $\tr R^{\dot{G}}_{\disc,\dot{\phi}^N}(\dot{f})$ is tempered as a linear form in $\dot{f}_u$. Since $\dot{\pi} = \dot{\pi}_u \otimes \dot{\pi}^u$ occurs in $R^{\dot{G}}_{\disc,\dot{\phi}^N}$ this gives the temperedness of $\dot{\pi}_u$. From this one concludes that $\dot{\pi}_u \cong \pi$.
\end{proof}

\begin{corollary}
For the global generic parameter $\dot{\phi}^N$ constructed in the proof of lemma 7.2.3, we have $\dot{\phi}^N_v$ is spherical for $v \notin S_{\infty}(u)$. For $v \in S_{\infty}^u$, we have $\dot{\phi}_v^N = \dot{\xi}_*\phi_v$ (recall that $\phi_v$ are parameters in general position that was chosen in the beginning of proof).
\end{corollary}
\begin{proof}
We have already seen in the course of the proof of lemma 2.2.3 that $\dot{\phi}^N_v$ is spherical for $v \notin S_{\infty}(u)$. As for $v \in S_{\infty}^u$, the same argument that was applied to the place $u$ in the last part of the proof, applies equally well to $v \in S_{\infty}^u$. In particular this gives $\dot{\phi}_v^N \in \widetilde{\Phi}_{v,\bdd}(N)$. We have seen that $\dot{\widetilde{f}}_{v,N}(\dot{\phi}^N_v) \neq 0$ for any $\dot{\widetilde{f}}_v$ such that $\dot{\widetilde{f}}_v^{\dot{G}_v} = \dot{f}_v^{\dot{G}_v}$ (recall that $\dot{f}_v$ a stable pseudo-coefficient of the parameter $\phi_v$), and that $\dot{\widetilde{f}}_v^{\dot{G}_v^*}=0$ for any $\dot{G}^*_v \neq \dot{G}_v$. Hence it follows from the spectral transfer results of \cite{Me} and \cite{Sh3} that $\dot{\phi}^N_v \in \dot{\xi}_{v,*} \Phi_{\bdd}(\dot{G}_v)$. This in turn implies again by their results that $\dot{\phi}^N_v = \dot{\xi}_{v,*} \phi_v$, as required.
\end{proof}

\begin{rem}
\end{rem}
\noindent Lemma 7.2.3 and corollary 7.2.4 completes the proof of lemma 3.3.2.

\bigskip

\begin{rem}
\end{rem}
\noindent There are obvious variants of lemma 7.2.3. For example if $V$ is a finite set of non-archimedean places of $\dot{F}$ disjoint from $\{u\}$, and $\pi_v$ is a discrete series representation of $\dot{G}(\dot{F}_v)$ for each $v \in V$, then we can find a $\dot{\pi}$ in the discrete automorphic spectrum of $\dot{G}(\dot{\mathbf{A}})$, that in addition to satisfying $\dot{\pi}_u \cong \pi$ at the place $u$, also satisfies $\dot{\pi}_v \cong \pi_v$ for $v \in V$.

We state one more corollary to lemma 7.2.3, which will be needed to globalize local parameters in the next subsection. To emphasize that this is the only result of this subsection that depends on the induction hypothesis, we state this explicitly in the statement:

\begin{corollary} 
Suppose $N_1 \leq N$, and assume that all the local and global theorems are valid for parameters of degree up to $N_1$. Let $G =(G,\xi) \in \widetilde{\mathcal{E}}_{\simp}(N_1)$ over $F$ and $\dot{G} =(\dot{G},\dot{\xi}) \in \dot{\widetilde{\mathcal{E}}}_{\simp}(N_1)$ over $\dot{F}$ be constructed as in lemma 7.2.1 and 7.2.2 (so that $\dot{G}_u=G$ for a place $u$ of $\dot{F}$). Then for any simple local parameter $\phi \in \Phi_{\simp}(G)$, there exists a simple global parameter $\dot{\phi} \in \dot{\Phi}_{\simp}(\dot{G})$, with the following property:
\begin{enumerate}

\item $\dot{\phi}_u = \phi$. 

\item $\dot{\phi}_v$ is a spherical parameter for any $v \notin S_{\infty}(u)$.

\item For any $v \in S_{\infty}^u$, the parameter $\dot{\phi}_v$ belongs to $\Phi_2(\dot{G}_v)$ and is in general position. 
\end{enumerate} 
\end{corollary}
\begin{proof}
Since we are assuming the local theorems are valid for parameters of degree up to $N_1$, the packet $\Pi_{\phi}$ associated to $\phi$ exists. Pick any $\pi \in \Pi_{\phi}$. By lemma 7.2.3, we can globalize $\pi$ to a discrete automorphic representation $\dot{\pi}$ on $\dot{G}(\dot{\mathbf{A}})$ such that $\dot{\pi}_u=\pi$. Furthermore, we have constructed a $\dot{\psi}^{N_1} \in \dot{\widetilde{\Psi}}(N_1)$ such that $\dot{\pi}$ occurs in $R^{\dot{G}}_{\disc,\dot{\psi}^{N_1}}$. We have seen in the course of the proof of lemma 7.2.3 that $\dot{\psi}^{N_1}=\dot{\phi}^{N_1}$ is actually a generic parameter. From the induction hypothesis concerning the global theorems, we have $\dot{\phi}^{N_1} = \xi_* \dot{\phi}$ for $\dot{\phi} \in \Phi_2(\dot{G})$. Furthermore, the global theorem concerning the decomposition of the global discrete spectrum of $\dot{G}$ also gives the result that $\pi=\dot{\pi}_u$ belongs to the packet corresponding to the parameter $\dot{\phi}_u$. Thus by the disjointness of the packets we have $\dot{\phi}_u=\phi$, and similarly that property (2) and (3) holds for $\dot{\phi}_v$ for $v \neq u$. Furthermore, since $\phi^{N_1}=\dot{\phi}^{N_1}_u$ is a simple parameter, the same must be true for $\dot{\phi}^{N_1}$. Thus $\dot{\phi}$ gives what we want.  
\end{proof}
\begin{rem}
\end{rem}
\noindent As we have seen in the proof of lemma 7.2.3, for $v \in S_{\infty}^u$ the parameter $\dot{\phi}_v$ can be chosen to be any preassigned parameter $\phi_v \in \Phi_2(\dot{G}_v)$ that is in general position, subject to the only condition that the product of the central characters of the parameters $\dot{\phi}_{v}$ over $v \in S_{\infty}(u)$ has to be trivial on $\dot{Z}^{\infty,u}$.

\subsection{Construction of global parameter}
As in section 7.2, $F$ is either real or $p$-adic, and $E/F$ quadratic as before. We let $G=(G,\xi) \in \widetilde{\mathcal{E}}_{\ellip}(N)$. The endoscopic datum $G$ need not be simple. In this subsection we fix a local parameter $\phi \in \Phi_{\bdd}(G) - \Phi_{\simp}(G)$ that is not simple. As usual put $\phi^N:=\xi_* \phi \in \widetilde{\Phi}_{\bdd}(N)$. We assume that all the irreducible components of $\phi^N$ are conjugate self-dual. We can then decompose
\begin{eqnarray}
\phi^N = l_1 \phi_1^{N_1} \oplus \cdots \oplus l_r \phi_r^{N_r}
\end{eqnarray}
where $\phi_i^{N_i} \in \widetilde{\Phi}_{\simp}(N_i)$, and $N_i<N$ for all $i$. Recall that we are under the inductive hypothesis that the local (and global) theorems hold for parameters of degree less than $N$. Hence there exists unique endoscopic datum $G_i=(G_i,\xi_i) \in \widetilde{\mathcal{E}}_{\simp}(N_i)$, and parameters $\phi_i \in \Phi_{\simp}(G_i)$, such that $\phi_i^{N_i}=\xi_{i,*} \phi_i$. In the case where $l_i>1$ for some $i$ in (7.3.1), we assume that $G$ itself is simple. This is to ensure that $\phi$ is uniquely determined by $\phi^N$.

In the case where $l_i>1$ for some $i$, i.e. that $\phi \notin \Phi_2(G)$, then there is a proper Levi subgroup $M$ of $G$, unique up to conjugation, such that $\Phi_2(M,\phi)$ is non-empty, i.e. there exists $\phi_M \in \Phi_2(M)$ that maps to $\phi$. We equip $M$ with the $L$-embedding $\xi$ so that $M=(M,\xi) \in \widetilde{\mathcal{E}}(N)$ (which is not elliptic). In the situation (7.3.1), if we put:
\begin{eqnarray}
N_- = \sum_{l_i  \odd} N_i
\end{eqnarray}
then we have
\begin{eqnarray}
M \cong G_- \times G_{E/F}(N_1)^{l_1^{\prime}} \times \cdots \times G_{E/F}(N_r)^{l_r^{\prime}}
\end{eqnarray}
with $l_i^{\prime}$ being the floor of  $l_i/2$, and $G_- = (G_-,\xi_-) \in \widetilde{\mathcal{E}}_{\ellip}(N_-)$ (here $\xi_-:\leftexp{L}{G_-} \hookrightarrow \leftexp{L}{G_{E/F}(N_-)}$ is the restriction of $\xi$ to $\leftexp{L}{G_-}$). We have
\begin{eqnarray}
\xi_{-,*} \phi_- = \bigoplus_{l_i \odd} \phi_i^{N_i}.
\end{eqnarray}
for uniquely determined $\phi_- \in \Phi_2(G_-)$. By the results of section 7.2, we can globalize all these data. More precisley, from lemma 7.2.1-7.2.2, and corollary 7.2.4, we can choose $\dot{E}/\dot{F}$, and a place $u$ of $\dot{F}$ such that $\dot{E}_u/\dot{F}_u=E/F$, and such that the data
\[
(G,\phi,M,\phi_M,G_i,\phi_i)
\]
can be globalized to data:
\[
(\dot{G},\dot{\phi},\dot{M},\dot{\phi}_M,\dot{G}_i,\dot{\phi}_i)
\]
satisfying the anaogue of (7.3.1), (7.3.3) and (7.3.4). For example, if $\dot{\phi}^N = \dot{\xi}_* \dot{\phi}$, then 
\[
\dot{\phi}^N = l_1 \dot{\phi}_1^{N_1} \boxplus \cdots \boxplus l_r \dot{\phi}_r^{N_r}
\]
with $\dot{\phi}^{N_i}_i = \dot{\xi}_{i,*} \dot{\phi}_i$ for $\dot{\phi} \in \Phi_{\simp}(\dot{G}_i)$; the datum $\dot{G}_i = (\dot{G}_i,\dot{\xi}_i) \in \dot{\widetilde{\mathcal{E}}}_{\simp}(N_i)$ being a globalization of $G_i=(G_i,\xi_i) \in \widetilde{\mathcal{E}}_{\simp}(N_i)$. Similarly for $\dot{\phi}_M$.

In order to apply the global results from section 6 we need to choose $\dot{\phi}$ so that the localization $\dot{\phi}_v$ for $v \in S_{\infty}^u$ has further properties. This is the point of the following:

\begin{proposition}
We can choose the global data:
\[
(\dot{G},\dot{\phi},\dot{M},\dot{\phi}_M,\dot{G}_i,\dot{\phi}_i)
\]
so that the following condition holds:
\bigskip

\noindent (i) $(\dot{G}_u,\dot{\phi}_u,\dot{M}_u,\dot{\phi}_{M,u},\dot{G}_{i,u},\dot{\phi}_{i,u}) =(G,\phi,M,\phi_M,G_i,\phi_i)$, and such that the canonical maps:
\begin{eqnarray}
\mathcal{S}_{\dot{\phi}} &\rightarrow& \mathcal{S}_{\phi} \\
\mathcal{S}_{\dot{\phi}_M} &\rightarrow& \mathcal{S}_{\phi_M} \nonumber
\end{eqnarray}
are isomorphisms.

\bigskip

\noindent (ii) For $v \notin S_{\infty}(u)$, the parameter $\dot{\phi}_v$ is spherical. For $v \in S_{\infty}^u$, the parameters $\dot{\phi}_{i,v}$ belongs to $\Phi_2(\dot{G}_{i,v})$, and the irreducible components of the family of parameters $\{\dot{\phi}^{N_i}_{i,v} \}_{i=1}^r$ are {\it distinct} one-dimensional characters.  

\bigskip

\noindent (iii)(a) Put $V=S_{\infty}^u$. Then the mappings
\begin{eqnarray}
\Pi_{\dot{\phi}_V}  := \bigotimes_{v \in V} \Pi_{\dot{\phi}_v} &\rightarrow & \widehat{\mathcal{S}}_{\dot{\phi}} \cong \widehat{\mathcal{S}}_{\phi} \\
\Pi_{\dot{\phi}_{M,V}} := \bigotimes_{v \in V} \Pi_{\dot{\phi}_{M,v}}  &\rightarrow & \widehat{\mathcal{S}}_{\dot{\phi}_M} \cong \widehat{\mathcal{S}}_{\phi_M} 
\end{eqnarray}
obtained from the combined places in $V$ are surjective.

\bigskip

\noindent (iii)(b) If $l_i=1$ for all $i$, then for every $v \in V$ the following is satisfied: If $\dot{\phi}_v^N$ lies in $\xi_{v,*}^* \Phi(G^*_v)$ for some $G_v^*=(G_v^*,\xi_v^*) \in \widetilde{\mathcal{E}}_{\simp,v}(N)$, then we have $(G_v^*,\xi^*_v)=(\dot{G}_v,\dot{\xi}_v)$ as (equivalences classes of) endoscopic datum in $\widetilde{\mathcal{E}}_{\ellip,v}(N)$.

\bigskip

\noindent (iii)(c) If $l_i >1$ for some $i$ (so that $\dot{G}$ is simple) then for every $v \in V$ the following is satisfied: the kernel of the mapping:
\[
\mathcal{S}_{\dot{\phi}} \rightarrow \mathcal{S}_{\dot{\phi}_v} \rightarrow R_{\dot{\phi}_v}
\]
contains no elements whose image in $R_{\dot{\phi}}$ belongs to $R_{\dot{\phi},\reg}$.
\end{proposition}
\begin{proof}
Conditions (i) and (ii) are clear from the way we apply corollary 7.2.7 to construct the global datum. We only need to show that they can be chosen so that the conditions in (iii) are satisfied.

As for (iii)(a), we treat the case of (7.3.6); the case for (7.3.7) will be similar. The first step is to show that the canonical map
\begin{eqnarray}
\mathcal{S}_{\dot{\phi}} \rightarrow \prod_{v \in V} \mathcal{S}_{\dot{\phi}_v}
\end{eqnarray}
is an injection. Since
\[
Z(\widehat{\dot{G}})^{\Gamma_{\dot{F}}} = Z(\widehat{\dot{G}}_v)^{\Gamma_{\dot{F}_v}} 
\]
for $v \in V=S^u_{\infty}$, it suffices to show that if $s \in S_{\dot{\phi}}$ maps to $S_{\dot{\phi}_v}^0$ for each $v \in V$, then $s \in S_{\dot{\phi}}^0$. We have to allow $G$ and hence $\dot{G}$ not being simple, i.e. $\dot{G}=\dot{G}_O \times \dot{G}_S$, with $\dot{G}_O,\dot{G}_S$ being simple {\it endoscopic datum} of opposite parity.  Similar to section 2.4, denote by $I^{+}_O$ the set of indices $i$ such that $\dot{\phi}_i^{N_i} $ is conjugate self-dual of the same parity as $\dot{G}_O$, and $I^{+}_S$ the set of indices $i$ such that $\dot{\phi}_i^{N_i} $ is conjugate self-dual of the same parity as $\dot{G}_S$. Then we have
\[
S_{\dot{\phi}} = \prod_{i \in I^{+}_O} O(l_i,\mathbf{C}) \times \prod_{i \in I^{+}_S} O(l_i,\mathbf{C}) \times  \mbox{symplectic factors}.
\]
On the other hand, from the induction hypothesis, we have theorem 2.4.10 being valid for the localization of each $\dot{\phi}^{N_i}$ at $v$. Thus if $i \in I^{+}_O$, then for each $v \in V$, the localization $\dot{\phi}^{N_i}_{i,v}$ is conjugate self-dual of the same parity as that of $\dot{\phi}_i^{N_i}$. The same applies to the set of indices $I^{+}_S$. Since the constituents of the set of local parameters $\{ \dot{\phi}^{N_i}_{i,v}\}_{i=1}^r$ are {\it distinct} one-dimensional characters, we see that
\[
S_{\dot{\phi}_v} = \prod_{i \in I^{+}_O} O(l_i,\mathbf{C})^{N_i} \times \prod_{i \in I^{+}_S} O(l_i,\mathbf{C})^{N_i} \times \mbox{symplectic factors}
\]
and the map $S_{\dot{\phi}}  \rightarrow S_{\dot{\phi}_v}$ is given on the orthogonal factors as the obvious diagonal embedding $O(l_i,\mathbf{C}) \hookrightarrow O(l_i,\mathbf{C})^{N_i}$. It follows that we already have the injectivity of $\mathcal{S}_{\dot{\phi}} \rightarrow \mathcal{S}_{\dot{\phi}_v}$, and hence the injectivity of (7.3.8).

Thus dualizing the injectivity of (7.3.8), we have the surjectivity of
\[
\prod_{v \in V} \widehat{\mathcal{S}}_{\dot{\phi}_v} \rightarrow \widehat{\mathcal{S}}_{\dot{\phi}}.
\]
On the other hand, from proposition 7.1.4, the map (for $v \in V$)
\[
\Pi_{\dot{\phi}_v} \rightarrow \widehat{\mathcal{S}}_{\dot{\phi}_v}
\]
has the property that the image generates $\widehat{\mathcal{S}}_{\dot{\phi}_v}$. Thus we need to show that we can choose the global data so that the combined map:
\begin{eqnarray}
\Pi_{\dot{\phi}_V} \rightarrow \prod_{v \in V}      \widehat{\mathcal{S}}_{\dot{\phi}_v}         \rightarrow \widehat{\mathcal{S}}_{\dot{\phi}}
\end{eqnarray}
is surjective. This can be insured simply by ``enlarging $V$"; more precisely, we can replace $\dot{F}$ by a totally real extension $\dot{F}^{\prime}$ of $\dot{F}$, whose set of archimedean places we denote as $S_{\dot{F}^{\prime},\infty}$ (and a corresponding extension $\dot{E}^{\prime}$ of $\dot{E}$ such that $\dot{E}^{\prime}/\dot{F}^{\prime}$ is totally imaginary). We choose correspondingly a global parameter $\dot{\phi}^{\prime}$, such that $\dot{\phi}^{\prime}_{u^{\prime}} = \dot{\phi}_u = \phi$ for some place $u^{\prime}$ of $\dot{F}^{\prime}$ above $u$, and that $\dot{\phi}^{\prime}_{v^{\prime}} = \dot{\phi}_v$ for each $v^{\prime} \in S_{\dot{F}^{\prime},\infty}$ above $v \in S_{\infty}^u$ ({\it c.f.} remark 7.2.8; the condition on the central characters of the parameters $\dot{\phi}^{\prime}_{v^{\prime}}$ over the places $S_{\dot{F}^{\prime},\infty}(u^{\prime})$ can easily be achieved by considering the cases that $u \in S_{\infty}$ and $u \notin S_{\infty}$ separately). Then if $[\dot{F}^{\prime}:\dot{F}]$ is large enough this would insure that (7.3.9), with $\dot{\phi}$ replaced by $\dot{\phi}^{\prime}$ etc., is surjective. 

For (iii)(b), we have $l_i=1$ for all $i$. Hence for $v \in V$, the parameter $\dot{\phi}_v^N$ is multiplicity free by condition (ii). Suppose that $\dot{\phi}^N_v \in \xi^*_{v,*} \Phi(G^*_v)$ for some $G^*_v =(G^*_v,\xi^*_v) \in \widetilde{\mathcal{E}}_{\simp,v}(N)$. Then for each index $i$ the parameter $\dot{\phi}^{N_i}_{i,v}$ must have the same parity as that of $G^*_v$. By theorem 2.4.10 applied to the global parameter $\dot{\phi}_i^{N_i}$, we see that $\dot{G}_i$ must have the same parity as $G^*_v$ for each index $i$. This implies that $\dot{G}$ is simple, and that $\dot{G}_v = G^*_v$ as (equivalence classes of) endoscopic datum.

For (iii)(c), we have $l_j>1$ for some $j$, and $\dot{G}$ is simple. Suppose that $x \in \mathcal{S}_{\dot{\phi}}$ whose image in $R_{\dot{\phi}}$ lies in $R_{\dot{\phi}_{\reg}}$. We need to show that the image $x_v$ of $x$ in $\mathcal{S}_{\dot{\phi}_v}$ does not vanish in $R_{\dot{\phi}_v}$. From the fact that $R_{\dot{\phi},\reg}$ is non-empty, the group $S_{\dot{\phi}}$ must take the form:
\[
S_{\dot{\phi}} = \prod_{i=1}^q O(2,\mathbf{C}) \times \prod^r_{i=q+1} O(1,\mathbf{C}) 
\]
(in particular $l_i \leq 2$ for all $i$) and the (unique) element of $R_{\dot{\phi},\reg}$ is represented by an element of $S_{\dot{\phi}}$ that lies in the non-identity component of each $O(2,\mathbf{C})$ factor. As in the discussion in (iii)(a), for each $v \in V$, we then have 
\[
S_{\dot{\phi}_v} = \prod_{i=1}^q O(2,\mathbf{C})^{N_i} \times \prod^r_{i=q+1} O(1,\mathbf{C})^{N_i} 
\]
and the map $\mathcal{S}_{\dot{\phi}} \rightarrow \mathcal{S}_{\dot{\phi}_v}$ is given by the diagonal map on the orthogonal factors. In particular the image $x_v$ of $x$ in $\mathcal{S}_{\dot{\phi}_v}$ lies in the non-identity components of each $O(2,\mathbf{C})$-factor, and thus $x_v$ does not vanish in $R_{\dot{\phi}_v}$.
\end{proof}

\begin{rem}
\end{rem}
More generally, given the irreducible components $\phi_1^{N_1},\cdots,\phi_r^{N_r}$ of $\phi^N$ as in (7.3.1), and their globalization $\dot{\phi}_1^{N_1}, \cdots, \dot{\phi}_r^{N_r}$. Then we can consider more generally any parameter $\phi^{N^{\prime}} \in \widetilde{\Phi}(N^{\prime})$ generated by the simple parameters $\phi_1^{N_1},\cdots,\phi_r^{N_r}$, i.e. 
\[
\phi^{N^{\prime}} = l_1^{\prime} \phi_1^{N_1} \oplus \cdots \oplus l_r^{\prime} \phi_r^{N_r}
\] 
for {\it any} positive integers $l_1^{\prime}, \cdots,l_r^{\prime}$ (thus $N^{\prime} = l_1^{\prime} N_1 + \cdots + l_r^{\prime} N_r$; in particular we allow $N^{\prime}$ to be larger than $N$), and the corresponding global parameter
\[
\dot{\phi}^{N^{\prime}} = l_1^{\prime} \dot{\phi}_1^{N_1} \boxplus \cdots \boxplus l_r^{\prime} \dot{\phi}_r^{N_r} \in \dot{\widetilde{\Phi}}(N^{\prime})
\]
{\it c.f.} the discussions in section 6.1. Suppose that $\phi^{N^{\prime}} = \xi_*^{\prime} \phi^{\prime}$ for $\phi^{\prime} \in \Phi(G^{\prime})$ and $G^{\prime}=(G^{\prime},\xi^{\prime}) \in \widetilde{\mathcal{E}}_{\simp}(N^{\prime})$. Then with a corresponding globalization $\dot{G}^{\prime}=(\dot{G}^{\prime},\dot{\xi}^{\prime}) \in \dot{\widetilde{\mathcal{E}}}_{\simp}(N^{\prime})$, and $\dot{\phi}^{\prime}$ of $\phi^{\prime}$ such that $\dot{\phi}^{N^{\prime}} = \dot{\xi}^{\prime}_* \dot{\phi}^{\prime}$, we have again have the validity of proposition 7.3.1 with respect to $\dot{G}^{\prime}$ and $\dot{\phi}^{\prime}$, etc.

\bigskip

Finally in the case where $F=\mathbf{R}$ we also need a variant of proposition 7.3.1:
\begin{lemma}
Suppose that $F=\mathbf{R}$, and $\phi \in \Phi_{\bdd}(G)$ as in (7.3.1). Assume that the infinitesimal characters of the distinct irreducible components of $\phi^N$ are in general position. Then we can choose the global data:
\[
(\dot{G},\dot{\phi},\dot{M},\dot{\phi}_M,\dot{G}_i,\dot{\phi}_i)
\]
so that we have $\dot{\phi}_v=\phi$ for each $v \in S_{\infty}$, and such that all the conditions of proposition 7.3.1 are satisfied.
\end{lemma}
\begin{proof}
Since $E=\mathbf{C}$ and $L_E=\mathbf{C}^{\times}$, the irreducible components of $\phi^N$ are just one-dimensional characters of $\mathbf{C}^{\times}$, and so this reduces to the globalization result for one-dimensional characters, which is just the elementary abelian case of proposition 7.3.1. 
\end{proof}

\subsection{The local intertwining relation}

In this subsection we prove the main technical result, the local intertwining relation for generic parameters, stated as theorem 3.4.3. It is based on the partial result in the archimedean case obtained by Shelstad (proposition 7.1.2), combined with the global results from trace formula comparisons in section 6. 

Thus let $\phi^N \in \widetilde{\Phi}_{\bdd}(N)$. The descent argument used in the proof of proposition 5.7.4, which applies equally well to the current local setting, shows that for any $G^*=(G^*,\xi^*) \in \widetilde{\mathcal{E}}_{\simp}(N)$ such that $\phi^N = \xi_* \phi$ for $\phi \in \Phi_{\bdd}(G) \smallsetminus \Phi_2(G)$, the local intertwining relation is valid for $\phi$ (with respect to $G^*$), {\it unless} the parameter $\phi^N$ belongs to one of the following three cases:

\begin{enumerate}
\item  $\phi^N$ belongs to $\xi_* \phi$, with $\phi \in \Phi_{\ellip}(G)$ for $(G,\xi)$ among one of the two simple (equivalence classes of) endoscopic data $\widetilde{\mathcal{E}}_{\simp}(N)$.
\bigskip

\item The local analogue of the case (5.7.12).
\bigskip

\item The local analogue of the case (5.7.13).
\end{enumerate}

In case (1) we have
\begin{eqnarray}
& & \phi^N = 2 \phi_1^{N_1} \oplus \cdots \oplus 2 \phi_q^{N_q} \oplus \phi_{q+1}^{N_{q+1}} \oplus \cdots \oplus \phi_r^{N_r}, \,\ q \geq 1 \\
& &  S_{\phi} = \prod_{i=1}^q O(2,\mathbf{C}) \times \prod_{i=q+1}^r O(1,\mathbf{C}).    \nonumber
\end{eqnarray}

While cases (2) and (3) corresponds to the following situation: there is a $G=(G,\xi) \in \widetilde{\mathcal{E}}_{\simp}(N)$ such that $\phi^N \in \xi_* \phi$ for $\phi \in \Phi_{\bdd}(G)$, and such that
\begin{eqnarray}
& & \phi^N = 2 \phi_1^{N_1} \oplus \phi_2^{N_2} \oplus \cdots \oplus \phi_r^{N_r} \\
& & S_{\phi} = Sp(2,\mathbf{C}) \times \prod_{i=2}^r O(1,\mathbf{C}) \nonumber
\end{eqnarray}
in the case (2), or
\begin{eqnarray}
& & \phi^N = 3 \phi_1^{N_1} \oplus \phi_2^{N_2} \oplus \cdots \oplus \phi_r^{N_r} \\
& & S_{\phi} = O(3,\mathbf{C}) \times \prod_{i=2}^r O(1,\mathbf{C}) \nonumber
\end{eqnarray}
in the case (3).

Thus we must treat these remaining three cases. We fix $G =(G,\xi) \in \widetilde{\mathcal{E}}_{\simp}(N)$, and $\phi \in \Phi_{\bdd}(G)$ such that $\phi^N=\xi_* \phi$, as in one of the three cases above. In particular $\phi^N$ is as in (7.3.1) of the previous subsection. We assume that $\phi \notin \Phi_2(G)$, i.e. there is a proper Levi $M$ of $G$ and $\phi_M \in \Phi_2(M)$ mapping to $\phi$. Recall that we equip $M$ with the $L$-embedding $\xi$ and so regard $M=(M,\xi)$ as a twisted endoscopic datum in $\widetilde{\mathcal{E}}(N)$. Recall that the local intertwining relation asserts:
\begin{eqnarray*}
f^{\prime}_G(\phi,s) = f_G(\phi,u), f \in \mathcal{H}(G)
\end{eqnarray*} 
for $s \in \overline{S}_{\phi}$ and $u \in \mathcal{N}_{\phi}$ having the same image in $\mathcal{S}_{\phi}$. In case (1) a descent argument shows that it suffices to establish the case where the common image of $s$ and $u$ in $\mathcal{S}_{\phi}$ lies in $\mathcal{S}_{\phi,\ellip}$.

By proposition 7.3.1, we can globalize the data
\[
(G,\phi,M,\phi_M,G_i,\phi_i)
\]
to data:
\[
(\dot{G},\dot{\phi},\dot{M},\dot{\phi}_{\dot{M}},\dot{G}_i,\dot{\phi}_i)
\]
with the specific properties as listed in (i)-(iii) of the proposition. In the notation of (7.4.1) - (7.4.3):
\[
\phi^N =\xi_* \phi = l_1 \phi_1^{N_1} \oplus \cdots \oplus l_r \phi_r^{N_r}
\]
and
\[
\dot{\phi}^N = \dot{\xi}_* \dot{\phi}= l_1 \dot{\phi}_1^{N_1} \oplus \cdots \oplus l_r \dot{\phi}_r^{N_r}
\]
and we can then form the family of global parameters:
\begin{eqnarray}
\dot{\widetilde{\mathcal{F}}} = \widetilde{\mathcal{F}}(\dot{\phi}_1^{N_1},\cdots,\dot{\phi}_r^{N_r})
\end{eqnarray}
generated by the simple parameters $\dot{\phi}_1^{N_1},\cdots,\dot{\phi}_r^{N_r}$, and also the family $\dot{\widetilde{\mathcal{F}}}(\dot{G})$, as in section 6.1. In particular $\dot{\phi}^N \in \dot{\widetilde{\mathcal{F}}}$ and $\dot{\phi} \in \dot{\widetilde{\mathcal{F}}}(\dot{G})$. The family $\dot{\widetilde{\mathcal{F}}}$ provides the global input we need to establish the local intertwining relation, in accordance with the results from section 6.4. We state this as:

\begin{proposition}
For $x \in \mathcal{S}_{\phi}$, we denote by $\dot{x} \in \mathcal{S}_{\dot{\phi}}$ the corresponding element of $\mathcal{S}_{\dot{\phi}}$ under the isomorphism (7.3.5). Then:

\noindent In case (1), we have the identity:
\begin{eqnarray}
\sum_{x \in \mathcal{S}_{\phi,\ellip}} (\dot{f}_{\dot{G}}^{\prime}(\dot{\phi},\dot{x}) - \dot{f}_{\dot{G}}(\dot{\phi},\dot{x}) ) =0.
\end{eqnarray}

\noindent In cases (2) and (3), we have the identity:
\begin{eqnarray}
\sum_{x \in \mathcal{S}_{\phi}} (\dot{f}_{\dot{G}}^{\prime}(\dot{\phi},\dot{x}) - \dot{f}_{\dot{G}}(\dot{\phi},\dot{x}) ) =0.
\end{eqnarray}
\end{proposition}
\begin{proof}
The validity of (7.4.6) in case (2) and (3) is exactly the content of corollary 5.7.5, on noting that since we are dealing with generic parameters the $\epsilon$ sign character and the element $s_{\phi}$ are trivial. For case (1), the validity of (7.4.5) is exactly the content of Proposition 6.4.4, the part on the vanishing of the right hand side of (6.2.2), and together with Proposition 6.4.5, the part on the vanishing of the right hand side of (6.2.10). In applying Proposition 6.4.4 and 6.4.5 we are using the fact that Assumption 6.4.1 of for the family (7.4.4) of global parameters $\dot{\widetilde{\mathcal{F}}}$, for which the results of section 6.4 are contingent upon, are satisifed by part (iii) of Proposition 7.3.1.
\end{proof}

\begin{rem}
\end{rem}
\noindent Recall that $\mathcal{S}_{\phi_M}$ is naturally a subgroup of $\mathcal{S}_{\phi}$. We note that in case (1) the set $\mathcal{S}_{\phi,\ellip}$ is a torsor under $\mathcal{S}_{\phi_M}$; in case (2) and (3) $\mathcal{S}_{\phi,\ellip}$ is empty and $\mathcal{S}_{\phi_M} = \mathcal{S}_{\phi}$.

\bigskip
We can now begin to extract the local intertwining relation from the global identity of Proposition 7.4.1. We first establish:

\begin{proposition}
Suppose $F$ is non-archimedean and $E/F$ is unramified. Then the local intertwining relation 
\[
f^{\prime}_G(\phi,x) = f_G(\phi,x), \,\ f \in \mathcal{H}(G)
\]
is valid if $\phi \in \Phi_{\bdd}(G)$ is a spherical parameter. In addition we have
\[
f^{\prime}_G(\phi,x) = f_G(\phi,x) =1
\]
for any $x \in \mathcal{S}_{\phi}$, if $f$ is the characteristic function of the (standard) special maximal compact subgroup of $G(F)$.
\end{proposition}
\begin{proof}
Since $\phi$ is a spherical parameter, we have $N_i=1$ for all index $i$, and $M=M_0$ is the minimal Levi subgroup, which is just the standard diagonal maximal torus. By the previous discussion, we need to consider only cases (1), (2) and (3) above. In is in fact easy to check that in the unramified case we must have $r=q=1$ in case (1), $r \leq 2$ in case (2), and $r=1$ in case (3). Note that $\mathcal{S}_{\phi}$ is of order two in case (1), and is trivial in cases (2) and (3). The the underlying group $G$ is just $U_{E/F}(2)$ or $U_{E/F}$(3). These cases should be known in the literature. In any case we give an argument along the lines of lemma 6.4.1 of \cite{A1}.

The treatment of these three cases are similar, so we only treat case (1). By applying proposition 7.4.1, we have, for $x$ the unique element of $\mathcal{S}_{\phi,\ellip}$:
\[
 \dot{f}_{\dot{G}}^{\prime}(\dot{\phi},\dot{x}) - \dot{f}_{\dot{G}}(\dot{\phi},\dot{x}) =0, \,\ \dot{f} \in \mathcal{H}(\dot{G})
\]
and on choosing $\dot{f}$ to be decomposable we then have:
\begin{eqnarray}
\end{eqnarray}
\begin{eqnarray*}
  \dot{f}^{\prime}_{\infty,\dot{G}} (\dot{\phi}_{\infty},\dot{x}_{\infty}) \prod_{v \notin S_{\infty}}  \dot{f}^{\prime}_{v,\dot{G}}(\dot{\phi}_v,\dot{x}_v)  - \dot{f}_{\infty,\dot{G}}(\dot{\phi}_{\infty},\dot{x}_{\infty}) \prod_{v \notin S_{\infty}}  \dot{f}_{v,\dot{G}}(\dot{\phi}_v,\dot{x}_v) =0. 
\end{eqnarray*}

\noindent For $v \in S_{\infty}$, the local intertwining relation is satisfied for the archimedean parameter $\dot{\phi}_v$. Indeed, the parameter $\dot{\phi}_{\dot{M},v}$ just corresponds to an unitary character $\dot{\pi}_{\dot{M},v}$ of $\dot{M}_0(\dot{F}_v)$ (and the packet $\Pi_{\dot{\phi}_{\dot{M},v}}$ is course just the singleton $\{\dot{\pi}_{\dot{M},v}\}$, and is generic trivially), and so by the results of Shelstad quoted as proposition 7.1.3, the character $\epsilon_{\dot{\pi}_{\dot{M},v}}$ is trivial, and hence the local intertwining relation holds for $\dot{\phi}_v$, by (7.1.5).

Hence we can write (7.4.7) as
\begin{eqnarray}
 \dot{f}_{\infty,\dot{G}} (\dot{\phi}_{\infty},\dot{x}_{\infty}) \Big(\prod_{v \notin S_{\infty}}  \dot{f}^{\prime}_{v,\dot{G}}(\dot{\phi}_v,\dot{x}_v)  -  \prod_{v \notin S_{\infty}}  \dot{f}_{v,\dot{G}}(\dot{\phi}_v,\dot{x}_v)\Big) =0. 
\end{eqnarray}
We can certainly choose $\dot{f}_{\infty}$ so that $ \dot{f}_{\infty,\dot{G}} \neq 0$. Hence we have
\[
\prod_{v \notin S_{\infty}}  \dot{f}^{\prime}_{v,\dot{G}}(\dot{\phi}_v,\dot{x}_v)  =  \prod_{v \notin S_{\infty}}  \dot{f}_{v,\dot{G}}(\dot{\phi}_v,\dot{x}_v).
\]
Similarly we can choose $\dot{f}_v$ for $v \notin S_{\infty}(u)$ so that $\dot{f}_v(\dot{\phi}_v,\dot{x}_v) \neq 0$. Hence we obtain:
\begin{eqnarray}
f_G(\phi,x) = e(x) f^{\prime}_G (\phi,x), \,\ f \in \mathcal{H}(G)
\end{eqnarray} 
for a non-zero constant $e(x)$ that is independent of $f$. We can now apply Lemma 2.5.5 of \cite{A1} to our situation. There are three conditions in Lemma 2.5.5 of {\it loc. cit.} Condition (iii) is our (7.4.9), while condition (ii) is trivially satisfied as noted above. Condition (i) is that the characteristic of the residue field of $F$ is not two. Thus whenever the residue field of $F$ is not of characteristic two, we can apply Lemma 2.5.5 of \cite{A1} to conclude the validity of the local intertwining relation for $\phi$. 

Thus it remain to deduce the case where the residue field of $F$ is of characteristic two from the results already established for odd residue characteristic. At this point, the argument is then the same as that of lemma 6.4.1 of \cite{A1}, so we refer to {\it loc. cit.} for completing the argument in the case of residue characteristic two.

For the final assertion, we only need to know that (with the notation as in section 7.1 concerning spherical parameters) for the unique spherical constituent $\pi$ of the induced representation $\mathcal{I}_{P_0}(\pi_{M_0})$, the action of the normalized intertwining operator $R_{P_0}(w,\widetilde{\pi}_{M_0},\phi_{M_0})$ on $\pi$ is trivial; this follows from the result of Shahidi \cite{S,S2}, {\it c.f.} the discussion of Corollary 2.5.2 in \cite{A1}. From this it follows that if $f$ is the characteristic function of the standard special maximal compact subgroup of $G(F)$, then
\[
f_G(\phi,x) =1
\]     
for any $x \in \mathcal{S}_{\phi}$. 
\end{proof}

We now treat the case where $F$ is archimedean. For archimedean $F$, case (2) and (3) already follows from Shelstad's results, namely proposition 7.2.1. Indeed, the $R$-group $R_{\phi}$ is trivial in cases (2) and (3), so the character $\epsilon_{\pi_M}$ for $\pi_M \in \Pi_{\phi_M}$ occuring in (7.1.5) is trivial.

\begin{lemma}
Suppose $F$ is archimedean, and that $\phi \in \Phi(G)$ is as in case (1) is in relative general position (in the sense that each component $\phi_i$ is in general position). Then there exists $\epsilon_1 \in \mathcal{S}_{\phi_M}$ such that for any $x \in \mathcal{S}_{\phi,\ellip}$, we have
\begin{eqnarray}
f^{\prime}_G(\phi, x ) = f_G(\phi,x \epsilon_1), \,\ f \in \mathcal{H}(G).
\end{eqnarray}
\end{lemma}  
\begin{proof}
This is the same as in the proof of lemma 6.4.2 of \cite{A1}. For this proposition we need to apply the variant 7.3.3 of the globalization result. Thus this time we choose the global parameter $\dot{\phi}$ such that $\dot{\phi}_v=\phi$ for every $v \in S_{\infty}$, and such that the conditions of proposition 7.3.1 still holds. We then apply the global identity (7.4.5) to $\dot{\phi}$ and choosing $\dot{f}$ to be decomposable:
\begin{eqnarray}
\end{eqnarray}
\begin{eqnarray*}
\sum_{x \in \mathcal{S}_{\phi,\ellip}} \Big( \dot{f}^{\prime}_{\infty,\dot{G}}(\dot{\phi},\dot{x}) \prod_{v \notin S_{\infty}}  \dot{f}^{\prime}_{v,\dot{G}_v}(\dot{\phi}_v,\dot{x}_v) -   \dot{f}_{\infty,\dot{G}}(\dot{\phi},\dot{x}) \prod_{v \notin S_{\infty}} \dot{f}_{v,\dot{G}_v}(\dot{\phi}_v,\dot{x}_v)       \Big) = 0.
\end{eqnarray*}  
For $v \notin S_{\infty}$ the parameter $\dot{\phi}_v$ is spherical, and so we have the local intertwining relation $\dot{f}^{\prime}_{v,\dot{G}_v}(\dot{\phi}_v,\dot{x}_v)  = \dot{f}_{v,\dot{G}_v}(\dot{\phi}_v,\dot{x}_v)$ by proposition 2.4.3. Furthermore, we have seen above that for spherical $\dot{\phi}_v$ we have $\dot{f}_{v,\dot{G}_v}(\dot{\phi}_v,\dot{x}_v)=1$ for $\dot{f}_v$ the characteristic function of the (standard) special maximal compact subgroup of $\dot{G}(\dot{F}_v)$. Hence we have
\begin{eqnarray}
\sum_{x \in \mathcal{S}_{\phi,\ellip}} \Big( \dot{f}^{\prime}_{\infty,\dot{G}}(\dot{\phi},\dot{x}) -   \dot{f}_{\infty,\dot{G}}(\dot{\phi},\dot{x})        \Big) = 0.
\end{eqnarray} 
In {\it loc. cit.} a Fourier tranform argument on the finite abelian group $\mathcal{S}_{\phi_M}$ is applied to (7.4.12) to obtain (7.4.10), and so we refer to {\it loc. cit.} for the proof (we remark that the proof also uses the result of Shelstad stated as proposition 7.1.3). 
\end{proof}

\begin{proposition}
Suppose $F$ is archimedean and that $\phi$ is in relative general position. Then the local intertwining relation is valid for $\phi$.
\end{proposition}
\begin{proof}
We have already seen that for archimedean $F$, we only need to treat case (1), and where $x \in \mathcal{S}_{\phi,\ellip}$. Thus we have
\begin{eqnarray*}
& & \phi^N = 2 \phi_1^{N_1} \oplus \cdots \oplus 2 \phi_q^{N_q} \oplus \phi_{q+1}^{N_{q+1}} \oplus \cdots \oplus \phi_r^{N_r}  \\
& & S_{\phi} = \prod_{i=1}^q O(2,\mathbf{C}) \times \prod_{i=q+1}^r O(1,\mathbf{C})
\end{eqnarray*}
(in the archimedean case we of course have $N_i=1$ for all indices $i$.) In particular since we are in case (1), all the components $\phi_i^{N_i}$ are conjugate self-dual with the same parity as $G=(G,\xi)$. 

Choose global endoscopic datum $\dot{G} \in \dot{\widetilde{\mathcal{E}}}_{\simp}(N)$ over $\dot{F}$ as in lemma 7.2.1 and 7.2.2; thus $(\dot{F}_u,\dot{G}_u)=(F,G)$ for some $u \in S_{\infty}$. Define new degrees $N_j^{\#}$ for $j=1,\cdots,q+1$ as follows:
\begin{eqnarray*}
& & N_j^{\#} = N_j=1,\,\ j=1,\cdots,q \\
& & N_j^{\#} = N_{q+1} + \cdots, + N_r = r-q, \,\ j=q+1.
\end{eqnarray*} 
We now use the following trick of Arthur ({\it c.f.} proof of lemma 6.4.3 in \cite{A1}): arbitrarily choose a non-archimedean place $u^{\#}$ of $\dot{F}$ that does not split in $\dot{E}$, and consider the localization $(G^{\#},\xi^{\#})= (\dot{G}_{u^{\#}},\dot{\xi}_{u^{\#}}) \in \dot{\widetilde{\mathcal{E}}}_{\simp,u^{\#}}(N)$. We can choose parameters:
\begin{eqnarray}
\phi_j^{\#, N_j^{\#}} \in \widetilde{\Phi}_{\simp,u^{\#}}(N_j^{\#}), \,\ j =1,\cdots,q+1
\end{eqnarray}
with the following requirement: $\phi^{\#,N_j^{\#}} = \xi_{j,*}^{\#} \phi_j^{\#}$ for $\phi_j^{\#} \in \Phi_{\simp}(G^{\#}_j)$, with $(G^{\#}_j,\xi_j^{\#}) \in \dot{\widetilde{\mathcal{E}}}_{\simp,u^{\#}}(N_j^{\#})$ that have the same parity as $(G^{\#},\xi^{\#})$. The existence of the local parameters $\phi_j^{\#, N_j^{\#}}$ is easily seen by considering parameters that is an irreducible $N_j^{\#}$-dimensional representation of $\SU(2)$, for the $\SU(2)$ component of $L_{\dot{E}_{u^{\#}}}$, and a character $\chi \in \mathcal{Z}^{\pm}_{\dot{E}_{u^{\#}}}$ on the $W_{\dot{E}_{u^{\#}}}$ component of $L_{\dot{E}_{u^{\#}}}$ (note that since $N_j^{\#} <N$ the local classification theorems are valid for parameters of degree $N_j^{\#}$, by the induction hypothesis). 

Note that if we put we 
\[
\phi^{\#,N} :=2 \phi_1^{N_1^{\#}} \oplus \cdots \oplus 2 \phi_q^{N_q^{\#}} \oplus \phi_{q+1}^{N_{q+1}^{\#}}
\]
then $\phi^{\#,N} = \xi^{\#}_* \phi^{\#}$, for a parameter $\phi^{\#} \in \Phi_{\ellip}(G^{\#})$.

We now globalize the parameter $\phi^{\#}$ to a global parameter $\dot{\phi}^{\#}$ over $\dot{F}$ such that $\dot{\phi}^{\#}_{u^{\#}} = \phi^{\#}$; in addition we also impose the conditions:
\begin{enumerate}
\item $\dot{\phi}_u^{\#} = \phi$

\bigskip

\item $\dot{\phi}^{\#}_v$ is spherical for $v \notin S_{\infty}(u^{\#})$

\bigskip

\item Condition (i) of proposition 7.3.1 is valid for $\dot{\phi}^{\#}$ with respect to the place $u^{\#}$, while conditions (ii) and (iii) of proposition 7.3.1 for $\dot{\phi}^{\#}$ are valid with respect to the set of places $V= S_{\infty}^{u^{\#}}=S_{\infty}$.
\end{enumerate}
That such $\phi^{\#}$ can be chosen is a simple variant of the globalization constructions leading up to proposition 7.3.1 ({\it c.f.} remark 7.2.6). 

Once we have the global parameter $\dot{\phi}^{\#}$, we can apply (7.4.5) to $\dot{\phi}^{\#}$ and obtain:
\begin{eqnarray}
\sum_{\dot{x} \in \mathcal{S}_{ \dot{\phi}^{\#},\ellip}} \Big( \dot{f}^{\prime}_{\dot{G}}(\dot{\phi}^{\#},\dot{x}) - \dot{f}_{\dot{G}}(\dot{\phi}^{\#},\dot{x}) \Big) =0.
\end{eqnarray}

\noindent The key point is that by the chocie of $\dot{\phi}^{\#}$, the set $\mathcal{S}_{\dot{\phi}^{\#},\ellip}$ consists of a single point $\dot{x}_1$, and it can be identitfied with a base point $x_1$ in the larger set $\mathcal{S}_{\phi,\ellip}$. Hence on choosing $\dot{f}$ to be a decomposable function in (7.4.14), we obtain:
\[
f^{\prime}_G(\phi,x_1) \prod_{v \neq u} \dot{f}^{\prime}_{v ,\dot{G}_v}( \dot{\phi}^{\#}_v, \dot{x}_{1,v}) = f_G(\phi,x_1) \prod_{v \neq u} \dot{f}_{v ,\dot{G}_v}( \dot{\phi}^{\#}_v, \dot{x}_{1,v}).
\]
We can choose $\dot{f}_v$ for $v \neq u$ such that $\prod_{v \neq u} \dot{f}_{v ,\dot{G}_v}( \dot{\phi}^{\#}_v, \dot{x}_{1,v}) \neq 0$. Hence we can write
\begin{eqnarray}
e(\phi) f^{\prime}_G(\phi,x_1)  = f_G(\phi,x_1), \,\ f \in \mathcal{H}(G)
\end{eqnarray}
for a scalar $e(\phi)$ that is independent of $f$. Combining (7.4.15) with lemma 7.4.4, we obtain:
\[
 f_G(\phi,x_1)= e(\phi)  f_G(\phi,x_1 \epsilon_1)
\]   
for an element $\epsilon_1 \in \mathcal{S}_{\phi_M}$. But it is easy to see that the two linear forms $ f_G(\phi,x_1)$ and $ f_G(\phi,x_1 \epsilon_1)$ are linearly independent, unless $\epsilon_1 =1$. Hence lemma 7.4.4 combined with $\epsilon_1=1$ give the local intertwining relation for $\phi$.  
\end{proof}

We can finally prove the main result of this subsection:
\begin{proposition}
The local intertwining relation is valid for $\phi$ over any $F$.
\end{proposition}
\begin{proof}
Again it remains to treat the cases (1), (2) and (3). As before we illustrate for case (1), as the treatment of (2), (3) are similar. It suffices to prove:
\[
f^{\prime}_G(\phi,x) = f_G(\phi,x), \,\ f \in \mathcal{H}(G)
\] 
for $x \in \mathcal{S}_{\phi,\ellip}$. From the global identity (7.4.5), and on choosing $\dot{f}$ to be decomposable, and such that $\dot{f}_u=f \in \mathcal{H}(G)$, we have:
\begin{eqnarray}
\end{eqnarray}
\begin{eqnarray*}
 \sum_{x \in \mathcal{S}_{\phi,\ellip}} \Big( f^{\prime}_G (\phi,x) \prod_{v \neq u}  \dot{f}^{\prime}_{v,\dot{G}}(\dot{\phi}_v,\dot{x}_v)  - f_G(\phi,x) \prod_{v \neq u}  \dot{f}_{v,\dot{G}}(\dot{\phi}_v,\dot{x}_v) \Big)=0. 
\end{eqnarray*}
For $v \notin S_{\infty}(u)$, the parameter $\dot{\phi}_v$ is spherical, and so the local intertwining relation is valid by proposition 7.4.3 (if $v \notin S_{\infty}(u)$ splits in $\dot{E}$, then $\dot{G}_v = \GL_{N/\dot{F}_v}$, and the local intertwining relation for $\dot{\phi}_v$ becomes trivial); for $v \in S_{\infty}^u$, the archimedean parameter $\dot{\phi}_v$ is in relative general position, and so the local intertwining relation is valid by proposition 7.4.5. Thus we have:
\begin{eqnarray}
\end{eqnarray}
\begin{eqnarray*}
\sum_{x \in \mathcal{S}_{\phi,\ellip}} \prod_{v \in S^u_{\infty}}  \dot{f}_{v,\dot{G}}(\dot{\phi}_v,\dot{x}_v) \prod_{v \notin S_{\infty}(u)}  \dot{f}_{v,\dot{G}}(\dot{\phi}_v,\dot{x}_v) \Big( f^{\prime}_G (\phi,x)- f_G(\phi,x)  \Big)=0. 
\end{eqnarray*}
Recall that in case (1) $\mathcal{S}_{\phi,\ellip}$ is a torsor under $\mathcal{S}_{\phi_M}$. By part (iii)(a) of proposition 7.3.1, specifically the surjectivity of (7.3.7), we can isolate the term corresponding to any $x \in \mathcal{S}_{\phi,\ellip}$, and so obtain
\begin{eqnarray}
\prod_{v \notin S_{\infty}(u)}  \dot{f}_{v,\dot{G}}(\dot{\phi}_v,\dot{x}_v) \Big( f^{\prime}_G (\phi,x)- f_G(\phi,x)  \Big) =0.
\end{eqnarray}
Since we can choose $\dot{f}_v$ for $v \notin S_{\infty}(u)$ so that $\prod_{v \notin S_{\infty}(u)}  \dot{f}_{v,\dot{G}}(\dot{\phi}_v,\dot{x}_v) \neq 0$ we thus finally obtain:
\[
f^{\prime}_G (\phi,x) = f_G(\phi,x).
\]
\end{proof}

With the proof of the local intertwining relation for generic parameters given by proposition 7.4.6, we also obtain:
\begin{corollary}
For $w^0 \in W^0_{\phi}$, we have
\[
R_P(w^0,\widetilde{\pi}_M,\phi_M) =1.
\]
\end{corollary}
\begin{proof}
Same as the way Corollary 6.4.5 is deduced from Propositin 6.4.4 in \cite{A1}.
\end{proof}

\subsection{Elliptic orthogonality relation}
With the completion of the local intertwining relation for generic parameters, we see in particular that all the local theorems are established for generic parameters in the archimedean case. Thus we can assume that $F$ is non-archimedean in the rest of section 7. In this subsection $G \in \widetilde{\mathcal{E}}_{\simp}(N)$.

In order to complete the local classification in the non-archimedean case, the global inputs have to be supplemented by local information from the elliptic orthogonality relation. There is no essential difference with section 6.5 of \cite{A1}, so we will be brief. 

Recall the subspace
\[
\mathcal{I}_{\cusp}(G) \subset \mathcal{I}(G)
\] 
consisting of invariant orbital integrals $f_G \in \mathcal{I}(G)$ which, when considered as a function on the set of regular semi-simple elements of $G(F)$, are supported on the set of elliptic elements. The space $\mathcal{I}_{\cusp}(G)$ has a spectral interpretation in terms of the sets $T_{\ellip}(G)$, which we recall.

First we have the set $T_{\temp}(G)$ consisting of $W_0^G$-orbits of triplets
\[
\tau = (M,\pi_M,r)
\]
where $M$ is a Levi subgroup of $G$, $\pi_M \in \Pi_{2,\temp}(M)$, and $r \in R(\sigma)$, where $R(\pi_M)$ is the representation theoretic $R$-group of $\pi_M$ given by $R(\pi_M) = W(\pi_M)/W(\pi_M)^0$ ({\it c.f.} the discussion in section 3.5 of \cite{A1} for general discussion of the representation theoretic $R$-group). The subset 
\[
T_{\ellip}(G) \subset T_{\temp}(G)
\]
consists of triplets $\tau = (M,\pi_M,r)$ with $r$ in the set $R_{\reg}(\pi_M)$ of regular elements, i.e.
\[
d(\tau):=\det(r-1)_{\mathfrak{a}_M^G} \neq 0.
\] 
Then there is a linear isomorphism from $\mathcal{I}_{\cusp}(G)$ to the space of functions of the set $T_{\ellip}(G)$ of finite support, given by sending $f_G \in \mathcal{I}_{\cusp}(G)$ to the function:
\[
\tau \mapsto f_G(\tau). 
\]
Here if $M=G$, then $\pi_G  \in \Pi_{2,\temp}(G)$ (and $r$ is trivial), and $f_G(\tau)$ is the usual character $f_G(\pi_G)$. If $M \neq G$, then we can apply our induction hypothesis on local theorems, and see that $\pi_M \in \Pi_{2,\temp}(M)$ belongs to a unique packet $\Pi_{\phi_M}$ associated to a parameter $\phi_M \in \Phi_{2,\bdd}(M)$. Then
\[
f_G(\tau) = \tr( R_P(r,\widetilde{\pi}_M,\phi_M) \mathcal{I}_P(\pi_M,f)).
\]
For future reference we also recall the subspace $\mathcal{H}_{\cusp}(G) \subset \mathcal{H}(G)$ consisting of functions $f \in \mathcal{H}(G)$ such that $f_G \in \mathcal{I}_{\cusp}(G)$.

\bigskip

Recall the endoscopic correspondence for parameters:
\begin{eqnarray}
(G^{\prime},\phi^{\prime}) \leftrightarrow (\phi,s).
\end{eqnarray}
We will need the correspondence between the two sets:
\[
\mathfrak{X}_{\ellip}(G)=\{ G^{\prime} \in \mathcal{E}_{\ellip}(G), \phi^{\prime} \in \Phi_2(G^{\prime})   \}
\]
and
\[
\mathcal{Y}_{\ellip}(G)=\{ \phi \in \Phi_{\ellip}(G), s \in \overline{S}_{\phi,\ellip}\}
\] 
for which (7.5.1) restricts to a bijective correspondence between $\mathfrak{X}_{\ellip}(G)$ and $\mathcal{Y}_{\ellip}(G)$. 

Given $(\phi,s) \in \mathcal{Y}_{\ellip}(G)$, the linear form
\begin{eqnarray}
f \mapsto f^{\prime}(\phi^{\prime}), \,\ f \in \mathcal{H}(G)
\end{eqnarray}
(where $(G^{\prime},\phi^{\prime}) \in \mathfrak{X}_{\ellip}(G)$ corresponds to $(\phi,s)$) is defined from our induction hypothesis when $ \phi \in \Phi_{\ellip}(G) \smallsetminus \Phi_2(G)$. In the case where $\phi \in \Phi_2(G)$, we still need to establish the existence of (7.5.2). In any case, assuming the existence of (7.5.2) (which we recall depends only on the image $x$ of $s$ in $\mathcal{S}_{\phi,\ellip}$ if it is defined), we can write its restriction to $\mathcal{H}_{\cusp}(G)$ as follows (by using the spectral intepretation for $\mathcal{I}_{\cusp}(G)$ above):
\begin{eqnarray}
f^{\prime}(\phi^{\prime}) = \sum_{\tau \in T_{\ellip}(G)} c_{\phi,x}(\tau) f_G(\tau), \,\ f \in \mathcal{H}_{\cusp}(G)
\end{eqnarray}
for uniquely determined scalar coefficients $c_{\phi,x}(\tau)$. 

Suppose that $\phi \in \Phi_{\bdd}(G) \smallsetminus \Phi_2(G)$. Then the local intertwining relation established in the last subsection allows us to express the coefficients $c_{\phi,x}(\tau)$ as follows. For this purpose, recall the notation used before: for $\phi \in \Phi_{\bdd}(G)$ and $M$ a proper Levi subgroup of $G$ such that there exists $\phi_M \in \Phi_2(M)$ mapping to $\phi$, put
\begin{eqnarray}
T_{\phi,\ellip}(G) =\{ \tau = (M,\pi_M,r):\,\ \pi_M \in \Pi_{\phi_M}, r \in R_{\reg}(\pi_M)    \}.
\end{eqnarray}
\begin{lemma}
$T_{\phi,\ellip}(G)$ is non-empty only when $\phi \in \Phi_{\ellip}(G)$.
\end{lemma}
\begin{proof}
By corollary 7.4.7, the action of any element $w^0 \in W^0_{\phi}$ on $\mathcal{I}_P(\pi_M)$ given by the normalized intertwining operator
\[
R_P(w^0,\widetilde{\pi}_M,\phi_M)
\]
is trivial. It follows that $W^0_{\phi} \subset W(\pi_M)^0$, while as we have already noted in section 3.4, for our current setting where $G$ is a unitary group, we have $W_{\phi} = W(\pi_M)$. Hence we have a surjection:
\begin{eqnarray}
R_{\phi} = W_{\phi}/W^0_{\phi} \rightarrow R(\pi_M) = W(\pi_M)/W^0(\pi_M). 
\end{eqnarray}
In particular, if $R_{\reg}(\pi_M)$ is non-empty, then so is $R_{\phi,\reg}$, in which case $\phi$ must lie in $\Phi_{\ellip}(G)$ (and that in fact $R_{\phi,\reg}$ is a singleton in this case, and so (7.5.5) is actually a bijection, and also a bijection on the set of regular elements).
\end{proof}
We will establish in the next subsection that $R_{\phi} \cong R(\pi_M)$, so that $T_{\phi,\ellip}(G)$ is indeed non-empty when $\phi \in \Phi_{\ellip}(G)$.

\bigskip

Thus assume that $\phi \in \Phi_{\ellip}(G) \smallsetminus \Phi_2(G)$. From our induction hypothesis, we have a pairing:
\begin{eqnarray}
\langle \cdot , \cdot \rangle : \mathcal{S}_{\phi_M} \times \Pi_{\phi_M} \rightarrow \{ \pm 1\}
\end{eqnarray}
which is perfect since $F$ is non-archimedean. On the other hand, since $\phi \in \Phi_{\ellip}(G)$, the elliptic subset $\mathcal{S}_{\phi,\ellip}$ is a torsor under $\mathcal{S}_{\phi_M}$. And we then have the extension of the pairing to:
\begin{eqnarray}
\langle \cdot, \cdot, \rangle: \mathcal{S}_{\phi,\ellip} \times T_{\phi,\ellip}
\end{eqnarray}
given by, for $\tau=(M,\pi_M,r)$:
\[
\langle x,\tau \rangle = \langle x,\widetilde{\pi}_M \rangle
\]
with $\langle x,\widetilde{\pi}_M \rangle$ being given by the extension of the pairing (7.5.6) to the $\mathcal{S}_{\phi_M}$-torsor $\mathcal{S}_{\phi,\ellip}$, as in section 3.4. 

We then have:
\begin{proposition}
For $\phi \in \Phi_{\ellip}(G) \smallsetminus \Phi_2(G)$ and $x \in \mathcal{S}_{\phi,\ellip}$, we have
\begin{eqnarray}
c_{\phi,x}(\tau) = \left \{ \begin{array}{c} \langle x,\tau \rangle, \mbox{  if } \tau \in T_{\phi,\ellip}(G) \\ 0 \mbox{  otherwise.} \end{array} \right.
\end{eqnarray}
\end{proposition}
\begin{proof}
This follows immediately on comparing (7.5.3) with the local intertwining relation.
\end{proof}

We can now state the elliptic orthogonality relations, which form the extra input we need to establish the local classification in the next subsections (in addition to the global inputs from section 6):
\begin{proposition} ({\it c.f.} Corollary 6.5.2 of \cite{A1})
Suppose that 
\[
y_i=(\phi_i,x_i), \,\ i=1,2
\]
are two pairs in $\mathcal{Y}_{\ellip}(G)$, such that the associated linear form (7.5.3) is defined. Then we have
\begin{eqnarray}
\sum_{\tau \in T_{\ellip}(G)} b(\tau) c_{\phi_1,x_1}(\tau) \overline{c_{\phi_2,x_2}(\tau)} = \left \{ \begin{array}{c} |\mathcal{S}_{\phi_1}|, \mbox{  if } y_1=y_2 \\ 0 \mbox{  otherwise.} \end{array} \right.
\end{eqnarray} 
\end{proposition}
\noindent Here in (7.5.9), the coefficient $b(\tau)$ for $\tau \in T_{\ellip}(G)$ is given by:
\[
b(\tau )= |d(\tau)| \cdot |R(\tau)|
\]
with $R(\tau):=R(\pi_M)$ if $\tau=(M,\pi_M,r)$.

The proof of proposition 7.5.3 are exactly the same as that in \cite{A1}, so we will not repeat the details (in {\it loc. cit.} it is obtained as corollary of Proposition 6.5.1) We content to mention that as a consequence of the local trace formulas for $G$ and the twisted group $\widetilde{G}_{E/F}(N)$, we have spectral expansions of the elliptic inner product on $G$, and respectively on $\widetilde{G}_{E/F}(N)$ (the local trace formula in the twisted case has been established by Waldspurger \cite{W6}). One then obtains (7.5.9) as a consequence of combining the spectral expansions with the endoscopic expansions of the elliptic inner products on $G$ and $\widetilde{G}_{E/F}(N)$, {\it c.f.} the proof of Proposition 6.5.1 of \cite{A1}. In our case this is the local and more elementary version of the discussion in section 5.6. (We remark that, similar to \cite{A1} we are taking for granted the stabilization of the elliptic inner product on the twisted group $\widetilde{G}_{E/F}(N)$, which is implicit in our assumption on the stabilization of the twisted trace formula.)

\bigskip

From proposition 7.5.3 we can deduce:
\begin{proposition} 

\noindent (a) Suppose that $(\phi,x) \in \mathcal{Y}_{\ellip}(G)$ such that $\phi \in \Phi_2(G)$, and such that the associated linear form (7.5.3) is defined. Then we have
\begin{eqnarray}
f^{\prime}(\phi^{\prime}) = \sum_{\Pi_2(G)} c_{\phi,x}(\pi) f_G(\pi), \,\ f \in \mathcal{H}_{\cusp}(G)
\end{eqnarray}
i.e. $c_{\phi,x}(\tau)=0$ for $\tau \in T_{\ellip}(G) \smallsetminus \Pi_2(G)$.

\bigskip
\noindent (b) Suppose that 
\[
y_i = (\phi_i,x_i), \,\ i=1,2
\]
are two pairs in $\mathcal{Y}_{\ellip}(G)$ that both satisfy the conditions of (a). Then
\begin{eqnarray}
\sum_{\pi \in \Pi_2(G)} c_{\phi_1,x_1}(\pi) \overline{c_{\phi_2,x_2}(\pi)} = \left \{ \begin{array}{c} |\mathcal{S}_{\phi_1}|, \mbox{  if } y_1=y_2 \\ 0 \mbox{  otherwise.} \end{array} \right.
\end{eqnarray}
\end{proposition}
\begin{proof}
This is the same as proof of lemma 6.5.3 of \cite{A1}. For part (a), suppose that $\tau_1 \in T_{\ellip}(G) \smallsetminus \Pi_2(G)$. Then $\tau_1 \in T_{\phi_1,\ellip}(G)$ for a unique $\phi_1 \in \Phi_{\bdd}(G) - \Phi_2(G)$, which by lemma 7.5.1, has to be in $\Phi_{\ellip}(G) - \Phi_2(G)$. Now by inversion of the formula (7.5.8), we have for $\tau \in T_{\ellip}(G)$:
\begin{eqnarray}
\frac{1}{|\mathcal{S}_{\phi_1,\ellip}|} \sum_{x_1 \in \mathcal{S}_{\phi_1,\ellip}} \langle x_1,\tau_1 \rangle \,\ \overline{c_{\phi_1,x_1}(\tau)} = \left \{ \begin{array}{c} 1, \mbox{  if } \tau=\tau_1 \\ 0 \mbox{  otherwise.} \end{array} \right.
\end{eqnarray}
Hence we have for any $x_1 \in \mathcal{S}_{\phi_1,\ellip}$:
\begin{eqnarray*}
 b(\tau_1) c_{\phi,x}(\tau_1) &=& \sum_{\tau \in T_{\ellip}(G)} \frac{1}{|\mathcal{S}_{\phi_1,\ellip}|} \sum_{x_1 \in \mathcal{S}_{\phi_1,\ellip}} \langle x_1 , \tau\rangle \,\ b(\tau) c_{\phi,x}(\tau) \overline{c_{\phi_1,x_1}(\tau)} \\
&=& \frac{1}{|\mathcal{S}_{\phi_1,\ellip}|} \sum_{x_1 \in \mathcal{S}_{\phi_1,\ellip}}  \langle x_1 , \tau\rangle \sum_{\tau \in T_{\ellip}(G)}   b(\tau) c_{\phi,x}(\tau) \overline{c_{\phi_1,x_1}(\tau)}. 
\end{eqnarray*}
The last inner sum vanishes by (7.5.9), since $\phi_1 \neq \phi$ (as $\phi \in \Phi_2(G)$). This gives the vanishing of $c_{\phi,x}(\tau_1)$. 

Part (b) follows on applying part (a) to (7.5.9), on noting that $b(\tau)=1$ for $\tau =\pi \in \Pi_2(G)$.
\end{proof}

\subsection{Local packets for non square-integrable parameters}

With the proof of the local intertwining relation in section 7.4, we can complete the proof of the local theorems for $\phi \in \Phi_{\bdd}(G) \smallsetminus \Phi_2(G)$. 

First, for part (a) of theorem 3.2.1, the existence of the stable linear form 
\[
f \mapsto f^G(\phi), \,\ f \in \mathcal{H}(G)
\]
for $\phi \in \Phi_{\bdd}(G) \smallsetminus \Phi_2(G)$ satisfying part (a) of theorem 3.2.1 is of course already implicit in our discussions in the previous subsection, and as we have seen it follows from descent and our local induction hypothesis. For exactly the same reason we also see that the linear form $f^G(\phi)$ depends only on $G$ and not on the $L$-embedding $\xi : \leftexp{L}{G} \hookrightarrow \leftexp{L}{G_{E/F}(N)}$ (note that part (a) of theorem 3.2.1 includes the case where $G$ is composite, but the descent argument is the same).

In section 3.4 we constructed the packet $\Pi_{\phi}$ from the irreducible constituents of $\mathcal{I}_P(\pi_M)$ with $\pi_M$ ranges over $\Pi_{\phi_M}$. In a moment we will see that the packet $\Pi_{\phi}$ is multiplicity free. Then by proposition 3.4.4, the packet $\Pi_{\phi}$ satisfies the endoscopic character relations, i.e. part (b) of theorem 3.2.1 is valid.

We also observe that in part (a) of theorem 2.5.1 we also have the assertion that if $\phi$ is a spherical parameter, then $\langle \cdot, \pi \rangle =1$ for the unique spherical $\pi \in \Pi_{\phi}$. This follows from the intertwining relation. Indeed for $\phi$ spherical we have $M=M_0$, and $\pi_{M_0}$ is a spherical unitary character of $M_0(F)$. By the local intertwining relation, the assertion is equivalent to the that the action of $R_{P_0}(w,\widetilde{\pi}_{M_0},\phi_{M_0})$ on $\pi$ is trivial for any $w \in W(M_0)$. But this is already known from Shahidi's results \cite{S,S2}.

It thus remains to verifiy that the packet $\Pi_{\phi}$ is multiplicity free, as representation of $\mathcal{S}_{\phi} \times G(F)$. We already know from induction hypothesis that the packet $\Pi_{\phi_M}$ is multiplicity free, so we must show that the the constituents of $\mathcal{I}_P(\pi_M)$ is multiplicity free as representation of $R_{\phi} \times G(F)$, for $\pi_M \in \Pi_{\phi_M}$.  

Recall that the decomposition of $\mathcal{I}_P(\pi_M)$ as a representation of $R(\pi_M) \times G(F)$ is multiplicity free, by the theory of representation theoretic $R$-group. Hence it suffices to show that $R_{\phi} \cong R(\pi_M)$. We have already seen that (7.5.5) is a surjection, and so we need to show that it is a bijection.

For this we first consider the case that $\phi \in \Phi_{\ellip}(G)$. Then (7.5.3) and the local intertwining relation implies the surjection
\begin{eqnarray}
R_{\phi} \rightarrow R(\pi_M)
\end{eqnarray}
restricts to a surjection
\begin{eqnarray}
R_{\phi,\reg} \rightarrow R_{\reg}(\pi_M)
\end{eqnarray}
and hence (7.6.2) is bijective ($R_{\phi,\reg}$ is either empty or singleton). 

Now let $\phi \in \Phi_{\bdd}(G) \smallsetminus \Phi_2(G)$ be general. Then $R_{\phi}$ is the disjoint union of $R_{\phi_L,\reg}$, with $L$ ranges over all the Levi subgroups $L$ of $G$ containing $M$ such that $\phi$ is the image of $\phi_L \in \Phi_{\ellip}(L)$. Similarly $R(\pi_M)$ is the disjoint union over all such Levi subgroups $L$ of $G$ containing $M$, of the sets $R^L_{\reg}(\pi_M)$; here we have denoted by $R^L(\pi_M)$ the representation theoretic $R$-group of $\pi_M$ with respect to $L$, and $R^L_{\reg}(\pi_M)$ its subset of regular elements ({\it c.f.} the discussion on page 530-531 of \cite{A11}). Thus by (7.6.2) applied to $\phi_L$ for all such $L$ we see that (7.6.1) must be a bijection.

With this we conclude the proof of all the local theorems for $\phi \in \Phi_{\bdd}(G) \smallsetminus \Phi_2(G)$.

\subsection{Local packets for square-integrable composite parameters}

We continue to assume $F$ is non-archimedean. In this subsection we construct the packets associated to generic square-integrable composite parameters $\phi \in \Phi_2(G)$. Thus we have
\[
\phi^N = \xi_* \phi = \phi_1^{N_1} \oplus \cdots \oplus \phi_r^{N_r} 
\]
with $\phi_i^{N_i} \in \widetilde{\Phi}_{\simp}(N_i)$, and $r > 1$. 

As before we globalize the data $(G,\phi)$, $\{\phi_i^{N_i}\}$ and $\phi^N$ to the global data:
\[
(\dot{G},\dot{\phi}),  \{\dot{\phi}_i^{N_i}\}, \dot{\phi}^N
\]
such that the conditions of proposition 7.3.1 are satisfied. We then consdered the family of global parameters:
\begin{eqnarray}
\dot{\widetilde{\mathcal{F}}} = \dot{\widetilde{\mathcal{F}}}(\dot{\phi}_1^{N_1},\cdots,\dot{\phi}_r^{N_r})=\{l_1 \dot{\phi}_1^{N_1}        \boxplus \cdots \boxplus l_r \dot{\phi}_r^{N_r}     , \,\ l_i \geq 0  \} 
\end{eqnarray}
and similarly the family $\dot{\widetilde{\mathcal{F}}}(\dot{G})$. As noted before, the conditions of proposition 7.3.1 implies that Assumption 6.4.1 are satisfied, and in particular we are in a position to apply the global results of section 6.4.

The first thing is the construction of the stable linear form:
\begin{eqnarray}
f \mapsto f^G(\phi), \,\ f \in \mathcal{H}(G).
\end{eqnarray}
i.e. we need to show that the linear form $\widetilde{f}_N(\phi^N)$ on $\widetilde{\mathcal{H}}(N)$ descends to the stable linear form (7.7.2) on $G=(G,\xi)$, i.e.
\begin{eqnarray}
\widetilde{f}^G(\phi) = \widetilde{f}_N(\phi^N), \,\ \widetilde{f} \in \widetilde{\mathcal{H}}(N).
\end{eqnarray} 
In fact we need to show the existence of (7.7.2) for $G^* = (G^*,\xi^*) \in \widetilde{\mathcal{E}}_{\ellip}(N)$, not just the simple ones. In the next proposition, we allow $G$ to be an element of $\widetilde{\mathcal{E}}_{\ellip}(N)$. 

\begin{proposition}
Suppose that $G \in \widetilde{\mathcal{E}}_{\ellip}(N)$, $\phi \in \Phi_2(G)$ such that $\phi^N$ is not a simple parameter. Then the stable linear form (7.7.2) exists. Furthermore, if $G$ is composite:
\[
G = G_O \times G_S, \,\ G_O,G_S \in \widetilde{\mathcal{E}}_{\simp}(N_i)
\]
\[
\phi = \phi_O \times \phi_S
\] 
\noindent then we have
\begin{eqnarray}
f^G(\phi)=f^{G_O}_O(\phi_O) f^{G_S}_S(\phi_S), \,\ f=f_O \times f_S.
\end{eqnarray}
\end{proposition}
\begin{proof}
We apply the global results of section 6.4 to the global datum $(\dot{G},\dot{\phi})$ and $\{\dot{\phi}_i^{N_i}\},\dot{\phi}^N$ above. By lemma 6.4.2, the global linear form 
\[
\dot{\widetilde{f}}_N(\dot{\phi}^N), \,\ \dot{\widetilde{f}} \in \dot{\widetilde{\mathcal{H}}}(N)
\]
descends to a stable linear form $\dot{f}^{\dot{G}}(\dot{\phi})$ on $\dot{G}=(\dot{G},\dot{\xi})$, i.e. 
\begin{eqnarray}
\dot{\widetilde{f}}^{\dot{G}}(\dot{\phi}) = \dot{\widetilde{f}}_N(\dot{\phi}^N), \,\ \dot{\widetilde{f}} \in \dot{\widetilde{\mathcal{H}}}(N)
\end{eqnarray}
and such that the global analogue of (7.7.4) is satisfied when $\dot{G}=\dot{G}_O \times \dot{G}_S$ is composite. We may choose $\dot{\widetilde{f}}$ and $\dot{f}$ to be decomposable. For $v \notin S_{\infty}(u)$, the parameter $\dot{\phi}_v$ is spherical, and hence $\dot{\phi}_v$ is not square-integrable, so the corresponding local assertion for $\dot{\phi}_v$ is known as we have seen in section 7.6. For $v \in S_{\infty}$, the assertions hold for the archimedean parameter $\dot{\phi}_v$, by the results of Mezo \cite{Me} and Shelstad \cite{Sh3} as discussed before. Hence we may cancel the contributions for places $v \neq u$ in (2.7.5) to obtain the corresponding assertion for the place $u$ and hence the corresponding local assertion holds for the parameter $\dot{\phi}_u=\phi$.
\end{proof}

In particular, we may assume from now on that $G \in \widetilde{\mathcal{E}}_{\simp}(N)$, and similarly in the global setup for $\dot{G}$. Note that at this point we have not shown that the stable linear form of proposition 7.7.1 depends on $G$ {\it only as a group and not as endoscopic datum in} $\widetilde{\mathcal{E}}_{\simp}(N)$. We will show this in section 7.9.

As before in the global setup we have the stabilization:
\begin{eqnarray}
I^{\dot{G}}_{\disc,\dot{\phi}^N}(\dot{f}) = \sum_{\dot{G}^{\prime} \in \dot{\widetilde{\mathcal{E}}}_{\ellip}(N) } \iota(\dot{G},\dot{G}^{\prime}) \widehat{S}^{\dot{G}^{\prime}}_{\disc,\dot{\phi}^N}(\dot{f}^{\prime}), \,\ \dot{f} \in \mathcal{H}(\dot{G}).
\end{eqnarray}

\noindent On the left hand side of (7.7.6), there is no contribution from the proper Levi subgroups of $\dot{G}$ to $I^{\dot{G}}_{\disc,\dot{\phi}^N}$, as follolws from our global induction hypothesis and the fact that $\dot{\phi}^N \in \dot{\widetilde{\Phi}}_{\ellip}(N)$ (we have already noted this on several occasions in section 4 and 5). Hence we have:
\begin{eqnarray}
I^{\dot{G}}_{\disc,\dot{\phi}^N}(\dot{f}) &=&  \tr R^{\dot{G}}_{\disc,\dot{\phi}^N}(\dot{f}) \\
&=& \sum_{\dot{\pi}} n(\dot{\pi}) \dot{f}_{\dot{G}}(\dot{\pi})               \nonumber
\end{eqnarray}
with $\dot{\pi}$ runs over the irreducible unitary representations of $\dot{G}(\dot{\mathbf{A}})$, and $n(\dot{\pi})$ is the multiplicity for $\dot{\pi}$ to occur in $ R^{\dot{G}}_{\disc,\dot{\phi}^N}$; in particular $n(\dot{\pi})$ {\it is a non-negative integer}.

On the other hand, by proposition 6.4.6, we have the validity of the stable multiplicity formula for the stable distributions $S_{\disc,\dot{\phi}^N}^{\dot{G}^{\prime}}$: 
\begin{eqnarray}
S^{\dot{G}^{\prime}}_{\disc,\dot{\phi}^N} = \sum_{\dot{\phi}^{\prime} \rightarrow \dot{\phi}^N} \frac{1}{|\mathcal{S}_{\dot{\phi}^{\prime}}|} \dot{f}^{\prime}(\dot{\phi}^{\prime}).
\end{eqnarray}

\noindent We then have
\begin{proposition}
We have
\begin{eqnarray}
\sum_{\dot{\pi}} n(\dot{\pi}) \dot{f}_{\dot{G}}(\dot{\pi}) = \frac{1}{|\mathcal{S}_{\dot{\phi}}|} \sum_{\dot{x} \in \mathcal{S}_{\dot{\phi}}} \dot{f}^{\prime}(\dot{\phi}^{\prime}).
\end{eqnarray}
\end{proposition}
\begin{proof}
We just need to apply (7.7.8) to the right hand side of (7.7.6):
\begin{eqnarray}
& & \sum_{\dot{G}^{\prime} \in \dot{\widetilde{\mathcal{E}}}_{\ellip}(N) } \iota(\dot{G},\dot{G}^{\prime})  \widehat{S}^{\dot{G}^{\prime}}_{\disc,\dot{\phi}^N}(\dot{f}^{\prime}) \\
&=& \sum_{\dot{G}^{\prime} \in \dot{\widetilde{\mathcal{E}}}_{\ellip}(N) }   \sum_{\dot{\phi}^{\prime} \rightarrow \dot{\phi}^N}    \iota(\dot{G},\dot{G}^{\prime})  \frac{1}{|\mathcal{S}_{\dot{\phi}^{\prime}}|} \dot{f}^{\prime}(\dot{\phi}^{\prime}). \nonumber
\end{eqnarray}
As an elementary case of the discussion in section 5.6, we have
\[
 \iota(\dot{G},\dot{G}^{\prime})  |\mathcal{S}_{\dot{\phi}^{\prime}}|^{-1} = |\mathcal{S}_{\dot{\phi}}|^{-1}
\]
and the double sum 
\[
\sum_{\dot{G}^{\prime}} \sum_{\dot{\phi}^{\prime}}
\]
in (7.7.10) can be indexed by a single sum 
\[
\sum_{\dot{x} \in \mathcal{S}_{\dot{\phi}}}
\]
with $(\dot{G}^{\prime},\dot{\phi}^{\prime}) \leftrightarrow (\dot{\phi},\dot{x})$. This gives (7.7.9).
\end{proof}

We now apply proposition 7.7.2 with a decomposable $\dot{f}$:
\[
\dot{f} = \dot{f}_u \cdot \dot{f}_{\infty} \cdot  \dot{f}^{\infty,u}
\]
with $\dot{f}_u = f \in \mathcal{H}_{\cusp}(G)$. For the places in $ S_{\infty}$, we have the endoscopic character identity:
\begin{eqnarray}
\dot{f}^{\prime}_{\infty}(\dot{\phi}^{\prime}_{\infty}) = \sum_{\dot{\pi}_{\infty} \in \Pi_{\dot{\phi}_{\infty}}} \langle \dot{x}_{\infty},\dot{\pi}_{\infty} \rangle \dot{f}_{\infty,\dot{G}_{\infty}}(\dot{\pi}_{\infty}), \,\ \dot{f}_{\infty} \in \mathcal{H}(\dot{G}_{\infty})
\end{eqnarray}
by the results of Shelstad; while for $v \notin S_{\infty}(u)$, the endoscopic character identity for the spherical parameter $\dot{\phi}_v$ is already known, as seen in section 7.6 (of course if $v$ splits in $\dot{E}$ then the assertion is trivial):
\begin{eqnarray}
\dot{f}^{\prime}_v(\dot{\phi}^{\prime}_v) = \sum_{\dot{\pi}_v \in \Pi_{\dot{\phi}_v}} \langle \dot{x}_v, \dot{\pi}_v \rangle \dot{f}_{v,\dot{G}_v}(\dot{\pi}_v), \,\ \dot{f}_v \in \mathcal{H}(\dot{G}_v).
\end{eqnarray}
Finally at the place $u$, we can apply part (a) of proposition 7.5.4 (proposition 7.7.1 shows that the hypothesis for proposition 7.5.4 is satisfied for $\phi = \dot{\phi}_u$): for $f \in \mathcal{H}_{\cusp}(G)$ we have:
\begin{eqnarray}
f^{\prime}(\phi^{\prime})=f^{\prime}(\dot{\phi}^{\prime}_u) = \sum_{\pi \in \Pi_2(G)} c_{\phi,x}(\pi) f_G(\pi), \,\ f \in \mathcal{H}_{\cusp}(G).
\end{eqnarray}
Hence combining (7.7.9) and (7.7.11) - (7.7.13) we obtain
\begin{eqnarray}
& & \sum_{\dot{\pi}} n(\dot{\pi}) \dot{f}_{\dot{G}}(\dot{\pi}) \\
&=& \frac{1}{|\mathcal{S}_{\dot{\phi}}|} \sum_{\dot{\pi}_{\infty}} \sum_{\dot{\pi}^{\infty,u}  } \sum_{\pi \in \Pi_2(G)} \sum_{\dot{x} \in \mathcal{S}_{\dot{\phi}}} c_{\phi,x}(\pi) \cdot \langle  \dot{x}_{\infty} , \dot{\pi}_{\infty} \rangle \cdot \langle \dot{x}^{\infty,u},\dot{\pi}^{\infty,u} \rangle \dot{f}_{\dot{G}}(\dot{\pi}) \nonumber
\end{eqnarray}
here on the right hand side of (7.7.14) $\dot{\pi}_{\infty}$ ranges over the elements of the packet $\Pi_{\dot{\phi}_{\infty}}$, and similarly $\dot{\pi}^{\infty,u}$ ranges over elements of the packet
\begin{eqnarray*}
& & \Pi_{\dot{\phi}^{\infty,u}} = \bigotimes_{v \notin S_{\infty}(u)} \Pi_{\dot{\phi}_v} \\
&:=& \Big\{ \dot{\pi}^{\infty,u} = \bigotimes_{v \notin S_{\infty}(u)}  \dot{\pi}_v, \,\ \dot{\pi}_v \in \Pi_{\dot{\phi}_v}, \,\ \langle \cdot , \dot{\pi}_v \rangle =1 \mbox{ for almost all } v  \Big\}.
\end{eqnarray*}

Recall that by part (i) of proposition 7.3.1 that $\mathcal{S}_{\dot{\phi}} \cong \mathcal{S}_{\phi}$. Fix a character $\xi \in \widehat{\mathcal{S}}_{\dot{\phi}} \cong \widehat{\mathcal{S}}_{\phi}$. By part (iii)(a) of proposition 7.3.1, we can pick a $\dot{\pi}_{\infty} \in \Pi_{\dot{\phi}_{\infty}}$ such that
\[
\xi(\dot{x})^{-1} = \langle \dot{x}_{\infty}, \dot{\pi}_{\infty} \rangle, \,\ \dot{x} \in \mathcal{S}_{\dot{\phi}}
\]
and for $\dot{\pi}^{\infty,u}$ we pick 
\[
\dot{\pi}^{\infty,u}(1) = \bigotimes_{v \notin S_{\infty}(u)} \dot{\pi}_v(1)
\]
with $\dot{\pi}_v(1) \in \Pi_{\dot{\phi}_v}$ corresponding to the trivial character of $\mathcal{S}_{\dot{\phi}_v}$. 

\noindent We define, for $\pi \in \Pi_2(G)$, the non-negative integer:
\[
n_{\phi}(\xi,\pi) := n(\pi \otimes \dot{\pi}_{\infty} \otimes \dot{\pi}^{\infty,u}).
\]

\begin{proposition}
The expression
\begin{eqnarray}
\frac{1}{|\mathcal{S}_{\phi}|} \sum_{x \in \mathcal{S}_{\phi}} c_{\phi,x}(\pi) \xi(x)^{-1}
\end{eqnarray}
is a non-negative integer, given by the integer $n_{\phi}(\xi,\pi)$. In particular the integer $n_{\phi}(\xi,\pi)$ defined above depends only on $\xi,\phi, \pi$.
\end{proposition}
\begin{proof}
With our choice of $\dot{\pi}_{\infty}$ and $\dot{\pi}^{\infty,u}$, by linear independence of characters of representations of $\dot{G}(\dot{\mathbf{A}}^u)$, it follows from (7.7.14) that
\[
\sum_{\pi \in \Pi_2(G)} n_{\phi}(\xi,\pi)  f_G(\pi) = \frac{1}{|\mathcal{S}_{\phi}|} \sum_{x \in \mathcal{S}_{\phi}} \sum_{\pi \in \Pi_2(G)} c_{\phi,x} (\pi) \xi(x)^{-1} f_G(\pi), \,\ f \in \mathcal{H}_{\cusp}(G).
\]
Then since characters of representation in $\Pi_2(G)$ are linear independent on $\mathcal{H}_{\cusp}(G)$, it follows that (7.7.15) holds, as required.
\end{proof}

\begin{proposition}
For $\phi,\phi^{*} \in \Phi_2(G)$ such that $\phi^N,\phi^{*,N} \notin \widetilde{\Phi}_{\simp}(N)$, we have for $\xi \in \widehat{\mathcal{S}}_{\phi}$, $\xi^* \in \widehat{\mathcal{S}}_{\phi^{*}}$:
\begin{eqnarray}
\sum_{\pi \in \Pi_2(G)} n_{\phi} (\xi,\pi) n_{\phi^{*}}(\xi^{*},\pi) = \left \{ \begin{array}{c} 1, \mbox{  if } (\phi,\xi) = (\phi^{*},\xi^{*}) \\ 0 \mbox{  otherwise.} \end{array} \right.
\end{eqnarray}
\end{proposition}
\begin{proof}
Compute:
\begin{eqnarray*}
& & \sum_{\pi \in \Pi_2(G)} n_{\phi} (\xi,\pi) n_{\phi^{*}}(\xi^{*},\pi) \\
&=& \sum_{\pi \in \Pi_2(G)} n_{\phi}(\xi,\pi) \overline{n_{\phi^*}(\xi^{*},\pi)}\\
&=& \frac{1}{|\mathcal{S}_{\phi}| \cdot |\mathcal{S}_{\phi^{*}}|} \sum_{x \in \mathcal{S}_{\phi}} \sum_{x^{*} \in \mathcal{S}_{\phi^{*}}} \xi(x)^{-1} \xi^{*}(x^{*}) \sum_{\pi \in \Pi_2(G)} c_{\phi,x}(\pi)\overline{ c_{\phi^{*},x^{*}}(\pi)}.
\end{eqnarray*}
By the orthogonality relation (7.5.11), the last inner sum is non-zero unless $(\phi,x)=(\phi^{*},x^{*})$, in which case the inner sum is $|\mathcal{S}_{\phi}|$. The assertion then follows.
\end{proof}

From (7.7.16), and the fact that $n_{\phi}(\xi,\pi)$ is a non-negative integer, it follows that given $\xi \in \widehat{\mathcal{S}}_{\phi}$, there is an unique $\pi=\pi(\xi) \in \Pi_2(G)$ such that $n(\xi,\pi)=1$, with $n(\xi,\pi^*)=0$ for $\pi^* \ncong \pi$, and also that the assignment
\[
\xi \mapsto \pi(\xi)
\]
is injective. Define the packet
\begin{eqnarray}
\Pi_{\phi} = \{ \pi(\xi) \in \Pi_2(G), \,\ \xi \in \widehat{\mathcal{S}}_{\phi}  \}
\end{eqnarray}
and for $\pi = \pi(\xi) \in \Pi_{\phi}$, define the character $\langle \cdot, \pi \rangle \in \widehat{\mathcal{S}}_{\phi}$:
\[
\langle \cdot , \pi \rangle = \xi(\cdot)
\]
This gives the perfect pairing between $\Pi_{\phi}$ and $\mathcal{S}_{\phi}$. It also folllows from (7.7.16) above that $\Pi_{\phi}$ and $\Pi_{\phi^*}$ are disjoint for distinct elements $\phi$ and $\phi^*$ in $\Phi_2^{\simp}(G)$.  

We have
\begin{proposition}
For $x \in \mathcal{S}_{\phi}$, we have
\begin{eqnarray}
f^{\prime}(\phi^{\prime}) = \sum_{\pi \in \Pi_{\phi}} \langle x,\pi \rangle f_G(\pi), \,\ f \in \mathcal{H}_{\cusp}(G).
\end{eqnarray}
\end{proposition}
\begin{proof}
By construction we have
\begin{eqnarray}
\frac{1}{|\mathcal{S}_{\phi}|} \sum_{x \in \mathcal{S}_{\phi}} c_{\phi,x}(\pi) \xi(x)^{-1} =  \left \{ \begin{array}{c} 1, \mbox{  if } \pi = \pi(\xi) \\ 0 \mbox{  otherwise.} \end{array} \right.
\end{eqnarray}
Hence inverting (7.7.19) we have
\begin{eqnarray}
c_{\phi,x}(\pi) =  \left \{ \begin{array}{c} \xi(x) = \langle x,\pi \rangle, \mbox{  if } \pi = \pi(\xi) \\ 0 \mbox{  otherwise.} \end{array} \right.
\end{eqnarray}
Applying (7.7.20) to (7.7.13), we obtain the proposition.
\end{proof}

\subsection{Local packets for simple parameters}
We now tackle the packets associated to simple parameters of $G$ in this section, which we eventually show to be singleton. For the purpose of carrying out the argument, it will be convenient to instead work with a temporary substitute of the set of simple parameters $\Phi_{\simp}(G)$, in a way similar to section 6.1.

We first recall some setup as in section 6.1 of \cite{A1}. First for our $G \in \widetilde{\mathcal{E}}_{\simp}(N)$ we have the linear isomorphism
\begin{eqnarray}
\mathcal{I}_{\cusp}(G) & \stackrel{\sim}{\rightarrow} & \bigoplus_{G^{\prime} \in \mathcal{E}_{\ellip}(G)} \mathcal{S}_{\cusp}(G^{\prime})^{\Out_G(G^{\prime})} \\
 f_G & \rightarrow &  \bigoplus_{G^{\prime}} f^{G^{\prime}}\nonumber
\end{eqnarray}
($\mathcal{S}_{\cusp}(G^{\prime})$ is the image of $\mathcal{H}_{\cusp}(G^{\prime})$ in $\mathcal{S}(G^{\prime})$). The twisted analogue of (7.8.1) is:
\begin{eqnarray}
\widetilde{\mathcal{I}}_{\cusp}(N) &\stackrel{\sim}{\rightarrow} & \bigoplus_{G \in \widetilde{\mathcal{E}}_{\ellip}(N)} \mathcal{S}_{\cusp}(G)^{\widetilde{\Out}_N(G)}  \\
\widetilde{f}_N & \rightarrow & \bigoplus_G \widetilde{f}^G \nonumber.
\end{eqnarray}

The linear isomorphism (7.8.1) is given by Proposition 3.5 of \cite{A11} for general $G$. The twisted case (7.8.2) in our current setting is given in \cite{Moe}.

Thus given any $\phi^N \in \widetilde{\Phi}_{\ellip}(N)$, we have a decomposition of the linear form:
\begin{eqnarray}
\widetilde{f}_N(\phi^N) = \sum_{G^* \in \widetilde{\mathcal{E}}_{\ellip}(N)} \widetilde{f}^{G^*}(\phi^*), \,\ \widetilde{f} \in \widetilde{\mathcal{H}}_{\cusp}(N)
\end{eqnarray}
here 
\[
f \mapsto f^{G^*}(\phi^*), \,\ f \in \mathcal{H}(G^*)
\]
is a stable linear form on $\mathcal{H}_{\cusp}(G^*)$. At this point $\phi^*$ is only a symbol standing for the pair $(G^*,\phi^N)$ in (7.8.3). As in section 6.5 of \cite{A1}, we introduce the following terminology:
 
\begin{definition}
\end{definition}
\noindent We say that $\phi^N$ is a {\it cuspidal lift}, if there is a $G^{\#}=(G^{\#},\xi^{\#}) \in \widetilde{\mathcal{E}}_{\ellip}(N)$, such that 
\[
\widetilde{f}_N(\phi^N) = \widetilde{f}^{G^{\#}}(\phi^{\#})
\]
i.e. there is only one term in the sum (7.8.3) corresponding to $G^{\#}$. 

Thus by proposition 7.7.1 if $\phi^N \in \widetilde{\Phi}_{\ellip}(N) \smallsetminus \widetilde{\Phi}_{\simp}(N)$ is such that $\phi^N = \xi^{\#}_* \phi^{\#}$, for $\phi^{\#} \in \Phi_{2}(G^{\#})$ with $G^{\#}=(G^{\#},\xi^{\#}) \in \widetilde{\mathcal{E}}_{\ellip}(N)$, then $\phi^N$ is a cuspidal lift.

With this preparation, we define, for $G=(G,\xi) \in \widetilde{\mathcal{E}}_{\simp}(N)$, the following temporary set of simple parameters for $G$:
\begin{eqnarray}
& & \Phi_{\simp}(G) \\
&=&\{\phi= (G,\phi^N), \,\ f^G(\phi) \neq 0 \mbox{ for some } f \in \mathcal{H}_{\cusp}(G)   \}. \nonumber
\end{eqnarray}
We shall show in the next subsection that this coincides with the usual definition in tems of Langlands parameters. We also simply define $\mathcal{S}_{\phi}$ to be trivial for $\phi \in \Phi_{\simp}(G)$.

Note that at this point we do not know that if $\phi=(G,\phi^N) \in \Phi_{\simp}(G)$, then $\phi^N$ is a cuspidal lift. Thus we define
\[
\Phi^c_{\simp}(G) \subset \Phi_{\simp}(G)
\]
to be the subset of simple parameters in $\Phi_{\simp}(G)$ that are cuspidal lifts. Similarly, we redefine for the temporary purpose of this subsection:
\begin{eqnarray}
\Phi_2(G) &:=& \Phi_{\simp}(G) \coprod \Phi_2^{\simp}(G)\\
\Phi_2^c(G) &:=& \Phi_{\simp}^c(G) \coprod \Phi_2^{\simp}(G) \nonumber 
\end{eqnarray}
where $\Phi_2^{\simp}(G)$ has the same meaning as in the last subsection.

\begin{definition}
\end{definition}
\noindent Define $\Pi_{\simp}(G)$ to be the subset of $\Pi_2(G)$ given by the complement of $\Pi_{\phi}$ over $\phi \in \Phi_2^{\simp}(G)$. 

\bigskip

With our setup and the results established in the previous subsections, the arguments in section 6.7 of \cite{A1} for constructing the packets associated to simple parameters apply to our case {\it verbatim}. So we will be brief.
 
Using the spectral interpretation of the spaces $\mathcal{I}_{\cusp}(G)$ and $\mathcal{S}_{\cusp}(G)$, we define $\mathcal{I}_{\simp}(G)$ to be the subspace of $\mathcal{I}_{\cusp}(G)$ consisting of those $f_G$ such that $f_G(\tau)=0$ for $\tau \in T_{\ellip}(G) - \Pi_{\simp}(G)$. Similarly define $\mathcal{S}_{\simp}(G)$ to be the subspace of $\mathcal{S}_{\cusp}(G)$ consisting of those $f^G$ such that $f^G(\phi)=0$ for any $\phi \in \Phi_2^{\simp}(G)$. As in \cite{A1}, we have: 

\begin{proposition} (Lemma 6.7.1 of \cite{A1})
Under the isomorphism (7.8.1), the subspace $\mathcal{I}_{\simp}(G) \subset \mathcal{I}_{\cusp}(G)$ is mapped isomorphically to the subspace $\mathcal{S}_{\simp}(G) \subset \mathcal{S}_{\cusp}(G)$. 
\end{proposition}
\noindent Here $\mathcal{S}_{\cusp}(G)$ is identified with the subspace of the right hand side of (7.8.1) consisting of elements $\bigoplus_{G^{\prime}} f^{G^{\prime}}$ such that $f^{G^{\prime}}=0$ for $G^{\prime} \neq G$. In particular for $f_G \in \mathcal{I}_{\simp}(G)$, we have $f^{G^{\prime}}=0$ for $G^{\prime} \neq G$.

\bigskip

\noindent In order to complete the local classification, the first step is to associate to $\pi \in \Pi_{\simp}(G)$ a parameter $\phi \in \Phi_{\simp}^c(G)$. Thus let $\pi \in \Pi_{\simp}(G)$. Then by definition, if $f_{\pi}$ is a pseudo-coefficient of $\pi$, we have $f_{\pi,G} \in \mathcal{I}_{\simp}(G)$. As in lemma 7.2.3, we globalize the data $(F,G,\pi)$ to $(\dot{F},\dot{G},\dot{\pi})$, such that $\dot{\pi}_u \cong \pi$ for a nonarchimedean place $u$ of $\dot{F}$. In the proof of lemma 7.2.3, we have shown that $\dot{\pi}$ belongs to $R^{\dot{G}}_{\disc,\dot{\phi}^N}$ for $\dot{\phi}^N \in \dot{\widetilde{\Phi}}(N)$, such that $\phi^N:=\dot{\phi}^N_u$ lies in $\widetilde{\Phi}_{\ellip}(N)$. 

Furthermore, as seen from the argument that led to (7.2.7), we have the following: for $\widetilde{f} \in \widetilde{\mathcal{H}}(N)$, the value $\widetilde{f}(\phi^N)$ depends only on the transfer $\widetilde{f}^G$, and the value is non-zero when $\widetilde{f}^{G}=f^G_{\pi}$ (as before $f_{\pi}$ is a pseudo-coefficient for $\pi$). In particular, the linear form
\[
\widetilde{f} \mapsto \widetilde{f}_N(\phi^N), \,\ \widetilde{f} \in \widetilde{\mathcal{H}}_{\cusp}(N)
\]     
depends only on the transfer $\widetilde{f}^G$, i.e. in (7.8.3) all the terms corresponding to $G^* \neq G$ vanish. Thus we have
\begin{eqnarray}
\widetilde{f}_N(\phi^N) = \widetilde{f}^G(\phi)
\end{eqnarray}
for a stable linear form $f \mapsto f^G(\phi)$ on $\mathcal{H}_{\cusp}(G)$.

From this it follows that $\phi^N \in \widetilde{\Phi}_{\simp}(N)$. Indeed, supppose not. Then $\phi^N \in \widetilde{\Phi}_{\ellip}^{\simp}(N)$, and hence $\phi^N \in \xi^{\#}_* \Phi_2(G^{\#})$, for a unique $G^{\#}=(G^{\#},\xi^{\#}) \in \widetilde{\mathcal{E}}_{\ellip}(N)$. Thus we have 
\[
\phi^N = \xi^{\#}_* \phi^{\#}, \,\ \phi^{\#} \in \Phi_2(G^{\#}).
\]
By proposition 7.7.1, we have
\[
\widetilde{f}_N(\phi^N) = \widetilde{f}^{G^{\#}}(\phi^{\#}), \,\ \widetilde{f} \in \widetilde{\mathcal{H}}_{\cusp}(N)
\]
for a stable linear form $f \mapsto f^{\#}(\phi^{\#})$ on $\mathcal{H}_{\cusp}(G^{\#})$. Thus by the uniqueness of the decomposition (7.8.3), we must have $G^{\#}=G$ and $\phi^{\#}=\phi$. Hence $\phi \in \Phi_2^{\simp}(G)$. In particular from proposition 7.7.5, we have
\begin{eqnarray}
f^G(\phi) = \sum_{\pi^* \in \Pi_{\phi}} f_G(\pi^*).
\end{eqnarray}
Now $\pi \in \Pi_{\simp}(G)$, so $\pi$ does not lie in $\Pi_{\phi}$. In particular if follows from (7.8.7) that
\begin{eqnarray}
f_{\pi}^G(\phi) =0.
\end{eqnarray}
But if $\widetilde{f} \in \widetilde{\mathcal{H}}_{\cusp}(N)$ is such that $\widetilde{f}^G=f_{\pi}^G$, then from (7.8.6) and (7.8.8) that
\[
\widetilde{f}_N(\phi^N) = \widetilde{f}^G(\phi) =f_{\pi}^G(\phi)=0
\]
contradicting what we have observed above. Thus we conclude that $\phi^N \in \widetilde{\Phi}_{\simp}(N)$, and by (7.8.6) again, $\phi^N$ is a cuspidal lift. Thus $\phi \in \Phi_{\simp}^c(G)$.  

Having shown that $\phi^N \in \widetilde{\Phi}_{\simp}(N)$, we have in particular that $\dot{\phi}^N \in \dot{\widetilde{\Phi}}_{\simp}(N)$, i.e. $\dot{\phi}^N$ is just a conjugate self-dual cuspidal automorphic representation on $G_{\dot{E}/\dot{F}}(\dot{\mathbf{A}})$. 

We can thus consider the global family of parameters
\[
\dot{\widetilde{\mathcal{F}}} = \dot{\widetilde{\mathcal{F}}}(\dot{\phi}^N)=\{l \dot{\phi}^N, \,\ l \geq 0 \}
\]
and the same argument as in the proof of proposition 7.3.1 shows that Assumption 6.4.1 is satisfied for the family $\dot{\widetilde{\mathcal{F}}}$, and hence we can apply the results of section 6.4 to $\dot{\widetilde{\mathcal{F}}}$. 

For this, recall that $\dot{\pi}$ occurs in $R^{\dot{G}}_{\disc,\dot{\phi}^N}$, and we have seen in the proof of lemma 7.2.3 (the equation after (7.2.5)) that
\[
\tr R^{\dot{G}}_{\disc,\dot{\phi}^N}(\dot{f}) = S^{\dot{G}}_{\disc,\dot{\phi}^N}(\dot{f})
\]
and both sides being non-zero, for a suitable choice of $\dot{f} \in \dot{\mathcal{H}}(\dot{G})$. The condition that $S^{\dot{G}}_{\disc,\dot{\phi}^N} \neq 0$ is exactly the (temporary) definition of the condition $\dot{\phi}^N \in \dot{\xi}_* \Phi_{\simp}(\dot{G})$ in section 6.1, and $(\dot{G},\dot{\phi}^N)$ defines the global simple parameter $\dot{\phi} \in \dot{\Phi}_{\simp}(\dot{G})$. In particular, by lemma 6.4.2, the global linear form 
\[
\dot{\widetilde{f}} \mapsto \dot{\widetilde{f}}_N(\dot{\phi}^N), \,\ \dot{\widetilde{f}} \in \dot{\widetilde{\mathcal{H}}}(N)
\]
descends to $\dot{G}$, i.e. there is a stable linear form
\[
\dot{f} \mapsto \dot{f}(\dot{\phi}), \,\ \dot{f} \in \dot{\mathcal{H}}(\dot{G})
\]
such that $\dot{\widetilde{f}}_N(\dot{\phi}^N)=\dot{\widetilde{f}}^{\dot{G}}(\dot{\phi})$ for $\dot{\widetilde{f}} \in \dot{\widetilde{\mathcal{H}}}(N)$. In particular by localizing at $u$, we see that the linear form
\[
\widetilde{f} \mapsto \widetilde{f}_N(\phi^N), \,\ \widetilde{f} \in \widetilde{\mathcal{H}}(N)
\] 
descends to $G$, i.e. there is a stable linear form 
\[
f \mapsto f^G(\phi), \,\ f \in \mathcal{H}(G)
\]
such that
\[
\widetilde{f}^G(\phi) = \widetilde{f}_N(\phi^N), \,\ \widetilde{f} \in \widetilde{\mathcal{H}}(N).
\]
The stable linear form $f^G(\phi)$ extends the one in (7.8.6) that is already defined for $f \in \mathcal{H}_{\cusp}(G)$. 

By proposition 6.4.7, the stable multiplicity formula is valid for $\dot{\phi}$, i.e. we have
\begin{eqnarray}
S_{\disc,\dot{\phi}^N}^{\dot{G}}(\dot{f}) = \dot{f}^{\dot{G}}(\dot{\phi})
\end{eqnarray}
As in proposition 7.7.2, the application of the stable multiplicity formula (7.8.9) gives:
\begin{eqnarray}
\sum_{\dot{\pi}^*} n(\dot{\pi}^*) \dot{f}_{\dot{G}}(\dot{\pi}^*) = \dot{f}^{\dot{G}}(\dot{\phi})
\end{eqnarray}
with $n(\dot{\pi}^{*})$ being the multiplicity of an irreducible unitary representation $\dot{\pi}^*$ of $\dot{G}(\dot{\mathbf{A}})$ to occur in $R^{\dot{G}}_{\disc,\dot{\phi}^N}$. Then we similarly define for $\pi^* \in \Pi_2(G)$:
\[
n_{\phi}(1,\pi^*) = n(\pi^* \otimes \dot{\pi}_{\infty} \otimes \dot{\pi}^{\infty,u})
\]
with $\dot{\pi}_{\infty} $ and $\dot{\pi}^{\infty,u}$ have similar meaning as in the last subsection (the ``1" appears simply because the group $\mathcal{S}_{\phi}$ is trivial). By the construction in the proof of lemma 7.2.3, we have $n_{\phi}(1,\pi)$ is a {\it positive} integer. 

By restricting $\dot{f}$ in (7.8.10) such that $\dot{f}_u \in \mathcal{H}_{\cusp}(G)$, we extract the following identity similar to proposition 7.7.3:
\begin{eqnarray}
c_{\phi,1}(\pi^*) = n_{\phi}(1,\pi^*)
\end{eqnarray}
with $c_{\phi,1}(\pi^*)$ being by definition the coefficients in the expansion:
\[
f^G(\phi) = \sum_{\pi^* \in \Pi_2(G)} c_{\phi,1}(\pi^*) f_G(\pi^*).
\]

\noindent Application of the orthogonality relation (7.5.11), combined with (7.8.11) gives the following: for any $\phi^{\#} \in \Phi_2^c(G)$, we have
\begin{eqnarray}
\sum_{\pi^* \in \Pi_2(G)} n_{\phi}(1,\pi^*) n_{\phi^{\#}}(1,\pi^*) =  \left \{ \begin{array}{c} 1, \mbox{  if } \phi = \phi^{\#} \\ 0 \mbox{  otherwise.} \end{array} \right.
\end{eqnarray}
Hence we have for any $\pi^* \in \Pi_2(G)$:
\begin{eqnarray}
n_{\phi}(1,\pi^*)=c_{\phi,1}(\pi^*)=  \left \{ \begin{array}{c} 1, \mbox{  if } \pi^{*} = \pi \\ 0 \mbox{  otherwise} \end{array} \right.
\end{eqnarray}
and the assignment
\begin{eqnarray}
\pi \mapsto \phi
\end{eqnarray}
defined by condition (7.8.13) gives a well-defined injection from $\Pi_{\simp}(G)$ to $\Phi^c_{\simp}(G)$. We also see that
\begin{eqnarray}
f^G(\phi) = f_G(\pi), \,\ f \in \mathcal{H}_{\cusp}(G).
\end{eqnarray}

To complete the argument, one needs to show:
\begin{proposition}
The map (7.8.14) gives a bijection between $\Pi_{\simp}(G)$ and $\Phi_{\simp}(G)$. In particular, $\Phi_{\simp}(G)=\Phi_{\simp}^c(G)$. Given $\phi \in \Phi_{\simp}(G)$, the associated representation $\pi_{\phi} \in \Pi_{\simp}(G)$ satisfies 
\begin{eqnarray}
f^G(\phi) = f_G(\pi_{\phi}), \,\ f \in \mathcal{H}_{\cusp}(G).
\end{eqnarray}
\end{proposition}
\begin{proof}
This is the same as proposition 6.7.2 of \cite{A1}.
\end{proof}

\noindent Given $\phi \in \Phi_{\simp}(G)$, we then define the packet $\Pi_{\phi}$ associated to $\phi$ to be the singleton $\{\pi_{\phi}\} \subset \Pi_{\simp}(G)$. These are disjoint for distinct $\phi$. Then we have:

\begin{proposition}
For $\phi \in \Phi_2(G)$, and $\phi^N=\xi_* \phi$ as before, we have for any $x \in \mathcal{S}_{\phi}$, the endoscopic character identity:
\begin{eqnarray}
f^{\prime}(\phi^{\prime}) = \sum_{\pi \in \Pi_{\phi}} \langle x,\pi \rangle f_G(\pi)
\end{eqnarray} 
where as usual $(G^{\prime},\phi^{\prime})\leftrightarrow (\phi,x)$.
\end{proposition}
\begin{proof}
Same as corollary 6.7.4 of \cite{A1}. We need to show that the identities (2.7.18), (2.8.15) established for $f \in \mathcal{H}_{\cusp}(G)$ remains valid for all $f \in \mathcal{H}(G)$. This entails using Theorem 6.1, 6.2 of \cite{A11}, and their expected twisted analogues (which we is implicit in our assumption on stabilization of the twisted trace formula).
\end{proof}

From the results of section 7.6 - 7.8, the following is clear:
\begin{corollary}
$\Pi_{\temp}(G)$ is the disjoint union of the packets $\Pi_{\phi}$ for $\phi \in \Phi_{\bdd}(G)$.
\end{corollary}

\noindent In order to complete the proof of all the local theorems for generic parameters, it remains to show that the temporary definition of $\Phi_{\simp}(G)$ given in (7.8.4) coincides with the usual one in terms of $L$-parameters, and to show that for $\phi \in \Phi_{\bdd}(G)$, the stable linear form $f \rightarrow f^G(\phi)$ depends on $G$ and not on the $L$-embedding $\xi: \leftexp{L}{G} \hookrightarrow \leftexp{L}{G_{E/F}(N)}$. We finish these in the next subsection.

\subsection{Resolution}

As before $F$ is non-archimedean, and $G=(G,\xi) \in \widetilde{\mathcal{E}}_{\simp}(N)$, with $\xi =\xi_{\chi}$ for $\chi \in \mathcal{Z}_E^{\kappa}$ ($\kappa= \pm 1$). We need to show that the temporary definition of the set of simple parameters $\Phi_{\simp}(G)$ coincides with the usual one in terms of the local Langlands group $L_F$. This is the content of the following:

\begin{proposition}
Let $\phi^N \in \widetilde{\Phi}_{\simp}(N)$ Then the following are equivalent:

\noindent (a) The parameter $\phi=(G,\phi^N)$ belongs to $\Phi_{\simp}(G)$ (in the sense of the definition in (7.8.4)).

\bigskip

\noindent (b) The Langlands parameter $\phi^N : L_F \rightarrow \leftexp{L}{G_{E/F}(N)}$ factors through the $L$-embedding $\xi: \leftexp{L}{G} \hookrightarrow \leftexp{L}{G_{E/F}(N)}$.

\bigskip

\noindent (c) The (Langlands-Shahidi) $L$-factor $L(s,\phi^N,\Asai^{(-1)^{N-1}\kappa})$ has a pole at $s=0$.
\end{proposition}
\begin{proof}
Condition (b) is equivalent to the statement that the Artin $L$-factor $L(s,\Asai^{(-1)^{N-1}\kappa}  \phi^N)$ has a pole at $s=0$. Thus the equivalence of (b) and (c) follows from the theorem of Henniart \cite{H2}, which gives the equality between the Artin and the Langlands-Shahidi $L$-factors:
\[
L(s,\Asai^{(-1)^{N-1}\kappa}  \phi^N) = L(s,\phi^N,\Asai^{(-1)^{N-1}\kappa}).
\]
Next we note that among the two $L$-factors $L(s,\phi^N,\Asai^+)$ and $L(s,\phi^N,\Asai^-)$, {\it exactly} one of them has a pole at $s=0$. Hence to prove the equivalences of (a) and (c), it only remains to prove that (a) implies (c).

To do this, we need to use the method of supplementary parameters as in section 6.3. Put $N_+=2N$, and consider the parameter
\[
\phi_+^{N_+} := \phi^N \oplus \phi^N \in \widetilde{\Phi}_{\ellip}(N_+).
\]
We choose the endoscopic datum $G_+ \in \widetilde{\mathcal{E}}_{\simp}(N_+)$ that is of the same parity as $G$ as a twisted endoscopic datum in $\widetilde{\mathcal{E}}_{\simp}(N)$ (recall that the parity of $G=(G,\xi)$ is given by the sign $(-1)^{N-1}\kappa$, and similarly for $G_+$). Thus for example, if $N$ is even, then the $L$-embedding $\xi_+=\xi_{+,\chi_+}$ in the datum for $G_+$ is given by a character $\chi_+ \in \mathcal{Z}_E^{\kappa}$; while if $N$ is odd, then $\chi_+ \in \mathcal{Z}_E^{-\kappa}$. Thus $\chi_+ \in \mathcal{Z}_E^{(-1)^N \cdot \kappa}$. 

With this choice we then have $\phi_+^{N_+} = \xi_{+,*} \phi_+$, for $\phi_+ \in \Phi_{\ellip}(G_+)$. Indeed we have $S_{\phi_+}=O(2,\mathbf{C})$. We also put $M_+=G_{E/F}(N)$, which is the Siegel Levi of $G_+$, and equip $M_+=(M_+,\xi_+) \in \widetilde{\mathcal{E}}(N_+)$ as a (twisted) endoscopic datum of $\widetilde{G}_{E/F}(N_+)$. Then we can take $\phi_{+,M_+} \in \Phi_2(M_+,\phi_+^{N_+})$ such that $\phi_+^{N_+}=\xi_{+,*} \phi_{+,M_+}$; in fact we have
\begin{eqnarray}
 \phi_{+,M_+} =    (\chi_+ \circ \det)^{-1} \otimes  \phi^N.
\end{eqnarray}
We denote by $\pi_{+,M_+}$ the irreducible admissible representation of $M_+(F) = G_{E/F}(F)$ that corresponds to $\phi_{+,M_+}$ under the usual local Langlands correspondence for $G_{E/F}(N)$ (in other words $\{\pi_{+,M_+}\}$ is the packet corresponding to $\phi_{+,M_+}$ ). 

\noindent With this setup, we claim that the induced representation $\mathcal{I}^{G_+}_{P_+}(\pi_{+,M_+})$ is reducible. Taking this for granted for the moment, we can complete the proof of proposition 7.9.1 as follows. The reducibility assertion is equivalent to the condition that the unnormalized intertwining operator $J_{w^{-1}P_+|P_+}(\pi_{+,M_+})$ being holomorphic, with $w$ being the unique non-trivial element of $W_{\phi_+}(G_+,M_+)$ (here the notation is as in section 3.4). Hence this is equivalent to the condition that the normalizing factor $r_{P_+}(w,\phi_{+,M_+})$ for $J_{w^{-1}P_+|P_+}(\pi_{+,M_+})$ being finite. From the expression of $r_{P_+}(w,\phi_{+,M_+})$ given by (3.3.44), this is equivalent to the condition that the Langlands-Shahidi $L$-factor $L(s,\phi_{+,M_+},\Asai^+)$ does not have a pole at $s=0$, i.e. that $L(s,\phi^N,\Asai^{(-1)^N \cdot \kappa})$ does not have a pole at $s=0$, by (7.9.1). This in turn is equivalent to the condition that $L(s,\phi^N,\Asai^{(-1)^{N-1} \kappa})$ has a pole at $s=0$.  
\end{proof}

\noindent It remains to show:
\begin{lemma}
In the setting above, the induced representation $\mathcal{I}^{G_+}_{P_+}(\pi_{+,M_+})$ is reducible.
\end{lemma}
\begin{proof}
In the setting above, we have $\phi=(G,\phi^N) \in \Phi_{\simp}(G)$ (again with $\Phi_{\simp}(G)$ being defined as in (7.8.4)). By proposition 7.8.4, the parameter $\phi$ corresponds to $\pi \in \Pi_{\simp}(G)$. We now apply global methods: as in section 7.8, we can as in the proof of lemma 7.2.3, globalize the data:
\[
(F,G,\phi,\phi^N,\pi)
\] 
to global data
\[
(\dot{F},\dot{G},\dot{\phi},\dot{\phi}^N,\dot{\pi})
\]
such that the conditions of proposition 7.3.1 holds. Similarly we globalize the data
\[
(G_+,\phi_+,\phi_+^{N_+},M_+)
\]
to
\[
(\dot{G}_+,\dot{\phi}_+,\dot{\phi}^{N_+},\dot{M}_+)
\]
compatibly with that of $(\dot{G},\dot{\phi},\dot{\phi}^N)$. Then as in section 7.8, we apply the results of section 6.4 to the global family
\[
\dot{\widetilde{\mathcal{F}}} = \dot{\widetilde{\mathcal{F}}}(\dot{\phi}^N)=\{l \dot{\phi}^N, \,\ l \geq 0 \}.
\]   
In particular, we apply proposition 6.4.7 concerning the vanishing of (6.4.26) which (when combined with (6.4.27) gives the global intertwining relation:
\begin{eqnarray}
\dot{f}_{\dot{G}_+}^{\prime}(\dot{\phi}_+,\dot{x}_+) = \dot{f}_{\dot{G}_+}(\dot{\phi}_+,\dot{x}_+), \,\ \dot{f} \in \dot{\mathcal{H}}(\dot{G}_+)
\end{eqnarray}
here $\dot{x}_+$ is the unique element of $\mathcal{S}_{\dot{\phi}_+,\ellip} \cong \mathcal{S}_{\phi_+,\ellip}$. Note that since we are dealing with parameters of size $N_+ >N$, (7.9.2) does not follow immediately from the local intertwining relation for parameters of size $N$ treated in section 7.4. In any case, our purpose is only to extract the reducibility of $\mathcal{I}_{P_+}(\pi_{+,M_+})$ from (7.9.2), and this can be done as follows. We choose decomposable $\dot{f} = \dot{f}_u \dot{f}^u$, with $\dot{f}_u=f$ being variable, and $\dot{f}^u$ is fixed such that $\dot{f}^u_{\dot{G}_+}(\dot{\phi}_+,\dot{x}_+) \neq 0$. Then (7.9.2) gives:
\begin{eqnarray}
c \cdot f^{\prime}_{G_+}(\phi_+,x_+) = f_{G_+}(\phi_+,x_+), \,\ f \in \mathcal{H}(G_+)
\end{eqnarray}  
for a non-zero constant $c$. On the other hand, under the endoscopic correspondence for parameters we have
\[
(G \times G,\phi \times \phi) \leftrightarrow (\phi_+,x_+)
\] 
hence the linear form $f^{\prime}_{G_+}(\phi_+,x_+)$ is given by the non-zero linear form:
\begin{eqnarray}
& & f^{\prime}_{G_+}(\phi_+,x_+) = f^{G \times G}(\phi \times \phi) \\
&=& f_1^G(\phi) \cdot f_2^G(\phi), \,\ \mbox{ if } f^{G \times G} = f_1^G \times f_2^G. \nonumber
\end{eqnarray}
On the other hand, we have
\begin{eqnarray}
& & f_{G_+}(\phi_+,x_+) \\
&=& \tr \big( R_{P_+}( w_{x_+} , \widetilde{\pi}_{+,M_+} , \phi_{+,M_+}) \mathcal{I}_{P_+}(\pi_{+,M_+},f) \big). \nonumber
\end{eqnarray}
If $\mathcal{I}_{P_+}(\pi_{+,M_+})$ were irreducible, then the action of $ R_{P_+}( w_{x_+} , \widetilde{\pi}_{+,M_+} , \phi_{+,M_+})$ on $\mathcal{I}_{P_+}(\pi_{+,M_+})$ would be given by a scalar, and hence for some non-zero constant $c^{\prime}$ we have:
\begin{eqnarray}
  f_{G_+}(\phi_+,x_+) &=& c^{\prime} \tr \mathcal{I}_{P_+}(\pi_{+,M_+},f) \\
&=& c^{\prime} f_{M_+}(\pi_{+,M_+}). \nonumber
\end{eqnarray}
However, (7.9.6) contradicts (7.9.4), since the linear form $f \mapsto f^{G \times G}(\phi \times \phi)$ on $\mathcal{H}(G_+)$ does not factor through the map $f \mapsto f_{M_+}$ (this is seen for example, by applying proposition 3.1.1(a). We conclude that $\mathcal{I}_{P_+}(\pi_{+,M_+})$ is reducible.

\end{proof}

\begin{rem}
\end{rem}
The relation between poles of local Asai $L$-function and reducibility of parabolic induction (with respect to the Siegel parabolic) is already considered in Goldberg \cite{G}.
\bigskip

With proposition 7.9.1 proved, we can thus use the original interpretation of $\Phi_{\simp}(G)$ in terms of Langlands parameters. In order to complete the induction step of the local argument for generic parameters, it only remains to show that the local linear form 
\[
f \mapsto f^G(\phi), \,\ f \in \mathcal{H}(G)
\]
for $\phi \in \Phi_{\bdd}(G)$ depends only on $G$ and not on the $L$-embedding $\xi:\leftexp{L}{G} \hookrightarrow \leftexp{L}{G_{E/F}(N)}$. More precisely, we fix representatives of the two elements of $\widetilde{\mathcal{E}}_{\simp}(N)$, given by 
\[
G=(U_{E/F}(N),\xi), \,\ G^{\vee} =(U_{E/F},\xi^{\vee})
\]
with $\xi = \xi_{\chi_+}$ for $\chi_+ \in \mathcal{Z}_E^+$, and similarly $\xi^{\vee}=\xi^{\vee}_{\chi_-}$ for $\chi_- \in \mathcal{Z}_E^-$. Given an $L$-parameter for $U_{E/F}(N)$:
\[
\phi:L_F \rightarrow \leftexp{L}{U_{E/F}(N)}
\]
we put $\phi^N :=\xi_* \phi$, and $\phi^{\vee,N}:=\xi^{\vee}_* \phi$. Then $\phi^N$ (resp. $\phi^{\vee,N}$) define stable linear forms:
\[
f \mapsto f^{G}(\phi), \,\ f \in \mathcal{H}(U_{E/F}(N))
\]
(resp.
\[
f  \mapsto f^{G^{\vee}}(\phi), \,\ f \in \mathcal{H}(U_{E/F}(N)) \,\ )
\]
such that
\[
\widetilde{f}^{G}(\phi) = \widetilde{f}_N(\phi^N), \,\ \widetilde{f} \in \widetilde{\mathcal{H}}(N)
\]
(resp.
\[
\widetilde{f}^{G^{\vee}}(\phi) = \widetilde{f}_N(\phi^{\vee,N}), \,\ \widetilde{f} \in \widetilde{\mathcal{H}}(N)\,\ )
\]
The stable linear forms $f^{G}(\phi)$, resp. $f^{G^{\vee}}(\phi)^{\vee}$ define packets $\Pi_{\phi} \subset \Pi_2(U_{E/F}(N))$, resp. $\Pi_{\phi^{\vee}} \subset \Pi_2(U_{E/F}(N))$ which might be different {\it a priori}. What we must show is that $f^{G}(\phi)=f^{G^{\vee}}(\phi)$ for $f \in \mathcal{H}(U_{E/F}(N))$; in other words, we need to show the following: if $\widetilde{f}_1,\widetilde{f}_2 \in \widetilde{\mathcal{H}}(N)$ are such that
\begin{eqnarray}
\widetilde{f}_1^{G} = \widetilde{f}_2^{G^{\vee}}
\end{eqnarray}
then we have
\begin{eqnarray}
\widetilde{f}_{1,N}(\phi^N) = \widetilde{f}_{2,N}(\phi^{\vee,N}).
\end{eqnarray}
We also note that (7.9.8) is already a consequence of the work of Mezo \cite{Me} and Shelstad \cite{Sh3} in the archimedean case.

By usual descent argument and the local induction hypothesis, it suffices to establish (7.9.8) in the case where $\phi \in \Phi_2(U_{E/F}(N))$. We then treat this case again by using global methods. Thus arguing as in corollary 7.2.7 we can globalize the data:
\[
(F,G,\phi^N)
\]
to global data:
\[
(\dot{F},\dot{G},\dot{\phi}^N)
\]
with $\dot{\phi}^N \in \dot{\widetilde{\Phi}}_{\ellip}(N)$ such that $\dot{\widetilde{\phi}}^N_u=\phi^N$ (with $u$ being the place of $\dot{F}$ such that $\dot{F}_u=F$ and $\dot{G}_u=G$). Similarly we globalize the data
\[
(F,G,\phi^{\vee,N})
\]
to global data
\[
(\dot{F},\dot{G}^{\vee},\dot{\phi}^{\vee,N}).
\]
In fact we only need to form the global datum $\dot{G}^{\vee}$, with $\dot{\phi}^{\vee,N}$ being determined as follows. Namely if
\[
\dot{G}=(U_{\dot{E}/\dot{F}}(N),\dot{\xi}_{\dot{\chi}_+}), \,\ \dot{G}^{\vee}=(U_{\dot{E}/\dot{F}}(N),\dot{\xi}^{\vee}_{\dot{\chi}_-})
\]
with $\dot{\chi}_+ \in \dot{\mathcal{Z}}^+_{\dot{E}}$ (resp. $\dot{\chi}_- \in \dot{\mathcal{Z}}^-_{\dot{E}}$) then we take
\begin{eqnarray}
& & \dot{\phi}^{\vee,N} = \dot{\eta}  \otimes\dot{\phi}^N, \\
& & \dot{\eta} = \dot{\chi}_- /\dot{\chi}_+. \nonumber
\end{eqnarray}
From (7.9.9) we have from definition
\begin{eqnarray}
R^{\dot{G}}_{\disc,\dot{\phi}^N}(\dot{f}) &=& R^{\dot{G}^{\vee}}_{\disc,\dot{\phi}^{\vee,N}}(\dot{f}), \,\  \\
S^{\dot{G}}_{\disc,\dot{\phi}^N}(\dot{f}) &=& S^{\dot{G}^{\vee}}_{\disc,\dot{\phi}^{\vee,N}}(\dot{f}), \,\ \dot{f} \in \dot{\mathcal{H}}(U_{\dot{E}/\dot{F}}(N)).
\end{eqnarray}

\noindent We now choose (decomposable) global test function $\dot{\widetilde{f}}_1 \in \dot{\mathcal{H}}(N)$ with $\widetilde{f}_1 = \dot{\widetilde{f}}_{1,u} \in \mathcal{H}(N)$ being the component at the place $u$ of $\dot{F}$. Then we have seen in the proof of lemma 7.2.3 that we can choose $\dot{\widetilde{f}}_1^u$ such that for some non-zero constant $c$, we have:
\begin{eqnarray}
c \cdot \widetilde{f}_1(\phi^N) = \widehat{S}^{\dot{G}}_{\disc,\dot{\phi}^N}(\widetilde{f}_1^{G} \cdot (\dot{\widetilde{f}}_1^u)^{\dot{G}}), \,\ \widetilde{f} \in \widetilde{\mathcal{H}}(N).
\end{eqnarray} 
Similarly we can choose $\dot{\widetilde{f}}_2^u$ such that the analogue of (7.9.12) holds for $\dot{G}^{\vee}$, i.e. there is a non-zero constant $c^{\prime}$ such that:
\begin{eqnarray}
c^{\prime} \cdot \widetilde{f}_2(\phi^{\vee,N}) = \widehat{S}^{\dot{G}^{\vee}}_{\disc,\dot{\phi}^{\vee,N}}(\widetilde{f}_2^{G^{\vee}} \cdot (\dot{\widetilde{f}}_2^u)^{\dot{G}^{\vee}}), \,\ \widetilde{f}_2 \in \widetilde{\mathcal{H}}(N).
\end{eqnarray} 
Hence by (7.9.11)-(7.9.13), we see that
\begin{eqnarray}
 \widetilde{f}_{1,N}(\phi^N) = \frac{c^{\prime}}{c} \cdot \widetilde{f}_{2,N}(\phi^{\vee,N})
\end{eqnarray}
whenever
\[
\widetilde{f}_1^G = \widetilde{f}_2^{G^{\vee}}
\]
holds. In other words, we have
\begin{eqnarray}
f^{G}(\phi) = \frac{c^{\prime}}{c} \cdot f^{G^{\vee}}(\phi), \,\ f \in \mathcal{H}(U_{E/F}(N)).
\end{eqnarray}

\noindent It thus suffices to show hat the constant $c^{\prime}/c$ in (7.9.15) is equal to one. But we know that 
\[
f^{G}(\phi) = \sum_{\pi \in \Pi_{\phi}} f_{U_{E/F}(N)}(\pi)
\]
and similarly
\[
f^{G^{\vee}}(\phi) = \sum_{\pi^{\prime} \in \Pi_{\phi^{\vee} }} f_{U_{E/F}(N)}(\pi^{\prime})
\]
so it follows that $c^{\prime}/c =1$ (and $\Pi_{\phi}=\Pi_{\phi^{\vee}}$), by linear independence of characters.

\bigskip

To complete the induction arguments for this section, it only remains to verify the global theorems for parameters in families $\widetilde{\mathcal{F}}$ of degree $N$. This is in fact already implicit in section 6.4, together with the results of this section.

First as in remark 6.4.9, the results of section 6.4 shows that the the definition of the set $\dot{\Phi}_{\simp}(\dot{G})$ of simple generic parameters of $\dot{G}$ coincides with the original definition given in section 2.4. In particular, the results of section 6.4 implies the seed theorem 2.4.2 is valid for the simple generic parameters $\dot{\widetilde{\mathcal{F}}}_{\gsimp}(N)$ of degree $N$.

For the second seed theorem 2.4.10, we simply note that for $\dot{G}=(\dot{G},\dot{\xi})$, if $\dot{\phi}^N \in \dot{\xi}_* \dot{\widetilde{\mathcal{F}}}_{\gsimp}(\dot{G})$, then the linear form 
\[
\dot{\widetilde{f}} \mapsto \dot{\widetilde{f}}_N(\dot{\phi}^N), \,\ \dot{\widetilde{f}} \in \dot{\widetilde{\mathcal{H}}}(N)
\]
descends under the Kottwitz-Shelstad transfer to a stable linear form on $\dot{G}$ (with respect to $\dot{\xi}$), according to the definition of the set $\dot{\widetilde{\mathcal{F}}}_{\gsimp}(\dot{G})$ given in section 6.1 (which we know is equivalent to the original definition in section 2.4 for parameters in $\dot{\widetilde{\mathcal{F}}}$). In particular for any place $v$ of $\dot{F}$ (which we may assume to be non-split in $\dot{E}$), the linear form 
\[
\dot{\widetilde{f}}_v \mapsto \dot{\widetilde{f}}_{v,N}(\dot{\phi}^N_v), \,\ \dot{\widetilde{f}}_v \in \dot{\widetilde{\mathcal{H}}}_v(N)
\]
descends under the Kottwitz-Shelstad transfer to $\dot{G}_v$ (with respect to $\dot{\xi}_v$). The local classification results established before in this section then implies that $\dot{\phi}_v \in \dot{\xi}_{v,*} \Phi(\dot{G}_v)$.

The validity of the stable multiplicity formula for parameters $\dot{\widetilde{\mathcal{F}}}(N)$ with respect to $\dot{G}^* \in \dot{\widetilde{\mathcal{E}}}_{\simp}(N)$ is a consequence of Proposition 5.7.4. and Propositions 6.4.4-6.4.7. With the stable multiplicity formula in hand, together with the local theorems for generic parameters of degree $N$ established in this section, the spectral multiplicity formula theorem 2.5.2 for parameters in the families $\dot{\widetilde{\mathcal{F}}}$ (which consists only of generic parameters) follows from lemma 5.7.6.

Theorem 5.2.1 (for our families $\dot{\widetilde{\mathcal{F}}}$ of generic global parameters) is then just follows from the local intertwining relation we established in this section (applied to each place $v$ of $\dot{F}$).

The only remaining global theorem is theorem 2.5.4(a) (recall that part (b) of theorem 2.5.4 really concerns non-generic parameters, {\it c.f.} remark 2.5.7). But this is just a consequence of Proposition 6.4.8.

\bigskip

We have thus completed the induction step for the global theorems for parameters in the families $\dot{\widetilde{\mathcal{F}}}$. With this also completes all the local theorems for generic parameters.

\section{\bf Nontempered representations}

In section seven we established all the local theorems for generic parameters. In this section we turn to the local theorems for non-generic parameters. As in section seven, we are going to establish all the local theorems for non-generic parameters for {\it all} degree $N$. Again the method is based on global inputs from section 6.  

The use of induction in section eight is similar to that given in the beginning of section seven. Thus we fix the integer $N$, and we assume that all the local theorems are valid for (non-generic) parameters of degree smaller than $N$. We also need to work with families $\dot{\widetilde{\mathcal{F}}}$ of global parameters to be introduced in section 8.3; all the global arguments are to be carried out within families $\dot{\widetilde{\mathcal{F}}}$, and we need to establish the global theorems for the families of parameters $\dot{\widetilde{\mathcal{F}}}$ of all degree in this section. Thus for the induction part of the argument we assume that all the global theorems are valid for global parameters in families $\dot{\widetilde{\mathcal{F}}}$ of degree less than $N$.

\subsection{Duality operator of Aubert-Schneider-Stuhler}

To establish the local theorems for general (non-generic) parameters, we have to construct the packets corresponding to non-generic parameters, and to verify the local intertwining relation for the packets. 

Given our results of section seven for the generic parameters, the constructions given in chapter 7 of \cite{A1} applies essentially without change for unitary groups. Thus we just give a summary of the construction, and mostly refer to \cite{A1} for proofs.

The construction is based again based on global method, combined with the duality operator of Aubert \cite{Au} and Schneider-Stuhler \cite{SS} (which exists when the local field $F$ is non-archimedean). 

Thus we first assume that $F$ is non-archimedean, and $G$ an arbitrary connected reductive group over $F$. The duality operator $D$ is an involution on the Grothendieck group $\mathcal{K}(G)$ of finite length $G(F)$-modules given by:
\[
D = \sum_{P \supset P_0} (-1)^{\dim A_{P_0}/A_P} i^G_P \circ r^G_P
\] 
with $r_P^G$ and $i_P^G$ being restriction and induction functors respectively (here $P_0$ is a fixed minimal parabloic subgroup of $G$ whose Levi component is noted as $M_0$). For such an object $\pi$ we denote by $[\pi]$ it image in $\mathcal{K}(G)$. If $\pi$ is irreducible then it is a theorem of Aubert \cite{Au} and Schneider-Stuhler \cite{SS} that $D \pi$ is again irreducible up to sign:
\begin{eqnarray}
D \pi = \beta(\pi) [\widehat{\pi}]
\end{eqnarray}
with $\widehat{\pi}$ irreducible. Here $\beta(\pi)$ is the sign
\begin{eqnarray}
\beta(\pi) = (-1)^{\dim A_{M_0}/A_{M_{\pi}}}
\end{eqnarray}
with $M_{\pi}$ being Levi subgroup of $G$, defined in terms of the Bernstein component to which $\pi$ belongs. The definition of the duality operator can be made similarly for a twisted group $\widetilde{G}$, and the twisted version of the results for the duality operator remain valid  \cite{A13}.

The complexification $\mathcal{K}(G)_{\mathbf{C}}$ can be identified with the space of finite invariant distributions on $G(F)$ (in the sense of being a finite sum of irreducible characters), given by the characters of representations. Thus $D$ extends as an operator on the space of invariant distributions. The operator $D$ is compatible with endoscopy \cite{Hi,A13}, in the following sense: if $G^{\prime}$ is a (twisted) endoscopic datum of $G$, then the dual of the Langlands-Shelstad-Kottwitz transfer gives a mapping:
\[
S^{\prime} \mapsto S^{\prime}_G
\]
with $S^{\prime}$ a stable distribution on $G^{\prime}(F)$ and $S^{\prime}_G$ an invariant distribution on $G(F)$. Then we have
\begin{eqnarray}
(D^{\prime} S^{\prime})_G = \alpha(G,G^{\prime}) D S^{\prime}_G
\end{eqnarray}
with $\alpha(G,G^{\prime})$ being the sign
\[
\alpha(G,G^{\prime}) = (-1)^{\dim A_{M_0}/A_{M_0^{\prime}}};
\]
here $D^{\prime}$ is the duality operator on $G^{\prime}$. The operator $D$ preserves the space of stable distributions.

The case of the twisted group $\widetilde{G}_{E/F}(N)$ again plays an important role (as before $E/F$ is a quadratic extension). Given $\phi^N \in \widetilde{\Phi}_{\bdd}(N)$ it can be identified as an irreducible conjugate self-dual tempered representation $\pi_{\phi^N}$ of $G_{E/F}(F)=GL_N(E)$. Then we have the following result \cite{A13} on the twisted version of the duality operator $\widetilde{D}$ on $\widetilde{G}_{E/F}(N)$, which is the twisted version of Aubert and Schneider-Stuhler's result \cite{Au,SS} on the Zelevinsky's conjecture:
\begin{eqnarray}
\widetilde{D} \phi^N = \widetilde{\beta}(\phi^N) [\psi^N]
\end{eqnarray}  
here $\psi^N \in \widetilde{\Psi}(N)$ is the parameter given by the dual $\widehat{\phi}^N$ of $\phi^N$, namely 
\[
\phi^N: L_F = W_F \times \SU_2 \rightarrow \leftexp{L}{G_{E/F}(N)}
\] 
and
\[
\psi^N : L_F \times \SU_2 = W_F \times \SU_2 \times \SU_2 \rightarrow \leftexp{L}{G_{E/F}(N)} 
\]
(with the $\SU_2$ factors being as ordered above) are related as
\begin{eqnarray}
& & \psi^N(w,u_1,u_2) =\widehat{\phi}^N(w,u_1,u_2)\\
&=& \phi^N(w,u_2), \,\ w \in W_F, u_1,u_2 \in \SU_2 \nonumber
\end{eqnarray}
(more symmetrically, we can write
\[
\widehat{\phi}^N(w,u_1,u_2) = \phi^N(w,u_2,u_1)
\]
if $\phi^N$ is identified as a parameter on $L_F \times \SU_2 = W_F \times \SU_2 \times \SU_2$ that is trivial on the second $\SU_2$-factor). The sign $\widetilde{\beta}(\phi^N)$ is given as on p.388 of \cite{A1}. Again we identify $\psi^N$ with the conjugate self-dual representation $\pi_{\psi^N}$ of $G_{E/F}(N)(F)$, which is obtained from the Langlands quotient of the standard representation associated to the parameter $\phi_{\psi^N} \in \widetilde{\Phi}(N)$.

\begin{rem}
\end{rem}
The Zelevinsky conjecture (proved by Aubert \cite{Au} and Schneider-Stuhler \cite{SS}) is the statement:
\begin{eqnarray*}
D \phi^N = \beta(\phi^N) \psi^N
\end{eqnarray*}
where $D,\phi^N,\psi^N$ and $\beta(\phi^N)$ are in the context of the untwisted group $G_{E/F}(F)(F)=\GL_N(E)$.

\subsection{Local parameters}

We return to the setting where $F$ is a general local field. We need to establish theorem 3.2.1 for non-generic parameters. Thus let $G = (G,\xi) \in \widetilde{\mathcal{E}}_{\ellip}(N)$. Then we have

\begin{proposition}
Let $\psi \in \Psi(G)$, and as before denote $\psi^N = \xi_* \psi \in \widetilde{\Psi}(N)$. Then there exists a (unique) stable linear form
\begin{eqnarray}
f \mapsto f^G(\psi), \,\ f \in \mathcal{H}(G)
\end{eqnarray}
that depends only on $G$ (and not on the $L$-embedding $\xi: \leftexp{L}{G} \hookrightarrow \leftexp{L}{G_{E/F}(N)}$), such that 
\begin{eqnarray}
\widetilde{f}^G(\psi)=\widetilde{f}_N(\psi^N),\,\ \widetilde{f} \in \widetilde{\mathcal{H}}(N).
\end{eqnarray}
Furthermore, if $G=G_O \times G_S$ is composite, and $\psi = \psi_O \times \psi_S$, then
\begin{eqnarray}
 f^G(\phi) = f_1^{G_O}(\psi_O) \cdot f^{G_S}(\psi_S)
\end{eqnarray}
for 
\begin{eqnarray*}
f^G=f^{G_O} \times f^{G_S}. 
\end{eqnarray*}
\end{proposition}
\begin{proof}
This is the same as lemma 2.2.2 of \cite{A1}. By usual descent argument it suffices to treat the case $\psi \in \Psi_2(G)$. Suppose first that $G \in \widetilde{\mathcal{E}}_{\simp}(N)$. Thus let $\psi \in \Psi(G)$, and put $\psi^N =\xi_* \psi \in \widetilde{\Psi}(N)$. We can then write ({\it c.f. loc. cit.}) the twisted character $\widetilde{f}_N(\psi^N)$ as a (finite) linear combination of twisted standard characters:
\begin{eqnarray*}
\widetilde{f}_N(\psi^N) = \sum_{\phi^N \in \widetilde{\Phi}(N)} \widetilde{n}(\psi^N,\phi^N) \cdot \widetilde{f}_N(\phi^N), \,\ \widetilde{f} \in \widetilde{\mathcal{H}}(N)
\end{eqnarray*}  
for coefficients $n(\psi^N,\phi^N)$, with the property that if $\widetilde{n}(\psi^N,\phi^N) \neq 0$ then $\phi^N$ has the same infinitesimal character (resp. cuspidal support) as that of $\phi_{\psi^N}$ for $F$ archimedean (resp. non-archimedean). Hence by considering the infinitesimal character (resp. cuspidal support) of $\phi^N$, we have $\phi^N \in \xi_* \Phi(G)$ for $\phi^N$ such that $\widetilde{n}(\psi^N,\phi^N) \neq 0$. Hence if we put
\begin{eqnarray*}
& & f^G(\psi) := \sum_{\phi \in \Phi(G)} n(\psi,\phi) f^G(\phi) \\
& & n(\psi,\phi):= \widetilde{n}(\psi^N,\phi^N)
\end{eqnarray*}
(note that $n(\psi,\phi)$ does not depend on the $L$-embedding $\xi$ which justifies the notation). Then the linear form $f^G(\psi)$ is stable, and satisfies (8.2.2). Furthermore, by the results of section 7 we know that the linear form $f^G(\phi)$ depends only on $G$ and not on the $L$-embedding $\xi$. Hence the same is true for $f^G(\psi)$.

\noindent For the case where $G = G_O \times G_S \in \widetilde{\mathcal{E}}_{\ellip}(N)$ and $\psi = \psi_O \times \psi_S$ is composite, we define the linear form $f^G(\psi)$ by (8.2.3), and then one needs to show that (8.2.2) is satisfied. The argument is the same as proof of lemma 2.2.2 of \cite{A1}. The key point being that, if $\widetilde{\pi}_{\psi_O^{N_O}}$ and $\widetilde{\pi}_{\psi_S^{N_S}}$ are the Langlands quotient of the (twisted) standard representations $\widetilde{\rho}_{\psi_O^{N_O}}$ and $\widetilde{\rho}_{\psi_S^{N_S}}$ of $\widetilde{G}^+_{E/F}(N)(F)$ respectively, then 
\[
\pi_{\psi^N}:=\mathcal{I}_{P_{N_O,N_S}} \big(\pi_{\psi_O^{N_O}} \boxtimes \pi_{\psi_S^{N_S}}\big) 
\]
is already irreducible, by Bernstein's theorem \cite{Be} (here $P_{N_O,N_S}$ is the standard parabolic subgroup of $G_{E/F}(N)$ of block size $(N_O,N_S)$, and $\psi^N:= \psi_O^{N_O} \oplus \psi_S^{N_S}$), and hence $\widetilde{\pi}_{\psi^N}$ is the Langlands quotient of $\widetilde{\rho}_{\psi^N}$.

\end{proof}

In particular, we may assume from now on that $G \in \widetilde{\mathcal{E}}_{\simp}(N)$. And in fact since we now know that the stable linear form depends only on $G$ and not on the $L$-embedding $\xi$, we may simply consider the case where $G=(G,\xi)$ with $\xi$ being the ``standard base change" $L$-embedding (i.e. $\xi=\xi_{\chi_{\triv}}$ with $\chi_{\triv} \in \mathcal{Z}_E^+$ being the trivial character).

Thus let $\psi \in \Psi(G)$, and $s \in \overline{S}_{\psi}$, with $x$ being the image of $s$ in $\mathcal{S}_{\psi}$; as usual we have the endoscopic correspondence for parameters:
\[
(G^{\prime},\psi^{\prime}) \leftrightarrow (\psi,s).
\]
The descent argument of lemma 5.3.1, which applies equally well in the local case, shows that the linear form 
\[
f \rightarrow f^{\prime}_G(\psi,s) :=f^{\prime}(\psi^{\prime}), \,\ f \in \mathcal{H}(G)
\]
depends only on the image $x$ of $s$ in $\mathcal{S}_{\psi}$. Recall the element
\[
s_{\psi} = \psi \Big(1, \begin{pmatrix}   -1 &  0 \\    0  &    -1    \\ \end{pmatrix}     \Big).
\]
We can then form an expansion as in \cite{A1}:
\begin{eqnarray}
f^{\prime}(\psi^{\prime}) = \sum_{\sigma \in \Sigma_{\psi} } \langle  s_{\psi} x , \sigma \rangle \cdot f_G(\sigma)
\end{eqnarray}
here $\Sigma_{\psi}$ is another name for the set of charcters $\widehat{\mathcal{S}}_{\psi}$ of $\mathcal{S}_{\psi}$; for $\sigma \in \Sigma_{\psi}$, the corresponding character of $\mathcal{S}_{\psi}$ is noted as $\langle \cdot,\sigma \rangle$, and $f_G(\sigma)$ is an invariant linear form on $\mathcal{H}(G)$. The set $\Sigma_{\psi}$ is temporarily called the packet associated to $\psi$, and we identity $\sigma \in \Sigma_{\psi}$ with the associated linear form $f_G(\sigma)$. Thus we have to show that the linear form $f_G(\sigma)$ is a sum of irreducible characters of unitary representations with non-negative integral coefficients. 

To establish this by global method, we need partial information supplied the non-generic parameters obtained from the dual of generic parameters. Thus let $F$ be non-archimedean. As in (8.1.5), for $\phi \in \Phi_{\bdd}(G)$, we denote $\psi=\widehat{\phi} \in \Psi(G)$ the dual of $\phi$:
\begin{eqnarray}
& & \psi(w,u_1,u_2) = \widehat{\phi}(w,u_2,u_1) \\
&=& \phi(w,u_2), \,\ w \in W_F,u_1,u_2 \in \SU_2. \nonumber
\end{eqnarray}
In particular $\psi^N = \widehat{\phi}^N$. We can identify $\mathcal{S}_{\phi} \cong \mathcal{S}_{\psi}$, and the corresponding identification $\Pi_{\phi} \leftrightarrow \Sigma_{\psi}$. Given $\pi \in \Pi_{\phi}$ we denote by $\sigma_{\pi}$ the corresponding element of $\Sigma_{\psi}$. 

As a consequnece of (8.1.4) and the twisted version of (8.1.3), we have
\begin{eqnarray}
D \phi =\beta(\phi) \psi
\end{eqnarray}
as stable linear forms on $G$, {\it c.f.} equation (7.1.6) of \cite{A1}. Here $\phi$ and $\psi$ are identified as their assocaited stable linear forms $f^G(\phi),f^G(\psi)$, and $\beta(\phi)$ is the sign associated to $\phi$ as on p.388 of \cite{A1}. 

From (8.1.3), (8.2.6) and the endoscopic character relation for $\phi$ one has:
\begin{lemma} (lemma 7.1.1 of \cite{A1})
With $\psi = \widehat{\phi}$ as above we have
\begin{eqnarray}
\langle s_{\psi},\pi \rangle \sigma_{\pi} = \beta(\pi) \beta(\phi) [\widehat{\pi}].
\end{eqnarray}
\end{lemma}

\noindent Thus we define
\begin{eqnarray}
\Pi^G_{\phi}=\{ \pi \in \Pi_{\phi},\,\  \sigma_{\pi} = [\widehat{\pi}] \}
\end{eqnarray}
which is equivalent by (8.2.7) to the condition $\langle s_{\psi},\pi \rangle=\beta(\pi)\beta(\phi)$. And put
\begin{eqnarray}
\Pi_{\psi}^G = \{ [\widehat{\pi}], \,\ \pi \in \Pi_{\phi}^G    \} \subset \Sigma_{\psi}.
\end{eqnarray}
Eventually we show via global methods that $\Pi_{\psi}^G=\Pi_{\psi}$; among other things, this shows that the representations $\widehat{\pi}$ for $\pi \in \Pi_{\phi}$ are unitary. Before one can apply global methods, one needs independent partial information about $\Pi_{\psi}^G$ for particular types of $\psi$; see proposition 8.2.4 below.

\noindent The situation for the local intertwining relation is similar. First we have spectral interpretations of the spaces $\mathcal{I}(G)$ and $\mathcal{S}(G)$ (for $\mathcal{S}(G)$ this is of course the consequence of the local classification for generic parameters in section 7). Then one has the adjoint action of $D$ on these spaces:
\begin{eqnarray*}
& & (D f_G)(\pi)= f_G(D \pi) \\
& & (D f^G)(\phi) = f^G(D \phi).
\end{eqnarray*}

\noindent Thus let $\psi \in \Psi(G)$ of the form $\psi = \widehat{\phi}$ for $\phi \in \Phi_{\bdd}(G)$. Assume that $\psi \notin \Psi_2(G)$. Thus we have a proper Levi $M$ of $G$, and $\psi_M \in \Psi_2(M,\psi)$ and a corresponding $\phi_M \in \Psi_2(M,\phi)$ such that $\psi_M = \widehat{\phi}_M$. From the local induction hypothesis one has $\Pi_{\psi_M} = \{\widehat{\pi}_M,\,\ \pi_M \in \Pi_{\phi_M}  \}$. Then from the property of duality operator, one has the following:  ({\it c.f.} p.392 of \cite{A1}) for every $u \in \mathcal{N}_{\psi}=\mathcal{N}_{\phi}$, we have
\begin{eqnarray*}
f^{\prime}_G(\psi,u) = (Df)_G^{\prime}(\phi,u), \,\ f \in \mathcal{H}(G)
\end{eqnarray*}

\noindent Hence by the local intertwining relation for the generic parameter $\phi$ (which was established in section 7) we have:
\begin{lemma}
For $\psi=\widehat{\phi}$ as above, the local intertwining relation for $\psi$:
\[
f_G^{\prime}(\psi,s_{\psi}x) = f_G(\psi,x), \,\ f \in \mathcal{H}(G)
\]
is equivalent to:
\begin{eqnarray}
f_G(\psi,u) = \beta(\phi) (D f)_G(\phi,s_{\psi}u), \,\ u \in \mathcal{N}_{\phi}, \,\  f \in \mathcal{H}(G).
\end{eqnarray}
\end{lemma}
More precisely, put for $u \in \mathcal{N}_{\phi} =\mathcal{N}_{\psi}$:
\begin{eqnarray*}
& & f_G(\phi,u,\pi_M) = \langle \widetilde{u}, \widetilde{\pi}_M \rangle \tr \big( R_P(w_u, \widetilde{\pi}_M,\phi_M) \mathcal{I}_P(\pi_M,f) \big) \\
& & f_G(\psi,u,\sigma_M) = \langle \widetilde{u}, \widetilde{\sigma}_M \rangle \tr \big( R_P(w_u, \widetilde{\sigma}_M,\psi_M) \mathcal{I}_P(\sigma_M,f) \big)
\end{eqnarray*}
for $\pi_M \in \Pi_{\phi_M}$ and $\sigma_M \in \Pi_{\psi_M}$, then (8.2.10) is equivalent to the condition:
\begin{eqnarray}
f_G(\psi,u,\widehat{\pi}_M) = \beta(\phi) (D f)_G(\phi,s_{\psi} u,\pi_M)
\end{eqnarray} 
for all $\pi_M \in \Pi_{\phi_M}$. Analogous to (8.2.8) we put
\begin{eqnarray*}
& & \Pi_{\phi_M}^G=\{\pi_M \in \Pi_{\phi_M}, \mbox{ such that }(8.2.11) \mbox{ is valid.} \} \\
& & \Pi_{\psi_M}^G=\{\widehat{\pi}_M, \,\ \pi_M \in \Pi^G_{\phi_M}  \}. 
\end{eqnarray*}

\noindent Then the extra input we need to establish the local theorems in the non-generic case is contained in the following proposition. Before stating that we need a terminology. We say that a parameter $\psi$ is {\it tamely ramified linear parameter}, if $E/F$ is tamely ramified, and the restriction of $\psi^N$ to $W_E$ decomposes as a sum of tamely ramified one-dimensional characters of $W_E$. Similar definition applies to $\phi$.
\begin{proposition} \cite{A13}
Suppose that $\psi = \widehat{\phi}$ is a tamely ramified linear parameter. Then

\bigskip
\noindent (a) The image of $\Pi_{\psi}^G$ in $\widehat{\mathcal{S}}_{\psi}$ generates $\widehat{\mathcal{S}}_{\psi}$, and contains the trivial character.

\bigskip
\noindent (b) The image of $\Pi_{\psi_M}^G$ in $\widehat{\mathcal{S}}_{\psi_M}$ generates $\widehat{\mathcal{S}}_{\psi_M}$, and contains the trivial character.
\end{proposition}

\noindent Proposition 8.2.4 plays a similar role as proposition 7.1.4 in using global methods to establish the local theorems. For parameters $\psi=\widehat{\phi}$ that are not tamely ramified linear, one has the result of Ban \cite{Ban} weaker than (8.2.11). It gives an isomorphism between the representation theoretic $R$-groups of $\pi_M$ and $\widehat{\pi}_M$. From the results of section seven, we already have the isomorphism between the representation theoretic $R$-group $R(\pi_M)$ and the endoscopic $R$-group $R_{\phi}=R_{\psi}$. Then in terms of the linear forms occuring in (8.2.11), the result of Ban can then be stated as:
\begin{proposition} \cite{Ban}
For $\psi=\widehat{\phi}$ as above and $\pi_M \in \Pi_{\phi_M}$, there is a sign character $\epsilon_{\pi_M}$ on $\mathcal{N}_{\psi}$ that is pulled back from $R_{\phi}=R_{\psi}$, such that
\begin{eqnarray}
f_G(\psi,u,\widehat{\pi}_M) = \epsilon_{\pi_M}(u) f_G(\phi, u,\pi_M)
\end{eqnarray}
\end{proposition}

\noindent We remark that the results in section 8 of \cite{Ban} are actually stated only for split orthogonal and symplectic groups, but with the formalism for quasi-split unitary groups of the present paper, the results of {\it loc. cit.} extends to this case also. The partial information from proposition 8.2.4 and 8.2.5 is then combined with the global method to establish the local intertwining relation in general, in section 8.5.

\subsection{Construction of global parameters with local constraints}

This is parallel to section 7.3. To apply global methods we need to globalize the local parameters to global parameters, with particular constraints at local places in order to apply the results of section 6. We can assume that the parameter in question is non-generic, as the generic case has already been dealt with in section 7. As in \cite{A1}, the difference with the generic case is that, in the present situation, the local constraints are imposed at a finite set of non-archimedean places (whereas in the generic case it is imposed at the archimedean places).

As before $G=(G,\xi) \in \widetilde{\mathcal{E}}_{\simp}(N)$, and $\psi \in \Psi(G)$. Similar to the generic case we write $\psi^N :=\xi_* \psi \in \widetilde{\Psi}(N)$ and assume that it has the form:
\begin{eqnarray}
\psi^N = l_1 \psi_1^{N_1} \oplus \cdots \oplus l_r \psi_r^{N_r}
\end{eqnarray}
with $\psi_i^{N_i} \in \widetilde{\Psi}_{\simp}(N_i)$ of the form
\begin{eqnarray}
\psi_i^{N_i} = \mu_i \otimes \nu_i, 
\end{eqnarray}
with $\mu_i$ being the generic component of $\psi_i^{N_i}$, and $\nu_i$ being a finite dimensional algebraic representation of $\SU_2$. We have $\psi_i^{N_i} = \xi_{i,*} \psi_i$ for $\psi_i \in \Psi_{\simp}(G_i)$ and certain $G_i=(G_i,\xi_i) \in \widetilde{\mathcal{E}}_{\simp}(N_i)$. Note that since we are assuming $\psi$ is not generic, we have $\mu_i \in \widetilde{\Phi}_{\simp}(m_i)$ with $m_i < N$ for every index $i$. In particular from our induction hypothesis we may assume that the global theorems are valid for generic parameters of degree $m_i$.

\noindent Finally as in the generic case we let $M$ be Levi subgroup of $G$ such that there is $\psi_M \in \Psi_2(M,\psi)$. We regard $M$ as the twisted endoscopic datum $M=(M,\xi) \in \widetilde{\mathcal{E}}(N)$. Then
\begin{eqnarray*}
\xi_* \psi_M = \bigoplus_{l_i  \odd} \psi_i^{N_i}.
\end{eqnarray*}

\noindent The globalization result that we need for the next two subsections is as follows:

\begin{proposition}
We can globalize the local data:
\[
(F,E,G,\psi, \{ \psi_i^{N_i} = \mu_i \otimes \nu_i \},\{G_i\},\{\psi_i\},M,\psi_M)
\]
to global data:
\[
(\dot{F}, \dot{E},\dot{G},\dot{\psi},\{ \dot{\psi}_i^{N_i} = \dot{\mu}_i \boxtimes \nu_i  \} , \{\dot{G}_i\},\{\dot{\psi}_i \},\dot{M},\dot{\psi}_{\dot{M}})
\]
with $\dot{E}$ a totally imaginary extension of a totally real field $\dot{F}$, $\dot{\mu}_i \in \dot{\widetilde{\Phi}}_{\simp}(m_i)$, $\dot{G}_i=(\dot{G}_i,\dot{\xi}_i) \in \dot{\widetilde{\mathcal{E}}}_{\simp}(N_i)$, and $\dot{\psi}_i^{N_i}=\dot{\xi}_{i,*} \dot{\psi}$ with $\dot{\psi} \in \dot{\Psi}_{\simp}(\dot{G}_i)$, and $u$ is a place of $\dot{F}$ that does not split in $\dot{E}$, such that
\begin{eqnarray*}
& & (\dot{F}_u, \dot{E}_u,\dot{G}_u,\dot{\psi}_u,\{ \dot{\psi}_{i,u}^{N_i} = \dot{\mu}_{i,u} \boxtimes \nu_i  \} , \{\dot{G}_{i,u}\},\dot{M}_u,\dot{\psi}_{\dot{M},u}) \\
&=& (F,E,G,\psi, \{ \psi_i^{N_i} = \mu_i \otimes \nu_i \},\{G_i\},M,\psi_M)
\end{eqnarray*}
and such that the following conditions hold:

\bigskip
\noindent (i) The natural morphisms:
\begin{eqnarray*}
\mathcal{S}_{\dot{\psi}} & \rightarrow & \mathcal{S}_{\psi} \\
 \mathcal{S}_{\dot{\psi_M}}    & \rightarrow & \mathcal{S}_{\psi_M}
\end{eqnarray*}
are isomorphisms

\bigskip
\noindent (ii)(a) Put $U:=S_{\infty}(u)$, where $S_{\infty}$ is the set of archimedean places of $\dot{F}$. Then for every $v \in U - \{u\}$, the infinitesimal characacter of the distinct irreducible generic constituents of $\dot{\psi}_v^N$ (which are just one-dimensional characters of $W_{\dot{E}_v} \cong \mathbf{C}^{\times}$) are in general position.

\bigskip

\noindent (ii)(b) There is a finite set of non-archimedean places $V$ of $\dot{F}$, such that each $v \in V$ does not split in $\dot{E}$ and is (at most) tamely ramified in $\dot{E}$, such that $\dot{\mu}_{i,v} \in \dot{\widetilde{\Phi}}_{\ellip,v}(m_i)$ is a direct sum of distinct tamely ramified characters of $W_{\dot{E}_v}$. Furthermore, for $v \notin U \cup V$, the generic components of $\dot{\psi}^N_v$ are unramified.

\bigskip

\noindent The finite set of places $V$ satisfy the following additional constraints

\bigskip
\noindent (iii)(a) The natural maps 
\begin{eqnarray}
\Pi_{\dot{\psi}_{V}}^{\dot{G}} := \bigotimes_{v \in V}  \Pi_{\dot{\psi}_{v}}^{\dot{G}_v} & \rightarrow & \widehat{\mathcal{S}}_{\dot{\psi}} \cong \widehat{\mathcal{S}}_{\psi} \\
\Pi^{\dot{G}}_{\dot{\psi}_{\dot{M},V}}  :=  \bigotimes_{v \in V}  \Pi_{\dot{\psi}_{\dot{M},v}}^{\dot{G}_v}  & \rightarrow &   \widehat{\mathcal{S}}_{\dot{\psi}_{\dot{M}}} \cong \widehat{\mathcal{S}}_{\psi_M}
\end{eqnarray} 
are surjective.

\bigskip

\noindent (iii)(b) If $l_i=1$ for all $i$, then for every $v \in V$ the following is satisfied: If $\dot{\psi}_v^N$ lies in $\xi_{v,*}^* \Psi(G^*_v)$ for some $G_v^*=(G_v^*,\xi_v^*) \in \widetilde{\mathcal{E}}_{\simp,v}(N)$, then we have $(G_v^*,\xi^*_v)=(\dot{G}_v,\dot{\xi}_v)$ as (equivalences classes of) endoscopic datum in $\widetilde{\mathcal{E}}_{\simp,v}(N)$.

\bigskip

\noindent (iii)(c) If $l_i >1$ for some $i$ then for every $v \in V$ the following is satisfied: the kernel of the mapping:
\[
\mathcal{S}_{\dot{\psi}} \rightarrow \mathcal{S}_{\dot{\psi}_v} \rightarrow R_{\dot{\psi}_v}
\]
contains no elements whose image in $R_{\dot{\psi}}$ belongs to $R_{\dot{\psi},\reg}$.
\end{proposition}

\begin{rem}
\end{rem}
\noindent For the application in the next two subsections, we also demand the following choice of $\dot{F}$ in the case when $F=\mathbf{R}$ and when the infinitesimal characters of distinct irreducible constituents of $\psi^N$ are in general position, namely that we simply take $\dot{F}=\mathbf{Q}$, and $\dot{E}$ to be an imaginary quadratic extension of $\mathbf{Q}$ that ramifies over $\mathbf{Q}$ only at odd primes. 

\begin{proof}
Essentially the same as proposition 7.3.1. The only difference is that in the generic case of proposition 7.3.1 the set $V$ consists of archimedean places, whereas in the present situation $V$ consists of non-archimedean places. Thus for instance condition (iii)(a) of proposition 7.3.1 is based on proposition 7.1.4, while for the present situation condition (iii)(a) is based on proposition 8.2.4; to ensure that one has sufficiently many tamely ramified characters, one simply need to take the residue characteristics of the set of places $V$ of $\dot{F}$ to be sufficiently large.
\end{proof}

\subsection{Local packets for square-integrable parameters}

We finish the proof of the local theorems in this and the next subsection. We emphasize again that the arguments in chapter 7 of \cite{A1} work without change for the unitary case. So we will only sketch the arguments. In this subsection, we deal with the case of square-integrable parameters.

Thus let $\psi \in \Psi_2(G)$. We are going to establish theorem 3.2.1 for the parameter $\psi$. We can assume that $\psi$ is non-generic (since the generic case is already dealt with in section 7). In particular then the generic components of $\psi^N$ have degree strictly less than $N$. We now again apply global methods. By proposition 8.3.1, we globalize the data:
\[
(F,G,\psi,\psi^N)
\]
to the data
\[
(\dot{F},\dot{G},\dot{\psi},\dot{\psi}^N)
\]
with respect to a place $u$ of $\dot{F}$ such that $\dot{F}_u=F$. As in the case of generic parameters we put ourselves in the formalism in the beginning of section 6.1, with the global family of parameters $\dot{\widetilde{\mathcal{F}}}$ being defined by the simple constituents of $\dot{\psi}^N$. In order to apply the global results from section 6, we need to verify Hypothesis 6.1.1 (we cannot apply the results of section 6.4 because the results of section 6.4 apply only to generic parameters). However, with proposition 8.2.1 in hand (applied to all the completions $\dot{F}_v$ and the parameters $\dot{\psi}_v$ of degree $N$), it follows that Hypothesis 6.1.1 is satisfied for the global parameter $\dot{\psi}$. 

We have the validity of the stable multiplicity formula:
\begin{eqnarray}
S^{\dot{G}}_{\disc,\dot{\psi}^N}(\dot{f}) = \frac{1}{|\mathcal{S}_{\dot{\psi}}|} \epsilon(\dot{\psi}) \dot{f}^{\dot{G}}(\dot{\psi}).
\end{eqnarray}
Indeed the proof is similar to the generic case treated in proposition 6.4.6 and 6.4.7; the only difference is that we utilize the local conditions given by proposition 8.3.1(iii), instead of the conditions given by proposition 7.3.1(iii) in the generic case; {\it c.f.} \cite{A1}, Lemma 7.3.2.

Then as in section 7.7, we have the following consequence from (8.4.1), which is an elementary version of the endoscopic expansion in section 5.6:
\begin{eqnarray}
I^{\dot{G}}_{\disc,\dot{\psi}^N}(\dot{f}) = \frac{1}{|\mathcal{S}_{\dot{\psi}}|} \sum_{\dot{x} \in \mathcal{S}_{\dot{\psi}}} \epsilon(\dot{\psi}^{\prime}) \dot{f}^{\prime}(\dot{\psi}^{\prime})
\end{eqnarray}
where as usual in (8.4.2) we have the endoscopic correspondence of parameters
\[
(\dot{G}^{\prime},\dot{\psi}^{\prime}) \leftrightarrow (\dot{\psi},\dot{x}).
\]
Again the global induction hypothesis gives, since $\psi \in \Psi_2(\dot{G})$:
\[
I^{\dot{G}}_{\disc,\dot{\psi}^N}(\dot{f}) = \tr R^{\dot{G}}_{\disc,\dot{\psi}^N}(\dot{f}).
\]
So we can write:
\begin{eqnarray}
\sum_{\dot{\pi}} n(\dot{\pi}) \dot{f}_{\dot{G}}(\dot{\pi}) &=& \frac{1}{|\mathcal{S}_{\dot{\psi}}|} \sum_{\dot{x} \in \mathcal{S}_{\dot{\psi}}} \epsilon(\dot{\psi}^{\prime}) \dot{f}^{\prime}(\dot{\psi}^{\prime}) \\
&=&  \frac{1}{|\mathcal{S}_{\dot{\psi}}|} \sum_{\dot{x} \in \mathcal{S}_{\dot{\psi}}} \epsilon_{\dot{\psi}}(s_{\dot{\psi}} \dot{x}) \dot{f}^{\prime}(\dot{\psi}^{\prime}) \nonumber
\end{eqnarray}
with $\dot{\pi}$ runs over irreducible representations of $\dot{G}(\dot{\mathbf{A}})$, and $n(\dot{\pi})$ is its multiplicity in $R^{\dot{G}}_{\disc,\dot{\psi}^N}$, in particular a non-negative integer (and also $n(\dot{\pi})=0$ if $\dot{\pi}$ is not unitary). Here in the second equality we have applied the endoscopic sign lemma 5.6.1.

We now specialize (8.4.3), and choose $\dot{f}$ to be decomposable:
\[
\dot{f} = \dot{f}_U \cdot \dot{f}_V \cdot \dot{f}^{U,V}
\]
where $U=S_{\infty}(u)$ and $V$ are the sets of places of $\dot{F}$ as in proposition 8.3.1. We can then apply the expansion (8.2.4) to places in $U$, $V$, and outside $U \cup V$: 
\begin{eqnarray}
\dot{f}^{\prime}_U(\dot{\psi}_U^{\prime}) = \sum_{\dot{\sigma}_U \in \Sigma_{\dot{\psi}_U}} \langle s_{\dot{\psi}} \dot{x}_U , \dot{\sigma}_U  \rangle \dot{f}_{U,\dot{G}}(\dot{\sigma}_U)
\end{eqnarray}
\begin{eqnarray}
\dot{f}^{\prime}_V(\dot{\psi}_V^{\prime}) = \sum_{\dot{\sigma}_V \in \Sigma_{\dot{\psi}_V}} \langle s_{\dot{\psi}} \dot{x}_V , \dot{\sigma}_V  \rangle \dot{f}_{V,\dot{G}}(\dot{\sigma}_V)
\end{eqnarray}
\begin{eqnarray}
(\dot{f}^{U,V})^{\prime}( (\dot{\psi}^{U,V})^{\prime}) = \sum_{\dot{\sigma}^{U,V} \in \Sigma_{\dot{\psi}^{U,V}}} \langle s_{\dot{\psi}} \dot{x}^{U,V} , \dot{\sigma}^{U,V}  \rangle \dot{f}^{U,V}_{\dot{G}}(\dot{\sigma}^{U,V}).
\end{eqnarray}
Making the substitution in (8.4.3) (and making the shift $\dot{x} \rightarrow s_{\dot{\psi}} \dot{x}$) we have
\begin{eqnarray}
& & \sum_{\dot{\pi}} n(\dot{\pi}) \dot{f}_{\dot{G}}(\dot{\pi}) \\
&=& \frac{1}{|\mathcal{S}_{\dot{\psi}}|} \sum_{\dot{x} \in \mathcal{S}_{\dot{\psi}}} \epsilon_{\dot{\psi}}(\dot{x})  \langle \dot{x},\dot{\sigma}_U \rangle \cdot \langle \dot{x},\dot{\sigma}_V \rangle \cdot \langle \dot{x}, \dot{\sigma}^{U,V } \rangle \dot{f}_{\dot{G}} (\dot{\sigma}). \nonumber
\end{eqnarray}

\noindent Now let $\xi \in \widehat{\mathcal{S}}_{\dot{\psi}} \cong \widehat{\mathcal{S}}_{\psi}$ be arbitrary. Then by the condition imposed on the set of places $V$ in proposition 8.3.1, we can pick a $\dot{\sigma}_{V,\xi}  \in \Sigma_{\dot{\psi}_V}$, such that $\dot{\sigma}_{V,\xi} = \dot{\pi}_{V,\xi}$ is an irreducible representation of $\dot{G}(\dot{\mathbf{A}}_V)$, and that
\[
\langle \cdot , \dot{\sigma}_{V,\xi} \rangle = \epsilon_{\dot{\psi}}^{-1} \xi(\cdot)^{-1}.
\]
Similarly, at places $v$ outside $U \cup V$, the parameter $\dot{\psi}_v$ is unramified, so in particular is a tamely ramified linear parameter, and part (a) of proposition 8.2.4 applies, and thus we can pick a $\dot{\sigma}^{U,V}(1) \in \Sigma^{U,V}_{\dot{\psi}}$ that is an irreducible representation $\dot{\pi}^{U,V}(1)$ such that the corresponding character 
\[
\langle \,\ \cdot \,\ , \dot{\pi}^{U,V}(1) \rangle
\]
is trivial.
Thus setting for $\dot{\pi}_U \in \Pi(\dot{G}_U)$:
\[
n_{\psi}(\xi,\dot{\pi}_U) = n(\dot{\pi}_U \otimes \dot{\pi}_{V,\xi} \otimes \dot{\pi}^{U,V}(1)  )
\]
which is a non-negative integer, we can extract, from (8.4.7) the identity:
\begin{eqnarray}
\sum_{\dot{\pi}_U \in \Pi(\dot{G}_U)} n_{\psi}(\xi,\dot{\pi}_U) \dot{f}_{U,\dot{G}}(\dot{\pi}_U) &=& \frac{1}{|\mathcal{S}_{\dot{\psi}}|} \sum_{\dot{x} \in \mathcal{S}_{\dot{\psi}}} \langle \dot{x},\dot{\sigma}_U \rangle \xi(\dot{x})^{-1} \dot{f}_{U,\dot{G}}(\dot{\sigma}_U) \\
&=& \sum_{\dot{\sigma}_U \in \Sigma_{\dot{\psi}_U}(\xi) } \dot{f}_{U,\dot{G}}(\dot{\sigma}_U) \nonumber
\end{eqnarray}
with $\Sigma_{\dot{\psi}_U}(\xi)$ being the subset consisting of $\dot{\sigma}_U$ such that the associated character satisfies:
\[
\langle \cdot , \dot{\sigma}_U \rangle = \xi(\cdot).
\]

\noindent At this point we first consider the case that $F=\mathbf{R}$, and that the infinitesimal characters of the distinct irreducible generic components of $\psi=\dot{\psi}_u$ are in sufficient general position, then in proposition 8.3.1 we can simply choose $\dot{F}=\mathbf{Q}$, and thus $U=\{u\}$, and (8.4.8) gives the desired result: for the unique $\sigma \in \Sigma_{\psi}$ such that $\langle \cdot, \sigma \rangle =\xi$, the linear form $f_G(\sigma)$ is a linear combination of characters of irreducible unitary representations with non-negative integral coefficients. Since $\xi \in \widehat{\mathcal{S}}_{\psi} \cong \widehat{\mathcal{S}}_{\dot{\psi}}$ is arbitrary this gives the desired result, and Theorem 3.2.1 thus hold in this case.

\noindent In the remaining cases, one uses the results just obtained, combined with a result of Arthur which is of separate importance: for a general parameter $\psi \in \Psi(G)$, recall the generic parameter $\phi_{\psi} \in \Phi(G)$ associated to $\psi$. One has $S_{\psi} \subset S_{\phi_{\psi}}$ and $\mathcal{S}_{\psi}$ surjects to $\mathcal{S}_{\phi_{\psi}}$. The packet $\Pi_{\phi_{\psi}}$ of irreducible representations associated to the generic parameter $\phi_{\psi}$ is constructed from Langlands quotient of the packet of standard representations associated to $\phi_{\psi}$ (thus given by the results of section 7).  

\begin{proposition} (Proposition 7.4.1 of \cite{A1})
Suppose that theorem 3.2.1 holds for a parameter $\psi \in \Psi(G)$. Then the elements of $\Pi_{\phi_{\psi}}$ belong to $\Pi_{\psi}$ (in particular are unitary representations) and they occur with multiplicity one in the packet $\Pi_{\psi}$, such that the following diagram commutes:
\begin{eqnarray*}
\xymatrix{ \Pi_{\phi_{\psi}} \ar@{^{(}->}[r] \ar@{^{(}->}[d]  & \Pi_{\psi}  \ar[d] \\  \widehat{\mathcal{S}}_{\phi_{\psi}} \ar@{^{(}->}[r] & \widehat{\mathcal{S}}_{\psi} }
\end{eqnarray*}
\end{proposition}

In the remaining cases, we specialize (8.4.8) as follows: for each $v \in S_{\infty}^u \subset U$, the parameter $\dot{\psi}_v$ is chosen as in proposition 8.3.1 so that the distinct irreducible generic constituents of $\dot{\psi}_v$ are in sufficient general position. Hence we can apply the result just established for $\dot{\psi}_v$. By proposition 8.4.1, we can choose $\dot{\pi}_v(1) \in \Sigma_{\dot{\psi}_v} = \Pi_{\dot{\psi}_v}$ corresponding to the trivial character in $\widehat{\mathcal{S}}_{\phi_{\dot{\psi}_v}} \hookrightarrow \widehat{\mathcal{S}}_{\dot{\psi}_v}$, which by the proposition occurs with multiplicity one in the packet $\Pi_{\dot{\psi}_v}$. Thus in (8.4.8) if we choose the function $\dot{f}_v \in \mathcal{H}(\dot{G}_v)$ for $v \in S_{\infty}^u$ such that for any $\dot{\sigma}_v \in \Pi_{\dot{\psi}_v}$:
\begin{eqnarray}
\dot{f}_{v,\dot{G}_v}(\dot{\sigma}_v) = \left \{ \begin{array}{c} 1 \mbox{  if }  \dot{\sigma}_v = \dot{\pi}_v(1) \\ 0 \mbox{ if } \dot{\sigma}_v \neq \dot{\pi}_v(1) \end{array} \right.
\end{eqnarray}
then (8.4.8) reduces to
\begin{eqnarray}
\sum_{\pi \in \Pi_{\unit}(G)} n_{\psi}(\xi, \pi   \bigotimes_{v \in S_{\infty}^u}\dot{\pi}_v(1) ) \cdot f_G(\pi) = f_G(\sigma(\xi)), \,\ f \in \mathcal{H}(G)
\end{eqnarray}
where $\sigma(\xi)$ is the element in the packet $\Sigma_{\psi}$ such that $\langle \cdot , \sigma \rangle = \xi$. Since $\xi \in \widehat{\mathcal{S}}_{\psi}$ is arbitrary this again proves that any elemet in the packet $\Sigma_{\psi}$ is a non-negative integral linear combination of irreducible unitary characters. In particular, the non-negative integer $n_{\psi}(\xi,\pi) := n_{\psi}(\xi, \pi   \bigotimes_{v \in S_{\infty}^u}\dot{\pi}_v(1) )$ does not depend on auxiliary choices.

\noindent Thus to conclude, we simply need to define the packet $\Pi_{\psi}$ as:
\begin{eqnarray*}
\Pi_{\psi} = \coprod_{\xi \in \widehat{\mathcal{S}}_{\psi}} \Pi_{\psi}(\xi)
\end{eqnarray*}
with $\Pi_{\psi}(\xi)$ being the (finite) multi-set consisting of $n_{\psi}(\xi,\pi)$ copies of $\pi$, as $\pi$ ranges over $\Pi_{\unit}(G)$, then the endoscopic character identity of part (b) of Theorem 3.2.1 is satisfied, with $\pi \in \Pi_{\psi}(\xi)$ being sent to the character: 
\[
\langle \cdot, \pi \rangle = \xi(\cdot)
\]
of $\mathcal{S}_{\psi}$.

\noindent The following result which comes out of the construction is worth emphasizing:
\begin{proposition}
Suppose that $\psi \in \Psi_2(G)$. Then all the elements of the packet $\Pi_{\psi}$ are automorphic, in the sense that they occur as local components of discrete automorphic representations.
\end{proposition}

\bigskip

\noindent Thus theorem 3.2.1 holds for the case where $\psi \in \Psi_2(G)$. It remains to treat the case of non-square-integrable parameters. By proposition 3.4.4, this is a consequence of the local intertwining relation, and it remains to establish this to complete the inductive step for the proof of all the local theorems.

\subsection{The local intertwining relation}
It remains to prove the local intertwining relation. Thus we let $\psi \in \Psi(G) \smallsetminus \Psi_2(G)$, and $M$ be a proper Levi subgroup of $G$ such that $\psi_M \in \Psi_2(M,\psi)$. Again the local intertwining relation is established via global methods, similar to the argument in section 7.4.

As in the case of generic parameters the local version of the descent argument in the proof of proposition 5.7.4 reduces the proof of the local intertwining relation to the following three cases:

\begin{eqnarray}
& & \psi^N = 2 \psi_1^{N_1} \oplus \cdots \oplus 2 \psi_q^{N_q} \oplus \psi_{q+1}^{N_{q+1}} \oplus \cdots \oplus \psi_r^{N_r}, \,\ q \geq 1 \\
& &  S_{\psi} = \prod_{i=1}^q O(2,\mathbf{C}) \times \prod_{i=q+1}^r O(1,\mathbf{C})    \nonumber
\end{eqnarray}

\begin{eqnarray}
& & \psi^N = 2 \psi_1^{N_1} \oplus \psi_2^{N_2} \oplus \cdots \oplus \psi_r^{N_r} \\
& & S_{\psi} = Sp(2,\mathbf{C}) \times \prod_{i=2}^r O(1,\mathbf{C}) \nonumber
\end{eqnarray}

\begin{eqnarray}
& & \psi^N = 3 \psi_1^{N_1} \oplus \psi_2^{N_2} \oplus \cdots \oplus \psi_r^{N_r} \\
& & S_{\psi} = O(3,\mathbf{C}) \times \prod_{i=2}^r O(1,\mathbf{C}). \nonumber
\end{eqnarray}

\noindent For  (8.5.1) we have $\psi \in \Psi_{\ellip}(G)$, and thus $\mathcal{S}_{\psi,\ellip}$ is a torsor under $\mathcal{S}_{\psi_M}$, while for (8.5.2) and (8.5.3) we have $\mathcal{S}_{\psi_M} = \mathcal{S}_{\psi}$.

\bigskip
\noindent Using proposition 8.3.1 once more, we globalize the data:
\[
(F,G,\psi,\psi_M, \{\psi_i^{N_i}\})
\]
to global data:
\[
(\dot{F},\dot{G},\dot{\psi},\dot{\psi}_M,\{\dot{\psi}_i^{N_i}\})
\]
with respect to a place $u$ of $\dot{F}$ such that $\dot{F}_u=F$. As in section 7 we then form the global family:
\[
\dot{\widetilde{\mathcal{F}}} = \dot{\widetilde{\mathcal{F}}}(\dot{\psi}_1^{N_1},\cdots,\dot{\psi}_r^{N_r}).
\]

\bigskip
\noindent We have $\mathcal{S}_{\dot{\psi}} \cong \mathcal{S}_{\psi} = \mathcal{S}_{\psi_M} \cong \mathcal{S}_{\dot{\psi}_M}$. In cases (8.5.2) and (8.5.3), we have by proposition 5.7.5:
\begin{eqnarray}
\sum_{x \in \mathcal{S}_{\psi}} \epsilon^{\dot{G}}_{\dot{\psi}}(\dot{x}) (\dot{f}^{\prime}_{\dot{G}}(\dot{\psi},s_{\psi} \dot{x}) - \dot{f}_{\dot{G}}(\dot{\psi},\dot{x}) ) =0, \,\ \dot{f} \in \mathcal{H}(\dot{G}).
\end{eqnarray} 

\noindent We have a similar result for the case of elliptic parameter:
\begin{proposition}
For $\psi \in \Psi_{\ellip}(G)$ as in (8.5.1), we have
\begin{eqnarray}
\sum_{x \in \mathcal{S}_{\psi,\ellip}} \epsilon_{\dot{\psi}}^{\dot{G}}(\dot{x})(\dot{f}^{\prime}_{\dot{G}}(\dot{\psi},s_{\psi} \dot{x}) - \dot{f}_{\dot{G}}(\dot{\psi},\dot{x})  ) =0.
\end{eqnarray}
\end{proposition}
\begin{proof}
This is the assertion given as part of the statement of Lemma 7.3.1 of \cite{A1}. We again apply the global result of section 6. As observed in the previous subsection, Hypothesis 6.1.1 for the global family $\dot{\widetilde{\mathcal{F}}}$ is satisfied, as it  already follows from Proposition 8.2.1 (applied to the localization of parameters for each place $v$ of $\dot{F}$), and in particular we can apply Proposition 6.2.1 (when $r >1$) or Proposition 6.2.2 (when $r=1$) to the global parameter $\dot{\psi}$ (we cannot apply results of section 6.4 for the same reason, namely that applies only to generic parameters). When treat the case $r >1$, as the case $r=1$ is similar. Applied to our setting, Proposition 6.2.1 asserts that there is a non-zero constant $c$ such that the following identity holds:
\begin{eqnarray}
& & \sum_{\dot{G}^* \in \dot{\widetilde{\mathcal{E}}}_{\simp}(N)} \dot{\widetilde{\iota}}(N,\dot{G^*}) \tr R_{\disc,\dot{\psi}^N}^{\dot{G}^*}(\dot{f}^*)\\
&=& c   \sum_{\dot{x} \in \mathcal{S}_{\psi,\ellip}} \epsilon^{\dot{G}}_{\dot{\psi}}(\dot{x}) \big(  \dot{f}^{\prime}_{\dot{G}}(\dot{\psi},s_{\psi} \dot{x}) - \dot{f}_{\dot{G}}(\dot{\psi},\dot{x})       \big)  \nonumber
\end{eqnarray}
for any compatible family of functions $\{\dot{f}^* \}$, with $\dot{f}$ being the function corresponding to $\dot{G} \in \dot{\widetilde{\mathcal{E}}}_{\simp}(N)$. Then the main point is that we can write (8.5.6) in the form as (4.3.31), so that Proposition 4.3.9 (which is Corollary 3.5.3 of \cite{A1}) on vanishing of coefficients can be applied to yield the vanishing of both sides of (8.5.6). With the condition on the global parameter $\dot{\psi}$ imposed at the set of places $V$, the proof of lemma 7.3.1 of \cite{A1} applies without change to the present situation. Here we note that it is in the proof of this proposition that the result of Ban \cite{Ban}, that we stated as proposition 8.2.5, is needed to write the right hand side of (8.5.6) in the form of the right hand side of (4.3.31).

\end{proof}

\noindent Similar to the case of generic parameters, we extract the local intertwining relation from the global identity (8.5.5) in the case (8.5.1), and from (8.5.4) in the case (8.5.2) and (8.5.3). We again just give a summary of the argument extracted from section 7.3 of \cite{A1}. We treat the case (8.5.1) as the other two cases are similar. For (8.5.1), it suffices (again by descent argument) to treat the local intertwining relation for $x \in \mathcal{S}_{\psi,\ellip}$.

\noindent We again choose decomposable function
\[
\dot{f} = \dot{f}_U \cdot \dot{f}_V \cdot \dot{f}^{U,V}
\] 
with $U=S_{\infty}(u)$ and $V$ being the set of places of $\dot{F}$ as in proposition 8.3.1. In particular, for $v \notin U \cup V$, the parameter $\dot{\psi}_v$ is unramified. The following result can be proved as in \cite{A1}, which allow us to remove the contribution to (8.5.5) from places $v \notin U \cup V$:

\begin{lemma} (Lemma 7.3.4, combined with Lemma 7.3.3 of \cite{A1})
Suppose that in general $\psi$ is an unramified parameter of $G(F)$. Then for $f$ the characteristic function of the (standard) hyperspecial maximal compact subgroup of $G(F)$, we have
\begin{eqnarray}
f^{\prime}_G(\psi,s_{\psi}x) =f_G(\psi,x)=1, \,\ x \in \mathcal{S}_{\psi}.
\end{eqnarray}
\end{lemma}  

\noindent Thus by choosing $\dot{f}^{U,V} = \prod_{v \notin U \cup V} \dot{f}_v$, with $\dot{f}_v$ being the characteristic function of the standard maximal compact subgroup of $\dot{G}(\dot{F}_v)$, we can apply (8.5.7), and hence the contributions from the places $v \notin U \cup V$ to (8.5.5) can be removed. For a place $v \in V$, the global parameter $\dot{\psi}$ is chosen as in proposition 8.3.1 so that $\dot{\psi}_v$ is a tamely ramified linear parameter, hence dual to a tamely ramified linear generic parameter $\phi_v$. We have the following:

\begin{lemma} (Lemma 7.3.3 of \cite{A1})
Suppose that $\psi = \widehat{\phi} \in \Psi(G)$ is a tamely ramified linear parameter, and that $f \in \mathcal{H}(G)$ is chosen such that the function
\[
\pi \mapsto f(\widehat{\pi}), \,\ \pi \in \Pi_{\phi}
\] 
is supported on the set of constituents of $\mathcal{I}_P(\pi_M)$ for $\pi_M \in \Pi^G_{\phi_M}$. Then 
\[
f^{\prime}_G(\psi,s_{\psi}x) =f_G(\psi,x), \,\ x \in \mathcal{S}_{\psi}.
\]
\end{lemma} 

\noindent Thus by choosing $\dot{f}_v$ for $v \in V$ such that the function $\dot{f}_v$ satisfies the condition as in lemma 8.5.3, with respect to the parameter $\dot{\psi}_v$, then we have:
\[
\dot{f}^{\prime}_{\dot{G}_v}(\dot{\psi}_v,s_{\psi} \dot{x}_v) = \dot{f}_{\dot{G}_v}(\dot{\psi}_v,\dot{x}_v)
\]   
and hence (8.5.5) simplifies to:
\begin{eqnarray}
\end{eqnarray}
\begin{eqnarray*}
\sum_{x \in \mathcal{S}_{\psi,\ellip}} \epsilon^{\dot{G}}_{\dot{\psi}}(\dot{x}) \dot{f}_{V,\dot{G}}(\dot{\psi}_V,\dot{x}_V) \big(  \dot{f}_{U,\dot{G}}^{\prime}(\dot{\psi}_U, s_{\psi} \dot{x}_U  ) -   \dot{f}_{U,\dot{G}}(\dot{\psi}_U,  \dot{x}_U  )  \big) =0.
\end{eqnarray*}

\noindent Recall that $\mathcal{S}_{\psi,\ellip}$ is a torsor under $\mathcal{S}_{\psi_M}$. Hence by the conditions imposed on $\dot{\psi}_V$ as in proposition 8.3.1, we can vary $\dot{f}_V$ under the constraint as in lemma 8.5.3, and still be able to isolate the contribution of the term in (8.5.7) from any given $x \in \mathcal{S}_{\psi,\ellip}$. Thus we have:
\begin{eqnarray}
 \dot{f}_{U,\dot{G}}^{\prime}(\dot{\psi}_U, s_{\psi} \dot{x}_U  ) =  \dot{f}_{U,\dot{G}}(\dot{\psi}_U,  \dot{x}_U  ). 
\end{eqnarray} 

\noindent The rest is then the same as in section 8.4. Namely when $F=\mathbf{R}$, and when the infinitesimal characters of the distinct irreducible constituents of $\psi^N$ are in general position, then we took $\dot{F}=\mathbf{Q}$ in proposition 8.3.1, and hence $U=\{u\}$. Thus (8.5.9) gives the local intertwining relation for $\dot{\psi}_u=\psi$. With this case established, the general case follows. Namely in accordance with proposition 8.3.1, the archimedean parameter $\dot{\psi}_v$ for $v \in S_{\infty}^u$ is of the type just treated. Hence choosing $\dot{f}_v$ for $v \in S_{\infty}^u$ such that $\dot{f}_{v,\dot{G}}(\dot{\psi}_v,\dot{x}_v) \neq 0$, we can cancel the contribution of the of the places from $v \in S_{\infty}^u$ in (8.5.9), and thus yield the local intertwining relation for $\dot{\psi}_u=\psi$.

\bigskip

With the local intertwining relation for general parameters established, we in particular completed the proof of Theorem 3.2.1, by the reduction given by Proposition 3.4.4. As in the case of generic parameters ({\it c.f.} corollary 7.4.7), we also obtain the following corollary from the general local intertwining relation:
\begin{corollary}
For $w^0 \in W_{\psi}^0$, we have
\[
R_P(w^0,\widetilde{\pi}_M,\psi_M) =1.
\]
\end{corollary} 
\begin{proof}
Same argument as in Corollary 6.4.5 of \cite{A1}.
\end{proof}

\bigskip

To complete the induction argument in this section, it remains to finish the induction argument for the global theorems for the families $\dot{\widetilde{\mathcal{F}}}$ of global non-generic parameters used in this section. The only relevant global theorems for non-generic parameters are theorem 2.5.3(b) ({\it c.f} remark 2.5.7), which we established in proposition 6.1.5. The other global theorems are the stable multiplicity formula, and the spectral multiplicity formula theorem 2.5.2 (for theorem 5.2.1 it simply follows from the local intertwining relation we established, applied to each place $v$ of $\dot{F}$).

The stable multiplicity formula follows from Proposition 5.7.4 in the case where the parameter $\dot{\psi}^N$ is ``degenerate", i.e. when $\dot{\psi}^N \notin \dot{\widetilde{\mathcal{F}}}_{\ellip}(N)$, and $\dot{\psi}^N \notin (\dot{\xi}^*)_* \dot{\widetilde{\mathcal{F}}}_{\ellip}(\dot{G}^*)$ for any $\dot{G}^*=(\dot{G}^*,\dot{\xi}^*) \in \dot{\widetilde{\mathcal{E}}}_{\simp}(N)$. 

In the case where $\dot{\psi}^N \in \dot{\xi}_* \dot{\widetilde{\mathcal{F}}}_{\ellip}(\dot{G})$ for $\dot{G}=(\dot{G},\dot{\xi}) \in \dot{\widetilde{\mathcal{E}}}_{\simp}(N)$, the proof of the stable multiplicity formula for $\dot{\psi}^N$ (with respect to any $\dot{G}^*=(\dot{G}^*,\dot{\xi}^*) \in \dot{\widetilde{\mathcal{E}}}_{\simp}(N)$) follows the same arguments as in the proof of proposition 6.4.4-6.4.7 (for proposition 6.4.7 we simply interpret the sign $\delta_{\psi}$ as being equal to $+1$ in the case of non-generic parameters); indeed as we observe in the beginning of section 8.4, the only difference being in the use of the local conditions given by proposition 8.3.1(iii), instead of the conditions given by proposition 7.3.1(iii) in treating the case for generic parameters.

With the stable multiplicity formula, the spectral multiplicity formula follows by application of lemma 5.7.6, together with the local results established in this section. 

\bigskip

Thus we have completed the proof of all the local theorems.

\section{\bf Global classification}

\subsection{Completion of induction arguments, part I}

In section 7 and 8 we complete the induction arguments concerning the local theorems. In this final section we complete the induction argumetns for the global classification theorems. Again we rely on section 6 to establish the results. Thus $F$ will now denote a global field.

We now simply take the family $\widetilde{\mathcal{F}}$ simply to be the maximal family $\widetilde{\Psi}$. In order to apply results of section 6, we need to verify that the Hypothesis 6.1.1 is satisfied:

\begin{proposition}
Hypothesis 6.1.1 is satisfied for the maximal family $\widetilde{\mathcal{F}}=\widetilde{\Psi}$.
\end{proposition}
\begin{proof}
Suppose first that $\psi^N \in \widetilde{\Psi}(N)$ is not a simple generic parameter. Then hypothesis 6.1.1 for $\psi^N$ is the assertion that if $\psi^N = \xi_* \psi$ for some $\psi \in \Psi(G)$ and $G =(G,\xi) \in \widetilde{\mathcal{E}}_{\ellip}(N)$, then the linear form
\[
\widetilde{f} \mapsto \widetilde{f}_N(\psi^N), \,\ \widetilde{f} \in \widetilde{\mathcal{H}}(N)
\]
transfers to a stable linear form $f^G(\psi)$ on $G(\mathbf{A}_F)$. Since we assumed that $\psi^N$ is not simple generic parameter, we see that the seed theorem 2.4.10 can be applied to the generic constituents of $\psi^N$. As a consequence one has $\psi^N_v = \xi_{v,*} \psi_v$ with $\psi_v \in \Psi^+(G_v)$ for any place $v$ of $F$. By the results of section 7 and 8, the local linear form
\[
\widetilde{f}_v \mapsto \widetilde{f}_{v,N}(\psi^N_v), \,\ \widetilde{f}_v \in \widetilde{\mathcal{H}}_v(N)
\]
transfers to $G_v$ as a stable linear form $f_v(\psi_v)$. It follows that the global linear form $\widetilde{f}_N(\psi^N)$ transfers to $G$, as required.

It thus remains to treat the class of simple generic parameters $\widetilde{\Phi}_{\simp}(N)$. In this case, Hypothesis 6.1.1 asserts the condition:
\begin{eqnarray}
\bigcup_{G^*  \in \widetilde{\mathcal{E}}_{\simp}(N)} \xi_* \Phi_{\simp}(G^*) = \Phi_{\simp}(N)
\end{eqnarray}
to hold. Here we recall that as in section 6.1, we are using the definition of $\Phi_{\simp}(G)$ as the set of pairs $(G,\phi^N)$ such that $S^G_{\disc,\phi^N} \nequiv 0$ and such that the linear form $\widetilde{f}_N(\phi^N)$ transfers to $G$. Thus we must show equality in (9.1.1) (we remark that at this point we still do not know that the union in (9.1.1) is disjoint; nevertheless it would follow from the stable multiplicity formula for generic parameters which we will establish till the end of this subsection). 

Thus let $\phi^N \in \widetilde{\Phi}_{\simp}(N)$. By the usual simplification of the trace formula, we have:
\[
\widetilde{I}^N_{\disc,\phi^N}(\widetilde{f}) = \tr \widetilde{R}^N_{\disc,\phi^N}(\widetilde{f}) = \widetilde{f}_N(\phi^N)
\]
and
\[
\widetilde{I}^N_{\disc,\phi^N}(\widetilde{f}) = \widetilde{\iota}(N,G) \widehat{S}^G_{\disc,\phi^N}(\widetilde{f}^G) + \widetilde{\iota}(N,G^{\vee}) \widehat{S}^{G^{\vee}}_{\disc,\phi^N}(\widetilde{f}^{G^{\vee}})
\]
where as usual we have noted $G=(G,\xi)$ and $G^{\vee}=(G^{\vee},\xi^{\vee})$ the two representatives of $\widetilde{\mathcal{E}}_{\simp}(N)$. Thus we have
\begin{eqnarray}
& & \\
& & \widetilde{f}_N(\phi^N)=\widetilde{\iota}(N,G) \widehat{S}^G_{\disc,\phi^N}(\widetilde{f}^G) + \widetilde{\iota}(N,G^{\vee}) \widehat{S}^{G^{\vee}}_{\disc,\phi^N}(\widetilde{f}^{G^{\vee}}). \nonumber
\end{eqnarray}  

\noindent In (9.1.2), if one of the distributions $S^G_{\disc,\phi^N}, \,\ S^{G^{\vee}}_{\disc,\phi^N}$ vanishes identifically, say $S^{G^{\vee}}_{\disc,\phi^N} \equiv 0$, then from (9.1.2) we must have $S^G_{\disc,\phi^N} \nequiv 0$, and that the linear form $\widetilde{f}_N(\phi^N)$ transfers to a stable linear form on $G$. 

To complete the proof we argue by contradiction, as in Lemma 8.1.1 of \cite{A1}: suppose that neither the distributions $S^G_{\disc,\phi^N}, \,\ S^{G^{\vee}}_{\disc,\phi^N}$ vanishes, and that the linear form $\widetilde{f}_N(\phi^N)$ does not transfer to $G$ or $G^{\vee}$. We shall obtain a contradiction as follows. From the local results of section 7, we see that if $\widetilde{f}_N(\phi^N)$ does not transfer to $G$, then there must be a place $v$ of $F$ (which we may assume does not split in $E$) such that the local linear form $\widetilde{f}_{v,N}(\phi^N_v)$ does not transfer to $G_v$. Similarly there is a place $v^{\vee}$ that does not split in $E$ such that the local linear form $\widetilde{f}_{v^{\vee},N}(\phi^N_{v^{\vee}})$ does not transfer to $G^{\vee}_{v^{\vee}}$.  

Suppose that $v=v^{\vee}$. Then the condition implies that we have $\phi^N_v \in \xi^{\#}_{v,*} \Phi(G_v^{\#})$ for some $G^{\#}_v=(G^{\#}_v,\xi^{\#}_v) \in \widetilde{\mathcal{E}}_{\ellip,v}(N) - \widetilde{\mathcal{E}}_{\simp_v}(N)$, and that $\phi^N_v \notin \xi_{v,*} \Phi(G_v)$, $\phi^N_v \notin \xi^{\vee}_{v,*} \Phi(G_v^{\vee})$. It follows that we can choose $\widetilde{f}_v$ such that $\widetilde{f}_{v,N}(\phi^N_v) \neq 0$, but such that $\widetilde{f}_v^{G_v}=0$ and $\widetilde{f}_v^{G^{\vee}_v} =0$ (the existence of such $\widetilde{f}_v$ follows from the by now familiar application of proposition 3.1.1(a). If we choose an arbitrary $\widetilde{f}^v$ such that $\widetilde{f}^v_N(\phi^{v,N}) \neq 0$, then the test function $\widetilde{f} = \widetilde{f}_v \widetilde{f}^v$ would give a contradiction in (9.1.2). Now if $v \neq v^{\vee}$, then we similarly choose $\widetilde{f}_v$ such that $\widetilde{f}_{v,N}(\phi^N_v) \neq 0$ and $\widetilde{f}_v^{G_v} =0$, and similarly $\widetilde{f}_{v^{\vee}}$ such that $\widetilde{f}_{v^{\vee},N}(\phi^N_{v^{\vee}}) \neq 0$ and $\widetilde{f}_{v^{\vee}}^{G^{\vee}_{v^{\vee}}} =0$. Then upon choosing an arbitrary $\widetilde{f}^{v,v^{\vee}}$ such that $\widetilde{f}_N^{v,v^{\vee}}(\phi^{v,v^{\vee},N}) \neq 0$, we again see that, with the test function $\widetilde{f} = \widetilde{f}_v \cdot \widetilde{f}_{v^{\vee}} \cdot \widetilde{f}^{v,v^{\vee}}$, we have a contradiction in (9.1.2).

\end{proof}

\bigskip

We can now begin to complete the induction arguments for the global theorems. Thus as before $G \in \widetilde{\mathcal{E}}_{\simp}(N)$. First recall that in Proposition 6.1.5 we have already completed the induction argument for the family $\widetilde{\mathcal{F}} = \widetilde{\Psi}$ for part (b) of Theorem 2.5.4 concerning the root number. We next turn to the stable multiplicity formula, which is stated as theorem 5.1.2. In proposition 5.7.4, we have already established the stable multiplicity formula for a parameter $\psi^N \in \widetilde{\Psi}(N)$ with respect to any $G^* \in \widetilde{\mathcal{E}}_{\simp}(N)$, whenever $\psi^N \notin \widetilde{\Psi}_{\ellip}(N)$, and that $\psi^N \notin \xi^*_* \Psi_{\ellip}(G^*)$ for any $G^*=(G^*,\xi^*) \in \widetilde{\mathcal{E}}_{\simp}(N)$. We also showed that $\psi^N$ does not contribute to the discrete automorphic spectrum of any $G^* \in \widetilde{\mathcal{E}}_{\simp}(N)$. 

Thus we suppose now that $\psi^N \in \xi_* \Psi_{\ellip}(G)$ for some $G=(G,\xi) \in \widetilde{\mathcal{E}}_{\simp}(N)$. Then $\psi^N=\xi_* \psi$ for $\psi \in \Psi_{\ellip}(G)$, and $G$ is uniquely determined by this condition (i.e. $\psi^N \notin \xi^{\vee}_* \Psi_{\ellip}(G^{\vee})$ for the other $G^{\vee} \in \widetilde{\mathcal{E}}_{\simp}(N)$). We first suppose that $\psi \in \Psi_{\ellip}(G) \smallsetminus \Psi_2(G)$. We can then apply Proposition 6.2.1 and 6.2.2. By the local intertwining relation that we have already established in section 7 and 8, we see that both propositions reduce to the statement:
\begin{eqnarray}
\sum_{G^* \in \widetilde{\mathcal{E}}_{\simp}(N)} \widetilde{\iota}(N,G^*) \tr R^{G^*}_{\disc,\psi^N}(f^*) =0 
\end{eqnarray}
for any compatible family of functions $\{f^*\}$. By the result on vanishing of coefficients (lemma 4.3.6), we have the vanishing:
\begin{eqnarray*}
R^{G^*}_{\disc,\psi^N} \equiv 0, \,\ G^* \in \widetilde{\mathcal{E}}_{\simp}(N).
\end{eqnarray*}
In particular $\psi^N$ does not contribute to the discrete spectrum of $G$. We also see from this that if $\psi^N$ were to lie in $\xi^{\vee}_{*}\Psi_{\ellip}(G^{\vee}) \smallsetminus \xi^{\vee}_* \Psi_{2}(G^{\vee})$, then $\psi^N$ also could not contribute to the discrete spectrum of $G$.   

Furthermore, from (6.2.8) and (6.2.11), we see, by combining the local intertwining relation and the vanishing of $R^G_{\disc,\psi^N}$ just proved, that
\begin{eqnarray}
\leftexp{0}{S^G_{\disc,\psi^N}} \equiv 0
\end{eqnarray}
i.e. the stable multiplicity formula is valid for $\psi^N$ with respect to $G$ (remark that since $\psi \in \Psi_{\ellip}(G) \smallsetminus \Psi_2(G)$, the stable multiplicity formula for $\psi^N$ with respect to $G$ reduces to the assertion $S_{\disc,\psi^N}^G$ vanishes). Similarly, from (6.2.4) and (6.2.12),  we can combine with the local intertwining relation and the vanishing of $R_{\disc,\psi^N}^{G^{\vee}}$ to obtain:
\begin{eqnarray}
\leftexp{0}{S^{G^{\vee}}_{\disc,\psi^N}} \equiv 0
\end{eqnarray}
thus again the stable multiplicity formula is valid for $\psi^N$ with respect to $G^{\vee}$. 

We now suppose that $\psi^N \in \widetilde{\Psi}_{\ellip}(N)$. If $\psi^N \in \xi_*^{\#} \Psi_2(G^{\#})$ for some $G^{\#}=(G^{\#},\xi^{\#}) \in \widetilde{\mathcal{E}}_{\ellip}(N) \smallsetminus \widetilde{\mathcal{E}}_{\simp}(N)$, the we have seen in proposition 6.1.3 and 6.1.5 that
\[
R^{G^*}_{\disc,\psi^N} \equiv 0, \,\ G^* \in \widetilde{\mathcal{E}}_{\simp}(N).
\]
Hence we now assume $\psi^N \in \xi_*\Psi_2(G)$. The case where $\psi^N \in \widetilde{\Psi}_{\simp}(N)$ is a simple parameter will be the most difficult case to treat, and which we will complete in the next subsection. Thus we assume that $\psi^N \notin \widetilde{\Psi}_{\simp}(N)$. We then apply the results of section 6.3 on supplementary parameters. Namely we apply Proposition 6.3.1 (noting that the hypothesis on stable multiplicity formula for in $\widetilde{\Psi}(N) \smallsetminus \widetilde{\Psi}_{\ellip}(N)$ has just been established). Now by the local intertwining relation, which is established in section 8.5 (for each place $v$ of $F$), we see that the right hand side of (6.3.3) vanishes. We hence obtain:
\begin{eqnarray}
& & \\
& & \sum_{G^* \in \widetilde{\mathcal{E}}_{\simp}(N_+)} \widetilde{\iota}(N_+,G^*) \tr R^{G^*}_{\disc,\psi_+^{N_+}}(f^*) + b_+ f_1^{L_+}(\psi_1 \times \Lambda) =0. \nonumber
\end{eqnarray}
Here as in the proof of proposition 6.4.6 (and with similar notations there) if $f_1$ is the function associated to $G_1^{\vee} = G_1 \times G^{\vee} \in \widetilde{\mathcal{E}}_{\ellip}(N_+)$ in the compatible family occuring in (9.1.6), then we have
\[
f^{L_+}(\psi_1 \times \Lambda) = f_1^{L_1^{\vee}}(\psi_1 \times \Lambda^{\vee}).
\]
And as seen in the proof of proposition 6.4.6, the linear form $f_1^{L_1^{\vee}}(\psi_1 \times \Lambda^{\vee})$ is a unitary character on $G_1^{\vee}$. The familiar application of lemma 4.3.6, to (9.1.6), gives the vanishing:
\begin{eqnarray*}
 R^{G^*}_{\disc,\psi_+^{N_+}} &\equiv& 0, \,\ G^* \in \widetilde{\mathcal{E}}_{\simp}(N_+),\\
 \psi_1 \times \Lambda^{\vee} &\equiv& 0  .
\end{eqnarray*}  
Since the stable linear form defined by $\psi_1$ does not vanish, it follows that we have the vanishing of $\Lambda^{\vee}$ and hence $\Lambda$. In other words we have established the stable multiplicity formula for $\psi^N$ ({\it c.f.} proposition 6.1.4). 

With the stable multiplicity formula in hand, we can now apply lemma 5.7.6 to yield the spectral multiplicity formula for $\psi^N$ with respect to $G$. Similarly lemma 5.7.6 also yields the assertion that if $\psi^N$ belongs to $\xi^{\vee}_*\Psi_2(G^{\vee})$, then $\psi^N$ does not contribute to the discrete spectrum of $G$. Thus we finally see that only $\psi^N \in \xi_*\Psi_2(G)$ contribute to $R^G_{\disc}$. This thus establishes theorem 2.5.2 with respect to composite square-integrable parameters.

It remains to treat the case of simple parameters, i.e. that $\psi^N \in \widetilde{\Psi}_{\simp}(N)$. Following Arthur in \cite{A1}, the completion of the global induction argument in the case of simple parameter $\psi^N$ is to analyze the supplementary parameter:
\[
\psi_{++}^{N_{++}} := \psi^N \boxplus \psi^N \boxplus \psi^N, \,\ N_{++}=3N
\]
which we will carry out in the next subsection. In the rest of this subsection we first make some reductions.

We first consider the terms appearing in proposition 6.3.3. By the local intertwining relation (applied to each place $v$ of $F$), we have the identities
\begin{eqnarray}
f^{\prime}_{G_+}(\psi_+,s_{\psi_+}x_{+.1}) = f_{G_+}(\psi_+,x_{+,1}),  \,\ f \in \mathcal{H}(G_+) 
\end{eqnarray}
\begin{eqnarray}
& & \\
& &  (f^{\vee})^{M_+^{\vee}}(\psi^{\vee}_{M_+^{\vee}}) = (f^{\vee}_{G_+^{\vee}})^{\prime}(\psi_+^{\vee},x_{+,1}^{\vee})  = f^{\vee}_{G_+^{\vee}}(\psi_+^{\vee},x_{+,1}^{\vee}), \,\ f^{\vee} \in \mathcal{H}(G_+^{\vee}). \nonumber
\end{eqnarray}
On the other hand we also have:
\begin{eqnarray}
f^{M_+}(\psi_{M_+}) = (f^{\vee})^{M_+^{\vee}}(\psi^{\vee}_{M_+^{\vee}})
\end{eqnarray}
with $f$ and $f^{\vee}$ being the functions in a compatible family associated to $G_+$ and resp. $G^{\vee}_+$. This follows from the relation (3.1.5) between $f^{M_+}$ and $f^{M_+^{\vee}}$ (applied to the present situation); namely, if $G_+=(U_{E/F}(N_+),\xi_{\chi_1})$ and $G^{\vee}_+=(U_{E/F}(N_+),\xi_{\chi_2})$, then we have
\[
f^{\vee,M_+^{\vee}} = ((\chi_2/\chi_1)^N \circ \det) \cdot f^{M_+}
\]
together with the relation:
\[
\psi^{\vee}_{M_+^{\vee}} = ((\chi_1/\chi_2 )^N\circ \det) \otimes \psi_{M_+}.
\]

Thus by applying (9.1.7)-(9.1.9) to (6.3.19), we see that proposition 6.3.3 asserts the vanishing of:
\begin{eqnarray}
& & \sum_{G^* \in \widetilde{\mathcal{E}}_{\simp}(N_+)} \widetilde{\iota}(N_+,G^*) \tr R^{G^*}_{\disc,\psi_+^{N_+}}(f^*) \\
& &  + \frac{1}{8}(1-\delta_{\psi})( f^{M_+}(\psi_+^{M_+})-f_{G_+}(\psi_+,x_{+,1}) ) + \frac{1}{2} f^{L_+}(\Gamma \times \Lambda). \nonumber
\end{eqnarray}
Furthermore as in the proof of proposition 6.3.3, we have
\[
f^{L_+}(\Gamma \times \Lambda) = f_1^{L_1^{\vee}}(\Gamma \times \Lambda^{\vee})
\]
with $f_1$ being the function in the compatible family associated to $G_1^{\vee} = G \times G^{\vee} \in \widetilde{\mathcal{E}}_{\ellip}(N_+)$, and as seen in the proof of proposition 6.3.3, the term $f_1^{L_1^{\vee}}(\Gamma \times \Lambda^{\vee})$ is a linear combination with non-negative coefficients of irreducible representations on $G_1^{\vee}(\mathbf{A}_F)$. The same is true for the term
\begin{eqnarray}
\frac{1}{8}(1-\delta_{\psi})( f^{M_+}(\psi_+^{M_+})-f_{G_+}(\psi_+,x_{+,1}) ). 
\end{eqnarray}
Indeed the coefficient $(1-\delta_{\psi})$ is either $0$ or $2$, and we can write:
\begin{eqnarray*}
& & f^{M_+}(\psi_+^{M_+}) = f(\pi_1)+f(\pi_2) \\
& & f_{G_+}(\psi_+,x_{+,1}) = f(\pi_1)-f(\pi_2)
\end{eqnarray*}
for irreducible admissible representations $\pi_1,\pi_2$ on $G_+(\mathbf{A}_F)$; here $\pi_2$ can be zero, i.e. at this point we do not know whether 
\begin{eqnarray}
\mathcal{I}_{P_+}^{G_+} (\pi_{\psi_+^{M_+}})
\end{eqnarray}
is reducible (but in fact, after the conclusion of the induction argument in the next subsection, we can conclude that (9.1.12) is reducible, {\it c.f.} the discussion in section 3.4). In any case we conclude that (9.1.11) is a linear combination with non-negative coefficients of irreducible admissible representations on $G_+(\mathbf{A}_F)$.

We thus conclude that (9.1.10) can be written in a form such that lemma 4.3.6 applies. We thus obtain the vanishing of the following quantities:
\begin{eqnarray}
& & R^{G^*}_{\disc,\psi_+^{N_+}}, \,\ G^* \in \widetilde{\mathcal{E}}_{\simp}(N_+) \\
& & \,\  \Gamma \times \Lambda^{\vee} \\
& & (1-\delta_{\psi})( f^{M_+}(\psi_+^{M_+})-f_{G_+}(\psi_+,x_{+,1}) ), \,\ f \in \mathcal{H}(G_+).
\end{eqnarray}

Now suppose that $\psi \in \Phi_{\simp}(G)$ is a simple generic parameter of $G$. Recall that we are still following the definition of the set $\Phi_{\simp}(G)$ of simple generic parameters of $G$ as in section 6.1, and in particular the stable linear form
\[
S^G_{\disc,\psi^N}(f) = f^G(\Gamma)
\]
does not vanish identically. Thus by (9.1.14) we have the vanishing of $\Lambda^{\vee}$ and hence $\Lambda$. In other words the stable multiplicity formula holds for $\psi^N$ with respect to $G$, and also with respect to $G^{\vee}$, which amounts to the vanishing of $S^{G^{\vee}}_{\disc,\psi^N}$. Thus the twisted endoscopic datum $G=(G,\xi) \in \widetilde{\mathcal{E}}_{\simp}(N)$ such that $\psi^N \in \xi_* \Phi_{\simp}(G)$ is uniquely determined by $\psi^N$. In other words (6.1.4) is a disjoint union.

Recall that we have 
\[
\tr R^G_{\disc,\psi^N} \equiv S^G_{\disc,\psi^N}, \,\ \tr R^{G^{\vee}}_{\disc,\psi^N} \equiv S^{G^{\vee}}_{\disc,\psi^N}
\]
since $\psi^N$ is simple. In particular if $\psi^N \in \xi_* \Phi_{\simp}(G)$, then $R^{G^{\vee}}_{\disc,\psi^N} \equiv 0$. Now we already know the vanishing of $R^{G^*}_{\disc,\psi^N}$ for $G^* \in \widetilde{\mathcal{E}}_{\ellip}(N) \smallsetminus \widetilde{\mathcal{E}}_{\simp}(N)$. We thus finally conclude the proof of the seed theorem 2.4.2 (which concerns only simple generic parameters). In particular we resolve with the original definition of the set of simple generic parameters of $G$ given as in section 2.4 based on the seed theorem 2.4.2. The only remaining assertion to be established for the simple generic parameter $\psi^N$ is theorem 2.5.4(a), in other words the assertion $\delta_{\psi}=1$.

The second seed theorem 2.4.10 (which again concerns only simple generic parameters) also follows. Indeed if $\phi^N \in \xi_* \Phi_{\simp}(G)$, then by the stable multiplicity formula we just established we have $S_{\disc,\phi^N}^{G^{\vee}} \equiv 0$. Hence (9.1.2) becomes
\[
\widetilde{f}_N(\phi^N) = \widetilde{\iota}(N,G) \widehat{S}^G_{\disc,\phi^N}(\widetilde{f}^G)
\] 
and then the argument given in the last two paragraphs of the proof of proposition 9.1.1 shows that we must have $\phi_v^N \in \xi_{v,*}\Phi(G_v)$ for every place $v$ of $F$.

With the proof of the two seed theorems 2.4.2 and 2.4.10, the proof of theorem 5.2.1 then just follows from the corresponding results in the local situation which we established in section 7 and 8.

Now if on the other hand $\psi^N$ is a non-generic parameter, then the only assertion that remains to be established is the stable multiplicity formula for $\psi^N$, i.e. the vanishing of $\Lambda$ (we may assume that $N$ is even in this case, since when $N$ is odd the stable multiplicity formula for $\psi^N$ is already established in proposition 6.1.3).

We will complete the proof in the next subsection, by using the results we obtained above on the vanishing of (9.1.13)-(9.1.15), together with the analysis of the supplementary parameter $\psi_{++}^{N_{++}}$.

\subsection{Completion of induction arguments, part II}

We now complete the remaining portion of the induction argument. 

Recall that there are two cases depending on whether $\psi^N \in \widetilde{\Psi}_{\simp}(N)$ is generic or not. In the case that $\psi^N$ is generic, we need to establish that $\delta_{\psi}=1$, while in the case where $\psi^N$ is non-generic we need to establish the vanishing of $\Lambda$. 

We argue by contradiction. Namely in the case that $\psi^N$ is generic we suppose that $\delta_{\psi}=-1$. Then by (9.1.15) we have:
\begin{eqnarray}
f^{M_+}(\psi_+^{M_+}) = f_{G_+}(\psi_+,x_{+,1}), \,\ f \in \mathcal{H}(G_+).
\end{eqnarray}

Similarly in the case that $\psi^N$ is non-generic, we suppose that the linear form $\Lambda$ (or equivalently $\Lambda^{\vee}$) does not vanish; recall that in this case we may assume that $N$ is even. Then by (9.1.14) we have
\begin{eqnarray}
\Gamma \equiv 0.
\end{eqnarray}

We will show that both (9.2.1) and (9.2.2) would lead to a contradiction, which thus complete the global induction argument. 

We first consider the generic case. In the previous subsection we have already established the stable multiplicity formula for $\psi^N$ (or equivalently the vanishing of $\Lambda$); thus 
\[
f_1^G(\psi) = f_1^G(\Gamma), \,\ f_1 \in \mathcal{H}(G).
\]
Hence we have:
\begin{eqnarray*}
& & S^G_{\disc,\psi^N}(f_1) = f_1^G(\psi), \,\ f_1 \in \mathcal{H}(G) \\
& & S^{G^{\vee}}_{\disc,\psi^N}(f_1^{\vee}) =0, \,\ f_1^{\vee} \in \mathcal{H}(G^{\vee}).
\end{eqnarray*}
For the stable distributions associated to the supplementary parameter $\psi_+^{N_+}$, we apply the identities from corollary 6.3.4. For instance for (6.3.26), we note that:
\begin{eqnarray}
f_2^{G \times G}(\Gamma \times \Gamma) &=& f_2^{G \times G}(\psi \times \psi) \\
&=& (f_2)^{\prime}_{G_+}(\psi_+,x_{+,1}) \nonumber \\
&=& (f_2)_{G_+}(\psi_+,x_{+,1}) \nonumber \\
&=& (f_2)^{M_+}(\psi_{M_+}) , \,\ f_2 \in \mathcal{H}(G_+). \nonumber
\end{eqnarray}
(here the third equality is by the local intertwining relation, and the fourth equality is by (9.2.1)). Hence (6.3.26) simplifies to (using the hypothesis $\delta_{\psi}=-1$ and the vanishing of (9.1.13)):
\begin{eqnarray}
S^{G_+}_{\disc,\psi_+^{N_+}}(f_2) = -\frac{1}{2} (f_2)^{M_+}(\psi_{M_+}), \,\ f_2 \in \mathcal{H}(G_+).
\end{eqnarray}

Similarly for (6.3.27); we have
\begin{eqnarray}
 (f_2^{\vee})_{G_+^{\vee}}(\psi_+^{\vee},x_{+,1}^{\vee})   &=&   (f_2^{\vee})^{\prime}_{G_+^{\vee}}(\psi_+^{\vee},x_{+,1}^{\vee}) \\
&=&   (f_2^{\vee})^{M_+^{\vee}}(\psi_{M_+^{\vee}}) \nonumber 
\end{eqnarray}
(here the first equality is again by the local intertwining relation). Thus (6.3.27) simplifies to:
\begin{eqnarray}
S^{G^{\vee}}_{\disc,\psi^N}(f_2^{\vee}) = \frac{1}{4} (f_2^{\vee})^{M_+^{\vee}}(\psi_{M_+^{\vee}}).
\end{eqnarray}

To summarize: in the generic case and under our hypothesis that $\delta_{\psi}=-1$ we have the following equalities: here we are taking $f_1$ and $f_1^{\vee}$ to be functions associated to $G$ and resp. $G^{\vee}$ in a compatible family, and similarly $f_2$ and $f_2^{\vee}$ to be functions associated to $G_+$ and $G_+^{\vee}$ in a compatible family: 
\begin{eqnarray}
& & S^G_{\disc,\psi^N}(f_1) = f_1^G(\psi) \\
& & S^{G^{\vee}}_{\disc,\psi^N}(f_1^{\vee}) =0 \nonumber \\
& & S^{G_+}_{\disc,\psi_+^{N_+}}(f_2) = -\frac{1}{2} (f_2)^{M_+}(\psi_{M_+}) \nonumber \\
& & S^{G^{\vee}_+}_{\disc,\psi^N}(f_2^{\vee}) = \frac{1}{4} (f_2)^{M_+}(\psi_{M_+}) \nonumber \\
& & \,\ \,\ \,\ \,\ \,\ \,\ \,\ \,\ \,\ \,\ \,\  =\frac{1}{4} (f_2^{\vee})^{M_+^{\vee}}(\psi^{\vee}_{M_+^{\vee}}). \nonumber
\end{eqnarray}
Here the fourth equation follows from (9.2.6) and (9.1.9). 

We next consider the case that $\psi^N$ is non-generic. Under the hypothesis that $\Lambda$ (or equivalently $\Lambda^{\vee}$) does not vanish, we have $\Gamma \equiv 0$ as in (9.2.2), and hence
\[
f_1^G(\psi) = f_1^L(\Lambda). 
\]
Analogous to (9.2.7) is a similar set of identities obtained as follows.

First from the vanishing of $\Gamma$ we have by proposition 6.1.3, the following (here as above $f_1$ and $f_1^{\vee}$ are functions associated to $G$ and $G^{\vee}$ in a compatible family):
\begin{eqnarray*}
S^G_{\disc,\psi^N}(f_1) = f_1^G(\psi) - f_1^L(\Lambda) =0, 
\end{eqnarray*}
\begin{eqnarray}
S^{G^{\vee}}_{\disc,\psi^N}(f^{\vee}_1) &=& (f_1^{\vee})^{L^{\vee}}(\Lambda^{\vee}) \\
&=& f_1^L(\Lambda) = f_1^G(\psi). \nonumber
\end{eqnarray}
(See proposition 6.1.4 for the second equality of (9.2.8).)

For the stable distributions associated to the supplementary parameter $\psi_+^{N_+}$ we need the following:
\begin{proposition}
Suppose we have $\Gamma \equiv 0$. Then for any place $v$ of $F$ that does not split in $E$, the localization $\psi_v^N$ of the parameter $\psi^N$ factors through both the $L$-embeddings $\leftexp{L}{L_v} \hookrightarrow \leftexp{L}{G_{E_v/F_v}(N)}$ and $\leftexp{L}{L^{\vee}_v} \hookrightarrow \leftexp{L}{G_{E_v/F_v}(N)}$. 
\end{proposition}
In order not to interrupt the line of reasoning we relegate the proof of proposition 9.2.1 to the appendix.  

We then draw the following consequence from proposition 9.2.1. For any place $v$ of $F$ that does not split in $E$, the localization $\psi_{+,v}$ of the supplementary parameter $\psi_+$ factors through the image of the $L$-embedding
\begin{eqnarray}
\leftexp{L}{(L_v \times L_v)} \hookrightarrow \leftexp{L}{(G_{E_v/F_v}(N) \times G_{E_v/F_v}(N))}.
\end{eqnarray}
The centralizer of the image of the $L$-embedding (9.2.9) in the dual group $\widehat{G}_{E_v/F_v}(N) \times \widehat{G}_{E_v/F_v}(N)$ can be identified with the connected group
\[
\GL(2,\mathbf{C}) \times \GL(2,\mathbf{C})
\]
and under this identification, the global centralizer $S_{\psi_+}$ can be identified as the diagonal image of $O(2,\mathbf{C})$ in $\GL(2,\mathbf{C}) \times\GL(2,\mathbf{C})$. It thus follows that for any $x_+ \in \mathcal{S}_{\psi_+}$ the image $x_{+,v}$ of $x_{+}$ in $\mathcal{S}_{\psi_{+,v}}$ under the localization map is trivial This last statement of course remains trivially true if $v$ splits in $E$. From this it follows that
\[
(f_2)_{G_+}^{\prime}(\psi_+,x_+) =(f_2)^{M_+}(\psi_{M_+}), \,\ f_2 \in \mathcal{H}(G_+).
\]
Hence on combining with the local intertwining relation we have for any $x_+ \in \mathcal{S}_{\psi_+}$:
\begin{eqnarray}
(f_2)_{G_+}(\psi_+,x_+) = (f_2)^{M_+}(\psi_{M_+}), \,\ f_2 \in \mathcal{H}(G_+).
\end{eqnarray}

Now as in the generic case we let $f_2$ and $f_2^{\vee}$ be functions in a compatible family associated to $G_+$ and resp. $G_+^{\vee}$. Thus by (9.2.10) and the vanishing of (9.1.13) we obtain from (6.3.26) the following (remember that since $\psi^N$ is non-generic we have $\delta_{\psi}=+1$):
\begin{eqnarray}
S^{G_+}_{\disc,\psi_+^{N_+}}(f_2) = \frac{1}{4}(f_2)^{M_+}(\psi_{M_+}).
\end{eqnarray}
Similarly for the terms occuring in (6.3.27) we have:
\begin{eqnarray}
f_2^{L \times L}(\Lambda \times \Lambda) &=& f_2^{G \times G}(\psi \times \psi) \\
&=& (f_2)^{\prime}_{G_+}(\psi_+,x_{+,1}) \nonumber \\
&=& (f_2)^{M_+}(\psi_{M_+}) \nonumber
\end{eqnarray}
(the last equality by (9.2.10)). Hence (6.3.27) gives:
\begin{eqnarray}
S_{\disc,\psi_+^{N_+}}^{G^{\vee}_+}(f_2^{\vee})= -\frac{1}{2} (f_2)^{M_+}(\psi_{M_+}).
\end{eqnarray}

Thus to summarize, in the non-generic case, under the hypothesis that $\Lambda $ does not vanish we have the identities parallel to (9.2.7):
\begin{eqnarray}
& & S^G_{\disc,\psi^N}(f_1) = 0 \\
& & S^{G^{\vee}}_{\disc,\psi^N}(f_1^{\vee}) = f_1(\psi) \nonumber \\
& & S^{G_+}_{\disc,\psi_+^{N_+}}(f_2) =  \frac{1}{4} (f_2)^{M_+}(\psi_{M_+}) \nonumber \\
& & S^{G^{\vee}_+}_{\disc,\psi^N}(f_2^{\vee}) = -\frac{1}{2} (f_2)^{M_+}(\psi_{M_+}) \nonumber \\
& & \,\ \,\ \,\ \,\ \,\ \,\ \,\ \,\ \,\ \,\ \,\  =-\frac{1}{2} (f_2^{\vee})^{M_+^{\vee}}(\psi^{\vee}_{M_+^{\vee}}). \nonumber
\end{eqnarray}

For comparison, the corresponding set of identities in the expected setting (namely $\delta_{\psi}=1$ and $\Lambda \equiv 0$, in both the generic and non-generic case) is given by:
\begin{eqnarray}
& & S^G_{\disc,\psi^N}(f_1) = f_1^G(\psi) \\
& & S^{G^{\vee}}_{\disc,\psi^N}(f_1^{\vee}) = 0 \nonumber \\
& & S^{G_+}_{\disc,\psi_+^{N_+}}(f_2) = 0 \nonumber \\
& & S^{G^{\vee}_+}_{\disc,\psi^N}(f_2^{\vee}) = -\frac{1}{4} (f_2)^{M_+}(\psi_{M_+}) \nonumber \\
& & \,\ \,\ \,\ \,\ \,\ \,\ \,\ \,\ \,\ \,\ \,\  =-\frac{1}{4} (f_2^{\vee})^{M_+^{\vee}}(\psi^{\vee}_{M_+^{\vee}}). \nonumber
\end{eqnarray}
which can be established by the similar reasoning as above.

In order to show that (9.2.7) or (9.2.14) leads to a contradiction, we follow the crucial idea of Arthur \cite{A1} by considering the supplementary parameter:
\[
\psi_{++}^{N_{++}} :=\psi^N \boxplus \psi^N \boxplus \psi^N, \,\ N_{++}=3N.
\] 
Denote by $G_{++}$ and $G_{++}^{\vee}$ the elements in $\widetilde{\mathcal{E}}_{\simp}(N_{++})$ with the same parity (resp. opposite parity) as $G \in \widetilde{\mathcal{E}}_{\simp}(N)$. Then $\psi_{++}^{N_{++}}$ defines the parameter $\psi_{++}$ of $G_{++}$. 

Similarly denote by $M_{++} = G \times M_+ $ (resp. $M_{++}^{\vee} =G^{\vee} \times  M_+^{\vee}$ the Levi sub-datum of $G_{++}$ (resp. $G_{++}^{\vee}$) whose underlying group is $\cong U_{E/F}(N) \times G_{E/F}(N) $. Then the parameter $\psi_{++}^{N_{++}}$ defines the parameter
\[
\psi \times \psi_{M_+} 
\]
with respect to $M_{++}$. We can now state the final result to be established by consideration of the standard model of trace formula comparison, which would allow us to obtain the desired contradiction:

\begin{proposition}
Under the assumption that $\delta_{\psi} \neq 1$ or $\Lambda$ does not vanish, we have the identity (as usual for a compatible family of functions with respect to $\widetilde{\mathcal{E}}_{\ellip}(N_{++})$):
\begin{eqnarray}
& & \\
& & \sum_{G^* \in \widetilde{\mathcal{E}}_{\simp}(N_{++})} \widetilde{\iota}(N_{++},G^*) \tr R^{G^*}_{\disc,\psi_{++}^{N_{++}}}(f^*) +\frac{1}{2} f^{M_{++}}(\psi \times \psi_{M_+}) \equiv 0. \nonumber
\end{eqnarray}
(Here $f$ is the function in the compatible family associated to $G_{++}$.)
\end{proposition}
\begin{proof}
We follow Arthur's method in the proof of lemma 8.2.2 in \cite{A1}. Thus we apply proposition 6.2.3 to the supplementary parameter $\psi_{++}^{N_{++}}$, which is applicable since condition (6.2.19) is valid. Then we write (6.2.20) as the following assertion: that for a compatible family of functions we have the equality between (9.2.17):
\begin{eqnarray}
\sum_{G^* \in \widetilde{\mathcal{E}}_{\simp}(N_{++})} \widetilde{\iota}(N_{++},G^*) \tr R^{G^*}_{\disc,\psi_{++}^{N_{++}}}(f^*) 
\end{eqnarray}
and the sum of (9.2.18) and (9.2.19):
\begin{eqnarray}
& &\\
& &
\sum_{G^* \in \widetilde{\mathcal{E}}_{\simp}(N_{++})} \widetilde{\iota}(N_{++},G^*) \big(\tr R^{G^*}_{\disc,\psi_{++}^{N_{++}}}(f^*)  - \leftexp{0}{S^{G^*}_{\disc,\psi_{++}^{N_{++}}} (f^*)} \big) \nonumber
\end{eqnarray}
\begin{eqnarray}
-\sum_{G^* \in \widetilde{\mathcal{E}}^0_{\ellip}(N_{++})} \widetilde{\iota}(N_{++},G^*) \leftexp{0}{S^{G^*}_{\disc, \psi_{++}^{N_{++}}}(f^*) }
\end{eqnarray}
here $\widetilde{\mathcal{E}}_{\ellip}^0(N_{++}) := \widetilde{\mathcal{E}}_{\ellip}(N_{++}) \smallsetminus \widetilde{\mathcal{E}}_{\simp}(N_{++})$. 

We first consider the case that $\psi^N$ is non-generic, and thus we are under the hypothesis that $\Lambda$ does not vanish. Now following the reasoning as in the proof of proposition 6.3.1 and 6.3.3, we see that the only possible non-vanishing contributions in (9.2.18) and (9.2.19) can be enlisted as follows.  

For (9.2.19), the only possible non-vanishing terms comes from $G^*=G \times G_+^{\vee}$ or $G^{\vee} \times G_+$. For $G^*=G \times G_+^{\vee}$, the actual value of 
\begin{eqnarray}
S^{G \times G_+^{\vee}}_{\disc,\psi_{++}^{N_{++}}}(f^*) = S^G_{\disc,\psi^N}(f_1) \cdot S^{G_+^{\vee}}_{\disc,\psi_+^{N_+}}(f_2^{\vee})
\end{eqnarray}
for $f^* =f_1 \times f_2^{\vee}$, vanishes, according to the first equation of (9.2.14). While its expected value is, according to the first and fourth equations of (9.2.15):
\begin{eqnarray}
-\frac{1}{4} f_1^G(\psi)  \cdot (f_2^{\vee})^{M_+^{\vee}}(\psi^{\vee}_{M_+^{\vee}}) = -\frac{1}{4} f^{M_{++}}(\psi \times \psi_{M_+}).
\end{eqnarray}
Hence we have
\begin{eqnarray}
\leftexp{0}{S^{G \times G_+^{\vee}}_{\disc,\psi_{++}^{N_{++}}}}(f^*) = \frac{1}{4} f^{M_{++}}(\psi \times \psi_{M_+}).
\end{eqnarray}
Similarly for $G^*=G^{\vee} \times G_+$, we have the actual value 
\begin{eqnarray}
S^{G^{\vee} \times G_+}_{\disc, \psi_{\psi_{++}^{N_{++}}}}(f^*) =\frac{1}{4} f^{M_{++}}(\psi \times \psi_{M_+})
\end{eqnarray}
while its expected value vanishes, by the second or third equation of (9.2.15). Hence we have
\begin{eqnarray}
\leftexp{0}{S^{G^{\vee} \times G_+}_{\disc, \psi_{\psi_{++}^{N_{++}}}}(f^*)   } =\frac{1}{4} f^{M_{++}}(\psi \times \psi_{M_+}).
\end{eqnarray}
Since we have $\widetilde{\iota}(N_{++},G \times G_+^{\vee})= \widetilde{\iota}(N_{++},G^{\vee} \times G_+)=1/4$, we see that 
\begin{eqnarray*}
(9.2.19)= -\frac{1}{8} f^{M_{++}}(\psi \times \psi_{M_+}).
\end{eqnarray*}

For (9.2.18) it can be analyzed as follows. For $G^* \in \widetilde{\mathcal{E}}_{\simp}(N_{++})$, we write the spectral and endoscopic expansions for $I^{G^*}_{\disc,\psi_{++}^{N_{++}}}$ as:
\begin{eqnarray*}
& & I^{G^*}_{\disc,\psi_{++}^{N_{++}}}(f^*) = \tr R^{G^*}_{\disc,\psi_{++}^{N_{++}}}(f^*) + (I^{G^*}_{\spec})^{\prime}(f^*) \\
& & I^{G^*}_{\disc,\psi_{++}^{N_{++}}}(f^*) = S^{G^*}_{\disc,\psi_{++}^{N_{++}}}(f^*) + (I^{G^*}_{\scopy})^{\prime}(f^*).
\end{eqnarray*}
Noting that the expected value of $R^{G^*}_{\disc,\psi_{++}^{N_{++}}}$ vanishes, we have the identity:
\begin{eqnarray}
& & \tr R^{G^*}_{\disc,\psi_{++}^{N_{++}}}(f^*) - \leftexp{0}{S^{G^*}_{\disc,\psi_{++}^{N_{++}}}(f^*)} \\
&=&  \leftexp{0}{(I^{G^*}_{\scopy})^{\prime}(f^*)}  - \leftexp{0}{(I^{G^*}_{\spec})^{\prime}(f^*)}   \nonumber
\end{eqnarray}
where $ \leftexp{0}{(I^{G^*}_{\spec})^{\prime}(f^*)}$ denotes the difference between the actual and expected value of $(I^{G^*}_{\spec})^{\prime}(f^*)$ (and similarly for $(I^{G^*}_{\scopy})^{\prime}(f^*)$). 

We now analyze the term $ \leftexp{0}{(I^{G^*}_{\scopy})^{\prime}}$. For $G^* =G_{++}$, the only possible non-vanishing contribution can be seen to come only from (similar to proof of proposition 6.3.1 and 6.3.3) $G^{\prime} = G \times G_+ \in \widetilde{\mathcal{E}}_{\ellip}(G_{++})$. For the distribution 
\[
S^{G \times G_+}_{\disc,\psi_{++}^{N_{++}}}(f^{\prime}) = S^{G}_{\disc,\psi^N}(f_1) \cdot S^{G_+}_{\disc,\psi_+^{N_+}}(f_2), \,\ f^{\prime}=f_1 \times f_2
\] 
both its actual and expected values vanishes, by the first equation of (9.2.14) and the third equation of (9.2.15). Thus 
\begin{eqnarray*}
\leftexp{0}{S^{G \times G_+}_{\disc,\psi_{++}^{N_{++}}}} \equiv 0.
\end{eqnarray*}
hence
\begin{eqnarray}
\leftexp{0}{ (I^{G_{++}}_{\scopy})^{\prime} } \equiv 0.
\end{eqnarray}
Similarly for $G^*=G_{++}^{\vee}$, the only possible non-vanishing contribution comes from $G^{\prime} =G^{\vee} \times G_+^{\vee} \in \widetilde{\mathcal{E}}_{\ellip}(G_{++}^{\vee})$. By the second and fourth equation of (9.2.14), we have the actual value of the distribution (for $f^{\prime}=f_1^{\vee} \times f_2^{\vee}$):
\begin{eqnarray}
S^{G \times G_+}_{\disc,\psi_{++}^{N_{++}}}(f^{\prime}) &=& S^G_{\disc,\psi^N}(f_1^{\vee}) \cdot S^{G_+}_{\disc,\psi_+^{N_+}}(f_2^{\vee}) \nonumber \\
&=& -\frac{1}{2} f^{M_{++}}(\psi \times \psi_{M_+}) \nonumber
\end{eqnarray} 
while its expected value vanishes, by the second equation of (9.2.15). Thus we have
\begin{eqnarray}
\leftexp{0}{S^{G^{\vee} \times G^{\vee}_+}_{\disc,\psi_{++}^{N_{++}}}}(f^*) = -\frac{1}{2} f^{M_{++}}(\psi \times \psi_{M_+}).
\end{eqnarray}
Hence 
\begin{eqnarray}
& & \leftexp{0}{(I^{G^{\vee}_{++}}_{\scopy})^{\prime}(f^*)} \\
&=&  \iota(G^{\vee}_{++},G^{\vee} \times G_+^{\vee})  \cdot ( -\frac{1}{2} f^{M_{++}}(\psi \times \psi_{M_+})) \nonumber\\
&=& \frac{1}{2}   \cdot ( -\frac{1}{2} f^{M_{++}}(\psi \times \psi_{M_+})) \nonumber \\
&= & - \frac{1}{4}  f^{M_{++}}(\psi \times \psi_{M_+}). \nonumber
\end{eqnarray}

For the distribution $\leftexp{0}{(I^{G^*}_{\spec})^{\prime}}$, we have to analyze the term:
\begin{eqnarray}
\tr M_{P^*,\psi_{++}^{N_{++}}}(w^*) \mathcal{I}_{P^*,\psi_{++}^{N_{++}}}(f^*).
\end{eqnarray}
The only possible non-vanishing contribution to $\leftexp{0}{(I^{G^*}_{\spec})^{\prime}}$ comes from $M^*=M_{++}$ when $G^*=G_{++}$, and similarly it comes from $M^*=M_{++}^{\vee}$ when $G^*=G_{++}^{\vee}$. In both cases $w^*$ is the unique element of $W(M^*)_{\reg}$.

Consider first the case for $G_{++}$. Then with $M^*=M_{++}$, we have the vanishing of $R^{M_{++}}_{\disc, \psi \times \psi_{M_+}}$, by the first equation of (9.2.14). Thus the actual value of (9.2.29) vanishes. To compute its expected value, we note that we have
\[
\overline{S}_{\psi_{++}} = SO(3,\mathbf{C})
\]
which is a connected group. Thus the distribution (9.2.29) is equal to the product of the expected value of the global normalizing factor, which is equal to $(-1)$ in this case (coming from the non-trivial Weyl element of $SO(3,\mathbf{C})$), and the distribution:
\[
f_{G_{++}}(\psi_{++},x_{++}) = f^{M_{++}}(\psi \times \psi_{M_+})
\] 
with equality being due to the intertwining relation, together with the fact that $\overline{S}_{\psi_{++}}$ is connected ($x_{++}$ being the only element in the trivial group $\mathcal{S}_{\psi_{++}}$). Hence:
\begin{eqnarray}
& &  \leftexp{0}{I_{\spec}^{G_{++}}(f) } \\
& =& |W(M_{++})|^{-1} |\det(w^* -1)|^{-1} (-1)\cdot (-1) f^{M_{++}}(\psi \times \psi_{M_+}) \nonumber \\
&=&  \frac{1}{2} \cdot \frac{1}{2} \cdot  f^{M_{++}}(\psi \times \psi_{M_+}) \nonumber \\
&=& \frac{1}{4}  f^{M_{++}}(\psi \times \psi_{M_+}). \nonumber
\end{eqnarray}
In particular substituting (9.2.26) and (9.2.30) in (9.2.25) we have
\begin{eqnarray}
& & \tr R^{G_{++}}_{\disc,\psi_{++}^{N_{++}}}(f) - \leftexp{0}{S^{G_{++}}_{\disc,\psi_{++}^{N_{++}}}(f)} \\
&=&  \leftexp{0}{(I^{G_{++}}_{\scopy})^{\prime}(f)}  - \leftexp{0}{(I^{G_{++}}_{\spec})^{\prime}(f)}   \nonumber \\
&=& -\frac{1}{4} f^{M_{++}}(\psi \times \psi_{M_+}). \nonumber
\end{eqnarray}

It remains to treat (9.2.29) for $G^*=G_{++}^{\vee}$. By the second equation of (9.2.15), the expected value of the distribution
\[
\tr R^{M_{++}^{\vee}}_{\disc,\psi^N \times \psi^N}
\]
vanishes, hence we have the vanishing of the expected value of (9.2.29). Thus it suffices to compute its actual value. 

We claim that the global intertwining operator occuring in (9.2.29) is equal to identity. This is seen as follows. First we compute the (actual value of the) global normalizing factor. Applying the discussion in section 5.8, the global normalizing factor is equal to 
\begin{eqnarray}
(-1)^{\ord_{s=1} L(s,\psi^N \times \psi^N)} \times (-1)^{\ord_{s=1} L(s,\psi^N,\Asai^{(-1)^{N-1}\kappa})}
\end{eqnarray}
here $\kappa$ is the sign such that $G=(G,\xi_{\chi})$ with $\chi \in \mathcal{Z}_E^{\kappa}$. We know that the Rankin-Selberg $L$-function $L(s,\psi^N \times \psi^N)$ has a pole of odd order at $s=1$ since $\psi^N$ is conjugate self-dual; one the other hand, since $\psi^N$ is non-generic, we have the validity of theorem 2.5.4(a) for the simple generic factor of $\psi^N$ (by induction hypothesis), and thus the Langlands-Shahidi $L$-function $L(s,\psi^N,\Asai^{(-1)^{N-1} \kappa})$ also has a pole of odd order at $s=1$. Thus (9.2.32) is equal to $(-1) \times (-1)=+1$.

Next we show that the normalized intertwining operator in (9.2.29) is equal to the identity. For this it suffices to check locally at every places $v$ of $F$. By proposition 9.2.1, if $v$ does not split in $E$, the localization $\psi_v^N$ of the global parameter $\psi^N$ factors through $\leftexp{L}{L^{\vee}_v} \hookrightarrow \leftexp{L}{G_{E_v/F_v}(N)}$, in particular factors through $\leftexp{L}{G_v^{\vee}} \hookrightarrow \leftexp{L}{G_{E_v/F_v}(N)}$. Hence the local parameter
\[
\psi^{N_{++}}_{++,v} = \psi^N_v \oplus \psi^N_v \oplus \psi^N_v
\]
defines a local parameter $\psi^{\vee}_{++,v}$ of $G^{\vee}_{++,v}$. The same conclusion holds trivially if $v$ splits in $E$.

Now if $v$ does not split in $E$, then noting that 
\[
\Cent(\leftexp{L}{(G^{\vee}_v \times G^{\vee}_v \times G^{\vee}_v})    ,\widehat{G}^{\vee}_{++})/Z(\widehat{G}_{++}^{\vee})^{\Gamma_{F_v}} \cong SO(3,\mathbf{C})
\]
a connected group, we see that the image of $w^*$ in the local $R$-group $R_{\psi^{\vee}_{++,v}}$ is trivial. Hence the normalized local intertwining operator defined by the local image of $w^*$ is trivial. When $v$ splits in $E$, then the same conclusion holds since the local $\mathcal{S}$-group $\mathcal{S}_{\psi^{\vee}_{++,v}}$ is already trivial.

From the triviality of the global intertwining operator in (9.2.29), we obtain:
\begin{eqnarray*}
(9.2.29) &=& \prod_v (f^{\vee}_v)^{G_{++,v}^{\vee}}(\psi^{\vee}_{++,v}) \\
&=&  \prod_v (f_v)^{G_{++,v}}(\psi_{++,v}) \\
&=& f^{G_{++}}(\psi_{++}) \\
&=& f^{M_{++}}(\psi \times \psi_{M_+}).
\end{eqnarray*}
And hence
\begin{eqnarray}
& & \leftexp{0}{(I^{G_{++}^{\vee}}_{\spec})^{\prime}  }(f^{\vee}) \\
&=& |W(M_{++}^{\vee})|^{-1}|\det (w^*-1)|^{-1} f^{M_{++}}(\psi \times \psi_{M_+})  \nonumber \\
&=& \frac{1}{4}  f^{M_{++}}(\psi \times \psi_{M_+}). \nonumber
\end{eqnarray}

Substituting (9.2.28) and (9.2.33) in (9.2.25) we obtain 
\begin{eqnarray}
& & \tr R^{G_{++}^{\vee}}_{\disc,\psi_{++}^{N_{++}}}(f^{\vee}) - \leftexp{0}{S^{G_{++}^{\vee}}_{\disc,\psi_{++}^{N_{++}}}(f^{\vee})} \\
&=&  \leftexp{0}{(I^{G^{\vee}_{++}}_{\scopy})^{\prime}(f^{\vee})}  - \leftexp{0}{(I^{G^{\vee}_{++}}_{\spec})^{\prime}(f^{\vee})}   \nonumber \\
&=& -\frac{1}{2} f^{M_{++}}(\psi \times \psi_{M_+}). \nonumber
\end{eqnarray}
Thus we obtain:
\begin{eqnarray*}
 (9.2.18) & =&  \widetilde{\iota}(N_{++},G_{++})   (-\frac{1}{4}  f^{M_{++}}(\psi \times \psi_{M_+}) )  \\ 
& & \,\ \,\ \,\  + \widetilde{\iota}(N_{++},G_{++}^{\vee}) (-\frac{1}{2}  f^{M_{++}}(\psi \times \psi_{M_+}) ) \nonumber \\
&=&  -\frac{3}{8}   f^{M_{++}}(\psi \times \psi_{M_+}) . \nonumber
\end{eqnarray*}
To conclude:
\begin{eqnarray*}
(9.2.17)&=& (9.2.18)+(9.2.19)\\
&=& -\frac{3}{8} f^{M_{++}}(\psi \times \psi_{M_+}) -\frac{1}{8} f^{M_{++}}(\psi \times \psi_{M_+})\\
&=& -\frac{1}{2} f^{M_{++}}(\psi \times \psi_{M_+})
\end{eqnarray*}
thus establishing (9.2.16), in the case where $\psi^N$ is non-generic (and under our hypothesis that $\Lambda$ does not vanish).

We now treat the case where $\psi^N$ is a simple generic parameter. Thus we assume that $\delta_{\psi}=-1$. The derivation of (9.2.16) is similar to the non-generic case, and so we will be brief.

For (9.2.19), we have the following: the actual value of the distribution
\[
S^{G \times G_+^{\vee}}_{\disc,\psi_{++}^{N_{++}}}(f^*)
\]
is given by: 
\[
\frac{1}{4} f^{M_{++}}(\psi \times \psi_{M_+})
\]
by the first and fourth equation of (9.2.7), while its expected value is given by: 
\[
-\frac{1}{4} f^{M_{++}}(\psi \times \psi_{M_+})
\]
by the first and fourth equation of (9.2.15). Hence
\begin{eqnarray}
\leftexp{0}{S^{G \times G_+^{\vee}}_{\disc,\psi_{++}^{N_{++}}}}(f^*) = \frac{1}{2} f^{M_{++}}(\psi \times \psi_{M_+}).
\end{eqnarray}
For the distribution $S^{G^{\vee} \times G_+}_{\disc,\psi_{++}^{N_{++}}}$, its actual value vanishes by the second equation of (9.2.7), while its expected value vanishes by the second or third equation of (9.2.15). Thus 
\begin{eqnarray}
\leftexp{0}{S^{G^{\vee} \times G_+}_{\disc,\psi_{++}^{N_{++}}}} \equiv 0
\end{eqnarray}
and hence
\begin{eqnarray*}
(9.2.19) &=&- \widetilde{\iota}(N_{++},G \times G_+^{\vee}) \cdot \frac{1}{2} f^{M_{++}}(\psi \times \psi_{M_+}) \\
&=&  -\frac{1}{8} f^{M_{++}}(\psi \times \psi_{M_+}) .  \nonumber
\end{eqnarray*}
Next consider (9.2.25). First consider the term $\leftexp{0}{ (I_{\scopy}^{G^*})^{\prime}(f^*)}$. For $G^*=G_{++}$, we have the actual value:
\begin{eqnarray}
S^{G \times G_+}_{\disc,\psi_{++}^{N_{++}}}(f) = -\frac{1}{2} f^{M_{++}}(\psi \times \psi_{M_+}). 
\end{eqnarray}
by the first and third equation of (9.2.7), while its expected value vanishes, by the third equation of (9.2.15). Hence 
\begin{eqnarray}
\leftexp{0}{ (I_{\scopy}^{G_{++}})^{\prime}(f)} &=& \iota(G_{++},G \times G_+) \leftexp{0}{S^{G \times G_+}_{\disc,\psi_{++}^{N_{++}}}(f)}\\
&=&  -\frac{1}{4} f^{M_{++}}(\psi \times \psi_{M_+}).  \nonumber
\end{eqnarray}

We next compute $\leftexp{0}{ (I_{\spec}^{G_{++}})^{\prime}(f)}$, and thus need to analyze (9.2.29). Still for $G^*=G_{++}$, we have
\[
\overline{S}_{\psi_{++}} = SO(3,\mathbf{C})
\]
a connected group. Hence the normalized global intertwining operator occuring in (9.2.29) is trivial. Thus we have:
\begin{eqnarray*}
& & \mbox{actual value of } (9.2.29) \\
&=& \mbox{actual value of global normalizing factor } \times f^{G_{++}}(\psi_{++}) \\
&=& \mbox{actual value of global normalizing factor } \times f^{M_{++}}(\psi \times \psi_{M_+}).
\end{eqnarray*}
Similarly we have:
\begin{eqnarray*}
& & \mbox{expected value of } (9.2.29) \\
&=& \mbox{expected value of global normalizing factor } \times f^{G_{++}}(\psi_{++}) \\
&=& \mbox{expected value of global normalizing factor } \times f^{M_{++}}(\psi \times \psi_{M_+}).
\end{eqnarray*}

Now the actual value of the global normalizing factor occuring in (9.2.29) is given by 
\begin{eqnarray}
(-1)^{\ord_{s=1} L(s,\psi^N \times \psi^N)} \cdot (-1)^{\ord_{s=1} L(s,\psi^N,\Asai^{(-1)^N \kappa })}
\end{eqnarray}
here $\kappa$ is as before the sign such that $G=(G,\xi_{\chi})$ with $\chi \in \mathcal{Z}_E^{\kappa}$. Now $L(s,\psi^N \times \psi^N)$ has a simple pole at $s=1$ as usual, while under our hypothesis that $\delta_{\psi}=-1$, we have $L(s,\psi^N,\Asai^{(-1)^N \kappa })$ has a simple pole at $s=1$. Thus
\[
(9.2.39)=+1.
\]
On the other hand the expected value of the global normalizing factor is given by $-1$, as is easily seen from the fact that $\overline{S}_{\psi_{++}} = SO(3,\mathbf{C})$.

It thus follows that
\begin{eqnarray}
& & \\
& & \leftexp{0}{(I^{G_{++}}_{\spec})^{\prime}(f)} = |W(M_{++})|^{-1}|\det( w^* -1)|^{-1} \cdot 2 \cdot  f^{M_{++}}(\psi \times \psi_{M_+}) \nonumber \\
&=& \frac{1}{2} \cdot \frac{1}{2} \cdot 2 \cdot f^{M_{++}}(\psi \times \psi_{M_+}) \nonumber \\
&=& \frac{1}{2} f^{M_{++}}(\psi \times \psi_{M_+}).\nonumber
\end{eqnarray}

We thus obtain 
\begin{eqnarray}
& & \tr R^{G_{++}}_{\disc,\psi_{++}^{N_{++}}}(f) - \leftexp{0}{S^{G_{++}}_{\disc,\psi_{++}^{N_{++}}}   }(f) \\
&=&   \leftexp{0}{ (I_{\scopy}^{G_{++}})^{\prime}(f)}  - \leftexp{0}{ (I_{\spec}^{G_{++}})^{\prime}(f)}    \nonumber \\
&=&  -\frac{1}{4}  f^{M_{++}}(\psi \times \psi_{M_+}) - \frac{1}{2} f^{M_{++}}(\psi \times \psi_{M_+}) \nonumber \\
&=& -\frac{3}{4}  f^{M_{++}}(\psi \times \psi_{M_+}) . \nonumber
\end{eqnarray}

We now turn to $G^*=G_{++}^{\vee}$. Both the actual and expected value of the distribution $S^{G^{\vee} \times G_+^{\vee}}_{\disc,\psi_{++}^{N_{++}}}$ vanishes by the second equation of (9.2.7) and (9.2.15). Thus we have
\begin{eqnarray}
\leftexp{0}{ (I_{\scopy}^{G^{\vee}_{++}})^{\prime}} \equiv 0.
\end{eqnarray}
Similarly both the actual and expected value of $R^{M_{++}^{\vee}}_{\disc,\psi^N \times \psi^N}$ and hence (9.2.29) vanishes, by the second equation of (9.2.7) and (9.2.15). Thus we obtain:
\begin{eqnarray}
\leftexp{0}{ (I_{\spec}^{G^{\vee}_{++}})^{\prime}} \equiv 0.
\end{eqnarray}

It thus follows that
\begin{eqnarray*}
(9.2.18)&=& \widetilde{\iota}(N_{++},G_{++}) \cdot ( -\frac{3}{4}  f^{M_{++}}(\psi \times \psi_{M_+}))\\
&=& -\frac{3}{8}  f^{M_{++}}(\psi \times \psi_{M_+}). \nonumber
\end{eqnarray*}
Hence as before:
\begin{eqnarray*}
(9.2.17)&=&(9.2.18)+(9.2.19) \\
&=& -\frac{3}{8}  f^{M_{++}}(\psi \times \psi_{M_+}) -\frac{1}{8}  f^{M_{++}}(\psi \times \psi_{M_+}) \\
&=& -\frac{1}{2}  f^{M_{++}}(\psi \times \psi_{M_+}).
\end{eqnarray*}

\end{proof}

With proposition 9.2.2 established, we can complete the proof of the induction argument. Indeed, (9.2.16) is clearly in a form where lemma 4.3.6 applies, which asserts the vanishing of all coefficients. Since not all coefficients in (9.2.16) vanishes, this gives the desired final contradiction. In other words we conclude that $\delta_{\psi}=+1$ and $\Lambda$ vanishes (in both the generic and non-generic case), thus finishing the induction step for theorem 2.5.4(a) (which concerns only simple generic parameters) and the stable multiplicity formula for non-generic simple parameters. Finally, the spectral multiplicity formula, namely theorem 2.5.2, with respect to simple parameters follows by application of lemma 5.7.6.

We have thus completed all the local and global classification theorems. In particular, we highlight the following result which came from our analysis of simple generic parameters:

\begin{corollary}
Let $\phi^N \in \widetilde{\Phi}_{\simp}(N)$ be a simple generic parameter. Let $G=(G,\xi_{\chi}) \in \widetilde{\mathcal{E}}_{\simp}(N)$ with $\chi \in \mathcal{Z}_E^{\kappa}$. Then the following are equivalent:

\bigskip

\noindent a) $S^G_{\disc,\phi^N} \nequiv 0$.

\bigskip
\noindent b) $R^{G}_{\disc,\phi^N} \nequiv 0$, i.e. $\phi^N$ contributes to the discrete spectrum of $G$ (with respect to $\xi_{\chi}$).

\bigskip
\noindent c) The Asai $L$-function $L(s,\phi^N,\Asai^{(-1)^{N-1}\kappa})$ has a pole at $s=1$.

\end{corollary}

\bigskip

Finally given our local and global results, together with the generic descent theorems of Ginzburg-Rallis-Soudry \cite{GRS} for unitary groups, the argument of Proposition 8.3.2 of \cite{A1} on the local and global generic packet conjecture applies without change for quasi-split unitary groups. We thus state the following as the final result of the paper. For this, recall that if $F$ is either local or global, and $\psi_F$ is a non-trivial additive character of $F$ (when $F$ is local) or $\mathbf{A}/F$ (when $F$ is global), we denote by $(B,\lambda)$ the Whittaker datum associated to the standard splitting of $U_{E/F}(N)$ and the additive character $\psi_F$. Then we have:

\begin{corollary} (Generic packet conjecture)

\noindent a) Suppose that $F$ is local, and $\phi \in \Phi_{\bdd}(G)$ is a bounded generic parameter. Then the representation $\pi$ of the packet $\Pi_{\phi}$ correspinding to the trivial character of $\mathcal{S}_{\phi}$ is $(B,\lambda)$-generic.

\bigskip

\noindent b) Suppose that $F$ is global, and $\phi \in \Phi_2(G)$ is a square-integrable generic parameter. Then there exists a globally $(B,\lambda)$-generic cuspidal automorphic representation $\pi$ in the global packet $\Pi_{\phi}$ such that the corresponding character on $\mathcal{S}_{\phi}$ is trivial. 
\end{corollary}

\subsection{Appendix}

In this appendix we establish proposition 9.2.1. As in \cite{A1} this is based on formulas for characters of irreducible representations on general linear groups.

Thus for the most part we are in the local situation, with $E$ being a local field. We first review the character formulas for irreducible representations of $\GL_N(E)$, following mostly the notations of section 7.5 of \cite{A1}.

We first consider the case where $F$ is archimedean. For our purpose we only need to consider the case $E=\mathbf{C}$. We denote for $k \in \mathbf{Z}$ and $\lambda \in \mathbf{C}$ the character $\theta(k,\lambda)$ on $\mathbf{C}^{\times}$ given by:
\begin{eqnarray*}
z \mapsto (z/|z|)^k |z|^{\lambda}, \,\ z \in \mathbf{C}^{\times}.
\end{eqnarray*}
It is identified with a character $\mu(k,\lambda)$ of $W_{\mathbf{C}}=\mathbf{C}^{\times}$.

Denote $\psi^n(k) :=\mu(k,0) \boxtimes \nu^n$, for $k \in \mathbf{Z}$ and $n \geq 1$. Define $\theta^n(k)$ to be the character associated to the irreducible representation of $\GL_n(\mathbf{C})$, whose Langlands parameter is given by $\phi_{\psi^n(k)}$. Then the character $\theta^n(k)$ has an expansion in terms of standard characters, due to Tadi\'c \cite{Ta1} ({\it c.f.} p. 427-428 of \cite{A1} for notations):
\begin{eqnarray}
\theta^n(k) = \sum_{w \in S_n} \sgn(w) \cdot \theta^w(k)
\end{eqnarray}
\begin{eqnarray}
\theta^w(k) = \bigboxplus_{i=1}^n \theta \big( k -(i-wi) , (n+1)-(i+wi) \big).
\end{eqnarray}

Now consider the case that $E$ is non-archimedean. For $r$ an irreducible unitary representation of $W_E$, and $\lambda \in \mathbf{C}$, we denote by $r_{\lambda}$ the twist given as:
\[
r_{\lambda}(w) = r(w) \otimes |w|^{\lambda}, \,\ w \in W_{E}
\]
and for $k \geq 0$, put:
\[
\mu_r(k,\lambda) :=r_{\lambda} \boxtimes \nu^{k+1}
\]
which is an irreducible representation of $L_E=W_E \times \SU(2)$ of dimension equal to $m=m_r (k+1)$ (here $m_r=\dim r$). Denote by $\theta_r(k,\lambda)$ the character associated to the irreducible representation of $\GL_m(E)$ whose Langlands parameter is given by $\mu_r(k,\lambda)$.  

For $n \geq 1$, define:
\[
\psi_r^n(k) := \mu_r(k,0) \boxtimes \nu^n
\]
which is a parameter for $\GL_N(E)$, with $N=m \cdot n$.

Similar to the archimedean case, define $\theta_r^n(k)$ to be the character of the irreducible representation of $\GL_N(E)$, whose Langlands parameter is given by $\phi_{\psi^n_r(k)}$. Then we have the expansion of $\theta_r^n(k)$ in terms of standard characters, again due to Tadi\'c \cite{Ta1} ({\it c.f.} p.429 of \cite{A1}), which is formally the same as in the archimedean case:
\begin{eqnarray*}
\theta_r^n(k) = \sum_{w \in S_n} \sgn(w) \cdot \theta^w_r(k)
\end{eqnarray*}
\begin{eqnarray*}
\theta^w_r(k) = \bigboxplus_{i=1}^n \theta_r \big( k-(i-wi),(n+1)-(i+wi)   \big).
\end{eqnarray*}
Here the character $\theta(-1,\lambda)$ is interpreted as the trivial character on the trivial group, and $\theta(k,\lambda)$ is interpreted as the zero representation for $k < -1$. The notations can be made to be consistent with the archimedean case, if we simply interpret $r$ as the trivial representation of $W_{E}$ in the case where $E$ is archimedean.

Given the parameter $\psi_r^n(k)$, the following Weyl element $w^* \in S_n$ plays a special role. For $E$ archimedean $w^*$ is the longest Weyl element, i.e.
\[
w^*(i)=n+1-i.
\]
For $E$ non-archimedean, the element $w^*$ is defined as:
\begin{eqnarray*}
w^*(i) = \left \{ \begin{array}{c} n+1 - i \mbox{ if } i  \leq k+1 \\  i -(k+1) \mbox{ if } i \geq k+2
\end{array} \right.
\end{eqnarray*}
(thus $w^*$ is the longest Weyl element if $n \leq k+1$).

Corresponding to the Weyl element $w^*$ we set, for $E=\mathbf{C}$ (here $k \in \mathbf{Z}$):
\[
\theta^{n,*}_r(k) = \theta^{w^*}_r(k)= \bigboxplus_{i=1}^{n} \theta_r(k+n+1-2i,0).
\]
while for $E$ non-archimedean:
\begin{eqnarray*}
\theta^{n,*}_r(k) &=& \theta^{w^*}_r(k) = \bigboxplus_{i=1}^{n^*} \theta_r(k+n+1-2i,0),\\
n^* &=&  \min(n,k+1).
\end{eqnarray*}

We then have the following:
\begin{proposition} Here $k \geq 0$ if $E$ is non-archimedean, and $k \in \mathbf{Z}$ if $E=\mathbf{C}$.

\noindent (i) The character $\theta^*_r(k)$ for $k \geq 0$ occurs with multiplicity $\pm 1$ in the expansion of $\theta_r^n(k)$ into standard characters, and is the unique tempered character in the expansion.

\bigskip
\noindent (ii) The essentially square-integrable tempered character $\theta_r(k+n-1,0)$ occurs in the decomposition of $\theta_r^{n,*}(k)$ with multiplicity one.
\end{proposition}
\begin{proof}
This is Lemma 7.5.2 of \cite{A1}, except that the case $E=\mathbf{C}$ is not treated in {\it loc. cit.} However given the character formulas for $\GL_N(\mathbf{C})$ of Tadi\'c stated above, the same argument applies. 
\end{proof}

Proposition 9.3.1(i) can be stated as follows: for $\psi  \in \Psi(\GL_N(E))$, we have the general expansion:
\begin{eqnarray}
& & \\
& & f_N(\psi) = \sum_{\phi \in \Phi(N,\psi)} n(\psi,\phi) \cdot f_N(\phi), \,\ f \in \mathcal{H}(N)=\mathcal{H}(\GL_N(E)) \nonumber
\end{eqnarray}
({\it c.f.} p.437 of \cite{A1} for the definition of the set $\Phi(N,\psi) \subset \Phi(\GL_N(E))$). Then proposition 9.3.1(i) asserts that $n(\psi_r^n(k),\phi_r^{n,*}(k)) = \pm 1$; here we have denoted by $\phi_r^{n,*}(k)$ the bounded generic parameter of $\GL_N(E)$ corresponding to the tempered character $\theta_r^{n,*}(k)$.

We need the twisted version of (9.3.3), in the setting where $E/F$ is a quadratic extension. Thus for $\widetilde{\psi} \in \widetilde{\Psi}(\GL_N(E)) = \widetilde{\Psi}(N)$, we have the general expansion:
\begin{eqnarray}
\widetilde{f}_N(\widetilde{\psi}) = \sum_{\widetilde{\phi} \in \widetilde{\Phi}(N,\psi)} \widetilde{n} (\widetilde{\psi},\widetilde{\phi} ) \cdot \widetilde{f}_N(\widetilde{\phi}) , \,\ \widetilde{f} \in \widetilde{\mathcal{H}}(N).
\end{eqnarray}
(with $\widetilde{\Phi}(N,\widetilde{\psi})$ being the set of conjugate self-dual elements of $\Phi(N,\widetilde{\psi})$). Then the following proposition can be proved by exactly the same argument as in Corollary 7.5.5 of \cite{A1} (which is a consequence of Lemma 7.5.4 of {\it loc. cit.}):

\begin{proposition}
For $\widetilde{\psi} \in \widetilde{\Psi}(N)$ we have 
\[
\widetilde{n}(\widetilde{\psi},\widetilde{\phi}) \equiv n(\widetilde{\psi},\widetilde{\phi}) \bmod{2}.
\]
In particular for $\widetilde{\psi} =\psi_r^n(k)$ ($r$ being conjugate self-dual), we have 
\[
\widetilde{n}(\psi_r^n(k),\phi_r^{n,*}(k)) \equiv 1 \bmod{2}.
\]
\end{proposition}

\bigskip

After these preliminaries, we can finally prove Proposition 9.2.1. The argument is a variant of the proof of Lemma 8.2.1 of \cite{A1}.

Since the argument is local, we revert to the local notations used in section 7 and 8. Thus $F$ is a local field, and $E$ being a quadratic extension of $F$. This arises from a place $v$ of $\dot{F}$ that does not split in $\dot{E}$, and we take take $E=\dot{E}_v$.

We are given $\widetilde{\psi} \in \xi_* \psi$ with $\psi \in \Psi^+(G)$. Such a $\widetilde{\psi}$ arises in the global context as the localization at a non-split place $v$ of $\dot{F}$ of a non-generic simple parameter of even degree $N$. Hence we assume that $\widetilde{\psi}$ is of the form:
\begin{eqnarray*}
\widetilde{\psi}  = \mu \otimes \nu^n = \bigoplus_{i \in I} l_i \widetilde{\psi}_i
\end{eqnarray*}
\begin{eqnarray*}
\mu =  \bigoplus_{i \in I}  l_i \mu_i, \,\ \mu_i \in \Phi_{\simp}(m_i),  \mbox{ mutually distinct} 
\end{eqnarray*}
\begin{eqnarray*}
\widetilde{\psi}_i = \mu_i \otimes \nu^n
\end{eqnarray*}
\begin{eqnarray*}
 N &= &m \cdot n = \sum_{i \in I} l_i \cdot N_i \\
m &=& \sum_{i \in I} l_i m_i, \,\ N_i = m_i \cdot n.
\end{eqnarray*}

We are given the condition that the linear form $f^G(\psi)$ transfers to the Siegel Levi subgroup $L$ of $G$. We are going to show that the parameter $\widetilde{\psi}$ factors through the two $L$-embeddings 
\begin{eqnarray*}
\xi|_{\leftexp{L}{L}} \mbox{ and } \xi^{\vee}|_{\leftexp{L}{L}}.
\end{eqnarray*}
Here we take $\xi$ and $\xi^{\vee}$ to be the two $L$-embeddings $\leftexp{L}{G} \hookrightarrow \leftexp{L}{G_{E/F}(N)}$ associated to the two (representatives of equivalence classes of) simple twisted endoscopic datum $G=(U_{E/F}(N),\xi)$ and $G^{\vee}=(U_{E/F}(N),\xi^{\vee})$ in $\widetilde{\mathcal{E}}_{\simp}(N)$. Any parameter $\widetilde{\psi}$ that factors through $\xi|_{\leftexp{L}{L}}$ also factors through $\xi^{\vee}|_{\leftexp{L}{L}}$ (and conversely), and so it suffices to treat the case concerning $\xi|_{\leftexp{L}{L}}$.

Now since $\widetilde{\psi} \in \widetilde{\Psi}^+(N)$ is conjugate self-dual, there is an involution (similar to the discussion in section 2.4):
\begin{eqnarray*}
i \leftrightarrow i^{*}
\end{eqnarray*}
such that
\begin{eqnarray*}
\mu_{i^*} &=& (\mu_i)^* \\
l_i &=& l_{i^*}.
\end{eqnarray*}

If $\mu_i$ is not conjugate self-dual, then the sub-parameter of $\widetilde{\psi}$:
\begin{eqnarray*}
& & l_i \psi_i \oplus l_{i^*} \psi_{i^*} \\
&=& l_i (\psi \oplus (\psi_i)^*)
\end{eqnarray*}
factors through $\xi|_{\leftexp{L}{L}}$. For similar reason, if $\mu_i$ is conjugate self-dual, then the sub-parameter of $\widetilde{\psi}$:
\begin{eqnarray*}
2 l_i^{\prime} \psi_i &=& l_i^{\prime}(\psi_i \oplus (\psi_i)^*) \\
l_i^{\prime} &=& \mbox{ integer part of } l_i/2
\end{eqnarray*}
again factors through $\xi|_{\leftexp{L}{L}}$. Thus it suffices to analyze the set of indices:
\begin{eqnarray*}
I_-:=\{i  \in I \,\ |\,\  (\mu_i)^*=\mu_i, \,\ l_i \mbox{ odd} \}
\end{eqnarray*}
and the corresponding sub-parameter $\widetilde{\psi}_-$ of $\widetilde{\psi}$:
\begin{eqnarray*}
\widetilde{\psi}_- =  \bigoplus_{i \in I_-} \widetilde{\psi}_i = \mu_- \otimes \nu^n
\end{eqnarray*}
\begin{eqnarray*}
\mu_- = \bigoplus_{i \in I_-} \mu_i.
\end{eqnarray*}

To establish the proposition we must show that $I_-$ is empty. We are going to show a contradiction by assuming otherwise.

Thus assuming $I_-$ is non-empty. Put
\begin{eqnarray*}
m_- &=& \sum_{i \in I_-} m_i \\
N_- &=& m_- \cdot n
\end{eqnarray*}
(note that $N_- \equiv N \bmod{2}$ and in particular $N_-$ is again even). Then in particular we have $\widetilde{\psi}_- \in \widetilde{\Psi}_{\ellip}(N_-)$. Now the given condition on $\widetilde{\psi}$ implies that the twisted character:
\begin{eqnarray*}
\widetilde{f}_- \mapsto \widetilde{f}_{-,N_-}(\widetilde{\psi}_-), \,\ \widetilde{f}_- \in \widetilde{\mathcal{H}}(N_-) = \mathcal{H}(\widetilde{G}_{E/F}(N_-))
\end{eqnarray*}
transfers to a linear form on the Levi-subset
\[
\widetilde{L}_- = ( G_{E/F}(N_-/2) \times G_{E/F}(N_-/2) ) \rtimes \theta(N_-)
\]
of $\widetilde{G}_{E/F}(N_-)$. On the other hand, we have the general expansion as in (9.3.4) for the linear form $\widetilde{f}_{-,N_-}(\widetilde{\psi}_-)$:
\begin{eqnarray}
\widetilde{f}_{-,N_-}(\widetilde{\psi}_-) = \sum_{\widetilde{\phi}_- \in \widetilde{\Phi}(N_-,\widetilde{\psi}_-) } \widetilde{n}(\widetilde{\psi}_-,\widetilde{\phi}_-) \cdot \widetilde{f}_{-,N_-}(\widetilde{\phi}_-)
\end{eqnarray}
and by a germ expansion argument we have as a consequence that the linear form $\widetilde{f}_{-,N_-}(\widetilde{\phi}_-)$ transfers to a linear form on the Levi subset $\widetilde{L}_-$, for any $\widetilde{\phi}_-$ such that $\widetilde{n}(\widetilde{\psi}_-,\widetilde{\phi}_-) \neq 0$. In particular, this implies that such parameters $\widetilde{\phi}_-$ factor through: 
\begin{eqnarray}
\xi_-|_{\leftexp{L}{\widetilde{L}^0_-}} : \leftexp{L}{\widetilde{L}^0_-} \hookrightarrow \leftexp{L}{G_{E/F}(N_-)}
\end{eqnarray}
here we are denoting by $G_-=(G_-,\xi_-)$ the simple twisted endoscopic datum in $\widetilde{\mathcal{E}}_{\simp}(N_-)$ that has the same parity as $G=(G,\xi)$ ({\it c.f.} discussion on p.456 of \cite{A1}).

In particular by using proposition 9.3.2, this can be applied to the (unique) tempered $\widetilde{\phi}_-$ occuring in the expansion (9.3.5). More precisely, in the decomposition:
\[
\widetilde{\psi}_- = \bigoplus_{i \in I_-} \widetilde{\psi}_i, \,\ \widetilde{\psi}_i \in \widetilde{\Psi}_{\simp}(N_i)
\]
each of the simple sub-parameter $\widetilde{\psi}_i$ is of the form discussed in the beginning of the subsection. In other words, we can identify the index set $I_-$ as a set $J_-$ consisting of pairs $(r,k)$ (in both the archimedean and non-archimedean case), and write the decomposition of $\widetilde{\psi}_-$ as:
\begin{eqnarray}
\widetilde{\psi}_- = \bigoplus_{(r,k) \in J_-} \psi^n_r(k).
\end{eqnarray}

Then from the discussion on character formulas of Tadi\'c (namely Proposition 9.3.1(i) in the case where $E=\mathbf{C}$, and in the non-archimedean case respectively), the tempered component of $\widetilde{\psi}_-$ in the expansion (9.3.5) is given by:
\begin{eqnarray}
\widetilde{\phi}_-^* := \bigoplus_{(r,k) \in J_-} \phi^{n,*}_r(k).
\end{eqnarray}

To obtain a contradiction, we need to show that $\widetilde{\phi}^*_-$ does not factor through (9.3.6). For this it suffices to exhibit a sub-parameter of $\widetilde{\phi}_-^*$ that occurs with multiplicity one. Indeed, let 
\begin{eqnarray*}
k_1 := \max_{(r,k) \in J_-} k.
\end{eqnarray*}
Also choose $r_1$ such that $(r_1,k_1) \in J_-$ and having $\dim r_1$ being maximum (the choice of $r_1$ is not unique in the non-archimedean case, but it does not matter). Then by Proposition 9.3.1(ii), the essentially square-integrable character 
\[
\theta_{r_1}(k_1+n-1,0)
\]
occurs in the decomposition of $\theta_{r_1}^{n,*}(k_1)$ with multiplicity one; furthermore by construction $\theta_{r_1}(k_1+n-1,0)$ does not occur in $\theta_r^{n,*}(k)$ for any $(r,k) \in J_-$ with $(r,k) \neq (r_1,k_1)$. It follows that $\widetilde{\phi}_-^*$ has a sub-parameter (corresponding to the character $\theta_{r_1}(k_1+n-1,0)$) that occurs with multiplicity one, as required. This concludes the proof of proposition 9.2.1.

\section{\bf{Addendum}}

In this paper we have been mainly concerned with the case of quasi-split unitary groups. However, with the results already established in this paper we can establish some results in the non quasi-split case. The results of this appendix are used for example in the work of W.Zhang \cite{Z} on the refined Gan-Gross-Prasad conjecture for unitary groups. More complete results would of course be obtained, if one has an extension of the results of chapter nine of \cite{A1} to the setting of unitary groups.

We consider the global setup. Thus $E/F$ is a quadratic extension of global number fields. For any non-degenerate Hermitian space $V$ over $E$ (with respect to the extension $E/F$), we denote by $U(V)$ the unitary group over $F$ defined by the Hermitian space $V$. If $G$ is a unitary group then we denote by $G^*$ the quasi-split inner form of $G$. The $L$-group of $G$ depends of course only on $G^*$. For $\xi: \leftexp{L}{G} = \leftexp{L}{G^*} \hookrightarrow \leftexp{L}{G_{E/F}(N)}$, the pair $(G^*,\xi)$ defines an element in $\widetilde{\mathcal{E}}_{\simp}(N)$ (here $N=\dim_E V$). Note that $G^*$ is the principal elliptic endoscopic group of $G$.

Firstly the results of section 4 can be extended without change:
\begin{proposition}
For $G=U(V)$ with $\dim_E V=N$, $\xi : \leftexp{L}{G} = \leftexp{L}{G^*}\hookrightarrow \leftexp{L}{G_{E/F}(N)}$, and $c^N \in \widetilde{\mathcal{C}}_{\mathbf{A}}(N)$, we have
\begin{eqnarray*}
I^G_{\disc,c^N,t,\xi}(f) =0, \,\ f \in \mathcal{H}(G)
\end{eqnarray*}
unless $(c^N,t)=(c(\psi^N,t(\psi^N))$ for $\psi^N \in \widetilde{\Psi}(N)$. 
\end{proposition}
\begin{proof}
Same argument as in the proof of proposition 4.3.4. In fact the argument is even simpler in the present case, since the corresponding arguments concerning the stable distributions, which pertains only to quasi-split groups, are already established. 
\end{proof}
Similar to Corollary 4.3.8 we have the following consequence of Proposition 10.0.1:
\begin{corollary}
In the above notation we have
\[
L^2_{\disc,c^N,t,\xi}(G(F) \backslash G(\mathbf{A})) = 0
\]
unless $(c^N,t)=(c(\psi^N),t(\psi^N))$ for some $\psi^N \in \widetilde{\Psi}(N)$. We have a decomposition:
\begin{eqnarray*}
& & \\
& & L^2_{\disc}(G(F) \backslash G(\mathbf{A})) = \bigoplus_{\psi^N \in \widetilde{\Psi}(N)} L^2_{\disc,\psi^N,\xi} (G(F)  \backslash G(\mathbf{A})). \nonumber
\end{eqnarray*}
\end{corollary}
In particular if $\pi$ is a discrete automorphic representation of $G(\mathbf{A}_F)$, and $\psi^N \in \widetilde{\Psi}(N)$ such that $L^2_{\disc,\psi^N,\xi} \neq 0$, then as is customary $\psi^N$ is called a weak base change of $\pi$.

For $v$ a place of $F$ that splits in $E$, we make the identification as in the Notation part of section 1:
\[
G_v \cong \GL_{N/F_v}
\]
and for $\pi$ an irreducible admissible representation of $G(\mathbf{A})$ we identify $\pi_v$ accordingly as an irreducible admissible representation of $\GL_N(F_v)$. Then we have:

\begin{proposition}
Suppose $\pi$ occurs in $L^2_{\disc,\psi^N,\xi}(G(F)\backslash G(\mathbf{A}))$, with $\psi^N \in \widetilde{\Psi}_{\ellip}(N)$. Then for any place $v$ of $F$ that splits in $E$, we have $\pi_v \cong \pi_{\psi^N_v}$.
\end{proposition}
\begin{proof}
First since $\psi^N \in \widetilde{\Psi}_{\ellip}(N)$, the usual argument for the spectral expansion of $I^G_{\disc,\psi^N,\xi}$ shows that the terms coming from proper Levi subgroups do not contribute:
\begin{eqnarray*}
I^G_{\disc,\psi^N,\xi}(f) = \tr R^G_{\disc,\psi^N,\xi}(f), \,\ f \in \mathcal{H}(G).
\end{eqnarray*}
On the other hand we have the endoscopic expansion:
\begin{eqnarray*}
I^G_{\disc,\psi^N,\xi}(f) = \sum_{(G^{\prime},\zeta^{\prime}) \in \widetilde{\mathcal{E}}_{\ellip}(G)} \iota(G,G^{\prime}) \widehat{S}^{G^{\prime}}_{\disc,\psi^N,\xi \circ \zeta^{\prime}}(f^{G^{\prime}}).
\end{eqnarray*}

Choose decomposable $f=f_v f^v$. For $(G^{\prime},\zeta^{\prime}) \in \widetilde{\mathcal{E}}_{\ellip}(G)$, the group $G^{\prime}$ is a quasi-split unitary group or a product of two quasi-split unitary groups, and hence the results established in this paper applies to $G^{\prime}$. In particular, from the stable multiplicity formula given by Theorem 5.1.2, we see that for such $G^{\prime}$ and each choice of $f^v$, the distribution:
\[
f_v \rightarrow \widehat{S}^{G^{\prime}}_{\disc,\psi^N,\xi \circ \zeta^{\prime}}((f_v f^v)^{G^{\prime}})
\]
is proportional to $f_v(\pi_{\psi^N_v})$. Thus there is a scalar $c$, depending on the choice of $f^v$, such that:
\begin{eqnarray*} 
 \tr R^G_{\disc,\psi^N,\xi}(f_v f^v) = c \cdot f_v(\pi_{\psi^N_v}), \,\ f_v \in \mathcal{H}(G_v).
\end{eqnarray*}
Now since $\tr R^G_{\disc,\psi^N,\xi}$ (as a distribution on $G(\mathbf{A}_F)$) is a linear combination with non-negative coefficients of irreducible characters of representations of $G(\mathbf{A}_F)$. Hence by a standard linear independence of character argument, we see that if $\pi$ occurs in $L^2_{\disc,\psi^N,\xi}$, then $\pi_v \cong \pi_{\psi^N_v}$, as required.
\end{proof}

\end{document}